\documentclass[11pt,letterpaper]{article}
\usepackage[square,comma,authoryear,sort]{natbib}
\usepackage{hyphenat}
\usepackage[colorlinks,linktocpage,breaklinks]{hyperref}
\makeatletter
\def\@linkcolor{blue}
\def\@citecolor{blue}
\gdef\NAT@aysep{}
\let\@old@citep\citep
\let\@old@citet\citet
\let\@old@citeauthor\citeauthor
\def\citep{\@old@citep*}
\def\citet{\@old@citet*}
\def\citeauthor{\@old@citeauthor*}
\let\cite\citep
\defcitealias{Stacks}{Stacks 2013}
\defcitealias{AB:71}{Baernstein [1971]}
\makeatother
\usepackage[top=1.5in,bottom=1.2in,left=1.2in,right=1.2in,letterpaper]%
  {geometry}
\makeatletter
\def\ps@headings{%
    \let\@mkboth\@gobbletwo
    \def\@oddhead{\hss\scshape\shorttitle\hss\reset@font\rmfamily\thepage}
    \def\@evenhead{\reset@font\rmfamily\thepage\hss\scshape\shortauthors\hss}
    \let\@oddfoot\@empty\let\@evenfoot\@empty}
\@twosidetrue
\makeatother
\usepackage[neveradjust]{paralist}
\makeatletter
\@definecounter{@ADL@savenum}
\newcommand\savenum{\setcounter{@ADL@savenum}%
  {\the\@nameuse{c@\@listctr}}}
\newcommand\resumenum{\setcounter{\@listctr}{\arabic{@ADL@savenum}}}
\makeatother
\usepackage[reqno,centertags]{amsmath}
\usepackage{amssymb,xspace,twoopt,xy,mathrsfs}
\xyoption{arrow}
\xyoption{matrix}
\xyoption{curve}
\usepackage[mathcal]{eucal}
\DeclareFontFamily{U}{mathx}{\hyphenchar\font45}
\DeclareFontShape{U}{mathx}{m}{n}{
      <5> <6> <7> <8> <9> <10>
      <10.95> <12> <14.4> <17.28> <20.74> <24.88>
      mathx10
      }{}
\DeclareSymbolFont{mathx}{U}{mathx}{m}{n}
\DeclareFontSubstitution{U}{mathx}{m}{n}
\DeclareMathAccent{\widecheck}{0}{mathx}{"71}
\input cyracc.def
\let\subsetneq\subset
\let\subset\subseteq

\let\supset\supseteq
\newcommand\real{\mathbb{R}}
\newcommand\ereal{\overline{\real}}
\newcommand\erealnn{\ereal_{\ge0}}
\newcommand\realp{\real_{>0}}
\newcommand\realnn{\real_{\ge0}}
\newcommand\tdomain{\mathbb{T}}
\newcommand\complex{\mathbb{C}}
\newcommand\integer{\mathbb{Z}}
\newcommand\integerp{\integer_{>0}}
\newcommand\integernn{\integer_{\ge0}}

\newcommand\eul{\textup{e}}
\newcommand\imag{\textup{i}}
\newcommand\sphere{\mathbb{S}}
\newcommand\scirc{\raise1pt\hbox{$\,\scriptstyle\circ\,$}}
\newcommand\eqdef{\triangleq}
\newcommand\subscr[2]{#1_{\textup{#2}}}
\newcommand\supscr[2]{#1^{\textup{#2}}}
\newcommand\ol[1]{\overline{#1}}
\newcommand\slnorm{\lvert}
\newcommand\srnorm{\rvert}
\newcommand\snorm[1]{\slnorm #1\srnorm}
\newcommand\Bigsnorm[1]{\Bigl\lvert #1\Bigr\rvert}
\newcommand\dlnorm{\lVert}
\newcommand\drnorm{\rVert}
\newcommand\dnorm[1]{\dlnorm #1\drnorm}

\newcommand\setdef[2]{\{#1\;|\enspace#2\}}
\newcommand\asetdef[2]{\left\{\left.#1\;\right|\enspace#2\right\}}
\newcommand\bigsetdef[2]{\big\{#1\;\big|\enspace#2\big\}}
\newcommand\Bigsetdef[2]{\Big\{#1\;\Big|\enspace#2\Big\}}
\newcommand\ifam[1]{(#1)}
\newcommand\inprod[2]{\langle#1,#2\rangle}
\newcommand\natpair[2]{\langle#1;#2\rangle}
\makeatletter
\newcommand\vecspan{\@ifnextchar[{\@ADL@normalvecspan}
  {\@ADL@@normalvecspan}}
\def\@ADL@normalvecspan[#1]#2{\operatorname{span}_{#1}(#2)}
\newcommand\@ADL@@normalvecspan[1]{\operatorname{span}(#1)}
\makeatother
\newcommand\ie{i.e.,}
\newcommand\eg{e.g.,}
\newcommand\cf{cf.}
\newcommand\etc{etc.}
\newcommand\resp{resp.}
\renewcommand\ae{\textrm{a.e.}}
\newcommand\gcs{tautological control system}
\newcommand\Gcs{Tautological control system}
\newcommand\mathupper[1]{\textup{#1}}
\newcommand\End{\mathupper{End}}
\newcommand\Hom{\mathupper{Hom}}
\renewcommand\d[1]{{\normalfont\textrm{d}}#1}
\newcommand\image{\operatorname{image}}
\newcommand\graph{\operatorname{graph}}
\newcommand\supp{\operatorname{supp}}
\newcommand\cohull{\operatorname{conv}}
\newcommand\coker{\operatorname{coker}}
\newcommand\id{\operatorname{id}}
\newcommand\interior{\operatorname{int}}
\newcommand\closure{\operatorname{cl}}
\newcommand\ver{\operatorname{ver}}
\newcommand\hor{\operatorname{hor}}
\newcommand\hol{\textup{hol}}
\newcommand\lip{\textup{lip}}
\newcommand\card{\mathupper{card}}
\newcommand\CO{\mathupper{CO}}
\newcommand\dil[1]{\mathupper{dil}\,#1}
\newcommand\Dil[1]{\mathupper{Dil}\,#1}
\newcommand\ev{\mathupper{ev}}

\renewcommand\L{\mathupper{L}}
\newcommand\Lloc{\subscr{\L}{loc}}
\newcommand\pr{\operatorname{pr}}
\newcommand\Sym{\operatorname{Sym}}
\newcommand\Eq{\mathupper{Eq}}
\newcommand\Sh{\mathupper{Sh}}
\newcommand\Rtraj{\mathupper{Rtraj}}
\newcommand\Traj{\mathupper{Traj}}
\newcommand\hlft{\mathupper{hlft}}
\newcommand\vlft{\mathupper{vlft}}
\newcommand\sA{\mathscr{A}}
\newcommand\sB{\mathscr{B}}
\newcommand\sC{\mathscr{C}}
\newcommand\sE{\mathscr{E}}
\newcommand\sF{\mathscr{F}}
\newcommand\sG{\mathscr{G}}
\newcommand\sI{\mathscr{I}}
\newcommand\sL{\mathscr{L}}
\newcommand\sM{\mathscr{M}}
\newcommand\sN{\mathscr{N}}
\newcommand\sO{\mathscr{O}}

\newcommand\sS{\mathscr{S}}
\newcommand\sX{\mathscr{X}}
\newcommand\fG{\mathfrak{G}}
\newcommand\fH{\mathfrak{H}}
\let\dual\stardual 
\let\topdual\primedual
\newcommand\deriv[2]{\frac{{\normalfont\mathupper{d}}#1}
  {{\normalfont\mathupper{d}}#2}}
\newcommand\derivatzero[2]{\frac{{\normalfont\mathupper{d}}#1}
  {{\normalfont\mathupper{d}}#2}\Big|_{#2=0}}
\newcommand\tderiv[2]{{\textstyle\frac{{\normalfont\mathupper{d}}#1}
  {{\normalfont\mathupper{d}}#2}}}
\newcommand\tderivatzero[2]{{\textstyle\frac{{\normalfont\mathupper{d}}#1}
  {{\normalfont\mathupper{d}}#2}}\big|_{#2=0}}
\newcommand\pderiv[2]{\frac{\partial#1}{\partial#2}}
\makeatletter
\newcommand\linder{\@ifnextchar[{\@ADL@rlinder}{\@ADL@linder}}
\def\@ADL@rlinder[#1]#2{\boldsymbol{D}^{#1}#2}
\newcommand\@ADL@linder[1]{\boldsymbol{D}#1}
\newcommand\plinder{\@ifnextchar[{\@ADL@rplinder}{\@ADL@plinder}}
\def\@ADL@rplinder[#1]#2#3{\boldsymbol{D}^{#1}_{#2}#3}
\newcommand\@ADL@plinder[2]{\boldsymbol{D}_{#1}#2}
\makeatother
\newcommand\map[3]{#1\colon#2\rightarrow#3}
\newcommand\setmap[3]{#1\colon#2\twoheadrightarrow#3}
\newcommand\mapdef[5]{\begin{aligned}
  #1\colon&\begin{aligned}[t]#2\end{aligned}\rightarrow
  \begin{aligned}[t]#3\end{aligned}\\&\begin{aligned}[t]#4\end{aligned}
  \mapsto\begin{aligned}[t]#5\end{aligned}\end{aligned}}
\newcommand\mapschar{\mathupper{C}}
\newcommand\C{\mapschar}
\renewcommand\c{\mathupper{c}}
\newcommand\csd{\c_{\downarrow0}}

\newcommand\mat[1]{\boldsymbol{#1}}
\newcommand\vect[1]{\boldsymbol{#1}}
\newcommand\alg[1]{\mathsf{#1}}
\newcommand\field{\mathbb{F}}
\newcommand\metric{\mathbb{G}}
\newcommand\hmetric{\mathbb{H}}
\newcommand\dirlim{\qopname\relax n{dir\,lim}}
\newcommand\symmgroup[1]{\mathfrak{S}_{#1}}
\makeatletter
\newcommand\GL{\@ifnextchar[{\@ADL@GL}{\@ADL@@GL}}
\def\@ADL@GL[#1]#2{\Lie{GL}(#1;#2)}
\newcommand\@ADL@@GL[1]{\Lie{GL}(#1)}
\newcommand\Symalg{\@ifnextchar[{\@ADL@symalg}{\@ADL@@symalg}}
\def\@ADL@symalg[#1]#2{\mathupper{S}^{#1}(#2)}
\newcommand\@ADL@@symalg[1]{\mathupper{S}(#1)}
\newcommand\tensor{\@ifnextchar[{\@ADL@ktensors}{\@ADL@noktensors}}
\def\@ADL@ktensors[#1]#2{\mathupper{T}^{#1}(#2)}
\newcommand\@ADL@noktensors[1]{\mathupper{T}(#1)}
\newcommand\lin{\@ifnextchar[{\@ADL@klinmap}{\@ADL@linmap}}
\def\@ADL@klinmap[#1]#2#3{\def\@tempa{#3}%
  \ifx\@tempa\@empty\mathupper{L}^{#1}(#2;#2)\else%
  \mathupper{L}^{#1}(#2;#3)\fi}
\def\@ADL@linmap#1#2{\def\@tempa{#2}\ifx\@tempa\@empty\mathupper{L}(#1;#1)\else
  \mathupper{L}(#1;#2)\fi}
\newcommand\symlin{\@ifnextchar[{\@ADL@ksymlinmap}{\@ADL@symlinmap}}
\def\@ADL@ksymlinmap[#1]#2#3{\def\@tempa{#3}\ifx\@tempa\@empty
  \subscr{\mathupper{L}}{sym}^{#1}(#2;#2)\else%
  \subscr{\mathupper{L}}{sym}^{#1}(#2;#3)\fi}
\def\@ADL@symlinmap#1#2{\def\@tempa{#2}%
  \ifx\@tempa\@empty\subscr{\mathupper{L}}{sym}(#1;#1)\else
  \subscr{\mathupper{L}}{sym}(#1;#2)\fi}
\makeatother
\newcommand\cs[1]{\mathcal{#1}}
\newcommand\ts[1]{\mathcal{#1}}
\newcommand\nbhd[1]{\mathcal{#1}}

\newcommand\lebmes{\lambda}
\newcommand\charfunc[1]{\chi_{#1}}
\newcommand\injotimes{\widecheck{\otimes}_e}
\newcommand\protimes{\widehat{\otimes}_\pi}
\makeatletter
\newcommand\oball{\@ifnextchar[{\@ADL@oballarg}{\@ADL@oballnoarg}}
\def\@ADL@oballarg[#1]#2#3{\mathsf{B}_{#1}(#2,#3)}
\newcommand\@ADL@oballnoarg[2]{\mathsf{B}(#1,#2)}
\newcommand\cball{\@ifnextchar[{\@ADL@cballarg}{\@ADL@cballnoarg}}
\def\@ADL@cballarg[#1]#2#3{\ol{\mathsf{B}}_{#1}(#2,#3)}
\newcommand\@ADL@cballnoarg[2]{\ol{\mathsf{B}}(#1,#2)}
\newcommand\odisk{\@ifnextchar[{\@ADL@odiskdim}{\@ADL@odisknodim}}
\def\@ADL@odiskdim[#1]#2#3{\mathsf{D}^{#1}(#2,#3)}
\newcommand\@ADL@odisknodim[2]{\mathsf{D}^n(#1,#2)}
\newcommand\cdisk{\@ifnextchar[{\@ADL@cdiskdim}{\@ADL@cdisknodim}}
\def\@ADL@cdiskdim[#1]#2#3{\ol{\mathsf{D}}\null^{#1}(#2,#3)}
\newcommand\@ADL@cdisknodim[2]{\ol{\mathsf{D}}\null^n(#1,#2)}
\makeatother
\newcommandtwoopt\sections[3][\sinfty][\null]{\Gamma^{#1}_{#2}(#3)}
\newcommand\secbdd[2][\infty]{\sections[#1][\textup{bdd}]{#2}}
\newcommand\Lie[1]{\mathsf{#1}}
\newcommand\man[1]{\mathsf{#1}}
\newcommand\dist[1]{\mathsf{#1}}
\newcommand\hlift[1]{#1^{H}}
\newcommand\tlift[1]{#1^{T}}
\newcommand\vlift[1]{#1^{V}}
\makeatletter
\newcommand\tb{\@ifnextchar[{\@ADL@tbarg}{\@ADL@tb}}
\def\@ADL@tbarg[#1]#2{\man{T}_{#1}#2}
\newcommand\@ADL@tb[1]{\man{T}#1}
\newcommand\tbproj[1]{\pi_{\tb{#1}}}
\newcommand\ctb{\@ifnextchar[{\@ADL@ctbarg}{\@ADL@ctb}}
\def\@ADL@ctbarg[#1]#2{\man{T}^*_{#1}#2}
\newcommand\@ADL@ctb[1]{\man{T}^*#1}
\newcommand\ctbproj[1]{\pi_{\ctb{#1}}}
\newcommand\vb{\@ifnextchar[{\@ADL@vbarg}{\@ADL@vb}}
\def\@ADL@vbarg[#1]#2{\man{V}_{#1}#2}
\newcommand\@ADL@vb[1]{\man{V}#1}
\newcommand\func{\@ifnextchar[{\@ADL@crfuncs}{\@ADL@cinftyfuncs}}
\def\@ADL@crfuncs[#1]#2{\mapschar^{#1}(#2)}
\newcommand\@ADL@cinftyfuncs[1]{\mapschar^\infty(#1)}
\newcommand\mappings{\@ifnextchar[{\@ADL@crmappings}{\@ADL@cinftymappings}}
\def\@ADL@crmappings[#1]#2#3{\mapschar^{#1}(#2;#3)}
\newcommand\@ADL@cinftymappings[2]{\mapschar^\infty(#1;#2)}
\newcommand\gfunc{\@ifnextchar[{\@ADL@crgfuncs}{\@ADL@cinftygfuncs}}
\def\@ADL@crgfuncs[#1]#2#3{\mathscr{C}^{#1}_{#2,#3}}
\newcommand\@ADL@cinftygfuncs[2]{\mathscr{C}^\infty_{#1,#2}}
\newcommand\sfunc{\@ifnextchar[{\@ADL@crsfuncs}{\@ADL@cinftysfuncs}}
\def\@ADL@crsfuncs[#1]#2{\mathscr{C}^{#1}_{#2}}
\newcommand\@ADL@cinftysfuncs[1]{\mathscr{C}^\infty_{#1}}
\newcommand\gsections{\@ifnextchar[{\@ADL@rgsections}{\@ADL@gsections}}
\def\@ADL@rgsections[#1]#2#3{\mathscr{G}^{#1}_{#2,#3}}
\newcommand\@ADL@gsections[2]{\mathscr{G}^\infty_{#1,#2}}
\newcommand\ssections{\@ifnextchar[{\@ADL@rssections}{\@ADL@ssections}}
\def\@ADL@rssections[#1]#2{\mathscr{G}^{#1}_{#2}}
\newcommand\@ADL@ssections[1]{\mathscr{G}^\infty_{#1}}
\newcommand\tf{\@ifnextchar[{\@ADL@tfarg}{\@ADL@tf}}
\def\@ADL@tfarg[#1]#2{\@ifnextchar[{\@ADL@@tfargr{#1}{#2}}
  {\@ADL@@tfarg{#1}{#2}}}
\def\@ADL@@tfargr#1#2[#3]{T^{#3}_{#1}#2}
\newcommand\@ADL@@tfarg[2]{T_{#1}#2}
\newcommand\@ADL@tf[1]{\@ifnextchar[{\@ADL@@tfr{#1}}{\@ADL@@tf{#1}}}
\def\@ADL@@tfr#1[#2]{T^{#2}#1}
\newcommand\@ADL@@tf[1]{T#1}
\newcommand\lieder[2]{\def\@tempa{#2}\ifx\@tempa\@empty%
  \boldsymbol{\mathscr{L}}_{#1}\else\boldsymbol{\mathscr{L}}_{#1}#2\fi}
\newcommand\SPC{\@ifnextchar[{\@SPCbrack}{\@SPCnobrack}}
\def\@SPCbrack[#1]#2{\mathupper{SPC}^{#1}(#2)}
\def\@SPCnobrack#1{\mathupper{SPC}^{\infty}(#1)}
\newcommand\SPsections{\@ifnextchar[{\@SPCsecbrack}{\@SPCsecnobrack}}
\def\@SPCsecbrack[#1]#2{\mathupper{SP}\Gamma^{#1}(#2)}
\def\@SPCsecnobrack#1{\mathupper{SP}\Gamma^{\infty}(#1)}
\newcommand\JPC{\@ifnextchar[{\@JPCbrack}{\@JPCnobrack}}
\def\@JPCbrack[#1]#2{\mathupper{JPC}^{#1}(#2)}
\def\@JPCnobrack#1{\mathupper{JPC}^{\infty}(#1)}
\newcommand\JPsections{\@ifnextchar[{\@JPCsecbrack}{\@JPCsecnobrack}}
\def\@JPCsecbrack[#1]#2{\mathupper{JP}\Gamma^{#1}(#2)}
\def\@JPCsecnobrack#1{\mathupper{JP}\Gamma^{\infty}(#1)}
\newcommand\LIC{\@ifnextchar[{\@LICbrack}{\@LICnobrack}}
\def\@LICbrack[#1]#2{\mathupper{LIC}^{#1}(#2)}
\def\@LICnobrack#1{\mathupper{LIC}^{\infty}(#1)}
\newcommand\LIsections{\@ifnextchar[{\@LICsecbrack}{\@LICsecnobrack}}
\def\@LICsecbrack[#1]#2{\mathupper{LI}\Gamma^{#1}(#2)}
\def\@LICsecnobrack#1{\mathupper{LI}\Gamma^{\infty}(#1)}
\newcommand\LBC{\@ifnextchar[{\@LBCbrack}{\@LBCnobrack}}
\def\@LBCbrack[#1]#2{\mathupper{LBC}^{#1}(#2)}
\def\@LBCnobrack#1{\mathupper{LBC}^{\infty}(#1)}
\newcommand\LBsections{\@ifnextchar[{\@LBCsecbrack}{\@LBCsecnobrack}}
\def\@LBCsecbrack[#1]#2{\mathupper{LB}\Gamma^{#1}(#2)}
\def\@LBCsecnobrack#1{\mathupper{LB}\Gamma^{\infty}(#1)}
\newcommand\CF{\@ifnextchar[{\@CFbrack}{\@CFnobrack}}
\def\@CFbrack[#1]#2{\mathupper{CF}^{#1}(#2)}
\def\@CFnobrack#1{\mathupper{CF}^{\infty}(#1)}
\newcommand\CFsections{\@ifnextchar[{\@CFsectionsbrack}{\@CFsectionsnobrack}}
\def\@CFsectionsbrack[#1]#2{\mathupper{CF}\Gamma^{#1}(#2)}
\def\@CFsectionsnobrack#1{\mathupper{CF}\Gamma^{\infty}(#1)}
\def\flow{\@ifstar{\@flow@star}{\@flow@nostar}}
\def\@flow@star#1{\Phi^{#1}}
\def\@flow@nostar#1#2{\@ifnextchar[{\@tflow{#1}{#2}}{\@flow{#1}{#2}}}
\def\@tflow#1#2[#3]{\Phi^{#1}_{#2,#3}}
\def\@flow#1#2{\Phi^{#1}_{#2}}
\def\jet{\@ifnextchar[{\@ADL@jetbase}{\@ADL@jetnobase}}
\def\@ADL@jetbase[#1]#2#3{\man{J}_{#1}^{#2}#3}
\def\@ADL@jetnobase#1#2{\man{J}^{#1}#2}
\def\jetalg{\@ifnextchar[{\@ADL@jetalgbase}{\@ADL@jetalgnobase}}
\def\@ADL@jetalgbase[#1]#2#3{\man{T}^{*#2}_{#1}#3}
\def\@ADL@jetalgnobase#1#2{\man{T}^{*#1}#2}
\makeatother
\makeatletter
\def\interval{\@ifnextchar({\@ADL@openleftint}{\@ADL@closedleftint}}
\def\@ADL@openleftint(#1,#2{(#1,#2%
  \@ifnextchar){\@ADL@openrightint}{\@ADL@closedrightint}}
\def\@ADL@closedleftint[#1,#2{[#1,#2%
  \@ifnextchar){\@ADL@openrightint}{\@ADL@closedrightint}}
\def\@ADL@openrightint){)}
\def\@ADL@closedrightint]{]}
\makeatother
\usepackage{theorem}
\newtheorem{theorem}{Theorem}[section]
\newtheorem{proposition}[theorem]{Proposition}
\newtheorem{lemma}[theorem]{Lemma}
\newtheorem{corollary}[theorem]{Corollary}
\newtheorem{prooflemma}{Lemma}[theorem]
\newtheorem{proofsublemma}{Sublemma}[theorem]

{\theorembodyfont{\normalfont\rmfamily}
\newtheorem{definition}[theorem]{Definition}
\newtheorem{remark}[theorem]{Remark}
\newtheorem{remarks}[theorem]{Remarks}
\newtheorem{example}[theorem]{Example}
\newtheorem{examples}[theorem]{Examples}
\newtheorem{notation}[theorem]{Notation}}
\makeatletter
\newcommand\pushright{\protect\@ADL@pushright}
\newcommand\@ADL@pushright[1]{{\ifvmode\null\hfill{#1}\par\else\ifmmode%
  \@ADLmaths@pushright{\hbox{#1}}\else\ifinner\@ADLhbox@pushright{#1}%
  \else\@ADLparag@pushright{#1}\fi\fi\fi}}
\newcommand\@ADLmaths@pushright[1]{{\ifinner\@ADLhbox@pushright{#1}\else%
  \tag*{$#1$}\fi}}
\newcommand\@ADLparag@pushright[1]{{\parfillskip=0pt\widowpenalty=10000%
  \displaywidowpenalty=10000\finalhyphendemerits=0\@ADLhbox@pushright#1\par}}
\newcommand\@ADLhbox@pushright{\unskip\nobreak\hfil\penalty50\hskip.2em%
  \null\hfill\hfill}
\newenvironment{proof}{\trivlist\item[\hskip\labelsep\textit{Proof:}\/]%
  \@ADLsave@set@qed\xspace\normalfont\rmfamily}
  {\qed\@ADLrestore@qed\endtrivlist}
\newif\if@ADL@qed\@ADL@qedfalse
\newcommand\qed{\protect\@ADL@qed{$\blacksquare$}}
\newcommand\@ADL@qed[1]{\if@ADL@qed\global\@ADL@qedfalse%
  \pushright{#1}\else\ifhmode\ifinner\else\par\fi\fi\fi}
\newcommand\@ADLrestore@qed{\global\let\if@ADL@qed\@ADLsaved@ifqed}
\newcommand\@ADLsave@set@qed{\let\@ADLsaved@ifqed
  \if@ADL@qed\global\@ADL@qedtrue}
\newenvironment{subproof}{\trivlist\item[\hskip\labelsep\textit{Proof:}\/]%
  \@ADLsave@set@subqed\normalfont\rmfamily}
  {\subqed\@ADLrestore@subqed\endtrivlist}
\newif\if@ADL@subqed\@ADL@subqedfalse
\newcommand\subqed{\protect\@ADL@subqed{$\blacktriangledown$}}
\newcommand\@ADL@subqed[1]{\if@ADL@subqed\global\@ADL@subqedfalse%
  \pushright{#1}\else\ifhmode\ifinner\else\par\fi\fi\fi}
\newcommand\@ADLrestore@subqed{\global\let\if@ADL@subqed\@ADLsaved@ifsubqed}
\newcommand\@ADLsave@set@subqed{\let\@ADLsaved@ifsubqed
  \if@ADL@subqed\global\@ADL@subqedtrue}
\newcommand\eqsubqed{\tag*{\subqed}}
\newif\if@ADL@oprocend\@ADL@oprocendfalse
\newcommand\oprocend{\@ADLsave@set@oprocend
  \protect\@ADL@oprocend{$\bullet$}\@ADLrestore@oprocend}
\newcommand\@ADL@oprocend[1]{\if@ADL@oprocend\global\@ADL@oprocendfalse%
  \pushright{#1}\else\ifhmode\ifinner\else\par\fi\fi\fi}
\newcommand\@ADLrestore@oprocend{\global
  \let\if@ADL@oprocend\@ADLsaved@ifoprocend}
\newcommand\@ADLsave@set@oprocend{\let\@ADLsaved@ifoprocend\if@ADL@oprocend%
  \global\@ADL@oprocendtrue}
\newcommand\eqoprocend{\tag*{\oprocend}}
\makeatother
\newenvironment{keywords}{\quote\small\textbf{Keywords.}}{\endquote}
\newenvironment{AMS}{\quote\small\textbf{AMS Subject Classifications (2010).}}
   {\endquote}
\newcommand\defn[1]{{\normalfont\bfseries\emph{\mathversion{bold}#1}}}
\makeatletter
\@namedef{procnonum*}{\trivlist \@topsep \theorempreskipamount
  \@topsepadd \theorempostskipamount
  \@ifnextchar[{\@ADL@xprocnonumstar}{\@ADL@yprocnonumstar}}
\@namedef{endprocnonum*}{\endtrivlist}
\def\@ADL@xprocnonumstar[#1]{\item[\hskip \labelsep{\theorem@headerfont #1}]
  \normalfont\rmfamily}
\def\@ADL@yprocnonumstar{\item[] \normalfont\rmfamily}
\makeatother
\numberwithin{equation}{section}
\newcommand\pldblref[2]{\mbox{\ref{#1}(\ref{#2})}}
\newcommand\enumdblref[2]{\mbox{\ref{#1}--\ref{#2}}}
\title{Mathematical models for geometric control theory\thanks{Research
supported in part by a grant from the Natural Sciences and Engineering
Research Council of Canada}}
\author{Saber Jafarpour\thanks{Graduate student, Department of Mathematics
and Statistics, Queen's University, Kingston, ON K6L 3N6, Canada,
email:~\texttt{saber.jafarpour@queensu.ca}} \and Andrew D.\
Lewis\thanks{Professor, Department of Mathematics and Statistics, Queen's
University, Kingston, ON K7L 3N6, Canada, email:~\texttt{andrew@mast.queensu.ca}}}
\date{12/06/2014}
\newcommand\shorttitle{Mathematical models for geometric control theory}
\newcommand\shortauthors{S.\ Jafarpour and A.\ D.\ Lewis}
\begin{document}
\maketitle

\begin{abstract}
Just as an explicit parameterisation of system dynamics by state, i.e., a
choice of coordinates, can impede the identification of general structure, so
it is too with an explicit parameterisation of system dynamics by control.
However, such explicit and fixed parameterisation by control is commonplace
in control theory, leading to definitions, methodologies, and results that
depend in unexpected ways on control parameterisation.  In this paper a
framework is presented for modelling systems in geometric control theory in a
manner that does not make any choice of parameterisation by control; the
systems are called ``tautological control systems.''  For the framework to be
coherent, it relies in a fundamental way on topologies for spaces of vector
fields.  As such, classes of systems are considered possessing a variety of
degrees of regularity: finitely differentiable; Lipschitz; smooth; real
analytic.  In each case, explicit geometric seminorms are provided for the
topologies of spaces of vector fields that enable straightforward
descriptions of time-varying vector fields and control systems.  As part of
the development, theorems are proved for regular (including real analytic)
dependence on initial conditions of flows of vector fields depending
measurably on time.  Classes of ``ordinary'' control systems are
characterised that interact with the regularity under consideration in a
comprehensive way.  In this framework, for example, the statement that ``a
smooth or real analytic control-affine system is a smooth or real analytic
control system'' becomes a theorem.  Correspondences between ordinary control
systems and tautological control systems are carefully examined, and
trajectory correspondence between the two classes is proved for
control-affine systems and for systems with general control dependence when
the control set is compact.
\end{abstract}


\tableofcontents

\section{Introduction}

One can study nonlinear control theory from the point of view of
applications, or from a more fundamental point of view, where system
structure is a key element.  From the practical point of view, questions that
arise are often of the form, ``How can we\ldots'', for example, ``How can we
steer a system from point $A$ to point $B$?''\ or, ``How can we stabilise
this unstable equilibrium point?''\ or, ``How can we manoeuvre this vehicle
in the most efficient manner?''  From a fundamental point of view, the
problems are often of a more existential nature, with, ``How can we''
replaced with, ``Can we''.  These existential questions are often very
difficult to answer in any sort of generality.

As one thinks about these fundamental existential questions and looks into
the quite extensive existing literature, one comes to understand that the
question, ``What is a control system?''\ is one whose answer must be decided
upon with some care.  One also begins to understand that structure coming
from common physical models can be an impediment to general understanding.
For example, in a real physical model, states are typically physical
quantities of interest,~\eg~position, current, quantity of reactant~X, and so
the explicit labelling of these is natural.  This labelling amounts to a
specific choice of coordinates, and it is now well understood that such
specific choices of coordinates obfuscate structure, and so are to be avoided
in any general treatment.  In like manner, in a real physical model, controls
are likely to have meaning that one would like to keep track of,~\eg~force,
voltage, flow.  The maintenance of these labels in a model provides a
specific parameterisation of the inputs to the system, completely akin to
providing a specific coordinate parameterisation for states.  However, while
specific coordinate parameterisations have come (by many) to be understood as
a bad idea in a general treatment, this is not the case for specific control
parameterisations; models with fixed control parameterisation are commonplace
in control theory.  In contrast to the situation with dependence of
\emph{state} on parameterisation, the problem of eliminating dependence of
\emph{control} on parameterisation is not straightforward.  In our discussion
below we shall overview some of the common models for control systems, and
some ways within these modelling frameworks for overcoming the problem of
dependence on control parameterisation.  As we shall see, the common models
all have some disadvantage or other that must be confronted when using these
models.  In this paper we provide a means for eliminating explicit
parameterisation of controls that, we believe, overcomes the problems with
existing techniques.  Our idea has some of its origins in the work on
``chronological calculus'' of \citet{AAA/RVG:78} (see
also~\cite{AAA/YS:04})\@, but the approach we describe here is more general
(in ways that we will describe below) and more fully developed as concerns
its relationship to control theory (chronological calculus is primarily a
device for understanding time-varying vector fields and flows).  There are
some ideas similar to ours in the approach of \citet{HJS:97b}\@, but there
are also some important differences,~\eg~our families of vector fields are
time-invariant (corresponding to vector fields with frozen control values)
while \citeauthor{HJS:97b} considers families of time-varying vector fields
(corresponding to selecting an open-loop control).  Also, the work of
\citeauthor{HJS:97b} does not touch on real analytic systems.

We are interested in models described by ordinary differential equations
whose states are in a finite-dimensional manifold.  Even within this quite
narrow class of control systems, there is a lot of room to vary the models
one might consider.  Let us now give a brief outline of the sorts of models
and methodologies of this type that are commonly present in the literature.

\subsection{Models for geometric control systems: pros and cons}\label{subsec:modelcomp}

By this time, it is well-understood that the language of systems such as we
are considering should be founded in differential geometry and vector fields
on manifolds~\cite{AAA/YS:04,AMB:02,FB/ADL:04,AI:95,VJ:97,HN/AJV:90}\@.  This
general principle can go in many directions, so let us discuss a few of
these.  Our presentation here is quite vague and not very careful.  In the
main body of the paper, we will be less vague and more careful.

\subsubsection{Family of vector field models}

Given that manifolds and vector fields are important, a first idea of what
might comprise a control system is that it is a family of vector fields.  For
these models, trajectories are concatenations of integral curves of vector
fields from the family.  This is the model used in the development of the
theory of accessibility of \citet{HJS/VJ:72} and in the early work of
\citet{HJS:78} on local controllability.  The work of \citet{RH/AJK:77}\@,
while taking place in the setting of systems parameterised by control (such
as we shall discuss in Section~\ref{subsubsec:parameterised}), uses the
machinery of families of vector fields to study controllability and
observability of nonlinear systems.  Indeed, a good deal of the early work in
control theory is developed in this sort of framework, and it is more or less
sufficient when dealing with questions where piecewise constant controls are
ample enough to handle the problems of interest.  The theory is also highly
satisfying in that it is very differential geometric, and the work utilising
this approach is often characterised by a certain elegance.

However, the approach does have the drawback of not handling well some of the
more important problems of control theory, such as feedback (where controls
are specified as functions of state) and optimal control (where piecewise
constant controls are often not a sufficiently rich
class~\cite[\cf][]{ATF:60}).

It is worth mentioning at this early stage in our presentation that one of
the ingredients of our approach is a sort of fusion of the ``family of vector
fields'' approach with the more common control parameterisation approach to
whose description we now turn.

\subsubsection{Models with control as a parameter}\label{subsubsec:parameterised}

Given the limitations of the ``family of vector fields'' models for physical
applications and also for a theory where merely measurable controls are
needed, one feels as if one has to have the control as a parameter in the
model, a parameter that one can vary in a quite general manner.  These sorts
of models are typically described by differential equations of the form
\begin{equation*}
\dot{x}(t)=F(x(t),u(t)),
\end{equation*}
where $t\mapsto u(t)$ is the control and $t\mapsto x(t)$ is a corresponding
trajectory.  For us, the trajectory is a curve on a differentiable manifold
$\man{M}$\@, but there can be some freedom in attributing properties to the
control set $\cs{C}$ in which $u$ takes its values, and on the properties of
the system dynamics $F$\@.  (In Section~\ref{sec:systems} we describe classes
of such models in differential geometric terms.)  This sort of model is
virtually synonymous with ``nonlinear control system'' in the existing
control literature.  A common class of systems that are studied are
control-affine systems, where
\begin{equation*}
F(x,\vect{u})=f_0(x)+\sum_{a=1}^ku^af_a(x),
\end{equation*}
for vector fields $f_0,f_1,\dots,f_k$ on $\man{M}$\@, and where the control
$\vect{u}$ takes values in a subset of $\real^k$\@.  For control-affine
systems, there is an extensively developed theory of controllability based on
free Lie algebras~\cite{RMB/GS:93,MK:90b,MK:99,MK:06,HJS:83a,HJS:87}\@.  We
will see in Section~\ref{subsec:Cr-systemsII} that control-affine systems fit
into our framework in a particularly satisfying way.

The above general model, and in particular the control-affine special case,
are all examples where there is an explicit parameterisation of the control
set,~\ie~the control $u$ lives in a particular set and the dynamics $F$ is
determined to depend on $u$ in some particular way.  It could certainly be
the case, for instance, that one could have two different systems
\begin{equation*}
\dot{x}(t)=F_1(x(t),u_1(t)),\quad\dot{x}(t)=F_2(x(t),u_2(t))
\end{equation*}
with exactly the same trajectories.  This has led to an understanding that
one should study equivalence classes of systems.  A little precisely, if one
has two systems
\begin{equation*}
\dot{x}_1(t)=F_1(x_1(t),u_1(t)),\quad\dot{x}_2(t)=F_2(x_2(t),u_2(t)),
\end{equation*}
with $x_a(t)\in\man{M}_a$ and $u_a(t)\in\cs{C}_a$\@, $a\in\{1,2\}$\@, then
there may exist a diffeomorphism $\map{\Phi}{\man{M}_1}{\man{M}_2}$ and a
mapping $\map{\kappa}{\man{M}_1\times\cs{C}_1}{\cs{C}_2}$ (with some sort of
regularity that we will not bother to mention) such that
\begin{compactenum}
\item $\tf[x_1]{\Phi}\scirc F_1(x_1,u_1)=F_2(\Phi(x_1),\kappa(x_1,u_1))$ and
\item the trajectories $t\mapsto x_1(t)$ for the first system are in 1--1
correspondence with those of the second system by $t\mapsto\Phi\scirc
x_1(t)$\@.\footnote{We understand that there are many ways of formulating
system equivalence.  But here we are content to be, not only vague, but far
from comprehensive.}
\end{compactenum}
Let us say a few words about this sort of ``feedback equivalence.''  One can imagine it being useful in at least two ways.
\begin{compactenum}
\item First of all, one might use it as a kind of ``acid test'' on the
viability of a control theoretic construction.  That is, a control theoretic
construction should make sense, not just for a system, but for the
equivalence class of that system.  This is somewhat akin to asking that
constructions in differential geometry should be independent of coordinates.
Indeed, in older presentations of differential geometry, this was often how
constructions were defined: they were given in coordinates, and then
demonstrated to behave properly under changes of coordinate.  We shall
illustrate in Example~\ref{eg:bad-linearise} below that many common
constructions in control theory do not pass the ``acid test'' for viability
as feedback-invariant constructions.
\item Feedback equivalence is also a device for classifying control systems, the prototypical example being ``feedback linearisation,'' the determination of those systems that are linear systems in disguise~\cite{BJ/WR:80}\@.  In differential geometry, this is akin to the classification of geometric structures on manifolds,~\eg~Riemannian, symplectic, \etc
\end{compactenum}
In Section~\ref{subsec:gcs-category} we shall consider a natural notion of
equivalence for systems of the sort we are introducing in this paper, and we
will show that ``feedback transformations'' are vacuous in that they amount
to being described by mappings between manifolds.  This is good news, since
the whole point of our framework is to eliminate control parameterisation
from the picture and so eliminate the need for considering the effects of
varying this parameterisation,~\cf~``coordinate-free'' versus
``coordinate-independent'' in differential geometry.  Thus the first of the
preceding uses of feedback transformations simply does not come up for
us:~our framework is naturally feedback-invariant.  The second use of
feedback transformations, as will be seen in
Section~\ref{subsec:gcs-category}\@, amounts to the classification of
families of vector fields under push-forward by diffeomorphisms.  This is
generally a completely hopeless undertaking, so we will have nothing to say
about this.  Studying this under severe restrictions using, for
example,~(1)~the Cartan method of
equivalence~\cite[\eg][]{RLB/RBG:93,RBG:89}\@,~(2)~the method of generalised
transformations~\cite[\eg][]{WK/AJK:98,WK/AJK:06}\@,~(3)~the study of
singularities of vector fields and
distributions~\cite[\eg][]{BJ/WR:80,WPL/WR:02}\@, one might expect that some
results are possible.

Let us consider an example that shows how a classical control-theoretic construction, linearisation, is not invariant under even the very weak notion of equivalence where equivalent systems are those with the same trajectories.
\begin{example}\label{eg:bad-linearise}
We consider two control-affine systems
\begin{equation*}
\begin{aligned}
\dot{x}_1(t)=&\;x_2(t),\\
\dot{x}_2(t)=&\;x_3(t)u_1(t),\\
\dot{x}_3(t)=&\;u_2(t),
\end{aligned}\qquad\quad
\begin{aligned}
\dot{x}_1(t)=&\;x_2(t),\\
\dot{x}_2(t)=&\;x_3(t)+x_3(t)u_1(t),\\
\dot{x}_3(t)=&\;u_2(t),
\end{aligned}
\end{equation*}
with $(x_1,x_2,x_3)\in\real^3$ and $(u_1,u_2)\in\real^2$\@.  One can readily verify that these two systems have the same trajectories.  If we linearise these two systems about the equilibrium point at $(0,0,0)$\textemdash{}in the usual sense of taking Jacobians with respect to state and control~\cite[page~172]{AI:95}\@,~\cite[\S12.2]{HKK:96}\@,%
~\cite[Proposition~3.3]{HN/AJV:90}\@,~\cite[page~236]{SS:99}\@, and~\cite[Definition~2.7.14]{EDS:98}\textemdash{}then we get the two linear systems
\begin{equation*}
\mat{A}_1=\begin{bmatrix}0&1&0\\0&0&0\\0&0&0\end{bmatrix},\enspace
\mat{B}_1=\begin{bmatrix}0&0\\0&0\\0&1\end{bmatrix},\qquad
\mat{A}_2=\begin{bmatrix}0&1&0\\0&0&1\\0&0&0\end{bmatrix},\enspace
\mat{B}_2=\begin{bmatrix}0&0\\0&0\\0&1\end{bmatrix},
\end{equation*}
respectively.  The linearisation on the left is not controllable, while that on the right is.

The example suggests that~(1)~classical linearisation is not independent of
parameterisation of controls and/or~(2)~the classical notion of linear
controllability is not independent of parameterisation of controls.  We shall
see in Section~\ref{subsec:equil-linearise} that both things, in fact, are
true: neither classical linearisation nor the classical linear
controllability test are feedback-invariant.  This may come as a surprise to
some.\oprocend
\end{example}

This example has been particularly chosen to provide probably the simplest illustration of the phenomenon of lack of feedback-invariance of common control theoretic constructions.  Therefore, it should not be a surprise that an astute reader will notice that linearising the ``uncontrollable'' system about the control $(1,0)$ rather than the control $(0,0)$ will square things away as concerns the discrepancy between the two linearisations.  But after doing this, the questions of, ``What are the proper \emph{definitions} of linearisation and linear controllability?''\ still remain.  Moreover, one might expect that as one moves to constructions in control theory more advanced than mere linearisation, the dependence of these constructions on the parameterisation of controls becomes more pronounced.  Thus the likelihood that a sophisticated construction, made using a specific control parameterisation, is feedback-invariant is quite small, and in any case would need proof to verify that it is.  Such verification is not typically part of the standard development of methodologies in control theory.  There are at least three reasons for this:~(1)~the importance of feedback-invariance is not universally recognised;~(2)~such verifications are generally extremely difficult, nearly impossible, in fact;~(3)~most methodologies will fail the verification, so it is hardly flattering to one's methodology to point this out.  Some discussion of this is made by \citet{ADL:12a}\@.

But the bottom line is that our framework simply eliminates the need for any of this sort of verification.  As long as one remains within the framework, feedback-invariance is guaranteed.  One of the central goals of the paper is to provide the means by which one does not have to leave the framework to get things done.  As we shall see, certain technical difficulties have to be overcome to achieve this.

\subsubsection{Fibred manifold models}

As we have tried to make clear in the discussion just preceding, the standard model for control theory has the unpleasant attribute of depending on parameterisation of controls.  A natural idea to overcome this unwanted dependence is to do with controls as one does with states: regard them as taking values in a differentiable manifold.  Moreover, the manner in which control enters the model should also be handled in an intrinsic manner.  This leads to the ``fibred manifold'' picture of a control system which, as far as we can tell, originated in the papers of \citet{RWB:77} and \citet{JCW:79}\@, and was further developed by \citet{HN/AJV:82}\@.  This idea has been pickup up on by many researchers in geometric control theory, and we point to the papers~\cite{MBL/MCML:09,JCPB:84,MDT/AI:03,BL:03} as illustrative examples.

The basic idea is this.  A control system is modelled by a fibred manifold $\map{\pi}{\man{C}}{\man{M}}$ and a bundle map $\map{F}{\man{C}}{\tb{\man{M}}}$ over $\id_{\man{M}}$\@:
\begin{equation*}
\xymatrix{{\man{C}}\ar[r]^{F}\ar[rd]_{\pi}&
{\tb{\man{M}}}\ar[d]^{\tbproj{\man{M}}}\\&{\man{M}}}
\end{equation*}
One says that $F$ is ``a vector field over the bundle map $\pi$\@.''
Trajectories are then curves $t\mapsto x(t)$ in $\man{M}$ satisfying
$\dot{x}(t)=F(u(t))$ for some $t\mapsto u(t)$ satisfying $x(t)=\pi\scirc
u(t)$\@.  When it is applicable, this is an elegant and profitable model for
control theory.  For example, for control models that arise in problems of
differential geometry or the calculus of variations, this can be a useful
model.

The difficulty with the model is that it is not always applicable, especially
in physical system models.  The problem that arises is the strong regularity
of the control set and, implicitly, the controls: $\man{C}$ is a manifold so
it is naturally the codomain for smooth curves.  In practice, control sets in
physical models are seldom manifolds, as bounds on controls lead to
boundaries of the control set.  Moreover, the boundary sets are seldom
smooth.  Also, as we have mentioned above, controls cannot be restricted to
be smooth or piecewise smooth; natural classes of controls are typically
merely measurable.  These matters become vital in optimal control theory
where bounds on control sets lead to bang-bang extremals.  When these
considerations are overlaid on the fibred manifold picture, it becomes
considerably less appealing and indeed problematic.  One might try to patch
up the model by generalising the structure, but at some point it ceases to be
worthwhile; the framework is simply not well suited to certain problems of
control theory.

\subsubsection{Differential inclusion models}\label{subsubsec:di-intro}

Another way to eliminate the control dependence seen in the models with fixed
control parameterisation is to instead work with differential inclusions.  A
differential inclusion, roughly (we will be precise about differential
inclusions in Section~\ref{subsec:di}), assigns to each $x\in\man{M}$ a
subset $\sX(x)\subset\tb[x]{\man{M}}$\@, and trajectories are curves
$t\mapsto x(t)$ satisfying $\dot{x}(t)\in\sX(x(t))$\@.  There is a
well-developed theory for differential inclusions, and we refer to the
literature for what is known,~\eg~\cite{JPA/AC:84,AFF:88,GVS:02}\@.  There
are many appealing aspects to differential inclusions as far as our
objectives here are concerned.  In particular, differential inclusions do
away with the explicit parameterisations of the admissible tangent vectors at
a state $x\in\man{M}$ by simply prescribing this set of admissible tangent
vectors with no additional structure.  Moreover, differential inclusions
generalise the control-parameterised systems described above.  Indeed, given
such a control-parameterised system with dynamics $F$\@, we associate the
differential inclusion
\begin{equation*}
\sX_F(x)=\setdef{F(x,u)}{u\in\cs{C}}.
\end{equation*}
The trouble with differential inclusions is that their theory is quite
difficult to understand if one just starts with differential inclusions
coming ``out of the blue.''  Indeed, it is immediately clear that one needs
some sort of conditions on a differential inclusion to ensure that
trajectories exist.  Such conditions normally come in the form of some
combination of compactness, convexity, and semicontinuity.  However, the
differential inclusions that arise in control theory are \emph{highly}
structured; certainly they are more regular than merely semicontinuous and
they automatically possess many trajectories.  Moreover, it is not clear how
to develop an \emph{independent} theory of differential inclusions,~\ie~one
not making reference to standard models for control theory, that captures the
desired structure (in Example~\enumdblref{eg:backandforth}{enum:geomdi} we
suggest a natural way of characterising a class of differential inclusions
useful in geometric control theory).  Also, differential inclusions do not
themselves,~\ie~without additional structure, capture the notion of a flow
that is often helpful in the standard control-parameterised models,~\eg~in
the Maximum Principle of optimal control theory,~\cf~\cite{HJS:02}\@.
However, differential inclusions \emph{are} a useful tool for studying
trajectories, and we include them in our development of our new framework in
Section~\ref{sec:gcs}\@.

\subsubsection{The ``behavioural'' approach}

Starting with a series of papers~\cite{JCW:86a,JCW:86b,JCW:87} and the often
cited review~\cite{JCW:91}\@, \citeauthor{JCW:91} provides a framework for
studying system theory, with an emphasis on linear systems.  The idea in this
approach is to provide a framework for dynamical systems as subsets of
general functions of generalised time taking values in a set.  The framework
is also intended to provide a mathematical notion of interconnection as
relations in a set.  In this framework, the most general formulation is quite
featureless,~\ie~maps between sets and relations in sets.  With this level of
generality, the basic questions have a computer science flavour to them, in
terms of formal languages.  When one comes to making things more concrete,
say by making the time-domain an interval in $\real$ for continuous-time
systems, one ends up with differential-algebraic equations describing the
behaviours and relations.  For the most part, these ideas seem to have been
only reasonably fully developed for linear models~\cite{JWP/JCW:98}\@; we are
not aware of substantial work on nonlinear systems in the behavioural
approach.  It is also the case that the considerations of
feedback-invariance, such as we discuss above, are not a part of the current
landscape in behavioural models, although this is possible within the context
of linear systems,~\cf~the beautiful book of \cite{WMW:85}\@.

Thus, while there are some idealogical similarities with our objectives and
those of the behavioural approach, our thinking in this paper is in a quite
specific and complementary direction to the existing work on the behavioural
point of view.

\subsection{Attributes of a modelling framework for geometric control
systems}\label{subsec:attributes}

The preceding sections are meant to illustrate some standard frameworks for modelling control systems and the motivation for consideration of these, as well as pointing out their limitations.  If one is going to propose a modelling framework, it is important to understand \emph{a priori} just what it is that one hopes to be able to do in this framework.  Here is a list of possible criteria, criteria that we propose to satisfy in our framework.
\begin{compactenum}
\item Models should provide for control parameterisation-independent constructions as discussed above.
\item We believe that being able to handle real analytic systems is essential
to a useful theory.  In practice, any smooth control system is also real
analytic, and one wants to be able to make use of real analyticity to both
strengthen conclusions,~\eg~the real analytic version of Frobenius's
Theorem~\cite{TN:66}\@, and to weaken hypotheses,~\eg~the infinitesimal
characterisation of invariant distributions~\cite[\eg][Lemma~5.2]{AAA/YS:04}\@.
\item \label{enum:consreg} The framework should be able to handle regularity
in an internally consistent manner.  This means, for example, that the
conclusions should be consistent with hypotheses,~\eg~smooth hypotheses with
continuous conclusions suggest that the framework may not be perfectly
natural or perfectly well-developed.  The pursuit of this internal
consistency in the real analytic case contributes to many of the difficulties
we encounter in the paper.
\item The modelling framework should seamlessly deal with distinctions between local and global.  Many notions in control theory are highly localised,~\eg~local controllability of real analytic control systems.  A satisfactory framework should include a systematic way of dealing with constructions in control theory that are of an inherently local nature.  Moreover, the framework should allow a systematic means of understanding the passage from local to global in cases where this is possible and/or interesting.  As we shall see, there are some simple instances of these phenomena that can easily go unnoticed if one is not looking for them.
\item Our interest is in \emph{geometric} control theory, as we believe this is the right framework for studying nonlinear systems in general.  A proper framework for geometric control theory should make it natural to use the tools of differential geometry.
\item While (we believe that) differential geometric methods are essential in nonlinear control theory, the quest for geometric elegance should not be carried out at the expense of a useful theory.
\end{compactenum}

\subsection{An outline of the paper}

Let us discuss briefly the contents of the paper.

One of the essential elements of the paper is a characterisation of seminorms
for the various topologies we use.  Our definitions of these seminorms unify
the presentation of the various degrees of regularity we
consider\textemdash{}finitely differentiable, Lipschitz, smooth, holomorphic,
and real analytic\textemdash{}making it so that, after the seminorms are in
place, these various cases can be treated in very similar ways in many cases.
The key to the construction of the seminorms that we use is the use of
connections to decompose jet bundles into direct sums.  In
Section~\ref{sec:connections} we present these constructions.  As we see in
Section~\ref{sec:analytic-topology}\@, in the real analytic case, some
careful estimates must be performed to ensure that the geometric seminorms we
use do, indeed, characterise the real analytic topology.

In Sections~\ref{sec:smooth-topology}\@,~\ref{sec:holomorphic-topology}\@,
and~\ref{sec:analytic-topology} we describe topologies for spaces of finitely
differentiable, Lipschitz, smooth, holomorphic, and real analytic vector
fields.  (While we do not have a \emph{per se} interest in holomorphic
systems, holomorphic geometry has an important part to play in real analytic
geometry.)  While these topologies are more or less classical in the smooth,
finitely differentiable, and holomorphic cases, in the real analytic case the
description we give is less well-known, and indeed many of our results here
are new, or provide new and useful ways of understanding existing results.

Time-varying vector fields feature prominently in geometric control theory.
In Section~\ref{sec:time-varying} we review some notions concerning such
vector fields and develop a few not quite standard constructions and results
for later use.  In the smooth case, the ideas we present are probably
contained in the work of \citet{AAA/RVG:78} (see also~\cite{AAA/YS:04}), but
our presentation of the real analytic case is novel.  For this reason, we
present a rather complete treatment of the smooth case (with the finitely
differentiable and Lipschitz cases following along similar lines) so as to
provide a context for the more complicated real analytic case.  We should
point out that, even in the smooth case, we use properties of the topology
that are not normally called upon, and we see that it is these deeper
properties that really tie together the various regularity hypotheses we use.
Indeed, what our presentation reveals is the connection between the standard
pointwise\textemdash{}in time and state\textemdash{}conditions placed on
time-varying vector fields and topological characterisations.  This is, we
believe, a fulfilling way of understanding the meaning of the usual pointwise
conditions.

In Section~\ref{sec:systems} we review quite precisely a fairly general
standard modelling framework in geometric control theory.  While ultimately
we wish to assert that there are some difficulties with this framework,
understanding it clearly will give us some context for what will be, frankly,
our rather abstract notion of a control system to follow.  Also, we do wish
to make sure that our proposed model does indeed generalise this more
concrete and standard notion, so to prove this we need precise definitions.
Additionally, as with time-varying vector fields, we show how natural
pointwise regularity conditions are equivalent to topological
characterisations of systems.  Thus, while we do generalise the standard
modelling framework for control theory, in doing so we arrive at a deeper
understanding of this framework.  For example, we introduce for the first
time the notion of a ``real analytic control system,'' which means that the
real analytic structure is fully integrated into the structure of the control
system; this is only made possible by understanding the topology for the
space of real analytic vector fields.  As a result, seemingly tautological
statements like, ``A real analytic control-affine system is a real analytic
control system,'' now are theorems in our framework.  Also, interestingly, we
will show that, in many cases, our more general modelling framework can be
cast in the standard framework, albeit in a non-obvious way; see
Example~\enumdblref{eg:gcs}{enum:gcs->cs}\@.

In Section~\ref{sec:gcs} we provide our modelling framework for geometric
control systems, defining what we shall call ``\gcs{}s.\footnote{The
terminology ``tautological'' arises from two different attributes of our
framework.  First of all, when one makes the natural connection from our
systems to standard control systems, we encounter the identity map
(Example~\enumdblref{eg:gcs}{enum:gcs->cs}).  Second, in our framework we
prove that the only pure feedback transformation is the identity
transformation (\cf~Proposition~\ref{prop:natural-isomorphism}).}''  After
developing the background needed, we provide the definitions and then give
the notion of a trajectory for these systems.  We also show that our
framework includes the standard framework of Section~\ref{sec:systems} as a
special case.  We carefully establish correspondences between our generalised
models, the standard models, and differential inclusion models.  Included in
this correspondence is a description of the relationships between
trajectories for these models.  One feature of our framework that will appear
strange initially is our use of presheaves and sheaves.  These are the
devices by which we can attempt to patch together local constructions to give
global constructions.  We understand that the use of this language will seem
unnecessarily complicated initially.  However, it will have its uses in the
paper,~\eg~our notion of transformations between \gcs{}s is based on a
standard construction in sheaf theory, and we will point out places where the
reader may have unwittingly encountered some shadows of sheaf theory, even in
familiar places in control theory.

We study the linearisation of \gcs{}s in Section~\ref{sec:linearisation}\@.
The theory here has many satisfying elements attached to it.  First of all,
the framework naturally suggests \emph{two} sorts of linearisation, one with
respect to a reference trajectory and another with respect to a reference
flow.  This is an interesting distinction, and one that is, as far as we
know, hitherto not made clear in the literature.  Also, of course, our theory
comprehends and rectifies the problems encountered in
Example~\ref{eg:bad-linearise}\@.

What is presented in this paper is the result of initial explorations of a
modelling framework for geometric control theory.  We certainly have not
fully fleshed out all parts of this framework ourselves, despite the
substantial length of the paper.  In the closing section of the paper,
Section~\ref{sec:future-work}\@, we outline places where there is obvious
further work to be done.

\subsection{Summary of contributions}

This is a long and complex paper with many results, some significant, and
some necessary for the foundations of the approach, but not necessarily
significant \emph{per se}\@.  In order to facilitate the reading of the
paper, we highlight the contributions that we feel are important.  First we
point out the more significant contributions.
\begin{compactenum}
\item The main contribution of the paper is the general feedback-invariant
framework.  This main contribution has with it a few novel components.
\begin{compactenum}[(a)]
\item Our framework generalises the standard formulation and has some
satisfying relationships with the standard theory and the theory of
differential inclusions; see Proposition~\ref{prop:gcs-cs-di} and the
trajectory equivalence results of Section~\ref{subsec:trajequiv}\@.  We
conclude, for example, that our generalised formulation agrees with the
standard formulation in two important cases:~(i)~for control-affine systems
with arbitrary control sets (Theorem~\ref{the:cs->gcs-equivII});~(ii)~for
systems depending generally on the control with compact control sets
(Theorem~\ref{the:cs->gcs-equivI}).
\item The framework relies in an essential and nontrivial way on topologies
for spaces of vector fields.  The full development of these topologies, and
their integration into a theory for control systems, is fully executed here
for the first time.
\item The framework relies in an essential and nontrivial way on topologies
for spaces of vector fields.  The full development of these topologies, and
their integration into a theory for control systems, is fully executed here
for the first time.
\item The formulation uses the theory of presheaves and sheaves in an essential way.
\item Using a notion of morphism borrowed from sheaf theory, we \emph{prove}
that equivalence for our systems is simply diffeomorphism equivalence of
vector fields; see Proposition~\ref{prop:natural-isomorphism}\@.  That is to
say, we \emph{prove} that our framework cannot involve any ``feedback
transformation'' in the usual sense.
\end{compactenum}
\item We provide, for the first time, a comprehensive treatment of real analytic time-varying vector fields and control systems.  In particular,
\begin{compactenum}[(a)]
\item we provide a concrete, usable, geometric characterisation of the real
analytic topology by specifying a family of geometric seminorms
(Theorem~\ref{the:Comega-seminorms}),
\item we provide conditions that ensure that a real analytic vector field
with measurable time dependence will have a flow depending on initial
conditions in a real analytic manner (Theorem~\ref{the:Comega-flows}),
\item we provide conditions that ensure that trajectories for a real analytic control system depend on initial conditions in a real analytic manner (Propositions~\ref{prop:open-loopcsI} and~\ref{prop:open-loopcsII}), and
\item we show that real analytic vector fields depending measurably on time
and/or continuously on a parameter can often be extended to holomorphic
vector fields depending on time or parameter (Theorems~\ref{the:Comega->Chol}
and~\ref{the:Comega->Cholparam}).
\end{compactenum}
The last three results rely, sometimes in highly nontrivial ways, on the properties of the real analytic topology for vector fields.
\item We fully develop various ``weak'' formulations of properties such as
continuity, boundedness, measurability, and integrability for spaces of
finitely differentiable, Lipschitz, smooth, and real analytic vector fields.
These weak formulations come in two forms, one for evaluations of vector
fields on functions by Lie differentiation, which we call the ``weak-$\sL$''
topology (see Theorems~\ref{the:COinfty-weak}\@,~\ref{the:COm-weak}\@,%
~\ref{the:COmm'-weak}\@, and~\ref{the:Comega-weak} and their corollaries),
and one for evaluations in time and space (see
Theorems~\ref{the:smooth-td-summary}\@,~\ref{the:mm'-td-summary}\@,
and~\ref{the:ra-td-summary}).  These results use deep properties of the
topologies for spaces of vector fields derived in
Sections~\ref{sec:smooth-topology} and~\ref{sec:analytic-topology}\@.  In the
existing literature, these weak formulations are often used without reference
to their ``strong'' counterparts; here we make the (unsurprising, but
sometimes nontrivial) link explicit.
\item In Section~\ref{sec:linearisation} we provide a coherent theory for
linearisation of systems in our framework.  The theory of linearisation that
we develop is necessarily feedback-invariant, and as a consequence reveals
some interesting structure that has previously been hidden by the standard
treatment of linearisation which is \emph{not} feedback-invariant, as we have
seen in Example~\ref{eg:bad-linearise}\@.\savenum
\end{compactenum}
Along the way to these substantial definitions and results, we uncover a few minor, but still interesting, results and constructions.
\begin{compactenum}\resumenum
\item We use to advantage some not entirely elementary geometric constructions to make elegant coordinate-free proofs.  Here are some instances of this.
\begin{compactenum}[(a)]
\item We provide a decomposition for jet bundles of sections of a vector
bundle using the theory of connections; see Lemma~\ref{lem:Jrdecomp}\@.  This
decomposition is used to provide a concrete and useful collection of
seminorms for the finitely differentiable, Lipschitz, and smooth compact-open
topologies, and the real analytic topology.  Indeed, without these seminorms,
our descriptions of these topologies would be incomprehensible, as opposed to
merely difficult as it already is in the real analytic case.
\item We use our seminorms in an essential way to prove the equivalence of ``weak-$\sL$'' and ``strong'' versions of the finitely differentiable, Lipschitz, and smooth compact-open topologies, and the real analytic topology for vector fields; see Theorems~\ref{the:COinfty-weak}\@,~\ref{the:COm-weak}\@,%
~\ref{the:COmm'-weak}\@, and~\ref{the:Comega-weak}\@.
\item These seminorms allow for relatively clean characterisations of the finitely differentiable, Lipschitz, and smooth compact-open, and real analytic topologies for vector fields on tangent bundles, using induced affine connections and Riemannian metrics on tangent bundles.  These constructions appear in the proofs concerning linearisation; see Lemmata~\ref{lem:XTflow} and~\ref{lem:open-loop-lin}\@.
\item The double vector bundle structure of the double tangent bundle
$\tb{\tb{\man{M}}}$ is used to provide a slick justification of our
definition of linearisation, culminating in the
formula~\eqref{eq:linear-derivation}\@.
\end{compactenum}
\item We provide a ``weak-$\sL$'' characterisation of the compact-open topology for holomorphic vector fields on a Stein manifold; see Theorem~\ref{the:COhol-weak}\@.
\end{compactenum}

\subsection{Notation, conventions, and background}\label{subsec:notation}

In this section we overview what is needed to read the paper.  We do use a lot of specialised material in essential ways, and we certainly do not review this comprehensively.  Instead, we simply provide a few facts, the notation we shall use, and recommended sources.  Throughout the paper we have tried to include precise references to material needed so that a reader possessing enthusiasm and lacking background can begin to chase down all of the ideas upon which we rely.

We shall use the slightly unconventional, but perfectly rational, notation of writing $A\subset B$ to denote set inclusion, and when we write $A\subsetneq B$ we mean that $A\subset B$ and $A\not=B$\@.  By $\id_A$ we denote the identity map on a set $A$\@.  For a product $\prod_{i\in I}X_i$ of sets, $\map{\pr_j}{\prod_{i\in I}X_i}{X_j}$ is the projection onto the $j$th component.  For a subset $A\subset X$\@, we denote by $\charfunc{A}$ the characteristic function of $A$\@,~\ie
\begin{equation*}
\charfunc{A}(x)=\begin{cases}1,&x\in A,\\0,&x\not\in A.\end{cases}
\end{equation*}
By $\card(A)$ we denote the cardinality of a set $A$\@.  By $\symmgroup{k}$ we denote the symmetric group on $k$ symbols.  We shall have occasion to talk about set-valued maps.  If $X$ and $Y$ are sets and $\Phi$ is a set-valued map from $X$ to $Y$\@,~\ie~$\Phi(x)$ is a subset of $Y$\@, we shall write $\setmap{\Phi}{X}{Y}$\@.  By $\integer$ we denote the set of integers, with $\integernn$ denoting the set of nonnegative integers and $\integerp$ denoting the set of positive integers.  We denote by $\real$ and $\complex$ the sets of real and complex numbers.  By $\realnn$ we denote the set of nonnegative real numbers and by $\realp$ the set of positive real numbers. By $\erealnn=\realnn\cup\{\infty\}$ we denote the extended nonnegative real numbers.  By $\delta_{jk}$\@, $j,k\in\{1,\dots,n\}$\@, we denote the Kronecker delta.

We shall use constructions from algebra and multilinear algebra, referring to~\cite{TWH:74}\@,~\cite[Chapter~III]{NB:89a}\@, and~\cite[\S{}IV.5]{NB:90}\@.  If $\alg{F}$ is a field (for us, typically $\alg{F}\in\{\real,\complex\}$), if $\alg{V}$ is an $\alg{F}$-vector space, and if $A\subset\alg{V}$\@, by $\vecspan[\alg{F}]{A}$ we denote the subspace generated by $A$\@.  If $\alg{F}$ is a field and if $\alg{U}$ and $\alg{V}$ are $\alg{F}$-vector spaces, by $\Hom_{\alg{F}}(\alg{U};\alg{V})$ we denote the set of linear maps from $\alg{U}$ to $\alg{V}$\@.  We denote $\End_{\alg{F}}(\alg{V})=\Hom_{\alg{F}}(\alg{V};\alg{V})$ and $\dual{\alg{V}}=\Hom_{\alg{F}}(\alg{V};\alg{F})$\@.  If $\alpha\in\dual{\alg{V}}$ and $v\in\alg{V}$\@, we may sometimes denote by $\natpair{\alpha}{v}\in\alg{F}$ the natural pairing.  The $k$-fold tensor product of $\alg{V}$ with itself is denoted by $\tensor[k]{\alg{V}}$\@.  Thus, if $\alg{V}$ is finite-dimensional, we identify $\tensor[k]{\dual{\alg{V}}}$ with the $k$-multilinear $\alg{F}$-valued functions on $\alg{V}^k$ by
\begin{equation*}
(\alpha^1\otimes\dots\otimes\alpha^k)(v_1,\dots,v_k)=
\alpha^1(v_1)\cdots\alpha^k(v_k).
\end{equation*}
By $\Symalg[k]{\dual{\alg{V}}}$ we denote the symmetric tensor algebra of degree $k$\@, which we identify with the symmetric $k$-multilinear $\alg{F}$-valued functions on $\alg{V}^k$\@, or polynomial functions of homogeneous degree $k$ on $\alg{V}$\@.

If $\metric$ is an inner product on a $\real$-vector space $\alg{V}$\@, we denote by $\metric^\flat\in\Hom_{\real}(\alg{V};\dual{\alg{V}})$ the associated mapping and by $\metric^\sharp\in\Hom_{\real}(\dual{\alg{V}};\alg{V})$ the inverse of $\metric^\flat$ when it is invertible.

For a topological space $\ts{X}$ and $A\subset\ts{X}$\@, $\interior(A)$
denotes the interior of $A$ and $\closure(A)$ denotes the closure of $A$\@.
Neighbourhoods will always be open sets.  The support of a continuous
function $f$ (or any other kind of object for which it makes sense to have a
value ``zero'') is denoted by $\supp(f)$\@.

By $\oball{r}{\vect{x}}\subset\real^n$ we denote the open ball of radius $r$ and centre $\vect{x}$\@.  In like manner, $\cball{r}{\vect{x}}$ denotes the closed ball.  If $r\in\realp$ and if $x\in\field$\@, $\field\in\{\real,\complex\}$\@, we denote by
\begin{equation*}
\odisk[]{r}{x}=\setdef{x'\in\field}{\snorm{x'-x}<r}
\end{equation*}
the disk of radius $r$ centred at $x$\@.  If $\vect{r}\in\realp^n$ and if $\vect{x}\in\field^n$\@, we denote by
\begin{equation*}
\odisk[]{\vect{r}}{\vect{x}}=\odisk[]{r_1}{x_1}\times\dots\times
\odisk[]{r_n}{x_n}
\end{equation*}
the polydisk with radius $\vect{r}$ centred at $\vect{x}$\@.  In like manner, $\cdisk[]{\vect{r}}{\vect{x}}$ denotes the closed polydisk.

Elements of $\field^n$\@, $\field\in\{\real,\complex\}$\@, are typically
denoted with a bold font,~\eg~``$\vect{x}$\@.''  The standard basis for
$\field^n$ is denoted by $(\vect{e}_1,\dots,\vect{e}_n)$\@.  By $\mat{I}_n$
we denote the $n\times n$ identity matrix.  We denote by
$\lin{\real^n}{\real^m}$ the set of linear maps from $\real^n$ to $\real^m$
(this is the same as $\Hom_{\real}(\real^n;\real^m)$\@, of course, but the
more compact notation is sometimes helpful).  The invertible linear maps on
$\real^n$ we denote by $\GL[n]{\real}$\@.  By
$\lin{\real^{n_1},\dots,\real^{n_k}}{\real^m}$ we denote the set of
multilinear mappings from $\prod_{j=1}^k\real^{n_j}$ to $\real^m$\@.  We
abbreviate by $\lin[k]{\real^n}{\real^m}$ the $k$-multilinear maps from
$(\real^n)^k$ to $\real^m$\@.  We denote by $\symlin[k]{\real^n}{\real^m}$
the set of symmetric $k$-multilinear maps from $(\real^n)^k$ to $\real^m$\@.
With our notation above, $\symlin[k]{\real^n}{\real^m}\simeq
\Symalg[k]{\dual{(\real^n)}}\otimes\real^m$\@, but, again, we prefer the
slightly more compact notation in this special case.

If $\nbhd{U}\subset\real^n$ is open and if $\map{\vect{\Phi}}{\nbhd{U}}{\real^m}$ is differentiable at $\vect{x}\in\nbhd{U}$\@, we denote its derivative by $\linder{\vect{\Phi}}(\vect{x})$\@.  Higher-order derivatives, when they exist, are denoted by $\linder[r]{\vect{\Phi}}(\vect{x})$\@, $r$ being the order of differentiation.  We will also use the following partial derivative notation.  Let $\nbhd{U}_j\subset\real^{n_j}$ be open, $j\in\{1,\dots,k\}$\@, and let $\map{\vect{\Phi}}{\nbhd{U}_1\times\dots\times\nbhd{U}_k}{\real^m}$ be continuously differentiable.  The derivative of the map
\begin{equation*}
\vect{x}_j\mapsto\vect{\Phi}(\vect{x}_{1,0},\dots,\vect{x}_j,
\dots,\vect{x}_{k,0})
\end{equation*}
at $\vect{x}_{j,0}$ is denoted by $\plinder{j}{\vect{\Phi}}(\vect{x}_{1,0},\dots,\vect{x}_{k,0})$\@.  Higher-order partial derivatives, when they exist, are denoted by $\plinder[r]{j}{\vect{\Phi}}(\vect{x}_{1,0},\dots,\vect{x}_{k,0})$\@, $r$ being the order of differentiation.  We recall that if $\map{\vect{\Phi}}{\nbhd{U}}{\real^m}$ is of class $\C^k$\@, $k\in\integerp$\@, then $\linder[k]{\vect{\Phi}}(\vect{x})$ is symmetric.  We shall sometimes find it convenient to use multi-index notation for derivatives.  A \defn{multi-index} with length $n$ is an element of $\integernn^n$\@,~\ie~an $n$-tuple $I=(i_1,\dots,i_n)$ of nonnegative integers.  If $\map{\vect{\Phi}}{\nbhd{U}}{\real^m}$ is a smooth function, then we denote
\begin{equation*}
\linder[I]{\vect{\Phi}}(\vect{x})=
\plinder[i_1]{1}{\cdots\plinder[i_n]{n}{\vect{\Phi}}}(\vect{x}).
\end{equation*}
We will use the symbol $\snorm{I}=i_1+\dots+i_n$ to denote the order of the derivative.  Another piece of multi-index notation we shall use is
\begin{equation*}
\vect{a}^I=a_1^{i_1}\cdots a_n^{i_n},
\end{equation*}
for $\vect{a}\in\real^n$ and $I\in\integernn^n$\@.  Also, we denote $I!=i_1!\cdots i_n!$\@.

If $\alg{V}$ is a $\real$-vector space and if $A\subset\alg{V}$\@, we denote by $\cohull(A)$ the convex hull of $A$\@, by which we mean the set of all convex combinations of elements of $A$\@.

Our differential geometric conventions mostly follow~\cite{RA/JEM/TSR:88}\@.
Whenever we write ``manifold,'' we mean ``second-countable Hausdorff
manifold.''  This implies, in particular, that manifolds are assumed to be
metrisable~\cite[Corollary~5.5.13]{RA/JEM/TSR:88}\@.  If we use the letter
``$n$'' without mentioning what it is, it is the dimension of the connected
component of the manifold $\man{M}$ with which we are working at that time.
The tangent bundle of a manifold $\man{M}$ is denoted by
$\map{\tbproj{\man{M}}}{\tb{\man{M}}}{\man{M}}$ and the cotangent bundle by
$\map{\ctbproj{\man{M}}}{\ctb{\man{M}}}{\man{M}}$\@.  The derivative of a
differentiable map $\map{\Phi}{\man{M}}{\man{N}}$ is denoted by
$\map{\tf{\Phi}}{\tb{\man{M}}}{\tb{\man{N}}}$\@, with
$\tf[x]{\Phi}=\tf{\Phi}|\tb[x]{\man{M}}$\@.  If $I\subset\real$ is an
interval and if $\map{\xi}{I}{\man{M}}$ is a curve that is differentiable at
$t\in I$\@, we denote the tangent vector field to the curve at $t$ by
$\xi'(t)=\tf[t]{\xi}(1)$\@.  We use the symbols $\Phi^*$ and $\Phi_*$ for
pull-back and push-forward.  Precisely, if $g$ is a function on $\man{N}$\@,
$\Phi^*g=g\scirc\Phi$\@, and if $\Phi$ is a diffeomorphism, if $f$ is a
function on $\man{M}$\@, if $X$ is a vector field on $\man{M}$\@, and if $Y$
is a vector field on $\man{N}$\@, we have $\Phi_*f=f\scirc\Phi^{-1}$\@,
$\Phi_*X=\tf{\Phi}\scirc X\scirc\Phi^{-1}$\@, and
$\Phi^*Y=\tf{\Phi^{-1}}\scirc Y\scirc\Phi$\@.  The flow of a vector field $X$
is denoted by $\flow{X}{t}$\@, so $t\mapsto\flow{X}{t}(x)$ is the integral
curve of $X$ passing through $x$ at $t=0$\@.  We shall also use time-varying
vector fields, but will develop the notation for the flows of these in the
text.

If $\map{\pi}{\man{E}}{\man{M}}$ is a vector bundle, we denote the fibre over
$x\in\man{M}$ by $\man{E}_x$ and we sometimes denote by $0_x$ the zero vector
in $E_x$\@.  If $\man{S}\subset\man{M}$ is a submanifold, we denote by
$\man{E}|\man{S}$ the restriction of $\man{E}$ to $\man{S}$ which we regard
as a vector bundle over $\man{S}$\@.  The \defn{vertical subbundle} of
$\man{E}$ is the subbundle of $\tb{\man{E}}$ defined by
$\vb{\man{E}}=\ker(\tf{\pi})$\@.  If $\metric$ is a fibre metric on
$\man{E}$\@,~\ie~a smooth assignment of an inner product to each of the
fibres of $\man{E}$\@, then $\dnorm{\cdot}_{\metric}$ denotes the norm
associated with the inner product on fibres.  If
$\map{\pi}{\man{E}}{\man{M}}$ is a vector bundle and if
$\map{\Phi}{\man{N}}{\man{M}}$ is a smooth map, then
$\map{\Phi^*\pi}{\Phi^*\man{E}}{\man{N}}$ denotes the pull-back of $\man{E}$
to $\man{N}$~\cite[\S{}III.9.5]{IK/PWM/JS:93}\@.  The dual of a vector bundle
$\map{\pi}{\man{E}}{\man{M}}$ is denoted by
$\map{\pi^*}{\dual{\man{E}}}{\man{M}}$\@.

Generally we will try hard to avoid coordinate computations.  However, they are sometimes unavoidable and we will use the Einstein summation convention when it is convenient to do so, but we will not do so slavishly.

We will work in both the smooth and real analytic categories, with occasional
forays into the holomorphic category.  We will also work with finitely
differentiable objects,~\ie~objects of class $\C^r$ for $r\in\integernn$\@.
(We will also work with Lipschitz objects, but will develop the notation for
these in the text.)  A good reference for basic real analytic analysis
is~\cite{SGK/HRP:02}\@, but we will need ideas going beyond those from this
text, or any other text.  Relatively recent work
of~\eg~\cite{PD:10}\@,~\cite{DV:13}\@, and~\cite{PD/DV:00} has shed a great
deal of light on real analytic analysis, and we shall take advantage of this
work.  An analytic manifold or mapping will be said to be of \defn{class
$\C^\omega$}\@.  Let $r\in\integernn\cup\{\infty,\omega\}$\@.  The set of
mappings of class $\C^r$ between manifolds $\man{M}$ and $\man{N}$ is denoted
by $\mappings[r]{\man{M}}{\man{N}}$\@.  We abbreviate
$\func[r]{\man{M}}=\mappings[r]{\man{M}}{\real}$\@.  The set of sections of a
vector bundle $\map{\pi}{\man{E}}{\man{M}}$ of class $\C^r$ is denoted by
$\sections[r]{\man{E}}$\@.  Thus, in particular, $\sections[r]{\tb{\man{M}}}$
denotes the set of vector fields of class $\C^r$\@.  We shall think of
$\sections[r]{\man{E}}$ as a $\real$-vector space with the natural pointwise
addition and scalar multiplication operations.  If $f\in\func[r]{\man{M}}$\@,
$\d{f}\in\sections[r]{\ctb{\man{M}}}$ denotes the differential of $f$\@.  If
$X\in\sections[r]{\tb{\man{M}}}$ and $f\in\func[r]{\man{M}}$\@, we denote the
Lie derivative of $f$ with respect to $X$ by $\lieder{X}{f}$\@.

We also work with holomorphic,~\ie~complex analytic, manifolds and associated geometric constructions; real analytic geometry, at some level, seems to unavoidably rely on holomorphic geometry.  A nice overview of holomorphic geometry, and some of its connections to real analytic geometry, is given in the book of \citet{KC/YE:12}\@.  There are many specialised texts on the subject of holomorphic geometry, including~\cite{JPD:12,KF/HG:02,RCG/HR:65,LH:73} and the three volumes of~\citet{RCG:90a,RCG:90b,RCG:90c}\@.  For our purposes, we shall just say the following things.  By $\tb{\man{M}}$ we denote the holomorphic tangent bundle of $\man{M}$\@.  This is the object which, in complex differential geometry, is commonly denoted by $\man{T}^{1,0}\man{M}$\@.  For holomorphic manifolds $\man{M}$ and $\man{N}$\@, we denote by $\mappings[\hol]{\man{M}}{\man{N}}$ the set of holomorphic mappings from $\man{M}$ to $\man{N}$\@, by $\func[\hol]{\man{M}}$ the set of holomorphic functions on $\man{M}$ (note that these functions are $\complex$-valued, not $\real$-valued, of course), and by $\sections[\hol]{\man{E}}$ the space of holomorphic sections of an holomorphic vector bundle $\map{\pi}{\man{E}}{\man{M}}$\@.  We shall use both the natural $\complex$- and, by restriction, $\real$-vector space structures for $\sections[\hol]{\man{E}}$\@.

We will make use of the notion of a ``Stein manifold.''  For practical purposes, these can be taken to be holomorphic manifolds admitting a proper holomorphic embedding in complex Euclidean space.\footnote{The equivalence of this to other characterisations of Stein manifolds is due to \citet{RR:55}\@.  A reader unfamiliar with holomorphic manifolds should note that, unlike in the smooth or real analytic cases, it is \emph{not} generally true that an holomorphic manifold can be embedded in complex Euclidean space, even after the usual elimination of topological pathologies such as non-paracompactness.  For example, compact holomorphic manifolds can never be holomorphically embedded in complex Euclidean space.}  Stein manifolds are characterised by having lots of holomorphic functions, distinguishing them from general holomorphic manifolds,~\eg~compact holomorphic manifolds whose only holomorphic functions are those that are locally constant.  There is a close connection between Stein manifolds and real analytic manifolds, and this explains our interest in Stein manifolds.  We shall point out these connections as they arise in the text.

We shall occasionally make use of Cartan's Theorems~A and~B for Stein
manifolds and real analytic manifolds; these are theorems about the
cohomology of certain sheaves.  In the holomorphic case, the original source
is~\cite{HC:51}\@, but there are many good treatments in textbooks, including
in~\cite{JLT:02}\@.  For the real analytic case, the only complete reference
seems to be the original work of \citet{HC:57}\@, although the short book of
\citet{FG/PM/AT:86} is also helpful.  In using these theorems (and sometimes
in other places where we use sheaves) we will use the following notation.
Let $r\in\integernn\cup\{\infty,\omega,\hol\}$ and let $\man{M}$ be a smooth,
real analytic, or holomorphic manifold, such as is demanded by $r$\@.  By
$\sfunc[r]{\man{M}}$ we denote the sheaf of functions of class $\C^r$ and by
$\gfunc[r]{x}{\man{M}}$ the set of germs of this sheaf at $x\in\man{M}$\@.
If $\map{\pi}{\man{E}}{\man{M}}$ is a $\C^r$-vector bundle, then
$\ssections[r]{\man{E}}$ denotes the sheaf of $\C^r$-sections of $\man{E}$
with $\gsections[r]{x}{\man{E}}$ the set of germs at $x$\@.  The germ of a
function (\resp~section) at $x$ will be denoted by $[f]_x$ (\resp~$[\xi]_x$).

We will make use of jet bundles, and a standard reference is~\cite{DJS:89}\@.
Appropriate sections of~\cite{IK/PWM/JS:93} (especially~\S{}12) are also
useful.  If $\map{\pi}{\man{E}}{\man{M}}$ is a vector bundle and if
$k\in\integernn$\@, we denote by $\jet{k}{\man{E}}$ the bundle of $k$-jets of
$\man{E}$\@.  For a section $\xi$ of $\man{E}$\@, we denote by $j_k\xi$ the
corresponding section of $\jet{k}{\man{E}}$\@.  The projection from
$\jet{k}{\man{E}}$ to $\jet{l}{\man{E}}$\@, $l\le k$\@, is denoted by
$\pi^k_l$\@.  If $\man{M}$ and $\man{N}$ are manifolds, we denote by
$\jet{k}{(\man{M};\man{N})}$ the bundle of $k$ jets of mappings from
$\man{M}$ to $\man{N}$\@.  If $\Phi\in\mappings[\infty]{\man{M}}{\man{N}}$\@,
$j_k\Phi$ denotes its $k$-jet, which is a mapping from $\man{M}$ to
$\jet{k}{(\man{M};\man{N})}$\@.  In the proof of
Theorem~\ref{the:Cinfty-flows} we will briefly make use of jets of sections
of fibred manifolds.  We shall introduce there the notation we require, and
the reader can refer to \cite{DJS:89} to fill in the details.

We shall make use of connections, and refer
to~\cite[\S{}11,~\S{}17]{IK/PWM/JS:93} for a comprehensive treatment of
these, or to~\cite{SK/KN:63a} for another comprehensive treatment and an
alternative point of view.

We shall make reference to elementary ideas from sheaf theory; indeed we have
already made reference to sheaves above.  It will not be necessary to
understand this theory deeply, at least not in the present paper.  In
particular, a comprehensive understanding of sheaf cohomology is not
required, although, as indicated above, we do make use of Cartan's Theorems~A
and~B in places.  A nice introduction to the use of sheaves in smooth
differential geometry can be found in the book of \citet{SR:05b}\@.  More
advanced and comprehensive treatments include~\cite{GEB:97,MK/PS:90}\@, and
the classic~\cite{RG:58}\@.  The discussion of sheaf theory
in~\citepalias{Stacks} is also useful.  For readers who are expert in sheaf
theory, we comment that our reasons for using sheaves are not always the
usual ones, so an adjustment of point of view may be required.

We shall make frequent and essential use of nontrivial facts about locally
convex topological vector spaces, and refer to
\cite{JBC:90,AG:73,JH:66,HJ:81,WR:91,HHS/MPW:99} for details.  We shall also
access the contemporary research literature on locally convex spaces, and
will indicate this as we go along.  We shall denote by
$\lin{\alg{U}}{\alg{V}}$ the set of continuous linear maps from a locally
convex space $\alg{U}$ to a locally convex space $\alg{V}$\@.  In particular,
$\topdual{\alg{U}}$ is the topological dual of $\alg{U}$\@, meaning the
continuous linear scalar-valued functions.  We will break with the usual
language one sees in the theory of locally convex spaces and call what are
commonly called ``inductive'' and ``projective'' limits, instead ``direct''
and ``inverse'' limits, in keeping with the rest of category theory.

By $\lebmes$ we denote the Lebesgue measure on $\real$\@.  We will talk about measurability of maps taking values in topological spaces.  If $(\ts{T},\sM)$ is a measurable space and if $\ts{X}$ is a topological space, a mapping $\map{\Psi}{\ts{T}}{\ts{X}}$ is \defn{Borel measurable} if $\Psi^{-1}(\nbhd{O})\in\sM$ for every open set $\nbhd{O}\subset\ts{X}$\@.  This is equivalent to requiring that $\Psi^{-1}(\nbhd{B})\in\sM$ for every Borel subset $\nbhd{B}\subset\ts{X}$\@.

One not completely standard topic we shall need to understand is integration
of functions with values in locally convex spaces.  There are multiple
theories here,\footnote{\label{fn:integration}Most of the theories of
integration in locally convex spaces coincide for the sorts of locally convex
spaces we deal with.} so let us outline what we mean,
following~\cite{RB/AD:11}\@.  We let $(\ts{T},\sM,\mu)$ be a finite measure
space, let $\alg{V}$ be a locally convex topological vector space, and let
$\map{\Psi}{\ts{T}}{\alg{V}}$\@.  Measurability of $\Psi$ is Borel
measurability mentioned above, and we note that there are other forms of
measurability that arise for locally convex spaces (the comment made in
footnote~\ref{fn:integration} applies to these multiple notions of
measurability as well).  The notion of the integral we use is the
\defn{Bochner integral}\@.  This is well understood for Banach
spaces~\cite{JD/JJUJr:77} and is often mentioned in an offhand manner as
being ``the same'' for locally convex
spaces~\cite[\eg][page~96]{HHS/MPW:99}\@.  A detailed textbook treatment does
not appear to exist, but fortunately this has been worked out in the note
of~\cite{RB/AD:11}\@, to which we shall refer for details as needed.  One has
a notion of simple functions, meaning functions that are finite linear
combinations, with coefficients in $\alg{V}$\@, of characteristic functions
of measurable sets.  The \defn{integral} of a simple function
$\sigma=\sum_{j=1}^kv_j\chi_{A_j}$ is
\begin{equation*}
\int_{\ts{T}}\sigma\,\d{\mu}=\sum_{j=1}^k\mu(A_j)v_j,
\end{equation*}
in the usual manner.  A measurable function $\Psi$ is \defn{Bochner
approximable} if it can be approximated with respect to any continuous
seminorm by a net of simple functions.  A Bochner approximable function
$\Psi$ is \defn{Bochner integrable} if there is a net of simple functions
approximating $\Psi$ whose integrals converge in $\alg{V}$ to a unique value,
which is called the \defn{integral} of $\Psi$\@.  If $\alg{V}$ is separable
and complete, as will be the case for us in this paper, then a measurable
function $\map{\Psi}{\ts{T}}{\alg{V}}$ is Bochner integrable if and only if
\begin{equation*}
\int_{\ts{T}}p\scirc\Psi\,\d{\mu}<\infty
\end{equation*}
for every continuous seminorm $p$ on $\alg{V}$\@~\cite[Theorems~3.2
and~3.3]{RB/AD:11}\@.  This construction of the integral clearly agrees with
the standard construction of the Lebesgue integral for functions taking
values in $\real$ or $\complex$ (or any finite-dimensional vector space over
$\real$ or $\complex$\@, for that matter).  If $A\subset\alg{V}$\@, by
$\L^1(\ts{T};A)$ we denote the space of Bochner integrable functions with
values in $A$\@.  The space $\L^1(\ts{T};\alg{V})$ is itself a locally convex
topological vector space with topology defined by the seminorms
\begin{equation*}
\hat{p}(\Psi)=\int_{\ts{T}}p\scirc\Psi\,\d{\mu},
\end{equation*}
where $p$ is a continuous seminorm for $\alg{V}$~\cite[page~96]{HHS/MPW:99}\@.
In the case where $\ts{T}=I$ is an interval in $\real$\@, $\Lloc^1(I;A)$
denotes the set of locally integrable functions,~\ie~those functions whose
restriction to any compact subinterval is integrable.

While it does not generally make sense to talk about integrability of
measurable functions with values in a topological space, one \emph{can}
sensibly talk about essentially bounded functions.  This means that one needs
a notion of boundedness, this being supplied by a ``bornology.''\footnote{A
\defn{bornology} on a set $\ts{S}$ is a family $\sB$ of subsets of
$\ts{S}$\@, called \defn{bounded sets}\@, and satisfying the axioms:
\begin{compactenum}
\item $\ts{S}$ is covered by bounded sets,~\ie~$\ts{S}=\cup_{B\in\sB}B$\@;
\item subsets of bounded sets are bounded,~\ie~if $B\in\sB$ and if $A\subset B$\@, then $A\in\sB$\@;
\item finite unions of bounded sets are bounded,~\ie~if $B_1,\dots,B_k\in\sB$\@, then $\cup_{j=1}^kB_j\in\sB$\@.
\end{compactenum}}
Bornologies are less popular than topologies, but a treatment in some
generality can be found in~\cite{HH-N:77}\@.  There are two bornologies we
consider in this paper.  One is the \defn{compact bornology} for a
topological space $\ts{X}$ whose bounded sets are the relatively compact
sets.  The other is the \defn{von Neumann bornology} for a locally convex
topological vector space $\alg{V}$ whose bounded sets are those subsets
$\nbhd{B}\subset\alg{V}$ for which, for any neighbourhood $\nbhd{N}$ of
$0\in\alg{V}$\@, there exists $\lambda\in\realp$ such that
$\nbhd{B}\subset\lambda\nbhd{N}$\@.  On any locally convex topological vector
space we thus have these two bornologies, and generally they are not the
same.  Indeed, if $\alg{V}$ is an infinite-dimensional normed vector space,
then the compact bornology is strictly contained in the von Neumann
bornology.  We will, in fact, have occasion to use both of these bornologies,
and shall make it clear which we mean.  Now, if $(\ts{T},\sM,\mu)$ is a
measure space and if $(\ts{X},\sB)$ is a bornological space,~\ie~a set
$\ts{X}$ with a bornology $\sB$\@, a measurable map
$\map{\Psi}{\ts{T}}{\ts{X}}$ is \defn{essentially bounded} if there exists a
bounded set $B\subset\ts{X}$ such that
\begin{equation*}
\mu(\setdef{t\in\ts{T}}{\Psi(t)\not\in B})=0.
\end{equation*}
By $\L^\infty(\ts{T};\ts{X})$ we denote the set of essentially bounded maps.  If $\ts{T}=I$ is an interval in $\real$\@, a measurable map $\map{\Psi}{I}{\ts{X}}$ is \defn{locally essentially bounded} in the bornology $\sB$ if $\Psi|J$ is essentially bounded in the bornology $\sB$ for every compact subinterval $J\subset I$\@.  By $\Lloc^\infty(I;\ts{X})$ we denote the set of locally essentially bounded maps; thus the bornology is to be understood when we write expressions such as this.

\subsection*{Apologia}

This is a paper about differential geometric control theory.  It is, therefore, a paper touching upon two things,~(1)~differential geometry and~(2)~control theory.

It is our view that differential geometry \emph{is} the language of nonlinear
control theory.  As such, our attitude toward the differential geometric
aspects of what we do is unflinching in that our presentation relies,
sometimes in nontrivial ways, on all of the tools of a differential geometer,
including some that are not always a part of the nonlinear control
theoretician's tool box,~\eg~jet bundles, connections, locally convex
topologies.  In this paper, apart from presenting a new framework for control
theory, we also hope to illustrate the value of differential geometric tools
in analysing these systems, and, for that matter, any sort of geometric model
in control theory.  We have, therefore, eschewed the use of coordinates
wherever possible, since it is our opinion that unfettered coordinate
calculations are dangerous; they can lead one astray if one forgets for too
long the necessity of developing definitions and results that do not depend
on specific choices of coordinates.  Also, overuse of coordinates has a
tendency to mask structure, and it is structure that we are emphasising in
this paper.  We accept that our approach will make the paper difficult
reading for some.

This is also a paper about control theory.  And, as such, we wish to make the paper as faithful to the discipline as possible, within the confines of what we are doing.  We are certainly not including in our modelling all of the elements that would be demanded by a practicing control engineer,~\eg~no uncertainty, no robustness, no adaptive control, \etc\  And we are only considering our very limited class of models with ordinary differential equations on finite-dimensional manifolds,~\eg~no partial differential equations, no discrete-time systems, no hybrid systems, \etc\  However, with respect to those elements of control theory that we \emph{do} touch upon, we have tried to be sincere in making a framework that captures what one is likely to encounter in practice.  This means, for example, that we assiduously refrain from imposing geometric structure that is not natural from the point of view of control theory.  This tends to be a weakness of some purely differential geometric approaches to control theory, and it is a weakness that we have avoided duplicating.

\subsection*{Acknowledgements}

This research was funded in part by a grant from the Natural Sciences and Engineering Research Council of Canada.  The second author was a Visiting Professor in the Department of Mathematics at University of Hawaii, Manoa, when the paper was written, and would like to acknowledge the hospitality of the department, particularly that of Monique Chyba and George Wilkens.  The second author would also like to thank his departmental colleague Mike Roth for numerous useful conversations over the years.  While conversations with Mike did not lead directly to results in this paper, Mike's willingness to chat about complex geometry and to answer ill-informed questions was always appreciated, and ultimately very helpful.

\section{Fibre metrics for jet bundles}\label{sec:connections}

One of the principal devices we use in the paper are convenient seminorms for
the various topologies we use for spaces of sections of vector bundles.
Since such topologies rely on placing suitable norms on derivatives of
sections,~\ie~on jet bundles of vector bundles, in this section we present a
means for defining such norms, using as our starting point a pair of
connections, one for the base manifold, and one for the vector bundle.  These
allow us to provide a direct sum decomposition of the jet bundle into its
component ``derivatives,'' and so then a natural means of defining a fibre
metric for jet bundles using metrics on the tangent bundle of the base
manifold and the fibres of the vector bundle.

As we shall see, in the smooth case, these constructions are a convenience,
whereas in the real analytic case, they provide a crucial ingredient in our
global, coordinate-free description of seminorms for the topology of the
space of real analytic sections of a vector bundle.  For this reason, in this
section we shall also consider the existence of, and some properties of, real
analytic connections in vector bundles.

\subsection{A decomposition for the jet bundles of a vector bundle}\label{subsec:Jmdecomp}

We let $\map{\pi}{\man{E}}{\man{M}}$ be a smooth vector bundle with $\map{\pi_m}{\jet{m}{\man{E}}}{\man{M}}$ its $m$th jet bundle.  In a local trivialisation of $\jet{m}{\man{E}}$\@, the fibres of this vector bundle are
\begin{equation*}
\oplus_{j=0}^m\symlin[j]{\real^n}{\real^k},
\end{equation*}
with $n$ the dimension of $\man{M}$ and $k$ the fibre dimension of $\man{E}$\@.  This decomposition of the derivatives, order-by-order, that we see in the local trivialisation has no global analogue, but such a decomposition can be provided with the use of connections, and we describe how to do this.

We suppose that we have a linear connection $\nabla^0$ on the vector bundle
$\man{E}$ and an affine connection $\nabla$ on $\man{M}$\@.  We then have a
connection, that we also denote by $\nabla$\@, on $\ctb{\man{M}}$ defined by
\begin{equation*}
\lieder{Y}{\natpair{\alpha}{X}}=\natpair{\nabla_Y\alpha}{X}+
\natpair{\alpha}{\nabla_YX}.
\end{equation*}
For $\xi\in\sections[\infty]{\man{E}}$ we then have
$\nabla^0\xi\in\sections[\infty]{\ctb{\man{M}}\otimes\man{E}}$ defined by
$\nabla^0\xi(X)=\nabla^0_X\xi$ for $X\in\sections[\infty]{\tb{\man{M}}}$\@.
The connections $\nabla^0$ and $\nabla$ extend naturally to a connection,
that we denote by $\nabla^m$\@, on
$\tensor[m]{\ctb{\man{M}}}\otimes\man{E}$\@, $m\in\integerp$\@, by the
requirement that
\begin{multline*}
\nabla^m_X(\alpha^1\otimes\dots\otimes\alpha^m\otimes\xi)\\
=\sum_{j=1}^m(\alpha^1\otimes\dots\otimes(\nabla_X\alpha_j)\otimes\dots
\otimes\alpha^m\otimes\xi)+
\alpha^1\otimes\dots\otimes\alpha^m\otimes(\nabla^0_X\xi)
\end{multline*}
for $\alpha^1,\dots,\alpha^m\in\sections[\infty]{\ctb{\man{M}}}$ and
$\xi\in\sections[\infty]{\man{E}}$\@.  Note that
\begin{equation}\label{eq:nabla(j)}
\nabla^{(m)}\xi\eqdef\nabla^m(\nabla^{m-1}\cdots(\nabla^1(\nabla^0\xi)))\in
\sections[\infty]{\tensor[m+1]{\ctb{\man{M}}}\otimes\man{E}}.
\end{equation}
Now, given $\xi\in\sections[\infty]{\man{E}}$ and $m\in\integernn$\@, we
define
\begin{equation*}
P^{m+1}_{\nabla,\nabla^0}(\xi)=\Sym_{m+1}\otimes\id_{\man{E}}
(\nabla^{(m)}\xi)\in\sections[\infty]{\Symalg[m+1]{\ctb{\man{M}}}
\otimes\man{E}},
\end{equation*}
where $\map{\Sym_m}{\tensor[m]{\alg{V}}}{\Symalg[m]{\alg{V}}}$ is defined
by
\begin{equation*}
\Sym_m(v_1\otimes\dots\otimes v_m)=
\frac{1}{m!}\sum_{\sigma\in\symmgroup{m}}v_{\sigma(1)}\otimes\dots\otimes
v_{\sigma(m)}.
\end{equation*}
We take the convention that $P^0_{\nabla,\nabla^0}(\xi)=\xi$\@.

The following lemma is then key for our presentation.  While this lemma
exists in the literature in various forms, often in the form of results
concerning the extension of connections by ``bundle
functors''~\cite[\eg][Chapter~X]{IK/PWM/JS:93}\@, we were unable to find the
succinct statement we give here.  \citet{WFP:66} gives existential results
dual to what we give here, but stops short of giving an explicit formula such
as we give below.  For this reason, we give a complete proof of the lemma.
\begin{lemma}\label{lem:Jrdecomp}
The map
\begin{equation*}
\mapdef{S_{\nabla,\nabla^0}^m}{\jet{m}{\man{E}}}
{\oplus_{j=0}^m(\Symalg[j]{\ctb{\man{M}}}\otimes\man{E})}
{j_m\xi(x)}{(\xi(x),P^1_{\nabla,\nabla^0}(\xi)(x),\dots,
P^m_{\nabla,\nabla^0}(\xi)(x))}
\end{equation*}
is an isomorphism of vector bundles, and, for each\/ $m\in\integerp$\@, the
diagram
\begin{equation*}
\xymatrix{{\jet{m+1}{\man{E}}}\ar[r]^(0.3){S_{\nabla,\nabla^0}^{m+1}}
\ar[d]_{\pi^{m+1}_m}&
{\oplus_{j=0}^{m+1}(\Symalg[j]{\ctb{\man{M}}}\otimes\man{E})}
\ar[d]^{\pr^{m+1}_m}\\
{\jet{m}{\man{E}}}\ar[r]_(0.3){S_{\nabla,\nabla^0}^m}&
{\oplus_{j=0}^m(\Symalg[j]{\ctb{\man{M}}}\otimes\man{E})}}
\end{equation*}
commutes, where\/ $\pr^{m+1}_m$ is the obvious projection, stripping off the
last component of the direct sum.
\begin{proof}
We prove the result by induction on $m$\@.  For $m=0$ the result is a
tautology.  For $m=1$\@, as in \cite[\S17.1]{IK/PWM/JS:93}\@, we have a
vector bundle mapping $\map{S_{\nabla^0}}{\man{E}}{\jet{1}{\man{E}}}$ over
$\id_{\man{M}}$ that determines the connection $\nabla^0$ by
\begin{equation}\label{eq:Snablar1}
\nabla^0\xi(x)=j_1\xi(x)-S_{\nabla^0}(\xi(x)).
\end{equation}
Let us show that $S_{\nabla,\nabla^0}^1$ is well-defined.  Thus let
$\xi,\eta\in\sections[\infty]{\man{E}}$ be such that
$j_1\xi(x)=j_1\eta(x)$\@.  Then, clearly, $\xi(x)=\eta(x)$\@, and the
formula~\eqref{eq:Snablar1} shows that $\nabla\xi(x)=\nabla\eta(x)$\@, and so
$S_{\nabla,\nabla^0}^1$ is indeed well defined.  It is clearly linear on
fibres, so it remains to show that it is an isomorphism.  This will follow
from dimension counting if it is injective.  However, if
$S_{\nabla,\nabla^0}^1(j_1\xi(x))=0$ then $j_1\xi(x)=0$
by~\eqref{eq:Snablar1}\@.

For the induction step, we begin with a sublemma.
\begin{proofsublemma}
Let\/ $\alg{F}$ be a field and consider the following commutative diagram of
finite-dimensional\/ $\alg{F}$-vector spaces with exact rows and columns:
\begin{equation*}
\xymatrix{&{0}\ar[d]&{0}\ar[d]&{0}\ar[d]\\
{0}\ar[r]&{\alg{A}_1}\ar[r]^{\phi_1}\ar[d]^{\iota_1}&
{\alg{C}_1}\ar[r]^{\psi_1}\ar[d]^{\iota_2}\ar@{-->}@/^2ex/[l]^{p_1}&
{\alg{B}}\ar[r]\ar@{=}[d]\ar@{-->}@/^2ex/[l]^{\gamma_1}&{0}\\
{0}\ar[r]&{\alg{A}_2}\ar[r]^{\phi_2}\ar@/^2ex/[u]^{\sigma_1}&
{\alg{C}_2}\ar[r]^{\psi_2}\ar@/^2ex/[l]^{p_2}&
{\alg{B}}\ar@/^2ex/[l]^{\gamma_2}\ar[r]&{0}}
\end{equation*}
If there exists a mapping\/ $\gamma_2\in\Hom_{\alg{F}}(\alg{B};\alg{C}_2)$
such that\/ $\psi_2\scirc\gamma_2=\id_{\alg{B}}$ (with\/
$p_2\in\Hom_{\alg{F}}(\alg{C}_2;\alg{A}_2)$ the corresponding projection),
then there exists a unique mapping\/
$\gamma_1\in\Hom_{\alg{F}}(\alg{B};\alg{C}_1)$ such that\/
$\psi_1\scirc\gamma_1=\id_{\alg{B}}$ and such that\/
$\gamma_2=\iota_2\scirc\gamma_1$\@.  There is also induced a projection\/
$p_1\in\Hom_{\alg{F}}(\alg{C}_1;\alg{A}_1)$\@.

Moreover, if there additionally exists a mapping\/
$\sigma_1\in\Hom_{\alg{F}}(\alg{A}_2;\alg{A}_1)$ such that\/
$\sigma_1\scirc\iota_1=\id_{\alg{A}_1}$\@, then the projection\/ $p_1$ is
uniquely determined by the condition\/ $p_1=\sigma_1\scirc
p_2\scirc\iota_2$\@.
\begin{subproof}
We begin by extending the diagram to one of the form
\begin{equation*}
\xymatrix{&{0}\ar[d]&{0}\ar[d]&{0}\ar[d]\\
{0}\ar[r]&{\alg{A}_1}\ar[r]^{\phi_1}\ar[d]^{\iota_1}&
{\alg{C}_1}\ar[r]^{\psi_1}\ar[d]^{\iota_2}&
{\alg{B}}\ar[r]\ar@{=}[d]&{0}\\
{0}\ar[r]&{\alg{A}_2}\ar[r]^{\phi_2}\ar[d]^{\kappa_1}&
{\alg{C}_2}\ar[r]^{\psi_2}\ar[d]^{\kappa_2}&
{\alg{B}}\ar[r]\ar[d]&{0}\\{0}\ar[r]&{\coker(\iota_1)}
\ar@{-->}[r]^{\phi_3}\ar[d]&{\coker(\iota_2)}\ar[r]\ar[d]&{0}&\\&{0}&{0}}
\end{equation*}
also with exact rows and columns.  We claim that there is a natural mapping
$\phi_3$ between the cokernels, as indicated by the dashed arrow in the
diagram, and that $\phi_3$ is, moreover, an isomorphism.  Suppose that
$u_2\in\image(\iota_1)$ and let $u_1\in\alg{A}_1$ be such that
$\iota_1(u_1)=u_2$\@.  By commutativity of the diagram, we have
\begin{equation*}
\phi_2(u_2)=\phi_2\scirc\iota_1(u_1)=
\iota_2\scirc\phi_1(u_1),
\end{equation*}
showing that $\phi_2(\image(\iota_1))\subset\image(\iota_2)$\@.  We thus have
a well-defined homomorphism
\begin{equation*}
\mapdef{\phi_3}{\coker(\iota_1)}{\coker(\iota_2)}
{u_2+\image(\iota_1)}{\phi_2(u_2)+\image(\iota_2).}
\end{equation*}

We now claim that $\phi_3$ is injective.  Indeed,
\begin{equation*}
\phi_3(u_2+\image(\iota_1))=0\enspace\implies\enspace
\phi_2(u_2)\in\image(\iota_2).
\end{equation*}
Thus let $v_1\in\alg{C}_1$ be such that $\phi_2(u_2)=\iota_2(v_1)$\@.  Thus
\begin{align*}
&0=\psi_2\scirc\phi_2(u_2)=
\psi_2\scirc\iota_2(v_1)=\psi_1(v_1)\\
\implies\enspace&v_1\in\ker(\psi_1)=\image(\phi_1).
\end{align*}
Thus $v_1=\phi_1(u'_1)$ for some $u'_1\in\alg{A}_1$\@.  Therefore,
\begin{equation*}
\phi_2(u_2)=\iota_2\scirc\phi_1(u'_1)=
\phi_2\scirc\iota_1(u'_1),
\end{equation*}
and injectivity of $\phi_2$ gives $u_2\in\image(\iota_1)$ and so
$u_2+\image(\iota_1)=0+\image(\iota_1)$\@, giving the desired injectivity of
$\phi_3$\@.

Now note that
\begin{equation*}
\dim(\coker(\iota_1))=\dim(\alg{A}_2)-\dim(\alg{A}_1)
\end{equation*}
by exactness of the left column.  Also,
\begin{equation*}
\dim(\coker(\iota_2))=\dim(\alg{C}_2)-\dim(\alg{C}_1)
\end{equation*}
by exactness of the middle column.  By exactness of the top and middle rows,
we have
\begin{equation*}
\dim(\alg{B})=\dim(\alg{C}_2)-\dim(\alg{A}_2)=
\dim(\alg{C}_1)-\dim(\alg{A}_1).
\end{equation*}
This proves that 
\begin{equation*}
\dim(\coker(\iota_1))=\dim(\coker(\iota_2)).
\end{equation*}
Thus the homomorphism $\phi_3$ is an isomorphism, as claimed.

Now we proceed with the proof, using the extended diagram, and identifying
the bottom cokernels with the isomorphism $\phi_3$\@.  The existence of the
stated homomorphism $\gamma_2$ means that the middle row in the diagram
splits.  Therefore, $\alg{C}_2=\image(\phi_2)\oplus\image(\gamma_2)$\@.  Thus
there exists a well-defined projection
$p_2\in\Hom_{\alg{F}}(\alg{C}_2;\alg{A}_2)$ such that
$p_2\scirc\phi_2=\id_{\alg{A}_2}$~\cite[Theorem~41.1]{PRH:74b}\@.

We will now prove that $\image(\gamma_2)\subset\image(\iota_2)$\@.  By
commutativity of the diagram and since $\psi_1$ is surjective, if
$w\in\alg{B}$ then there exists $v_1\in\alg{C}_1$ such that
$\psi_2\scirc\iota_2(v_1)=w$\@.  Since
$\psi_2\scirc\gamma_2=\id_{\alg{B}}$\@, we have
\begin{equation*}
\psi_2\scirc\iota_2(v_1)=\psi_2\scirc\gamma_2(w)\enspace\implies\enspace
\iota_2(v_1)-\gamma_2(w)\in\ker(\psi_2)=\image(\phi_2).
\end{equation*}
Let $u_2\in\alg{A}_2$ be such that $\phi_2(u_2)=\iota_2(v_1)-\gamma_2(w)$\@.
Since $p_2\scirc\phi_2=\id_{\alg{A}_2}$ we have
\begin{equation*}
u_2=p_2\scirc\iota_2(v_1)-p_2\scirc\gamma_2(w),
\end{equation*}
whence
\begin{equation*}
\kappa_1(u_2)=\kappa_1\scirc p_2\scirc\iota_2(v_1)-
\kappa_1\scirc p_2\scirc\gamma_2(w)=0,
\end{equation*}
noting that~(1)~$\kappa_1\scirc p_2=\kappa_2$ (by
commutativity),~(2)~$\kappa_2\scirc\iota_2=0$ (by exactness),
and~(3)~$p_2\scirc\gamma_2=0$ (by exactness).  Thus
$u_2\in\ker(\kappa_1)=\image(\iota_1)$\@.  Let $u_1\in\alg{A}_1$ be such that
$\iota_1(u_1)=u_2$\@.  We then have
\begin{equation*}
\iota_2(v_1)-\gamma_2(w)=\phi_2\scirc\iota_1(u_1)=\iota_2\scirc\phi_1(u_1),
\end{equation*}
which gives $\gamma_2(w)\in\image(\iota_2)$\@, as claimed.

Now we define $\gamma_1\in\Hom_{\alg{F}}(\alg{B};\alg{C}_1)$ by asking that
$\gamma_1(w)\in\alg{C}_1$ have the property that
$\iota_2\scirc\gamma_1(w)=\gamma_2(w)$\@, this making sense since we just
showed that $\image(\gamma_2)\subset\image(\iota_2)$\@.  Moreover, since
$\iota_2$ is injective, the definition uniquely prescribes $\gamma_1$\@.
Finally we note that
\begin{equation*}
\psi_1\scirc\gamma_1=\psi_2\scirc\iota_2\scirc\gamma_1=\psi_2\scirc\gamma_2
=\id_{\alg{B}},
\end{equation*}
as claimed.

To prove the final assertion, let us denote $\hat{p}_1=\sigma_1\scirc
p_2\scirc\iota_2$\@.  We then have
\begin{equation*}
\hat{p}_1\scirc\phi_1=\sigma_1\scirc p_2\scirc\iota_2\scirc\phi_1
=\sigma_1\scirc p_2\scirc\phi_2\scirc\iota_1
=\sigma_1\scirc\iota_1=\id_{\alg{A}_1},
\end{equation*}
using commutativity.  We also have
\begin{equation*}
\hat{p}_1\scirc\gamma_1=\sigma_1\scirc p_2\scirc\iota_2\scirc\gamma_1=
\sigma_1\scirc p_2\scirc\gamma_2=0.
\end{equation*}
The two preceding conclusions show that $\hat{p}_1$ is the projection defined
by the splitting of the top row of the diagram,~\ie~$\hat{p}_1=p_1$\@.
\end{subproof}
\end{proofsublemma}

Now suppose that the lemma is true for $m\in\integerp$\@.  For any
$k\in\integerp$ we have a short exact sequence
\begin{equation*}
\xymatrix{{0}\ar[r]&{\Symalg[k]{\ctb{\man{M}}}\otimes\man{E}}
\ar[r]^(0.65){\epsilon_k}&{\jet{k}{\man{E}}}\ar[r]^(0.45){\pi^k_{k-1}}&
{\jet{k-1}{\man{E}}}\ar[r]&{0}}
\end{equation*}
for which we refer to~\cite[Theorem~6.2.9]{DJS:89}\@.  Recall
from~\cite[Definition~6.2.25]{DJS:89} that we have an inclusion $\iota_{1,m}$
of $\jet{m+1}{\man{E}}$ in $\jet{1}{(\jet{m}{\man{E}})}$ by
$j_{m+1}\xi(x)\mapsto j_1(j_m\xi(x))$\@.  We also have an induced injection
\begin{equation*}
\map{\hat{\iota}_{1,m}}{\Symalg[m+1]{\ctb{\man{M}}}\otimes\man{E}}
{\ctb{\man{M}}\otimes\jet{m}{\man{E}}}
\end{equation*}
defined by the composition
\begin{equation*}
\xymatrix{{\Symalg[m+1]{\ctb{\man{M}}}\otimes\man{E}}\ar[r]&
{\ctb{\man{M}}\otimes\Symalg[m]{\ctb{\man{M}}}\otimes\man{E}}
\ar[r]^(0.6){\id\otimes\epsilon_m}&{\ctb{\man{M}}\otimes\jet{m}{\man{E}}}}
\end{equation*}
Explicitly, the left arrow is defined by
\begin{equation*}
\alpha^1\odot\dots\odot\alpha^{m+1}\otimes\xi\mapsto
\sum_{j=1}^{m+1}\alpha^j\otimes\alpha^1\odot\dots\odot
\alpha^{j-1}\odot\alpha^{j+1}\odot\dots\odot\alpha^{m+1}\otimes\xi,
\end{equation*}
$\odot$ denoting the symmetric tensor product defined by
\begin{equation}\label{eq:odotdef}
A\odot B=\sum_{\sigma\in\symmgroup{k,l}}\sigma(A\otimes B),
\end{equation}
for $A\in\Symalg[k]{\alg{V}}$ and $B\in\Symalg[l]{\alg{V}}$\@, and with
$\symmgroup{k,l}$ the subset of $\symmgroup{k+l}$ consisting of permutations
$\sigma$ satisfying
\begin{equation*}
\sigma(1)<\dots<\sigma(k),\quad\sigma(k+1)<\dots<\sigma(k+l).
\end{equation*}
We thus have the following commutative diagram with exact rows and columns:
\begin{equation}\label{eq:J1Jmdiagram}
\xymatrix{&{0}\ar[d]&{0}\ar[d]&{0}\ar[d]&\\
{0}\ar[r]&{\Symalg[m+1]{\ctb{\man{M}}}\otimes\man{E}}
\ar[r]^(0.6){\epsilon_{m+1}}\ar[d]^{\hat{\iota}_{1,m}}&
{\jet{m+1}{\man{E}}}\ar[r]^(0.55){\pi^{m+1}_m}\ar[d]^{\iota_{1,m}}
\ar@/^2ex/@{-->}[l]^(0.45){P_{m+1}}&
{\jet{m}{\man{E}}}\ar[r]\ar@{=}[d]\ar@/^2ex/@{-->}[l]^(0.45){\Gamma_{m+1}}
&{0}\\{0}\ar[r]&{\ctb{\man{M}}\otimes\jet{m}{\man{E}}}
\ar[r]^(0.55){\epsilon_{1,m}}\ar@/^2ex/@{-->}[u]^{\lambda_{1,m}}&
{\jet{1}{(\jet{m}{\man{E}})}}\ar[r]^(0.55){(\pi_m)_1}
\ar@/^2ex/@{-->}[l]^(0.45){P_{1,m}}&
{\jet{m}{\man{E}}}\ar[r]\ar@/^2ex/@{-->}[l]^(0.45){\Gamma_{1,m}}&{0}}
\end{equation}
We shall define a connection on
$\map{(\pi_m)_1}{\jet{1}{(\jet{m}{\man{E}})}}{\jet{m}{\man{E}}}$ which gives
a splitting $\Gamma_{1,m}$ and $P_{1,m}$ of the lower row in the diagram.  By
the sublemma, this will give a splitting $\Gamma_{m+1}$ and $P_{m+1}$ of the
upper row, and so give a projection from $\jet{m+1}{\man{E}}$ onto
$\Symalg[m+1]{\ctb{\man{M}}}\otimes\man{E}$\@, which will allow us to prove
the induction step.  To compute $P_{m+1}$ from the sublemma, we shall also
give a map $\lambda_{1,m}$ as in the diagram so that
$\lambda_{1,m}\scirc\hat{\iota}_{1,m}$ is the identity on
$\Symalg[m+1]{\ctb{\man{M}}}\otimes\man{E}$\@.

We start, under the induction hypothesis, by making the identification
\begin{equation*}
\jet{m}{\man{E}}\simeq\oplus_{j=0}^m\Symalg[j]{\ctb{\man{M}}}\otimes\man{E},
\end{equation*}
and consequently writing a section of $\jet{m}{\man{E}}$ as
\begin{equation*}
x\mapsto(\xi(x),P^1_{\nabla,\nabla^0}(\xi(x)),\dots,
P^m_{\nabla,\nabla^0}(\xi(x))).
\end{equation*}
We then have a connection $\ol{\nabla}^m$ on $\jet{m}{\man{E}}$ given by
\begin{equation*}
\ol{\nabla}^m_X(\xi,P^1_{\nabla,\nabla^0}(\xi),\dots,
P^m_{\nabla,\nabla^0}(\xi))
=(\nabla^0_X\xi,\nabla^1_XP^1_{\nabla,\nabla^0}(\xi),\dots,
\nabla^m_XP^m_{\nabla,\nabla^0}(\xi)).
\end{equation*}
Thus
\begin{equation*}
\ol{\nabla}^m(\xi,P^1_{\nabla,\nabla^0}(\xi),\dots,
P^m_{\nabla,\nabla^0}(\xi))=
(\nabla^0\xi,\nabla^1P^1_{\nabla,\nabla^0}(\xi),\dots,
\nabla^mP^m_{\nabla,\nabla^0}(\xi)),
\end{equation*}
which\textemdash{}according to the jet bundle characterisation of connections
from~\cite[\S17.1]{IK/PWM/JS:93} and which we have already employed
in~\eqref{eq:Snablar1}\textemdash{}gives the mapping $P_{1,m}$ in the
diagram~\eqref{eq:J1Jmdiagram} as
\begin{equation*}
P_{1,m}(j_1(\xi,P^1_{\nabla,\nabla^0}(\xi),\dots,
P^m_{\nabla,\nabla^0}(\xi)))=
(\nabla^0\xi,\nabla^1P^1_{\nabla,\nabla^0}(\xi),\dots,
\nabla^mP^m_{\nabla,\nabla^0}(\xi)).
\end{equation*}

Now we define a mapping $\lambda_{1,m}$ for which
$\lambda_{1,m}\scirc\hat{\iota}_{1,m}$ is the identity on
$\Symalg[m+1]{\ctb{\man{M}}}\otimes\man{E}$\@.  We continue to use the
induction hypothesis in writing elements of $\jet{m}{\man{E}}$\@, so that we
consider elements of $\ctb{\man{M}}\otimes\jet{m}{\man{E}}$ of the form
\begin{equation*}
(\alpha\otimes\xi,\alpha\otimes A_1,\dots,\alpha\otimes A_m),
\end{equation*}
for $\alpha\in\ctb{\man{M}}$ and
$A_k\in\Symalg[k]{\ctb{\man{M}}}\otimes\man{E}$\@, $k\in\{1,\dots,m\}$\@.
We then define $\lambda_{1,m}$ by
\begin{multline*}
\lambda_{1,m}(\alpha_0\otimes\xi,\alpha_0\otimes\alpha_1^1\otimes\xi,\dots,
\alpha_0\otimes\alpha_m^1\odot\dots\odot\alpha_m^m\otimes\xi)\\
=\Sym_{m+1}\otimes\id_{\man{E}}
(\alpha_0\otimes\alpha_m^1\odot\dots\odot\alpha_m^m\otimes\xi).
\end{multline*}
Note that, with the form of $\jet{m}{\man{E}}$ from the induction hypothesis,
we have
\begin{multline*}
\hat{\iota}_{1,m}(\alpha^1\odot\dots\odot\alpha^{m+1}\otimes\xi)\\
=\Bigl(0,\dots,0,\frac{1}{m+1}\sum_{j=1}^{m+1}
\alpha^j\otimes\alpha^1\odot\dots\odot
\alpha^{j-1}\odot\alpha^{j+1}\odot\dots\odot\alpha^{m+1}\otimes\xi\Bigr).
\end{multline*}
We then directly verify that $\lambda_{1,m}\scirc\hat{\iota}_{1,m}$ is indeed
the identity.

We finally claim that
\begin{equation}\label{eq:Pm+1def}
P_{m+1}(j_{m+1}\xi(x))=P^{m+1}_{\nabla,\nabla^0}(\xi),
\end{equation}
which will establish the lemma.  To see this, first note that it suffices to
define $P_{m+1}$ on $\image(\epsilon_{m+1})$ since
\begin{compactenum}
\item $\jet{m+1}{\man{E}}\simeq(\Symalg[m+1]{\ctb{\man{M}}}\otimes\man{E})
\oplus\jet{m}{\man{E}}$\@,
\item $P_{m+1}$ is zero on $\jet{m}{\man{E}}\subset\jet{m+1}{\man{E}}$
(thinking of the inclusion arising from the connection-induced isomorphism
from the preceding item), and
\item $P_{m+1}\scirc\epsilon_{m+1}$ is the identity map on
$\Symalg[m+1]{\ctb{\man{M}}}\otimes\man{E}$\@.
\end{compactenum}
In order to connect the algebra and the geometry, let us write elements of
$\Symalg[m+1]{\ctb{\man{M}}}\otimes\man{E}$ in a particular way.  We let
$x\in\man{M}$ and let $f^1,\dots,f^{m+1}$ be smooth functions contained in
the maximal ideal of $\func[\infty]{\man{M}}$ at $x$\@,~\ie~$f^j(x)=0$\@,
$j\in\{1,\dots,m+1\}$\@.  Let $\xi$ be a smooth section of $\man{E}$\@.  We
then can work with elements of $\Symalg[m+1]{\ctb{\man{M}}}\otimes\man{E}$
of the form
\begin{equation*}
\d{f^1}(x)\odot\dots\odot\d{f^{m+1}}(x)\otimes\xi(x).
\end{equation*}
We then have
\begin{equation*}
\epsilon_{m+1}(\d{f^1}(x)\odot\dots\odot\d{f^{m+1}}(x)\otimes\xi(x))=
j_{m+1}(f^1\cdots f^{m+1}\xi)(x);
\end{equation*}
this is easy to see using the Leibniz Rule~\cite[\cf][Lemma~2.1]{HLG:67a}\@.
(See~\cite[Supplement~2.4A]{RA/JEM/TSR:88} for a description of the
higher-order Leibniz Rule.)  Now, using the last part of the sublemma, we
compute
{\small\begin{align*}
P_{m+1}(&j_{m+1}(f^1\cdots f^{m+1}\xi)(x))\\
=&\;\lambda_{1,m}\scirc P_{1,m}\scirc\iota_{1,m}
(j_{m+1}(f^1\cdots f^{m+1}\xi)(x))\\
=&\;\lambda_{1,m}\scirc P_{1,m}(j_1(f^1\cdots f^{m+1}\xi,
P^1_{\nabla,\nabla^0}(f^1\cdots f^{m+1}\xi),\dots,
P^m_{\nabla,\nabla^0}(f^1\cdots f^{m+1}\xi))(x))\\
=&\;\lambda_{1,m}(\nabla^0(f^1\cdots f^{m+1}\xi)(x),
\nabla^1P^1_{\nabla,\nabla^0}(f^1\cdots f^{m+1}\xi)(x),\dots,
\nabla^mP^m_{\nabla,\nabla^0}(f^1\cdots f^{m+1}\xi)(x))\\
=&\;\Sym_{m+1}\otimes\id_{\man{E}} 
(\nabla^mP^m_{\nabla,\nabla^0}(f^1\cdots f^{m+1}\xi)(x))\\
=&\;P^{m+1}_{\nabla,\nabla_0}(f^1\cdots f^{m+1}\xi)(x),
\end{align*}}%
which shows that, with $P_{m+1}$ defined as in~\eqref{eq:Pm+1def}\@,
$P_{m+1}\scirc\epsilon_{m+1}$ is indeed the identity on
$\Symalg[m+1]{\ctb{\man{M}}}\otimes\man{E}$\@.

The commuting of the diagram in the statement of the lemma follows directly
from the recursive nature of the constructions.
\end{proof}
\end{lemma}

\subsection{Fibre metrics using jet bundle decompositions}\label{subsec:olGm}

We also require the following result concerning inner products on tensor
products.
\begin{lemma}\label{lem:inprodotimes}
Let\/ $\alg{U}$ and\/ $\alg{V}$ be finite-dimensional\/ $\real$-vector spaces
and let\/ $\metric$ and\/ $\hmetric$ be inner products on\/ $\alg{U}$ and\/
$\alg{V}$\@, respectively.  Then the element\/ $\metric\otimes\hmetric$ of\/
$\tensor[2]{\dual{\alg{U}}\otimes\dual{\alg{V}}}$ defined by
\begin{equation*}
\metric\otimes\hmetric(u_1\otimes v_1,u_2\otimes v_2)=
\metric(u_1,u_2)\hmetric(v_1,v_2)
\end{equation*}
is an inner product on\/ $\alg{U}\otimes\alg{V}$\@.
\begin{proof}
Let $\ifam{e_1,\dots,e_m}$ and $\ifam{f_1,\dots,f_n}$ be orthonormal bases
for $\alg{U}$ and $\alg{V}$\@, respectively.  Then
\begin{equation}\label{eq:GHbasis}
\setdef{e_a\otimes f_j}{a\in\{1,\dots,m\},\ j\in\{1,\dots,n\}}
\end{equation}
is a basis for $\alg{U}\otimes\alg{V}$\@.  Note that
\begin{equation*}
\metric\otimes\hmetric(e_a\otimes f_j,e_b\otimes f_k)=
\metric(e_a,e_b)\hmetric(f_j,f_k)=\delta_{ab}\delta_{jk},
\end{equation*}
which shows that $\metric\otimes\hmetric$ is indeed an inner product,
as~\eqref{eq:GHbasis} is an orthonormal basis.
\end{proof}
\end{lemma}

Now, we let $\metric_0$ be a fibre metric on $\man{E}$ and let $\metric$ be a
Riemannian metric on $\man{M}$\@.  Let us denote by $\metric^{-1}$ the
associated fibre metric on $\ctb{\man{M}}$ defined by
\begin{equation*}
\metric^{-1}(\alpha_x,\beta_x)=
\metric(\metric^\sharp(\alpha_x),\metric^\sharp(\beta_x)).
\end{equation*}
By induction using the preceding lemma, we have a fibre metric $\metric_j$ on
$\tensor[j]{\ctb{\man{M}}}\otimes\man{E}$ induced by $\metric^{-1}$ and
$\metric_0$\@.  By restriction, this gives a fibre metric on
$\Symalg[j]{\ctb{\man{M}}}\otimes\man{E}$\@.  We can thus define a fibre
metric $\ol{\metric}_m$ on $\jet{m}{\man{E}}$ given by
\begin{equation*}
\ol{\metric}_m(j_m\xi(x),j_m\eta(x))=
\sum_{j=0}^m\metric_j\Bigl(\frac{1}{j!}P^j_{\nabla,\nabla^0}(\xi)(x),
\frac{1}{j!}P^j_{\nabla,\nabla^0}(\eta)(x)\Bigr),
\end{equation*}
with the convention that $\nabla^{(-1)}\xi=\xi$\@.  Associated to this inner
product on fibres is the norm on fibres, which we denote by
$\dnorm{\cdot}_{\ol{\metric}_m}$\@.  We shall use these fibre norms
continually in our descriptions of our various topologies below.

\subsection{Real analytic connections}\label{subsec:analytic-conn}

The fibre metrics from the preceding section will be used to define seminorms
for spaces of sections of vector bundles.  In the finitely differentiable and
smooth cases, the particular fibre metrics we define above are not really
required to give seminorms for the associated topologies: any fibre metrics
on the jet bundles will suffice.  Indeed, as long as one is only working with
finitely many derivatives at one time, the choice of fibre norms on jet
bundles is of no consequence, since different choices will be equivalent on
compact subsets of $\man{M}$\@,~\cf~Section~\ref{subsec:COinfty-vb}\@.
However, when we work with the real analytic topology, we are no longer
working only with finitely many derivatives, but with the infinite jet of a
section.  For this reason, different choices of fibre metric for jet bundles
may give rise to different topologies for the space of real analytic
sections, unless the behaviour of the fibre metrics is compatible as the
order of derivatives goes to infinity.  In this section we give a fundamental
inequality for our fibre metrics of Section~\ref{subsec:olGm} in the real
analytic case that ensures that they, in fact, describe the real analytic
topology.

First let us deal with the matter of existence of real analytic data defining
these fibre metrics.
\begin{lemma}\label{lem:analytic-conn}
If\/ $\map{\pi}{\man{E}}{\man{M}}$ is a real analytic vector bundle, then
there exist
\begin{compactenum}[(i)]
\item a real analytic linear connection on\/ $\man{E}$\@,
\item a real analytic affine connection on\/ $\man{M}$\@,
\item a real analytic fibre metric on\/ $\man{E}$\@, and
\item a real analytic Riemannian metric on\/ $\man{M}$\@.
\end{compactenum}
\begin{proof}
By \cite[Theorem~3]{HG:58}\@, there exists a proper real analytic embedding
$\iota_{\man{E}}$ of $\man{E}$ in $\real^N$ for some $N\in\integerp$\@.
There is then an induced proper real analytic embedding $\iota_{\man{M}}$ of
$\man{M}$ in $\real^N$ by restricting $\iota_{\man{E}}$ to the zero section
of $\man{E}$\@.  Let us take the subbundle $\hat{\man{E}}$ of
$\tb{\real^N}|\iota_{\man{M}}(\man{M})$ whose fibre at
$\iota_{\man{M}}(x)\in\iota_{\man{M}}(\man{M})$ is
\begin{equation*}
\hat{\man{E}}_{\iota_{\man{M}}(x)}=
\tf[0_x]{\iota_{\man{E}}}(\vb[0_x]{\man{E}}).
\end{equation*}
Now recall that $\man{E}\simeq\zeta^*\vb{\man{E}}$\@, where
$\map{\zeta}{\man{M}}{\man{E}}$ is the zero
section~\cite[page~55]{IK/PWM/JS:93}\@.  Let us abbreviate
$\hat{\iota}_{\man{E}}=\tf{\iota_{\man{E}}}|\zeta^*\vb{\man{E}}$\@.  We then
have the following diagram
\begin{equation}\label{eq:EMembedding}
\xymatrix{{\man{E}\simeq\zeta^*\vb{\man{E}}}
\ar[d]_{\pi}\ar[r]^{\hat{\iota}_{\man{E}}}&
{\real^N\times\real^N}\ar[d]^{\pr_2}\\{\man{M}}\ar[r]_{\iota_{\man{M}}}&
{\real^N}}
\end{equation}
describing a monomorphism of real analytic vector bundles over the proper
embedding $\iota_{\man{M}}$\@, with the image of $\hat{\iota}_{\man{E}}$
being $\hat{\man{E}}$\@.

Among the many ways to prescribe a linear connection on the vector bundle
$\man{E}$\@, we will take the prescription whereby one defines a mapping
$\map{K}{\tb{\man{E}}}{\man{E}}$ such that the two diagrams
\begin{equation}\label{eq:Kdiagrams}
\xymatrix{{\tb{\man{E}}}\ar[r]^{K}\ar[d]_{\tf{\pi}}&{\man{E}}\ar[d]^{\pi}\\
{\tb{\man{M}}}\ar[r]_{\tbproj{\man{M}}}&{\man{M}}}\qquad
\xymatrix{{\tb{\man{E}}}\ar[r]^{K}\ar[d]_{\tbproj{\man{E}}}&
{\man{E}}\ar[d]^{\pi}\\{\man{E}}\ar[r]_{\pi}&{\man{M}}}\qquad
\end{equation}
define vector bundle mappings~\cite[\S11.11]{IK/PWM/JS:93}\@.  We define $K$
as follows.  For $e_x\in\man{E}_x$ and $X_{e_x}\in\tb[e_x]{\man{E}}$ we have
\begin{equation*}
\tf[e_x]{\hat{\iota}_{\man{E}}}(X_{e_x})\in
\tb[\hat{\iota}_{\man{E}}(e_x)]{(\real^N\times\real^N)}\simeq
\real^N\oplus\real^N,
\end{equation*}
and we define $K$ so that
\begin{equation*}
\hat{\iota}_{\man{E}}\scirc K(X_{e_x})=
\pr_2\scirc\tf[e_x]{\hat{\iota}_{\man{E}}}(X_{e_x});
\end{equation*}
this uniquely defines $K$ by injectivity of $\hat{\iota}_{\man{E}}$\@, and
amounts to using on $\man{E}$ the connection induced on
$\image(\hat{\iota}_{\man{E}})$ by the trivial connection on
$\real^N\times\real^N$\@.  In particular, this means that we think of
$\hat{\iota}_{\man{E}}\scirc K(X_{e_x})$ as being an element of the fibre of
the trivial bundle $\real^N\times\real^N$ at $\iota_{\man{M}}(x)$\@.

If $v_x\in\tb{\man{M}}$\@, if $e,e'\in\man{E}$\@, and if
$X\in\tb[e]{\man{E}}$ and $X'\in\man{E}_{e'}$ satisfy
$X,X'\in\tf{\pi}^{-1}(v_x)$\@, then note that
\begin{align*}
\tf[e]{\pi}(X)=\tf[e']{\pi}(X')\ \implies&\
\tf[e]{(\iota_{\man{M}}\scirc\pi)}(X)=\tf[e']{(\iota_{\man{M}}\scirc\pi)}(X')\\
\implies&\ \tf[e]{(\pr_2\scirc\hat{\iota}_{\man{E}})}(X)=
\tf[e']{(\pr_2\scirc\hat{\iota}_{\man{E}})}(X')\\
\implies&\ \tf[\iota_{\man{M}}(x)]{\pr_2}\scirc\tf[e]{\hat{\iota}_{\man{E}}}(X)=
\tf[\iota_{\man{M}}(x)]{\pr_2}\scirc\tf[e']{\hat{\iota}_{\man{E}}}(X').
\end{align*}
Thus we can write
\begin{equation*}
\tf[e]{\hat{\iota_{\man{E}}}}(X)=(\vect{x},\vect{e},\vect{u},\vect{v}),\quad
\tf[e']{\hat{\iota_{\man{E}}}}(X)=(\vect{x},\vect{e}',\vect{u},\vect{v}')
\end{equation*}
for suitable
$\vect{x},\vect{u},\vect{e},\vect{e}',\vect{v},\vect{v}'\in\real^N$\@.
Therefore,
\begin{equation*}
\hat{\iota}_{\man{E}}\scirc K(X)=(\vect{x},\vect{v}),\quad
\hat{\iota}_{\man{E}}\scirc K(X')=(\vect{x},\vect{v}'),\quad
\hat{\iota}_{\man{E}}\scirc K(X+X')=(\vect{x},\vect{v}+\vect{v}'),
\end{equation*}
from which we immediately conclude that, for addition in the vector bundle
$\map{\tf{\pi}}{\tb{\man{E}}}{\tb{\man{M}}}$\@, we have
\begin{equation*}
\hat{\iota}_{\man{E}}\scirc K(X+X')=
\hat{\iota}_{\man{E}}\scirc K(X)+\hat{\iota}_{\man{E}}\scirc K(X'),
\end{equation*}
showing that the diagram on the left in~\eqref{eq:Kdiagrams} makes $K$ a
vector bundle mapping.

On the other hand, if $e_x\in\man{E}$ and if $X,X'\in\tb[e_x]{\man{E}}$\@,
then we have, using vector bundle addition in $\map{\tbproj{\man{E}}}{\tb{\man{E}}}{\man{E}}$\@,
\begin{align*}
\hat{\iota}_{\man{E}}\scirc K(X+X')=&\;
\pr_2\scirc\tf[e_x]{\hat{\iota}_{\man{E}}}(X+X')\\
=&\;\pr_2\scirc\tf[e_x]{\hat{\iota_{\man{E}}}}(X)+
\pr_2\scirc\tf[e_x]{\hat{\iota_{\man{E}}}}(X')\\
=&\;\hat{\iota}_{\man{E}}\scirc K(X)+\hat{\iota}_{\man{E}}\scirc K(X'),
\end{align*}
giving that the diagram on the right in~\eqref{eq:Kdiagrams} make $K$ a
vector bundle mapping.  Since $K$ is real analytic, this defines a real
analytic linear connection $\nabla^0$ on $\man{E}$ as
in~\cite[\S11.11]{IK/PWM/JS:93}\@.

The existence of $\metric_0$\@, $\metric$\@, and $\nabla$ are
straightforward.  Indeed, we let $\metric_{\real^N}$ be the Euclidean metric
on $\real^N$\@, and define $\metric_0$ and $\metric$ by
\begin{equation*}
\metric_0(e_x,e'_x)=\metric_{\real^N}(\hat{\iota}_{\man{E}}(e_x),
\hat{\iota}_{\man{E}}(e'_x))
\end{equation*}
and
\begin{equation*}
\metric(v_x,v'_x)=\metric_{\real^N}(\tf[x]{\iota_{\man{M}}}(v_x),
\tf[x]{\iota_{\man{M}}}(v'_x)).
\end{equation*}
The affine connection $\nabla$ can be taken to be the Levi-Civita connection
of $\metric$\@.
\end{proof}
\end{lemma}

The existence of a real analytic linear connection in a real analytic vector
bundle is asserted at the bottom of page~302 in~\cite{AK/PWM:97}\@, and we
fill in the blanks in the preceding proof.

Now let us provide a fundamental relationship between the geometric fibre
norms of Section~\ref{subsec:olGm} and norms constructed in local coordinate
charts.
\begin{lemma}\label{lem:pissy-estimate}
Let\/ $\nbhd{U}\subset\real^n$ be open, denote\/
$\real^k_{\nbhd{U}}=\nbhd{U}\times\real^k$\@, let\/ $K\subset\nbhd{U}$ be
compact, and consider the trivial vector bundle\/
$\map{\pr_1}{\real^k_{\nbhd{U}}}{\nbhd{U}}$\@.  Let\/ $\metric$ be a
Riemannian metric on\/ $\nbhd{U}$\@, let\/ $\metric_0$ be a vector bundle
metric on\/ $\real^k_{\nbhd{U}}$\@, let\/ $\nabla$ be an affine connection
on\/ $\nbhd{U}$\@, and let\/ $\nabla^0$ be a vector bundle connection on\/
$\real^k_{\nbhd{U}}$\@, with all of these being real analytic.  Then there
exist\/ $C,\sigma\in\realp$ such that
\begin{equation*}
\frac{\sigma^m}{C}\dnorm{j_m\vect{\xi}(\vect{x})}_{\ol{\metric}_m}\le
\sup\Bigsetdef{\frac{1}{I!}\snorm{\linder[I]{\xi^a}(\vect{x})}}
{\snorm{I}\le m,\ a\in\{1,\dots,k\}}\le
\frac{C}{\sigma^m}\dnorm{j_m\vect{\xi}(\vect{x})}_{\ol{\metric}_m}
\end{equation*}
for every\/ $\vect{\xi}\in\sections[\infty]{\real^k_{\nbhd{U}}}$\@,\/
$\vect{x}\in K$\@, and\/ $m\in\integernn$\@.
\begin{proof}
We begin the proof with a series of sublemmata of a fairly technical nature.
From these the lemma will follow in a more or less routine manner.

Let us first prove a result which gives a useful local trivialisation of a
vector bundle and a corresponding Taylor expansion for real analytic
sections.
\begin{proofsublemma}\label{psublem:local-nabla}
Let\/ $\map{\pi}{\man{E}}{\man{M}}$ be a real analytic vector bundle, let\/
$\nabla^0$ be a real analytic linear connection on\/ $\man{E}$\@, and let\/
$\nabla$ be a real analytic affine connection on\/ $\man{M}$\@.  Let\/
$x\in\man{M}$\@, and let\/ $\nbhd{N}\subset\tb[x]{\man{M}}$ be a convex
neighbourhood of\/ $0_x$ and\/ $\nbhd{V}\subset\man{M}$ be a neighbourhood
of\/ $x$ such that the exponential map\/ $\exp_x$ corresponding to\/ $\nabla$
is a real analytic diffeomorphism from\/ $\nbhd{N}$ to\/ $\nbhd{V}$\@.
For\/ $y\in\nbhd{V}$\@, let\/ $\map{\tau_{xy}}{\man{E}_x}{\man{E}_y}$ be
parallel transport along the geodesic\/ $t\mapsto\exp_x(t\exp_x^{-1}(y))$\@.  Define
\begin{equation*}
\mapdef{\kappa_x}{\nbhd{N}\times\man{E}_x}{\man{E}|\nbhd{V}}
{(v,e_x)}{\tau_{x,\exp_x(v)}(e_x).}
\end{equation*}
Then
\begin{compactenum}[(i)]
\item \label{pl:kappax1} $\kappa_x$ is a real analytic vector bundle
isomorphism over\/ $\exp_x$ and
\item \label{pl:kappax2} if\/ $\xi\in\sections[\omega]{\man{E}|\nbhd{V}}$\@,
then
\begin{equation*}
\kappa_x^{-1}\scirc\xi\scirc\exp_x(v)=
\sum_{m=0}^\infty\frac{1}{m!}\nabla^{(m-1)}\xi(x)
(\underbrace{v,\cdots,v}_{m\ \textrm{times}})
\end{equation*}
for\/ $v$ in a sufficiently small neighbourhood of\/
$0_x\in\tb[x]{\man{M}}$\@.
\end{compactenum}
\begin{subproof}
\eqref{pl:kappax1} Consider the vector field $X_{\nabla,\nabla^0}$ on the
Whitney sum $\tb{\man{M}}\oplus\man{E}$ defined by
\begin{equation*}
X_{\nabla,\nabla^0}(v_x,e_x)=\hlft(v_x,v_x)\oplus\hlft_0(e_x,v_x),
\end{equation*}
where $\hlft(v_x,u_x)$ is the horizontal lift of $u_x\in\tb[x]{\man{M}}$ to
$\tb[v_x]{\tb{\man{M}}}$ and $\hlft_0(e_x,u_x)$ is the horizontal lift of
$u_x\in\tb[x]{\man{M}}$ to $\tb[e_x]{\man{E}}$\@.  Note that, since
\begin{equation*}
\tf{\tbproj{\man{M}}}(\hlft(v_x,v_x))=
\tf{\pi}(\hlft_0(e_x,v_x)),
\end{equation*}
this is indeed a vector field on $\tb{\man{M}}\oplus\man{E}$\@.  Moreover,
the integral curve of $X_{\nabla,\nabla^0}$ through $(v_x,e_x)$ is
$t\mapsto\gamma'(t)\oplus\tau(t)$\@, where $\gamma$ is the geodesic with
initial condition $\gamma'(0)=v_x$ and where $t\mapsto\tau(t)$ is parallel
transport of $e_x$ along $\gamma$\@.  This is a real analytic vector field,
and so the flow depends in a real analytic manner on initial
condition~\cite[Proposition~C.3.12]{EDS:98}\@.  In particular, it depends in
a real analytic manner on initial conditions lying in
$\nbhd{N}\times\man{E}_x$\@.  But, in this case, the map from initial
condition to value at $t=1$ is exactly $\kappa_x$\@.  This shows that
$\kappa_x$ is indeed real analytic.  Moreover, it is clearly fibre preserving
over $\exp_x$ and is linear on fibres, and so is a vector bundle map~\cite[\cf][Proposition~3.4.12(iii)]{RA/JEM/TSR:88}\@.

\eqref{pl:kappax2} For $v\in\nbhd{N}$\@, let $\gamma_v$ be the geodesic
satisfying $\gamma_v'(0)=v$\@.  Then, for $t\in\realp$ satisfying
$\snorm{t}\le1$\@, define
\begin{equation*}
\alpha_v(t)=\kappa_x^{-1}\scirc\xi(\gamma_v(t))=
\tau_{x,\gamma_v(t)}^{-1}(\xi(\gamma_v(t)).
\end{equation*}
We compute derivatives of $\alpha_v$ as follows, by induction and using the
fact that $\nabla_{\gamma_v'(t)}\gamma_v'(t)=0$\@:
\begin{align*}
\linder{\alpha_v}(t)=&\;\tau_{x,\gamma_v(t)}^{-1}(\nabla^0\xi(\gamma_v'(t)))\\
\linder[2]{\alpha_v}(t)=&\;\tau_{x,\gamma_v(t)}^{-1}
(\nabla^{(1)}\xi(\gamma_v'(t),\gamma_v'(t)))\\
\vdots\;&\\
\linder[m]{\alpha_v}(t)=&\;\tau_{x,\gamma_v(t)}^{-1}(\nabla^{(m-1)}\xi
(\underbrace{\gamma_v'(t),\dots,\gamma_v'(t)}_{m\ \textrm{times}})).
\end{align*}
By these computations, we have
\begin{equation*}
\frac{\d^m}{\d t^m}\Big|_{t=0}(\kappa_x^{-1}\scirc\xi(\exp_x(tv))=
\nabla^{(m-1)}\xi(\underbrace{v,\dots,v}_{m\ \textrm{times}}),
\end{equation*}
and so
\begin{equation*}
\kappa_x^{-1}\scirc\xi(\exp_x(tv))=\sum_{m=0}^\infty\frac{t^m}{m!}
\nabla^{(m-1)}\xi(\underbrace{v,\dots,v}_{m\ \textrm{times}}),
\end{equation*}
which is the desired result upon letting $t=1$ and supposing that $v$ is in a
sufficiently small neighbourhood of $0_x\in\tb[x]{\man{M}}$\@.
\end{subproof}
\end{proofsublemma}

Next we introduce some notation in the general setting of the preceding
sublemma that will be useful later.  We fix $x\in\man{M}$\@.  We let
$\nbhd{N}_x\subset\tb[x]{\man{M}}$ and $\nbhd{V}_x\subset\man{M}$ be
neighbourhoods of $0_x$ and $x$\@, respectively, such that
$\map{\exp_x}{\nbhd{N}_x}{\nbhd{V}_x}$ is a diffeomorphism.  For
$y\in\nbhd{V}_x$ we then define
\begin{equation*}
\mapdef{I'_{xy}}{\nbhd{N}'_{xy}\times\man{E}_x}{\man{E}|\nbhd{V}'_{xy}}
{(v,e_x)}{\tau_{x,\exp_x(v+\exp_x^{-1}(y))}(e_x)}
\end{equation*}
for neighbourhoods $\nbhd{N}'_{xy}\subset\tb[x]{\man{M}}$ of
$0_x\in\tb[x]{\man{M}}$ and $\nbhd{V}'_{xy}\subset\man{M}$ of $y$\@.  We note
that $I'_{xy}$ is a real analytic vector bundle isomorphism over the
diffeomorphism
\begin{equation*}
\mapdef{i'_{xy}}{\nbhd{N}'_{xy}}{\nbhd{V}'_{xy}}{v}
{\exp_x(v+\exp_x^{-1}(y)).}
\end{equation*}
Thus $I_{xy}\eqdef I'_{xy}\scirc\kappa_x^{-1}$ is a real analytic vector
bundle isomorphism from $\man{E}|\nbhd{U}'_{xy}$ to $\man{E}|\nbhd{V}'_{xy}$
for appropriate neighbourhoods $\nbhd{U}'_{xy}\subset\man{M}$ of $x$ and
$\nbhd{V}'_{xy}\subset\man{M}$ of $y$\@.  If we define
$\map{i_{xy}}{\nbhd{U}'_{xy}}{\nbhd{V}'_{xy}}$ by
$i_{xy}=i'_{xy}\scirc\exp_x^{-1}$\@, then $I_{xy}$ is a vector bundle mapping
over $i_{xy}$\@.  Along similar lines, $\hat{I}_{xy}\eqdef\kappa_y^{-1}\scirc
I'_{xy}$ is a vector bundle isomorphism between the trivial bundles
$\nbhd{O}'_{xy}\times\man{E}_x$ and $\nbhd{N}'_{xy}\times\man{E}_y$ for
appropriate neighbourhoods $\nbhd{O}'_{xy}\subset\tb[x]{\man{M}}$ and
$\nbhd{N}'_{xy}\subset\tb[y]{\man{M}}$ of the origin.  If we define
$\map{\hat{i}_{xy}}{\nbhd{O}'_{xy}}{\nbhd{N}'_{xy}}$ by
$\hat{i}_{xy}=\exp_y^{-1}\scirc i'_{xy}$\@, then $\hat{I}_{xy}$ is a vector
bundle map over $\hat{i}_{xy}$\@.

The next sublemma indicates that the neighbourhoods $\nbhd{U}'_{xy}$ of $x$
and $\nbhd{O}'_{xy}$ of $0_x$ can be uniformly bounded from below.
\begin{proofsublemma}\label{psublem:bddnbhds}
The neighbourhood\/ $\nbhd{V}_x$ and the neighbourhoods\/ $\nbhd{U}'_{xy}$
and\/ $\nbhd{O}'_{xy}$ above may be chosen so that
\begin{equation*}
\interior(\cap_{y\in\nbhd{V}_x}\nbhd{U}'_{xy})\not=\emptyset,\quad
\interior(\cap_{y\in\nbhd{V}_x}\nbhd{O}'_{xy})\not=\emptyset.
\end{equation*}
\begin{subproof}
By~\cite[Theorem~III.8.7]{SK/KN:63a} we can choose $\nbhd{V}_x$ so that, if
$y\in\nbhd{V}_x$\@, then there is a normal coordinate neighbourhood
$\nbhd{V}_y$ of $y$ containing $\nbhd{V}_x$\@.  Taking
$\nbhd{V}'_{xy}=\nbhd{V}_x\cap\nbhd{V}_y$ and $\nbhd{U}'_{xy}=\nbhd{V}_x$
gives the sublemma.
\end{subproof}
\end{proofsublemma}

We shall always assume $\nbhd{V}_x$ chosen as in the preceding sublemma, and
we let $\nbhd{U}'_x\subset\man{M}$ be a neighbourhood of $x$ and
$\nbhd{O}'_x\subset\tb[x]{\man{M}}$ be a neighbourhood of $0_x$ such that
\begin{equation*}
\nbhd{U}'_x\subset\interior(\cap_{y\in\nbhd{V}_x}\nbhd{U}'_{xy}),\quad
\nbhd{O}'_x\subset\interior(\cap_{y\in\nbhd{V}_x}\nbhd{O}'_{xy}).
\end{equation*}

These constructions can be ``bundled together'' as one to include the
dependence on $y\in\nbhd{V}_x$ in a clearer manner.  Since this will be
useful for us, we explain it here.  Let us denote
$\nbhd{D}_x=\nbhd{V}_x\times\nbhd{U}'_x$\@, let
$\map{\pr_2}{\nbhd{D}_x}{\nbhd{U}'_x}$ be the projection onto the second
factor, and denote
\begin{equation*}
\mapdef{i_x}{\nbhd{D}_x}{\man{M}}{(y,x')}{i_{xy}(x').}
\end{equation*}
Consider the pull-back bundle
$\map{\pr_2^*\pi}{\pr_2^*\man{E}|\nbhd{U}'_x}{\nbhd{D}_x}$\@.  Thus
\begin{equation*}
\pr_2^*\man{E}|\nbhd{U}'_x=
\setdef{((y,x'),e_{y'})\in\nbhd{D}_x\times\man{E}|\nbhd{U}'_x}{y'=x'}.
\end{equation*}
We then have a real analytic vector bundle mapping
\begin{equation*}
\mapdef{I_x}{\pr_2^*\man{E}|\nbhd{U}'_x}{\man{E}}{((y,x'),e_{x'})}
{I_{xy}(e_{x'})}
\end{equation*}
which is easily verified to be defined over $i_x$ and is isomorphic on
fibres.  Given $\xi\in\sections[\infty]{\man{E}}$\@, we define
$I_x^*\xi\in\sections[\infty]{\pr_2^*\man{E}|\nbhd{U}'_x}$ by
\begin{equation*}
I_x^*\xi(y,x')=(I_x)^{-1}_{(y,x')}\scirc\xi\scirc i_x(y,x')=
I_{xy}^{-1}\scirc\xi\scirc i_{xy}(x').
\end{equation*}
For $y\in\nbhd{V}_x$ fixed, we denote by
$I_{xy}^*\xi\in\sections[\infty]{\man{E}|\nbhd{U}'_x}$ the section given by
\begin{equation*}
I_{xy}^*\xi(x')=I_x^*\xi(y,x')=I_{xy}^{-1}\scirc\xi\scirc i_{xy}(x').
\end{equation*}

A similar construction can be made in the local trivialisations.  Here we
denote $\hat{\nbhd{D}}_x=\nbhd{V}_x\times\nbhd{O}_x$\@, let
$\map{\pr_2}{\hat{\nbhd{D}}_x}{\nbhd{O}_x}$ be the projection onto the second
factor, and consider the map
\begin{equation*}
\mapdef{\hat{i}_x}{\hat{\nbhd{D}}_x}{\tb{\man{M}}}
{(y,v_x)}{\hat{i}_{xy}(v_x).}
\end{equation*}
Denote by
$\map{\tbproj{\man{M}}^*\pi}{\tbproj{\man{M}}^*\man{E}}{\tb{\man{M}}}$ the
pull-back bundle and also define the pull-back bundle
\begin{equation*}
\map{\pr_2^*\tbproj{\man{M}}^*\pi}
{\pr_2^*\tbproj{\man{M}}^*\man{E}}{\hat{\nbhd{D}}_x}.
\end{equation*}
Note that
\begin{equation*}
\pr_2^*\tbproj{\man{M}}^*\man{E}=
\setdef{((y,v_x),(u_y,e_y))\in\hat{\nbhd{D}}_x\times
\tbproj{\man{M}}^*\man{E}}{x=y}.
\end{equation*}
We then define the real analytic vector bundle map
\begin{equation*}
\mapdef{\hat{I}_x}{\pr_2^*\tbproj{\man{M}}^*\man{E}}
{\tbproj{\man{M}}^*\man{E}}{((y,v_x),(u_x,e_x))}
{(v_x,\hat{I}_{xy}(v_x,e_x)).}
\end{equation*}
Given a local section $\eta\in\sections[\infty]{\tbproj{\man{M}}^*\man{E}}$
defined in a neighbourhood of the zero section, define a local section
$\hat{I}\null_x^*\eta\in
\sections[\infty]{\pr_2^*\tbproj{\man{M}}^*\man{E}}$ in a neighbourhood of
the zero section of $\hat{\nbhd{D}}_x$ by
\begin{equation*}
\hat{I}\null_x^*\eta(y,v_x)=(\hat{I}_x)^{-1}_{(y,v_x)}\scirc
\eta\scirc\hat{i}_x(y,v_x)=\hat{I}\null_{xy}^{-1}\scirc\eta\scirc i_x(y,v_x).
\end{equation*}
For $y\in\nbhd{V}_x$ fixed, we denote by $\eta_y$ the restriction of $\eta$
to a neighbourhood of $0_y\in\tb[y]{\man{M}}$\@.  We then denote by
\begin{equation*}
\hat{I}\null^*_{xy}\eta_y(v_x)=\hat{I}\null^*_x\eta(y,v_x)=
\hat{I}_{xy}^{-1}\scirc\eta_y\scirc\hat{i}_{xy}(v_x)
\end{equation*}
the element of $\sections[\infty]{\nbhd{O}'_x\times\man{E}_x}$\@.

The following simple lemma ties the preceding two constructions together.
\begin{proofsublemma}\label{psublem:IhatI}
Let\/ $\xi\in\sections[\infty]{\man{E}}$ and let\/
$\hat{\xi}\in\sections[\infty]{\tbproj{\man{M}}^*\man{E}}$ be defined in a
neighbourhood of the zero section by
\begin{equation*}
\hat{\xi}=\kappa_y^{-1}\scirc\xi\scirc\exp_y.
\end{equation*}
Then, for each\/ $y\in\nbhd{V}_x$\@,
\begin{equation*}
\hat{I}\null_{xy}^*\hat{\xi}_y=
\kappa_x^{-1}\scirc I_{xy}^*\xi\scirc\exp_x.
\end{equation*}
\begin{subproof}
We have
\begin{align*}
\hat{I}_{xy}^*\hat{\xi}(v_x)=&\;
\hat{I}\null^{-1}_{xy}\scirc\hat{\xi}\scirc\hat{i}_{xy}(v_x)\\
=&\;(I'_{xy})^{-1}\scirc\kappa_y\scirc\hat{\xi}\scirc
\exp_y^{-1}\scirc i'_{xy}(v_x)\\
=&\;\kappa_x^{-1}\scirc I_{xy}^{-1}\scirc\kappa_y\scirc
\hat{\xi}\scirc\exp_y^{-1}\scirc i_{xy}\scirc\exp_x(v_x)\\
=&\;\kappa_x^{-1}\scirc I_{xy}^{-1}\scirc\xi\scirc
i_{xy}\scirc\exp_x(v_x)\\
=&\;\kappa_x^{-1}I_{xy}^*\xi\scirc\exp_x(v_x).
\end{align*}
as claimed.
\end{subproof}
\end{proofsublemma}

Let us leave these general vector bundle considerations and proceed to local
estimates.  We shall consider estimates associated with local vector bundle
maps.  First we consider an estimate arising from multiplication.
\begin{proofsublemma}\label{psublem:multcont}
If\/ $\nbhd{U}\subset\real^n$ is open, if\/ $f\in\func[\omega]{\nbhd{U}}$\@,
and if\/ $K\subset\nbhd{U}$ is compact, then there exist\/
$C,\sigma\in\realp$ such that
\begin{equation*}
\sup\Bigsetdef{\frac{1}{I!}\linder[I]{(fg)}(\vect{x})}{\snorm{I}\le m}\le
C\sigma^{-m}\sup\Bigsetdef{\frac{1}{I!}\linder[I]{g}(\vect{x})}
{\snorm{I}\le m}
\end{equation*}
for every\/ $g\in\func[\infty]{\nbhd{U}}$\@,\/ $\vect{x}\in K$\@, and\/
$m\in\integernn$\@.
\begin{subproof}
For multi-indices $I,J\in\integernn^n$\@, let us write $J\le I$ if
$I-J\in\integernn^n$\@.  For $I\in\integernn^n$ we have
\begin{equation*}
\frac{1}{I!}\linder[I]{(fg)}(\vect{x})=\sum_{J\le I}
\frac{\linder[J]{g}(\vect{x})}{J!}\frac{\linder[I-J]{f}(\vect{x})}{(I-J)!},
\end{equation*}
by the Leibniz Rule.  By~\cite[Lemma~2.1.3]{SGK/HRP:02}\@, the number of
multi-indices in $n$ variables of order at most $\snorm{I}$ is
$\frac{(n+\snorm{I})!}{n!\snorm{I}!}$\@.  Note that, by the binomial theorem,
\begin{equation*}
(a_1+a_2)^{n+\snorm{I}}=\sum_{j=0}^{n+\snorm{I}}
\frac{(n+\snorm{I})!}{(n+\snorm{I}-j)!j!}a_1^ja_2^{n+\snorm{I}-j}.
\end{equation*}
Evaluating at $a_1=a_2=1$ and considering the summand corresponding to
$j=\snorm{I}$\@, this gives
\begin{equation*}
\frac{(n+\snorm{I})!}{n!\snorm{I}!}\le2^{n+\snorm{I}}.
\end{equation*}
Using this inequality we derive
\begin{align*}
\frac{1}{I!}\snorm{\linder[I]{(fg)}(\vect{x})}
\le&\;\sum_{\snorm{J}\le\snorm{I}}
\sup\Bigsetdef{\frac{\snorm{\linder[J]{f}(\vect{x})}}{J!}}
{\snorm{J}\le\snorm{I}}
\sup\Bigsetdef{\frac{\snorm{\linder[J]{g}(\vect{x})}}{J!}}
{\snorm{J}\le\snorm{I}}\\
\le&\;\frac{(n+\snorm{I})!}{n!\snorm{I}!}
\sup\Bigsetdef{\frac{\snorm{\linder[J]{f}(\vect{x})}}{J!}}
{\snorm{J}\le\snorm{I}}
\sup\Bigsetdef{\frac{\snorm{\linder[J]{g}(\vect{x})}}{J!}}
{\snorm{J}\le\snorm{I}}\\
\le&\;2^{n+\snorm{I}}\sup
\Bigsetdef{\frac{\snorm{\linder[J]{f}(\vect{x})}}{J!}}
{\snorm{J}\le\snorm{I}}
\sup\Bigsetdef{\frac{\snorm{\linder[J]{g}(\vect{x})}}{J!}}
{\snorm{J}\le\snorm{I}}.
\end{align*}
By~\cite[Proposition~2.2.10]{SGK/HRP:02}\@, there exist $B,r\in\realp$ such
that
\begin{equation*}
\frac{1}{J!}\snorm{\linder[J]{f}(\vect{x})}\le Br^{-\snorm{J}},\qquad
J\in\integernn^n,\ \vect{x}\in K.
\end{equation*}
We can suppose, without loss of generality, that $r<1$ so that we have
\begin{equation*}
\frac{1}{I!}\snorm{\linder[I]{(fg)}(\vect{x})}\le
2^nB\Bigl(\frac{2}{r}\Bigr)^{\snorm{I}}
\sup\Bigsetdef{\frac{\snorm{\linder[J]{g}(\vect{x})}}{J!}}
{\snorm{J}\le\snorm{I}},\qquad\vect{x}\in K.
\end{equation*}
We conclude, therefore, that if $\snorm{I}\le m$ we have
\begin{equation*}
\frac{1}{I!}\snorm{\linder[I]{(fg)}(\vect{x})}\le
2^nB\Bigl(\frac{2}{r}\Bigr)^{m}
\sup\Bigsetdef{\frac{\snorm{\linder[J]{g}(\vect{x})}}{J!}}
{\snorm{J}\le m},\qquad\vect{x}\in K,
\end{equation*}
which is the result upon taking $C=2^nB$ and $\sigma=\frac{2}{r}$\@.
\end{subproof}
\end{proofsublemma}

Next we give an estimate for derivatives of compositions of mappings, one of
which is real analytic.  Thus we have a real analytic mapping
$\map{\vect{\Phi}}{\nbhd{U}}{\nbhd{V}}$ between open sets
$\nbhd{U}\subset\real^n$ and $\nbhd{V}\subset\real^k$ and
$f\in\func[\infty]{\nbhd{V}}$\@.  By the higher-order Chain
Rule~\cite[\eg][]{GMC/THS:96}\@, we can write
\begin{equation*}
\linder[I]{(f\scirc\vect{\Phi})}(\vect{x})=
\sum_{\substack{H\in\integernn^m\\\snorm{H}\le\snorm{I}}}A_{I,H}(\vect{x})
\linder[H]{f}(\vect{\Phi}(\vect{x}))
\end{equation*}
for $\vect{x}\in\nbhd{U}$ and for some real analytic functions
$A_{I,H}\in\func[\omega]{\nbhd{U}}$\@.  The proof of the next sublemma gives
estimates for the $A_{I,H}$'s, and is based on computations of \citet{VT:97}
in the proof of his Proposition~2.5.
\begin{proofsublemma}\label{psublem:thilliez}
Let\/ $\nbhd{U}\subset\real^n$ and\/ $\nbhd{V}\subset\real^k$ be open, let\/
$\vect{\Phi}\in\mappings[\omega]{\nbhd{U}}{\nbhd{V}}$\@, and let\/
$K\subset\nbhd{U}$ be compact.  Then there exist\/ $C,\sigma\in\realp$ such
that
\begin{equation*}
\snorm{\linder[J]{A_{I,H}}(\vect{x})}\le
C\sigma^{-(\snorm{I}+\snorm{J})}(\snorm{I}+\snorm{J}-\snorm{H})!
\end{equation*}
for every\/ $\vect{x}\in K$\@,\/ $I,J\in\integernn^n$\@, and\/ $H\in\integernn^k$\@.
\begin{subproof}
First we claim that, for $j_1,\dots,j_r\in\{1,\dots,n\}$\@,
\begin{equation*}
\frac{\partial^r(f\scirc\vect{\Phi})}
{\partial x^{j_1}\cdots\partial x^{j_r}}(\vect{x})=
\sum_{s=1}^r\sum_{a_1,\dots,a_s=1}^kB^{a_1\cdots a_s}_{j_1\cdots j_r}(\vect{x})
\frac{\partial^sf}
{\partial y^{a_1}\cdots\partial y^{a_s}}(\vect{\Phi}(\vect{x})),
\end{equation*}
where the real analytic functions $B^{a_1\cdots a_s}_{j_1\cdots j_r}$\@,
$a_1,\dots,a_s\in\{1,\dots,k\}$\@, $j_1,\dots,j_r\in\{1,\dots,n\}$\@,
$r,s\in\integerp$\@, $s\le r$\@, are defined by the following recursion,
starting with $B^a_j=\pderiv{\Phi^a}{x^j}$\@:
\begin{compactenum}
\item $\displaystyle B^a_{j_1\cdots j_r}=
\frac{\partial B^a_{j_2\cdots j_r}}{\partial x^{j_1}}$\@;
\item $\displaystyle B^{a_1\cdots a_s}_{j_1\cdots j_r}=
\frac{\partial B^{a_1\cdots a_s}_{j_2\cdots j_r}}{\partial x^{j_1}}
+\frac{\partial\Phi^{a_1}}{\partial x^{j_1}}
B^{a_2\cdots a_s}_{j_2\cdots j_r}$\@, $r\ge2$\@, $s\in\{2,\dots,r-1\}$\@;
\item $\displaystyle B^{a_1\cdots a_r}_{j_1\cdots j_r}=
\frac{\partial\Phi^{a_1}}{\partial x^{j_1}}
B^{a_2\cdots a_r}_{j_2\cdots j_r}$\@.
\end{compactenum}
This claim we prove by induction on $r$\@.  It is clear for $r=1$\@, so
suppose the assertion true up to $r-1$\@.  By the induction hypothesis we
have
\begin{equation*}
\frac{\partial^{r-1}(f\scirc\vect{\Phi})}
{\partial x^{j_2}\cdots\partial x^{j_r}}(\vect{x})=
\sum_{s=1}^{r-1}\sum_{a_1,\dots,a_s=1}^k
B^{a_1\cdots a_s}_{j_2\cdots j_r}(\vect{x})
\frac{\partial^sf}{\partial y^{a_1}\cdots\partial y^{a_s}}
(\vect{\Phi}(\vect{x})).
\end{equation*}
We then compute
{\allowdisplaybreaks\begin{align*}
\frac{\partial}{\partial x^{j_1}}&
\frac{\partial^{r-1}(f\scirc\vect{\Phi})}
{\partial x^{j_2}\cdots\partial x^{j_r}}(\vect{x})\\
=&\;\sum_{s=1}^{r-1}\sum_{a_1,\dots,a_s=1}^k
\Bigl(\frac{\partial B^{a_1,\dots,a_s}_{j_2\cdots j_r}}
{\partial x^{j_1}}(\vect{x})
\frac{\partial^sf}{\partial y^{a_1}\cdots\partial y^{a_s}}
(\vect{\Phi}(\vect{x}))\\
&\;+\sum_{b=1}^kB^{a_1\cdots a_s}_{j_2\cdots j_r}(\vect{x})
\frac{\partial\Phi^b}{\partial x^{j_1}}(\vect{x})
\frac{\partial^{s+1}f}
{\partial y^b\partial y^{a_1}\cdots\partial y^{a_s}}
(\vect{\Phi}(\vect{x}))\Bigr)\\
=&\;\sum_{s=1}^{r-1}\sum_{a_1,\dots,a_s=1}^k
\frac{\partial B^{a_1\cdots a_s}_{j_2\cdots j_r}}
{\partial x^{j_1}}(\vect{x})
\frac{\partial^sf}{\partial y^{a_1}\cdots\partial y^{a_s}}
(\vect{\Phi}(\vect{x}))\\
&\;+\sum_{s=2}^r\sum_{a_1,\dots,a_s=1}^k
B^{a_2\cdots a_s}_{j_2\cdots j_r}(\vect{x})
\frac{\partial\Phi^{a_1}}{\partial x^{j_1}}(\vect{x})
\frac{\partial^sf}{\partial y^{a_1}\cdots\partial y^{a_s}}
(\vect{\Phi}(\vect{x}))\\
=&\;\sum_{a=1}^k
\frac{\partial B^a_{j_2\cdots j_r}}{\partial x^{j_1}}(\vect{x})
\frac{\partial f}{\partial y^a}(\vect{\Phi}(\vect{x}))\\
&\;+\sum_{s=2}^{r-1}\sum_{a_1,\dots,a_s=1}^k\Bigl(
\frac{\partial B^{a_1\cdots a_s}_{j_2\cdots j_r}}{\partial x^{j_1}}(\vect{x})
+\frac{\partial\Phi^{a_1}}{\partial x^{j_1}}(\vect{x})
B^{a_2\cdots a_s}_{j_2\cdots j_r}(\vect{x})\Bigl)
\frac{\partial^sf}{\partial y^{a_1}\cdots\partial y^{a_s}}
(\vect{\Phi}(\vect{x}))\\
&\;+\sum_{a_1,\dots,a_r=1}^k
\frac{\partial\Phi^{a_1}}{\partial x^{j_1}}(\vect{x})
B^{a_2\cdots a_r}_{j_2\cdots j_r}(\vect{x})
\frac{\partial^sf}{\partial y^{a_1}\cdots\partial y^{a_r}}
(\vect{\Phi}(\vect{x})),
\end{align*}}%
from which our claim follows.

Next we claim that there exist $A,\rho,\alpha,\beta\in\realp$ such that
\begin{equation*}
\snorm{\linder[J]{B^{a_1\cdots a_s}_{j_1\cdots j_r}}(\vect{x})}
\le(A\alpha)^r\Bigl(\frac{\beta}{\rho}\Bigr)^{r+\snorm{J}-s}(r+\snorm{J}-s)!
\end{equation*}
for every $\vect{x}\in K$\@, $J\in\integernn^n$\@,
$a_1,\dots,a_s\in\{1,\dots,k\}$\@, $j_1,\dots,j_k\in\{1,\dots,n\}$\@,
$r,s\in\integerp$\@, $s\le r$.  This we prove by induction on $r$ once again.
First let $\beta\in\realp$ be sufficiently large that
\begin{equation*}
\sum_{I\in\integernn^n}\beta^{-\snorm{I}}<\infty,
\end{equation*}
and denote this value of this sum by $S$\@.  Then let $\alpha=2S$\@.
By~\cite[Proposition~2.2.10]{SGK/HRP:02} there exist $A,\rho\in\realp$ such
that
\begin{equation*}
\snorm{\linder[J]{\linder[j]{\Phi^a}}(\vect{x})}\le AJ!\rho^{-\snorm{J}}
\end{equation*}
for every $\vect{x}\in K$\@, $J\in\integernn^n$\@, $j\in\{1,\dots,n\}$\@, and
$a\in\{1,\dots,k\}$\@.  This gives the claim for $r=1$\@.  So suppose the
claim true up to $r-1$\@.  Then, for any $a_1,\dots,a_s\in\{1,\dots,k\}$ and
$j_1,\dots,j_r\in\{1,\dots,n\}$\@, $s\le r$\@, $B^{a_1\cdots a_s}_{j_1\cdots
j_r}$ has one of the three forms listed above in the recurrent definition.
These three forms are themselves sums of terms of the form
\begin{equation*}
\underbrace{\pderiv{B^{a_1\cdots a_s}_{j_2\cdots j_r}}{{x^{j_1}}}}_P,\qquad
\underbrace{\pderiv{\Phi^{a_1}}{{x^{j_1}}}B^{a_2\cdots a_s}_{j_2\cdots j_r}}_Q.
\end{equation*}
Let us, therefore, estimate derivatives of these terms, abbreviated by $P$
and $Q$ as above.

We directly have, by the induction hypothesis,
\begin{align*}
\snorm{\linder[J]{P}(\vect{x})}\le&\;
(A\alpha)^r\Bigl(\frac{\beta}{\rho}\Bigr)^{r+\snorm{J}-s}(r+\snorm{J}-s)!\\
\le&A^r\alpha^{r-1}S\Bigl(\frac{\beta}{\rho}\Bigr)^{r+\snorm{J}-s}
(r+\snorm{J}-s)!,
\end{align*}
noting that $\alpha=2S$\@.  By the Leibniz Rule we have
\begin{equation*}
\linder[J]{Q}(\vect{x})=\sum_{J_1+J_2=J}\frac{J!}{J_1!J_2!}
\linder[J_1]{\linder[j_1]{\Phi^{a_1}}}(\vect{x})
\linder[J_2]{B^{a_2\cdots a_s}_{j_2\cdots j_r}}(\vect{x}).
\end{equation*}
By the induction hypothesis we have
\begin{equation*}
\snorm{\linder[J_2]{B^{a_2\cdots a_s}_{j_2\cdots j_r}}(\vect{x})}
\le(A\alpha)^{r-1}\Bigl(\frac{\beta}{\rho}\Bigr)^{r+\snorm{J_2}-s}
(r+\snorm{J_2}-s)!
\end{equation*}
for every $\vect{x}\in K$ and $J_2\in\integernn$\@.  Therefore,
\begin{equation*}
\snorm{\linder[J]{Q}(\vect{x})}\le\sum_{J_1+J_2=J}
\frac{J!}{J_2!}A(A\alpha)^{r-1}\Bigl(\frac{\beta}{\rho}\Bigr)^{r+\snorm{J}-s}
\beta^{-\snorm{J_1}}(r+\snorm{J_2}-s)!
\end{equation*}
for every $\vect{x}\in K$ and $J\in\integernn^n$\@.  Now note that, for any
$a,b,c\in\integerp$ with $b<c$\@, we have
\begin{equation*}
\frac{(a+b)!}{b!}=(1+b)\cdots(a+b)<(1+c)\cdots(a+c)=\frac{(a+c)!}{c!}.
\end{equation*}
Thus, if $L,J\in\integernn^n$ satisfy $L<J$ (meaning that
$J-L\in\integernn^n$), then we have
\begin{equation*}
l_k\le j_k\quad\implies\quad
\frac{(a+l_k)!}{l_k!}\le\frac{(a+j_k)!}{j_k!}\quad\implies\quad
\frac{j_k!}{l_k!}\le\frac{(a+j_k)!}{(a+l_k)!}
\end{equation*}
for every $a\in\integerp$ and $k\in\{1,\dots,n\}$\@.  Therefore,
\begin{equation*}
\frac{(j_1+\dots+j_{n-1}+j_n)!}{(j_1+\dots+j_{n-1}+l_n)!}
\ge\frac{j_n!}{l_n!}
\end{equation*}
and
\begin{multline*}
\frac{(j_1+\dots+j_{n-2}+j_{n-1}+j_n)!}
{(j_1+\dots+j_{n-2}+l_{n-1}+l_n)!}\\
=\frac{(j_1+\dots+j_{n-1}+j_n)!}{(j_1+\dots+j_{n-1}+l_n)!}
\frac{(j_1+\dots+j_{n-2}+j_{n-1}+l_n)!}{(j_1+\dots+j_{n-2}+l_{n-1}+l_n)!}
\ge\frac{j_{n-1}!}{l_{n-1}!}\frac{j_n!}{l_n!}.
\end{multline*}
Continuing in this way, we get
\begin{equation*}
\frac{J!}{L!}\le\frac{\snorm{J}!}{\snorm{L}!}.
\end{equation*}
We also have
\begin{equation*}
\frac{(r+\snorm{J_2}-s)!}{\snorm{J_2}!}\le\frac{(r+\snorm{J}-s)!}{\snorm{J}!}.
\end{equation*}
Thus we have
\begin{align*}
\snorm{\linder[J]{Q}(\vect{x})}\le&\;\sum_{J_1+J_2=J}
\frac{J!}{J_2!}A(A\alpha)^{r-1}\Bigl(\frac{\beta}{\rho}\Bigr)^{r+\snorm{J}-s}
\beta^{-\snorm{J_1}}(r+\snorm{J_2}-s)!\\
\le&\;A(A\alpha)^{r-1}\Bigl(\frac{\beta}{\rho}\Bigr)^{r+\snorm{J}-s}
(r+\snorm{J}-s)!\sum_{J_1+J_2=J}\beta^{-\snorm{J_1}}\\
\le&\;AS(A\alpha)^{r-1}\Bigl(\frac{\beta}{\rho}\Bigr)^{r+\snorm{J}-s}
(r+\snorm{J}-s)!
\end{align*}
Combining the estimates for $P$ and $Q$ to give an estimate for their sum,
and recalling that $\alpha=2S$\@, gives our claim that there exist $A,\rho,\alpha,\beta\in\realp$ such that
\begin{equation*}
\snorm{\linder[J]{B^{a_1\cdots a_s}_{j_1\cdots j_r}}(\vect{x})}
\le(A\alpha)^r\Bigl(\frac{\beta}{\rho}\Bigr)^{r+\snorm{J}-s}(r+\snorm{J}-s)!
\end{equation*}
for every $\vect{x}\in K$\@, $J\in\integernn^n$\@,
$a_1,\dots,a_s\in\{1,\dots,k\}$\@, and $j_1,\dots,j_r\in\{1,\dots,n\}$\@,
$r,s\in\integerp$\@, $s\le r$.

To conclude the proof of the lemma, note that given an index
$\vect{j}=(j_1,\dots,j_r)\in\{1,\dots,n\}^r$ we define a multi-index
$I(\vect{j})=(i_1,\dots,i_n)\in\integernn^n$ by asking that $i_l$ be the
number of times $l$ appears in the list $\vect{j}$\@.  Similarly an index
$\vect{a}=(a_1,\dots,a_s)\in\{1,\dots,k\}^s$ gives rise to a multi-index
$H(\vect{a})\in\integernn^k$\@.  Moreover, by construction we have
\begin{equation*}
B^{a_1\cdots a_s}_{j_1\cdots j_r}=A_{I(\vect{j}),H(\vect{a})}.
\end{equation*}
Let $C=1$ and $\sigma^{-1}=\max\{A\alpha,\frac{\beta}{\rho}\}$ and suppose,
without loss of generality, that $\sigma\le1$\@.  Then
\begin{equation*}
(A\alpha)^{\snorm{I}}\le\sigma^{-(\snorm{I}+\snorm{J})},\quad
\Bigl(\frac{\beta}{\rho}\Bigr)^{r+\snorm{J}-s}\le\sigma^{-(\snorm{I}+\snorm{J})}
\end{equation*}
for every $I,J\in\integernn^n$\@.  Thus we have
\begin{equation*}
\snorm{\linder[J]{A_{I,H}}(\vect{x})}\le
C\sigma^{-(\snorm{I}+\snorm{J})}(\snorm{I}+\snorm{J}-\snorm{H})!
\end{equation*}
as claimed.
\end{subproof}
\end{proofsublemma}

Next we consider estimates for derivatives arising from composition.
\begin{proofsublemma}\label{psublem:compcont}
Let\/ $\nbhd{U}\subset\real^n$ and\/ $\nbhd{V}\subset\real^k$ be open, let\/
$\vect{\Phi}\in\mappings[\omega]{\nbhd{U}}{\nbhd{V}}$\@, and let\/
$K\subset\nbhd{U}$ be compact.  Then there exist\/ $C,\sigma\in\realp$ such
that
\begin{equation*}
\sup\Bigsetdef{\frac{1}{I!}
\snorm{\linder[I]{(f\scirc\vect{\Phi})}(\vect{x})}}{\snorm{I}\le m}\\
\le C\sigma^{-m}\sup\Bigsetdef{\frac{1}{I!}
\snorm{\linder[H]{f}(\vect{\Phi}(\vect{x}))}}{\snorm{H}\le m}
\end{equation*}
for every\/ $f\in\func[\infty]{\nbhd{V}}$\@,\/ $\vect{x}\in K$\@, and\/
$m\in\integernn$\@.
\begin{subproof}
As we denoted preceding the statement of Sublemma~\ref{psublem:thilliez}
above, let us write
\begin{equation*}
\linder[I]{(f\scirc\vect{\Phi})}(\vect{x})=
\sum_{\substack{H\in\integernn^m\\\snorm{H}\le\snorm{I}}}A_{I,H}(\vect{x})
\linder[H]{f}(\vect{\Phi}(\vect{x}))
\end{equation*}
for $\vect{x}\in\nbhd{U}$ and for some real analytic functions
$A_{I,H}\in\func[\omega]{\nbhd{U}}$\@.  By Sublemma~\ref{psublem:thilliez}\@,
let $A,r\in\realp$ be such that
\begin{equation*}
\snorm{\linder[J]{A_{I,H}}(\vect{x})}\le Ar^{-(\snorm{I}+\snorm{J})}(\snorm{I}+\snorm{J}-\snorm{H})!
\end{equation*}
for $\vect{x}\in K$\@.  By the multinomial
theorem~\cite[Theorem~1.3.1]{SGK/HRP:02} we can write
\begin{equation*}
(a_1+\dots+a_n)^{\snorm{I}}=
\sum_{\snorm{J}=\snorm{I}}\frac{\snorm{J}!}{J!}\vect{a}^J
\end{equation*}
for every $I\in\integernn^n$\@.  Setting $a_1=\dots=a_n=1$ gives
$\frac{\snorm{I}!}{I!}\le n^{\snorm{I}}$ for every $I\in\integernn^n$\@.  As
in the proof of Sublemma~\ref{psublem:multcont} we have that the number of
multi-indices of length $k$ and degree at most $\snorm{I}$ is bounded above
by $2^{k+\snorm{I}}$\@.  Also, by a similar binomial theorem argument, if
$\snorm{H}\le\snorm{I}$\@, then we have
\begin{equation*}
\frac{(\snorm{I}-\snorm{H})!\snorm{H}!}{\snorm{I}!}\le2^{\snorm{I}}.
\end{equation*}
Putting this together yields
\begin{align*}
\frac{1}{I!}\snorm{\linder[I]{(f\scirc\vect{\Phi})}(\vect{x})}\le&\;
An^{\snorm{I}}r^{-\snorm{I}}
\sum_{\snorm{H}\le\snorm{I}}\frac{(\snorm{I}-\snorm{H})!\snorm{H}!}{\snorm{I}!}
\frac{1}{H!}\snorm{\linder[H]{f}(\vect{\Phi}(\vect{x}))}\\
\le&\;An^{\snorm{I}}2^{k+\snorm{I}}2^{\snorm{I}}r^{-\snorm{I}}
\sup\Bigsetdef{\frac{1}{H!}\snorm{\linder[H]{f}(\vect{\Phi}(\vect{x}))}}
{\snorm{H}\le\snorm{I}}\\
=&\;2^kA(4nr^{-1})^{\snorm{I}}
\sup\Bigsetdef{\frac{1}{H!}\snorm{\linder[H]{f}(\vect{\Phi}(\vect{x}))}}
{\snorm{H}\le\snorm{I}}
\end{align*}
whenever $\vect{x}\in K$\@.  Let us denote $C=2^kA$ and
$\sigma^{-1}=4nr^{-1}$ and take $r$ so that $4nr^{-1}\ge1$\@, without loss of
generality.  We then have
\begin{equation*}
\sup\Bigsetdef{\frac{1}{I!}
\snorm{\linder[I]{(f\scirc\vect{\Phi})}(\vect{x})}}{\snorm{I}\le m}
\le C\sigma^{-1}
\sup\Bigsetdef{\frac{1}{H!}\snorm{\linder[H]{f}(\vect{\Phi}(\vect{x}))}}
{\snorm{H}\le m}
\end{equation*}
for every $f\in\func[\infty]{\nbhd{U}_2}$\@,
$\vect{x}\in K$\@, and $m\in\integernn$\@, as claimed.
\end{subproof}
\end{proofsublemma}

Now we can state the following estimate for vector bundle mappings which is
essential for our proof.
\begin{proofsublemma}\label{psublem:localvbest}
Let\/ $\nbhd{U}\subset\real^n$ and\/ $\nbhd{V}\subset\real^k$ be open, let\/
$l\in\integerp$\@, and consider the trivial vector bundles\/
$\real^l_{\nbhd{U}}$ and\/ $\real^l_{\nbhd{V}}$\@.  Let\/
$\vect{\Phi}\in\mappings[\omega]{\nbhd{U}}{\nbhd{V}}$\@, let\/
$\mat{A}\in\mappings[\omega]{\nbhd{U}}{\GL[l]{\real}}$\@, and let\/
$K\subset\nbhd{U}$ be compact.  Then there exist\/ $C,\sigma\in\realp$ such
that
\begin{multline*}
\sup\Bigsetdef{\frac{1}{I!}
\snorm{\linder[I]{(\mat{A}^{-1}\cdot(\vect{\xi}\scirc
\vect{\Phi}))^b}(\vect{x})}}
{\snorm{I}\le m,\ b\in\{1,\dots,l\}}\\
\le C\sigma^{-m}\sup\Bigsetdef{\frac{1}{H!}
\snorm{\linder[H]{\xi^a}(\vect{\Phi}(\vect{x}))}}
{\snorm{H}\le m,\ a\in\{1,\dots,l\}},
\end{multline*}
for every\/ $\vect{\xi}\in\sections[\infty]{\real^l_{\nbhd{V}}}$\@,\/
$\vect{x}\in K$\@, and\/ $m\in\integernn$\@.
\begin{subproof}
By Sublemma~\ref{psublem:compcont} there exist $C_1,\sigma_1\in\realp$ such
that
\begin{multline*}
\sup\Bigsetdef{\frac{1}{I!}
\snorm{\linder[I]{(\vect{\xi}\scirc\vect{\Phi})^a}
(\vect{x})}}{\snorm{I}\le m,\ a\in\{1,\dots,l\}}\\
\le C_1\sigma_1^{-m}\sup\Bigsetdef{\frac{1}{H!}
\snorm{\linder[H]{\xi^a}(\vect{\Phi}(\vect{x}))}}
{\snorm{H}\le m,\ a\in\{1,\dots,l\}}
\end{multline*}
for every $\vect{\xi}\in\sections[\infty]{\real^l_{\nbhd{V}}}$\@,
$\vect{x}\in K$\@, and $m\in\integernn$\@.

Now let $\vect{\eta}\in\sections[\infty]{\real^l_{\nbhd{U}}}$\@.  Let
$B^b_a\in\func[\omega]{\nbhd{U}}$\@, $a\in\{1,\dots,l\}$\@,
$b\in\{1,\dots,l\}$\@, be the components of $\mat{A}^{-1}$\@.  By
Sublemma~\ref{psublem:multcont}\@, there exist $C_2,\sigma_2\in\realp$ such
that
\begin{multline*}
\sup\Bigsetdef{\frac{1}{I!}
\snorm{\linder[I]{(B^b_a(\vect{x})\eta^a(\vect{x}))}}}
{\snorm{I}\le m,\ a,b\in\{1,\dots,l\}}\\
\le C_2\sigma_2^{-m}\sup
\Bigsetdef{\frac{1}{I!}\snorm{\linder[I]{\eta^a}(\vect{x})}}
{\snorm{I}\le m,\ a\in\{1,\dots,l\}}
\end{multline*}
for every $\vect{x}\in K$ and $m\in\integernn$\@.  (There is no implied sum
over ``$a$'' in the preceding formula.)  Therefore, by the triangle
inequality,
\begin{multline*}
\sup\Bigsetdef{\frac{1}{I!}
\snorm{\linder[I]{(\mat{A}^{-1}\cdot\vect{\eta})^b}(\vect{x})}}
{\snorm{I}\le m,\ b\in\{1,\dots,l\}}\\
\le lC_2\sigma_2^{-m}
\sup\Bigsetdef{\frac{1}{I!}\snorm{\linder[I]{\eta^a}(\vect{x})}}
{\snorm{I}\le m,\ a\in\{1,\dots,l\}}
\end{multline*}
for every $\vect{x}\in K$ and $m\in\integernn$\@.

Combining the estimates from the preceding two paragraphs gives
\begin{multline*}
\sup\Bigsetdef{\frac{1}{I!}
\snorm{\linder[I]{(\vect{\xi}\scirc\vect{\Phi})^b}(\vect{x})}}
{\snorm{I}\le m,\ b\in\{1,\dots,l\}}\\
\le lC_1C_2(\sigma_1\sigma_2)^{-m}
\sup\Bigsetdef{\frac{1}{H!}\snorm{\linder[H]{\xi^a}(\vect{\Phi}(\vect{x}))}}
{\snorm{H}\le m,\ a\in\{1,\dots,l\}}
\end{multline*}
for every $\vect{\xi}\in\sections[\infty]{\real^l_{\nbhd{V}}}$\@,
$\vect{x}\in K$\@, and $m\in\integernn$\@, which is the desired result after
taking $C=lC_1C_2$ and $\sigma=\sigma_1\sigma_2$\@.
\end{subproof}
\end{proofsublemma}

Now we begin to provide some estimates that closely resemble those in the
statement of the lemma.  We begin by establishing an estimate resembling that
of the required form for a fixed $\vect{x}\in\nbhd{U}$\@.
\begin{proofsublemma}\label{psublem:fixedx}
Let\/ $\nbhd{U}\subset\real^n$ be open, denote\/
$\real^k_{\nbhd{U}}=\nbhd{U}\times\real^k$\@, and consider the trivial vector
bundle\/ $\map{\pr_1}{\real^k_{\nbhd{U}}}{\nbhd{U}}$\@.  Let\/ $\metric$ be a
Riemannian metric on\/ $\nbhd{U}$\@, let\/ $\metric_0$ be a vector bundle
metric on\/ $\real^k_{\nbhd{U}}$\@, let\/ $\nabla$ be an affine connection
on\/ $\nbhd{U}$\@, and let\/ $\nabla^0$ be a vector bundle connection on\/
$\real^k_{\nbhd{U}}$\@, with all of these being real analytic.  For\/
$\vect{\xi}\in\sections[\infty]{\real^k_{\nbhd{U}}}$ and\/
$\vect{x}\in\nbhd{U}$\@, denote by\/ $\hat{\vect{\xi}}_{\vect{x}}$ the
corresponding section of\/ $\nbhd{N}_{\vect{x}}\times\real^k$ defined by the
isomorphism\/ $\kappa_{\vect{x}}$ of Sublemma~\ref{psublem:local-nabla}\@.
For\/ $K\subset\nbhd{U}$ compact, there exist\/ $C,\sigma\in\realp$ such that
the following inequalities hold for each\/
$\vect{\xi}\in\sections[\infty]{\real^k_{\nbhd{U}}}$\@,\/ $\vect{x}\in K$\@,
and\/ $m\in\integernn$\@:
\begin{compactenum}[(i)]
\item $\dnorm{j_m\vect{\xi}(\vect{x})}_{\ol{\metric}_m}\le
C\sigma^{-m}\sup
\Bigsetdef{\frac{1}{I!}
\snorm{\linder[I]{\hat{\xi}\null_{\vect{x}}^a}(\vect{0})}}
{\snorm{I}\le m,\ a\in\{1,\dots,k\}}$\@;
\item $\Bigsetdef{\frac{1}{I!}
\snorm{\linder[I]{\hat{\xi}\null_{\vect{x}}^a}(\vect{0})}}
{\snorm{I}\le m,\ a\in\{1,\dots,k\}}\le C\sigma^{-m}
\dnorm{j_m\vect{\xi}(\vect{x})}_{\ol{\metric}_m}$\@.
\end{compactenum}
\begin{subproof}
By Sublemma~\ref{psublem:local-nabla} we have
\begin{equation*}
\hat{\vect{\xi}}_{\vect{x}}(\vect{v})=\sum_{m=0}^\infty\frac{1}{m!}
\nabla^{(m-1)}\vect{\xi}(\vect{x})
(\underbrace{\vect{v},\dots,\vect{v}}_{m\ \textrm{times}})
\end{equation*}
for $\vect{v}$ in some neighbourhood of $\vect{0}\in\real^n$\@.  We also have
\begin{equation*}
\hat{\vect{\xi}}_{\vect{x}}(\vect{v})=\sum_{m\in\integernn}\frac{1}{m!}
\linder[m]{\hat{\vect{\xi}}_{\vect{x}}}(\vect{0})
(\underbrace{\vect{v},\dots,\vect{v}}_{m\ \textrm{times}})
\end{equation*}
for every $\vect{v}$ in some neighbourhood of $\vect{0}\in\real^n$\@.  As
the relation
\begin{equation*}
\sum_{m=0}^\infty\frac{1}{m!}\nabla^{(m-1)}\vect{\xi}(\vect{x})
(\underbrace{\vect{v},\dots,\vect{v}}_{m\ \textrm{times}})=
\sum_{m\in\integernn}\frac{1}{m!}
\linder[m]{\hat{\vect{\xi}}_{\vect{x}}}(\vect{0})
(\underbrace{\vect{v},\dots,\vect{v}}_{m\ \textrm{times}})
\end{equation*}
holds for every $\vect{v}\in\real^n$\@, it follows that
\begin{equation*}
P^m_{\nabla,\nabla^0}(\vect{\xi})(\vect{x})=
\linder[m]{\hat{\vect{\xi}}_{\vect{x}}}(\vect{0})
\end{equation*}
for every $m\in\integernn$\@.  Take $m\in\integernn$\@.  We have
\begin{equation*}
\sum_{r=0}^m\frac{1}{(r!)^2}
\dnorm{P^r_{\nabla,\nabla^0}(\vect{\xi})(\vect{x})}^2_{\metric_r}
\le\sum_{r=0}^m\frac{A'A^r}{(r!)^2}
\dnorm{\linder[r]{\hat{\vect{\xi}}_{\vect{x}}}(\vect{0})}^2,
\end{equation*}
where $A'\in\realp$ depends on $\metric_0$\@, $A\in\realp$ depends on
$\metric$\@, and where $\dnorm{\cdot}$ denotes the $2$-norm,~\ie~the square
root of the sum of squares of components.  We can, moreover, assume without
loss of generality that $A\ge1$ so that we have
\begin{equation*}
\sum_{r=0}^m\frac{1}{(r!)^2}
\dnorm{P^r_{\nabla,\nabla^0}(\vect{\xi})(\vect{x})}^2_{\metric_r}
\le A'A^m\sum_{r=0}^m\frac{1}{(r!)^2}
\dnorm{\linder[r]{\hat{\vect{\xi}}_{\vect{x}}}(\vect{0})}^2.
\end{equation*}
By~\cite[Lemma~2.1.3]{SGK/HRP:02}\@,
\begin{equation*}
\card(\setdef{I\in\integernn^n}{\snorm{I}\le m})=\frac{(n+m)!}{n!m!}.
\end{equation*}
Note that the $2$-norm for $\real^N$ is related to the $\infty$-norm for
$\real^N$ by $\dnorm{\vect{a}}_2\le\sqrt{N}\dnorm{\vect{a}}_\infty$ so that
\begin{equation*}
\sum_{r=0}^m\frac{1}{(r!)^2}
\dnorm{\linder[r]{\hat{\vect{\xi}}_{\vect{x}}}(\vect{0})}^2\le 
k\frac{(n+m)!}{n!m!}\Bigl(\sup\Bigsetdef{\frac{1}{r!}
\snorm{\linder[I]{\hat{\xi}\null_{\vect{x}}^a}(\vect{0})}}
{\snorm{I}\le m,\ a\in\{1,\dots,k\}}\Bigr)^2.
\end{equation*}
By the binomial theorem, as in the proof of
Sublemma~\ref{psublem:multcont}\@,
\begin{equation*}
\frac{(n+m)!}{n!m!}\le 2^{n+m}.
\end{equation*}
Thus
\begin{equation}\label{eq:pointest1}
\dnorm{j_m\vect{\xi}(\vect{x})}_{\ol{\metric}_m}\le
\sqrt{kA'2^n}(\sqrt{2A})^m\sup
\Bigsetdef{\frac{1}{I!}
\snorm{\linder[I]{\hat{\xi}\null_{\vect{x}}^a}(\vect{0})}}
{\snorm{I}\le m,\ a\in\{1,\dots,k\}}
\end{equation}
for every $m\in\integernn$\@.  The above computations show that this
inequality is satisfied for a real analytic section $\vect{\xi}$\@.  However,
it also is satisfied if $\vect{\xi}$ is a smooth section.  This we argue as
follows.  Let $\vect{\xi}\in\sections[\infty]{\real^k_{\nbhd{U}}}$ and, for
$m\in\integerp$\@, let $\vect{\xi}_m\in\sections[\omega]{\real^k_{\nbhd{U}}}$
be the section whose coefficients are polynomial functions of degree at most
$m$ and such that $j_m\vect{\xi}_m(\vect{x})=j_m\vect{\xi}(\vect{x})$\@.
Also let $\hat{\vect{\xi}}_{\vect{x},m}$ be the corresponding section of
$\nbhd{N}_{\vect{x}}\times\real^k$\@.  We then have
\begin{equation*}
j_m\vect{\xi}_m(\vect{x})=j_m\vect{\xi}(\vect{x}),\quad
\linder[I]{\hat{\vect{\xi}}_{\vect{x},m}}(\vect{0})=
\linder[I]{\hat{\vect{\xi}}_{\vect{x}}}(\vect{0}),
\end{equation*}
for every $I\in\integernn^n$ satisfying $\snorm{I}\le m$\@, the latter by the
formula for the higher-order Chain
Rule~\cite[Supplement~2.4A]{RA/JEM/TSR:88}\@.  Since $\vect{\xi}_m$ is real
analytic, this shows that~\eqref{eq:pointest1} is also satisfied for every
$m\in\integernn$ if $\vect{\xi}$ is smooth.

To establish the other estimate asserted in the sublemma, let $\vect{x}\in K$
and, using the notation of Sublemma~\ref{psublem:local-nabla}\@, let
$\nbhd{N}_{\vect{x}}$ be a relatively compact neighbourhood of
$\vect{0}\in\real^n\simeq\tb[\vect{x}]{\real^n}$ and
$\nbhd{V}_{\vect{x}}\subset\nbhd{U}$ be a relatively compact neighbourhood of
$\vect{x}$ such that
$\map{\kappa_{\vect{x}}}{\nbhd{N}_{\vect{x}}\times\real^k}
{\nbhd{V}_{\vect{x}}\times\real^k}$ is a real analytic vector bundle
isomorphism.  Let
$\vect{\xi}\in\sections[\omega]{\real^k_{\nbhd{V}_{\vect{x}}}}$ and let
$\hat{\vect{\xi}}_{\vect{x}}\in\sections[\omega]{\real^k_{\nbhd{N}_{\vect{x}}}}$
be defined by
$\hat{\vect{\xi}}_{\vect{x}}(\vect{v})=\kappa_{\vect{x}}^{-1}\scirc\vect{\xi}
(\exp_{\vect{x}}(\vect{v}))$\@.  As in the first part of the estimate, we
have
\begin{equation*}
\linder[m]{\hat{\vect{\xi}}_{\vect{x}}}(\vect{0})=
P^m_{\nabla,\nabla^0}(\vect{\xi})(\vect{x})
\end{equation*}
for every $m\in\integernn$\@.  For indices
$\vect{j}=(j_1,\dots,j_m)\in\{1,\dots,n\}^m$ we define
$I(\vect{j})=(i_1,\dots,i_n)\in\integernn^n$ by asking that $i_j$ be the
number of times ``$j$'' appears in the list $\vect{j}$\@.  We then have
\begin{multline*}
\sup\Bigsetdef{\frac{1}{I!}
\snorm{\linder[I]{\hat{\xi}_{\vect{x}}^a}(\vect{0})}}
{\snorm{I}\le m,\ a\in\{1,\dots,k\}}=\\
\sup\Bigl\{\frac{1}{I(\vect{j})!}
\snorm{(P^r_{\nabla,\nabla^0}(\vect{\xi})(\vect{x}))^a_{j_1\cdots j_r}}
\Big|\\j_1,\dots,j_r\in\{1,\dots,n\},\ r\in\{0,1,\dots,m\},\
a\in\{1,\dots,m\}\Bigr\}.
\end{multline*}
By an application of the multinomial theorem as in the proof of
Sublemma~\ref{psublem:compcont}\@, we have $\frac{\snorm{I}!}{I!}\le
n^{\snorm{I}}$ for every $I\in\integernn^n$\@.  We then have
\begin{equation*}
\frac{1}{I(\vect{j})!}
\snorm{(P^r_{\nabla,\nabla^0}(\vect{\xi})(\vect{x}))^a_{j_1\cdots j_r}}\le
\frac{n^r}{r!}\snorm{(P^r_{\nabla,\nabla^0}(\vect{\xi})
(\vect{x}))^a_{j_1\cdots j_r}}
\end{equation*}
for every $j_1,\dots,j_r\in\{1,\dots,n\}$ and $a\in\{1,\dots,k\}$\@.  Using
the fact that the $\infty$-norm for $\real^N$ is related to the $2$-norm for
$\real^N$ by $\dnorm{\vect{a}}_\infty\le\dnorm{\vect{a}}_2$\@, we have
\begin{equation*}
\sup\Bigsetdef{\frac{1}{I!}
\snorm{\linder[I]{\hat{\xi}_{\vect{x}}^a}(\vect{0})}}
{\snorm{I}\le m,\ a\in\{1,\dots,k\}}\le\left(\sum_{r=0}^m
\Bigl(\frac{n^r}{r!}\Bigr)^2B'B^r\dnorm{P^r_{\nabla,\nabla^0}
(\vect{\xi})(\vect{x})}_{\metric_r}^2\right)^{1/2},
\end{equation*}
where $B'\in\realp$ depends on $\metric_0$ and $B\in\realp$ depends on
$\metric$\@.  We may, without loss of generality, suppose that $B\ge1$ so
that we have
\begin{equation*}
\sup\Bigsetdef{\frac{1}{I!}
\snorm{\linder[I]{\hat{\xi}\null_{\vect{x}}^a}(\vect{0})}}
{\snorm{I}\le m,\ a\in\{1,\dots,k\}}\le
\sqrt{B'}(n\sqrt{B})^m\dnorm{j_m\vect{\xi}(\vect{x})}_{\ol{\metric}_m}
\end{equation*}
for every $m\in\integernn$\@.  As in the first part of the proof, while we
have demonstrated the preceding inequality for $\vect{\xi}$ real analytic, it
can also be demonstrated to hold for $\vect{\xi}$ smooth.

The sublemma follows by taking
\begin{equation*}\eqsubqed
C=\max\{\sqrt{kA'2^n},\sqrt{B'}\},\quad
\sigma^{-1}=\max\{\sqrt{2}A,n\sqrt{B}\}.
\end{equation*}
\end{subproof}
\end{proofsublemma}

The next estimates we consider will allow us to expand the pointwise estimate
from the preceding sublemma to a local estimate of the same form.  The
construction makes use of the vector bundle isomorphisms $I_{xy}$ and
$\hat{I}_{xy}$ defined after Sublemma~\ref{psublem:local-nabla}\@.  In the
statement and proof of the following sublemma, we make free use of the
notation we introduced where these mappings were defined.
\begin{proofsublemma}\label{psublem:Ixyest}
Let\/ $\nbhd{U}\subset\real^n$ be open, denote\/
$\real^k_{\nbhd{U}}=\nbhd{U}\times\real^k$\@, and consider the trivial vector
bundle\/ $\map{\pr_1}{\real^k_{\nbhd{U}}}{\nbhd{U}}$\@.  Let\/ $\metric$ be a
Riemannian metric on\/ $\nbhd{U}$\@, let\/ $\metric_0$ be a vector bundle
metric on\/ $\real^k_{\nbhd{U}}$\@, let\/ $\nabla$ be an affine connection
on\/ $\nbhd{U}$\@, and let\/ $\nabla^0$ be a vector bundle connection on\/
$\real^k_{\nbhd{U}}$\@, with all of these being real analytic.  For each\/
$\vect{x}\in\nbhd{U}$ there exist a neighbourhood\/ $\nbhd{V}_{\vect{x}}$
and\/ $C_{\vect{x}},\sigma_{\vect{x}}\in\realp$ such that we have the
following inequalities for each\/
$\vect{\xi}\in\sections[\infty]{\real^k_{\nbhd{U}}}$\@, $m\in\integernn$\@,
and\/ $\vect{y}\in\nbhd{V}_{\vect{x}}$\@:
\begin{compactenum}[(i)]
\item \label{pl:Ixyest1} $\sup\bigsetdef{\frac{1}{I!}
\snorm{\linder[I]{\hat{\xi}\null_{\vect{y}}^a}(\vect{0})}}
{\snorm{I}\le m,\ a\in\{1,\dots,k\}}$\\
$\le C_{\vect{x}}\sigma_{\vect{x}}^{-1}
\sup\bigsetdef{\frac{1}{I!}
\snorm{\linder[I]{((\hat{I}\null^*_{\vect{x}\vect{y}})^{-1}
\hat{\vect{\xi}}_{\vect{y}})^a}(\vect{0})}}
{\snorm{I}\le m,\ a\in\{1,\dots,k\}}$\@;
\item \label{pl:Ixyest2} $\sup\bigsetdef{\frac{1}{I!}
\snorm{\linder[I]{(\hat{I}\null^*_{\vect{x}\vect{y}}}
\hat{\vect{\xi}}_{\vect{y}})^a}(\vect{0})}
{\snorm{I}\le m,\ a\in\{1,\dots,k\}}$\\
$\le C_{\vect{x}}\sigma_{\vect{x}}^{-1}
\sup\bigsetdef{\frac{1}{I!}
\snorm{\linder[I]{\hat{\xi}\null_{\vect{y}}^a}(\vect{0})}}
{\snorm{I}\le m,\ a\in\{1,\dots,k\}}$\@;
\item \label{pl:Ixyest3} $\dnorm{j_m\vect{\xi}(\vect{y})}_{\ol{\metric}_m}\le
C_{\vect{x}}\sigma_{\vect{x}}^{-1}
\dnorm{j_m((I_{\vect{x}\vect{y}}^*)^{-1}
\vect{\xi})(\vect{x})}_{\ol{\metric}_m}$\@;
\item \label{pl:Ixyest4} $\dnorm{j_m(I_{\vect{x}\vect{y}}^*\vect{\xi})(\vect{x})}_{\ol{\metric}_m}\le
C_{\vect{x}}\sigma_{\vect{x}}^{-1}
\dnorm{j_m\vect{\xi}(\vect{y})}_{\ol{\metric}_m}$\@.
\end{compactenum}
\begin{subproof}
We begin the proof with an observation.  Suppose that we have an open subset
$\nbhd{U}\subset\real^n\times\real^k$ and $f\in\func[\omega]{\nbhd{U}}$\@.
We wish to think of $f$ as a function of $\vect{x}\in\real^n$ depending on a
parameter $\vect{p}\in\real^k$ in a jointly real analytic manner.  We note
that, for $K\subset\nbhd{U}$ compact, we have $C,\sigma\in\realp$ such that
the partial derivatives satisfy a bound
\begin{equation*}
\snorm{\plinder[I]{1}{f}(\vect{x},\vect{p})}\le CI!\sigma^{-\snorm{I}}
\end{equation*}
for every $(\vect{x},\vect{p})\in K$ and $I\in\integernn^n$\@.  This is a
mere specialisation of~\cite[Proposition~2.2.10]{SGK/HRP:02} to partial
derivatives.  The point is that the bound for the partial derivatives is
uniform in the parameter $\vect{p}$\@.  With this in mind, we note that the
following are easily checked:
\begin{compactenum}
\item the estimate of Sublemma~\ref{psublem:multcont} can be extended to the
case where $f$ depends in a jointly real analytic manner on a parameter, and
the estimate is uniform in the parameter over compact sets;
\item the estimate of Sublemma~\ref{psublem:thilliez} can be extended to the
case where $\vect{\Phi}$ depends in a jointly real analytic manner on a
parameter, and the estimate is uniform in the parameter over compact sets;
\item as a consequence of the preceding fact, the estimate of
Sublemma~\ref{psublem:compcont} can be extended to the case where
$\vect{\Phi}$ depends in a jointly real analytic manner on a parameter, and
the estimate is uniform in the parameter over compact sets;
\item as a consequence of the preceding three facts, the estimate of
Sublemma~\ref{psublem:localvbest} can be extended to the case where
$\vect{\Phi}$ and $\mat{A}$ depend in a jointly real analytic manner on a
parameter, and the estimate is uniform in the parameter over compact sets.
\end{compactenum}
Now let us proceed with the proof.

We take $\nbhd{V}_{\vect{x}}$ as in the discussion following
Sublemma~\ref{psublem:local-nabla}\@.  Let us introduce coordinate notation
for all maps needed.  We have
\begin{align*}
\hat{\vect{\xi}}_{\vect{y}}(\vect{u})=&\;\hat{\vect{\xi}}(\vect{y},\vect{u})=
\vect{\xi}\scirc\exp_{\vect{y}}(\vect{u}),\\
I_{\vect{x}\vect{y}}^*\vect{\xi}(\vect{x}')=&\;
\mat{A}(\vect{y},\vect{x}')\cdot
(\vect{\xi}\scirc i_{\vect{x}\vect{y}}(\vect{x}')),\\
\hat{I}\null_{\vect{x}\vect{y}}^*\hat{\vect{\xi}}_{\vect{y}}(\vect{v})=&\;
\hat{\mat{A}}(\vect{y},\vect{v})\cdot
(\hat{\vect{\xi}}_{\vect{y}}\scirc\hat{i}_{\vect{x}\vect{y}}(\vect{v})),\\
(I_{\vect{x}\vect{y}}^*)^{-1}\vect{\xi}(\vect{y}')=&\;
\mat{A}^{-1}(\vect{y},i_{\vect{x}\vect{y}}^{-1}(\vect{y}')\cdot
(\vect{\xi}\scirc i^{-1}_{\vect{x}\vect{y}}(\vect{y}')),\\
(\hat{I}\null_{\vect{x}\vect{y}}^*)^{-1}\hat{\vect{\xi}}_{\vect{y}}(\vect{v})=&\;
\hat{\mat{A}}\null^{-1}(\vect{y},\hat{i}_{\vect{x}\vect{y}}^{-1}(\vect{u}))
\cdot(\hat{\vect{\xi}}_{\vect{y}}\scirc
\hat{i}\null^{-1}_{\vect{x}\vect{y}}(\vect{v})),
\end{align*}
for appropriate real analytic mappings $\mat{A}$ and $\hat{\mat{A}}$ taking
values in $\GL[k]{\real}$\@.  Note that, for every $I\in\integernn^n$\@,
\begin{equation*}
\linder[I]{(\hat{I}\null^*_{\vect{x}\vect{y}}\hat{\vect{\xi}}_{\vect{y}})}
(\vect{0})=\plinder[I]{2}{(\hat{I}\null^*_{\vect{x}}
\hat{\vect{\xi}})}(\vect{y},\vect{0}),
\end{equation*}
and similarly for $\linder[I]{((\hat{I}\null^*_{\vect{x}\vect{y}})^{-1}
\hat{\vect{\xi}}_{\vect{y}})}(\vect{0})$\@.  The observation made at the
beginning of the proof shows that parts~\eqref{pl:Ixyest1}
and~\eqref{pl:Ixyest2} follow immediately from
Sublemma~\ref{psublem:localvbest}\@.  Parts~\eqref{pl:Ixyest3}
and~\eqref{pl:Ixyest4} follow from the first two parts after an application
of Sublemma~\ref{psublem:fixedx}\@.
\end{subproof}
\end{proofsublemma}

By applications
of~(a)~Sublemma~\ref{psublem:Ixyest}\@,~(b)~Sublemmata~\ref{psublem:IhatI}
and~\ref{psublem:fixedx}\@,~(c)~Sublemma~\ref{psublem:Ixyest} again,
and~(d)~Sublemma~\ref{psublem:localvbest}\@, there exist
\begin{equation*}
A_{1,\vect{x}},A_{2,\vect{x}},A_{3,\vect{x}},A_{4,\vect{x}},
r_{1,\vect{x}},r_{2,\vect{x}},r_{3,\vect{x}},r_{4,\vect{x}}\in\realp
\end{equation*}
and a relatively compact neighbourhood $\nbhd{V}_{\vect{x}}\subset\nbhd{U}$
of $\vect{x}$ such that
\begin{align*}
\dnorm{j_m\vect{\xi}(\vect{y})}_{\ol{\metric}_m}\le&\;
A_{1,\vect{x}}r_{1,\vect{x}}^{-m}
\dnorm{j_m((I_{\vect{x}\vect{y}}^*)^{-1}\vect{\xi})(\vect{x})}_{\ol{\metric}_m}\\
\le&\;A_{2,\vect{x}}r_{2,\vect{x}}^{-m}
\sup\Bigsetdef{\frac{1}{I!}
\snorm{\linder[I]{((\hat{I}\null_{\vect{x}\vect{y}}^*)^{-1}
\hat{\vect{\xi}}_{\vect{y}})^a}(\vect{0})}}
{\snorm{I}\le m,\ a\in\{1,\dots,k\}}\\
\le&\;A_{3,\vect{x}}r_{3,\vect{x}}^{-m}
\sup\Bigsetdef{\frac{1}{I!}
\snorm{\linder[I]{\hat{\xi}\null_{\vect{y}}^a}(\vect{0})}}
{\snorm{I}\le m,\ a\in\{1,\dots,k\}}\\
\le&\;A_{4,\vect{x}}r_{4,\vect{x}}^{-m}
\sup\Bigsetdef{\frac{1}{I!}\snorm{\linder[I]{\xi^a}(\vect{y})}}
{\snorm{I}\le m,\ a\in\{1,\dots,k\}}
\end{align*}
for every $\vect{\xi}\in\sections[\infty]{\real^k_{\nbhd{U}}}$\@,
$m\in\integernn$\@, and $\vect{y}\in\nbhd{V}_{\vect{x}}$\@.  Take
$\vect{x}_1,\dots,\vect{x}_k\in K$ such that
$K\subset\cup_{j=1}^k\nbhd{V}_{\vect{x}_j}$ and define
\begin{equation*}
C_1=\max\{A_{4,\vect{x}_1},\dots,A_{4,\vect{x}_k}\},\quad
\sigma_1=\min\{r_{4,\vect{x}_1},\dots,r_{4,\vect{x}_k}\},
\end{equation*}
so that
\begin{equation*}
\dnorm{j_m\vect{\xi}(\vect{x})}_{\ol{\metric}_m}\le
C_1\sigma_1^{-m}\sup\Bigsetdef{\frac{1}{I!}
\snorm{\linder[I]{\xi^a}(\vect{x})}}{\snorm{I}\le m,\ a\in\{1,\dots,k\}}
\end{equation*}
for every $\vect{\xi}\in\sections[\infty]{\real^k_{\nbhd{U}}}$\@,
$m\in\integernn$\@, and $\vect{x}\in K$\@.  This gives one half of the
estimate in the lemma.

For the other half of the estimate in the lemma, we
apply~(a)~Sublemma~\ref{psublem:localvbest}\@,~(b)~%
Sublemma~\ref{psublem:Ixyest}\@,~(c)~Sublemmata~\ref{psublem:IhatI}
and~\ref{psublem:fixedx}\@, and~(d)~Sublemma~\ref{psublem:Ixyest} again to
assert the existence of
\begin{equation*}
A_{1,\vect{x}},A_{2,\vect{x}},A_{3,\vect{x}},A_{4,\vect{x}},
r_{1,\vect{x}},r_{2,\vect{x}},r_{3,\vect{x}},r_{4,\vect{x}}\in\realp
\end{equation*}
and a relatively compact neighbourhood $\nbhd{V}_{\vect{x}}\subset\nbhd{U}$
of $\vect{x}$ such that
\begin{align*}
\sup\Bigl\{\frac{1}{I!}\snorm{\linder[I]{\xi^a}(\vect{y})}&\Big|\enspace
\snorm{I}\le m,\ a\in\{1,\dots,m\}\Bigr\}\\
\le&\;A_{1,\vect{x}}r_{1,\vect{x}}^{-m}\sup\Bigsetdef{\frac{1}{I!}
\snorm{\linder[I]{\hat{\xi}\null^a_{\vect{y}}}(\vect{0})}}
{\snorm{I}\le m,\ a\in\{1,\dots,k\}}\\
\le&\;A_{2,\vect{x}}r_{2,\vect{x}}^{-m}\sup\Bigsetdef{\frac{1}{I!}
\snorm{\linder[I]{((\hat{I}\null^*_{\vect{x}\vect{y}})^{-1}
\hat{\vect{\xi}}_{\vect{y}})^a}(\vect{0})}}
{\snorm{I}\le m,\ a\in\{1,\dots,k\}}\\
\le&\;A_{3,\vect{x}}r_{3,\vect{x}}^{-m}
\dnorm{j_m((I_{\vect{x}\vect{y}}^*)^{-1}\vect{\xi})(\vect{x})}_{\ol{\metric}_m}
\le A_{4,\vect{x}}r_{4,\vect{x}}^{-m}
\dnorm{j_m\vect{\xi}(\vect{y})}_{\ol{\metric}_m}
\end{align*}
for every $\vect{\xi}\in\sections[\infty]{\real^k_{\nbhd{U}}}$\@,
$m\in\integernn$\@, and $\vect{y}\in\nbhd{V}_{\vect{x}}$\@.  As we argued
above using a standard compactness argument, there exist
$C_2,\sigma_2\in\realp$ such that
\begin{equation*}
\sup\Bigsetdef{\frac{1}{I!}
\snorm{\linder[I]{\xi^a}(\vect{x})}}
{\snorm{I}\le m,\ a\in\{1,\dots,k\}}\le C_2\sigma_2^{-m}
\dnorm{j_m\vect{\xi}(\vect{x})}_{\ol{\metric}_m}
\end{equation*}
for every $\vect{\xi}\in\sections[\infty]{\real^k_{\nbhd{U}}}$\@,
$m\in\integernn$\@, and $\vect{x}\in K$\@.  Taking $C=\max\{C_1,C_2\}$ and
$\sigma=\min\{\sigma_1,\sigma_2\}$ gives the lemma.
\end{proof}
\end{lemma}

The preceding lemma will come in handy on a few crucial occasions.  To
illustrate how it can be used, we give the following characterisation of real
analytic sections, referring to Section~\ref{sec:smooth-topology} below for
the definition of the seminorm $p_{K,m}^\infty$ used in the statement.
\begin{lemma}[Characterisation of real analytic sections]\label{lem:rabound}
Let\/ $\map{\pi}{\man{E}}{\man{M}}$ be a real analytic vector bundle and
let\/ $\xi\in\sections[\infty]{\man{E}}$\@.  Then the following statements
hold:
\begin{compactenum}[(i)]
\item \label{pl:rabound1} $\xi\in\sections[\omega]{\man{E}}$\@;
\item \label{pl:rabound2} for every compact set\/ $K\subset\man{M}$\@,
there exist\/ $C,r\in\realp$ such that\/ $p^\infty_{K,m}(\xi)\le Cr^{-m}$
for every\/ $m\in\integernn$\@.
\end{compactenum}
\begin{proof}
\eqref{pl:rabound1}$\implies$\eqref{pl:rabound2} Let $K\subset\man{M}$ be
compact, let $x\in K$\@, and let $(\nbhd{V}_x,\psi_x)$ be a vector bundle
chart for $\man{E}$ with $(\nbhd{U}_x,\phi_x)$ the corresponding chart for
$\man{M}$\@.  Let $\map{\vect{\xi}}{\phi(\nbhd{U}_x)}{\real^k}$ be the local
representative of $\xi$\@.  By~\cite[Proposition~2.2.10]{SGK/HRP:02}\@, there
exist a neighbourhood $\nbhd{U}'_x\subset\nbhd{U}_x$ of $x$ and
$B_x,\sigma_x\in\realp$ such that
\begin{equation*}
\snorm{\linder[I]{\xi^a}(\vect{x}')}\le B_xI!\sigma_x^{-\snorm{I}}
\end{equation*}
for every $a\in\{1,\dots,k\}$\@, $\vect{x}'\in\closure(\nbhd{U}'_x)$\@, and
$I\in\integernn^n$\@.  We can suppose, without loss of generality, that
$\sigma_x\in\interval(0,1)$\@.  In this case, if $\snorm{I}\le m$\@,
\begin{equation*}
\frac{1}{I!}\snorm{\linder[I]{\xi^a}(\vect{x}')}\le B_x\sigma_x^{-m}
\end{equation*}
for every $a\in\{1,\dots,k\}$ and $\vect{x}'\in\closure(\nbhd{U}'_x)$\@.  By
Lemma~\ref{lem:pissy-estimate}\@, there exist $C_x,r_x\in\realp$ such that
\begin{equation*}
\dnorm{j_m\xi(x')}_{\ol{\metric}_m}\le C_xr_x^{-m},\qquad
x'\in\closure(\nbhd{U}'_x),\ m\in\integernn.
\end{equation*}
Let $x_1,\dots,x_k\in K$ be such that $K\subset\cup_{j=1}^k\nbhd{U}'_{x_j}$
and let $C=\max\{C_{x_1},\dots,C_{x_k}\}$ and
$r=\min\{r_{x_1},\dots,r_{x_k}\}$\@.  Then, if $x\in K$\@, we have
$x\in\nbhd{U}'_{x_j}$ for some $j\in\{1,\dots,k\}$ and so
\begin{equation*}
\dnorm{j_m\xi(x)}_{\ol{\metric}_m}\le C_{x_j}r_{x_j}^{-m}\le Cr^{-m},
\end{equation*}
as desired.

\eqref{pl:rabound2}$\implies$\eqref{pl:rabound2} Let $x\in\man{M}$ and let
$(\nbhd{V},\psi)$ be a vector bundle chart for $\man{E}$ such that the
associated chart $(\nbhd{U},\phi)$ for $\man{M}$ is a relatively compact
coordinate chart about $x$\@.  Let
$\map{\vect{\xi}}{\phi(\nbhd{U})}{\real^k}$ be the local representative of
$\xi$\@.  By hypothesis, there exist $C,r\in\realp$ such that
$\dnorm{j_m\xi(x')}_{\ol{\metric}_m}\le Cr^{-m}$ for every $m\in\integernn$
and $x'\in\nbhd{U}$\@.  Let $\nbhd{U}'$ be a relatively compact neighbourhood
of $x$ such that $\closure(\nbhd{U}')\subset\nbhd{U}$\@.  By
Lemma~\ref{lem:pissy-estimate}\@, there exist $B,\sigma\in\realp$ such that
\begin{equation*}
\snorm{\linder[I]{\xi^a}(\vect{x}')}\le BI!\sigma^{-\snorm{I}}
\end{equation*}
for every $a\in\{1,\dots,k\}$\@, $\vect{x}'\in\closure(\nbhd{U}')$\@, and
$I\in\integernn^n$\@.  We conclude real analyticity of $\xi$ in a
neighbourhood of $x$ by~\cite[Proposition~2.2.10]{SGK/HRP:02}\@.
\end{proof}
\end{lemma}

\section{The compact-open topologies for the spaces of finitely
differentiable, Lipschitz, and smooth vector fields}\label{sec:smooth-topology}

In Sections~\ref{sec:time-varying} and~\ref{sec:systems} we will look
carefully at two related things:~(1)~time-varying vector fields
and~(2)~control systems.  In doing so, we focus on structure that allows us
to prove useful properties such as regular dependence of flows on initial
conditions.  Also, in our framework of \gcs{}s in Section~\ref{sec:gcs}\@, we
will need to impose structure on systems where we have carefully eliminated
the usual structure of a control parameterisation.  To do this, we use the
topological structure of sets of vector fields in an essential way.  In this
and the subsequent two sections we describe appropriate topologies for
finitely differentiable, Lipschitz, smooth, holomorphic, and real analytic
vector fields.  The topology we use in this section in the smooth case (and
the easily deduced finitely differentiable case) is classical, and is
described, for example, in \cite[\S2.2]{AAA/YS:04}\@; see
also~\cite[Chapter~4]{PWM:80}\@.  What we do that is original is provide a
characterisation of the seminorms for this topology using the jet bundle
fibre metrics from Section~\ref{subsec:olGm}\@.  The fruits of the effort
expended in the next three sections is harvested in the remainder of the
paper, where our concrete definitions of seminorms permit a relatively
unified analysis in Sections~\ref{sec:time-varying} and~\ref{sec:systems} of
time-varying vector fields and control systems.  Also, the treatment of our
new class of systems in Section~\ref{sec:gcs} is made relatively simple by
our descriptions of topologies for spaces of vector fields.

One facet of our presentation that is novel is that we flesh out completely
the ``weak-$\sL$'' characterisations of topologies for vector fields.  These
topologies characterise vector fields by how they act on functions through
Lie differentiation.  The use of such ``weak'' characterisations is
commonplace~\cite[\eg][]{AAA/YS:04,HJS:97b}\@, although the equivalence with
strong characterisation is not typically proved; indeed, we know of no
existing proofs of our
Theorems~\ref{the:COinfty-weak}\@,~\ref{the:COm-weak}\@,%
~\ref{the:COmm'-weak}\@, and~\ref{the:Comega-weak}\@.  We show that, for the
issues that come up in this paper, the weak characterisations for vector
field topologies agree with the direct ``strong'' characterisations.  This
requires some detailed knowledge of the topologies we use.

While our primary interest is in vector fields,~\ie~sections of the tangent
bundle, it is advantageous to work instead with topologies for sections of
general vector bundles, and then specialise to vector fields.  We will also
work with topologies for functions, but this falls out easily from the
general vector bundle treatment.

\subsection{General smooth vector bundles}\label{subsec:COinfty-vb}

We let $\map{\pi}{\man{E}}{\man{M}}$ be a smooth vector bundle with
$\nabla^0$ a linear connection on $\man{E}$\@, $\nabla$ an affine connection
on $\man{M}$\@, $\metric_0$ a fibre metric on $\man{E}$\@, and $\metric$ a
Riemannian metric on $\man{M}$\@.  This gives us, as in Section~\ref{subsec:olGm}\@,
fibre metrics $\ol{\metric}_m$ on the jet bundles $\jet{m}{\man{E}}$\@,
$m\in\integernn$\@, and corresponding fibre norms
$\dnorm{\cdot}_{\ol{\metric}_m}$\@.

For a compact set $K\subset\man{M}$ we now define a seminorm $p^\infty_{K,m}$
on $\sections[\infty]{\man{E}}$ by
\begin{equation*}
p^\infty_{K,m}(\xi)=\sup\setdef{\dnorm{j_m\xi(x)}_{\ol{\metric}_m}}{x\in K}.
\end{equation*}
The locally convex topology on $\sections[\infty]{\tb{\man{M}}}$ defined by
the family of seminorms $p^\infty_{K,m}$\@, $K\subset\man{M}$ compact,
$m\in\integernn$\@, is called the \defn{smooth compact open} or
\defn{$\CO^\infty$-topology} for $\sections[\infty]{\man{E}}$\@.

We comment that the seminorms depend on the choices of $\nabla$\@,
$\nabla^0$\@, $\metric$\@, and $\metric_0$\@, but the $\CO^\infty$-topology
is independent of these choices.  We will constantly throughout the paper use
these seminorms, and in doing so we will automatically be assuming that we
have selected the linear connection $\nabla^0$\@, the affine connection
$\nabla$\@, the fibre metric $\metric_0$\@, and the Riemannian metric
$\metric$\@.  We will do this often without explicit mention of these objects
having been chosen.

\subsection{Properties of the $\CO^\infty$-topology}\label{subsec:COinfty-props}

Let us say a few words about the $\CO^\infty$-topology, referring to
references for details.  The locally convex $\CO^\infty$-topology has the
following attributes.
\begin{compactenum}[$\CO^\infty$-1.]
\item It is Hausdorff:~\cite[4.3.1]{PWM:80}\@.
\item \label{enum:COinfty-complete} It is complete:~\cite[4.3.2]{PWM:80}\@.
\item It is metrisable:~\cite[4.3.1]{PWM:80}\@.
\item \label{enum:COinfty-separable} It is separable: We could not find this
stated anywhere, but here's a sketch of a proof.  By embedding $\man{E}$ in
Euclidean space $\real^N$ and, using an argument like that for real analytic
vector bundles in the proof of Lemma~\ref{lem:analytic-conn}\@, we regard
$\man{E}$ as a subbundle of a trivial bundle over the submanifold
$\man{M}\subset\real^N$\@.  In this case, we can reduce our claim of
separability of the $\CO^\infty$-topology to that for smooth functions on
submanifolds of $\real^N$\@.  Here we can argue as follows.  If
$K\subset\man{M}$ is compact, it can be contained in a compact cube $C$ in
$\real^N$\@.  Then we can use a cutoff function to take any smooth function
on $\man{M}$ and leave it untouched on a neighbourhood of $K$\@, but have it
and all of its derivatives vanish outside a compact set contained in
$\interior(C)$\@.  Then we can use Fourier series to approximate in the
$\CO^\infty$-topology~\cite[Theorem~VII.2.11(b)]{EMS/GW:71}\@.  Since there
are countably many Fourier basis functions, this gives the desired
separability.
\item \label{enum:COinfty-nuclear} It is nuclear:\footnote{There are several
ways of characterising nuclear spaces.  Here is one.  A continuous linear
mapping $\map{L}{\alg{E}}{\alg{F}}$ between Banach spaces is \defn{nuclear}
if there exist sequences $\ifam{v_j}_{j\in\integerp}$ in $\alg{F}$ and
$\ifam{\alpha_j}_{j\in\integerp}$ in $\topdual{\alg{E}}$ such that
$\sum_{j\in\integerp}\dnorm{\alpha_j}\dnorm{v_j}<\infty$ and such that
\begin{equation*}
L(u)=\sum_{j=1}^\infty\alpha_j(u)v_j,
\end{equation*}
the sum converging in the topology of $\alg{V}$\@.  Now suppose that
$\alg{V}$ is a locally convex space and $p$ is a continuous seminorm on
$\alg{V}$\@.  We denote by $\ol{\alg{V}}_p$ the completion of
\begin{equation*}
\alg{V}/\setdef{v\in\alg{V}}{p(v)=0};
\end{equation*}
thus $\ol{\alg{V}}_p$ is a Banach space.  The space $\alg{V}$ is \defn{nuclear}
if, for any continuous seminorm $p$\@, there exists a continuous seminorm $q$
satisfying $q\le p$ such that the mapping
\begin{equation*}
\mapdef{i_{p,q}}{\ol{\alg{V}}_p}{\ol{\alg{V}}_q}
{v+\setdef{v'\in\alg{V}}{p(v)=0}}
{v+\setdef{v'\in\alg{V}}{q(v)=0}}
\end{equation*}
is nuclear.  It is to be understood that this definition is essentially
meaningless at a first encounter, so we refer to~\cite{HH-N/VBM:81,AP:69} and
relevant sections of~\cite{HJ:81} to begin understanding the notion of a
nuclear space.  The only attribute of nuclear spaces of interest to us here
is that their relatively compact subsets are exactly the von Neumann bounded
subsets~\cite[Proposition~4.47]{AP:69}\@.}~\cite[Theorem~21.6.6]{HJ:81}\@.
\item \label{enum:COinfty-suslin} It is Suslin:\footnote{\label{fn:suslin}A
\defn{Polish space} is a complete separable metrisable space.  A \defn{Suslin
space} is a continuous image of a Polish space.  A good reference for the
basic properties of Suslin spaces is~\cite[Chapter~6]{VIB:07b}\@.} This
follows since $\sections[\infty]{\tb{\man{M}}}$ is a Polish space (see
footnote~\ref{fn:suslin}), as we have already seen.
\end{compactenum}
Some of these attributes perhaps seem obscure, but we will, in fact, use all
of them!

Since the $\CO^\infty$-topology is metrisable, it is exactly characterised by
its convergent sequences, so let us describe these.  A sequence
$\ifam{\xi_k}_{k\in\integerp}$ in $\sections[\infty]{\man{E}}$ converges to
$\xi\in\sections[\infty]{\man{E}}$ if and only if, for each compact set
$K\subset\man{M}$ and for each $m\in\integernn$\@, the sequence
$\ifam{j_m\xi_k|K}_{k\in\integerp}$ converges uniformly to
$j_m\xi|K$\@,~\cf~combining~\cite[Theorem~46.8]{JRM:00}
and~\cite[Lemma~4.2]{PWM:80}\@.

Since the topology is nuclear, it follows that subsets of
$\sections[\infty]{\tb{\man{M}}}$ are compact if and only if they are closed
and von Neumann bounded~\cite[Proposition~4.47]{AP:69}\@.  That is to say, in
a nuclear locally convex space, the compact bornology and the von Neumann
bornology agree, according to the terminology introduced in
Section~\ref{subsec:notation}\@.  It is then interesting to characterise von
Neumann bounded subsets of $\sections[\infty]{\man{E}}$\@.  One can show that
a subset $\nbhd{B}$ is bounded in the von Neumann bornology if and only if
every continuous seminorm on $\alg{V}$ is a bounded function when restricted
to $\nbhd{B}$~\cite[Theorem~1.37(b)]{WR:91}\@.  Therefore, to characterise
von Neumann bounded subsets, we need only characterise subsets on which each
of the seminorms $p^\infty_{K,m}$ is a bounded function.  This obviously
gives the following characterisation.
\begin{lemma}\label{lem:COinftybdd}
A subset\/ $\nbhd{B}\subset\sections[\infty]{\man{E}}$ is bounded in the von
Neumann bornology if and only if the following property holds: for any
compact set\/ $K\subset\man{M}$ and any\/ $m\in\integernn$\@, there exists\/
$C\in\realp$ such that\/ $p^\infty_{K,m}(\xi)\le C$ for every\/
$\xi\in\nbhd{B}$\@.
\end{lemma}

Let us give a coordinate characterisation of the smooth compact-open
topology, just for concreteness and so that the reader can see that our
constructions agree with perhaps more familiar things.  If we have a smooth
vector bundle $\map{\pi}{\man{E}}{\man{M}}$\@, we let $(\nbhd{V},\psi)$ be a
vector bundle chart for $\man{E}$ inducing a chart $(\nbhd{U},\phi)$ for
$\man{M}$\@.  For $\xi\in\sections[\infty]{\man{E}}$\@, the local
representative of $\xi$ has the form
\begin{equation*}
\real^n\supset\phi(\nbhd{U})\ni\vect{x}\mapsto
(\vect{x},\vect{\xi}(\vect{x}))\in\phi(\nbhd{U})\times\real^k.
\end{equation*}
Thus we have an associated map $\map{\vect{\xi}}{\phi(\nbhd{U})}{\real^k}$
that describes the section locally.  A \defn{$\CO^\infty$-subbasic
neighbourhood} is a subset $\nbhd{B}_\infty(\xi,\nbhd{V},K,\epsilon,m)$ of
$\sections[\infty]{\man{E}}$\@, where
\begin{compactenum}
\item $\xi\in\sections[\infty]{\man{E}}$\@,
\item $(\nbhd{V},\psi)$ is a vector bundle chart for $\man{E}$ with
associated chart $(\nbhd{U},\phi)$ for $\man{M}$\@,
\item $K\subset\nbhd{U}$ is compact,
\item $\epsilon\in\realp$\@,
\item $m\in\integernn$\@, and
\item $\eta\in\nbhd{B}_\infty(\xi,\nbhd{V},K,\epsilon,m)$ if and only if
\begin{equation*}
\dnorm{\linder[l]{\vect{\eta}}(\vect{x})-\linder[l]{\vect{\xi}(\vect{x})}}
<\epsilon,\qquad\vect{x}\in\phi(K),\ l\in\{0,1,\dots,m\},
\end{equation*}
where $\map{\vect{\xi},\vect{\eta}}{\phi(\nbhd{U})}{\real^k}$ are the local
representatives.
\end{compactenum}
One can show that the $\CO^\infty$-topology is that topology having as a
subbase the $\CO^\infty$-subbasic neighbourhoods.  This is the definition
used by~\cite{MWH:76}\@, for example.  To show that this topology agrees with
our intrinsic characterisation is a straightforward bookkeeping chore, and
the interested reader can refer to Lemma~\ref{lem:pissy-estimate} to see how
this is done in the more difficult real analytic case.  This more concrete
characterisation using vector bundle charts can be useful should one ever
wish to verify some properties in examples.  It can also be useful in general
arguments in emergencies when one does not have the time to flesh out
coordinate-free constructions.

\subsection{The weak-$\sL$ topology for smooth vector fields}\label{subsec:COinfty-weak}

The $\CO^\infty$-topology for smooth sections of a vector bundle, merely by
specialisation, gives a locally convex topology on the set
$\sections[\infty]{\tb{\man{M}}}$ of smooth vector fields and the set
$\func[\infty]{\man{M}}$ of smooth functions (noting that a smooth function
is obviously identified with a section of the trivial vector bundle
$\man{M}\times\real$).  The only mildly interesting thing in these cases is
that one does not need a separate linear connection in the vector bundles or
a separate fibre metric.  Indeed, $\tb{\man{M}}$ is already assumed to have a
linear connection (the affine connection on $\man{M}$) and a fibre metric
(the Riemannian metric on $\man{M}$), and the trivial bundle has the
canonical flat linear connection defined by $\nabla_Xf=\lieder{X}{f}$ and the
standard fibre metric induced by absolute value on the fibres.

We wish to see another way of describing the $\CO^\infty$-topology on
$\sections[\infty]{\tb{\man{M}}}$ by noting that a vector field defines a
linear map, indeed a derivation, on $\func[\infty]{\man{M}}$ by Lie
differentiation: $f\mapsto\lieder{X}{f}$\@.  The topology we describe for
$\sections[\infty]{\tb{\man{M}}}$ is a sort of weak topology arising from the
$\CO^\infty$-topology on $\func[\infty]{\man{M}}$ and Lie differentiation.
To properly set the stage for the fact that we will repeat this construction
for our other topologies, it is most clear to work in a general setting for a
moment, and then specialise in each subsequent case.

The general setup is provided by the next definition.
\begin{definition}\label{def:weakA}
Let $\field\in\{\real,\complex\}$ and let $\alg{U}$ and $\alg{V}$ be
$\field$-vector spaces with $\alg{V}$ locally convex.  Let
$\sA\subset\Hom_{\real}(\alg{U};\alg{V})$ and let the \defn{weak-$\sA$
topology} on $\alg{U}$ be the weakest topology for which $A$ is continuous
for every $A\in\sA$~\cite[\S2.11]{JH:66}\@.

Also let $(\ts{X},\sO)$ be a topological space, let $(\ts{T},\sM)$ be a
measurable space, and let $\map{\mu}{\sM}{\erealnn}$ be a finite measure.  We
have the following notions:
\begin{compactenum}[(i)]
\item a subset $\nbhd{B}\subset\alg{U}$ is \defn{weak-$\sA$ bounded in the
von Neumann bornology} if $A(\nbhd{B})$ is bounded in the von Neumann
bornology for every $A\in\sA$\@;
\item a map $\map{\Phi}{\ts{X}}{\alg{U}}$ is \defn{weak-$\sA$ continuous} if
$A\scirc\Phi$ is continuous for every $A\in\sA$\@;
\item a map $\map{\Psi}{\ts{T}}{\alg{U}}$ is \defn{weak-$\sA$ measurable} if
$A\scirc\Psi$ is measurable for every $A\in\sA$\@;
\item a map $\map{\Psi}{\ts{T}}{\alg{U}}$ is \defn{weak-$\sA$ Bochner
integrable} with respect to $\mu$ if $A\scirc\Psi$ is Bochner integrable with
respect to $\mu$ for every $A\in\sA$\@.\oprocend
\end{compactenum}
\end{definition}

As can be seen in Section~2.11 of~\cite{JH:66}\@, the weak-$\sA$ topology is
a locally convex topology, and a subbase for open sets in this topology is
\begin{equation*}
\setdef{A^{-1}(\nbhd{O})}{A\in\sA,\ \nbhd{O}\subset\alg{V}\ \textrm{open}}.
\end{equation*}
Equivalently, the weak-$\sA$ topology is defined by the seminorms
\begin{equation*}
u\mapsto q(A(u)),\qquad A\in\sA,\ q\ \textrm{a continuous seminorm for}\
\alg{V}.
\end{equation*}
This is a characterisation of the weak-$\sA$ topology we will use often.

We now have the following result which gives conditions for the equivalence
of ``weak-$\sA$'' notions with the usual notions.  We call a subset
$\sA\subset\Hom_{\field}(\alg{U};\alg{V})$ \defn{point separating} if, given
distinct $u_1,u_2\in\alg{U}$\@, there exists $A\in\sA$ such that
$A(u_1)\not=A(u_2)$\@.
\begin{lemma}\label{lem:weakA}
Let\/ $\field\in\{\real,\complex\}$ and let\/ $\alg{U}$ and\/ $\alg{V}$ be
locally convex\/ $\field$-vector spaces.  Let\/
$\sA\subset\Hom_{\real}(\alg{U};\alg{V})$ and suppose that the weak-$\sA$
topology agrees with the locally convex topology for\/ $\alg{U}$\@.  Let\/
$(\ts{X},\sO)$ be a topological space, let\/ $(\ts{T},\sM)$ be a measurable
space, and let\/ $\map{\mu}{\sM}{\erealnn}$ be a finite measure.  Then the
following statements hold:
\begin{compactenum}[(i)]
\item \label{pl:weakAbdd} a subset\/ $\nbhd{B}\subset\alg{U}$ is bounded in
the von Neumann bornology if and only if it is weak-$\sA$ bounded in the von
Neumann bornology;
\item \label{pl:weakAcont} a map\/ $\map{\Phi}{\ts{X}}{\alg{U}}$ is
continuous if and only if it is weak-$\sA$ continuous;
\item \label{pl:weakAmeas} for a map\/ $\map{\Psi}{\ts{T}}{\alg{U}}$\@,
\begin{compactenum}[(a)]
\item if\/ $\Psi$ is measurable, then it is weak-$\sA$ measurable;
\item if\/ $\alg{U}$ and\/ $\alg{V}$ are Hausdorff Suslin spaces, if\/ $\sA$
contains a countable point separating subset, and if\/ $\Psi$ is weak-$\sA$
measurable, then\/ $\Psi$ is measurable;
\end{compactenum}
\item \label{pl:weakAint} if\/ $\alg{U}$ is complete and separable, a map\/
$\map{\Psi}{\ts{T}}{\alg{U}}$ is Bochner integrable with respect to\/ $\mu$
if and only if it is weak-$\sA$ Bochner integrable with respect to\/ $\mu$\@.
\end{compactenum}
\begin{proof}
\eqref{pl:weakAbdd} and \eqref{pl:weakAcont}\@: Both of these assertions
follows directly from the fact that the locally convex topology of $\alg{U}$
agrees with the weak-$\sA$ topology.  Indeed, the equivalence of these
topologies implies that~(a)~if $p$ is a continuous seminorm for the locally
convex topology of $\alg{U}$\@, then there exist continuous seminorms
$q_1,\dots,q_k$ for $\alg{V}$ and $A_1,\dots,A_k\in\sA$ such that
\begin{equation}\label{eq:weakA1}
p(u)\le q_1(A_1(u))+\dots+q_k(A_k(u)),\qquad u\in\alg{U},
\end{equation}
and~(b)~if $q$ is a continuous seminorm for $\alg{V}$ and if $A\in\sA$\@,
then there exists a continuous seminorm $p$ for the locally convex topology
for $\alg{U}$ such that
\begin{equation}\label{eq:weakA2}
q(A(u))\le p(u),\qquad u\in\alg{U}.
\end{equation}

\eqref{pl:weakAmeas} First suppose that $\Psi$ is measurable and let
$A\in\sA$\@.  Since the locally convex topology of $\alg{U}$ agrees with the
weak-$\sA$ topology, $A$ is continuous in the locally convex topology of
$\alg{U}$\@.  Therefore, if $\Psi$ is measurable, it follows immediately by
continuity of $A$ that $A\scirc\Psi$ is measurable.

Next suppose that $\alg{U}$ and $\alg{V}$ are Suslin, that $\sA$ contains a
countable point separating subset, and that $\Psi$ is weak-$\sA$ measurable.
Without loss of generality, let us suppose that $\sA$ is itself countable.
By $\alg{V}^{\sA}$ we denote the mappings from $\sA$ to $\alg{V}$\@, with the
usual pointwise vector space structure.  A typical element of $\alg{V}^{\sA}$
we denote by $\phi$\@.  By~\cite[Lemma~6.6.5(iii)]{VIB:07b}\@, $\alg{V}^{\sA}$
is a Suslin space.  Let us define a mapping
$\map{\iota_{\sA}}{\alg{U}}{\alg{V}^{\sA}}$ by $\iota_{\sA}(u)(A)=A(u)$\@.
Since $\sA$ is point separating, we easily verify that $\iota_{\sA}$ is
injective, and so we have $\alg{U}$ as a subspace of the countable product
$\alg{V}^{\sA}$\@.  For $A\in\sA$ let $\map{\pr_A}{\alg{V}^{\sA}}{\alg{V}}$
be the projection defined by $\pr_A(\phi)=\phi(A)$\@.  Since $\alg{V}$ is
Suslin, it is hereditary Lindel\"of~\cite[Lemma~6.6.4]{VIB:07b}\@.  Thus the
Borel $\sigma$-algebra of $\alg{V}^{\sA}$ is the same as the initial Borel
$\sigma$-algebra defined by the projections $\pr_A$\@, $A\in\sA$\@,~\ie~the
smallest $\sigma$-algebra for which the projections are
measurable~\cite[Lemma~6.4.2]{VIB:07b}\@.  By hypothesis,
$(A\scirc\Psi)^{-1}(\nbhd{B})$ is measurable for every $A\in\sA$ and every
Borel set $\nbhd{B}\subset\alg{V}$\@.  Now we note that
$\pr_A\scirc\iota_{\sA}(v)=A(v)$\@, from which we deduce that
\begin{equation*}
(A\scirc\Psi)^{-1}(\nbhd{B})=
(\iota_{\sA}\scirc\Psi)^{-1}(\pr_A^{-1}(\nbhd{B}))
\end{equation*}
is measurable for every $A\in\sA$ and every Borel set
$\nbhd{B}\subset\alg{V}$\@.  Thus $\iota_{\sA}\scirc\Psi$ is measurable.

Since $\alg{U}$ is Suslin, by definition there is a Polish space $\ts{P}$ and
a continuous surjection $\map{\sigma}{\ts{P}}{\alg{U}}$\@.  If
$\nbhd{C}\subset\alg{U}$ is a Borel set, then
$\sigma^{-1}(\nbhd{C})\subset\ts{P}$ is a Borel set.  Note that $\iota_{\sA}$
is continuous (since $\pr_A\scirc\iota_{\sA}$ is continuous for every
$A\in\sA$) and so is a Borel mapping.  By~\cite[Theorem~423I]{DHF:06}\@, we
have that $\iota_{\sA}\scirc\sigma(\sigma^{-1}(\nbhd{C}))\subset\alg{V}$ is
Borel.  Since $\sigma$ is surjective, this means that
$\iota_{\sA}(\nbhd{C})\subset\alg{V}$ is Borel.  Finally, since
\begin{equation*}
\Psi^{-1}(\nbhd{C})=(\iota_{\sA}\scirc\Psi)^{-1}(\iota_{\sA}(\nbhd{C})),
\end{equation*}
measurability of $\Psi$ follows.

\eqref{pl:weakAint} Since $\alg{U}$ is separable and complete, by
\citet[Theorems~3.2 and~3.3]{RB/AD:11} Bochner integrability of $\Psi$ is
equivalent to integrability, in the sense of Lebesgue, of $t\mapsto
p\scirc\Psi(t)$ for any continuous seminorm $p$\@.  Thus, $\Psi$ is Bochner
integrable with respect to the locally convex topology of $\alg{U}$ if and
only if $t\mapsto p\scirc\Psi(t)$ is integrable, and $\Psi$ is weak-$\sA$
Bochner integrable if and only if $t\mapsto q_A(\Psi(t))$ is integrable for
every $A\in\sA$\@.  This part of the proof now follows from the
inequalities~\eqref{eq:weakA1} and~\eqref{eq:weakA2} that characterise the
equivalence of the locally convex and weak-$\sA$ topologies for $\alg{U}$\@.
\end{proof}
\end{lemma}

The proof of the harder direction in part~\eqref{pl:weakAmeas} is an
adaptation of~\cite[Theorem~1]{GEFT:75} to our more general setting.  We will
revisit this idea again when we talk about measurability of time-varying
vector fields in Section~\ref{sec:time-varying}\@.

For $f\in\func[\infty]{\man{M}}$\@, let us define
\begin{equation*}
\mapdef{\sL_f}{\sections[\infty]{\tb{\man{M}}}}{\func[\infty]{\man{M}}}
{X}{\lieder{X}{f}.}
\end{equation*}
The topology for $\sections[\infty]{\tb{\man{M}}}$ we now define corresponds
to the general case of Definition~\ref{def:weakA} by taking
$\alg{U}=\sections[\infty]{\tb{\man{M}}}$\@,
$\alg{V}=\func[\infty]{\man{M}}$\@, and
$\sA=\setdef{\sL_f}{f\in\func[\infty]{\man{M}}}$\@.  To this end, we make the
following definition.
\begin{definition}
For a smooth manifold $\man{M}$\@, the \defn{weak-$\sL$ topology} for
$\sections[\infty]{\tb{\man{M}}}$ is the weakest topology for which $\sL_f$
is continuous for every $f\in\func[\infty]{\man{M}}$\@, if
$\func[\infty]{\man{M}}$ has the $\CO^\infty$-topology.\oprocend
\end{definition}

We now have the following result.
\begin{theorem}\label{the:COinfty-weak}
For a smooth manifold, the following topologies for\/
$\sections[\infty]{\tb{\man{M}}}$ agree:
\begin{compactenum}[(i)]
\item \label{pl:COinfty-weak1} the\/ $\CO^\infty$-topology;
\item \label{pl:COinfty-weak2} the weak-$\sL$ topology.
\end{compactenum}
\begin{proof}
\eqref{pl:COinfty-weak1}$\subset$\eqref{pl:COinfty-weak2} For this part of
the proof, we assume that $\man{M}$ has a well-defined dimension.  The proof
is easily modified by additional notation to cover the case where this may
not hold.  Let $K\subset\man{M}$ be compact and let $m\in\integernn$\@.  Let
$x\in K$ and let $(\nbhd{U}_x,\phi_x)$ be a coordinate chart for $\man{M}$
about $x$ with coordinates denoted by $(x^1,\dots,x^n)$\@.  Let
$\map{\vect{X}}{\phi_x(\nbhd{U}_x)}{\real^n}$ be the local representative of
$X\in\sections[\infty]{\tb{\man{M}}}$\@.  For $j\in\{1,\dots,n\}$ let
$f^j_x\in\func[\infty]{\man{M}}$ have the property that, for some relatively
compact neighbourhood $\nbhd{V}_x$ of $x$ with
$\closure(\nbhd{V}_x)\subset\nbhd{U}_x$\@, $f^j_x=x^j$ for $y$ in some
neighbourhood of $\closure(\nbhd{V}_x)$\@.  (This is done using standard
extension arguments for smooth
functions,~\cf~\cite[Proposition~5.5.8]{RA/JEM/TSR:88}\@.)  Then, in a
neighbourhood of $\closure(\nbhd{V}_x)$ in $\nbhd{U}_x$\@, we have
$\lieder{X}{f^j_x}=X^j$\@.  Therefore, for each $y\in\closure(\nbhd{V}_x)$\@,
\begin{equation*}
j_mX(y)\mapsto\sum_{j=1}^n\dnorm{j_m(\lieder{X}{f^j_x})(y)}_{\ol{\metric}_m}
\end{equation*}
is a norm on the fibre $\jet[y]{m}{\man{E}}$\@.  Therefore, there exists
$C_x\in\realp$ such that
\begin{equation*}
\dnorm{j_mX(y)}_{\ol{\metric}_m}\le
C_x\sum_{j=1}^n\dnorm{j_m(\lieder{X}{f^j_x})(y)}_{\ol{\metric}_m},
\qquad y\in\closure(\nbhd{V}_x).
\end{equation*}
Since $K$ is compact, let $x_1,\dots,x_k\in K$ be such that
$K\subset\cup_{a=1}^k\nbhd{V}_{x_a}$\@.  Let
\begin{equation*}
C=\max\{C_{x_1},\dots,C_{x_r}\}.
\end{equation*}
Then, if $y\in K$ we have $y\in\nbhd{V}_{x_a}$ for some
$a\in\{1,\dots,r\}$\@, and so
\begin{equation*}
\dnorm{j_mX(y)}_{\ol{\metric}_m}\le
C\sum_{j=1}^n\dnorm{j_m(\lieder{X}{f^j_{x_a}})(y)}_{\ol{\metric}_m}\le
C\sum_{a=1}^r\sum_{j=1}^n\dnorm{j_m(\lieder{X}{f^j_{x_a}})(y)}_{\ol{\metric}_m}.
\end{equation*}
Taking supremums over $y\in K$ gives
\begin{equation*}
p^\infty_{K,m}(X)\le
C\sum_{a=1}^r\sum_{j=1}^np^\infty_{K,m}(\lieder{X}{f^j_{x_a}}),
\end{equation*}
This part of the theorem then follows since the weak-$\sL$ topology, as we
indicated following Definition~\ref{def:weakA} above, is defined by the
seminorms
\begin{equation*}
X\mapsto p^\infty_{K,m}(\lieder{X}{f}),\qquad K\subset\man{M}\
\textrm{compact},\ m\in\integernn,\ f\in\func[\infty]{\man{M}}.
\end{equation*}

\eqref{pl:COinfty-weak2}$\subset$\eqref{pl:COinfty-weak1} As per~\eqref{eq:nabla(j)}\@, let us abbreviate
\begin{equation*}
\nabla^j(\dots(\nabla^1(\nabla^0A)))=\nabla^{(j)}A,
\end{equation*}
where $A$ can be either a vector field or one-form, in what we will need.
Since covariant differentials commute with
contractions~\cite[Theorem~7.03(F)]{CTJD/TP:91}\@, an elementary induction
argument gives the formula
\begin{equation}\label{eq:nablafX}
\nabla^{(m-1)}(\d{f}(X))=\sum_{j=0}^m\binom{m}{j}
C_{1,m-j+1}((\nabla^{(m-j-1)}X)\otimes(\nabla^{(j-1)}\d{f})),
\end{equation}
where $C_{1,m-j+1}$ is the contraction defined by
\begin{multline*}
C_{1,m-j+1}(v\otimes\alpha^1\otimes\dots\otimes\alpha^{m-j}
\otimes\alpha^{m-j+1}\otimes\alpha^{m-j+2}\otimes\dots\otimes\alpha^{m+1})\\
=(\alpha^{m-j+1}(v))(\alpha^1\otimes\dots\otimes\alpha^{m-j}
\otimes\alpha^{m-j+2}\otimes\dots\otimes\alpha^{m+1}).
\end{multline*}
In writing~\eqref{eq:nablafX} we use the convention $\nabla^{(-1)}X=X$ and
$\nabla^{(-1)}(\d{f})=\d{f}$\@.  Next we claim that $\sL_f$ is continuous for
every $f\in\func[\infty]{\man{M}}$ if $\sections[\infty]{\tb{\man{M}}}$ is
provided with $\CO^\infty$-topology.  Indeed, let $K\subset\man{M}$\@, let
$m\in\integerp$\@, and let $f\in\func[\infty]{\man{M}}$\@.
By~\eqref{eq:nablafX} (after a few moments of thought), we have, for some
suitable $M_0,M_1\dots,M_m\in\realp$\@,
\begin{equation*}
p^\infty_{K,m}(\lieder{X}{f})\le
\sum_{j=0}^mM_{m-j}p^\infty_{K,m-j}(X)p^\infty_{K,j+1}(f)
\le\sum_{j=0}^mM'_jp^\infty_{K,j}(X).
\end{equation*}
This gives continuity of the identity map, if we provide the domain with the
$\CO^\infty$-topology and the codomain with the weak-$\sL$
topology,~\cf~\cite[\S{}III.1.1]{HHS/MPW:99}\@.  Thus open sets in the
weak-$\sL$ topology are contained in the $\CO^\infty$-topology.
\end{proof}
\end{theorem}

With respect to the concepts of interest to us, this gives the following
result.
\begin{corollary}\label{cor:COinftyweak}
Let\/ $\man{M}$ be a smooth manifold, let\/ $(\ts{X},\sO)$ be a topological
space, let\/ $(\ts{T},\sM)$ be a measurable space, and let\/
$\map{\mu}{\sM}{\erealnn}$ be a finite measure.  The following statements
hold:
\begin{compactenum}[(i)]
\item \label{pl:COinftyweakbdd} a subset\/
$\nbhd{B}\subset\sections[\infty]{\tb{\man{M}}}$ is bounded in the von
Neumann bornology if and only if it is weak-$\sL$ bounded in the von Neumann
bornology;
\item \label{pl:COinftyweakcont} a map\/
$\map{\Phi}{\ts{X}}{\sections[\infty]{\tb{\man{M}}}}$ is continuous if and
only if it is weak-$\sL$ continuous;
\item \label{pl:COinftyweakmeas} a map\/
$\map{\Psi}{\ts{T}}{\sections[\infty]{\tb{\man{M}}}}$ is measurable if and
only if it is weak-$\sL$ measurable;
\item \label{pl:COinftyweakint} a map\/
$\map{\Psi}{\ts{T}}{\sections[\infty]{\tb{\man{M}}}}$ is Bochner integrable
if and only if it is weak-$\sL$ Bochner integrable.
\end{compactenum}
\begin{proof}
We first claim that $\sA\eqdef\setdef{\sL_f}{f\in\func[\infty]{\man{M}}}$ has
a countable point separating subset.  This is easily proved as follows.  For
notational simplicity, suppose that $\man{M}$ has a well-defined dimension.
Let $x\in\man{M}$ and note that there exist a neighbourhood $\nbhd{U}_x$ of
$x$ and $f^1_x,\dots,f^n_x\in\func[\infty]{\man{M}}$ such that
\begin{equation*}
\ctb[y]{\man{M}}=\vecspan[\real]{\d{f^1}(y),\dots,\d{f^n}(y)},\qquad
y\in\nbhd{U}_x.
\end{equation*}
Since $\man{M}$ is second countable it is
Lindel\"of~\cite[Theorem~16.9]{SW:04}\@.  Therefore, there exists
$\ifam{x_j}_{j\in\integerp}$ such that
$\man{M}=\cup_{j\in\integerp}\nbhd{U}_{x_j}$\@.  The countable collection of
linear mappings $\sL_{f^k_{x_j}}$\@, $k\in\{1,\dots,n\}$\@,
$j\in\integerp$\@, is then point separating.  Indeed, if
$X,Y\in\sections[\infty]{\tb{\man{M}}}$ are distinct, then there exists
$x\in\man{M}$ such that $X(x)\not=Y(x)$\@.  Let $j\in\integerp$ be such that
$x\in\nbhd{U}_{x_j}$ and note that we must have
$\sL_{f^k_{x_j}}(X)(x)\not=\sL_{f^k_{x_j}}(Y)(x)$ for some
$k\in\{1,\dots,n\}$\@, giving our claim.

The result is now a direct consequence of Lemma~\ref{lem:weakA}\@, noting
that the $\CO^\infty$-topology on $\sections[\infty]{\tb{\man{M}}}$ is
complete, separable, and Suslin (we also need that the $\CO^\infty$-topology
on $\func[\infty]{\man{M}}$ is Suslin, which it is), as we have seen above in
properties~$\CO^\infty$-\ref{enum:COinfty-complete}\@,%
~$\CO^\infty$-\ref{enum:COinfty-separable}\@,
and~$\CO^\infty$-\ref{enum:COinfty-suslin}\@.
\end{proof}
\end{corollary}

\subsection{Topologies for finitely differentiable vector fields}\label{subsec:COm-topology}

The constructions of this section so far are easily adapted to the case where
objects are only finitely differentiable.  We sketch here how this can be
done.  We let $\map{\pi}{\man{E}}{\man{M}}$ be a smooth vector bundle, and we
suppose that we have a linear connection $\nabla^0$ on $\man{E}$\@, an affine
connection $\nabla$ on $\man{M}$\@, a fibre metric $\metric_0$ on
$\man{E}$\@, and a Riemannian metric $\metric$ on $\man{M}$\@.  Let
$r\in\integernn\cup\{\infty\}$ and let $m\in\integernn$ with $m\le r$\@.  By
$\sections[r]{\man{E}}$ we denote the space of $\C^r$-sections of
$\man{E}$\@.  We define seminorms $p^m_K$\@, $K\subset\man{M}$ compact, on
$\sections[r]{\man{E}}$ by
\begin{equation*}
p^m_K(\xi)=\sup\setdef{\dnorm{j_m\xi(x)}_{\ol{\metric}_m}}{x\in K},
\end{equation*}
and these seminorms define a locally convex topology that we call the
\defn{$\CO^m$-topology}\@.  Let us list some of the attributes of this
topology.
\begin{compactenum}[$\CO^m$-1.]
\item It is Hausdorff:~\cite[4.3.1]{PWM:80}\@.
\item \label{enum:COm-complete} It is complete if and only if
$m=r$\@:~\cite[4.3.2]{PWM:80}\@.
\item It is metrisable:~\cite[4.3.1]{PWM:80}\@.
\item \label{enum:COm-separable} It is separable: This can be shown to follow
by an argument similar to that given above for the $\CO^\infty$-topology.
\item It is probably not nuclear: In case $\man{M}$ is compact, note that
$p^m_{\man{M}}$ is a norm that characterises the $\CO^m$-topology.  A normed
vector space is nuclear if and only if it is
finite-dimensional~\cite[Theorem~4.4.14]{AP:69}\@, so the $\CO^m$-topology
cannot be nuclear when $\man{M}$ is compact except in cases of degenerate
dimension.  But, even when $\man{M}$ is not compact, the $\CO^m$-topology is
not likely nuclear, although we have neither found a reference nor proved
this.
\item \label{enum:COm-suslin} It is Suslin when $m=r$\@: This follows since
$\sections[m]{\tb{\man{M}}}$ is a Polish space, as we have already seen.
\item The $\CO^m$-topology is weaker than the $\CO^r$-topology: This is more
or less clear from the definitions.
\end{compactenum}

From the preceding, we point out two places where one must take care in using
the $\CO^m$-topology, $m\in\integernn$\@, contrasted with the
$\CO^\infty$-topology.  First of all, the topology, if used on
$\sections[r]{\man{E}}$\@, $r>m$\@, is not complete, so convergence arguments
must be modified appropriately.  Second, it is no longer the case that
bounded sets are relatively compact.  Instead, relatively compact subsets
will be described by an appropriate version of the
Arzel\`a\textendash{}Ascoli Theorem,~\cf~\cite[Theorem~5.21]{JJ:05}\@.
Therefore, we need to specify for these spaces whether we will be using the
von Neumann bornology or the compact bornology when we use the word
``bounded.''  These caveats notwithstanding, it is oftentimes appropriate to
use these weaker topologies.

Of course, the preceding can be specialised to vector fields and functions,
and one can define the weak-$\sL$ topologies corresponding to the topologies
for finitely differentiable sections.  In doing this, we apply the general
construction of Definition~\ref{def:weakA} with
$\alg{U}=\sections[r]{\tb{\man{M}}}$\@, $\alg{V}=\func[r]{\man{M}}$ (with the
$\CO^m$-topology), and $\sA=\setdef{\sL_f}{f\in\func[\infty]{\man{M}}}$\@,
where
\begin{equation*}
\mapdef{\sL_f}{\sections[r]{\tb{\man{M}}}}{\func[r]{\man{M}}}
{X}{\lieder{X}{f}.}
\end{equation*}
This gives the following definition.
\begin{definition}
Let $\man{M}$ be a smooth manifold, let $m\in\integernn$\@, and let
$r\in\integernn\cup\{\infty\}$ have the property that $r\ge m$\@.  The
\defn{weak-$(\sL,m)$ topology} for $\sections[r]{\tb{\man{M}}}$ is the
weakest topology for which $\sL_f$ is continuous for each
$f\in\func[\infty]{\man{M}}$\@, where $\func[r]{\man{M}}$ is given the
$\CO^m$-topology.\oprocend
\end{definition}

We can show that the weak-$(\sL,m)$ topology agrees with the
$\CO^m$-topology.
\begin{theorem}\label{the:COm-weak}
Let\/ $\man{M}$ be a smooth manifold, let\/ $m\in\integernn$\@, and let\/
$r\in\integernn\cup\{\infty\}$ have the property that\/ $r\ge m$\@.  Then the
following two topologies for\/ $\sections[r]{\tb{\man{M}}}$ agree:
\begin{compactenum}[(i)]
\item the\/ $\CO^m$-topology;
\item the weak-$(\sL,m)$-topology.
\end{compactenum}
\begin{proof}
Let us first show that the $\CO^m$-topology is weaker than the weak-$(\sL,m)$
topology.  Just as in the corresponding part of the proof of
Theorem~\ref{the:COinfty-weak}\@, we can show that, for $K\subset\man{M}$
compact, there exist $f^1,\dots,f^r\in\func[\infty]{\man{M}}$\@, compact
$K_1,\dots,K_r\subset\man{M}$\@, and $C_1,\dots,C_r\in\realp$ such that
\begin{equation*}
p^m_K(X)\le C_1p^m_{K_1}(\lieder{X}{f^1})+\dots+C_rp^m_{K_r}(\lieder{X}{f^r})
\end{equation*}
for every $X\in\sections[r]{\tb{\man{M}}}$\@.  This estimate gives this part
of the theorem.

To prove that the weak $(\sL,m)$-topology is weaker than the
$\CO^m$-topology, it suffices to show that $\sL_f$ is continuous if
$\sections[r]{\tb{\man{M}}}$ and $\func[r]{\man{M}}$ are given the
$\CO^m$-topology.  This can be done just as in
Theorem~\ref{the:COinfty-weak}\@, with suitable modifications since we only
have to account for $m$ derivatives.
\end{proof}
\end{theorem}

We also have the corresponding relationships between various attributes and
their weak counterparts.
\begin{corollary}\label{cor:COmweak}
Let\/ $\man{M}$ be a smooth manifold, let\/ $m\in\integernn$\@, and let\/
$r\in\integernn\cup\{\infty\}$ have the property that\/ $r\ge m$\@.  Let\/
$(\ts{X},\sO)$ be a topological space, let\/ $(\ts{T},\sM)$ be a measurable
space, and let\/ $\map{\mu}{\sM}{\erealnn}$ be a finite measure.  The
following statements hold:
\begin{compactenum}[(i)]
\item a subset\/ $\nbhd{B}\subset\sections[r]{\tb{\man{M}}}$ is\/
$\CO^m$-bounded in the von Neumann bornology if and only if it is
weak-$(\sL,m)$ bounded in the von Neumann bornology;
\item \label{pl:COmweakcont} a map\/
$\map{\Phi}{\ts{X}}{\sections[r]{\tb{\man{M}}}}$ is\/ $\CO^m$-continuous if
and only if it is weak-$(\sL,m)$ continuous;
\item a map\/ $\map{\Psi}{\ts{T}}{\sections[m]{\tb{\man{M}}}}$ is\/
$\CO^m$-measurable if and only if it is weak-$(\sL,m)$ measurable;
\item a map\/ $\map{\Psi}{\ts{T}}{\sections[m]{\tb{\man{M}}}}$ is Bochner
integrable if and only if it is weak-$(\sL,m)$ Bochner integrable.
\end{compactenum}
\begin{proof}
In the proof of Corollary~\ref{cor:COinftyweak} we established that
$\setdef{\sL_f}{f\in\func[\infty]{\man{M}}}$ was point separating as a family
of linear mappings with domain $\sections[\infty]{\tb{\man{M}}}$\@.  The same
proof is valid if the domain is $\sections[m]{\tb{\man{M}}}$\@.  The result
is then a direct consequence of Lemma~\ref{lem:weakA}\@, taking care to note
that the $\CO^m$-topology on $\sections[r]{\tb{\man{M}}}$ is separable, and
is also complete and Suslin when $r=m$ (and $\func[r]{\man{M}}$ is Suslin
when $r=m$), as we have seen in
properties~$\CO^m$-\ref{enum:COm-complete}\@,%
~$\CO^m$-\ref{enum:COm-separable}\@, and~$\CO^m$-\ref{enum:COm-suslin} above.
\end{proof}
\end{corollary}

\subsection{Topologies for Lipschitz vector fields}\label{subsec:COlip-topology}

It is also possible to characterise Lipschitz sections, so let us indicate
how this is done in geometric terms.  Throughout our discussion of the
Lipschitz case, we make the assumption that the affine connection $\nabla$ on
$\man{M}$ is the Levi-Civita connection for $\metric$ and that the linear
connection $\nabla^0$ on $\man{E}$ is $\metric_0$-orthogonal, by which we
mean that parallel translation consists of isometries.  The existence of such
a connection is ensured by the reasoning of \citet{SK/KN:63a} following the
proof of their Proposition~III.1.5.  We suppose that $\man{M}$ is connected,
for simplicity.  If it is not, then one has to allow the metric we are about
to define to take infinite values.  This is not
problematic~\cite[Exercise~1.1.2]{DB/YB/SI:01}\@, but we wish to avoid the
more complicated accounting procedures.  The \defn{length} of a piecewise
differentiable curve $\map{\gamma}{\interval[a,b]}{\man{M}}$ is
\begin{equation*}
\ell_{\metric}(\gamma)=\int_a^b\sqrt{\metric(\gamma'(t),\gamma'(t))}\,\d{t}.
\end{equation*}
One easily shows that the length of the curve $\gamma$ depends only on
$\image(\gamma)$\@, and not on the particular parameterisation.  We can,
therefore, restrict ourselves to curves defined on $\interval[0,1]$\@.  In
this case, for $x_1,x_2\in\man{M}$\@, we define the \defn{distance} between
$x_1$ and $x_2$ to be
\begin{multline*}
\d_{\metric}(x_1,x_2)=\inf\{\ell_{\metric}(\gamma)|\enspace
\map{\gamma}{\interval[0,1]}{\man{M}}\ \textrm{is a piecewise}\\
\textrm{differentiable curve for which}\ \gamma(0)=x_1\ \textrm{and}\
\gamma(1)=x_2\}.
\end{multline*}
It is relatively easy to show that $(\man{M},\d_{\metric})$ is a metric
space~\cite[Proposition~5.5.10]{RA/JEM/TSR:88}\@.

Now we define a canonical Riemannian metric on the total space $\man{E}$ of a
vector bundle $\map{\pi}{\man{E}}{\man{M}}$\@, following the construction
of~\citet{SS:58} for tangent bundles.  The linear connection $\nabla^0$ gives
a splitting $\tb{\man{E}}\simeq\pi^*\tb{\man{M}}\oplus
\pi^*\man{E}$~\cite[\S11.11]{IK/PWM/JS:93}\@.  The second component of this
decomposition is the vertical component so $\tf[e_x]{\pi}$ restricted to the
first component is an isomorphism onto $\tb[x]{\man{M}}$\@,~\ie~the first
component is ``horizontal.''  Let us denote by
$\map{\hor}{\tb{\man{E}}}{\pi^*\tb{\man{M}}}$ and
$\map{\ver}{\tb{\man{E}}}{\pi^*\man{E}}$ the projections onto the first and
second components of the direct sum decomposition.  This then gives the
Riemannian metric $\metric_{\man{E}}$ on $\man{E}$ defined by
\begin{equation*}
\metric_{\man{E}}(X_{e_x},Y_{e_x})=\metric(\hor(X_{e_x}),\hor(Y_{e_x}))+
\metric_0(\ver(X_{e_x}),\ver(Y_{e_x})).
\end{equation*}

Now let us consider various ways of characterising Lipschitz sections.  To
this end, we let $\map{\xi}{\man{M}}{\man{E}}$ be such that
$\xi(x)\in\man{E}_x$ for every $x\in\man{M}$\@.  For compact
$K\subset\man{M}$ we then define
\begin{equation*}
L_K(\xi)=\sup\asetdef{\frac{\d_{\metric_{\man{E}}}(\xi(x_1),\xi(x_2))}
{\d_{\metric}(x_1,x_2)}}{x_1,x_2\in K,\ x_1\not=x_2}.
\end{equation*}
This is the \defn{$K$-dilatation} of $\xi$\@.  For a piecewise differentiable
curve $\map{\gamma}{\interval[0,T]}{\man{M}}$\@, we denote by
$\map{\tau_{\gamma,t}}{\man{E}_{\gamma(0)}}{\man{E}_{\gamma(t)}}$ the
isomorphism of parallel translation along $\gamma$ for each
$t\in\interval[0,T]$\@.  We then define
\begin{equation}\label{eq:lK}
l_K(\xi)=\sup\asetdef{\frac{\dnorm{\tau^{-1}_{\gamma,1}(\xi\scirc\gamma(1))-
\xi\scirc\gamma(0)}_{\metric_0}}{\ell_\metric(\gamma)}}
{\map{\gamma}{\interval[0,1]}{\man{M}},\ \gamma(0),\gamma(1)\in K,\
\gamma(0)\not=\gamma(1)},
\end{equation}
which is the \defn{$K$-sectional dilatation} of $\xi$\@.  Finally, we define
\begin{equation*}
\mapdef{\Dil{\xi}}{\man{M}}{\realnn}{x}
{\inf\setdef{L_{\closure(\nbhd{U})}(\xi)}{\nbhd{U}\ \textrm{is a
relatively compact neighbourhood of}\ x},}
\end{equation*}
and
\begin{equation*}
\mapdef{\dil{\xi}}{\man{M}}{\realnn}{x}
{\inf\setdef{l_{\closure(\nbhd{U})}(\xi)}{\nbhd{U}\ \textrm{is a
relatively compact neighbourhood of}\ x},}
\end{equation*}
which are the \defn{local dilatation} and \defn{local sectional
dilatation}\@, respectively, of $\xi$\@.
Following~\cite[Proposition~1.5.2]{NW:99} one can show that
\begin{equation*}
L_K(\xi+\eta)\le L_K(\xi)+L_K(\eta),\quad
l_K(\xi+\eta)\le l_K(\xi)+l_K(\eta),\qquad K\subset\man{M}\ \textrm{compact},
\end{equation*}
and
\begin{equation*}
\Dil{(\xi+\eta)}(x)\le\Dil{\xi}(x)+\Dil{\eta}(x),\quad
\dil{(\xi+\eta)}(x)\le\dil{\xi}(x)+\dil{\eta}(x),\qquad x\in\man{M}.
\end{equation*}

The following lemma connects the preceding notions.
\begin{lemma}\label{lem:loclip}
Let\/ $\map{\pi}{\man{E}}{\man{M}}$ be a smooth vector bundle and let\/
$\map{\xi}{\man{M}}{\man{E}}$ be such that\/ $\xi(x)\in\man{E}_x$ for every\/
$x\in\man{M}$\@.  Then the following statements are equivalent:
\begin{compactenum}[(i)]
\item \label{pl:lipchar1} $L_K(\xi)<\infty$ for every compact\/
$K\subset\man{M}$\@;
\item \label{pl:lipchar2} $l_K(\xi)<\infty$ for every compact\/
$K\subset\man{M}$\@;
\item \label{pl:lipchar3} $\Dil{\xi}(x)<\infty$ for every\/ $x\in\man{M}$\@;
\item \label{pl:lipchar4} $\dil{\xi}(x)<\infty$ for every\/ $x\in\man{M}$\@.
\end{compactenum}
Moreover, we have the equalities
\begin{equation*}
L_K(\xi)=\sqrt{l_K(\xi)^2+1},\quad \Dil{\xi}(x)=\sqrt{\dil{\xi}(x)^2+1}
\end{equation*}
for every compact\/ $K\subset\man{M}$ and every\/ $x\in\man{M}$\@.
\begin{proof}
The equivalence of~\eqref{pl:lipchar1} and~\eqref{pl:lipchar2}\@, along with
the equality $L_K=\sqrt{l_K^2+1}$\@, follows from the arguments of
\citet[Lemma~II.A.2.4]{RDC/DBAE/AM:06}\@.  This also implies the equality
$\Dil{\xi}(x)=\sqrt{\dil{\xi}(x)^2+1}$ when both $\Dil{\xi}(x)$ and
$\dil{\xi}(x)$ are finite.

\eqref{pl:lipchar1}$\implies$\eqref{pl:lipchar3} If $x\in\man{M}$ and if
$\nbhd{U}$ is a relatively compact neighbourhood of $x$\@, then
$L_{\closure(\nbhd{U})}(\xi)<\infty$ and so $\Dil{\xi}(x)<\infty$\@.

\eqref{pl:lipchar2}$\implies$\eqref{pl:lipchar4} This follows just as does the
preceding part of the proof.

\eqref{pl:lipchar3}$\implies$\eqref{pl:lipchar1} Suppose that
$\Dil{\xi}(x)<\infty$ for every $x\in\man{M}$ and that there exists a compact
set $K\subset\man{M}$ such that $L_K(\xi)\not<\infty$\@.  Then there exist
sequences $\ifam{x_j}_{j\in\integerp}$ and $\ifam{y_j}_{j\in\integerp}$ in
$K$ such that $x_j\not=y_j$\@, $j\in\integerp$\@, and
\begin{equation*}
\lim_{j\to\infty}\frac{\d_{\metric_{\man{E}}}(\xi(x_j),\xi(y_j))}
{\d_{\metric}(x_j,y_j)}=\infty.
\end{equation*}
Since $\Dil{\xi}(x)<\infty$ for every $x\in\man{M}$\@, it follows directly
that $\xi$ is continuous and so $\xi(K)$ is bounded in the metric
$\metric_{\man{E}}$\@.  Therefore, there exists $C\in\realp$ such that
\begin{equation*}
\d_{\metric_{\man{E}}}(\xi(x_j),\xi(y_j))\le C,\qquad j\in\integerp,
\end{equation*}
and so we must have $\lim_{j\to\infty}\d_{\metric}(x_j,y_j)=0$\@.  Let
$\ifam{x_{j_k}}_{k\in\integerp}$ be a subsequence converging to $x\in K$ and
note that $\ifam{y_{j_k}}_{k\in\integerp}$ then also converges to $x$\@.  This
implies that $\Dil{\xi}(x)\not<\infty$\@, which proves the result.

\eqref{pl:lipchar4}$\implies$\eqref{pl:lipchar2} This follows just as the
preceding part of the proof.
\end{proof}
\end{lemma}

With the preceding, we can define what we mean by a locally Lipschitz section
of a vector bundle, noting that, if $\dil{\xi}(x)<\infty$ for every
$x\in\man{M}$\@, $\xi$ is continuous.  Our definition is in the general
situation where sections are of class $\C^m$ with the $m$th derivative being,
not just continuous, but Lipschitz.
\begin{definition}
For a smooth vector bundle $\map{\pi}{\man{E}}{\man{M}}$ and for
$m\in\integernn$\@, $\xi\in\sections[m]{\man{E}}$ is of \defn{class
$\C^{m+\lip}$} if $\map{j_m\xi}{\man{M}}{\jet{m}{\man{E}}}$ satisfies any of
the four equivalent conditions of Lemma~\ref{lem:loclip}\@.  If $\xi$ is of
class $\C^{0+\lip}$ then we say it is \defn{locally Lipschitz}\@.  By
$\sections[\lip]{\man{E}}$ we denote the space of locally Lipschitz sections
of $\man{E}$\@.  For $m\in\integernn$\@, by $\sections[m+\lip]{\man{E}}$ we
denote the space of sections of $\man{E}$ of class $\C^{m+\lip}$\@.\oprocend
\end{definition}

It is straightforward, if tedious, to show that a section is of class
$\C^{m+\lip}$ if and only if, in any coordinate chart, the section is
$m$-times continuously differentiable with the $m$th derivative being locally
Lipschitz in the usual Euclidean sense.  The essence of the argument is that,
in any sufficiently small neighbourhood of a point in $\man{M}$\@, the
distance functions $\d_{\metric}$ and $\d_{\metric_{\man{E}}}$ are equivalent
to the Euclidean distance functions defined in coordinates.

The following characterisation of the local sectional dilatation is useful.
\begin{lemma}\label{lem:dil-deriv}
For a smooth vector bundle\/ $\map{\pi}{\man{E}}{\man{M}}$ and for\/
$\xi\in\sections[\lip]{\man{E}}$\@, we have
\begin{multline*}
\dil{\xi}(x)=\inf\{\sup\setdef{\dnorm{\nabla_{v_y}\xi}_{\metric_0}}
{y\in\closure(\nbhd{U}),\ \dnorm{v_y}_{\metric}=1,\ \xi\
\textrm{differentiable at}\ y}|\\
\nbhd{U}\ \textrm{is a relatively compact neighbourhood of}\ x\}.
\end{multline*}
\begin{proof}
As per~\cite[Proposition~IV.3.4]{SK/KN:63a}\@, let $\nbhd{U}$ be a
geodesically convex, relatively compact open set.  We claim that
\begin{equation*}
l_{\closure(\nbhd{U})}(\xi)=\sup\setdef{\dnorm{\nabla^0_{v_y}\xi}_{\metric_0}}
{y\in\closure(\nbhd{U}),\ \dnorm{v_y}_{\metric}=1,\ \xi\
\textrm{differentiable at}\ y}.
\end{equation*}
By~\cite[Lemma~II.A.2.4]{RDC/DBAE/AM:06}\@, to determine
$l_{\closure(\nbhd{U})}(\xi)$\@, it suffices in the formula~\eqref{eq:lK} to
use only length minimising geodesics whose images are contained in
$\closure(\nbhd{U})$\@.  Let $x\in\nbhd{U}$\@, let $v_x\in\tb[x]{\man{M}}$
have unit length, and let $\map{\gamma}{\interval[0,T]}{\closure(\nbhd{U})}$
be a minimal length geodesic such that $\gamma'(0)=v_x$\@.  If $x$ is a point
of differentiability for $\xi$\@, then
\begin{equation*}
\lim_{t\to0}\frac{\dnorm{\tau_{\gamma,t}^{-1}(\xi\scirc\gamma(t))-
\xi\scirc\gamma(0)}_{\metric_0}}{t}=\dnorm{\nabla^0_{v_y}\xi}_{\metric_0}.
\end{equation*}
From this we conclude that
\begin{equation*}
l_{\closure(\nbhd{U})}(\xi)\ge\sup\setdef{\dnorm{\nabla^0_{v_x}\xi}_{\metric_0}}
{x\in\closure(\nbhd{U}),\ \dnorm{v_x}_{\metric}=1,\ \xi\
\textrm{differentiable at}\ y}.
\end{equation*}
Suppose the opposite inequality does not hold.  Then there exist
$x_1,x_2\in\closure(\nbhd{U})$ such that, if
$\map{\gamma}{\interval[0,T]}{\man{M}}$ is the arc-length parameterised
minimal length geodesic from $x_1$ to $x_2$\@, then
\begin{equation}\label{eq:lKderiv}
\frac{\dnorm{\tau_{\gamma,T}^{-1}(\xi\scirc\gamma(T))-\xi\scirc\gamma(0)}}{T}
>\dnorm{\nabla^0_{v_x}\xi}_{\metric_0}
\end{equation}
for every $x\in\closure(\nbhd{U})$ for which $\xi$ is differentiable at $x$
and every $v_x\in\tb[x]{\man{M}}$ of unit length.  Note that $\alpha\colon
t\mapsto\tau_{\gamma,t}^{-1}(\xi\scirc\gamma(t))$ is a Lipschitz curve in
$\tb[x_1]{\man{M}}$\@.  By Rademacher's
Theorem~\cite[Theorem~3.1.5]{HF:96}\@, this curve is almost everywhere
differentiable.  If $\alpha$ is differentiable at $t$ we have
\begin{equation*}
\alpha'(t)=\tau_{\gamma,t}^{-1}(\nabla^0_{\gamma'(t)}\xi).
\end{equation*}
Therefore, also by Rademacher's Theorem and since $\nabla^0$ is $\metric_0$-orthogonal, we have
\begin{multline*}
\sup\asetdef{\frac{\dnorm{\tau_{\gamma,t}^{-1}(\xi\scirc\gamma(t))-
\xi\scirc\gamma(0)}_{\metric_0}}{t}}{t\in\interval[0,T]}\\
=\sup\setdef{\dnorm{\nabla^0_{\gamma'(t)}\xi}_{\metric_0}}
{t\in\interval[0,T],\ \xi\ \textrm{is differentiable at}\ \gamma(t)}.
\end{multline*}
This, however, contradicts~\eqref{eq:lKderiv}\@, and so our claim holds.

Now let $x\in\man{M}$ and let $\ifam{\nbhd{U}_j}_{j\in\integerp}$ be a
sequence of relatively compact, geodesically convex neighbourhood of $x$ such
that $\cap_{j\in\integerp}\nbhd{U}_j=\{x\}$\@.  Then
\begin{equation*}
\dil{\xi}(x)=\lim_{j\to\infty}l_{\closure(\nbhd{U}_j)}(\xi)
\end{equation*}
and
\begin{multline*}
\inf\{\sup\setdef{\dnorm{\nabla^0_{v_y}\xi}_{\metric_0}}
{y\in\closure(\nbhd{U}),\ \dnorm{v_y}_{\metric}=1,\ \xi\
\textrm{differentiable at}\ y}|\\
\nbhd{U}\ \textrm{is a relatively compact neighbourhood of}\ x\}\\
=\lim_{j\to\infty}\sup\setdef{\dnorm{\nabla^0_{v_y}\xi}_{\metric_0}}
{y\in\closure(\nbhd{U}_j),\ \dnorm{v_y}_{\metric}=1,\ \xi\
\textrm{differentiable at}\ y}.
\end{multline*}
The lemma now follows from the claim in the opening paragraph.
\end{proof}
\end{lemma}

Let us see how to topologise spaces of locally Lipschitz sections.
Lemma~\ref{lem:loclip} gives us four possibilities for doing this.  In order
to be as consistent as possible with our other definitions of seminorms, we
use the ``locally sectional'' characterisation of Lipschitz seminorms.  Thus,
for $\xi\in\sections[\lip]{\man{E}}$ and $K\subset\man{M}$ compact, let us
define
\begin{equation*}
\lambda_K(\xi)=\sup\setdef{\dil{\xi}(x)}{x\in K}
\end{equation*}
and then define a seminorm $p^{\lip}_K$\@, $K\subset\man{M}$ compact, on
$\sections[\lip]{\man{E}}$ by
\begin{equation*}
p^{\lip}_K(\xi)=\max\{\lambda_K(\xi),p_K^0(\xi)\}.
\end{equation*}
The seminorms $p_K^{\lip}$\@, $K\subset\man{M}$ compact, give the
\defn{$\CO^{\lip}$-topology} on $\sections[r]{\man{E}}$ for
$r\in\integerp\cup\{\infty\}$\@.  To topologise
$\sections[m+\lip]{\man{E}}$\@, note that the $\CO^{\lip}$-topology on
$\sections[\lip]{\jet{m}{\man{E}}}$ induces a topology on
$\sections[m+\lip]{\man{E}}$ that we call the
\defn{$\CO^{m+\lip}$-topology}\@.  The seminorms for this locally convex
topology are
\begin{equation*}
p^{m+\lip}_K(\xi)=\max\{\lambda_K^m(\xi),p^m_K(\xi)\},\qquad
K\subset\man{M}\ \textrm{compact},
\end{equation*}
where
\begin{equation*}
\lambda_K^m(\xi)=\sup\setdef{\dil{j_m\xi(x)}}{x\in K}.
\end{equation*}
Note that $\dil{j_m\xi}$ is unambiguously defined.  Let us briefly explain
why.  If the connections $\nabla$ and $\nabla^0$ are metric connections for
$\metric$ and $\metric_0$\@, as we are assuming, then the induced connection
$\nabla^m$ on $\tensor[k]{\ctb{\man{M}}}\otimes\man{E}$ is also metric with
respect to the induced metric determined from Lemma~\ref{lem:inprodotimes}\@.
It then follows from Lemma~\ref{lem:Jrdecomp} that the dilatation for
sections of $\jet{m}{\man{E}}$ can be defined just as for sections of
$\man{E}$\@.

Note that $\sections[\lip]{\man{E}}\subset\sections[0]{\man{E}}$ and
$\sections[r]{\man{E}}\subset\sections[\lip]{\man{E}}$ for $r\in\integerp$\@.
Thus we adopt the convention that $0<\lip<1$ for the purposes of ordering
degrees of regularity.  Let $m\in\integernn$\@, and let
$r\in\integernn\cup\{\infty\}$ and $r'\in\{0,\lip\}$ be such that $r+r'\ge
m+\lip$\@.  We adopt the obvious convention that $\infty+\lip=\infty$\@.  The
seminorms $p^{m+\lip}_K$\@, $K\subset\man{M}$ compact, can then be defined on
$\sections[r+r']{\man{E}}$\@.

Let us record some properties of the $\CO^{m+\lip}$-topology for
$\sections[r+r']{\man{E}}$\@.  This topology is not extensively studied like
the other differentiable topologies, but we can nonetheless enumerate its
essential properties.
\begin{compactenum}[$\CO^{m+\lip}$-1.]
\item It is Hausdorff: This is clear.
\item \label{enum:COmm'-complete} It is complete if and only if
$r+r'=m+\lip$\@: This is more or less because, for a compact metric space,
the space of Lipschitz functions is a Banach
space~\cite[Proposition~1.5.2]{NW:99}\@.  Since $\sections[m+\lip]{\man{E}}$
is the inverse limit of the Banach spaces
$\sections[m+\lip]{\man{E}|K_j}$\@,\footnote{To be clear, by
$\sections[m+\lip]{\man{E}|K}$ we denote the space of sections of class
$m+\lip$ defined on a neighbourhood of $K$\@.} $j\in\integerp$\@, for a
compact exhaustion $\ifam{K_j}_{j\in\integerp}$ of $\man{M}$\@, and since the
inverse limit of complete locally convex spaces is
complete~\cite[Proposition~2.11.3]{JH:66}\@, we conclude the stated
assertion.
\item It is metrisable: This is argued as follows.  First of all, it is a
countable inverse limit of Banach spaces.  Inverse limits are closed
subspaces of the direct product~\cite[Proposition~V.19]{APR/WR:80}\@.  The
direct product of metrisable spaces, in particular Banach spaces, is
metrisable~\cite[Theorem~22.3]{SW:04}\@.
\item \label{enum:COmm'-separable} It is separable: This is a consequence of
the result of \citet[Theorem~1.2$'$]{REG/HW:79} which says that Lipschitz
functions on Riemannian manifolds can be approximated in the
$\CO^{\lip}$-topology by smooth functions, and by the separability of the
space of smooth functions.
\item It is probably not nuclear: For compact base manifolds,
$\sections[m+\lip]{\man{E}}$ is an infinite-dimensional normed space,
and so not nuclear~\cite[Theorem~4.4.14]{AP:69}\@.  But, even when $\man{M}$
is not compact, the $\CO^{m+\lip}$-topology is not likely nuclear,
although we have neither found a reference nor proved this.
\item \label{enum:COmm'-suslin} It is Suslin when $m+\lip=r+r'$\@: This
follows since $\sections[m+\lip]{\man{E}}$ is a Polish space, as we have
already seen.
\end{compactenum}

Of course, the preceding can be specialised to vector fields and functions,
and one can define the weak-$\sL$ topologies corresponding to the above
topologies.  To do this, we apply the general construction of
Definition~\ref{def:weakA} with $\alg{U}=\sections[r+r']{\tb{\man{M}}}$\@,
$\alg{V}=\func[r+r']{\man{M}}$ (with the $\CO^m$-topology), and
$\sA=\setdef{\sL_f}{f\in\func[\infty]{\man{M}}}$\@, where
\begin{equation*}
\mapdef{\sL_f}{\sections[r+r']{\tb{\man{M}}}}{\func[r+r']{\man{M}}}
{X}{\lieder{X}{f}.}
\end{equation*}
We then have the following definition.
\begin{definition}
Let $\man{M}$ be a smooth manifold, let $m\in\integernn$\@, and let
$r\in\integernn\cup\{\infty\}$ and $r'\in\{0,\lip\}$ have the property that
$r+r'\ge m+\lip$\@.  The \defn{weak-$(\sL,m+\lip)$ topology} for
$\sections[r+r']{\tb{\man{M}}}$ is the weakest topology for which $\sL_f$ is
continuous for each $f\in\func[\infty]{\man{M}}$\@, where
$\func[r+r']{\man{M}}$ is given the $\CO^{m+\lip}$-topology.\oprocend
\end{definition}

We can show that the weak-$(\sL,m+\lip)$ topology agrees with the
$\CO^{m+\lip}$-topology.
\begin{theorem}\label{the:COmm'-weak}
Let\/ $\man{M}$ be a smooth manifold, let\/ $m\in\integernn$\@, and let\/
$r\in\integernn\cup\{\infty\}$ and\/ $r'\in\{0,\lip\}$ have the property
that\/ $r+r'\ge m+\lip$\@.  Then the following two topologies for\/
$\sections[r+r']{\man{E}}$ agree:
\begin{compactenum}[(i)]
\item the\/ $\CO^{m+\lip}$-topology;
\item the weak-$(\sL,m+\lip)$-topology.
\end{compactenum}
\begin{proof}
We prove the theorem only for the case $m=0$\@, since the general case
follows from this in combination with Theorem~\ref{the:COm-weak}\@.

Let us first show that the $\CO^{\lip}$-topology is weaker than the
weak-$(\sL,\lip)$ topology.  Let $K\subset\man{M}$ be compact and for
$x\in\man{M}$ choose a coordinate chart $(\nbhd{U}_x,\phi_x)$ and functions
$f^1_x,\dots,f^n_x\in\func[\infty]{\man{M}}$ agreeing with the coordinate
functions in a neighbourhood of a geodesically convex relatively compact
neighbourhood $\nbhd{V}_x$ of $x$~\cite[Proposition~IV.3.4]{SK/KN:63a}\@.  We
denote by $\map{\vect{X}}{\phi_x(\nbhd{U}_x)}{\real^n}$ the local
representative of $X$\@.  Since $\lieder{X}{f^j_x}=X^j$ on a neighbourhood of
$\nbhd{V}_x$\@, there exists $C_x\in\realp$ such that
\begin{equation*}
\dnorm{\tau_{\gamma,1}^{-1}(X(x_1))-X(x_2)}_{\metric}\le C_x\sum_{j=1}^n
\snorm{\lieder{X}{f^j_x}(x_1)-\lieder{X}{f^j_x}(x_2)}
\end{equation*}
for every distinct $x_1,x_2\in\closure(\nbhd{V}_x)$\@, where $\gamma$ is the
unique minimal length geodesic from $x_2$ to $x_1$ (the inequality is a
consequence of the fact that the $\ell^1$ norm for $\real^n$ is equivalent to
any other norm).  This gives an inequality
\begin{equation*}
\dil{X}(y)\le C_x(\dil{\lieder{X}{f^1_x}}(y)+\dots+\dil{\lieder{X}{f^n_x}}(y))
\end{equation*}
for every $y\in\nbhd{V}_x$\@.  Now let $x_1,\dots,x_k\in K$ be such that
$K\subset\cup_{j=1}^k\nbhd{V}_{x_j}$\@.  From this point, it is a bookkeeping
exercise, exactly like that in the corresponding part of the proof of
Theorem~\ref{the:COinfty-weak}\@, to arrive at the inequality
\begin{equation*}
\lambda_K(X)\le C_1\lambda_K(\lieder{X}{f^1})+\dots+C_r\lambda_K(\lieder{X}{f^r}).
\end{equation*}
From the proof of Theorem~\ref{the:COm-weak} we also have
\begin{equation*}
p_K^0(X)\le C'_1p_K^0(\lieder{X}{f^1})+\dots+C'_rp_K^0(\lieder{X}{f^r}),
\end{equation*}
and this gives the result.

To prove that the weak $(\sL,\lip)$-topology is weaker than the
$\CO^{\lip}$-topology, it suffices to show that $\sL_f$ is continuous for
every $f\in\func[\infty]{\man{M}}$ if $\sections[r+r']{\tb{\man{M}}}$ and
$\func[r+r']{\man{M}}$ are given the $\CO^{m+\lip}$-topology.  Thus let
$K\subset\man{M}$ be compact and let $f\in\func[\infty]{\man{M}}$\@.  We
choose a relatively compact geodesically convex chart $(\nbhd{U}_x,\phi_x)$
about $x\in K$ and compute, for distinct $x_1,x_2\in\nbhd{U}_x$\@,
\begin{align*}
|\lieder{X}{f}&(x_1)-\lieder{X}{f}(x_2)|\\
\le&\;\sum_{j=1}^n\Bigsnorm{X^j(x_1)\pderiv{f}{x^j}(x_1)-
X^j(x_2)\pderiv{f}{x^j}(x_2)}\\
\le&\;\sum_{j=1}^n\Bigl(\snorm{X^j(x_1)}
\Bigsnorm{\pderiv{f}{x^j}(x_1)-\pderiv{f}{x^j}(x_2)}+
\snorm{X^j(x_1)-X^j(x_2)}\Bigsnorm{\pderiv{f}{x^j}(x_2)}\Bigr)\\
\le&\;\sum_{j=1}^n\Bigl(A_xp_{\closure(\nbhd{U}_x)}^0(X)\pderiv{f}{x^j}(y)
\d_{\metric}(x_1,x_2)\Bigr)+B_x\dnorm{\tau_{\gamma,1}^{-1}X(x_1)-X(x_2)}_{\metric},
\end{align*}
for some $y\in\nbhd{U}_x$\@, using the mean value theorem~\cite[Proposition~2.4.8]{RA/JEM/TSR:88}\@, and where $\gamma$
is the unique length minimising geodesic from $x_2$ to $x_1$\@.  Thus we have
an inequality
\begin{equation*}
\lambda_{\closure(\nbhd{U}_x)}(\lieder{X}{f})\le
A_xp_{\closure(\nbhd{U}_x)}^0(X)+B_x\lambda_{\closure(\nbhd{U}_x)}(X),
\end{equation*}
for a possibly different $A_x$\@.  Letting $x_1,\dots,x_k\in K$ be such that
$K\subset\cup_{j=1}^k\nbhd{U}_x$\@, some more bookkeeping like that in the
first part of the proof of Theorem~\ref{the:COinfty-weak} gives
\begin{equation*}
\lambda_K(\lieder{X}{f})\le\sum_{j=1}^r(A_jp_{\closure(\nbhd{U}_{x_j})}^0(X)+
B_j\lambda_{\closure(\nbhd{U}_{x_j})}(X))
\end{equation*}
for suitable constants $A_j,B_j\in\realp$\@, $j\in\{1,\dots,r\}$\@.  Since,
from the proof of Theorem~\ref{the:COm-weak}\@, we also have
\begin{equation*}
p_K^0(\lieder{X}{f})\le\sum_{j=1}^rC_jp_K^0(X)
\end{equation*}
for suitable constants $C_1,\dots,C_r\in\realp$\@, the result follows.
\end{proof}
\end{theorem}

We also have the corresponding relationships between various attributes and
their weak counterparts.
\begin{corollary}\label{cor:COmm'weak}
Let\/ $\man{M}$ be a smooth manifold, let\/ $m\in\integernn$\@, and let\/
$r\in\integernn\cup\{\infty,\lip\}$ and $r'\in\{0,\lip\}$ have the property
that\/ $r+r'\ge m+\lip$\@.  Let\/ $(\ts{X},\sO)$ be a topological space,
let\/ $(\ts{T},\sM)$ be a measurable space, and let\/
$\map{\mu}{\sM}{\erealnn}$ be a finite measure.  The following statements
hold:
\begin{compactenum}[(i)]
\item a subset\/ $\nbhd{B}\subset\sections[r+r']{\tb{\man{M}}}$ is\/
$\CO^{m+\lip}$-bounded in the von Neumann bornology if and only if it is
weak-$(\sL,m+\lip)$ bounded in the von Neumann bornology;
\item \label{pl:COmm'weakcont} a map\/
$\map{\Phi}{\ts{X}}{\sections[r+r']{\tb{\man{M}}}}$ is\/
$\CO^{m+\lip}$-continuous if and only if it is weak-$(\sL,m+\lip)$
continuous;
\item a map\/ $\map{\Psi}{\ts{T}}{\sections[m+\lip]{\tb{\man{M}}}}$ is\/
$\CO^{m+\lip}$-measurable if and only if it is weak-$(\sL,m+\lip)$ measurable;
\item a map\/ $\map{\Psi}{\ts{T}}{\sections[m+\lip]{\tb{\man{M}}}}$ is
Bochner integrable if and only if it is weak-$(\sL,m+\lip)$ Bochner
integrable.
\end{compactenum}
\begin{proof}
In the proof of Corollary~\ref{cor:COinftyweak} we established that
$\setdef{\sL_f}{f\in\func[\infty]{\man{M}}}$ was point separating as a family
of linear mappings with domain $\sections[\infty]{\tb{\man{M}}}$\@.  The same
proof is valid if the domain is $\sections[m+\lip]{\tb{\man{M}}}$\@.  The
result is then a direct consequence of Lemma~\ref{lem:weakA}\@, noting that
the $\CO^{m+\lip}$-topology on $\sections[r+r']{\tb{\man{M}}}$ is separable,
and is also complete and Suslin when $r+r'=m+\lip$ (and
$\func[r+r']{\man{M}}$ is Suslin when $r+r'=m+\lip$), as we have seen above
in properties~$\CO^{m+\lip}$-\ref{enum:COmm'-complete}\@,%
~$\CO^{m+\lip}$-\ref{enum:COmm'-separable}\@,
and~$\CO^{m+\lip}$-\ref{enum:COmm'-suslin}\@.
\end{proof}
\end{corollary}

\begin{notation}
In order to try to compactify the presentation of the various degrees of
regularity we consider, we will frequently speak of the class ``$m+m'$''
where $m\in\integernn$ and $m'\in\{0,\lip\}$\@.  This allows us to include
the various Lipschitz cases alongside the finitely differentiable cases.
Thus, whenever the reader sees ``$m+m'$\@,'' this is what they should have in
mind.\oprocend
\end{notation}

\section{The $\CO^{\hol}$-topology for the space of holomorphic vector
fields}\label{sec:holomorphic-topology}

While in this paper we have no \emph{per se} interest in holomorphic vector
fields, it is the case that an understanding of certain constructions for
real analytic vector fields rely in an essential way on their holomorphic
extensions.  Also, as we shall see, we will arrive at a description of the
real analytic topology that, while often easy to use in general arguments, is
not well suited for verifying hypotheses in examples.  In these cases, it is
often most convenient to extend from real analytic to holomorphic, where
things are easier to verify.

Thus in this section we overview the holomorphic case.  We begin with vector
bundles, as in the smooth case.

\subsection{General holomorphic vector bundles}\label{subsec:COhol-vb}

We let $\map{\pi}{\man{E}}{\man{M}}$ be an holomorphic vector bundle with
$\sections[\hol]{\man{E}}$ the set of holomorphic sections.  We let $\metric$
be an Hermitian fibre metric on $\man{E}$\@, and, for $K\subset\man{M}$
compact, define a seminorm $p^{\hol}_K$ on $\sections[\hol]{\man{E}}$ by
\begin{equation*}
p^{\hol}_K(\xi)=\sup\setdef{\dnorm{\xi(z)}_{\metric}}{z\in K}.
\end{equation*}
The \defn{$\CO^{\hol}$-topology} for $\sections[\hol]{\man{E}}$ is the
locally convex topology defined by the family of seminorms $p^{\hol}_K$\@,
$K\subset\man{M}$ compact.

We shall have occasion to make use of bounded holomorphic sections.  Thus we
let $\map{\pi}{\man{E}}{\man{M}}$ be an holomorphic vector bundle with
Hermitian fibre metric $\metric$\@.  We denote by $\secbdd[\hol]{\man{E}}$
the sections of $\man{E}$ that are bounded, and on $\secbdd[\hol]{\man{E}}$
we define a norm
\begin{equation*}
p^{\hol}_\infty(\xi)=\sup\setdef{\dnorm{\xi(z)}_{\metric}}{z\in\man{M}}.
\end{equation*}
If we wish to draw attention to the domain of the section, we will write the
norm as $p^{\hol}_{\man{M},\infty}$\@.  This will occur when we have sections
defined on an open subset of the manifold.

The following lemma makes an assertion of which we shall make use.
\begin{lemma}\label{lem:secbddhol}
Let\/ $\map{\pi}{\man{E}}{\man{M}}$ be an holomorphic vector bundle.  The
subspace topology on\/ $\secbdd[\hol]{\man{E}}$\@, induced from the\/
$\CO^{\hol}$-topology, is weaker than the norm topology induced by the norm\/
$p^{\hol}_\infty$\@.  Moreover,\/ $\secbdd[\hol]{\man{E}}$ is a Banach space.
Also, if\/ $\nbhd{U}\subset\man{M}$ is a relatively compact open set with\/
$\closure(\nbhd{U})\subsetneq\man{M}$\@, then the restriction map from\/
$\sections[\hol]{\man{E}}$ to\/ $\secbdd[\hol]{\man{E}|\nbhd{U}}$ is
continuous.
\begin{proof}
It suffices to show that a sequence $\ifam{\xi_j}_{j\in\integerp}$ in
$\secbdd[\hol]{\man{E}}$ converges to $\xi\in\secbdd[\hol]{\man{E}}$
uniformly on compact subsets of $\man{M}$ if it converges in norm.  This,
however, is obvious.  It remains to prove completeness of
$\secbdd[\hol]{\man{E}}$ in the norm topology.  By
\cite[Theorem~7.9]{EH/KS:75}\@, a Cauchy sequence
$\ifam{\xi_j}_{j\in\integerp}$ in $\secbdd[\hol]{\man{E}}$ converges to a
bounded continuous section $\xi$ of $\man{E}$\@.  That $\xi$ is also
holomorphic follows since uniform limits of holomorphic sections are
holomorphic~\cite[page~5]{RCG:90a}\@.  For the final assertion, since the
topology of $\sections[\hol]{\man{E}}$ is metrisable (see
$\CO^\hol$-\ref{enum:COhol-metrisable} below), it suffices to show that the
restriction of a convergent sequence in $\sections[\hol]{\man{E}}$ to
$\nbhd{U}$ converges uniformly.  This, however, follows since
$\closure(\nbhd{U})$ is compact.
\end{proof}
\end{lemma}

One of the useful attributes of holomorphic geometry is that properties of
higher derivatives can be deduced from the mapping itself.  To make this
precise, we first make the following observations.
\begin{compactenum}
\item Hermitian inner products on $\complex$-vector spaces give inner
products on the underlying $\real$-vector space.
\item By Lemma~\ref{lem:analytic-conn}\@, there exist a real analytic affine
connection $\nabla$ on $\man{M}$ and a real analytic vector bundle connection
$\nabla^0$ on $\man{E}$\@.
\end{compactenum}
Therefore, the seminorms defined in Section~\ref{subsec:COinfty-vb} can be
made sense of for holomorphic sections.
\begin{proposition}\label{prop:Cauchyest}
Let\/ $\map{\pi}{\man{E}}{\man{M}}$ be an holomorphic vector bundle, let\/
$K\subset\man{M}$ be compact, and let\/ $\nbhd{U}$ be a relatively compact
neighbourhood of\/ $K$\@.  Then there exist\/ $C,r\in\realp$ such that
\begin{equation*}
p^{\infty}_{K,m}(\xi)\le Cr^{-m}p^{\hol}_{\nbhd{U},\infty}(\xi)
\end{equation*}
for every\/ $m\in\integernn$ and\/ $\xi\in\secbdd[\hol]{\man{E}|\nbhd{U}}$\@.

Moreover, if\/ $\ifam{\nbhd{U}_j}_{j\in\integerp}$ is a sequence of
relatively compact neighbourhoods of\/ $K$ such
that~(i)~$\closure(\nbhd{U}_j)\subset\nbhd{U}_{j+1}$
and~(ii)~$K=\cap_{j\in\integerp}\nbhd{U}_j$\@, and if\/ $C_j,r_j\in\realp$
are such that
\begin{equation*}
p^{\infty}_{K,m}(\xi)\le
C_jr_j^{-m}p^{\hol}_{\nbhd{U}_j,\infty}(\xi),
\qquad m\in\integernn,\ \xi\in\secbdd[\hol]{\man{E}|\nbhd{U}_j},
\end{equation*}
then\/ $\lim_{j\to\infty}r_j=0$\@.
\begin{proof}
Let $z\in K$ and let $(\nbhd{W}_z,\psi_z)$ be an holomorphic vector bundle
chart about $z$ with $(\nbhd{U}_z,\phi_z)$ the associated chart for
$\man{M}$\@, supposing that $\nbhd{U}_z\subset\nbhd{U}$\@.  Let
$k\in\integerp$ be such that
$\psi_z(\nbhd{W}_z)=\phi_z(\nbhd{U}_z)\times\complex^k$\@.  Let
$\vect{z}=\phi_z(z)$ and let
$\map{\vect{\xi}}{\phi_z(\nbhd{U}_z)}{\complex^k}$ be the local
representative of $\xi\in\secbdd[\hol]{\man{E}|\nbhd{U}}$\@.  Note that when
taking real derivatives of $\vect{\xi}$ with respect to coordinates, we can
think of taking derivatives with respect to
\begin{equation*}
\pderiv{}{z^j}=\frac{1}{2}\Bigl(\pderiv{}{x^j}-\imag\pderiv{}{y^j}\Bigr),
\quad\pderiv{}{\bar{z}\null^j}=\frac{1}{2}\Bigl(\pderiv{}{x^j}+
\imag\pderiv{}{y^j}\Bigr),\qquad j\in\{1,\dots,n\}.
\end{equation*}
Since $\vect{\xi}$ is holomorphic, the $\pderiv{}{\bar{z}\null^j}$
derivatives will vanish~\cite[page~27]{SGK:92}\@.  Thus, for the purposes of
the multi-index calculations, we consider multi-indices of length $n$ (not
$2n$).  In any case, applying the usual Cauchy
estimates~\cite[Lemma~2.3.9]{SGK:92}\@, there exists $r\in\realp$ such that
\begin{equation*}
\snorm{\linder[I]{\xi^a}(\vect{z})}\le I!r^{-\snorm{I}}
\sup\setdef{\snorm{\xi^a(\vect{\zeta})}}
{\vect{\zeta}\in\cdisk[]{\vect{r}}{\vect{z}}}
\end{equation*}
for every $a\in\{1,\dots,k\}$\@, $I\in\integernn^n$\@, and
$\xi\in\secbdd[\hol]{\man{E}|\nbhd{U}}$\@.  We may choose
$r\in\interval(0,1)$ such that $\cdisk[]{\vect{r}}{\vect{z}}$ is contained in
$\phi_z(\nbhd{U}_z)$\@, where $\vect{r}=(r,\dots,r)$\@.  Denote
$\nbhd{V}_z=\phi_z^{-1}(\odisk[]{\vect{r}}{\vect{z}})$\@.  There exists a
neighbourhood $\nbhd{V}'_z$ of $z$ such that
$\closure(\nbhd{V}'_z)\subset\nbhd{V}_z$ and such that
\begin{equation*}
\snorm{\linder[I]{\xi^a}(\vect{z}')}\le 2I!r^{-\snorm{I}}
\sup\setdef{\snorm{\xi^a(\vect{\zeta})}}
{\vect{\zeta}\in\cdisk[]{\vect{r}}{\vect{z}}}
\end{equation*}
for every $\vect{z}'\in\phi_z(\nbhd{V}'_z)$\@,
$\xi\in\secbdd[\hol]{\man{E}|\nbhd{U}}$\@, $a\in\{1,\dots,k\}$\@, and
$I\in\integernn^n$\@.  If $\snorm{I}\le m$ then, since we are assuming that
$r<1$\@, we have
\begin{equation*}
\frac{1}{I!}\snorm{\linder[I]{\xi^a}(\vect{z}')}\le 2r^{-m}
\sup\setdef{\snorm{\xi^a(\vect{\zeta})}}
{\vect{\zeta}\in\cdisk[]{\vect{r}}{\vect{z}}}
\end{equation*}
for every $a\in\{1,\dots,k\}$\@, $\vect{z}'\in\phi_z(\nbhd{V}'_z)$\@, and
$\xi\in\secbdd[\hol]{\man{E}|\nbhd{U}}$\@.  By
Lemma~\ref{lem:pissy-estimate}\@, it follows that there exist
$C_z,r_z\in\realp$ such that
\begin{equation*}
\dnorm{j_m\xi(z)}_{\ol{\metric}_m}
\le C_zr_z^{-m}p^{\hol}_{\nbhd{V}_z,\infty}(\xi)
\end{equation*}
for all $z\in\nbhd{V}'_z$\@, $m\in\integernn$\@, and
$\xi\in\secbdd[\hol]{\man{E}|\nbhd{U}}$\@.  Let $z_1,\dots,z_k\in K$ be such
that $K\subset\cup_{j=1}^k\nbhd{V}'_{z_j}$\@, and let
$C=\max\{C_{z_1},\dots,C_{z_k}\}$ and $r=\min\{r_{z_1},\dots,r_{z_k}\}$\@.
If $z\in K$\@, then $z\in\nbhd{V}'_{z_j}$ for some $j\in\{1,\dots,k\}$ and so
we have
\begin{equation*}
\dnorm{j_m\xi(z)}_{\ol{\metric}_m}
\le C_{z_j}r_{z_j}^{-m}p^{\hol}_{\nbhd{V}_{z_j},\infty}(\xi)\le
Cr^{-m}p^{\hol}_{\nbhd{U},\infty}(\xi),
\end{equation*}
and taking supremums over $z\in K$ on the left gives the result.

The final assertion of the proposition immediately follows by observing in
the preceding construction how ``$r$'' was defined, namely that it had to be
chosen so that polydisks of radius $r$ in the coordinate charts remained in
$\nbhd{U}$\@.
\end{proof}
\end{proposition}

\subsection{Properties of the $\CO^{\hol}$-topology}\label{subsec:COhol-props}

The $\CO^{\hol}$-topology for $\sections[\hol]{\man{E}}$ has the following
attributes.
\begin{compactenum}[$\CO^{\hol}$-1.]
\item It is Hausdorff:~\cite[Theorem~8.2]{AK/PWM:97}\@.
\item \label{enum:COhol-complete} It is
complete:~\cite[Theorem~8.2]{AK/PWM:97}\@.
\item \label{enum:COhol-metrisable} It is metrisable:~\cite[Theorem~8.2]{AK/PWM:97}\@.
\item \label{enum:COhol-separable} It is separable: This follows since
$\sections[\hol]{\man{E}}$ is a closed subspace of
$\sections[\infty]{\man{E}}$ by~\cite[Theorem~8.2]{AK/PWM:97} and since
subspaces of separable metric spaces are separable~\cite[Theorems~16.2,~16.9
and~16.11]{SW:04}\@.
\item It is nuclear:~\cite[Theorem~8.2]{AK/PWM:97}\@.  Note that, when
$\man{M}$ is compact, $p^{\hol}_{\man{M}}$ is a norm for the
$\C^{\hol}$-topology.  A consequence of this is that
$\sections[\hol]{\man{E}}$ must be finite-dimensional in these cases since the
only nuclear normed vector spaces are those that are
finite-dimensional~\cite[Theorem~4.4.14]{AP:69}\@.
\item \label{enum:COhol-suslin} It is Suslin: This follows since
$\sections[\hol]{\man{E}}$ is a Polish space, as we have seen above, at least
when the base manifold is Stein.
\end{compactenum}

Being metrisable, it suffices to describe the $\CO^{\hol}$-topology by
describing its convergent sequences; these are more or less obviously the
sequences that converge uniformly on every compact set.

As with spaces of smooth sections, we are interested in the fact that
nuclearity of $\sections[\hol]{\man{E}}$ implies that compact sets are
exactly those sets that are closed and von Neumann bounded.  The following
result is obvious in the same way that Lemma~\ref{lem:COinftybdd} is obvious
once one understands Theorem~1.37(b) from~\cite{WR:91}\@.
\begin{lemma}\label{lem:CObdd}
A subset\/ $\nbhd{B}\subset\sections[\hol]{\man{E}}$ is bounded in
the von Neumann bornology if and only if the following property holds: for
any compact set\/ $K\subset\man{M}$\@, there exists\/ $C\in\realp$ such
that\/ $p^{\hol}_K(\xi)\le C$ for every\/ $\xi\in\nbhd{B}$\@.
\end{lemma}

\subsection{The weak-$\sL$ topology for holomorphic vector fields}

As in the smooth case, one simply specialises the constructions for general
vector bundles to get the $\CO^{\hol}$-topology for the space
$\sections[\hol]{\tb{\man{M}}}$ of holomorphic vector fields and the space
$\func[\hol]{\man{M}}$ of holomorphic functions, noting that an holomorphic
function is obviously identified with a section of the trivial holomorphic
vector bundle $\man{M}\times\complex$\@.

As with smooth vector fields, for holomorphic vector fields we can seek a
weak-$\sL$ characterisation of the $\CO^{\hol}$-topology.  To begin, we need
to understand the Lie derivative in the holomorphic case.  Thinking of
$\func[\hol]{\man{M}}\subset\func[\infty]{\man{M}}\otimes\complex$ and using
the Wirtinger formulae,
\begin{equation*}
\pderiv{}{z^j}=\frac{1}{2}\Bigl(\pderiv{}{x^j}-\imag\pderiv{}{y^j}\Bigr),
\quad\pderiv{}{\bar{z}\null^j}=\frac{1}{2}\Bigl(\pderiv{}{x^j}+
\imag\pderiv{}{y^j}\Bigr),\qquad j\in\{1,\dots,n\},
\end{equation*}
in an holomorphic chart, one sees that the usual differential of a
$\complex$-valued function can be decomposed as $\d_{\complex}f=\partial
f+\bar{\partial}f$\@, the first term on the right corresponding to
``$\pderiv{}{z}$'' and the second to ``$\pderiv{}{\bar{z}}$\@.''  For
holomorphic functions, the Cauchy\textendash{}Riemann
equations~\cite[page~27]{SGK:92} imply that $\d_{\complex}f=\partial f$\@.
Thus we define the Lie derivative of an holomorphic function $f$ with respect
to an holomorphic vector field $X$ by $\lieder{X}{f}=\natpair{\partial
f}{X}$\@.  Fortunately, in coordinates this assumes the expected form:
\begin{equation*}
\lieder{X}{f}=\sum_{j=1}^nX^j\pderiv{f}{z^j}.
\end{equation*}
It is \emph{not} the case that on a general holomorphic manifold there is a
correspondence between derivations of the $\complex$-algebra
$\func[\hol]{\man{M}}$ and holomorphic vector fields by Lie
differentiation.\footnote{For example, on a compact holomorphic manifold, the
only holomorphic functions are locally
constant~\cite[Corollary~IV.1.3]{KF/HG:02}\@, and so the only derivation is
the zero derivation.  However, the $\complex$-vector space of holomorphic
vector fields, while not large, may have positive dimension.  For example,
the space of holomorphic vector fields on the Riemann sphere has
$\complex$-dimension three~\cite[Problem~17.9]{YI/SY:08}\@.}  However, for a
certain class of holomorphic manifolds, those known as ``Stein manifolds,''
the exact correspondence between derivations of the $\complex$-algebra
$\func[\hol]{\man{M}}$ and holomorphic vector fields under Lie
differentiation \emph{does} hold~\cite{JG:81}\@.  This is good news for us,
since Stein manifolds are intimately connected with real analytic manifolds,
as we shall see in the next section.

With the preceding discussion in mind, we can move ahead with
Definition~\ref{def:weakA} with $\alg{U}=\sections[\hol]{\tb{\man{M}}}$\@,
$\alg{V}=\func[\hol]{\man{M}}$ (with the $\CO^{\hol}$-topology), and
$\sA=\setdef{\sL_f}{f\in\func[\hol]{\man{M}}}$\@, where
\begin{equation*}
\mapdef{\sL_f}{\sections[\hol]{\tb{\man{M}}}}{\func[\hol]{\man{M}}}
{X}{\lieder{X}{f}.}
\end{equation*}
We make the following definition.
\begin{definition}\label{def:COhol-weak}
For an holomorphic manifold $\man{M}$\@, the \defn{weak-$\sL$ topology} for
$\sections[\hol]{\tb{\man{M}}}$ is the weakest topology for which $\sL_f$ is
continuous for every $f\in\func[\hol]{\man{M}}$\@, if $\func[\hol]{\man{M}}$
has the $\CO^{\hol}$-topology.\oprocend
\end{definition}

We then have the following result.
\begin{theorem}\label{the:COhol-weak}
For a Stein manifold\/ $\man{M}$\@, the following topologies for\/
$\sections[\hol]{\tb{\man{M}}}$ agree:
\begin{compactenum}[(i)]
\item \label{pl:COhol-weak1} the $\CO^{\hol}$-topology;
\item \label{pl:COhol-weak2} the weak-$\sL$ topology.
\end{compactenum}
\begin{proof}
\eqref{pl:COhol-weak1}$\subset$\eqref{pl:COhol-weak2} As we argued in the
proof of the corresponding assertion of Theorem~\ref{the:COinfty-weak}\@, it
suffices to show that
\begin{equation*}
p^{\hol}_K(X)\le C_1p^{\hol}_{K_1}(\lieder{X}{f^1})+\dots+
C_rp^{\hol}_{K_r}(\lieder{X}{f^r})
\end{equation*}
for some $C_1,\dots,C_r\in\realp$\@, some $K_1,\dots,K_r\subset\man{M}$
compact, and some $f^1,\dots,f^r\in\func[\hol]{\man{M}}$\@.

Let $K\subset\man{M}$ be compact.  For simplicity, we assume that $\man{M}$
is connected and so has a well-defined dimension $n$\@.  If not, then the
arguments are easily modified by change of notation to account for this.
Since $\man{M}$ is a Stein manifold, for every $z\in K$ there exists a
coordinate chart $(\nbhd{U}_z,\phi_z)$ with coordinate functions
$\map{z^1,\dots,z^n}{\nbhd{U}_z}{\complex}$ that are restrictions to
$\nbhd{U}_z$ of globally defined holomorphic functions on $\man{M}$\@.
Depending on your source, this is either a theorem or part of the definition
of a Stein manifold~\cite{KF/HG:02,LH:73}\@.  Thus, for
$j\in\{1,\dots,n\}$\@, let $f^j_z\in\func[\hol]{\man{M}}$ be the holomorphic
function which, when restricted to $\nbhd{U}_z$\@, gives the coordinate
function $z^j$\@.  Clearly, $\lieder{X}{f^j_z}=X^j$ on $\nbhd{U}_z$\@.  Also,
there exists $C_z\in\realp$ such that
\begin{equation*}
\dnorm{X(\zeta)}_{\metric}\le
C_z(\snorm{X^1(\zeta)}+\dots+\snorm{X^n(\zeta)}),
\qquad\zeta\in\closure(\nbhd{V}_z),
\end{equation*}
for some relatively compact neighbourhood $\nbhd{V}_z\subset\nbhd{U}_z$ of
$z$ (this follows from the fact that all norms are equivalent to the $\ell^1$
norm for $\complex^n$).  Thus
\begin{equation*}
\dnorm{X(\zeta)}_{\metric}\le
C_z(\snorm{\lieder{X}{f_z^1}(\zeta)}+\dots+
\snorm{\lieder{X}{f_z^n}(\zeta)}),
\qquad\zeta\in\closure(\nbhd{V}_z).
\end{equation*}
Let $z_1,\dots,z_k\in K$ be such that $K\subset\cup_{j=1}^k\nbhd{V}_{z_j}$\@.
Let $f^1,\dots,f^{kn}$ be the list of globally defined holomorphic functions
\begin{equation*}
f^1_{z_1},\dots,f^n_{z_1},\dots,f^1_{z_k},\dots,f^n_{z_k}
\end{equation*}
and let $C_1,\dots,C_{kn}$ be the list of coefficients
\begin{equation*}
\underbrace{C_{z_1},\dots,C_{z_1}}_{n\ \text{times}},\dots,
\underbrace{C_{z_k},\dots,C_{z_k}}_{n\ \text{times}}.
\end{equation*}
If $z\in K$\@, then $z\in\nbhd{V}_{z_j}$ for some $j\in\{1,\dots,k\}$ and so
\begin{equation*}
\dnorm{X(z)}_{\metric}\le
C_1\snorm{\lieder{X}{f^1}(z)}+\dots+C_{kn}\snorm{\lieder{X}{f^{kn}}(z)},
\end{equation*}
which gives
\begin{equation*}
p^{\hol}_K(X)\le C_1p^{\hol}_K(\lieder{X}{f^1})+\dots+
C_{kn}p^{\hol}_K(\lieder{X}{f^{kn}}),
\end{equation*}
as needed.

\eqref{pl:COhol-weak2}$\subset$\eqref{pl:COhol-weak1} We claim that $\sL_f$
is continuous for every $f\in\func[\hol]{\man{M}}$ if
$\sections[\hol]{\tb{\man{M}}}$ has the $\CO^{\hol}$-topology.  Let
$K\subset\man{M}$ be compact and let $\nbhd{U}$ be a relatively compact
neighbourhood of $K$ in $\man{M}$\@.  Note that, for
$f\in\func[\hol]{\man{M}}$\@,
\begin{equation*}
p^{\hol}_K(\lieder{X}{f})\le Cp^{\hol}_{K,1}(f)p^{\hol}_K(X)\le
C'p^{\hol}_K(X),
\end{equation*}
using Proposition~\ref{prop:Cauchyest}\@, giving continuity of the identity
map if we provide the domain with the $\CO^{\hol}$-topology and the codomain
with the weak-$\sL$ topology,~\cf~\cite[\S{}III.1.1]{HHS/MPW:99}\@.  Thus
open sets in the weak-$\sL$ topology are contained in the
$\CO^{\hol}$-topology.
\end{proof}
\end{theorem}

As in the smooth case, we shall use the theorem according to the following
result.
\begin{corollary}\label{cor:COweak}
Let\/ $\man{M}$ be a Stein manifold, let\/ $(\ts{X},\sO)$ be a topological
space, let\/ $(\ts{T},\sM)$ be a measurable space, and let\/
$\map{\mu}{\sM}{\erealnn}$ be a finite measure.  The following statements hold:
\begin{compactenum}[(i)]
\item \label{pl:COweakbdd} a subset\/
$\nbhd{B}\subset\sections[\hol]{\tb{\man{M}}}$ is bounded in the von Neumann
bornology if and only if it is weak-$\sL$ bounded in the von Neumann
bornology;
\item \label{pl:COweakcont} a map\/
$\map{\Phi}{\ts{X}}{\sections[\hol]{\tb{\man{M}}}}$ is continuous if and only
if it is weak-$\sL$ continuous;
\item \label{pl:COweakmeas} a map\/
$\map{\Psi}{\ts{T}}{\sections[\hol]{\tb{\man{M}}}}$ is measurable if and only
if it is weak-$\sL$ measurable;
\item \label{pl:COweakint} a map\/
$\map{\Psi}{\ts{T}}{\sections[\hol]{\tb{\man{M}}}}$ is Bochner integrable if
and only if it is weak-$\sL$ Bochner integrable.
\end{compactenum}
\begin{proof}
As in the proof of Corollary~\ref{cor:COinftyweak}\@, we need to show that
$\setdef{\sL_f}{f\in\func[\hol]{\man{M}}}$ has a countable point separating
subset.  The argument here follows that in the smooth case, except that here
we have to use the properties of Stein manifolds,~\cf~the proof of the first
part of Theorem~\ref{the:COhol-weak} above, to assert the existence, for each
$z\in\man{M}$\@, of a neighbourhood on which there are globally defined
holomorphic functions whose differentials span the cotangent space at each
point.  Since $\sections[\hol]{\tb{\man{M}}}$ is complete, separable, and
Suslin, and since $\func[\hol]{\man{M}}$ is Suslin by
properties~$\CO^{\hol}$-\ref{enum:COhol-complete}\@,%
~$\CO^{\hol}$-\ref{enum:COhol-separable}
and~$\CO^{\hol}$-\ref{enum:COhol-suslin} above, the corollary follows from
Lemma~\ref{lem:weakA}\@.
\end{proof}
\end{corollary}

\section{The $\C^\omega$-topology for the space of real analytic vector fields}\label{sec:analytic-topology}

In this section we examine a topology on the set of real analytic vector
fields.  As we shall see, this requires some considerable effort.
\citet{AAA/RVG:78} consider the real analytic case by considering bounded
holomorphic extensions to neighbourhoods of $\real^n$ of fixed width
$\complex^n$\@.  Our approach is more general, more geometric, and global,
using a natural real analytic topology described, for example, in the work of
\citet{AM:66}\@.  This allows us to dramatically broaden the class of real
analytic systems that we can handle to include ``all'' analytic systems.

The first observation we make is that $\sections[\omega]{\man{E}}$ is not a
closed subspace of $\sections[\infty]{\man{E}}$ in the $\CO^\infty$-topology.
To see this, consider the following.  Take a smooth but not real analytic
function on $\sphere^1$\@.  The Fourier series of this function gives rise,
by taking partial sums, to a sequence of real analytic functions.  Standard
harmonic analysis~\cite[Theorem~VII.2.11(b)]{EMS/GW:71} shows that this
sequence and all of its derivatives converge uniformly, and so in the
$\CO^\infty$-topology, to the original function.  Thus we have a Cauchy
sequence in $\func[\omega]{\sphere^1}$ that does not converge, with respect
to the $\CO^\infty$-topology, in $\func[\omega]{\sphere^1}$\@.

The second observation we make is that a plain restriction of the topology
for holomorphic objects is not sufficient.  The reason for this is that, upon
complexification (a process we describe in detail below) there will not be a
uniform neighbourhood to which all real analytic objects can be extended.
Let us look at this for an example, where ``object'' is ``function.''  For
$r\in\realp$ we consider the real analytic function $\map{f_r}{\real}{\real}$
defined by $f_r(x)=\frac{r^2}{r^2+x^2}$\@.  We claim that there is no
neighbourhood $\ol{\nbhd{U}}$ of $\real$ in $\complex$ to which all of the
functions $f_r$\@, $r\in\realp$\@, can be extended.  Indeed, take some such
neighbourhood $\ol{\nbhd{U}}$ and let $r\in\realp$ be sufficiently small that
$\cdisk[]{r}{0}\subset\ol{\nbhd{U}}$\@.  To see that $f_r$ cannot be extended
to an holomorphic function $\ol{f}_r$ on $\ol{\nbhd{U}}$\@, let $\ol{f}_r$ be
such an holomorphic extension.  Then $\ol{f}_r(z)$ must be equal to
$\frac{r^2}{r^2+z^2}$ for $z\in\odisk[]{r}{0}$ by uniqueness of holomorphic
extensions~\cite[Lemma~5.40]{KC/YE:12}\@.  But this immediately prohibits
$\ol{f}_r$ from being holomorphic on any neighbourhood of $\cdisk[]{r}{0}$\@,
giving our claim.

Therefore, to topologise the space of real analytic vector fields, we will
need to do more than either~(1)~restrict the $\CO^\infty$-topology or~(2)~use
the $\CO^\hol$-topology in an ``obvious'' way.  Note that it is the
``obvious'' use of the $\CO^\hol$-topology for holomorphic objects that is
employed by \citet{AAA/RVG:78} in their study of time-varying real analytic
vector fields.  Moreover, \citet{AAA/RVG:78} also restrict to \emph{bounded}
holomorphic extensions.  What we propose is an improvement on this in that it
works far more generally, and is also more natural to a geometric treatment
of the real analytic setting.  We comment at this point that we shall see in
Theorems~\ref{the:Comega->Chol} and~\ref{the:Comega->Cholparam} below that
the consideration of bounded holomorphic extensions to fixed neighbourhoods
in the complexification is sometimes sufficient locally.  But conclusions
such as this become hard theorems with precise hypotheses in our approach,
not starting points for the theory.

As in the smooth and holomorphic cases, we begin by considering a general
vector bundle.

\subsection{A natural direct limit topology}\label{subsec:Comega-direct}

We let $\map{\pi}{\man{E}}{\man{M}}$ be a real analytic vector bundle.  We
shall extend $\man{E}$ to an holomorphic vector bundle that will serve an an
important device for all of our constructions.

\subsubsection{Complexifications}\label{subsubsec:complexify}

Let us take some time to explain how holomorphic extensions can be
constructed.  The following two paragraphs distill out important parts of
about forty years of intensive development of complex analysis, culminating
in the paper of \citet{HG:58}\@.

For simplicity, let us assume that $\man{M}$ is connected and so has pure
dimension, and so the fibres of $\man{E}$ also have a fixed dimension.  As in
Section~\ref{subsec:analytic-conn}\@, we suppose that we have a real analytic
affine connection $\nabla$ on $\man{M}$\@, a real analytic vector bundle
connection $\nabla^0$ on $\man{E}$\@, a real analytic Riemannian metric
$\metric$ on $\man{M}$\@, and a real analytic fibre metric $\metric_0$ on
$\man{E}$\@.  We also assume the data required to make the
diagram~\eqref{eq:EMembedding} giving $\map{\pi}{\man{E}}{\man{M}}$ as the
image of a real analytic vector bundle monomorphism in the trivial vector
bundle $\real^N\times\real^N$ for some suitable $N\in\integerp$\@.

Now we complexify.  Recall that, if $\alg{V}$ is a $\complex$-vector space,
then multiplication by $\sqrt{-1}$ induces a $\real$-linear map
$J\in\End_{\real}(\alg{V})$\@.  A $\real$-subspace $\alg{U}$ of $\alg{V}$ is
\defn{totally real} if $\alg{U}\cap J(\alg{U})=\{0\}$\@.  A submanifold of an
holomorphic manifold, thinking of the latter as a smooth manifold, is
\defn{totally real} if its tangent spaces are totally real subspaces.  By
\cite[Proposition~1]{HW/FB:59}\@, for a real analytic manifold $\man{M}$
there exists a complexification $\ol{\man{M}}$ of $\man{M}$\@,~\ie~an
holomorphic manifold having $\man{M}$ as a totally real submanifold and where
$\ol{\man{M}}$ has the same $\complex$-dimension as the $\real$-dimension of
$\man{M}$\@.  As shown by \citet[\S3.4]{HG:58}\@, for any neighbourhood
$\ol{\nbhd{U}}$ of $\man{M}$ in $\ol{\man{M}}$\@, there exists a Stein
neighbourhood $\ol{\nbhd{S}}$ of $\man{M}$ contained in $\ol{\nbhd{U}}$\@.
By arguments involving extending convergent real power series to convergent
complex power series (the conditions on coefficients for convergence are the
same for both real and complex power series), one can show that there is an
holomorphic extension of $\iota_{\man{M}}$ to
$\map{\iota_{\ol{\man{M}}}}{\ol{\man{M}}}{\complex^N}$\@, possibly after
shrinking $\ol{\man{M}}$~\cite[Lemma~5.40]{KC/YE:12}\@.  By applying similar
reasoning to the transition maps for the real analytic vector bundle
$\man{E}$\@, one obtains an holomorphic vector bundle
$\map{\ol{\pi}}{\ol{\man{E}}}{\ol{\man{M}}}$ for which the diagram
\begin{equation*}
\xymatrix{{\ol{\man{E}}}\ar[rrr]\ar[ddd]_{\ol{\pi}}&&&
{\complex^N\times\complex^N}\ar[ddd]^{\pr_2}\\
&{\man{E}}\ar[r]^(0.35){\hat{\iota}_{\man{E}}}\ar[d]_{\pi}\ar[lu]&
{\real^N\times\real^N}\ar[d]^{\pr_2}\ar[ru]&\\
&{\man{M}}\ar[r]_{\iota_{\man{M}}}\ar[ld]&{\real^N}\ar[rd]\\
{\ol{\man{M}}}\ar[rrr]_{\iota_{\ol{\man{M}}}}&&&{\complex^N}}
\end{equation*}
commutes, where all diagonal arrows are complexification and where the inner
diagram is as defined in the proof of Lemma~\ref{lem:analytic-conn}\@.  One
can then define an Hermitian fibre metric $\ol{\metric}_0$ on $\ol{\man{E}}$
induced from the standard Hermitian metric on the fibres of the vector bundle
$\complex^N\times\complex^N$ and an Hermitian metric $\ol{\metric}$ on
$\ol{\man{M}}$ induced from the standard Hermitian metric on $\complex^N$\@.

In the remainder of this section, we assume that the preceding constructions
have been done and fixed once and for all.

\subsubsection{Germs of holomorphic sections over subsets of a real analytic
manifold}\label{subsubsec:hologerms}

In two different places, we will need to consider germs of holomorphic
sections.  In this section we organise the methodology for doing this to
unify the notation.

Let $A\subset\man{M}$ and let $\sN_A$ be the set of neighbourhoods of $A$ in
the complexification $\ol{\man{M}}$\@.  For
$\ol{\nbhd{U}},\ol{\nbhd{V}}\in\sN_A$\@, and for
$\ol{\xi}\in\sections[\hol]{\ol{\man{E}}|\ol{\nbhd{U}}}$ and
$\ol{\eta}\in\sections[\hol]{\ol{\man{E}}|\ol{\nbhd{V}}}$\@, we say that
$\ol{\xi}$ is \defn{equivalent} to $\ol{\eta}$ if there exist
$\ol{\nbhd{W}}\in\sN_A$ and
$\ol{\zeta}\in\sections[\hol]{\ol{\man{E}}|\ol{\nbhd{W}}}$ such that
$\ol{\nbhd{W}}\subset\ol{\nbhd{U}}\cap\ol{\nbhd{V}}$ and such that
\begin{equation*}
\ol{\xi}|\ol{\nbhd{W}}=\ol{\eta}|\ol{\nbhd{W}}=\ol{\zeta}.
\end{equation*}
By $\gsections[\hol]{A}{\ol{\man{E}}}$ we denote the set of equivalence
classes, which we call the set of germs of sections of $\ol{\man{E}}$ over
$A$\@.  By $[\ol{\xi}]_A$ we denote the equivalence class of
$\ol{\xi}\in\sections[\hol]{\ol{\man{E}}|\ol{\nbhd{U}}}$ for some
$\ol{\nbhd{U}}\in\sN_A$\@.

Now, for $x\in\man{M}$\@, $\man{E}_x$ is a totally real subspace of
$\ol{\man{E}}_x$ with half the real dimension, and so it follows that
\begin{equation*}
\ol{\man{E}}_x=\man{E}_x\oplus J(\man{E}_x),
\end{equation*}
where $J$ is the complex structure on the fibres of $\ol{\man{E}}$\@.  For
$\ol{\nbhd{U}}\in\sN_A$\@, denote by
$\sections[\hol,\real]{\ol{\man{E}}|\ol{\nbhd{U}}}$ those holomorphic
sections $\ol{\xi}$ of $\ol{\man{E}}|\ol{\nbhd{U}}$ such that
$\ol{\xi}(x)\in\man{E}_x$ for $x\in\ol{\nbhd{U}}\cap\man{M}$\@.  We think of
this as being a locally convex topological $\real$-vector space with the
seminorms $p^{\hol}_{\ol{K}}$\@, $\ol{K}\subset\ol{\nbhd{U}}$ compact,
defined by
\begin{equation*}
p^{\hol}_{\ol{K}}(\ol{\xi})=
\sup\setdef{\dnorm{\ol{\xi}(\ol{x})}_{\ol{\metric}_0}}{\ol{x}\in\ol{K}},
\end{equation*}
\ie~we use the locally convex structure induced from the usual
$\CO^\hol$-topology on $\sections[\hol]{\ol{\man{E}}|\ol{\nbhd{U}}}$\@.
\begin{remark}\label{rem:real-closed}
We note that $\sections[\hol,\real]{\ol{\man{E}}|\ol{\nbhd{U}}}$ is a closed
$\real$-subspace of $\sections[\hol]{\ol{\man{E}}}$ in the
$\CO^\hol$-topology,~\ie~the restriction of requiring ``realness'' on
$\man{M}$ is a closed condition.  This is easily shown, and we often assume
it often without mention.\oprocend
\end{remark}

Denote by $\gsections[\hol,\real]{A}{\ol{\man{E}}}$ the set of germs of
sections from $\sections[\hol,\real]{\ol{\man{E}}|\ol{\nbhd{U}}}$\@,
$\ol{\nbhd{U}}\in\sN_A$\@.  If $\ol{\nbhd{U}}_1,\ol{\nbhd{U}}_2\in\sN_A$
satisfy $\ol{\nbhd{U}}_1\subset\ol{\nbhd{U}}_2$\@, then we have the
restriction mapping
\begin{equation*}
\mapdef{r_{\ol{\nbhd{U}}_2,\ol{\nbhd{U}}_1}}
{\sections[\hol,\real]{\ol{\man{E}}|\ol{\nbhd{U}}_2}}
{\sections[\hol,\real]{\ol{\man{E}}|\ol{\nbhd{U}}_1}}
{\ol{\xi}}{\ol{\xi}|\ol{\nbhd{U}}_1.}
\end{equation*}
This restriction is continuous since, for any compact set
$\ol{K}\subset\ol{\nbhd{U}}_1\subset\ol{\nbhd{U}}_2$ and any
$\ol{\xi}\in\sections[\hol,\real]{\ol{\man{E}}|\ol{\nbhd{U}}_2}$\@,
we have
$p^{\hol}_{\ol{K}}(r_{\ol{\nbhd{U}}_2,\ol{\nbhd{U}}_1}(\ol{\xi}))\le
p^{\hol}_{\ol{K}}(\ol{\xi})$ (in fact we have equality, but the
inequality emphasises what is required for our assertion to be
true~\cite[\S{}III.1.1]{HHS/MPW:99}).  We also have maps
\begin{equation*}
\mapdef{r_{\ol{\nbhd{U}},A}}
{\sections[\hol,\real]{\ol{\man{E}}|\ol{\nbhd{U}}}}
{\gsections[\hol,\real]{A}{\ol{\man{E}}}}
{\ol{\xi}}{[\ol{\xi}]_A.}
\end{equation*}
Note that $\sN_A$ is a directed set by inclusion; that is,
$\ol{\nbhd{U}}_2\preceq\ol{\nbhd{U}}_1$ if
$\ol{\nbhd{U}}_1\subset\ol{\nbhd{U}}_2$\@.  Thus we have the directed system
$\ifam{\sections[\hol,\real]
{\tb{\ol{\nbhd{U}}}}}_{\ol{\nbhd{U}}\in\sN_A}$\@, along with the mappings
$r_{\ol{\nbhd{U}}_2,\ol{\nbhd{U}}_1}$\@, in the category of locally convex
topological $\real$-vector spaces.  The usual notion of direct limit in the
category of $\real$-vector spaces gives
$\gsections[\hol,\real]{A}{\ol{\man{E}}}$\@, along with the linear mappings
$r_{\ol{\nbhd{U}},A}$\@, $\ol{\nbhd{U}}\in\sN_A$\@, as the direct limit of
this directed system~\cite[\cf][Theorem~III.10.1]{SL:02b}\@.  This vector
space then has the finest locally convex topology making the maps
$r_{\ol{\nbhd{U}},A}$\@, $\ol{\nbhd{U}}\in\sN_A$\@, continuous,~\ie~the
direct limit in the category of locally convex topological vector spaces.  We
refer to this as the \defn{direct limit topology} for
$\gsections[\hol,\real]{A}{\ol{\man{E}}}$\@.

\subsubsection{The direct limit topology}\label{subsubsec:Comega-direct}

We shall describe four topologies (or more, depending on which descriptions
you regard as being distinct) for the space of real analytic sections of a
real analytic vector bundle.  The first is quite direct, involving an
application of the construction above to the case of $A=\man{M}$\@.  In this
case, the following lemma is key to our constructions.
\begin{lemma}\label{lem:section<->germ}
There is a natural\/ $\real$-vector space isomorphism between\/
$\sections[\omega]{\man{E}}$ and\/
$\gsections[\hol,\real]{\man{M}}{\ol{\man{E}}}$\@.
\begin{proof}
Let $\xi\in\sections[\omega]{\man{E}}$\@.  As
in~\cite[Lemma~5.40]{KC/YE:12}\@, there is an extension of $\xi$ to a section
$\ol{\xi}\in\sections[\hol,\real]{\ol{\man{E}}|\ol{\nbhd{U}}}$ for some
$\ol{\nbhd{U}}\in\sN_{\man{M}}$\@.  We claim that the map
$\map{i_{\man{M}}}{\sections[\omega]{\man{E}}}
{\gsections[\hol,\real]{\man{M}}{\ol{\man{E}}}}$ defined by
$i_{\man{M}}(\xi)=[\ol{\xi}]_{\man{M}}$ is the desired isomorphism.  That
$i_{\man{M}}$ is independent of the choice of extension $\ol{\xi}$ is a
consequence of the fact that the extension to $\ol{\xi}$ is unique inasmuch
as any two such extensions agree on some neighbourhood contained in their
intersection; this is the uniqueness assertion
of~\cite[Lemma~5.40]{KC/YE:12}\@.  This fact also ensures that $i_{\man{M}}$
is injective.  For surjectivity, let $[\ol{\xi}]_{\man{M}}\in
\gsections[\hol,\real]{\man{M}}{\ol{\man{E}}}$ and let us define
$\map{\xi}{\man{M}}{\man{E}}$ by $\xi(x)=\ol{\xi}(x)$ for $x\in\man{M}$\@.
Note that the restriction of $\ol{\xi}$ to $\man{M}$ is real analytic because
the values of $\ol{\xi}|\man{M}$ at points in a neighbourhood of
$x\in\man{M}$ are given by the restriction of the (necessarily convergent)
$\complex$-Taylor series of $\ol{\xi}$ to $\man{M}$\@.  Obviously,
$i_{\man{M}}(\xi)=[\ol{\xi}]_{\man{M}}$\@.
\end{proof}
\end{lemma}

Now we use the direct limit topology on
$\gsections[\hol,\real]{\man{M}}{\ol{\man{E}}}$ described above, along with
the preceding lemma, to immediately give a locally convex topology for
$\sections[\omega]{\man{E}}$ that we refer to as the \defn{direct
$\C^\omega$-topology}\@.

Let us make an important observation about the direct $\C^\omega$-topology.
Let us denote by $\sS_{\man{M}}$ the set of all Stein neighbourhoods of
$\man{M}$ in $\ol{\man{M}}$\@.  As shown by \citet[\S3.4]{HG:58}\@, if
$\ol{\nbhd{U}}\in\sN_{\man{M}}$ then there exists
$\ol{\nbhd{S}}\in\sS_{\man{M}}$ with $\ol{\nbhd{S}}\subset\ol{\nbhd{U}}$\@.
Therefore, $\sS_{\man{M}}$ is cofinal in $\sN_{\man{M}}$ and so the directed
systems $\ifam{\sections[\hol]{\man{E}|
\ol{\nbhd{U}}}}_{\ol{\nbhd{U}}\in\sN_{\man{M}}}$ and
$\ifam{\sections[\hol]{\man{E}|
\ol{\nbhd{S}}}}_{\ol{\nbhd{S}}\in\sS_{\man{M}}}$ induce the same final
topology on $\sections[\omega]{\man{E}}$~\cite[page~137]{AG:73}\@.

\subsection{Topologies for germs of holomorphic functions about compact sets}\label{subsec:GK-topology}

In the preceding section, we gave a more or less direct description of a
topology for the space of real analytic sections.  This description has a
benefit of being the one that one might naturally arrive at after some
thought.  However, there is not a lot that one can do with this description
of the topology.  In this section we develop the means by which one can
consider alternative descriptions of this topology that, for example, lead to
explicit seminorms for the topology on the space of real analytic sections.
These seminorms will be an essential part of our developing a useful theory
for time-varying real analytic vector fields and real analytic control
systems.

\subsubsection{The direct limit topology for the space of germs about a
compact set}

We continue with the notation from Section~\ref{subsubsec:hologerms}\@.  For
$K\subset\man{M}$ compact, we have the direct limit topology, described above
for general subsets $A\subset\man{M}$\@, on
$\gsections[\hol,\real]{K}{\ol{\man{E}}}$\@.  We seem to have gained nothing,
since we have yet another direct limit topology.  However, the direct limit
can be shown to be of a friendly sort as follows.  Unlike the general
situation, since $K$ is compact there is a \emph{countable} family
$\ifam{\ol{\nbhd{U}}_{K,j}}_{j\in\integerp}$ from $\sN_K$ with the property
that $\closure(\ol{\nbhd{U}}_{K,j+1})\subset\ol{\nbhd{U}}_{K,j}$ and
$K=\cap_{j\in\integerp}\ol{\nbhd{U}}_{K,j}$\@.  Moreover, the sequence
$\ifam{\ol{\nbhd{U}}_{K,j}}_{j\in\integerp}$ is cofinal in $\sN_K$\@,~\ie~if
$\ol{\nbhd{U}}\in\sN_K$\@, then there exists $j\in\integerp$ with
$\ol{\nbhd{U}}_{K,j}\subset\ol{\nbhd{U}}$\@.  Let us fix such a family of
neighbourhoods.  Let us fix $j\in\integerp$ for a moment.  Let
$\secbdd[\hol,\real]{\ol{\man{E}}|\ol{\nbhd{U}}_{K,j}}$ be the set of bounded
sections from $\sections[\hol,\real]{\ol{\man{E}}|\ol{\nbhd{U}}_{K,j}}$\@,
boundedness being taken relative to the Hermitian fibre metric
$\ol{\metric}_0$\@.  As we have seen in Lemma~\ref{lem:secbddhol}\@, if we
define a norm on $\secbdd[\hol,\real]{\ol{\man{E}}|\ol{\nbhd{U}}_{K,j}}$ by
\begin{equation*}
p^{\hol}_{\ol{\nbhd{U}}_{K,j},\infty}(\ol{\xi})=
\sup\setdef{\dnorm{\ol{\xi}(\ol{x})}_{\ol{\metric}_0}}
{\ol{x}\in\ol{\nbhd{U}}_{K,j}},
\end{equation*}
then this makes $\secbdd[\hol,\real]{\ol{\nbhd{U}}_{K,j}}$ into a Banach
space, a closed subspace of the Banach space of bounded continuous sections
of $\ol{\man{E}}|\ol{\nbhd{U}}_{K,j}$\@.  Now, no longer fixing $j$\@, we
have a sequence of inclusions
\begin{multline*}
\secbdd[\hol,\real]{\ol{\man{E}}|\ol{\nbhd{U}}_{K,1}}\subset
\sections[\hol,\real]{\ol{\man{E}}|\ol{\nbhd{U}}_{K,1}}\subset
\secbdd[\hol,\real]{\ol{\man{E}}|\ol{\nbhd{U}}_{K,2}}\subset\\
\dots\subset
\sections[\hol,\real]{\ol{\man{E}}|\ol{\nbhd{U}}_{K,j}}\subset
\secbdd[\hol,\real]{\ol{\man{E}}|\ol{\nbhd{U}}_{K,j+1}}\subset
\cdots.
\end{multline*}
The inclusion $\sections[\hol,\real]{\ol{\nbhd{U}}_{K,j}}\subset
\secbdd[\hol,\real]{\ol{\nbhd{U}}_{K,j+1}}$\@, $j\in\integerp$\@, is by
restriction from $\ol{\nbhd{U}}_{K,j}$ to the smaller
$\ol{\nbhd{U}}_{K,j+1}$\@, keeping in mind that
$\closure(\ol{\nbhd{U}}_{K,j+1})\subset\ol{\nbhd{U}}_{K,j}$\@.  By
Lemma~\ref{lem:secbddhol}\@, all inclusions are continuous.  For
$j\in\integerp$ define
\begin{equation}\label{eq:Gammahol->Ghol}
\mapdef{r_{K,j}}{\secbdd[\hol,\real]
{\ol{\man{E}}|\ol{\nbhd{U}}_{K,j}}}
{\gsections[\hol,\real]{K}{\ol{\man{E}}}}{\ol{\xi}}{[\ol{\xi}]_K.}
\end{equation}
Now one can show that the direct limit topologies induced on
$\gsections[\hol,\real]{K}{\ol{\man{E}}}$ by the directed system
$\ifam{\sections[\hol,\real]{\ol{\man{E}}|
\ol{\nbhd{U}}}}_{\ol{\nbhd{U}}\in\sN_K}$ of Fr\'echet spaces and by the
directed system $\ifam{\secbdd[\hol,\real]{\ol{\man{E}}|
\ol{\nbhd{U}}_{K,j}}}_{j\in\integerp}$ of Banach spaces
agree~\cite[Theorem~8.4]{AK/PWM:97}\@.  We refer to~\cite{KDB:88}\@, starting
on page~63, for a fairly comprehensive discussion of the topology we have
just described in the context of germs of holomorphic functions about a compact
subset $K\subset\complex^n$\@.

\subsubsection{A weighted direct limit topology for sections of bundles of infinite jets}\label{subsubsec:E(K)}

Here we provide a direct limit topology for a subspace of the space of
continuous sections of the infinite jet bundle of a vector bundle.  Below we
shall connect this direct limit topology to the direct limit topology
described above for germs of holomorphic sections about a compact set.  The
topology we give here has the advantage of providing explicit seminorms for
the topology of germs, and subsequently for the space of real analytic
sections.

For this description, we work with infinite jets, so let us introduce the
notation we will use for this, referring to~\cite[Chapter~7]{DJS:89} for
details.  Let us denote by $\jet{\infty}{\man{E}}$ the bundle of infinite
jets of a vector bundle $\map{\pi}{\man{E}}{\man{M}}$\@, this being the
inverse limit (in the category of sets, for the moment) of the inverse system
$\ifam{\jet{m}{\man{E}}}_{m\in\integernn}$ with mappings $\pi^{m+1}_m$\@,
$m\in\integernn$\@.  Precisely,
\begin{equation*}
\jet{\infty}{\man{E}}=\Bigsetdef{\phi\in\prod_{m\in\integernn}\jet{m}{\man{E}}}
{\pi^k_l\scirc\phi(k)=\phi(l),\ k,l\in\integernn,\ k\ge l}.
\end{equation*}
We let $\map{\pi^\infty_m}{\jet{\infty}{\man{E}}}{\jet{m}{\man{E}}}$ be the
projection defined by $\pi^\infty_m(\phi)=\phi(m)$\@.  For
$\xi\in\sections[\infty]{\man{E}}$ we let
$\map{j_\infty\xi}{\man{M}}{\jet{\infty}{\man{E}}}$ be defined by
$\pi^\infty_m\scirc j_\infty\xi(x)=j_m\xi(x)$\@.  By a theorem of
\citet{EB:95}\@, if $\phi\in\jet{\infty}{\man{E}}$\@, there exist
$\xi\in\sections[\infty]{\man{E}}$ and $x\in\man{M}$ such that
$j_\infty\xi(x)=\phi$\@.  We can define sections of $\jet{\infty}{\man{E}}$ in
the usual manner: a section is a map
$\map{\Xi}{\man{M}}{\jet{\infty}{\man{E}}}$ satisfying
$\pi^\infty_0\scirc\Xi(x)=x$ for every $x\in\man{M}$\@.  We shall equip
$\jet{\infty}{\man{E}}$ with the initial topology so that a section $\Xi$ is
continuous if and only if $\pi^\infty_m\scirc\Xi$ is continuous for every
$m\in\integernn$\@.  We denote the space of continuous sections of
$\jet{\infty}{\man{E}}$ by $\sections[0]{\jet{\infty}{\man{E}}}$\@.  Since we
are only dealing with continuous sections, we can talk about sections defined
on any subset $A\subset\man{M}$\@, using the relative topology on $A$\@.  The
continuous sections defined on $A\subset\man{M}$ will be denoted by
$\sections[0]{\jet{\infty}{\man{E}}|A}$\@.

Now let $K\subset\man{M}$ be compact and, for $j\in\integerp$\@, denote
\begin{equation*}
\sE_j(K)=\setdef{\Xi\in\sections[0]{\jet{\infty}{\man{E}}|K}}
{\sup\setdef{j^{-m}\dnorm{\pi^\infty_m\scirc
\Xi(x)}_{\ol{\metric}_m}}{m\in\integernn,\ x\in K}<\infty},
\end{equation*}
and on $\sE_j(K)$ we define a norm $p_{K,j}$ by
\begin{equation*}
p_{K,j}(\Xi)=\sup\setdef{j^{-m}\dnorm{\pi^\infty_m\scirc
\Xi(x)}_{\ol{\metric}_m}}{m\in\integernn,\ x\in K}.
\end{equation*}
One readily verifies that, for each $j\in\integerp$\@, $(\sE_j(K),p_{K,j})$
is a Banach space.  Note that $\sE_j(K)\subset\sE_{j+1}(K)$ and that
$p_{K,j+1}(\Xi)\le p_{K,j}(\Xi)$ for $\Xi\in\sE_j(K)$\@, and so the inclusion
of $\sE_j(K)$ in $\sE_{j+1}(K)$ is continuous.  We let $\sE(K)$ be the direct
limit of the directed system $\ifam{\sE_j(K)}_{j\in\integerp}$\@.

We shall subsequently explore more closely the relationship between the
direct limit topology for $\sE(K)$ and the topology for
$\gsections[\hol,\real]{K}{\ol{\man{E}}}$\@.  For now, we merely observe that
the direct limit topology for $\sE(K)$ admits a characterisation by
seminorms.  To state the result, let us denote by $\csd(\integernn;\realp)$
the set of nonincreasing sequences $\ifam{a_m}_{m\in\integernn}$ in $\realp$
that converge to $0$\@.  Let us abbreviate such a sequence by
$\vect{a}=\ifam{a_m}_{m\in\integernn}$\@.  The following result is modelled
after~\cite[Lemma~1]{DV:13}\@.
\begin{lemma}\label{lem:EKseminorms}
The direct limit topology for\/ $\sE(K)$ is defined by the seminorms
\begin{equation*}
p_{K,\vect{a}}=\sup\setdef{a_0a_1\cdots a_m
\dnorm{\pi^\infty_m\scirc\Xi(x)}_{\ol{\metric}_m}}
{m\in\integernn,\ x\in K},
\end{equation*}
for\/ $\vect{a}\in\csd(\integernn;\realp)$\@.
\begin{proof}
First we show that the seminorms $p_{K,\vect{a}}$\@,
$\vect{a}\in\csd(\integernn;\realp)$\@, are continuous on $\sE(K)$\@.  It
suffices to show that $p_{K,\vect{a}}|\sE_j(K)$ is continuous for each
$j\in\integerp$~\cite[Proposition~IV.5.7]{JBC:90}\@.  Thus, since $\sE_j(K)$
is a Banach space, it suffices to show that, if
$\ifam{\Xi_k}_{k\in\integerp}$ is a sequence in $\sE_j(K)$ converging to
zero, then $\lim_{k\to\infty}p_{K,\vect{a}}(\Xi_k)=0$\@.  Let
$N\in\integernn$ be such that $a_N<\frac{1}{j}$\@.  Let $C\ge1$ be such that
\begin{equation*}
a_0a_1\cdots a_m\le Cj^{-m},\qquad m\in\{0,1,\dots,N\},
\end{equation*}
this being possible since there are only finitely many inequalities to
satisfy.  Therefore, for any $m\in\integernn$\@, we have $a_0a_1\cdots a_m\le
Cj^{-m}$\@.  Then, for any $\Xi\in\sections[0]{\jet{\infty}{\man{E}}|K}$\@,
\begin{equation*}
a_0a_1\cdots a_m\dnorm{\pi^\infty_m\scirc\Xi(x)}_{\ol{\metric}_m}
\le Cj^{-m}\dnorm{\pi^\infty_m\scirc\Xi(x)}_{\ol{\metric}_m}
\end{equation*}
for every $x\in K$ and $m\in\integernn$\@.  From this we immediately have
$\lim_{k\to\infty}p_{K,\vect{a}}(\Xi_k)=0$\@, as desired.  This shows that
the direct limit topology on $\sE(K)$ is stronger than the topology defined
by the family of seminorms $p_{K,\vect{a}}$\@,
$\vect{a}\in\csd(\integernn;\realp)$\@.

For the converse, we show that every neighbourhood of $0\in\sE(K)$ in the
direct limit topology contains a neighbourhood of zero in the topology
defined by the seminorms $p_{K,\vect{a}}$\@,
$\vect{a}\in\csd(\integernn;\realp)$\@.  Let $\nbhd{B}_j$ denote the unit
ball in $\sE_j(K)$\@.  A neighbourhood of $0$ in the direct limit topology
contains a union of balls $\epsilon_j\nbhd{B}_j$ for some
$\epsilon_j\in\realp$\@, $j\in\integerp$\@, (see~\cite[page~54]{HHS/MPW:99})
and we can assume, without loss of generality, that
$\epsilon_j\in\interval(0,1)$ for each $j\in\integerp$\@.  We define an
increasing sequence $\ifam{m_j}_{j\in\integerp}$ in $\integernn$ as follows.
Let $m_1=0$\@.  Having defined $m_1,\dots,m_j$\@, define $m_{j+1}>m_j$ by
requiring that $j<\epsilon_{j+1}^{1/m_{j+1}}(j+1)$\@.  For
$m\in\{m_j,\dots,m_{j+1}-1\}$\@, define $a_m\in\realp$ by
$a_m^{-1}=\epsilon^{1/m_j}_jj$\@.  Note that, for
$m\in\{m_j,\dots,m_{j+1}-1\}$\@, we have
\begin{equation*}
a_m^{-m}=\epsilon_j^{m/m_j}j^m\le\epsilon_jj^m.
\end{equation*}
Note that $\lim_{m\to\infty}a_m=0$\@.  If
$\Xi\in\sections[0]{\jet{\infty}{\man{E}}|K}$ satisfies
$p_{K,\vect{a}}(\Xi)\le1$ then, for $m\in\{m_j,\dots,m_{j+1}-1\}$\@, we have
\begin{align*}
j^{-m}\dnorm{\pi^\infty_m\scirc\Xi(x)}_{\ol{\metric}_m}\le&\;
a_m^m\epsilon_j\dnorm{\pi^\infty_m\scirc\Xi(x)}_{\ol{\metric}_m}\\
\le&\;a_0a_1\cdots a_m\epsilon_j
\dnorm{\pi^\infty_m\scirc\Xi(x)}_{\ol{\metric}_m}\le\epsilon_j
\end{align*}
for $x\in K$\@.  Thus, if $\Xi\in\sections[0]{\jet{\infty}{\man{E}}|K}$
satisfies $p_{K,\vect{a}}(\Xi)\le1$ then, for
$m\in\{m_j,\dots,m_{j+1}-1\}$\@, we have
$\pi^\infty_m\scirc\Xi\in\epsilon_j\nbhd{B}_j$\@.  Therefore,
$\Xi\in\cup_{j\in\integerp}\epsilon_j\nbhd{B}_j$\@, and this shows that, for
$\vect{a}$ as constructed above,
\begin{equation*}
\setdef{\Xi\in\sections[0]{\jet{\infty}{\man{E}}|K}}
{p_{K,\vect{a}}(\Xi)\le1}\subset\cup_{j\in\integerp}\epsilon_j\nbhd{B}_j,
\end{equation*}
giving the desired conclusion.
\end{proof}
\end{lemma}

The following attribute of the direct limit topology for $\sE(K)$ will also
be useful.
\begin{lemma}\label{lem:E(K)-regular}
The direct limit topology for\/ $\sE(K)$ is regular,~\ie~if\/
$\nbhd{B}\subset\sE(K)$ is von Neumann bounded, then there exists\/
$j\in\integerp$ such that\/ $\nbhd{B}$ is contained in and von Neumann
bounded in\/ $\sE_j(K)$\@.
\begin{proof}
Let $\nbhd{B}_j\subset\sE_j(K)$\@, $j\in\integerp$\@, be the closed unit ball
with respect to the norm topology.  We claim that $\nbhd{B}_j$ is closed in
the direct limit topology of $\sE(K)$\@.  To prove this, we shall prove that
$\nbhd{B}_j$ is closed in a topology that is weaker than the direct limit
topology.

The weaker topology we use is the topology induced by the topology of
pointwise convergence in $\sections[0]{\jet{\infty}{\man{E}}|K}$\@.  To be
precise, let $\sE'_j(K)$ be the vector space $\sE_j(K)$ with the topology
defined by the seminorms
\begin{equation*}
p_{x,j}(\Xi)=\sup\setdef{j^{-m}\dnorm{\pi^\infty_m\scirc
\Xi(x)}_{\ol{\metric}_m}}{m\in\integernn},\qquad x\in K.
\end{equation*}
Clearly the identity map from $\sE_j(K)$ to $\sE'_j(K)$ is continuous, and so
the topology of $\sE'_j(K)$ is weaker than the usual topology of $\sE(K)$\@.
Now let $\sE'(K)$ be the direct limit of the directed system
$\ifam{\sE'_j(K)}_{j\in\integerp}$\@.  Note that, algebraically,
$\sE'(K)=\sE(K)$\@, but the spaces have different topologies, the topology
for $\sE'(K)$ being weaker than that for $\sE(K)$\@.

We will show that $\nbhd{B}_j$ is closed in $\sE'(K)$\@.  Let $(I,\preceq)$
be a directed set and let $\ifam{\Xi_i}_{i\in I}$ be a convergent net in
$\nbhd{B}_j$ in the topology of $\sE'(K)$\@.  Thus we have a map
$\map{\Xi}{K}{\jet{\infty}{\man{E}}|K}$ such that, for each $x\in K$\@,
$\lim_{i\in I}\Xi_i(x)=\Xi(x)$\@.  If $\Xi\not\in\nbhd{B}_j$ then there
exists $x\in K$ such that
\begin{equation*}
\sup\setdef{j^{-m}
\dnorm{\pi^\infty_m\scirc\Xi(x)}_{\ol{\metric}_m}}{m\in\integernn}>1.
\end{equation*}
Let $\epsilon\in\realp$ be such that
\begin{equation*}
\sup\setdef{j^{-m}
\dnorm{\pi^\infty_m\scirc\Xi(x)}_{\ol{\metric}_m}}{m\in\integernn}>1+\epsilon
\end{equation*}
and let $i_0\in I$ be such that
\begin{equation*}
\sup\setdef{j^{-m}\dnorm{\pi^\infty\scirc\Xi_i(x)-
\pi^\infty_m\scirc\Xi(x)}_{\ol{\metric}_m}}{m\in\integernn}<\epsilon
\end{equation*}
for $i_0\preceq i$\@, this by pointwise convergence.  We thus have, for all
$i_0\preceq i$\@,
\begin{align*}
\epsilon<&\;\sup\setdef{j^{-m}
\dnorm{\pi^\infty_m\scirc\Xi(x)}_{\ol{\metric}_m}}{m\in\integernn}-
\sup\setdef{j^{-m}
\dnorm{\pi^\infty_m\scirc\Xi_i(x)}_{\ol{\metric}_m}}{m\in\integernn}\\
\le&\;\sup\setdef{j^{-m}\dnorm{\pi^\infty_m\scirc\Xi_i(x)-
\pi^\infty_m\scirc\Xi(x)}_{\ol{\metric}_m}}{m\in\integernn}<\epsilon,
\end{align*}
which contradiction gives the conclusion that $\Xi\in\nbhd{B}_j$\@.

Since $\nbhd{B}_j$ has been shown to be closed in $\sE(K)$\@, the lemma now
follows from~\cite[Corollary~7]{KDB:88}\@.
\end{proof}
\end{lemma}

\subsubsection{Seminorms for the topology of spaces of holomorphic germs}

Let us define seminorms $p^\omega_{K,\vect{a}}$\@, $K\subset\man{M}$ compact,
$\vect{a}\in\csd(\integernn;\realp)$\@, for
$\gsections[\hol,\real]{K}{\ol{\man{E}}}$ by
\begin{equation*}
p^\omega_{K,\vect{a}}([\ol{\xi}]_K)=\sup\setdef{a_0a_1\cdots a_m
\dnorm{j_m\ol{\xi}(x)}_{\ol{\metric}_m}}{x\in K,\ m\in\integernn}.
\end{equation*}
We can (and will) also think of $p^\omega_{K,\vect{a}}$ as being a seminorm
on $\sections[\omega]{\man{E}}$ defined by the same formula.

Let us prove that the seminorms $p^\omega_{K,\vect{a}}$\@, $K\subset\man{M}$
compact, $\vect{a}\in\csd(\integernn;\realp)$\@, can be used to define the
direct limit topology on $\gsections[\hol,\real]{K}{\ol{\man{E}}}$\@.
\begin{theorem}\label{the:Comega-seminorms}
Let\/ $\map{\pi}{\man{E}}{\man{M}}$ be a real analytic vector bundle and
let\/ $K\subset\man{M}$ be compact.  Then the family of seminorms\/
$p^\omega_{K,\vect{a}}$\@,\/ $\vect{a}\in\csd(\integernn;\realp)$\@, defines
a locally convex topology on\/ $\gsections[\hol,\real]{K}{\ol{\man{E}}}$
agreeing with the direct limit topology.
\begin{proof}
Let $K\subset\man{M}$ be compact and let
$\ifam{\ol{\nbhd{U}}_j}_{j\in\integerp}$ be a sequence of neighbourhoods of
$K$ in $\ol{\man{M}}$ such that
$\closure(\ol{\nbhd{U}}_{j+1})\subset\ol{\nbhd{U}}_j$\@, $j\in\integerp$\@,
and such that $K=\cap_{j\in\integerp}\ol{\nbhd{U}}_j$\@.  We have mappings
\begin{equation*}
\mapdef{r_{\ol{\nbhd{U}}_j,K}}
{\secbdd[\hol,\real]{\ol{\man{E}}|\ol{\nbhd{U}}_j}}
{\gsections[\hol,\real]{K}{\ol{\man{E}}}}{\ol{\xi}}{[\ol{\xi}]_K.}
\end{equation*}
The maps $r_{\ol{\nbhd{U}}_j,K}$ can be assumed to be injective without loss
of generality, by making sure that each open set $\ol{\nbhd{U}}_j$ consists
of disconnected neighbourhoods of the connected components of $K$\@.  Since
$\man{M}$ is Hausdorff and the connected components of $K$ are compact, this
can always be done by choosing the initial open set $\ol{\nbhd{U}}_1$
sufficiently small.  In this way,
$\secbdd[\hol,\real]{\ol{\man{E}}|\ol{\nbhd{U}}_j}$\@, $j\in\integerp$\@, are
regarded as subspaces of $\gsections[\hol,\real]{K}{\ol{\man{E}}}$\@.  It is
convenient to be able to do this.

We will work with the locally convex space $\sE(K)$ introduced in
Section~\ref{subsubsec:E(K)}\@, and define a mapping
$\map{L_K}{\gsections[\hol,\real]{K}{\ol{\man{E}}}}{\sE(K)}$ by
$L_K([\ol{\xi}]_K)=j_\infty\xi|K$\@.  Let us prove that this mapping is
well-defined,~\ie~show that, if
$[\ol{\xi}]_K\in\gsections[\hol,\real]{K}{\ol{\man{E}}}$\@, then
$L_K([\ol{\xi}]_K)\in\sE_j(K)$ for some $j\in\integerp$\@.  Let
$\ol{\nbhd{U}}$ be a neighbourhood of $K$ in $\ol{\man{M}}$ on which the
section $\ol{\xi}$ is defined, holomorphic, and bounded.  Then
$\xi|(\man{M}\cap\ol{\nbhd{U}})$ is real analytic and so, by
Lemma~\ref{lem:rabound}\@, there exist $C,r\in\realp$ such that
\begin{equation*}
\dnorm{j_m\xi(x)}_{\ol{\metric}_m}\le Cr^{-m},\qquad x\in K,\ m\in\integernn.
\end{equation*}
If $j>r^{-1}$ it immediately follows that
\begin{equation*}
\sup\setdef{j^{-m}\dnorm{j_m\xi(x)}_{\ol{\metric}_m}}
{x\in K,\ m\in\integernn}<\infty,
\end{equation*}
\ie~$L_K([\ol{\xi}]_K)\in\sE_j(K)$\@.

The following lemma records the essential feature of $L_K$\@.
\begin{prooflemma}
The mapping\/ $L_K$ is a continuous, injective, open mapping, and so an
homeomorphism onto its image.
\begin{subproof}
To show that $L_K$ is continuous, it suffices to show that
$L_K|\secbdd[\hol,\real]{\ol{\man{E}}|\ol{\nbhd{U}}_j}$ is continuous for
each $j\in\integernn$\@.  We will show this by showing that, for each
$j\in\integerp$\@, there exists $j'\in\integerp$ such that
$L_K(\secbdd[\hol]{\ol{\man{E}}|\ol{\nbhd{U}}_j})\subset\sE_{j'}(K)$ and such
that $L_K$ is continuous as a map from
$\secbdd[\hol]{\ol{\man{E}}|\ol{\nbhd{U}}_j}$ to $\sE_{j'}(K)$\@.  Since
$\sE_{j'}(K)$ is continuously included in $\sE(K)$\@, this will give the
continuity of $L_K$\@.  First let us show that
$L_K(\secbdd[\hol]{\ol{\man{E}}|\ol{\nbhd{U}}_j})\subset\sE_{j'}(K)$ for some
$j'\in\integerp$\@.  By Proposition~\ref{prop:Cauchyest}\@, there exist
$C,r\in\realp$ such that
\begin{equation*}
\dnorm{j_m\xi(x)}_{\ol{\metric}_m}\le
Cr^{-m}p^\hol_{\ol{\nbhd{U}}_j,\infty}(\ol{\xi})
\end{equation*}
for every $m\in\integernn$ and
$\ol{\xi}\in\secbdd[\hol]{\ol{\man{E}}|\ol{\nbhd{U}}_j}$\@.  Taking
$j'\in\integerp$ such that $j'\ge r^{-1}$ we have
$L_K(\secbdd[\hol]{\ol{\man{E}}|\ol{\nbhd{U}}_j})\subset\sE_{j'}(K)$\@, as
claimed.  To show that $L_K$ is continuous as a map from
$\secbdd[\hol]{\ol{\man{E}}|\ol{\nbhd{U}}_j}$ to $\sE_{j'}(K)$\@, let
$\ifam{[\xi_k]_K}_{k\in\integerp}$ be a sequence in
$\secbdd[\hol]{\ol{\man{E}}|\ol{\nbhd{U}}_j}$ converging to zero.  We then
have
\begin{equation*}
\lim_{k\to\infty}\sup\setdef{(j')^{-m}
\dnorm{j_m\xi_k(x)}_{\ol{\metric}}}{x\in K,\ m\in\integernn}\le
\lim_{k\to\infty}C\sup\setdef{\dnorm{\ol{\xi}_k(z)}_{\ol{\metric}}}
{z\in\ol{\nbhd{U}}_j}=0,
\end{equation*}
giving the desired continuity.

Since germs of holomorphic sections are uniquely determined by their infinite
jets, injectivity of $L_K$ follows.

We claim that, if $\nbhd{B}\subset\sE(K)$ is von Neumann bounded, then
$L_K^{-1}(\nbhd{B})$ is also von Neumann bounded.  By
Lemma~\ref{lem:E(K)-regular}\@, if $\nbhd{B}\subset\sE(K)$ is bounded, then
$\nbhd{B}$ is contained and bounded in $\sE_j(K)$ for some $j\in\integerp$\@.
Therefore, there exists $C\in\realp$ such that, if
$L_K([\ol{\xi}]_K)\subset\nbhd{B}$\@, then
\begin{equation*}
\dnorm{j_m\xi(x)}_{\ol{\metric}_m}\le Cj^{-m},\qquad x\in K,\ m\in\integernn.
\end{equation*}
Let $x\in K$ and let $(\nbhd{V}_x,\psi_x)$ be a vector bundle chart for
$\man{E}$ about $x$ with corresponding chart $(\nbhd{U}_x,\phi_x)$ for
$\man{M}$\@.  Suppose the fibre dimension of $\man{E}$ over $\nbhd{U}_x$ is
$k$ and that $\phi_x$ takes values in $\real^n$\@.  Let
$\nbhd{U}'_x\subset\nbhd{U}_x$ be a relatively compact neighbourhood of $x$
such that $\closure(\nbhd{U}'_x)\subset\nbhd{U}_x$\@.  Denote
$K_x=K\cap\closure(\nbhd{U}'_x)$\@.  By Lemma~\ref{lem:pissy-estimate}\@,
there exist $C_x,r_x\in\realp$ such that, if
$L_K([\ol{\xi}]_K)\subset\nbhd{B}$\@, then
\begin{equation*}
\snorm{\linder[I]{\xi^a}(\vect{x})}\le C_xI!r_x^{-\snorm{I}},
\qquad\vect{x}\in\phi_x(K_x),\ I\in\integernn^n,\ a\in\{1,\dots,k\},
\end{equation*}
where $\vect{\xi}$ is the local representative of $\xi$\@.  Note that this
implies the following for each $[\ol{\xi}]_K$ such that
$L_K([\ol{\xi}]_K)\subset\nbhd{B}$ and for each $a\in\{1,\dots,k\}$\@:
\begin{compactenum}
\item $\ol{\xi}\null^a$ admits a convergent power series expansion to an
holomorphic function on the polydisk $\odisk[]{\vect{\sigma}_x}{\phi_x(x)}$
for $\sigma_x<r_x$\@;
\item on the polydisk $\odisk[]{\vect{\sigma}_x}{\phi_x(x)}$\@,
$\ol{\xi}\null^a$ satisfies
$\snorm{\ol{\xi}\null^a}\le(\frac{1}{1-\sigma_x})^n$\@.
\end{compactenum}
It follows that, if $L_K([\ol{\xi}]_K)\in\nbhd{B}$\@, then $\ol{\xi}$ has a
bounded holomorphic extension in some coordinate polydisk around each $x\in
K$\@.  By a standard compactness argument and since
$\cap_{j\in\integerp}\ol{\nbhd{U}}_j=K$\@, there exists $j'\in\integerp$ such
that $\ol{\xi}\in\secbdd[\hol,\real]{\ol{\man{E}}|\ol{\nbhd{U}}_{j'}}$ for
each $[\ol{\xi}]_K$ such that $L_K([\ol{\xi}]_K)\in\nbhd{B}$\@, and that the
set of such sections of $\ol{\man{E}}|\ol{\nbhd{U}}_{j'}$ is von Neumann
bounded,~\ie~norm bounded.  Thus $L_K^{-1}(\nbhd{B})$ is von Neumann bounded,
as claimed.

Note also that $\sE(K)$ is a DF-space since Banach spaces are
DF-spaces~\cite[Corollary~12.4.4]{HJ:81} and countable direct limits of
DF-spaces are DF-spaces~\cite[Theorem~12.4.8]{HJ:81}\@.  Therefore, by the
open mapping lemma from \S2\ of \citetalias{AB:71}\@, the result follows.
\end{subproof}
\end{prooflemma}

From the lemma, it follows that the direct limit topology of
$\gsections[\hol,\real]{K}{\ol{\man{E}}}$ agrees with that induced by its
image in $\sE(K)$\@.  Since the seminorms $p_{K,\vect{a}}$\@,
$\vect{a}\in\csd(\integernn;\realp)$\@, define the locally convex topology of
$\sE(K)$ by Lemma~\ref{lem:EKseminorms}\@, it follows that the seminorms
$p^\omega_{K,\vect{a}}$\@, $\vect{a}\in\csd(\integernn;\realp)$\@, define the
direct limit topology of $\gsections[\hol,\real]{K}{\ol{\man{E}}}$\@.
\end{proof}
\end{theorem}

The problem of providing seminorms for the direct limit topology of
$\gsections[\hol,\real]{K}{\ol{\man{E}}}$ is a nontrivial one, so let us
provide a little history for what led to the preceding theorem.  First of
all, the first concrete characterisation of seminorms for germs of
holomorphic functions about compact subsets of $\complex^n$ comes
in~\cite{JM:84}\@.  \citeauthor{JM:84} provides seminorms having two parts,
one very much resembling the seminorms we use, and another part that is more
complicated.  These seminorms specialise to the case where the compact set
lies in $\real^n\subset\complex^n$\@, and the first mention of this we have
seen in the research literature is in the notes of \citet{PD:10}\@.  The
first full proof that the seminorms analogous to those we define are, in
fact, the seminorms for the space of real analytic functions on open subsets
of $\real^n$ appears in the recent note of \citet{DV:13}\@.  Our presentation
is an adaptation, not quite trivial as it turns out, of \citeauthor{DV:13}'s
constructions.  One of the principal difficulties is
Lemma~\ref{lem:pissy-estimate} which is essential in showing that our jet
bundle fibre metrics $\dnorm{\cdot}_{\ol{\metric}_m}$ are suitable for
defining the seminorms for the real analytic topology.  Note that one cannot
use arbitrary fibre metrics, since one needs to have the behaviour of these
metrics be regulated to the real analytic topology as the order of jets goes
to infinity.  Because our fibre metrics are constructed by differentiating
objects defined at low order,~\ie~the connections $\nabla$ and $\nabla^0$\@,
we can ensure that the fibre metrics are compatible with real analytic growth
conditions on derivatives.

\subsubsection{An inverse limit topology for the space of real analytic
sections}

In the preceding three sections we provided three topologies for the space
$\gsections[\hol,\real]{K}{\ol{\man{E}}}$ of holomorphic sections about a
compact subset $K$ of a real analytic manifold:~(1)~the ``standard'' direct
limit topology;~(2)~the topology induced by the direct limit topology on
$\sE(K)$\@;~(3)~the topology defined by the seminorms
$p^\omega_{K,\vect{a}}$\@, $K\subset\man{M}$ compact,
$\vect{a}\in\csd(\integernn;\realp)$\@.  We showed in
Lemma~\ref{lem:EKseminorms} and Theorem~\ref{the:Comega-seminorms} that these
three topologies agree.  Now we shall use these constructions to easily
arrive at~(1)~a topology on $\sections[\omega]{\man{E}}$ induced by the
locally convex topologies on the spaces
$\gsections[\hol,\real]{K}{\ol{\man{E}}}$\@, $K\subset\man{M}$ compact,
and~(2)~seminorms for the topology of $\sections[\omega]{\man{E}}$\@.

For a compact set $K\subset\man{M}$ we have an inclusion
$\map{i_K}{\sections[\omega]{\man{E}}}
{\gsections[\hol,\real]{K}{\ol{\man{E}}}}$ defined as follows.  If
$\xi\in\sections[\omega]{\man{E}}$\@, then $\xi$ admits an holomorphic
extension $\ol{\xi}$ defined on a neighbourhood
$\ol{\nbhd{U}}\subset\ol{\man{M}}$ of
$\man{M}$~\cite[Lemma~5.40]{KC/YE:12}\@.  Since $\ol{\nbhd{U}}\in\sN_K$ we
define $i_K(\xi)=[\ol{\xi}]_K$\@.  Now we have a compact exhaustion
$\ifam{K_j}_{j\in\integerp}$ of $\man{M}$\@.  Since
$\sN_{K_{j+1}}\subset\sN_{K_j}$ we have a projection
\begin{equation*}
\mapdef{\pi_j}{\gsections[\hol,\real]{K_{j+1}}{\ol{\man{E}}}}
{\gsections[\hol,\real]{K_j}{\ol{\man{E}}}}
{[\ol{\xi}]_{K_{j+1}}}{[\ol{\xi}]_{K_j}.}
\end{equation*}
One can check that, as $\real$-vector spaces, the inverse limit of the
inverse family
$\ifam{\gsections[\hol,\real]{K_j}{\ol{\man{E}}}}_{j\in\integerp}$ is
isomorphic to $\gsections[\hol,\real]{\man{M}}{\ol{\man{E}}}$\@, the
isomorphism being given explicitly by the inclusions
\begin{equation*}
\mapdef{i_j}{\gsections[\hol,\real]{\man{M}}{\ol{\man{E}}}}
{\gsections[\hol,\real]{K_j}{\ol{\man{E}}}}
{[\ol{\xi}]_{\man{M}}}{[\ol{\xi}]_{K_j}.}
\end{equation*}
Keeping in mind Lemma~\ref{lem:section<->germ}\@, we then have the inverse
limit topology on $\sections[\omega]{\man{E}}$ induced by the mappings
$i_j$\@, $j\in\integerp$\@.  The topology so defined we call the
\defn{inverse $\C^\omega$-topology} for $\sections[\omega]{\man{E}}$\@.

It is now a difficult theorem of \citet[Theorem~1.2(a)]{AM:66} that the
direct $\C^\omega$-topology of Section~\ref{subsubsec:Comega-direct} agrees
with the inverse $\C^\omega$-topology.  Therefore, we call the resulting
topology the \defn{$\C^\omega$-topology}\@.  It is clear from
Theorem~\ref{the:Comega-seminorms} and the preceding inverse limit
construction that the seminorms $p^\omega_{K,\vect{a}}$\@, $K\subset\man{M}$
compact, $\vect{a}\in\csd(\integernn;\realp)$\@, define the
$\C^\omega$-topology.

\subsection{Properties of the $\C^\omega$-topology}\label{subsec:Comega-props}

To say some relevant things about the $\C^\omega$-topology, let us first
consider the direct limit topology for
$\gsections[\hol,\real]{K}{\ol{\man{E}}}$\@, $K\subset\man{M}$ compact, as
this is an important building block for the $\C^\omega$-topology.  First, we
recall that a \defn{strict direct limit} of locally convex spaces consists of
a sequence $\ifam{\alg{V}_j}_{j\in\integerp}$ of locally convex spaces that
are subspaces of some vector space $\alg{V}$\@, and which have the nesting
property $\alg{V}_j\subset\alg{V}_{j+1}$\@, $j\in\integerp$\@.  In defining
the direct limit topology for $\gsections[\hol,\real]{K}{\ol{\man{E}}}$ we
defined it as a strict direct limit of Banach spaces.  Moreover, the
restriction mappings from $\secbdd[\hol,\real]{\ol{\man{E}}|\ol{\nbhd{U}}_j}$
to $\secbdd[\hol,\real]{\ol{\man{E}}|\ol{\nbhd{U}}_{j+1}}$ can be shown to be
compact~\cite[Theorem~8.4]{AK/PWM:97}\@.  Direct limits such as these are
known as ``Silva spaces'' or ``DFS spaces.''  Silva spaces have some nice
properties, and these provide some of the following attributes for the direct
limit topology for $\gsections[\hol,\real]{K}{\ol{\man{E}}}$\@.
\begin{compactenum}[$\mathscr{G}^{\hol,\real}$-1.]
\item It is Hausdorff:~\cite[Theorem~12.1.3]{LN/EB:10}\@.
\item It is complete:~\cite[Theorem~12.1.10]{LN/EB:10}\@.
\item It is not metrisable:~\cite[Theorem~12.1.8]{LN/EB:10}\@.
\item It is regular:~\cite[Theorem~8.4]{AK/PWM:97}\@.  This means that every
von Neumann bounded subset of $\gsections[\hol,\real]{K}{\ol{\man{E}}}$ is
contained and von Neumann bounded in
$\sections[\hol,\real]{\ol{\man{E}}|\ol{\nbhd{U}}_j}$ for some
$j\in\integerp$\@.
\item \label{enum:Ghol-reflexive} It is
reflexive:~\cite[Theorem~8.4]{AK/PWM:97}\@.
\item \label{enum:Ghol-dual} Its strong dual is a nuclear Fr\'echet
space:~\cite[Theorem~8.4]{AK/PWM:97}\@.  Combined with reflexivity, this
means that $\gsections[\hol,\real]{K}{\ol{\man{E}}}$ is the strong dual of a
nuclear Fr\'echet space.
\item It is nuclear:~\cite[Corollary~III.7.4]{HHS/MPW:99}\@.
\item \label{enum:Ghol-suslin} It is Suslin: This follows
from~\cite[Th\'eor\`eme~I.5.1(b)]{XF:67} since
$\gsections[\hol,\real]{K}{\ol{\man{E}}}$ is a strict direct limit of
separable Fr\'echet spaces.
\end{compactenum}
These attributes for the spaces $\gsections[\hol,\real]{K}{\ol{\man{E}}}$
lead, more or less, to the following attributes of
$\sections[\omega]{\man{E}}$\@.
\begin{compactenum}[$\C^\omega$-1.]
\item It is Hausdorff:~It is a union of Hausdorff topologies.
\item \label{enum:Comega-complete} It is complete:~\cite[Corollary
to~Proposition~2.11.3]{JH:66}\@.
\item It is not metrisable:~It is a union of non-metrisable topologies.
\item \label{enum:Comega-separable} It is separable:~\cite[Theorem~16]{PD:10}\@.
\item \label{enum:Comega-nuclear} It is
nuclear:~\cite[Corollary~III.7.4]{HHS/MPW:99}\@.
\item \label{enum:Comega-suslin} It is Suslin: Here we note that a countable
direct product of Suslin spaces is Suslin~\cite[Lemma~6.6.5(iii)]{VIB:07b}\@.
Next we note that the inverse limit is a closed subspace of the direct
product~\cite[Proposition~V.19]{APR/WR:80}\@.  Next, closed subspaces of
Suslin spaces are Suslin spaces~\cite[Lemma~6.6.5(ii)]{VIB:07b}\@.
Therefore, since $\sections[\omega]{\man{E}}$ is the inverse limit of the
Suslin spaces $\gsections[\hol,\real]{K_j}{\ol{\man{E}}}$\@,
$j\in\integerp$\@, we conclude that $\sections[\omega]{\man{E}}$ is Suslin.
\end{compactenum}

As we have seen with the $\CO^\infty$- and $\CO^\hol$-topologies for
$\sections[\infty]{\man{E}}$ and $\sections[\hol]{\man{E}}$\@, nuclearity of
the $\C^\omega$-topology implies that compact subsets of
$\sections[\omega]{\man{E}}$ are exactly those that are closed and von
Neumann bounded.  For von Neumann boundedness, we have the following
characterisation.
\begin{lemma}\label{lem:Comegabdd}
A subset\/ $\nbhd{B}\subset\sections[\omega]{\man{E}}$ is bounded in the von
Neumann bornology if and only if the following property holds: for any
compact set\/ $K\subset\man{M}$ and any\/
$\vect{a}\in\csd(\integernn;\realp)$\@, there exists\/ $C\in\realp$ such
that\/ $p^\omega_{K,\vect{a}}(\xi)\le C$ for every\/ $\xi\in\nbhd{B}$\@.
\end{lemma}

\subsection{The weak-$\sL$ topology for real analytic vector fields}

As in the finitely differentiable, Lipschitz, smooth, and holomorphic cases,
the above constructions for general vector bundles can be applied to the
tangent bundle and the trivial vector bundle $\man{M}\times\real$ to give the
\defn{$\C^\omega$-topology} on the space $\sections[\omega]{\tb{\man{M}}}$ of
real analytic vector fields and the space $\func[\omega]{\man{M}}$ of real
analytic functions.  As we have already done in these other cases, we wish to
provide a weak characterisation of the $\C^\omega$-topology for
$\sections[\omega]{\tb{\man{M}}}$\@.  First of all, if
$X\in\sections[\omega]{\tb{\man{M}}}$\@, then $f\mapsto\lieder{X}{f}$ is a
derivation of $\func[\omega]{\man{M}}$\@.  As we have seen, in the
holomorphic case this does not generally establish a correspondence between
vector fields and derivations, but it does for Stein manifolds.  In the real
analytic case, \citet{JG:81} shows that the map $X\mapsto\lieder{X}{}$ is
indeed an isomorphism of the $\real$-vector spaces of real analytic vector
fields and derivations of real analytic functions.  Thus the pursuit of a
weak description of the $\C^\omega$-topology for vector fields does not seem
to be out of line.

The definition of the weak-$\sL$ topology proceeds much as in the smooth and
holomorphic cases.
\begin{definition}
For a real analytic manifold $\man{M}$\@, the \defn{weak-$\sL$ topology} for
$\sections[\omega]{\tb{\man{M}}}$ is the weakest topology for which the map
$X\mapsto\lieder{X}{f}$ is continuous for every
$f\in\func[\omega]{\man{M}}$\@, if $\func[\omega]{\man{M}}$ has the
$\C^\omega$-topology.\oprocend
\end{definition}

We now have the following result.
\begin{theorem}\label{the:Comega-weak}
For a real analytic manifold\/ $\man{M}$\@, the following topologies for\/
$\sections[\omega]{\tb{\man{M}}}$ agree:
\begin{compactenum}[(i)]
\item \label{pl:Comega-weak1} the\/ $\C^\omega$-topology;
\item \label{pl:Comega-weak2} the weak-$\sL$ topology.
\end{compactenum}
\begin{proof}
\eqref{pl:Comega-weak1}$\subset$\eqref{pl:Comega-weak2} As we argued in the
corresponding part of the proof of Theorem~\ref{the:COinfty-weak}\@, it
suffices to show that, for $K\subset\man{M}$ compact and for
$\vect{a}\in\csd(\integernn;\realp)$\@, there exist compact sets
$K_1,\dots,K_r\subset\man{M}$\@,
$\vect{a}_1,\dots,\vect{a}_r\in\csd(\integernn;\realp)$\@,
$f^1,\dots,f^r\in\func[\omega]{\man{M}}$\@, and $C_1,\dots,C_r\in\realp$ such
that
\begin{equation*}
p^\omega_{K,\vect{a}}(X)\le C_1p^\omega_{K_1,\vect{a}_1}(\lieder{X}{f^1})+\dots+
C_rp^\omega_{K_r,\vect{a}_r}(\lieder{X}{f^r}),\qquad
X\in\sections[\omega]{\tb{\man{M}}}.
\end{equation*}

We begin with a simple technical lemma.
\begin{prooflemma}
For each\/ $x\in\man{M}$ there exist\/
$f^1,\dots,f^n\in\func[\omega]{\man{M}}$ such that\/
$\ifam{\d{f^1}(x),\dots,\d{f^n}(x)}$ is a basis for\/ $\ctb[x]{\man{M}}$\@.
\begin{subproof}
We are supposing, of course, that the connected component of $\man{M}$
containing $x$ has dimension $n$\@.  There are many ways to prove this lemma,
including applying Cartan's Theorem~A to the sheaf of real analytic functions
on $\man{M}$\@.  We shall prove the lemma by embedding $\man{M}$ in $\real^N$
by the embedding theorem of \citet{HG:58}\@.  Thus we have a proper real
analytic embedding $\map{\iota_{\man{M}}}{\man{M}}{\real^N}$\@.  Let
$g^1,\dots,g^N\in\func[\omega]{\real^N}$ be the coordinate functions.  Then
we have a surjective linear map
\begin{equation*}
\mapdef{\sigma_x}{\real^N}{\ctb[x]{\man{M}}}{(c_1,\dots,c_N)}
{\sum_{j=1}^Nc_j\d{(\iota_{\man{M}}^*g^j)}(x).}
\end{equation*}
Let $\vect{c}^1,\dots,\vect{c}^n\in\real^N$ be a basis for a complement of
$\ker(\sigma_x)$\@.  Then the functions
\begin{equation*}
f^j=\sum_{k=1}^Nc_k^jg^k
\end{equation*}
have the desired property.
\end{subproof}
\end{prooflemma}

We assume that $\man{M}$ has a well-defined dimension $n$\@.  This assumption
can easily be relaxed.  We use the notation
\begin{equation*}
p^{\prime\omega}_{K,\vect{a}}(f)=\sup\Bigsetdef{\frac{a_0a_1\cdots a_{\snorm{I}}}{I!}
\snorm{\linder[I]{f}(\vect{x})}}{\vect{x}\in K,\ I\in\integernn^n}
\end{equation*}
for a function $f\in\func[\omega]{\nbhd{U}}$ defined on an open subset of
$\real^n$ and with $K\subset\nbhd{U}$ compact.  We shall also use this local
coordinate notation for seminorms of local representatives of vector fields.
Let $K\subset\man{M}$ be compact and let
$\vect{a}\in\csd(\integernn;\realp)$\@.  Let $x\in K$ and let
$(\nbhd{U}_x,\phi_x)$ be a chart for $\man{M}$ about $x$ with the property
that the coordinate functions $x^j$\@, $j\in\{1,\dots,n\}$\@, are
restrictions to $\nbhd{U}_x$ of globally defined real analytic functions
$f^j_x$\@, $j\in\{1,\dots,n\}$\@, on $\man{M}$\@.  This is possible by the
lemma above.  Let $\map{\vect{X}}{\phi_x(\nbhd{U}_x)}{\real^n}$ be the local
representative of $X\in\sections[\omega]{\man{M}}$\@.  Then, in a
neighbourhood of the closure of a relatively compact neighbourhood
$\nbhd{V}_x\subset\nbhd{U}_x$ of $x$\@, we have $\lieder{X}{f^j_x}=X^j$\@,
the $j$th component of $X$\@.  By Lemma~\ref{lem:pissy-estimate}\@, there
exist $C_x,\sigma_x\in\realp$ such that
\begin{equation*}
\dnorm{j_mX(y)}_{\ol{\metric}_m}\le
C_x\sigma_x^{-m}\sup\Bigsetdef{\frac{1}{I!}
\snorm{\linder[I]{X^j}(\phi_x(y))}}{\snorm{I}\le m,\ j\in\{1,\dots,n\}}
\end{equation*}
for $m\in\integernn$ and $y\in\closure(\nbhd{V}_x)$\@.  By equivalence of the
$\ell^1$ and $\ell^\infty$-norms for $\real^n$\@, there exists $C\in\realp$
such that
\begin{multline*}
\sup\Bigsetdef{\frac{1}{I!}\snorm{\linder[I]{X^j}(\phi_x(y))}}
{\snorm{I}\le m,\ j\in\{1,\dots,n\}}\\
\le C\sum_{j=1}^n\sup\Bigsetdef{\frac{1}{I!}
\snorm{\linder[I]{(\lieder{X}{f^j_x})}(\phi_x(y))}}{\snorm{I}\le m}
\end{multline*}
for $m\in\integernn$ and $y\in\closure(\nbhd{V}_x)$\@.  Another application
of Lemma~\ref{lem:pissy-estimate} gives $B_x,r_x\in\realp$ such that
\begin{equation*}
\sup\Bigsetdef{\frac{1}{I!}\snorm{\linder[I]{(\lieder{X}{f_x^j})}(\phi_x(y))}}
{\snorm{I}\le m,\ j\in\{1,\dots,n\}}\le
B_xr_x^{-m}\dnorm{j_m(\lieder{X}{f^j_x})(y)}
\end{equation*}
for $m\in\integernn$\@, $j\in\{1,\dots,n\}$\@, and
$y\in\closure(\nbhd{V}_x)$\@.  Combining the preceding three estimates and
renaming constants gives
\begin{equation*}
\dnorm{j_mX(y)}_{\ol{\metric}_m}\le
\sum_{j=1}^nC_x\sigma_x^{-m}
\dnorm{j_m(\lieder{X}{f^j_x}(\phi_x(y)))}_{\ol{\metric}_m}
\end{equation*}
for $m\in\integernn$ and $y\in\closure(\nbhd{V}_x)$\@.  Define
\begin{equation*}
\vect{b}_x=\ifam{b_m}_{m\in\integernn}\in\csd(\integernn;\realp)
\end{equation*}
by $b_0=C_xa_0$ and $b_m=\sigma_x^{-1}a_m$\@, $m\in\integerp$\@.  Therefore,
\begin{equation*}
a_0a_1\cdots a_m\dnorm{j_mX(y)}_{\ol{\metric}_m}\le
\sum_{j=1}^nb_0b_1\cdots b_m
\dnorm{j_m(\lieder{X}{f^j_x}(\phi_x(y)))}_{\ol{\metric}_m}
\end{equation*}
for $m\in\integernn$ and $y\in\closure(\nbhd{V}_x)$\@.  Supping over
$y\in\closure(\nbhd{V}_x)$ and $m\in\integernn$ on the right gives
\begin{equation*}
a_0a_1\cdots a_m\dnorm{j_mX(y)}_{\ol{\metric}_m}\le
C_x\sum_{j=1}^np^\omega_{\closure(\nbhd{V}_x),\vect{b}_x}(\lieder{X}{f^j_x}),
\qquad m\in\integernn,\ y\in\closure(\nbhd{V}_x).
\end{equation*}
Let $x_1,\dots,x_k\in K$ be such that $K\subset\cup_{j=1}^k\nbhd{V}_{x_j}$\@,
let $f^1,\dots,f^{kn}$ be the list of functions
\begin{equation*}
f^1_{x_1},\dots,f^n_{x_1},\dots,f^1_{x_k},\dots,f^n_{x_k},
\end{equation*}
and let $\vect{a}_1,\dots,\vect{a}_{kn}\in\csd(\integernn;\realp)$ be the
list of sequences
\begin{equation*}
\underbrace{\vect{b}_{x_1},\dots,\vect{b}_{x_1}}_{n\ \textrm{times}},\dots,
\underbrace{\vect{b}_{x_k},\dots,\vect{b}_{x_k}}_{n\ \textrm{times}}.
\end{equation*}
If $x\in K$\@, then $x\in\nbhd{V}_{x_j}$ for some $j\in\{1,\dots,k\}$ and so
\begin{equation*}
a_0a_1\cdots a_m\dnorm{j_mX(x)}_{\ol{\metric}_m}\le
\sum_{j=1}^{kn}p^\omega_{K,\vect{b}_j}(\lieder{X}{f^j}),
\end{equation*}
and this part of the lemma follows upon taking the supremum over $x\in K$ and $m\in\integernn$\@.

\eqref{pl:Comega-weak2}$\subset$\eqref{pl:Comega-weak1} Here, as in the proof
of the corresponding part of Theorem~\ref{the:COinfty-weak}\@, it suffices to
show that, for every $f\in\func[\omega]{\man{M}}$\@, the map $\sL_f\colon
X\mapsto\lieder{X}{f}$ is continuous from $\sections[\omega]{\tb{\man{M}}}$
with the $\C^\omega$-topology to $\func[\omega]{\man{M}}$ with the
$\C^\omega$-topology.

We shall use the direct $\C^\omega$-topology to show this.  Thus we work with
an holomorphic manifold $\ol{\man{M}}$ that is a complexification of
$\man{M}$\@, as described in Section~\ref{subsubsec:complexify}\@.  We recall
that $\sN_{\man{M}}$ denotes the directed set of neighbourhoods of $\man{M}$
in $\ol{\man{M}}$\@, and that the set $\sS_{\man{M}}$ of Stein neighbourhoods
is cofinal in $\sN_{\man{M}}$\@.  As we saw in
Section~\ref{subsubsec:Comega-direct}\@, for
$\ol{\nbhd{U}}\in\sN_{\man{M}}$\@, we have mappings
\begin{equation*}
\mapdef{r_{\ol{\nbhd{U}},\man{M}}}
{\sections[\hol,\real]{\tb{\ol{\nbhd{U}}}}}
{\sections[\omega]{\tb{\man{M}}}}{\ol{X}}{\ol{X}|\man{M}}
\end{equation*}
and
\begin{equation*}
\mapdef{r_{\ol{\nbhd{U}},\man{M}}}
{\func[\hol,\real]{\ol{\nbhd{U}}}}
{\func[\omega]{\man{M}}}{\ol{f}}{\ol{f}|\man{M},}
\end{equation*}
making an abuse of notation by using $r_{\ol{\nbhd{U}},\man{M}}$ for two
different things, noting that context will make it clear which we mean.  For
$K\subset\man{M}$ compact, we also have the mapping
\begin{equation*}
\mapdef{i_{\man{M},K}}{\func[\omega]{\man{M}}}
{\gfunc[\hol,\real]{K}{\ol{\man{M}}}}{f}{[\ol{f}]_K,}
\end{equation*}
The $\C^\omega$-topology is the final topology induced by the mappings
$r_{\ol{\nbhd{U}},\man{M}}$\@.  As such,
by~\cite[Proposition~2.12.1]{JH:66}\@, the map $\sL_f$ is continuous if and
only if $\sL_f\scirc r_{\ol{\nbhd{U}},\man{M}}$ for every
$\ol{\nbhd{U}}\in\sN_{\man{M}}$\@.  Thus let
$\ol{\nbhd{U}}\in\sN_{\man{M}}$\@.  To show that $\sL_f\scirc
r_{\ol{\nbhd{U}},\man{M}}$ is continuous, it suffices by \cite[\S2.11]{JH:66}
to show that $i_{\man{M},K}\scirc\sL_f\scirc r_{\ol{\nbhd{U}},\man{M}}$ is
continuous for every compact $K\subset\man{M}$\@.  Next, there is
$\ol{\nbhd{U}}\supset\ol{\nbhd{S}}\in\sS_{\man{M}}$ so that $f$ admits an
holomorphic extension $\ol{f}$ to $\ol{\nbhd{S}}$\@.  The following diagram
shows how this all fits together.
\begin{equation*}
\xymatrix{{\sections[\hol,\real]{\tb{\ol{\nbhd{U}}}}}
\ar@{-->}[rrd]\ar[r]^{r_{\ol{\nbhd{U}},\ol{\nbhd{S}}}}&
{\sections[\hol,\real]{\tb{\ol{\nbhd{S}}}}}
\ar[r]^{r_{\ol{\nbhd{S}},\man{M}}}\ar[d]_(0.3){\lieder{\bar{f}}{}}&
{\sections[\omega]{\tb{\man{M}}}}\ar@{-->}[d]^{\sL_f}&\\
&{\func[\hol,\real]{\ol{\nbhd{S}}}}\ar[r]_{r_{\ol{\nbhd{S}},\man{M}}}&
{\func[\omega]{\man{M}}}\ar[r]_{i_{\man{M},K}}&
{\gfunc[\hol,\real]{K}{\ol{\man{M}}}}}
\end{equation*}
The dashed arrows signify maps whose continuity is \emph{a priori} unknown to
us.  The diagonal dashed arrow is the one whose continuity we must verify to
ascertain the continuity of the vertical dashed arrow.  It is a simple matter
of checking definitions to see that the diagram commutes.  By
Theorem~\ref{the:COhol-weak}\@, we have that $\map{\lieder{\ol{f}}{}}
{\sections[\hol,\real]{\tb{\ol{\nbhd{S}}}}}
{\func[\hol,\real]{\ol{\nbhd{S}}}}$ is continuous (keeping
Remark~\ref{rem:real-closed} in mind).  We deduce that, since
\begin{equation*}
i_{\man{M},K}\scirc\sL_f\scirc r_{\ol{\nbhd{U}},\man{M}}=
i_{\man{M},K}\scirc r_{\ol{\nbhd{S}},\man{M}}\scirc\lieder{\ol{f}}
\scirc r_{\ol{\nbhd{U}},\ol{\nbhd{S}}},
\end{equation*}
$i_{\man{M},K}\scirc\sL_f\scirc r_{\ol{\nbhd{U}},\man{M}}$ is
continuous for every $\ol{\nbhd{U}}\in\sN_{\man{M}}$ and for every compact
$K\subset\man{M}$\@, as desired.
\end{proof}
\end{theorem}

As in the smooth and holomorphic cases, we can prove the equivalence of
various topological notions between the weak-$\sL$ and usual topologies.
\begin{corollary}\label{cor:Comegaweak}
Let\/ $\man{M}$ be a real analytic manifold, let\/ $(\ts{X},\sO)$ be a
topological space, let\/ $(\ts{T},\sM)$ be a measurable space, and let\/
$\map{\mu}{\sM}{\erealnn}$ be a finite measure.  The following statements hold:
\begin{compactenum}[(i)]
\item \label{pl:Comegaweakbdd} a subset\/
$\nbhd{B}\subset\sections[\omega]{\tb{\man{M}}}$ is bounded in the von Neumann
bornology if and only if it is weak-$\sL$ bounded in the von Neumann
bornology;
\item \label{pl:Comegaweakcont} a map\/
$\map{\Phi}{\ts{X}}{\sections[\omega]{\tb{\man{M}}}}$ is continuous if and
only if it is weak-$\sL$ continuous;
\item \label{pl:Comegaweakmeas} a map\/
$\map{\Psi}{\ts{T}}{\sections[\omega]{\tb{\man{M}}}}$ is measurable if and
only if it is weak-$\sL$ measurable;
\item \label{pl:Comegaweakint} a map\/
$\map{\Psi}{\ts{T}}{\sections[\omega]{\tb{\man{M}}}}$ is Bochner integrable
if and only if it is weak-$\sL$ Bochner integrable.
\end{compactenum}
\begin{proof}
The fact that $\setdef{\sL_f}{f\in\func[\omega]{\man{M}}}$ contains a
countable point separating subset follows from combining the lemma from the
proof of Theorem~\ref{the:Comega-weak} with the proof of the corresponding
assertion in Corollary~\ref{cor:COinftyweak}\@.  Since
$\sections[\omega]{\tb{\man{M}}}$ is complete, separable, and Suslin, and
since $\func[\omega]{\man{M}}$ is Suslin by
properties~$\C^\omega$-\ref{enum:Comega-complete}\@,~%
$\C^\omega$-\ref{enum:Comega-separable}\@,
and~$\C^\omega$-\ref{enum:Comega-suslin} above, the corollary follows from
Lemma~\ref{lem:weakA}\@, taking
``$\alg{U}=\sections[\omega]{\tb{\man{M}}}$\@,''
``$\alg{V}=\func[\omega]{\man{M}}$\@,'' and
``$\sA=\setdef{\sL_f}{f\in\func[\omega]{\man{M}}}$\@.''
\end{proof}
\end{corollary}

\section{Time-varying vector fields}\label{sec:time-varying}

In this section we consider time-varying vector fields.  The ideas in this
section originate (for us) with the paper of \citet{AAA/RVG:78}\@, and are
nicely summarised in the more recent book of \citet{AAA/YS:04}\@, at least in
the smooth case.  A geometric presentation of some of the constructions can
be found in the paper of \citet{HJS:97b}\@, again in the smooth case, and
\citeauthor{HJS:97b} also considers regularity less than smooth,~\eg~finitely
differentiable or Lipschitz.  There is some consideration of the real
analytic case in \cite{AAA/RVG:78}\@, but this consideration is restricted to
real analytic vector fields admitting a bounded holomorphic extension to a
fixed-width neighbourhood of $\real^n$ in $\complex^n$\@.  One of our
results, the rather nontrivial Theorem~\ref{the:Comega->Chol}\@, is that this
framework of \citet{AAA/RVG:78} is sufficient for the purposes of local
analysis.  However, our treatment of the real analytic case is global,
general, and comprehensive.  To provide some context for our novel treatment
of the real analytic case, we treat the smooth case in some detail, even
though the results are probably mostly known.  (However, we should say that,
even in the smooth case, we could not find precise statements with proofs of
some of the results we give.)  We also treat the finitely differentiable and
Lipschitz cases, so our theory also covers the ``standard'' Carath\'eodory
existence and uniqueness theorem for time-varying ordinary differential
equations,~\cite[\eg][Theorem~54]{EDS:98}\@.  We also consider holomorphic
time-varying vector fields, as these have a relationship to real analytic
time-varying vector fields that is sometimes useful to exploit.

One of the unique facets of our presentation is that we fully explain the
r\^ole of the topologies developed in
Sections~\ref{sec:smooth-topology}\@,~\ref{sec:holomorphic-topology}\@,
and~\ref{sec:analytic-topology}\@.  Indeed, one way to understand the
principal results of this section is that they show that the usual
pointwise\textemdash{}in state and time\textemdash{}conditions placed on
vector fields to regulate the character of their flows can be profitably
phrased in terms of topologies for spaces of vector fields.  While this idea
is not entirely new\textemdash{}it is implicit in the approach of
\cite{AAA/RVG:78}\@\textemdash{}we do develop it comprehensively and in new
directions.

While our principal interest is in vector fields, and also in functions, it
is convenient to conduct much of the development for general vector bundles,
subsequently specialising to vector fields and functions.

\subsection{The smooth case}\label{subsec:smooth-tvf}

Throughout this section we will work with a smooth vector bundle
$\map{\pi}{\man{E}}{\man{M}}$ with a linear connection $\nabla^0$ on
$\man{E}$\@, an affine connection $\nabla$ on $\man{M}$\@, a fibre metric
$\metric_0$ on $\man{E}$\@, and a Riemannian metric $\metric$ on $\man{M}$\@.
This defines the fibre norms $\dnorm{\cdot}_{\ol{\metric}_m}$ on
$\jet{m}{\man{E}}$ and seminorms $p^\infty_{K,m}$\@, $K\subset\man{M}$
compact, $m\in\integernn$\@, on $\sections[\infty]{\man{E}}$ as in
Section~\ref{subsec:COinfty-vb}\@.
\begin{definition}
Let $\map{\pi}{\man{E}}{\man{M}}$ be a smooth vector bundle and let
$\tdomain\subset\real$ be an interval.  A \defn{Carath\'eodory section of
class $\C^\infty$} of $\man{E}$ is a map
$\map{\xi}{\tdomain\times\man{M}}{\man{E}}$ with the following properties:
\begin{compactenum}[(i)]
\item $\xi(t,x)\in\man{E}_x$ for each $(t,x)\in\tdomain\times\man{M}$\@;
\item for each $t\in\tdomain$\@, the map $\map{\xi_t}{\man{M}}{\man{E}}$
defined by $\xi_t(x)=\xi(t,x)$ is of class $\C^\infty$\@;
\item for each $x\in\man{M}$\@, the map $\map{\xi^x}{\tdomain}{\man{E}}$
defined by $\xi^x(t)=\xi(t,x)$ is Lebesgue measurable.
\end{compactenum}
We shall call $\tdomain$ the \defn{time-domain} for the section.  By
$\CFsections[\infty]{\tdomain;\man{E}}$ we denote the set of Carath\'eodory
sections of class $\C^\infty$ of $\man{E}$\@.\oprocend
\end{definition}

Note that the curve $t\mapsto\xi(t,x)$ is in the finite-dimensional vector
space $\man{E}_x$\@, and so Lebesgue measurability of this is unambiguously
defined,~\eg~by choosing a basis and asking for Lebesgue measurability of the
components with respect to this basis.

Now we put some conditions on the time dependence of the derivatives of the section.
\begin{definition}
Let $\map{\pi}{\man{E}}{\man{M}}$ be a smooth vector bundle and let
$\tdomain\subset\real$ be an interval.  A Carath\'eodory section
$\map{\xi}{\tdomain\times\man{M}}{\man{E}}$ of class $\C^\infty$ is
\begin{compactenum}[(i)]
\item \defn{locally integrally $\C^\infty$-bounded} if, for every compact set
$K\subset\man{M}$ and every $m\in\integernn$\@, there exists
$g\in\Lloc^1(\tdomain;\realnn)$ such that
\begin{equation*}
\dnorm{j_m\xi_t(x)}_{\ol{\metric}_m}\le g(t),\qquad(t,x)\in\tdomain\times K,
\end{equation*}
and is
\item \defn{locally essentially $\C^\infty$-bounded} if, for every compact
set $K\subset\man{M}$ and every $m\in\integernn$\@, there exists
$g\in\Lloc^\infty(\tdomain;\realnn)$ such that
\begin{equation*}
\dnorm{j_m\xi_t(x)}_{\ol{\metric}_m}\le g(t),\qquad(t,x)\in\tdomain\times K.
\end{equation*}
\end{compactenum}
The set of locally integrally $\C^\infty$-bounded sections of $\man{E}$ with
time-domain $\tdomain$ is denoted by $\LIsections[\infty]{\tdomain,\man{E}}$
and the set of locally essentially $\C^\infty$-bounded sections of $\man{E}$
with time-domain $\tdomain$ is denoted by
$\LBsections[\infty]{\tdomain;\man{E}}$\@.\oprocend
\end{definition}

Note that $\LBsections[\infty]{\tdomain;\man{M}}\subset
\LIsections[\infty]{\tdomain;\man{M}}$\@, precisely because locally
essentially bounded functions (in the usual sense) are locally integrable (in
the usual sense).

We note that our definitions differ from those
in~\cite{AAA/RVG:78,AAA/YS:04,HJS:97b}\@.  The form of the difference is our
use of connections and jet bundles, aided by Lemma~\ref{lem:Jrdecomp}\@.
In~\cite{AAA/RVG:78} the presentation is developed on Euclidean spaces, and
so the geometric treatment we give here is not necessary.  (One way of
understanding why it is not necessary is that Euclidean space has a canonical
flat connection in which the decomposition of Lemma~\ref{lem:Jrdecomp}
becomes the usual decomposition of derivatives by their order.)  In
\cite{AAA/YS:04} the treatment is on manifolds, and the seminorms are defined
by an embedding of the manifold in Euclidean space by Whitney's Embedding
Theorem~\cite{HW:36}\@.  Also, \citet{AAA/YS:04} use the weak-$\sL$ topology
in the case of vector fields, but we have seen that this is the same as the
usual topology (Theorem~\ref{the:COinfty-weak}).  In \cite{HJS:97b} the
characterisation of Carath\'eodory functions uses Lie differentiation by
smooth vector fields, and the locally convex topology for
$\sections[\infty]{\tb{\man{M}}}$ is not explicitly considered, although it
is implicit in \citeauthor{HJS:97b}'s constructions.  \citeauthor{HJS:97b}
also takes a weak-$\sL$ approach to characterising properties of time-varying
vector fields.  In any case, all approaches can be tediously shown to be
equivalent once the relationships are understood.  An advantage of the
approach we use here is that it does not require coordinate charts or
embeddings to write the seminorms, and it makes the seminorms explicit,
rather than implicitly present.  The disadvantage of our approach is the
added machinery and complication of connections and our jet bundle
decomposition.

The following characterisation of Carath\'eodory sections and their
relatives is also useful and insightful.
\begin{theorem}\label{the:Cinfty-tvsec}
Let\/ $\map{\pi}{\man{E}}{\man{M}}$ be a smooth vector bundle and let\/
$\tdomain\subset\real$ be an interval.  For a map\/
$\map{\xi}{\tdomain\times\man{M}}{\man{E}}$ satisfying\/
$\xi(t,x)\in\man{E}_x$ for each\/ $(t,x)\in\tdomain\times\man{M}$\@, the
following two statements are equivalent:
\begin{compactenum}[(i)]
\item \label{pl:CinftyCF1} $\xi\in\CFsections[\infty]{\tdomain;\man{E}}$\@;
\item \label{pl:CinftyCF2} the map\/ $\tdomain\ni
t\mapsto\xi_t\in\sections[\infty]{\man{E}}$ is measurable,\savenum
\end{compactenum}
the following two statements are equivalent:
\begin{compactenum}[(i)]\resumenum
\item \label{pl:CinftyCF3} $\xi\in\LIsections[\infty]{\tdomain;\man{E}}$\@;
\item \label{pl:CinftyCF4} the map\/ $\tdomain\ni
t\mapsto\xi_t\in\sections[\infty]{\man{E}}$ is measurable and locally Bochner
integrable,\savenum
\end{compactenum}
and the following two statements are equivalent:
\begin{compactenum}[(i)]\resumenum
\item \label{pl:CinftyCF5} $\xi\in\LBsections[\infty]{\tdomain;\man{E}}$\@;
\item \label{pl:CinftyCF6} the map\/ $\tdomain\ni
t\mapsto\xi_t\in\sections[\infty]{\man{E}}$ is measurable and locally
essentially von Neumann bounded.
\end{compactenum}
\begin{proof}
It is illustrative, especially since we will refer to this proof at least
three times subsequently, to understand the general framework of the proof.
Much of the argument has already been carried out in a more general setting
in Lemma~\ref{lem:weakA}\@.

So we let $\alg{V}$ be a locally convex topological vector space over
$\field\in\{\real,\complex\}$\@, let $(\ts{T},\sM)$ be a measurable space,
and let $\map{\Psi}{\ts{T}}{\alg{V}}$\@.  Let us first characterise
measurability of $\Psi$\@.  We use here the results of \citet{GEFT:75} who
studies integrability for functions taking values in locally convex Suslin
spaces.  Thus we assume that $\alg{V}$ is a Hausdorff Suslin space (as is the
case for all spaces of interest to us in this paper).  We let
$\topdual{\alg{V}}$ denote the topological dual of $\alg{V}$\@.  A subset
$S\subset\topdual{\alg{V}}$ is \defn{point separating} if, for distinct
$v_1,v_2\in\alg{V}$\@, there exists $\alpha\in\topdual{\alg{V}}$ such that
$\alpha(v_1)\not=\alpha(v_2)$\@.  \citet{GEFT:75} proves the following result
as his Theorem~1, and whose proof we provide, as it is straightforward and
shows where the (not so straightforward) properties of Suslin spaces are
used.
\begin{prooflemma}
Let\/ $\alg{V}$ be a Hausdorff, Suslin, locally convex topological vector
space over\/ $\field\in\{\real,\complex\}$\@, let $(\ts{T},\sM)$ be a
measurable space, and let $\map{\Psi}{\ts{T}}{\alg{V}}$\@.  If\/
$S\subset\topdual{\alg{V}}$ is point separating, then\/ $\Psi$ is measurable
if and only if\/ $\alpha\scirc\Psi$ is measurable for every\/ $\alpha\in
S$\@.
\begin{subproof}
If $\Psi$ is measurable, then it is obvious that $\alpha\scirc\Psi$ is
measurable for every $\alpha\in\topdual{\alg{V}}$ since such $\alpha$ are
continuous.

Conversely, suppose that $\alpha\scirc\Psi$ is measurable for every
$\alpha\in S$\@.  First of all, locally convex topological vector spaces are
completely regular if they are Hausdorff~\cite[page~16]{HHS/MPW:99}\@.
Therefore, by~\cite[Theorem~6.7.7]{VIB:07b}\@, there is a countable subset of
$S$ that is point separating, so we may as well suppose that $S$ is
countable.  We are now in the same framework as
Lemma~\pldblref{lem:weakA}{pl:weakAmeas}\@, and the proof there applies by
taking ``$\alg{U}=\alg{V}$\@,'' ``$\alg{V}=\field$\@,'' and ``$\sA=S$\@.''
\end{subproof}
\end{prooflemma}

The preceding lemma will allow us to characterise measurability.  Let us now
consider integrability.
\begin{prooflemma}
Let\/ $\alg{V}$ be a complete separable locally convex topological vector
space over\/ $\field\in\{\real,\complex\}$ and let\/ $(\ts{T},\sM,\mu)$ be a
finite measure space.  A measurable function\/ $\map{\Psi}{\ts{T}}{\alg{V}}$
is Bochner integrable if and only if\/ $p\scirc\Psi$ is integrable for every
continuous seminorm\/ $p$ for\/ $\alg{V}$\@.
\begin{subproof}
It follows from~\cite[Theorems~3.2,~3.3]{RB/AD:11} that $\Psi$ is integrable
if $p\scirc\Psi$ is integrable for every continuous seminorm $p$\@.
Conversely, if $\Psi$ is integrable, it is implied that $\Psi$ is Bochner
approximable, and so, by~\cite[Theorem~3.2]{RB/AD:11}\@, we have that
$p\scirc\Psi$ is integrable for every continuous seminorm $p$\@.
\end{subproof}
\end{prooflemma}

\eqref{pl:CinftyCF1}$\iff$\eqref{pl:CinftyCF2} For $x\in\man{M}$ and
$\alpha_x\in\dual{\man{E}}_x$\@, define
$\map{\ev_{\alpha_x}}{\sections[\infty]{\man{E}}}{\real}$ by
$\ev_{\alpha_x}(\xi)=\natpair{\alpha_x}{\xi(x)}$\@.  Clearly $\ev_{\alpha_x}$
is $\real$-linear.  We claim that $\ev_{\alpha_x}$ is continuous.  Indeed,
for a directed set $(I,\preceq)$ and a net $\ifam{\xi}_{i\in I}$ converging
to $\xi$\@,\footnote{Since $\sections[\infty]{\man{E}}$ is metrisable, it
suffices to use sequences.  However, we shall refer to this argument when we
do not use metrisable spaces, so it is convenient to have the general
argument here.} we have
\begin{equation*}
\lim_{i\in I}\ev_{\alpha_x}(\xi_i)=\lim_{i\in I}\alpha_x(\xi_i(x))=
\alpha_x\Bigl(\lim_{i\in I}\xi_i(x)\Bigr)=\alpha_x(\xi(x))=\ev_{\alpha_x}(\xi),
\end{equation*}
using the fact that convergence in the $\CO^\infty$-topology implies
pointwise convergence.  It is obvious that the continuous linear functions
$\ev_{\alpha_x}$\@, $\alpha_x\in\dual{\man{E}}$\@, are point separating.  We
now recall from property~$\CO^\infty$-\ref{enum:COinfty-suslin} for the
smooth $\CO^\infty$-topology that $\sections[\infty]{\man{E}}$ is a Suslin
space with the $\CO^\infty$-topology.  Therefore, by the first lemma above,
it follows that $t\mapsto\xi_t$ is measurable if and only if
$t\mapsto\ev_{\alpha_x}(\xi_t)=\natpair{\alpha_x}{\xi_t(x)}$ is measurable
for every $\alpha_x\in\dual{\man{E}}$\@.  On the other hand, this is
equivalent to $t\mapsto\xi_t(x)$ being measurable for every $x\in\man{M}$
since $t\mapsto\xi_t(x)$ is a curve in the finite-dimensional vector space
$\man{E}_x$\@.  Finally, note that it is implicit in the statement
of~\eqref{pl:CinftyCF2} that $\xi_t$ is smooth, and this part of the
proposition follows easily from these observations.

\eqref{pl:CinftyCF3}$\iff$\eqref{pl:CinftyCF4} Let $\tdomain'\subset\tdomain$
be compact.

First suppose that $\xi\in\LIsections[\infty]{\tdomain;\man{E}}$\@.  By
definition of locally integrally $\C^\infty$-bounded, for each compact
$K\subset\man{M}$ and $m\in\integernn$\@, there exists
$g\in\L^1(\tdomain';\realnn)$ such that
\begin{equation*}
\dnorm{j_m\xi_t(x)}_{\ol{\metric}_m}\le g(t),\qquad(t,x)\in\tdomain'\times K
\quad\implies\quad p^\infty_{K,m}(\xi_t)\le g(t),\qquad t\in\tdomain'.
\end{equation*}
Note that continuity of $p^\infty_{K,m}$ implies that $t\mapsto
p^\infty_{K,m}(\xi_t)$ is measurable.  Therefore,
\begin{equation*}
\int_{\tdomain'}p^\infty_{K,m}(\xi_t)\,\d{t}<\infty,\qquad K\subset\man{M}\
\textrm{compact},\ m\in\integernn.
\end{equation*}
Since $\sections[\infty]{\man{E}}$ is complete and separable, it now follows
from the second lemma above that $t\mapsto\xi_t$ is Bochner integrable on
$\tdomain'$\@.  That is, since $\tdomain'$ is arbitrary, $t\mapsto\xi_t$ is
locally Bochner integrable.

Next suppose that $t\mapsto\xi_t$ is Bochner integrable on $\tdomain$\@.  By
the second lemma above,
\begin{equation*}
\int_{\tdomain'}p^\infty_{K,m}(\xi_t)\,\d{t}<\infty,\qquad K\subset\man{M}\
\textrm{compact},\ m\in\integernn.
\end{equation*}
Therefore, since
\begin{equation*}
\dnorm{j_m\xi_t(x)}_{\ol{\metric}_m}\le
p^\infty_{K,m}(\xi_t),\qquad(t,x)\in\tdomain'\times K,
\end{equation*}
we conclude that $\xi$ is locally integrally $\C^\infty$-bounded since
$\tdomain'$ is arbitrary.

\eqref{pl:CinftyCF5}$\iff$\eqref{pl:CinftyCF6} We recall our discussion of
von Neumann bounded sets in locally convex topological vector spaces
preceding Lemma~\ref{lem:COinftybdd} above.  With this in mind and using
Lemma~\ref{lem:COinftybdd}\@, this part of the theorem follows immediately.
\end{proof}
\end{theorem}

Note that Theorem~\ref{the:Cinfty-tvsec} applies, in particular, to vector
fields and functions, giving the classes $\CF[\infty]{\tdomain;\man{M}}$\@,
$\LIC[\infty]{\tdomain;\man{M}}$\@, and $\LBC[\infty]{\tdomain;\man{M}}$ of
functions, and the classes $\CFsections[\infty]{\tdomain;\tb{\man{M}}}$\@,
$\LIsections[\infty]{\tdomain;\tb{\man{M}}}$\@, and
$\LBsections[\infty]{\tdomain;\tb{\man{M}}}$ of vector fields.  Noting that
we have the alternative weak-$\sL$ characterisation of the
$\CO^\infty$-topology, we can summarise the various sorts of measurability,
integrability, and boundedness for smooth time-varying vector fields as
follows.  In the statement of the result, $\ev_x$ is the ``evaluate at $x$''
map for both functions and vector fields.
\begin{theorem}\label{the:smooth-td-summary}
Let\/ $\man{M}$ be a smooth manifold, let\/ $\tdomain\subset\real$ be a
time-domain, and let\/ $\map{X}{\tdomain\times\man{M}}{\tb{\man{M}}}$ have
the property that\/ $X_t$ is a smooth vector field for each\/
$t\in\tdomain$\@.  Then the following four statements are equivalent:
\begin{compactenum}[(i)]
\item $t\mapsto X_t$ is measurable;
\item $t\mapsto\lieder{X_t}{f}$ is measurable for every\/ $f\in\func[\infty]{\man{M}}$\@;
\item $t\mapsto\ev_x\scirc X_t$ is measurable for every\/ $x\in\man{M}$\@;
\item $t\mapsto\ev_x\scirc\lieder{X_t}{f}$ is measurable for every\/
$f\in\func[\infty]{\man{M}}$ and every\/ $x\in\man{M}$\@,\savenum
\end{compactenum}
the following two statements are equivalent:
\begin{compactenum}[(i)]\resumenum
\item $t\mapsto X_t$ is locally Bochner integrable;
\item $t\mapsto\lieder{X_t}{f}$ is locally Bochner integrable for every\/
$f\in\func[\infty]{\man{M}}$\@,\savenum
\end{compactenum}
and the following two statements are equivalent:
\begin{compactenum}[(i)]\resumenum
\item $t\mapsto X_t$ is locally essentially von Neumann bounded;
\item $t\mapsto\lieder{X_t}{f}$ is locally essentially von Neumann bounded
for every\/ $f\in\func[\infty]{\man{M}}$\@.
\end{compactenum}
\begin{proof}
This follows from Theorem~\ref{the:Cinfty-tvsec}\@, along with
Corollary~\ref{cor:COinftyweak}\@.
\end{proof}
\end{theorem}

Let us now discuss flows of vector fields from
$\LIsections[\infty]{\tdomain;\tb{\man{M}}}$\@.  To do so, let us provide the
definition of the usual attribute of integral curves, but on manifolds.
\begin{definition}
Let\/ $\man{M}$ be a smooth manifold and let\/ $\tdomain\subset\real$ be an
interval.
\begin{compactenum}[(i)]
\item A function $\map{f}{\interval[a,b]}{\real}$ is \defn{absolutely
continuous} if there exists $g\in\L^1(\interval[a,b];\real)$ such that
\begin{equation*}
f(t)=f(a)+\int_a^tg(\tau)\,\d{\tau},\qquad t\in\interval[a,b].
\end{equation*}
\item A function $\map{f}{\tdomain}{\real}$ is \defn{locally absolutely
continuous} if $f|\tdomain'$ is absolutely continuous for every compact
subinterval $\tdomain'\subset\tdomain$\@.
\item A curve $\map{\gamma}{\tdomain}{\man{M}}$ is \defn{locally absolutely
continuous} if $\phi\scirc\gamma$ is locally absolutely continuous for every $\phi\in\func[\infty]{\man{M}}$\@.\oprocend
\end{compactenum}
\end{definition}

One easily verifies that a curve is locally absolutely continuous according
to our definition if and only if its local representative is locally
absolutely continuous in any coordinate chart.

We then have the following existence, uniqueness, and regularity result for
locally integrally bounded vector fields.  In the statement of the result,
we use the notation
\begin{equation*}
|a,b|=\begin{cases}\interval[a,b],&a\le b,\\\interval[b,a],&
b<a.\end{cases}
\end{equation*}
In the following result, we do not provide the comprehensive list of
properties of the flow, but only those required to make sense of its
regularity with respect to initial conditions, as per our
specification~\ref{enum:consreg} for our theory in
Section~\ref{subsec:attributes}\@.
\begin{theorem}\label{the:Cinfty-flows}
Let\/ $\man{M}$ be a smooth manifold, let\/ $\tdomain$ be an interval, and
let\/ $X\in\LIsections[\infty]{\tdomain;\tb{\man{M}}}$\@.  Then there exist a
subset\/ $D_X\subset\tdomain\times\tdomain\times\man{M}$ and a map\/
$\map{\flow*{X}}{D_X}{\man{M}}$ with the following properties for each\/
$(t_0,x_0)\in\tdomain\times\man{M}$\@:
\begin{compactenum}[(i)]
\item the set
\begin{equation*}
\tdomain_X(t_0,x_0)=\setdef{t\in\tdomain}{(t,t_0,x_0)\in D_X}
\end{equation*}
is an interval;
\item there exists a locally absolutely continuous curve\/ $t\mapsto\xi(t)$
satisfying
\begin{equation*}
\xi'(t)=X(t,\xi(t)),\quad\xi(t_0)=x_0,
\end{equation*}
for almost all\/ $t\in|t_0,t_1|$ if and only if\/
$t_1\in\tdomain_X(t_0,x_0)$\@;
\item $\deriv{}{t}\flow*{X}(t,t_0,x_0)=X(t,\flow*{X}(t,t_0,x_0))$ for almost
all\/ $t\in\tdomain_X(t_0,x_0)$\@;
\item for each\/ $t\in\tdomain$ for which\/ $(t,t_0,x_0)\in D_X$\@, there
exists a neighbourhood\/ $\nbhd{U}$ of\/ $x_0$ such that the mapping\/
$x\mapsto\flow*{X}(t,t_0,x)$ is defined and of class\/ $\C^\infty$ on\/
$\nbhd{U}$\@.
\end{compactenum}
\begin{proof}
We observe that the requirement that
$X\in\LIsections[\infty]{\tdomain;\tb{\man{M}}}$ implies that, in any
coordinate chart, the components of $X$ and their derivatives are all bounded
by a locally integrable function.  This, in particular, implies that, in any
coordinate chart for $\man{M}$\@, the ordinary differential equation
associated to the vector field $X$ satisfies the usual conditions for
existence and uniqueness of solutions as per, for
example,~\cite[Theorem~54]{EDS:98}\@.  Of course, the differential equation
satisfies conditions much stronger than this, and we shall see how to use
these in our argument below.

The first three assertions are now part of the standard existence theorem for
solutions of ordinary differential equations, along with the usual Zorn's
Lemma argument for the existence of a maximal interval on which integral
curves is defined.

In the sequel we denote $\flow{X}{t}[t_0](x)=\flow*{X}(t,t_0,x_0)$\@.

For the fourth assertion we first make some constructions with vector fields
on jet bundles, more or less following~\cite[\S4.4]{DJS:89}\@.  We let
$\man{M}^2=\man{M}\times\man{M}$ and we consider $\man{M}^2$ as a fibred
manifold, indeed a trivial fibre bundle, over $\man{M}$ by
$\map{\pr_1}{\man{M}^2}{\man{M}}$\@,~\ie~by projection onto the first factor.
A section of this fibred manifold is naturally identified with a smooth map
$\map{\Phi}{\man{M}}{\man{M}}$ by $x\mapsto(x,\Phi(x))$\@.  We introduce the
following notation:
\begin{compactenum}
\item $\jet{m}{\pr_1}$\@: the bundle of $m$-jets of sections of the fibred
manifold $\map{\pr_1}{\man{M}^2}{\man{M}}$\@;
\item $\vb{\pr_{1,m}}$\@: the vertical bundle of the fibred manifold $\map{\pr_{1,m}}{\jet{m}{\pr_1}}{\man{M}}$\@;
\item $\vb{\pr_1}$\@: the vertical bundle of the fibred manifold $\map{\pr_1}{\man{M}^2}{\man{M}}$\@;
\item $\nu$\@: the projection $\pr_1\scirc(\tbproj{\man{M}^2}|\vb{\pr_1})$\@;
\item $\jet{m}{\nu}$\@: the bundle of $m$-jets of sections of the fibred
manifold $\map{\nu}{\vb{\pr_1}}{\man{M}}$\@.
\end{compactenum}
With this notation, we have the following lemma.
\begin{prooflemma}
There is a canonical diffeomorphism\/
$\map{\alpha_m}{\jet{m}{\nu}}{\vb{\pr_{1,m}}}$\@.
\begin{subproof}
We describe the diffeomorphism, and then note that the verification that it
is, in fact, a diffeomorphism is a fact easily checked in jet bundle
coordinates.

Let $I\subset\real$ be an interval with $0\in\interior(I)$ and consider a
smooth map $\map{\phi}{I\times\man{M}}{\man{M}\times\man{M}}$ of the form
$\phi(t,x)=(x,\phi_1(t,x))$ for a smooth map $\phi_1$\@.  We let
$\phi_t(x)=\phi^x(t)=\phi(t,x)$\@.  We then have maps
\begin{equation*}
\mapdef{j^x_m\phi}{I}{\jet{m}{\pr_1}}{t}{j_m\phi_t(x)}
\end{equation*}
and
\begin{equation*}
\mapdef{\phi'}{\man{M}}{\vb{\pr_1}}{x}{\derivatzero{}{t}\phi^x(t).}
\end{equation*}
Note that the curve $j^x_m\phi$ is a curve in the fibre of
$\map{\pr_{1,m}}{\jet{m}{\pr_1}}{\man{M}}$\@.  Thus we can sensibly define
$\alpha_m$ by
\begin{equation*}
\alpha_m(j_m\phi'(x))=\derivatzero{}{t}j^x_m\phi(t).
\end{equation*}
In jet bundle coordinates, one can check that $\alpha_m$ has the local
representative
\begin{equation*}
((\vect{x}_1,(\vect{x}_2,\vect{A}_0)),(\mat{B}_1,\mat{A}_1,\dots,
\mat{B}_m,\mat{A}_m))\mapsto
((\vect{x}_1,(\vect{x}_2,\mat{B}_1,\dots,\mat{B}_m)),
(\mat{A}_0,\mat{A}_1,\dots,\mat{A}_m)),
\end{equation*}
showing that $\alpha_m$ is indeed a diffeomorphism.
\end{subproof}
\end{prooflemma}

Given a smooth vector field $Y$ on $\man{M}$\@, we define a vector field
$\tilde{Y}$ on $\man{M}^2$ by $\tilde{Y}(x_1,x_2)=(0_{x_1},Y(x_2))$\@.  Note
that we have the following commutative diagram
\begin{equation*}
\xymatrix{{\man{M}^2}\ar[r]^{\tilde{Y}}\ar[d]_{\pr_1}&
{\vb{\pr_1}}\ar[d]^{\nu}\\{\man{M}}\ar@{=}[r]&{\man{M}}}
\end{equation*}
giving $\tilde{Y}$ as a morphism of fibred manifolds.  It is thus a candidate
to have its $m$-jet taken, giving a morphism of fibred manifolds
$\map{j_m\tilde{Y}}{\jet{m}{\pr_1}}{\jet{m}{\nu}}$\@.  By the lemma,
$\alpha_m\scirc j_m\tilde{Y}$ is a vertical vector field on $\jet{m}{\pr_1}$
that we denote by $\nu_mY$\@, the \defn{$m$th vertical prolongation} of
$Y$\@.  Let us verify that this is a vector field.  First of all, for a
section $\tilde{\Phi}$ of $\pr_1$ given by $x\mapsto(x,\Phi(x))$\@, note that
$\tilde{Y}\scirc\tilde{\Phi}(x)=(0_x,Y(\Phi(x)))$\@, and so
$j_m(\tilde{Y}\scirc\tilde{\Phi})(x)$ is vertical.  By the notation from the
proof of the lemma, we can write
$j_m(\tilde{Y}\scirc\tilde{\Phi})(x)=j_m\phi'(x)$ for some suitable map $\phi$
as in the lemma.  We then have
\begin{equation*}
\alpha_m\scirc j_m(\tilde{Y}\scirc\tilde{\Phi})(x)=
\alpha_m(j_m\phi'(x))\in\vb[j_m\tilde{\Phi}(x)]{\pr_{1,m}}.
\end{equation*}
Therefore,
\begin{equation*}
\tbproj{\jet{m}{\pr_1}}(\alpha_m\scirc j_m\tilde{Y}(j_m\tilde{\Phi}(x)))=
\tbproj{\jet{m}{\pr_1}}(\alpha_m\scirc
j_m(\tilde{Y}\scirc\tilde{\Phi})(x))=j_m\tilde{\Phi}(x).
\end{equation*}

Note that since $\jet{m}{\pr_1}$ is naturally identified with
$\jet{m}{(\man{M};\man{M})}$ via the identification
\begin{equation*}
j_m\tilde{\Phi}(x)\mapsto j_m\Phi(x)
\end{equation*}
if $\tilde{\Phi}(x)=(x,\Phi(x))$\@, we can as well think of $\nu_mY$ as being
a vector field on the latter space.  Sorting through all the definitions
gives the form of $\nu_mY$ in coordinates as
\begin{equation}\label{eq:numY}
((\vect{x}_1,\vect{x}_2),\mat{A}_1,\dots,\mat{A}_m)\mapsto
(((\vect{x}_1,\vect{x}_2),\mat{A}_1,\dots,\mat{A}_m),
\vect{0},\vect{Y},\linder{\vect{Y}},\dots,\linder[m]{\vect{Y}}).
\end{equation}

We now apply the above constructions, for each fixed $t\in\tdomain$\@, to get
the vector field $\nu_mX_t$\@, and so the time-varying vector field $\nu_mX$
defined by $\nu_mX(t,j_m\Phi(x))=\nu_mX_t(j_m\Phi(x))$ on
$\jet{m}{(\man{M};\man{M})}$\@.  The definition of
$\LIsections[\infty]{\tdomain;\tb{\man{M}}}$\@, along with the coordinate
formula~\eqref{eq:numY}\@, shows that $\nu_mX$ satisfies the standard
conditions for existence and uniqueness of integral curves, and so its flow
depends continuously on initial condition~\cite[Theorem~55]{EDS:98}\@.

The fourth part of the theorem, therefore, will follow if we can show that
\begin{compactenum}
\item \label{enum:smoothflow1} for each $m\in\integernn$\@, the flow of
$\nu_mX$ depends on the initial condition in $\man{M}$ in a $\C^m$ way,
\item \label{enum:smoothflow2} $\flow{\nu_mX}{t}[t_0](j_m\flow{X}{t_0}[t_0](x_0))=
j_m\flow{X}{t}[t_0](j_m\flow{X}{t_0}[t_0](x_0))$\@, and
\item \label{enum:smoothflow3} if $\{t\}\times\{t_0\}\times\nbhd{U}\subset
D_X$\@, then $\{t\}\times\{t_0\}\times\pr_{1,m}^{-1}(\nbhd{U})\subset
D_{\nu_mX}$\@.
\end{compactenum}
We ask for property~\ref{enum:smoothflow3} to ensure that the domain of
differentiability does not get too small as the order of the derivatives gets
large.

To prove these assertions, it suffices to work locally.  According
to~\eqref{eq:numY}\@, we have the time-dependent differential equation defined
on
\begin{equation*}
\nbhd{U}\times\lin{\real^n}{\real^n}\times\dots\times
\symlin[m]{\real^n}{\real^n},
\end{equation*}
where $\nbhd{U}$ is an open subset of $\real^n$\@, and given by
\begin{align*}
\dot{\vect{\gamma}}(t)=&\;\vect{X}(t,\vect{\gamma}(t)),\\
\dot{\mat{A}}_1(t)=&\;\linder{\vect{X}}(t,\vect{\gamma}(t)),\\
\dot{\mat{A}}_2(t)=&\;\linder[2]{\vect{X}}(t,\vect{\gamma}(t)),\\
\vdots\;&\\
\dot{\mat{A}}_m(t)=&\;\linder[m]{\vect{X}}(t,\vect{\gamma}(t)),
\end{align*}
$(t,\vect{x})\mapsto(\vect{x},\vect{X}(t,\vect{x}))$ being the local
representative of $X$\@.  The initial conditions of interest for the vector
field $\nu_mX$ are of the form $j_m\flow{X}{t_0}[t_0](x)$\@.  In coordinates,
keeping in mind that $\flow{X}{t_0}[t_0]=\id_{\man{M}}$\@, this gives
\begin{equation}\label{eq:smoothflowic}
\vect{\gamma}(t_0)=\vect{x}_0,\ \mat{A}_1(t_0)=\mat{I}_n,\
\mat{A}_j(t_0)=\mat{0},\qquad j\ge2.
\end{equation}
Let us denote by $t\mapsto\vect{\gamma}(t,t_0,\vect{x})$ and
$t\mapsto\mat{A}_j(t,t_0,\vect{x})$\@, $j\in\{1,\dots,m\}$\@, the solutions
of the differential equations above with these initial conditions.

We will show that assertions~\ref{enum:smoothflow1}--\ref{enum:smoothflow3}
hold by induction on $m$\@.  In doing this, we will need to understand how
differential equations depending differentiably on state also have solutions
depending differentiably on initial condition.  Such a result is not readily
found in the textbook literature, as this latter is typically concerned with
continuous dependence on initial conditions for cases with measurable
time-dependence, and on differentiable dependence when the dependence on time
is also differentiable. However, the general case (much more general than we
need here) is worked out by \citet{FS/HvdM:00}\@.

For $m=0$\@, the assertions are simply the result of the usual continuous
dependence on initial conditions~\cite[\eg][Theorem~55]{EDS:98}\@.  Let us
consider the case $m=1$\@.  In this case, the properties of
$\LIsections[\infty]{\tdomain;\tb{\man{M}}}$ ensure that the hypotheses
required to apply Theorem~2.1 of \cite{FS/HvdM:00} hold for the
differential equation
\begin{align*}
\dot{\vect{\gamma}}(t)=&\;\vect{X}(t,\vect{\gamma}(t)),\\
\dot{\mat{A}}_1(t)=&\;\linder{\vect{X}}(t,\vect{\gamma}(t)).
\end{align*}
This allows us to conclude that
$\vect{x}\mapsto\vect{\gamma}(t,t_0,\vect{x})$ is of class $\C^1$\@.  This
establishes the assertion~\ref{enum:smoothflow1} in this case.  Therefore, on
a suitable domain, $j_1\flow{X}{t}[t_0]$ is well-defined.  In coordinates the
map $\map{j_1\flow{X}{t}[t_0]}{\jet{1}{(\man{M};\man{M})}}
{\jet{1}{(\man{M};\man{M})}}$ is given by
\begin{equation}\label{eq:smoothflowm=1}
(\vect{x},\vect{y},\mat{B}_1)\mapsto(\vect{x},\vect{\gamma}(t,t_0,\vect{x}),
\plinder{3}{\vect{\gamma}}(t,t_0,\vect{x})\scirc\mat{B}_1),
\end{equation}
this by the Chain Rule.  We have
\begin{equation*}
\deriv{}{t}\plinder{3}{\vect{\gamma}}(t,t_0,\vect{x})=
\plinder{3}{(\tderiv{}{t}\vect{\gamma}(t,t_0,\vect{x}))}=
\linder{\vect{X}}(t,\vect{\gamma}(t,t_0,\vect{x})),
\end{equation*}
the swapping of the time and spatial derivatives being valid
by~\cite[Corollary~2.2]{FS/HvdM:00}\@.  Combining this
with~\eqref{eq:smoothflowm=1} and the initial
conditions~\eqref{eq:smoothflowic} shows that
assertion~\ref{enum:smoothflow2} holds for $m=1$\@.  Moreover, since
$\mat{A}_1(t,t_0,\vect{x})$ is obtained by merely integrating a continuous
function of $t$ from $t_0$ to $t$\@, we also conclude that
assertion~\ref{enum:smoothflow3} holds.

Now suppose that assertions~\ref{enum:smoothflow1}--\ref{enum:smoothflow3}
hold for $m$\@.  Again, the properties of
$\LIsections[\infty]{\tdomain;\tb{\man{M}}}$ imply that the hypotheses of
Theorem~2.1 of \cite{FS/HvdM:00} hold, and so solutions of the differential
equation
\begin{align*}
\dot{\vect{\gamma}}(t)=&\;\vect{X}(t,\vect{\gamma}(t)),\\
\dot{\mat{A}}_1(t)=&\;\linder{\vect{X}}(t,\vect{\gamma}(t)),\\
\dot{\mat{A}}_2(t)=&\;\linder[2]{\vect{X}}(t,\vect{\gamma}(t)),\\
\vdots\;&\\
\dot{\mat{A}}_m(t)=&\;\linder[m]{\vect{X}}(t,\vect{\gamma}(t))
\end{align*}
depend continuously differentiably on initial condition.  By the induction
hypothesis applied to the assertion~\ref{enum:smoothflow2}\@, this means that
\begin{equation*}
(t,x)\mapsto\flow{\nu_mX}{t}[t_0](j_m\flow{X}{t_0}[t_0](x))=
j_m\flow{X}{t}[t_0](x)
\end{equation*}
depends continuously differentiably on $x$\@, and so we conclude that
$(t,x)\mapsto\flow{X}{t}[t_0](x)$ depends on $x$ in a $\C^{m+1}$ manner.
This establishes assertion~\ref{enum:smoothflow1} for $m+1$\@.  After an
application of the Chain Rule for high-order derivatives
(see~\cite[Supplement~2.4A]{RA/JEM/TSR:88}) we can, admittedly after just a
few moments thought, see that the local representative of
$j_{m+1}\flow{X}{t}[t_0](j_{m+1}\flow{X}{t_0}[t_0](x))$ is
\begin{equation*}
(\vect{x},\vect{\gamma}(t,t_0,\vect{x}),
\plinder{3}{\vect{\gamma}}(t,t_0,\vect{x}),\dots,
\plinder[m+1]{3}{\vect{\gamma}}(t,t_0,\vect{x})),
\end{equation*}
keeping in mind the initial conditions~\eqref{eq:smoothflowic} in coordinates.

By the induction hypothesis,
\begin{equation*}
\deriv{}{t}\plinder[j]{3}{\vect{\gamma}}(t)=
\linder[j]{\vect{X}}(t,\vect{\gamma}(t,t_0,\vect{x})),\qquad
j\in\{1,\dots,m\}.
\end{equation*}
Using Corollary~2.2 of~\cite{FS/HvdM:00} we compute
\begin{equation*}
\deriv{}{t}\plinder[m+1]{3}{\vect{\gamma}}(t,t_0,\vect{x})=
\linder{(\tderiv{}{t}\plinder[m]{3}{\vect{\gamma}}(t,t_0,\vect{x}))}=
\linder[m+1]{\vect{X}}(t,\vect{\gamma}(t,t_0,\vect{x})),
\end{equation*}
giving assertion~\ref{enum:smoothflow2} for $m+1$\@.  Finally, by the
induction hypothesis and since $\mat{A}_{m+1}(t,t_0,\vect{x})$ is obtained by
simple integration from $t_0$ to $t$\@, we conclude that
assertion~\ref{enum:smoothflow3} holds for $m+1$\@.
\end{proof}
\end{theorem}

\subsection{The finitely differentiable or Lipschitz case}

The requirement that the flow depends smoothly on initial conditions is not
always essential, even when the vector field itself depends smoothly on the
state.  In such cases as this, one may want to consider classes of vector
fields characterised by one of the weaker topologies described in
Section~\ref{subsec:COm-topology}\@.  Let us see how to do this.  In this
section, so as to be consistent with our definition of Lipschitz norms in
Section~\ref{subsec:COlip-topology}\@, we suppose that the affine connection
$\nabla$ on $\man{M}$ is the Levi-Civita connection for the Riemannian metric
$\metric$ and that the vector bundle connection $\nabla^0$ in $\man{E}$ is
$\metric_0$-orthogonal.
\begin{definition}
Let $\map{\pi}{\man{E}}{\man{M}}$ be a smooth vector bundle and let
$\tdomain\subset\real$ be an interval.  Let $m\in\integernn$ and let
$m'\in\{0,\lip\}$\@.  A \defn{Carath\'eodory section of class $\C^{m+m'}$} of
$\man{E}$ is a map $\map{\xi}{\tdomain\times\man{M}}{\man{E}}$ with the
following properties:
\begin{compactenum}[(i)]
\item $\xi(t,x)\in\man{E}_x$ for each $(t,x)\in\tdomain\times\man{M}$\@;
\item for each $t\in\tdomain$\@, the map $\map{\xi_t}{\man{M}}{\man{E}}$
defined by $\xi_t(x)=\xi(t,x)$ is of class $\C^{m+m'}$\@;
\item for each $x\in\man{M}$\@, the map $\map{\xi^x}{\tdomain}{\man{E}}$
defined by $\xi^x(t)=\xi(t,x)$ is Lebesgue measurable.
\end{compactenum}
We shall call $\tdomain$ the \defn{time-domain} for the section.  By
$\CFsections[m+m']{\tdomain;\man{E}}$ we denote the set of Carath\'eodory
sections of class $\C^{m+m'}$ of $\man{E}$\@.\oprocend
\end{definition}

Now we put some conditions on the time dependence of the derivatives of the
section.
\begin{definition}
Let $\map{\pi}{\man{E}}{\man{M}}$ be a smooth vector bundle and let
$\tdomain\subset\real$ be an interval.  Let $m\in\integernn$ and let
$m'\in\{0,\lip\}$\@.  A Carath\'eodory section
$\map{\xi}{\tdomain\times\man{M}}{\man{E}}$ of class $\C^{m+m'}$ is
\begin{compactenum}[(i)]
\item \defn{locally integrally $\C^{m+m'}$-bounded} if:
\begin{compactenum}[(a)]
\item $m'=0$\@: for every compact set $K\subset\man{M}$\@, there exists
$g\in\Lloc^1(\tdomain;\realnn)$ such that
\begin{equation*}
\dnorm{j_m\xi_t(x)}_{\ol{\metric}_m}\le
g(t),\qquad(t,x)\in\tdomain\times K;
\end{equation*}
\item $m'=\lip$\@: for every compact set $K\subset\man{M}$\@, there exists
$g\in\Lloc^1(\tdomain;\realnn)$ such that
\begin{equation*}
\dil{j_m\xi_t}(x),\dnorm{j_m\xi_t(x)}_{\ol{\metric}_m}\le
g(t),\qquad(t,x)\in\tdomain\times K,
\end{equation*}
\end{compactenum}
and is
\item \defn{locally essentially $\C^{m+m'}$-bounded} if:
\begin{compactenum}[(a)]
\item $m'=0$\@: for every compact set $K\subset\man{M}$\@, there exists
$g\in\Lloc^\infty(\tdomain;\realnn)$ such that
\begin{equation*}
\dnorm{j_m\xi_t(x)}_{\ol{\metric}_m}\le
g(t),\qquad(t,x)\in\tdomain\times K;
\end{equation*}
\item $m'=\lip$\@: for every compact set $K\subset\man{M}$\@, there exists
$g\in\Lloc^\infty(\tdomain;\realnn)$ such that
\begin{equation*}
\dil{j_m\xi_t}(x),\dnorm{j_m\xi_t(x)}_{\ol{\metric}_m}\le
g(t),\qquad(t,x)\in\tdomain\times K.
\end{equation*}
\end{compactenum}
\end{compactenum}
The set of locally integrally $\C^{m+m'}$-bounded sections of $\man{E}$ with
time-domain $\tdomain$ is denoted by $\LIsections[m+m']{\tdomain,\man{E}}$
and the set of locally essentially $\C^{m+m'}$-bounded sections of $\man{E}$
with time-domain $\tdomain$ is denoted by
$\LBsections[m+m']{\tdomain;\man{E}}$\@.\oprocend
\end{definition}

\begin{theorem}\label{the:Cmm'-tvsec}
Let\/ $\map{\pi}{\man{E}}{\man{M}}$ be a smooth vector bundle and let\/
$\tdomain\subset\real$ be an interval.  Let\/ $m\in\integernn$ and let\/
$m'\in\{0,\lip\}$\@.  For a map\/ $\map{\xi}{\tdomain\times\man{M}}{\man{E}}$
satisfying\/ $\xi(t,x)\in\man{E}_x$ for each\/
$(t,x)\in\tdomain\times\man{M}$\@, the following two statements are
equivalent:
\begin{compactenum}[(i)]
\item \label{pl:Cmm'CF1} $\xi\in\CFsections[m+m']{\tdomain;\man{E}}$\@;
\item \label{pl:Cmm'CF2} the map\/ $\tdomain\ni
t\mapsto\xi_t\in\sections[m+m']{\man{E}}$ is measurable,\savenum
\end{compactenum}
the following two statements are equivalent:
\begin{compactenum}[(i)]\resumenum
\item \label{pl:Cmm'CF3} $\xi\in\LIsections[m+m']{\tdomain;\man{E}}$\@;
\item \label{pl:Cmm'CF4} the map\/ $\tdomain\ni
t\mapsto\xi_t\in\sections[m+m']{\man{E}}$ is measurable and locally Bochner
integrable,\savenum
\end{compactenum}
and the following two statements are equivalent:
\begin{compactenum}[(i)]\resumenum
\item \label{pl:Cmm'CF5} $\xi\in\LBsections[m+m']{\tdomain;\man{E}}$\@;
\item \label{pl:Cmm'CF6} the map\/ $\tdomain\ni
t\mapsto\xi_t\in\sections[m+m']{\man{E}}$ is measurable and locally
essentially von Neumann bounded.
\end{compactenum}
\begin{proof}
\eqref{pl:Cmm'CF1}$\iff$\eqref{pl:Cmm'CF2} For $x\in\man{M}$ and
$\alpha_x\in\dual{\man{E}}_x$\@, define
$\map{\ev_{\alpha_x}}{\sections[m+m']{\man{E}}}{\real}$ by
$\ev_{\alpha_x}(\xi)=\natpair{\alpha_x}{\xi(x)}$\@.  It is easy to show that
$\ev_{\alpha_x}$ is continuous and that the set of continuous functionals
$\ev_{\alpha_x}$\@, $\alpha_x\in\dual{\man{E}}_x$\@, is point separating.
Since $\sections[m+m']{\man{E}}$ is a Suslin space (properties
$\CO^m$-\ref{enum:COm-suslin} and~$\CO^{m+\lip}$-\ref{enum:COmm'-suslin}),
this part of the theorem follows in the same manner as the corresponding part
of Theorem~\ref{the:Cinfty-tvsec}\@.

\eqref{pl:Cmm'CF3}$\iff$\eqref{pl:Cmm'CF4} Since $\sections[m+m']{\man{E}}$
is complete and separable (by properties~$\CO^m$-\ref{enum:COm-complete}
and~$\CO^m$-\ref{enum:COm-separable}\@,
and~$\CO^{m+\lip}$-\ref{enum:COmm'-complete}
and~$\CO^{m+\lip}$-\ref{enum:COmm'-separable}), the arguments from the
corresponding part of Theorem~\ref{the:Cinfty-tvsec} apply here, taking note
of the definition of the seminorms $p^{\lip}_K(\xi)$ in case $m'=\lip$\@.

\eqref{pl:Cmm'CF5}$\iff$\eqref{pl:Cmm'CF6} We recall our discussion of von
Neumann bounded sets in locally convex topological vector spaces preceding
Lemma~\ref{lem:COinftybdd} above.  With this in mind and using
Lemma~\ref{lem:CObdd}\@, this part of the proposition follows immediately.
\end{proof}
\end{theorem}

Note that Theorem~\ref{the:Cmm'-tvsec} applies, in particular, to vector
fields and functions, giving the classes $\CF[m+m']{\tdomain;\man{M}}$\@,
$\LIC[m+m']{\tdomain;\man{M}}$\@, and $\LBC[m+m']{\tdomain;\man{M}}$ of
functions, and the classes $\CFsections[m+m']{\tdomain;\tb{\man{M}}}$\@,
$\LIsections[m+m']{\tdomain;\tb{\man{M}}}$\@, and
$\LBsections[m+m']{\tdomain;\tb{\man{M}}}$ of vector fields.  Noting that we
have the alternative weak-$\sL$ characterisation of the
$\CO^{m+m'}$-topology, we can summarise the various sorts of measurability,
integrability, and boundedness for smooth time-varying vector fields as
follows.  In the statement of the result, $\ev_x$ is the ``evaluate at $x$''
map for both functions and vector fields.
\begin{theorem}\label{the:mm'-td-summary}
Let\/ $\man{M}$ be a smooth manifold, let\/ $\tdomain\subset\real$ be a
time-domain, let\/ $m\in\integernn$\@, let\/ $m'\in\{0,\lip\}$\@, and let\/
$\map{X}{\tdomain\times\man{M}}{\tb{\man{M}}}$ have the property that\/ $X_t$
is a vector field of class\/ $\C^{m+m'}$ for each\/ $t\in\tdomain$\@.  Then
the following four statements are equivalent:
\begin{compactenum}[(i)]
\item $t\mapsto X_t$ is measurable;
\item $t\mapsto\lieder{X_t}{f}$ is measurable for every\/ $f\in\func[\infty]{\man{M}}$\@;
\item $t\mapsto\ev_x\scirc X_t$ is measurable for every\/ $x\in\man{M}$\@;
\item $t\mapsto\ev_x\scirc\lieder{X_t}{f}$ is measurable for every\/
$f\in\func[\infty]{\man{M}}$ and every\/ $x\in\man{M}$\@,\savenum
\end{compactenum}
the following two statements are equivalent:
\begin{compactenum}[(i)]\resumenum
\item $t\mapsto X_t$ is locally Bochner integrable;
\item $t\mapsto\lieder{X_t}{f}$ is locally Bochner integrable for every\/
$f\in\func[\infty]{\man{M}}$\@,\savenum
\end{compactenum}
and the following two statements are equivalent:
\begin{compactenum}[(i)]\resumenum
\item $t\mapsto X_t$ is locally essentially von Neumann bounded;
\item $t\mapsto\lieder{X_t}{f}$ is locally essentially von Neumann bounded
for every\/ $f\in\func[\infty]{\man{M}}$\@.
\end{compactenum}
\begin{proof}
This follows from Theorem~\ref{the:Cmm'-tvsec}\@, along with
Corollaries~\ref{cor:COmweak} and~\ref{cor:COmm'weak}\@.
\end{proof}
\end{theorem}

It is also possible to state an existence, uniqueness, and regularity theorem
for flows of vector fields that depend on state in a finitely differentiable
or Lipschitz manner.
\begin{theorem}\label{the:Cmm'-flows}
Let\/ $\man{M}$ be a smooth manifold, let\/ $\tdomain$ be an interval, let\/
$m\in\integernn$\@, and let\/
$X\in\LIsections[m+\lip]{\tdomain;\tb{\man{M}}}$\@.  Then there exist
a subset\/ $D_X\subset\tdomain\times\tdomain\times\man{M}$ and a map\/
$\map{\flow*{X}}{D_X}{\man{M}}$ with the following properties for each\/ $(t_0,x_0)\in\tdomain\times\man{M}$\@:
\begin{compactenum}[(i)]
\item the set
\begin{equation*}
\tdomain_X(t_0,x_0)=\setdef{t\in\tdomain}{(t,t_0,x_0)\in D_X}
\end{equation*}
is an interval;
\item there exists a locally absolutely continuous curve\/ $t\mapsto\xi(t)$
satisfying
\begin{equation*}
\xi'(t)=X(t,\xi(t)),\quad\xi(t_0)=x_0,
\end{equation*}
for almost all\/ $t\in|t_0,t_1|$ if and only if\/
$t_1\in\tdomain_X(t_0,x_0)$\@;
\item $\deriv{}{t}\flow*{X}(t,t_0,x_0)=X(t,\flow*{X}(t,t_0,x_0))$ for almost
all\/ $t\in\tdomain_X(t_0,x_0)$\@;
\item for each\/ $t\in\tdomain$ for which\/ $(t,t_0,x_0)\in D_X$\@, there
exists a neighbourhood\/ $\nbhd{U}$ of\/ $x_0$ such that the mapping\/
$x\mapsto\flow*{X}(t,t_0,x)$ is defined and of class\/ $\C^m$ on\/
$\nbhd{U}$\@.
\end{compactenum}
\begin{proof}
The proof here is by truncation of the proof of
Theorem~\ref{the:Cinfty-flows} from ``$\infty$'' to ``$m$\@.''
\end{proof}
\end{theorem}

\subsection{The holomorphic case}

While we are not \emph{per se} interested in time-varying holomorphic vector
fields, our understanding of time-varying real analytic vector
fields\textemdash{}in which we are most definitely
interested\textemdash{}is connected with an understanding of the
holomorphic case,~\cf~Theorem~\ref{the:Comega->Chol}\@.

We begin with definitions that are similar to the smooth case, but which rely
on the holomorphic topologies introduced in Section~\ref{subsec:COhol-vb}\@.
We will consider an holomorphic vector bundle $\map{\pi}{\man{E}}{\man{M}}$
with an Hermitian fibre metric $\metric$\@.  This defines the seminorms
$p^{\hol}_K$\@, $K\subset\man{M}$ compact, describing the
$\CO^{\hol}$-topology for $\sections[\hol]{\man{E}}$ as in
Section~\ref{subsec:COhol-vb}\@.

Let us get started with the definitions.
\begin{definition}
Let $\map{\pi}{\man{E}}{\man{M}}$ be an holomorphic vector bundle and let
$\tdomain\subset\real$ be an interval.  A \defn{Carath\'eodory section of
class $\C^{\hol}$} of $\man{E}$ is a map
$\map{\xi}{\tdomain\times\man{M}}{\man{E}}$ with the following properties:
\begin{compactenum}[(i)]
\item $\xi(t,z)\in\man{E}_z$ for each $(t,z)\in\tdomain\times\man{M}$\@;
\item for each $t\in\tdomain$\@, the map $\map{\xi_t}{\man{M}}{\man{E}}$
defined by $\xi_t(z)$ is of class $\C^{\hol}$\@;
\item for each $z\in\man{M}$\@, the map $\map{\xi^z}{\tdomain}{\man{E}}$
defined by $\xi^z(t)=\xi(t,z)$ is Lebesgue measurable.
\end{compactenum}
We shall call $\tdomain$ the \defn{time-domain} for the section.  By
$\CFsections[\hol]{\tdomain;\man{E}}$ we denote the set of Carath\'eodory
sections of class $\C^{\hol}$ of $\man{E}$\@.\oprocend
\end{definition}

The associated notions for time-dependent sections compatible with the
$\CO^{\hol}$-topology are as follows.
\begin{definition}
Let $\map{\pi}{\man{E}}{\man{M}}$ be an holomorphic vector bundle and let
$\tdomain\subset\real$ be an interval.  A Carath\'eodory section
$\map{\xi}{\tdomain\times\man{M}}{\man{E}}$ of class $\C^{\hol}$ is
\begin{compactenum}[(i)]
\item \defn{locally integrally $\C^{\hol}$-bounded} if, for every compact set
$K\subset\man{M}$\@, there exists $g\in\Lloc^1(\tdomain;\realnn)$ such that
\begin{equation*}
\dnorm{\xi(t,z)}_{\metric}\le g(t),\qquad(t,z)\in\tdomain\times K
\end{equation*}
and is
\item \defn{locally essentially $\C^{\hol}$-bounded} if, for every compact
set $K\subset\man{M}$\@, there exists $g\in\Lloc^\infty(\tdomain;\realnn)$
such that
\begin{equation*}
\dnorm{\xi(t,z)}_{\metric}\le g(t),\qquad (t,z)\in\tdomain\times K.
\end{equation*}
\end{compactenum}
The set of locally integrally $\C^{\hol}$-bounded sections of $\man{E}$ with
time-domain $\tdomain$ is denoted by $\LIsections[\hol]{\tdomain,\man{E}}$
with time-domain $\tdomain$ is denoted by
and the set of locally essentially $\C^{\hol}$-bounded sections of $\man{E}$
$\LBsections[\hol]{\tdomain;\man{E}}$\@.\oprocend
\end{definition}

As with smooth sections, the preceding definitions admit topological
characterisations, now using the $\CO^{\hol}$-topology for
$\sections[\hol]{\man{E}}$\@.
\begin{theorem}\label{the:Chol-tvsec}
Let\/ $\map{\pi}{\man{E}}{\man{M}}$ be an holomorphic vector bundle and let\/
$\tdomain\subset\real$ be an interval.  For a map\/
$\map{\xi}{\tdomain\times\man{M}}{\man{E}}$ satisfying\/
$\xi(t,z)\in\man{E}_z$ for each\/ $(t,z)\in\tdomain\times\man{M}$\@, the
following two statements are equivalent:
\begin{compactenum}[(i)]
\item \label{pl:CholCF1} $\xi\in\CFsections[\hol]{\tdomain;\man{E}}$\@;
\item \label{pl:CholCF2} the map\/ $\tdomain\ni
t\mapsto\xi_t\in\sections[\hol]{\man{E}}$ is measurable,\savenum
\end{compactenum}
the following two statements are equivalent:
\begin{compactenum}[(i)]\resumenum
\item \label{pl:CholCF3} $\xi\in\LIsections[\hol]{\tdomain;\man{E}}$\@;
\item \label{pl:CholCF4} the map\/ $\tdomain\ni
t\mapsto\xi_t\in\sections[\hol]{\man{E}}$ is measurable and locally Bochner
integrable,\savenum
\end{compactenum}
and the following two statements are equivalent:
\begin{compactenum}[(i)]\resumenum
\item \label{pl:CholCF5} $\xi\in\LBsections[\hol]{\tdomain;\man{E}}$\@;
\item \label{pl:CholCF6} the map\/ $\tdomain\ni
t\mapsto\xi_t\in\sections[\hol]{\man{E}}$ is measurable and locally
essentially von Neumann bounded.
\end{compactenum}
\begin{proof}
\eqref{pl:CholCF1}$\iff$\eqref{pl:CholCF2} For $z\in\man{M}$ and
$\alpha_z\in\dual{\man{E}}_z$\@, define
$\map{\ev_{\alpha_z}}{\sections[\hol]{\man{E}}}{\complex}$ by
$\ev_{\alpha_z}(\xi)=\natpair{\alpha_z}{\xi(z)}$\@.  It is easy to show that
$\ev_{\alpha_z}$ is continuous and that the set of continuous functionals
$\ev_{\alpha_z}$\@, $\alpha_z\in\dual{\man{E}}_z$\@, is point separating.
Since $\sections[\hol]{\man{E}}$ is a Suslin space
by~$\CO^{\hol}$-\ref{enum:COhol-suslin}\@, this part of the theorem follows
in the same manner as the corresponding part of Theorem~\ref{the:Cinfty-tvsec}\@.

\eqref{pl:CholCF3}$\iff$\eqref{pl:CholCF4} Since $\sections[\hol]{\man{E}}$ is
complete and separable (by properties~$\CO^{\hol}$-\ref{enum:COhol-complete}
and~$\CO^{\hol}$-\ref{enum:COhol-separable}), the arguments from the
corresponding part of Theorem~\ref{the:Cinfty-tvsec} apply here.

\eqref{pl:CholCF5}$\iff$\eqref{pl:CholCF6} We recall our discussion of von
Neumann bounded sets in locally convex topological vector spaces preceding
Lemma~\ref{lem:COinftybdd} above.  With this in mind and using
Lemma~\ref{lem:CObdd}\@, this part of the proposition follows immediately.
\end{proof}
\end{theorem}

Since holomorphic vector bundles are smooth vector bundles (indeed, real
analytic vector bundles), we have natural inclusions
\begin{equation}\label{eq:CholsubsetCinfty}
\LIsections[\hol]{\tdomain;\man{E}}\subset
\CFsections[\infty]{\tdomain;\man{E}},\qquad
\LBsections[\hol]{\tdomain;\man{E}}\subset
\CFsections[\infty]{\tdomain;\man{E}}.
\end{equation}
Moreover, by Proposition~\ref{prop:Cauchyest} we have the following.
\begin{proposition}\label{prop:CholsubsetCinfty}
For an holomorphic vector bundle\/ $\map{\pi}{\man{E}}{\man{M}}$ and an
interval\/ $\tdomain$\@, the inclusions~\eqref{eq:CholsubsetCinfty} actually
induce inclusions
\begin{equation*}
\LIsections[\hol]{\tdomain;\man{E}}\subset
\LIsections[\infty]{\tdomain;\man{E}},
\qquad\LBsections[\hol]{\tdomain;\man{E}}\subset
\LBsections[\infty]{\tdomain;\man{E}}.
\end{equation*}
\end{proposition}

Note that Theorem~\ref{the:Chol-tvsec} applies, in particular, to vector
fields and functions, giving the classes $\CF[\hol]{\tdomain;\man{M}}$\@,
$\LIC[\hol]{\tdomain;\man{M}}$\@, and $\LBC[\hol]{\tdomain;\man{M}}$ of
functions and the classes $\CFsections[\hol]{\tdomain;\tb{\man{M}}}$\@,
$\LIsections[\hol]{\tdomain;\tb{\man{M}}}$\@, and
$\LBsections[\hol]{\tdomain;\tb{\man{M}}}$ of vector fields.  Unlike in the
smooth case preceding and the real analytic case following, there is, in
general, not an equivalent weak-$\sL$ version of the preceding definitions
and results.  This is because our Theorem~\ref{the:COhol-weak} on the
equivalence of the $\CO^{\hol}$-topology and the corresponding weak-$\sL$
topology holds only on Stein manifolds.  Let us understand the consequences
of this with what we are doing here via an example.
\begin{example}\label{eg:LIhol!stein}
Let $\man{M}$ be a compact holomorphic manifold.
By~\cite[Corollary~IV.1.3]{KF/HG:02}\@, the only holomorphic functions on
$\man{M}$ are the locally constant functions.  Therefore, since $\partial
f=0$ for every $f\in\func[\hol]{\man{M}}$\@, a literal application of the
definition shows that, were we to make weak-$\sL$ characterisations of vector
fields,~\ie~give their properties by ascribing those properties to the
functions obtained after Lie differentiation, we would have $\CFsections[\hol]{\tdomain;\tb{\man{M}}}$\@, and,
therefore, also $\LIsections[\hol]{\tdomain;\tb{\man{M}}}$ and
$\LBsections[\hol]{\tdomain;\tb{\man{M}}}$\@, consisting of \emph{all} maps
$\map{X}{\tdomain\times\man{M}}{\tb{\man{M}}}$ satisfying
$X(t,z)\in\tb[z]{\man{M}}$ for all $z\in\man{M}$\@.  This is not a very
useful class of vector fields.\oprocend
\end{example}

The following result summarises the various ways of verifying the
measurability, integrability, and boundedness of holomorphic time-varying
vector fields, taking into account that the preceding example necessitates
that we restrict our consideration to Stein manifolds.
\begin{theorem}
Let\/ $\man{M}$ be a Stein manifold, let\/ $\tdomain\subset\real$ be a
time-domain, and let\/ $\map{X}{\tdomain\times\man{M}}{\tb{\man{M}}}$ have
the property that\/ $X_t$ is an holomorphic vector field for each\/
$t\in\tdomain$\@.  Then the following statements are equivalent:
\begin{compactenum}[(i)]
\item $t\mapsto X_t$ is measurable;
\item $t\mapsto\lieder{X_t}{f}$ is measurable for every\/
$f\in\func[\hol]{\man{M}}$\@;
\item $t\mapsto\ev_z\scirc X_t$ is measurable for every\/ $z\in\man{M}$\@;
\item $t\mapsto\ev_z\scirc\lieder{X_t}{f}$ is measurable for every\/
$f\in\func[\hol]{\man{M}}$ and every\/ $z\in\man{M}$\@,\savenum
\end{compactenum}
the following two statements are equivalent:
\begin{compactenum}[(i)]\resumenum
\item $t\mapsto X_t$ is locally Bochner integrable;
\item $t\mapsto\lieder{X_t}{f}$ is locally Bochner integrable for every\/
$f\in\func[\hol]{\man{M}}$\@,\savenum
\end{compactenum}
and the following two statements are equivalent:
\begin{compactenum}[(i)]\resumenum
\item $t\mapsto X_t$ is locally essentially von Neumann bounded;
\item $t\mapsto\lieder{X_t}{f}$ is locally essentially von Neumann bounded
for every\/ $f\in\func[\hol]{\man{M}}$\@.
\end{compactenum}
\begin{proof}
This follows from Theorem~\ref{the:Chol-tvsec}\@, along with Corollary~\ref{cor:COweak}\@.
\end{proof}
\end{theorem}

Now we consider flows for the class of time-varying holomorphic vector fields
defined above.  Let $X\in\LIsections[\hol]{\tdomain;\tb{\man{M}}}$\@.
According to Proposition~\ref{prop:CholsubsetCinfty}\@, we can define the
flow of $X$ just as in the real case, and we shall continue to use the
notation $D_X\subset\tdomain\times\tdomain\times\man{M}$\@,
$\flow{X}{t}[t_0]$\@, and $\map{\flow*{X}}{D_X}{\man{M}}$ as in the smooth
case.  The following result provides the attributes of the flow in the
holomorphic case.  This result follows easily from the constructions in the
usual existence and uniqueness theorem for ordinary differential equations,
but we could not find the result explicitly in the literature for measurable
time-dependence.  Thus we provide the details here.
\begin{theorem}\label{the:Chol-flows}
Let\/ $\man{M}$ be an holomorphic manifold, let\/ $\tdomain$ be an interval,
and let\/ $X\in\LIsections[\hol]{\tdomain;\tb{\man{M}}}$\@.  Then there
exist a subset\/ $D_X\subset\tdomain\times\tdomain\times\man{M}$ and a map\/
$\map{\flow*{X}}{D_X}{\man{M}}$ with the following properties for each\/ $(t_0,z_0)\in\tdomain\times\man{M}$\@:
\begin{compactenum}[(i)]
\item the set
\begin{equation*}
\tdomain_X(t_0,z_0)=\setdef{t\in\tdomain}{(t,t_0,z_0)\in D_X}
\end{equation*}
is an interval;
\item there exists a locally absolutely continuous curve\/ $t\mapsto\xi(t)$
satisfying
\begin{equation*}
\xi'(t)=X(t,\xi(t)),\quad\xi(t_0)=z_0,
\end{equation*}
for almost all\/ $t\in|t_0,t_1|$ if and only if\/
$t_1\in\tdomain_X(t_0,z_0)$\@;
\item $\deriv{}{t}\flow*{X}(t,t_0,z_0)=X(t,\flow*{X}(t,t_0,z_0))$ for almost all\/ $t\in\tdomain_X(t_0,z_0)$\@;
\item for each\/ $t\in\tdomain$ for which\/ $(t,t_0,z_0)\in D_X$\@, there
exists a neighbourhood\/ $\nbhd{U}$ of\/ $z_0$ such that the mapping\/
$z\mapsto\flow*{X}(t,t_0,z)$ is defined and of class\/ $\C^{\hol}$ on\/
$\nbhd{U}$\@.
\end{compactenum}
\begin{proof}
Given Proposition~\ref{prop:CholsubsetCinfty}\@, the only part of the theorem
that does not follow from Theorem~\ref{the:Cinfty-flows} is the holomorphic
dependence on initial conditions.  This is a local assertion, so we let
$(\nbhd{U},\phi)$ be an holomorphic chart for $\man{M}$ with coordinates
denoted by $(z^1,\dots,z^n)$\@.  We denote by
$\map{\vect{X}}{\tdomain\times\phi(\nbhd{U})}{\complex^n}$ the local
representative of $X$\@.  By Proposition~\ref{prop:CholsubsetCinfty}\@, this
local representative is locally integrally $\C^\infty$-bounded.  To prove
holomorphicity of the flow, we recall the construction for the existence and
uniqueness theorem for the solutions of the initial value problem
\begin{equation*}
\dot{\vect{\gamma}}(t)=\vect{X}(t,\vect{\gamma}(t)),\qquad
\vect{\gamma}(t_0)=\vect{z},
\end{equation*}
see~\cite[\eg][\S1.2]{FS/HvdM:00}\@.  On some suitable product domain
$\tdomain'\times\oball{r}{\vect{z}_0}$ (the ball being contained in
$\phi(\nbhd{U})\subset\complex^n$) we denote by
$\mappings[0]{\tdomain'\times\oball{r}{\vect{z}_0}}{\complex^n}$ the Banach
space of continuous mappings with the
$\infty$-norm~\cite[Theorem~7.9]{EH/KS:75}\@.  We define an operator
\begin{equation*}
\map{\Phi}{\mappings[0]{\tdomain'\times\oball{r}{\vect{z}_0}}{\complex^n}}
{\mappings[0]{\tdomain'\times\oball{r}{\vect{z}_0}}{\complex^n}}
\end{equation*}
by
\begin{equation*}
\Phi(\vect{\gamma})(t,\vect{z})=\vect{z}+
\int_{t_0}^t\vect{X}(s,\vect{\gamma}(s,\vect{z}))\,\d{s}.
\end{equation*}
One shows that this mapping, with domains suitably defined, is a contraction
mapping, and so, by iterating the mapping, one constructs a sequence in
$\mappings[0]{\tdomain'\times\oball{r}{\vect{z}_0}}{\complex^n}$ converging
to a fixed point, and the fixed point, necessarily satisfying
\begin{equation*}
\vect{\gamma}(t,\vect{z})=\vect{z}+
\int_{t_0}^t\vect{X}(s,\vect{\gamma}(s,\vect{z}))\,\d{s}
\end{equation*}
and $\vect{\gamma}(t_0,\vect{z})=\vect{z}$\@, has the property that
$\vect{\gamma}(t,\vect{z})=\flow*{\vect{X}}(t,t_0,\vect{z})$\@.

Let us consider the sequence one constructs in this procedure.  We define
$\vect{\gamma}_0\in
\mappings[0]{\tdomain'\times\oball{r}{\vect{z}_0}}{\complex^n}$ by
$\vect{\gamma}_0(t,\vect{z})=\vect{z}$\@.  Certainly $\vect{\gamma}_0$ is
holomorphic in $\vect{z}$\@.  Now define $\vect{\gamma}_1\in
\mappings[0]{\tdomain'\times\oball{r}{\vect{z}_0}}{\complex^n}$ by
\begin{equation*}
\vect{\gamma}_1(t,\vect{z})=\Phi(\vect{\gamma}_0)=
\vect{z}+\int_{t_0}^t\vect{X}(s,\vect{z})\,\d{s}.
\end{equation*}
Since $\vect{X}\in\LIsections[\hol]{\tdomain';\tb{\oball{r}{\vect{z}_0}}}$\@,
we have
\begin{equation*}
\pderiv{}{\bar{z}^j}\vect{\gamma}_1(t,\vect{z})=
\pderiv{}{\bar{z}^j}\vect{z}+\int_{t_0}^t\pderiv{}{\bar{z}^j}
\vect{X}(s,\vect{\gamma}_0(s,\vect{z}))\,\d{s}=0,\qquad j\in\{1,\dots,n\},
\end{equation*}
swapping the derivative and the integral by the Dominated Convergence
Theorem~\cite[Theorem~16.11]{JJ:05} (also noting by
Proposition~\ref{prop:CholsubsetCinfty} that derivatives of $\vect{X}$ are
bounded by an integrable function).  Thus $\vect{\gamma}_1$ is holomorphic
for each fixed $t\in\tdomain'$\@.  By iterating with $t$ fixed, we have a
sequence $\ifam{\vect{\gamma}_{j,t}}_{j\in\integernn}$ of holomorphic
mappings from $\oball{r}{\vect{z}_0}$ converging uniformly to the function
$\vect{\gamma}$ that describes how the solution at time $t$ depends on the
initial condition $\vect{z}$\@.  The limit function is necessarily
holomorphic~\cite[page~5]{RCG:90a}\@.
\end{proof}
\end{theorem}

\subsection{The real analytic case}\label{subsec:analytic-tvvf}

Let us now turn to describing real analytic time-varying sections.  We thus
will consider a real analytic vector bundle $\map{\pi}{\man{E}}{\man{M}}$
with $\nabla^0$ a real analytic linear connection on $\man{E}$\@, $\nabla$ a
real analytic affine connection on $\man{M}$\@, $\metric_0$ a real analytic
fibre metric on $\man{E}$\@, and $\metric$ a real analytic Riemannian metric
on $\man{M}$\@.  This defines the seminorms $p^\omega_{K,\vect{a}}$\@,
$K\subset\man{M}$ compact, $\vect{a}\in\csd(\integernn;\realp)$\@, describing
the $\C^\omega$-topology as in Theorem~\ref{the:Comega-seminorms}\@.
\begin{definition}
Let $\map{\pi}{\man{E}}{\man{M}}$ be a real analytic vector bundle and let
$\tdomain\subset\real$ be an interval.  A \defn{Carath\'eodory section of
class $\C^\omega$} of $\man{E}$ is a map
$\map{\xi}{\tdomain\times\man{M}}{\man{E}}$ with the following properties:
\begin{compactenum}[(i)]
\item $\xi(t,x)\in\man{E}_x$ for each $(t,x)\in\tdomain\times\man{M}$\@;
\item for each $t\in\tdomain$\@, the map $\map{\xi_t}{\man{M}}{\man{E}}$
defined by $\xi_t(x)$ is of class $\C^\omega$\@;
\item for each $x\in\man{M}$\@, the map $\map{\xi^x}{\tdomain}{\man{E}}$
defined by $\xi^x(t)=\xi(t,x)$ is Lebesgue measurable.
\end{compactenum}
We shall call $\tdomain$ the \defn{time-domain} for the section.  By
$\CFsections[\omega]{\tdomain;\man{E}}$ we denote the set of Carath\'eodory
sections of class $\C^\omega$ of $\man{E}$\@.\oprocend
\end{definition}

Now we turn to placing restrictions on the time-dependence to allow us to do
useful things.
\begin{definition}\label{def:LILBomega}
Let $\map{\pi}{\man{E}}{\man{M}}$ be a real analytic vector bundle and let
$\tdomain\subset\real$ be an interval.  A Carath\'eodory section
$\map{\xi}{\tdomain\times\man{M}}{\man{E}}$ of class $\C^\omega$ is
\begin{compactenum}[(i)]
\item \defn{locally integrally $\C^\omega$-bounded} if, for every compact set
$K\subset\man{M}$ and every $\vect{a}\in\csd(\integernn;\realp)$\@, there
exists $g\in\Lloc^1(\tdomain;\realnn)$ such that
\begin{equation*}
a_0a_1\cdots a_m\dnorm{j_m\xi_t(x)}_{\ol{\metric}_m}\le
g(t),\qquad(t,x)\in\tdomain\times K,\ m\in\integernn,
\end{equation*}
and is
\item \defn{locally essentially $\C^\omega$-bounded} if, for every compact
set $K\subset\man{M}$ and every $\vect{a}\in\csd(\integernn;\realp)$\@, there
exists $g\in\Lloc^\infty(\tdomain;\realnn)$ such that
\begin{equation*}
a_0a_1\cdots a_m\dnorm{j_m\xi_t(x)}_{\ol{\metric}_m}\le g(t),
\qquad(t,x)\in\tdomain\times K,\ m\in\integernn.
\end{equation*}
\end{compactenum}
The set of locally integrally $\C^\omega$-bounded sections of $\man{E}$ with
time-domain $\tdomain$ is denoted by $\LIsections[\omega]{\tdomain,\man{E}}$
and the set of locally essentially $\C^\omega$-bounded sections of $\man{E}$
with time-domain $\tdomain$ is denoted by
$\LBsections[\omega]{\tdomain;\man{E}}$\@.\oprocend
\end{definition}

As with smooth and holomorphic sections, the preceding definitions admit
topological characterisations.
\begin{theorem}\label{the:Comega-tvsec}
Let\/ $\map{\pi}{\man{E}}{\man{M}}$ be a real analytic manifold and let\/
$\tdomain\subset\real$ be an interval.  For a map\/
$\map{\xi}{\tdomain\times\man{M}}{\man{E}}$ satisfying\/
$\xi(t,x)\in\man{E}_x$ for each\/ $(t,x)\in\tdomain\times\man{M}$\@, the
following two statements are equivalent:
\begin{compactenum}[(i)]
\item \label{pl:ComegaCF1} $\xi\in\CFsections[\omega]{\tdomain;\man{E}}$\@;
\item \label{pl:ComegaCF2} the map\/ $\tdomain\ni t\mapsto\xi_t\in\sections[\omega]{\man{E}}$ is measurable,\savenum
\end{compactenum}
the following two statements are equivalent:
\begin{compactenum}[(i)]\resumenum
\item \label{pl:ComegaCF3} $\xi\in\LIsections[\omega]{\tdomain;\man{E}}$\@;
\item \label{pl:ComegaCF4} the map\/ $\tdomain\ni
t\mapsto\xi_t\in\sections[\omega]{\man{E}}$ is measurable and locally Bochner
integrable,\savenum
\end{compactenum}
and the following two statements are equivalent:
\begin{compactenum}[(i)]\resumenum
\item \label{pl:ComegaCF5} $\xi\in\LBsections[\omega]{\tdomain;\man{E}}$\@;
\item \label{pl:ComegaCF6} the map\/ $\tdomain\ni
t\mapsto\xi_t\in\sections[\omega]{\man{E}}$ is measurable and locally
essentially von Neumann bounded.
\end{compactenum}
\begin{proof}
Just as in the smooth case in Theorem~\ref{the:Cinfty-tvsec}\@, this is
deduced from the following facts:~(1)~evaluation maps $\ev_{\alpha_x}$\@,
$\alpha_x\in\dual{\man{E}}$\@, are continuous and point
separating;~(2)~$\sections[\omega]{\man{E}}$ is a Suslin space
(property~$\C^\omega$-\ref{enum:Comega-suslin});%
~(3)~$\sections[\omega]{\man{E}}$ is complete and separable
(properties~$\C^\omega$-\ref{enum:Comega-complete}
and~$\C^\omega$-\ref{enum:Comega-separable};~(4)~we understand von Neumann
bounded subsets of $\sections[\omega]{\man{E}}$ by
Lemma~\ref{lem:Comegabdd}\@.
\end{proof}
\end{theorem}

Note that Theorem~\ref{the:Comega-tvsec} applies, in particular, to vector
fields and functions, giving the classes $\CF[\omega]{\tdomain;\man{M}}$\@,
$\LIC[\omega]{\tdomain;\man{M}}$\@, and $\LBC[\omega]{\tdomain;\man{M}}$ of
functions, and the classes $\CFsections[\omega]{\tdomain;\tb{\man{M}}}$\@,
$\LIsections[\omega]{\tdomain;\tb{\man{M}}}$\@, and
$\LBsections[\omega]{\tdomain;\tb{\man{M}}}$ of vector fields.  The following
result then summarises the various ways of verifying the measurability,
integrability, and boundedness of real analytic time-varying vector fields.
\begin{theorem}\label{the:ra-td-summary}
Let\/ $\man{M}$ be a real analytic manifold, let\/ $\tdomain\subset\real$ be
a time-domain, and let\/ $\map{X}{\tdomain\times\man{M}}{\tb{\man{M}}}$ have
the property that\/ $X_t$ is a real analytic vector field for each\/
$t\in\tdomain$\@.  Then the following statements are equivalent:
\begin{compactenum}[(i)]
\item $t\mapsto X_t$ is measurable;
\item $t\mapsto\lieder{X_t}{f}$ is measurable for every\/ $f\in\func[\omega]{\man{M}}$\@;
\item $t\mapsto\ev_x\scirc X_t$ is measurable for every\/ $x\in\man{M}$\@;
\item $t\mapsto\ev_x\scirc\lieder{X_t}{f}$ is measurable for every\/
$f\in\func[\omega]{\man{M}}$ and every\/ $x\in\man{M}$\@,\savenum
\end{compactenum}
the following two statements are equivalent:
\begin{compactenum}[(i)]\resumenum
\item $t\mapsto X_t$ is locally Bochner integrable;
\item $t\mapsto\lieder{X_t}{f}$ is locally Bochner integrable for every\/
$f\in\func[\omega]{\man{M}}$\@,\savenum
\end{compactenum}
and the following two statements are equivalent:
\begin{compactenum}[(i)]\resumenum
\item $t\mapsto X_t$ is locally essentially bounded;
\item $t\mapsto\lieder{X_t}{f}$ is locally essentially bounded in the von
Neumann bornology for every\/ $f\in\func[\omega]{\man{M}}$\@.
\end{compactenum}
\begin{proof}
This follows from Theorem~\ref{the:Comega-tvsec}\@, along with
Corollary~\ref{cor:Comegaweak}\@.
\end{proof}
\end{theorem}

Let us verify that real analytic time-varying sections have the expected
relationship to their smooth brethren.
\begin{proposition}\label{prop:ComegasubsetCinfty}
For a real analytic vector bundle\/ $\map{\pi}{\man{E}}{\man{M}}$ and an
interval\/ $\tdomain$\@, we have
\begin{gather*}
\LIsections[\omega]{\tdomain;\man{E}}\subset
\LIsections[\infty]{\tdomain;\man{E}},\qquad
\LBsections[\omega]{\tdomain;\man{M}}\subset
\LBsections[\infty]{\tdomain;\man{M}}.
\end{gather*}
\begin{proof}
It is obvious that real analytic Carath\'eodory sections are smooth
Carath\'eodory sections.

Let us verify only that $\LIsections[\omega]{\tdomain;\man{E}}\subset
\LIsections[\infty]{\tdomain;\man{E}}$\@, as the essentially bounded case
follows in the same manner.  We let $K\subset\man{M}$ be compact and let
$m\in\integernn$\@.  Choose (arbitrarily)
$\vect{a}\in\csd(\integernn;\realp)$\@.  Then, if
$\xi\in\LIsections[\omega]{\tdomain;\man{E}}$\@, there exists
$g\in\Lloc^1(\tdomain;\realnn)$ such that
\begin{equation*}
a_0a_1\cdots a_m\dnorm{j_m\xi_t(x)}_{\ol{\metric}_m}\le g(t),
\qquad x\in K,\ t\in\tdomain,\ m\in\integernn.
\end{equation*}
Thus, taking $g_{\vect{a},m}\in\Lloc^1(\tdomain;\realnn)$ defined by
\begin{equation*}
g_{\vect{a},m}(t)=\frac{1}{a_0a_1\cdots a_m}g(t),
\end{equation*}
we have
\begin{equation*}
\dnorm{j_m\xi_t(x)}_{\ol{\metric}_m}\le g_{\vect{a},m}(t),\qquad
x\in K,\ t\in\tdomain
\end{equation*}
showing that $\xi\in\LIsections[\infty]{\tdomain;\man{E}}$\@.
\end{proof}
\end{proposition}

Having understood the comparatively simple relationship between real analytic
and smooth time-varying sections, let us consider the correspondence between
real analytic and holomorphic time-varying sections.  First, note that if
$\tdomain\subset\real$ is an interval and if $\ol{\nbhd{U}}\in\sN_{\man{M}}$
is a neighbourhood of $\man{M}$ in a complexification $\ol{\man{M}}$\@, then
we have an inclusion
\begin{equation*}
\mapdef{\rho_{\ol{\nbhd{U}},\man{M}}}
{\CFsections[\hol,\real]{\tdomain;\ol{\man{E}}|\ol{\nbhd{U}}}}
{\CFsections[\omega]{\tdomain;\man{E}}}
{\ol{\xi}}{\ol{\xi}|\man{M}.}
\end{equation*}
(Here the notation
$\CFsections[\hol,\real]{\tdomain;\ol{\man{E}}|\ol{\nbhd{U}}}$ refers to
those Carath\'eodory sections that are real when restricted to
$\man{M}$\@,~\cf~the constructions of Section~\ref{subsubsec:hologerms}\@.)
However, this inclusion does not characterise all real analytic
Carath\'eodory sections, as the following example shows.
\begin{example}\label{eg:CFomeganothol}
Let $\tdomain$ be any interval for which $0\in\interior(\tdomain)$\@.  We
consider the real analytic Carath\'edory function on $\real$ with time-domain
$\tdomain$ defined by
\begin{equation*}
f(t,x)=\begin{cases}\frac{t^2}{t^2+x^2},&t\not=0,\\0,&t=0.\end{cases}
\end{equation*}
It is clear that $x\mapsto f(t,x)$ is real analytic for every $t\in\tdomain$
and that $t\mapsto f(t,x)$ is measurable for every $x\in\real$\@.  We claim
that there is no neighbourhood $\ol{\nbhd{U}}\subset\complex$ of
$\real\subset\complex$ such that $f$ is the restriction to $\real$ of an
holomorphic Carath\'eodory function on $\ol{\nbhd{U}}$\@.  Indeed, let
$\ol{\nbhd{U}}\subset\complex$ be a neighbourhood of $\real$ and choose
$t\in\realp$ sufficiently small that $\cdisk[]{t}{0}\subset\ol{\nbhd{U}}$\@.
Note that $f_t\colon x\mapsto\frac{1}{1+(x/t)^2}$ does not admit an
holomorphic extension to any open set containing $\cdisk[]{t}{0}$ since the
radius of convergence of $z\mapsto\frac{1}{1+(z/t)^2}$ is $t$\@,~\cf~the
discussion at the beginning of Section~\ref{sec:analytic-topology}\@.  Note
that our construction actually shows that in no neighbourhood of
$(0,0)\in\real\times\real$ is there an holomorphic extension of
$f$\@.\oprocend
\end{example}

Fortunately, the example will not bother us, although it does serve to
illustrate that the following result is not immediate.
\begin{theorem}\label{the:Comega->Chol}
Let\/ $\map{\pi}{\man{E}}{\man{M}}$ be a real analytic vector bundle with
complexification\/ $\map{\ol{\pi}}{\ol{\man{E}}}{\ol{\man{M}}}$\@, and let\/
$\tdomain$ be a time-domain.  For a map\/
$\map{\xi}{\tdomain\times\man{M}}{\man{E}}$ satisfying\/
$\xi(t,x)\in\man{E}_x$ for every\/ $(t,x)\in\tdomain\times\man{M}$\@, the
following statements hold:
\begin{compactenum}[(i)]
\item \label{pl:Comega->Chol1} if\/
$\xi\in\LIsections[\omega]{\tdomain;\man{E}}$\@, then, for each\/
$(t_0,x_0)\in\tdomain\times\man{M}$ and each bounded subinterval\/
$\tdomain'\subset\tdomain$ containing\/ $t_0$\@, there exist a
neighbourhood\/ $\ol{\nbhd{U}}$ of\/ $x_0$ in\/ $\ol{\man{M}}$ and\/
$\ol{\xi}\in\sections[\hol]{\tdomain';\ol{\man{E}}|\ol{\nbhd{U}}}$ such
that\/ $\ol{\xi}(t,x)=\xi(t,x)$ for each\/ $t\in\tdomain'$ and\/
$x\in\ol{\nbhd{U}}\cap\man{M}$\@;
\item \label{pl:Comega->Chol2} if, for each\/ $x_0\in\man{M}$\@, there exist
a neighbourhood\/ $\ol{\nbhd{U}}$ of\/ $x_0$ in\/ $\ol{\man{M}}$ and\/
$\ol{\xi}\in\sections[\hol]{\tdomain;\ol{\man{E}}|\ol{\nbhd{U}}}$ such
that\/ $\ol{\xi}(t,x)=\xi(t,x)$ for each\/ $t\in\tdomain$ and\/
$x\in\ol{\nbhd{U}}\cap\man{M}$\@, then\/ $\xi\in\LIsections[\omega]{\tdomain;\man{E}}$\@.
\end{compactenum}
\begin{proof}
\eqref{pl:Comega->Chol1} We let $\tdomain'\subset\tdomain$ be a bounded
subinterval containing $t_0$ and let $\nbhd{U}$ be a relatively compact
neighbourhood of $x_0$\@.  Let $\ifam{\ol{\nbhd{U}}_j}_{j\in\integerp}$ be a
sequence of neighbourhoods of $\closure(\nbhd{U})$ in $\ol{\man{M}}$ with the
properties that $\closure(\ol{\nbhd{U}}_j)\subset\ol{\nbhd{U}}_{j+1}$ and
that $\cap_{j\in\integerp}\ol{\nbhd{U}}_j=\closure(\nbhd{U})$\@.  We first
note that
\begin{equation*}
\L^1(\tdomain';\sections[\hol,\real]{\ol{\man{E}}|\ol{\nbhd{U}}_j})
\simeq\L^1(\tdomain';\real)\protimes
\sections[\hol,\real]{\ol{\man{E}}|\ol{\nbhd{U}}_j},
\end{equation*}
with $\protimes$ denoting the completed projective tensor
product~\cite[Theorem~III.6.5]{HHS/MPW:99}\@. The theorem of
\citeauthor{HHS/MPW:99} is given for Banach spaces, and they also assert the
validity of this for locally convex spaces; thus we also have
\begin{equation*}
\L^1(\tdomain';\gsections[\hol,\real]
{\closure(\nbhd{U})}{\ol{\man{E}}})\simeq
\L^1(\tdomain';\real)\protimes
\gsections[\hol,\real]{\closure(\nbhd{U})}{\ol{\man{E}}}.
\end{equation*}
In both cases, the isomorphisms are in the category of locally convex
topological vector spaces.  We claim that, with these identifications,
\begin{equation*}
\L^1(\tdomain';\real)\otimes_\pi
\gsections[\hol,\real]{\closure(\nbhd{U})}{\ol{\man{E}}}
\end{equation*}
is the direct limit of the directed system
\begin{equation*}
\ifam{\L^1(\tdomain';\real)\otimes_\pi
\sections[\hol,\real]{\ol{\man{E}}|\ol{\nbhd{U}}_j})}_{j\in\integerp}
\end{equation*}
with the associated mappings $\id\otimes_\pi r_{\closure(\nbhd{U}),j}$\@,
$j\in\integerp$\@, where $r_{\closure(\nbhd{U}),j}$ is defined as
in~\eqref{eq:Gammahol->Ghol}\@.  (Here $\otimes_\pi$ is the uncompleted
projective tensor product).  We, moreover, claim that the direct limit
topology is boundedly retractive, meaning that bounded sets in the direct
limit are contained in and bounded in a single component of the directed
system and, moreover, the topology on the bounded set induced by the
component is the same as that induced by the direct limit.

Results of this sort have been the subject of research in the area of locally
convex topologies, with the aim being to deduce conditions on the structure
of the spaces comprising the directed system, and on the corresponding
mappings (for us, the inclusion mappings and their tensor products with the
identity on $\L^1(\tdomain';\real)$), that ensure that direct limits commute
with tensor product, and that the associated direct limit topology is
boundedly retractive.  We shall make principal use of the results given by
\citet{EMM:97}\@.  To state the arguments with at least a little context, let
us reproduce two conditions used by \citeauthor{EMM:97}\@.
\begin{procnonum*}[Condition (M) of \citet{VSR:70}]
Let $\ifam{\alg{V}_j}_{j\in\integerp}$ be a directed system of locally convex
spaces with strict direct limit $\alg{V}$\@.  The direct limit topology of
$\alg{V}$ satisfies \defn{condition~(M)} if there exists a sequence
$\ifam{\nbhd{O}_j}_{j\in\integerp}$ for which
\begin{compactenum}[(i)]
\item $\nbhd{O}_j$ is a balanced convex neighbourhood of $0\in\alg{V}_j$\@,
\item $\nbhd{O}_j\subset\nbhd{O}_{j+1}$ for each $j\in\integerp$\@, and
\item for every $j\in\integerp$\@, there exists $k\ge j$ such that the
topology induced on $\nbhd{O}_j$ by its inclusion in $\alg{V}_k$ and its
inclusion in $\alg{V}$ agree.\oprocend
\end{compactenum}
\end{procnonum*}

\begin{procnonum*}[Condition~(MO) of \citet{EMM:97}]
Let $\ifam{\alg{V}_j}_{j\in\integerp}$ be a directed system of metrisable
locally convex spaces with strict direct limit $\alg{V}$\@.  Let
$\map{i_{j,k}}{\alg{V}_j}{\alg{V}_k}$ be the inclusion for $k\ge j$ and let
$\map{i_j}{\alg{V}_j}{\alg{V}}$ be the induced map into the direct limit.

Suppose that, for each $j\in\integerp$\@, we have a sequence
$\ifam{p_{j,l}}_{l\in\integerp}$ of seminorms defining the topology of
$\alg{V}_j$ such that $p_{j,l_1}\ge p_{j,l_2}$ if $l_1\ge l_2$\@.  Let
\begin{equation*}
\alg{V}_{j,l}=\alg{V}_j/\setdef{v\in\alg{V}_j}{p_{j,l}(v)=0}
\end{equation*}
and denote by $\hat{p}_{j,l}$ the norm on $\alg{V}_{j,l}$ induced by
$p_{j,l}$~\cite[page~97]{HHS/MPW:99}\@.  Let
$\map{\pi_{j,l}}{\alg{V}_j}{\alg{V}_{j,l}}$ be the canonical projection.  Let
$\ol{\alg{V}}_{j,l}$ be the completion of $\alg{V}_{j,l}$\@.  The family
$\ifam{\ol{\alg{V}}_{j,l}}_{j,l\in\integerp}$ is called a \defn{projective
spectrum} for $\alg{V}_j$\@.  Denote
\begin{equation*}
\nbhd{O}_{j,l}=\setdef{v\in\alg{V}_j}{p_{j,l}(v)\le1}.
\end{equation*}

The direct limit topology of $\alg{V}$ satisfies \defn{condition~(MO)} if
there exists a sequence $\ifam{\nbhd{O}_j}_{j\in\integerp}$ and if, for every
$j\in\integerp$\@, there exists a projective spectrum
$\ifam{\ol{\alg{V}}_{j,l}}_{j,l\in\integerp}$ for $\alg{V}_j$ for which
\begin{compactenum}[(i)]
\item $\nbhd{O}_j$ is a balanced convex neighbourhood of $0\in\alg{V}_j$\@,
\item $\nbhd{O}_j\subset\nbhd{O}_{j+1}$ for each $j\in\integerp$\@, and
\item for every $j\in\integerp$\@, there exists $k\ge j$ such that, for every
$l\in\integerp$\@, there exists\/ $A\in\lin{\alg{V}}{\ol{\alg{V}}_{k,l}}$
satisfying
\begin{equation*}
(\pi_{k,l}\scirc i_{jk}-A\scirc i_j)(\nbhd{O}_j)\subset
\closure(\pi_{k,l}(\alg{O}_{k,l})),
\end{equation*}
the closure on the right being taken in the norm topology of $\ol{\alg{V}}_{k,l}$\@.\oprocend
\end{compactenum}
\end{procnonum*}

With these concepts, we have the following statements.  We let
$\ifam{\alg{V}_j}_{j\in\integerp}$ be a directed system of metrisable locally
convex spaces with strict direct limit $\alg{V}$\@.
\begin{compactenum}
\item \label{enum:M-MO1} If the direct limit topology on $\alg{V}$ satisfies
condition~(MO), then, for any Banach space $\alg{U}$\@,
$\alg{U}\otimes_\pi\alg{V}$ is the direct limit of the directed system
$\ifam{\alg{U}\otimes_\pi\alg{V}_j}_{j\in\integerp}$\@, and the direct limit
topology on $\alg{U}\otimes_\pi\alg{V}$ satisfies
condition~(M)~\cite[Theorem~1.3]{EMM:97}\@.
\item \label{enum:M-MO2} If the spaces $\alg{V}_j$\@, $j\in\integerp$\@, are
nuclear and if the direct limit topology on $\alg{V}$ is regular, then the
direct limit topology on $\alg{V}$ satisfies
condition~(MO)~\cite[Theorem~1.3]{EMM:97}\@.
\item \label{enum:M-MO3} If the direct limit topology on $\alg{V}$ satisfies
condition~(M), then this direct limit topology is boundedly retractive~\cite{JW:95}\@.\savenum
\end{compactenum}
Using these arguments we make the following conclusions.
\begin{compactenum}\resumenum
\item The direct limit topology on
$\gsections[\hol,\real]{\closure(\nbhd{U})}{\ol{\man{E}}}$ satisfies
condition~(MO) (by virtue of assertion~\ref{enum:M-MO2} above and by the
properties of the direct limit topology enunciated in
Section~\ref{subsec:Comega-props}\@, specifically that the direct limit is a
regular direct limit of nuclear Fr\'echet spaces).
\item \label{enum:M-MO5} The space $\L^1(\tdomain';\real)\otimes_\pi
\gsections[\hol,\real]{\closure(\nbhd{U})}{\ol{\man{E}}}$ is the direct limit
of the directed sequence
$\ifam{\L^1(\tdomain';\real)\otimes_\pi\sections[\hol,\real]
{\ol{\man{E}}|\ol{\nbhd{U}}_j}}_{j\in\integerp}$ (by virtue of
assertion~\ref{enum:M-MO1} above).
\item The direct limit topology on
$\L^1(\tdomain';\real)\otimes_\pi\gsections[\hol,\real]
{\closure(\nbhd{U})}{\ol{\man{E}}}$ satisfies condition~(M) (by virtue of
assertion~\ref{enum:M-MO1} above).
\item \label{enum:M-MO7} The direct limit topology on
$\L^1(\tdomain';\real)\otimes_\pi\gsections[\hol,\real]
{\closure(\nbhd{U})}{\ol{\man{E}}}$ is boundedly retractive (by virtue of
assertion~\ref{enum:M-MO3} above).
\end{compactenum}

We shall also need the following lemma.
\begin{prooflemma}
Let\/ $K\subset\man{M}$ be compact.  If\/
$[\ol{\xi}]_K\in\L^1(\tdomain';\gsections[\hol,\real]{K}{\man{E}})$ then
there exists a sequence\/ $\ifam{[\ol{\xi}_k]_K}_{k\in\integerp}$ in\/
$\L^1(\tdomain';\real)\otimes\gsections[\hol,\real]{K}{\man{E}}$ converging
to\/ $[\ol{\xi}]_K$ in the topology of\/
$\L^1(\tdomain';\gsections[\hol]{K}{\ol{\man{E}}})$\@.
\begin{subproof}
Since $\L^1(\tdomain';\gsections[\hol,\real]{K}{\ol{\man{E}}})$ is the
completion of $\L^1(\tdomain';\real)\otimes_\pi\gsections[\hol,\real]
{K}{\ol{\man{E}}}$\@, there exists a \emph{net} $\ifam{[\ol{\xi}_i]_K}_{i\in
I}$ converging to $[\ol{\xi}]$\@, so the conclusion here is that we can
actually find a converging \emph{sequence}\@.

To prove this we argue as follows.  Recall
properties~$\mathscr{G}^{\hol,\real}$-\ref{enum:Ghol-reflexive}
and~$\mathscr{G}^{\hol,\real}$-\ref{enum:Ghol-dual} of
$\gsections[\hol,\real]{K}{\ol{\man{E}}}$\@, indicating that it is reflexive
and its dual is a nuclear Fr\'echet space.  Thus
$\gsections[\hol,\real]{K}{\ol{\man{E}}}$ is the dual of a nuclear Fr\'echet
space.  Also recall from
property~$\mathscr{G}^{\hol,\real}$-\ref{enum:Ghol-suslin} that
$\gsections[\hol,\real]{K}{\ol{\man{E}}}$ is a Suslin space.  Now, by
combining \cite[Theorem~7]{GEFT:75} with remark~(1) at the bottom of page~76
of~\cite{GEFT:75} (and being aware that Bochner integrability as defined by
\citeauthor{GEFT:75} is not \emph{a priori} the same as Bochner integrability
as we mean it), there exists a sequence
$\ifam{[\ol{\xi}_k]_K}_{k\in\integerp}$ of simple functions,~\ie~elements of
$\L^1(\tdomain';\real)\otimes \gsections[\hol,\real]{K}{\ol{\man{E}}}$\@,
such that
\begin{equation*}
\lim_{k\to\infty}[\ol{\xi}_k(t)]_K=[\ol{\xi}(t)]_K,\qquad\ae\ t\in\tdomain',
\end{equation*}
(this limit being in the topology of
$\gsections[\hol,\real]{K}{\ol{\man{E}}}$) and
\begin{equation*}
\lim_{k\to\infty}\int_{\tdomain'}([\ol{\xi}(t)]_K-[\ol{\xi}_k(t)]_K)\,\d{t}=0.
\end{equation*}
This implies, by the Dominated Convergence Theorem, that
\begin{equation*}
\lim_{k\to\infty}\int_{\tdomain'}p^\omega_{K,\vect{a}}
([\ol{\xi}(t)]_K-[\ol{\xi}_k(t)]_K)\,\d{t}=0
\end{equation*}
for every $\vect{a}\in\csd(\integernn;\realp)$\@, giving convergence in
\begin{equation*}
\L^1(\tdomain';\gsections[\hol,\real]{K}{\ol{\man{E}}})\simeq
\L^1(\tdomain';\real)\protimes\gsections[\hol,\real]{K}{\man{E}},
\end{equation*}
as desired.
\end{subproof}
\end{prooflemma}

The remainder of the proof is straightforward.  Since
$\xi\in\LIsections[\omega]{\tdomain;\man{E}}$\@, the map
\begin{equation*}
\tdomain'\ni t\mapsto\xi_t\in\sections[\omega]{\man{E}}
\end{equation*}
is an element of $\L^1(\tdomain';\sections[\omega]{\man{E}})$ by
Theorem~\ref{the:Comega-tvsec}\@.  Therefore, if
$[\ol{\xi}]_{\closure(\nbhd{U})}$ is the image of $\xi$ under the natural
mapping from $\sections[\omega]{\man{E}}$ to
$\gsections[\hol,\real]{\closure(\nbhd{U})}{\ol{\man{E}}}$\@, the map
\begin{equation*}
\tdomain'\ni t\mapsto[\ol{\xi(t)}]_{\closure(\nbhd{U})}\in
\gsections[\hol,\real]{\closure(\nbhd{U})}{\ol{\man{E}}}
\end{equation*}
is an element of $\L^1(\tdomain';\gsections[\hol,\real]
{\closure(\nbhd{U})}{\ol{\man{E}}})$\@, since continuous linear maps commute
with integration~\cite[Lemma~1.2]{RB/AD:11}\@.  Therefore, by the Lemma
above, there exists a sequence
$\ifam{[\ol{\xi}_k]_{\closure(\nbhd{U})}}_{k\in\integerp}$ in
$\L^1(\tdomain';\real)\otimes\gsections[\hol,\real]
{\closure(\nbhd{U})}{\ol{\man{E}}}$ that converges to
$[\ol{\xi}]_{\closure(\nbhd{U})}$\@.  By our conclusion~\ref{enum:M-MO5}
above, the topology in which this convergence takes place is the completion
of the direct limit topology associated to the directed system
$\ifam{\L^1(\tdomain';\real)\otimes_\pi\sections[\hol,\real]
{\ol{\man{E}}|\ol{\nbhd{U}}_j}}_{j\in\integerp}$\@.  The direct limit
topology on $\L^1(\tdomain';\real)\otimes_\pi\gsections[\hol,\real]
{\closure(\nbhd{U})}{\ol{\man{E}}}$ is boundedly retractive by our
conclusion~\ref{enum:M-MO7} above.  This is easily seen to imply that the
direct limit topology is sequentially retractive,~\ie~that convergent
sequences are contained in, and convergent in, a component of the direct
limit~\cite{CF:90}\@.  This implies that there exists $j\in\integerp$ such
that the sequence $\ifam{\ol{\xi}_k}_{k\in\integerp}$ converges in
$\L^1(\tdomain';\sections[\hol,\real]{\ol{\man{E}}|\ol{\nbhd{U}}_j})$ and so
converges to a limit $\ol{\eta}$ satisfying
$[\ol{\eta}]_{\closure(\nbhd{U}_j)}=[\ol{\xi}]_{\closure(\nbhd{U}_j)}$\@.
Thus $\ol{\xi}$ can be holomorphically extended to $\ol{\nbhd{U}}_j$\@.  This
completes this part of the proof.

\eqref{pl:Comega->Chol2} Let $K\subset\man{M}$ be compact and let
$\vect{a}\in\csd(\integernn;\realp)$\@.  Let
$\ifam{\ol{\nbhd{U}}_j}_{j\in\integerp}$ be a sequence of neighbourhoods of
$K$ in $\ol{\man{M}}$ such that
$\closure(\ol{\nbhd{U}}_j)\subset\ol{\nbhd{U}}_{j+1}$ and
$K=\cap_{j\in\integerp}\ol{\nbhd{U}}_j$\@.  By hypothesis, for $x\in K$\@,
there is a relatively compact neighbourhood
$\ol{\nbhd{U}}_x\subset\ol{\man{M}}$ of $x$ in $\ol{\man{M}}$ such that there
is an extension
$\ol{\xi}_x\in\LIsections[\hol,\real]{\tdomain;\ol{\man{E}}|\ol{\nbhd{U}}_x}$
of $\xi|(\tdomain\times(\ol{\nbhd{U}}_x\cap\man{M}))$\@.  Let
$x_1,\dots,x_k\in K$ be such that $K\subset\cup_{j=1}^k\ol{\nbhd{U}}_{x_j}$
and let $l\in\integerp$ be sufficiently large that
$\ol{\nbhd{U}}_l\subset\cup_{j=1}^k\ol{\nbhd{U}}_{x_j}$\@, so $\xi$ admits an
holomorphic extension $\ol{\xi}\in\LIsections[\hol,\real]
{\tdomain;\ol{\man{E}}|\ol{\nbhd{U}}_l}$\@.

Now we show that the above constructions imply that
$\xi\in\LIsections[\omega]{\tdomain;\tb{\man{M}}}$\@.  Let
$\ol{g}\in\Lloc^1(\tdomain;\realnn)$ be such that
\begin{equation*}
\dnorm{\ol{\xi}(t,z)}_{\ol{\metric}}\le\ol{g}(t),\qquad
(t,z)\in\tdomain\times\ol{\nbhd{U}}_l.
\end{equation*}
By Proposition~\ref{prop:Cauchyest}\@, there exist $C,r\in\realp$ such that
\begin{equation*}
\dnorm{j_m\xi(t,x)}\le Cr^{-m}\ol{g}(t)
\end{equation*}
for all $m\in\integernn$\@, $t\in\tdomain$\@, and $x\in K$\@.  Now let
$N\in\integernn$ be such that $a_{N+1}<r$ and let
$g\in\Lloc^1(\tdomain;\realnn)$ be such that
\begin{equation*}
Ca_0a_1\cdots a_mr^{-m}\ol{g}(t)\le g(t)
\end{equation*}
for $m\in\{0,1,\dots,N\}$\@.  Now, if $m\in\{0,1,\dots,N\}$\@, we have
\begin{equation*}
a_0a_1\cdots a_m\dnorm{j_m\xi(t,x)}_{\ol{\metric}_m}
\le a_0a_1\cdots a_mCr^{-m}\ol{g}(t)\le g(t)
\end{equation*}
for $(t,x)\in\tdomain\times K$\@.  If $m>N$ we also have
\begin{align*}
a_0a_1\cdots a_m\dnorm{j_m\xi(t,x)}_{\ol{\metric}_m}
\le&\;a_0a_1\cdots a_Nr^{-N}r^m\dnorm{j_m\xi(t,x)}_{\ol{\metric}_m}\\
\le&\;a_0a_1\cdots a_Nr^{-N}C\ol{g}(t)\le g(t),
\end{align*}
for $(t,x)\in\tdomain\times K$\@, as desired.
\end{proof}
\end{theorem}

Finally, let us show that, according to our definitions, real analytic
time-varying vector fields possess flows depending in a real analytic way on
initial condition.
\begin{theorem}\label{the:Comega-flows}
Let\/ $\man{M}$ be a real analytic manifold, let\/ $\tdomain$ be an interval,
and let\/ $X\in\LIsections[\omega]{\tdomain;\tb{\man{M}}}$\@.  Then there
exist a subset\/ $D_X\subset\tdomain\times\tdomain\times\man{M}$ and a map\/
$\map{\flow*{X}}{D_X}{\man{M}}$ with the following properties for each\/
$(t_0,x_0)\in\tdomain\times\man{M}$\@:
\begin{compactenum}[(i)]
\item the set
\begin{equation*}
\tdomain_X(t_0,x_0)=\setdef{t\in\tdomain}{(t,t_0,x_0)\in D_X}
\end{equation*}
is an interval;
\item there exists a locally absolutely continuous curve\/ $t\mapsto\xi(t)$
satisfying
\begin{equation*}
\xi'(t)=X(t,\xi(t)),\quad\xi(t_0)=x_0,
\end{equation*}
for almost all\/ $t\in|t_0,t_1|$ if and only if\/
$t_1\in\tdomain_X(t_0,x_0)$\@;
\item $\deriv{}{t}\flow*{X}(t,t_0,x_0)=X(t,\flow*{X}(t,t_0,x_0))$ for almost all\/ $t\in\tdomain_X(t_0,x_0)$\@;
\item for each\/ $t\in\tdomain$ for which\/ $(t,t_0,x_0)\in D_X$\@, there
exists a neighbourhood\/ $\nbhd{U}$ of\/ $x_0$ such that the mapping\/
$x\mapsto\flow*{X}(t,t_0,x)$ is defined and of class\/ $\C^\omega$ on\/
$\nbhd{U}$\@.
\end{compactenum}
\begin{proof}
The theorem follows from Theorems~\ref{the:Chol-flows}
and~\ref{the:Comega->Chol}\@, noting that the flow of an holomorphic
extension will leave invariant the real analytic manifold.
\end{proof}
\end{theorem}

\subsection{Mixing regularity hypotheses}\label{subsec:time-mixup}

It is possible to mix regularity conditions for vector fields.  By this we
mean that one can consider vector fields whose dependence on state is more
regular than their joint state/time dependence.  This can be done by
considering $m\in\integernn$\@, $m'\in\{0,\lip\}$\@,
$r\in\integernn\cup\{\infty,\omega\}$\@, and $r'\in\{0,\lip\}$ satisfying
$m+m'<r+r'$\@, and considering vector fields in
\begin{equation*}
\CFsections[r+r']{\tdomain;\tb{\man{M}}}\cap
\LIsections[m+m']{\tdomain;\tb{\man{M}}}\quad\textrm{or}\quad
\CFsections[r+r']{\tdomain;\tb{\man{M}}}\cap
\LBsections[m+m']{\tdomain;\tb{\man{M}}},
\end{equation*}
using the obvious convention that $\infty+\lip=\infty$ and
$\omega+\lip=\omega$\@.  This does come across as quite unnatural in our
framework, and perhaps it is right that it should.  Moreover, because the
$\CO^{m+m'}$-topology for $\sections[r+r']{\tb{\man{M}}}$ will be complete if
and only if $m+m'=r+r'$\@, some of the results above will not translate to
this mixed class of time-varying vector fields: particularly, the results on
Bochner integrability require completeness.  Nonetheless, this mixing of
regularity assumptions is quite common in the literature.  Indeed, this has
\emph{always} been done in the real analytic case, since the notions of
``locally integrally $\C^\omega$-bounded'' and ``locally essentially
$\C^\omega$-bounded'' given in Definition~\ref{def:LILBomega} are being given
for the first time in this paper.

\section{Control systems}\label{sec:systems}

Now, having at hand a thorough accounting of time-varying vector fields, we
turn to the characterisation of classes of control systems.  These classes of
systems will provide us with a precise point of comparison between our
general development of Section~\ref{sec:gcs} and the more common notion of a
control system.  Our system definitions are designed so that the act of
``substituting in a control'' leads to a time-varying vector field of the
sort considered in Section~\ref{sec:time-varying}\@.  This essentially means
that we need for our system vector fields to depend continuously on control
in the appropriate topology.  We note that, in practice, this is generally
not a limitation,~\eg~we show in Example~\ref{eg:control-linear} that
control-affine systems satisfy our conditions.  In cases where it is a
limitation, the definitions and results here can be replaced with suitably
modified versions with less smoothness, and we say a few words about this at
the end of the section.

As we have been doing all along so far, we initially consider separately the
finitely differentiable, Lipschitz, smooth, holomorphic, and real analytic
cases.  Also, the initial part of our discussion is carried out for
parameterised sections of vector bundles (control systems are parameterised
vector fields), as this allows us to handle vector fields and functions
simultaneously, just as we did in
Sections~\ref{sec:smooth-topology}\@,~\ref{sec:holomorphic-topology}\@,
and~\ref{sec:analytic-topology}\@.

When we turn to control systems starting in
Section~\ref{subsec:Cr-systemsI}\@, we merge as much as possible the
consideration of varying degrees of regularity to make clear the fact that,
once the general framework is in place, much of the analysis proceeds along
very similar lines, regardless of regularity.

We also include a brief discussion of differential inclusions since we shall
use these, as well as usual control systems, in understanding the position of
our ``\gcs{}s'' from Section~\ref{sec:gcs} in the existing order of things.

\subsection{Parameterised vector fields}

One can think of a control system as a family of vector fields parameterised
by control, as discussed in Section~\ref{subsubsec:parameterised}\@.  It is
the exact nature of this dependence on the parameter that we discuss in this
section.

\subsubsection{The smooth case}

We begin by discussing parameter dependent smooth sections.  Throughout this
section we will work with a smooth vector bundle
$\map{\pi}{\man{E}}{\man{M}}$ with a linear connection $\nabla^0$ on
$\man{E}$\@, an affine connection $\nabla$ on $\man{M}$\@, a fibre metric
$\metric_0$ on $\man{E}$\@, and a Riemannian metric $\metric$ on $\man{M}$\@.
These define the fibre metrics $\dnorm{\cdot}_{\ol{\metric}_m}$ and the
seminorms $p^\infty_{K,m}$\@, $K\subset\man{M}$ compact, $m\in\integernn$\@,
on $\sections[\infty]{\man{E}}$ as in Section~\ref{subsec:COinfty-vb}\@.
\begin{definition}\label{def:Cinfty-paramsec}
Let $\map{\pi}{\man{E}}{\man{M}}$ be a smooth vector bundle and let $\ts{P}$
be a topological space.  A map $\map{\xi}{\man{M}\times\ts{P}}{\man{E}}$ such
that $\xi(x,p)\in\man{E}_x$ for every $(x,p)\in\man{M}\times\ts{P}$
\begin{compactenum}[(i)]
\item is a \defn{separately parameterised section of class $\C^\infty$} if
\begin{compactenum}[(a)]
\item for each $x\in\man{M}$\@, the map $\map{\xi_x}{\ts{P}}{\man{E}}$ defined
by $\xi_x(p)=\xi(x,p)$ is continuous and
\item for each $p\in\ts{P}$\@, the map $\map{\xi^p}{\man{M}}{\man{E}}$
defined by $\xi^p(x)=\xi(x,p)$ is of class $\C^\infty$\@,
\end{compactenum}
and
\item is a \defn{jointly parameterised section of class $\C^\infty$} if it is
a separately parameterised section of class $\C^\infty$ and if the map
$(x,p)\mapsto j_m\xi^p(x)$ is continuous for every $m\in\integernn$\@.
\end{compactenum}
By $\SPsections[\infty]{\ts{P};\man{E}}$ we denote the set of separately
parameterised sections of $\man{E}$ of class $\C^\infty$ and by
$\JPsections[\infty]{\ts{P};\man{E}}$ we denote the set of jointly
parameterised sections of $\man{E}$ of class $\C^\infty$\@.\oprocend
\end{definition}

It is possible to give purely topological characterisations of this class of sections.
\begin{proposition}\label{prop:paramCinfty}
Let\/ $\map{\pi}{\man{E}}{\man{M}}$ be a smooth vector bundle, let\/ $\ts{P}$
be a topological space, and let\/ $\map{\xi}{\man{M}\times\ts{P}}{\man{E}}$
satisfy\/ $\xi(x,p)\in\man{E}_x$ for every\/ $(x,p)\in\man{M}\times\ts{P}$\@.
Then\/ $\xi\in\JPsections[\infty]{\ts{P};\man{E}}$ if and only if the map\/
$p\mapsto\xi^p\in\sections[\infty]{\man{E}}$ is continuous, where\/
$\sections[\infty]{\man{E}}$ has the\/ $\CO^\infty$-topology.
\begin{proof}
Given $\map{\xi}{\man{M}\times\ts{P}}{\man{E}}$ we let
$\map{\xi_m}{\man{M}\times\ts{P}}{\jet{m}{\man{E}}}$ be the map
$\xi_m(x,p)=j_m\xi^p(x)$\@.  We also denote by
$\map{\sigma_\xi}{\ts{P}}{\sections[\infty]{\man{E}}}$ the map given by
$\sigma_\xi(p)=\xi^p$\@.

First suppose that $\xi_m$ is continuous for every $m\in\integernn$\@.  Let
$K\subset\man{M}$ be compact, let $m\in\integernn$\@, let
$\epsilon\in\realp$\@, and let $p_0\in\ts{P}$\@.  Let $x\in K$ and let
$\nbhd{W}_x$ be a neighbourhood of $\xi_m(x,p_0)$ in $\jet{m}{\man{E}}$ for
which
\begin{equation*}
\nbhd{W}_x\subset\setdef{j_m\eta(x')\in\jet{m}{\man{E}}}
{\dnorm{j_m\eta(x')-\xi_m(x',p_0)}_{\ol{\metric}_m}<\epsilon}.
\end{equation*}
By continuity of $\xi_m$\@, there exist a neighbourhood
$\nbhd{U}_x\subset\man{M}$ of $x$ and a neighbourhood
$\nbhd{O}_x\subset\ts{P}$ of $p_0$ such that
$\xi_m(\nbhd{U}_x\times\nbhd{O}_x)\subset\nbhd{W}_x$\@.  Now let
$x_1,\dots,x_k\in K$ be such that $K\subset\cup_{j=1}^k\nbhd{U}_{x_j}$ and
let $\nbhd{O}=\cap_{j=1}^k\nbhd{O}_{x_j}$\@.  Then, if $p\in\nbhd{O}$ and
$x\in K$\@, we have $x\in\nbhd{U}_{x_j}$ for some $j\in\{1,\dots,k\}$\@.
Thus $\xi_m(x,p)\in\nbhd{W}_{x_j}$\@.  Thus
\begin{equation*}
\dnorm{\xi_m(x,p)-\xi_m(x,p_0)}_{\ol{\metric}_m}<\epsilon.
\end{equation*}
Therefore, taking supremums over $x\in K$\@,
$p^\infty_{K,m}(\sigma_\xi(p)-\sigma_\xi(p_0))\le\epsilon$\@.  As this can be
done for every compact $K\subset\man{M}$ and every $m\in\integernn$\@, we
conclude that $\sigma_\xi$ is continuous.

Next suppose that $\sigma_\xi$ is continuous and let $m\in\integernn$\@.  Let
$(x_0,p_0)\in\man{M}\times\ts{P}$ and let $\nbhd{W}\subset\jet{m}{\man{E}}$
be a neighbourhood of $\xi_m(x_0,p_0)$\@.  Let $\nbhd{U}\subset\man{M}$ be a
relatively compact neighbourhood of $x_0$ and let $\epsilon\in\realp$ be such
that
\begin{equation*}
\pi_m^{-1}(\nbhd{U})\cap\setdef{j_m\eta(x)\in\jet{m}{\man{E}}}
{\dnorm{j_m\eta(x)-\xi_m(x,p_0)}_{\ol{\metric}_m}<\epsilon}\subset\nbhd{W},
\end{equation*}
where $\map{\pi_m}{\jet{m}{\man{E}}}{\man{M}}$ is the projection.  By
continuity of $\sigma_\xi$\@, let $\nbhd{O}\subset\ts{P}$ be a neighbourhood
of $p_0$ such that
$p^\infty_{\closure(\nbhd{U}),m}(\sigma_\xi(p)-\sigma_\xi(p_0))<\epsilon$ for
$p\in\nbhd{O}$\@.  Therefore,
\begin{equation*}
\dnorm{j_m\sigma_\xi(p)(x)-j_m\sigma_\xi(p_0)(x)}_{\ol{\metric}_m}<\epsilon,
\qquad(x,p)\in\closure(\nbhd{U})\times\nbhd{O}.
\end{equation*}
Therefore, if $(x,p)\in\nbhd{U}\times\nbhd{O}$\@, then
$\pi_m(\xi_m(x,p))=x\in\nbhd{U}$ and so $\xi_m(x,p)\in\nbhd{W}$\@, showing
that $\xi_m$ is continuous at $(x_0,p_0)$\@.
\end{proof}
\end{proposition}

Of course, the preceding discussion applies, in particular, to give vector
fields of parameterised class $\C^\infty$ and functions of parameterised
class $\C^\infty$\@.  This gives the spaces $\SPC[\infty]{\ts{P};\man{M}}$
and $\JPC[\infty]{\man{M}}$ of parameterised functions, and the spaces
$\SPsections[\infty]{\ts{P};\tb{\man{M}}}$ and
$\JPsections[\infty]{\ts{P};\tb{\man{M}}}$ of parameterised vector fields.
Let us verify that we can as well use a weak-$\sL$ version of this
characterisation for jointly parameterised vector fields.
\begin{proposition}
Let\/ $\man{M}$ be a smooth manifold, let\/ $\ts{P}$ be a topological space,
and let\/ $\map{X}{\man{M}\times\ts{P}}{\tb{\man{M}}}$ satisfy\/
$X(x,p)\in\tb[x]{\man{M}}$ for every\/ $(x,p)\in\man{M}\times\ts{P}$\@.
Then\/ $X\in\JPsections[\infty]{\ts{P};\tb{\man{M}}}$ if and only if\/
$(x,p)\mapsto\lieder{X^p}{f}$ is a jointly parameterised function of class\/
$\C^\infty$ for every\/ $f\in\func[\infty]{\man{M}}$\@.
\begin{proof}
This follows from Corollary~\pldblref{cor:COinftyweak}{pl:COinftyweakcont}\@.
\end{proof}
\end{proposition}

\subsubsection{The finitely differentiable or Lipschitz case}

The preceding development in the smooth case is easily extended to the
finitely differentiable and Lipschitz cases, and we quickly give the results
and definitions here.  In this section, when considering the Lipschitz case,
we assume that $\nabla$ is the Levi-Civita connection associated to $\metric$
and we assume that $\nabla^0$ is $\metric_0$-orthogonal.
\begin{definition}\label{def:Cmm'-paramsec}
Let $\map{\pi}{\man{E}}{\man{M}}$ be a smooth vector bundle and let $\ts{P}$
be a topological space.  A map $\map{\xi}{\man{M}\times\ts{P}}{\man{E}}$ such
that $\xi(x,p)\in\man{E}_x$ for every $(x,p)\in\man{M}\times\ts{P}$
\begin{compactenum}[(i)]
\item is a \defn{separately parameterised section of class $\C^{m+m'}$} if
\begin{compactenum}[(a)]
\item for each $x\in\man{M}$\@, the map $\map{\xi_x}{\ts{P}}{\man{E}}$ defined
by $\xi_x(p)=\xi(x,p)$ is continuous and
\item for each $p\in\ts{P}$\@, the map $\map{\xi^p}{\man{M}}{\man{E}}$
defined by $\xi^p(x)=\xi(x,p)$ is of class $\C^{m+m'}$\@,
\end{compactenum}
and
\item is a \defn{jointly parameterised section of class $\C^{m+m'}$} if it is
a separately parameterised section of class $\C^{m+m'}$ and
\begin{compactenum}[(a)]
\item $m'=0$\@: the map $(x,p)\mapsto j_m\xi^p(x)$ is continuous;
\item $m'=\lip$\@: the map $(x,p)\mapsto j_m\xi^p(x)$ is continuous and, for
each $(x_0,p_0)\in\man{M}\times\ts{P}$ and each $\epsilon\in\realp$\@, there
exist a neighbourhood $\nbhd{U}\subset\man{M}$ of $x_0$ and a neighbourhood
$\nbhd{O}\subset\ts{P}$ of $p_0$ such that
\begin{equation*}
j_m\xi(\nbhd{U}\times\nbhd{O})\subset\setdef{j_m\eta(x)\in\jet{m}{\man{E}}}
{\dil{(j_m\eta-j_m\xi^{p_0})}(x)<\epsilon},
\end{equation*}
where, of course, $j_m\xi(x,p)=j_m\xi^p(x)$\@.
\end{compactenum}
\end{compactenum}
By $\SPsections[m+m']{\ts{P};\man{E}}$ we denote the set of separately
parameterised sections of $\man{E}$ of class $\C^{m+m'}$ and by
$\JPsections[m+m']{\ts{P};\man{E}}$ we denote the set of jointly
parameterised sections of $\man{E}$ of class $\C^{m+m'}$\@.\oprocend
\end{definition}

Let us give the purely topological characterisation of this class of
sections.
\begin{proposition}\label{prop:paramCmm'}
Let\/ $\map{\pi}{\man{E}}{\man{M}}$ be a smooth vector bundle, let\/ $\ts{P}$
be a topological space, and let\/ $\map{\xi}{\man{M}\times\ts{P}}{\man{E}}$
satisfy\/ $\xi(x,p)\in\man{E}_x$ for every\/ $(x,p)\in\man{M}\times\ts{P}$\@.
Then\/ $\xi\in\JPsections[m+m']{\ts{P};\man{E}}$ if and only if the map\/
$p\mapsto\xi^p\in\sections[m+m']{\man{E}}$ is continuous, where\/
$\sections[m+m']{\man{E}}$ has the\/ $\CO^{m+m'}$-topology.
\begin{proof}
We will prove the result only in the case that $m=0$ and $m'=\lip$\@, as the
general case follows by combining this case with the computations from the
proof of Proposition~\ref{prop:paramCinfty}\@.  We denote
$\sigma_\xi(p)=\xi(x,p)$\@.

Suppose that $(x,p)\mapsto\xi(x,p)$ is continuous and that, for every
$(x_0,p_0)\in\man{M}\times\ts{P}$ and for every $\epsilon\in\realp$\@, there
exist a neighbourhood $\nbhd{U}\subset\man{M}$ of $x_0$ and a neighbourhood
$\nbhd{O}\subset\ts{P}$ of $p_0$ such that, if
$(x,p)\in\nbhd{U}\times\nbhd{O}$\@, then
$\dil{(\xi^p-\xi^{p_0})}(x)<\epsilon$\@.  Let $K\subset\man{M}$ be compact,
let $\epsilon\in\realp$\@, and let $p_0\in\ts{P}$\@.  Let $x\in K$\@.  By
hypothesis, there exist a neighbourhood $\nbhd{U}_x\subset\man{M}$ of $x$
and a neighbourhood $\nbhd{O}_x\subset\ts{P}$ of $p_0$ such that
\begin{equation*}
\xi(\nbhd{U}_x\times\nbhd{O}_x)\subset
\setdef{\eta(x')\in\jet{m}{\man{E}}}{\dil{(\eta-\xi^{p_0})}(x')<\epsilon}.
\end{equation*}
Now let $x_1,\dots,x_k\in K$ be such that
$K\subset\cup_{j=1}^k\nbhd{U}_{x_j}$ and let
$\nbhd{O}=\cap_{j=1}^k\nbhd{O}_{x_j}$\@.  Then, if $p\in\nbhd{O}$ and $x\in
K$\@, we have $x\in\nbhd{U}_{x_j}$ for some $j\in\{1,\dots,k\}$\@.  Thus
\begin{equation*}
\dil{(\xi(x,p)-\xi(x,p_0))}_{\ol{\metric}_m}<\epsilon.
\end{equation*}
Therefore, taking supremums over $x\in K$\@, we have
$\lambda_K(\sigma_\xi(p)-\sigma_\xi(p_0))\le\epsilon$\@.  By choosing
$\nbhd{O}$ to be possibly smaller, the argument of
Proposition~\ref{prop:paramCinfty} ensures that
$p_K^0(\sigma_\xi(p)-\sigma_\xi(p_0))\le\epsilon$\@, and so
$p_K^{\lip}(\sigma_\xi(p)-\sigma_\xi(p_0))<\epsilon$ for $p\in\nbhd{O}$\@.
As this can be done for every compact $K\subset\man{M}$\@, we conclude that
$\sigma_\xi$ is continuous.

Next suppose that $\sigma_\xi$ is continuous, let
$(x_0,p_0)\in\man{M}\times\ts{P}$\@, and let $\epsilon\in\realp$\@.  Let
$\nbhd{U}$ be a relatively compact neighbourhood of $x_0$\@.  Since
$\sigma_\xi$ is continuous, let $\nbhd{O}$ be a neighbourhood of $p_0$ such
that
\begin{equation*}
p^{\lip}_{\closure(\nbhd{U})}(\sigma_\xi(p)-\sigma_\xi(p_0))<\epsilon,\qquad
p\in\nbhd{O}.
\end{equation*}
Thus, for every $(x,p)\in\nbhd{U}\times\nbhd{O}$\@,
$\dil{(\xi^p-\xi^{p_0})}(x)<\epsilon$\@.  Following the argument of
Proposition~\ref{prop:paramCinfty} one also shows that $\xi$ is continuous at
$(x_0,p_0)$\@, which shows that $\xi\in\JPsections[\lip]{\ts{P};\man{E}}$\@.
\end{proof}
\end{proposition}

Of course, the preceding discussion applies, in particular, to give vector
fields of jointly parameterised class $\C^{m+m'}$ and functions of jointly
parameterised class $\C^{m+m'}$\@.  This gives the spaces
$\SPC[m+m']{\ts{P};\man{M}}$ and $\JPC[m+m']{\man{M}}$ of parameterised
functions, and the spaces $\SPsections[m+m']{\ts{P};\tb{\man{M}}}$ and
$\JPsections[m+m']{\ts{P};\tb{\man{M}}}$ of parameterised vector fields.  Let
us verify that we can as well use a weak-$\sL$ version of this
characterisation for jointly parameterised vector fields.
\begin{proposition}
Let\/ $\man{M}$ be a smooth manifold, let\/ $\ts{P}$ be a topological space,
and let\/ $\map{X}{\man{M}\times\ts{P}}{\tb{\man{M}}}$ satisfy\/
$X(x,p)\in\tb[x]{\man{M}}$ for every\/ $(x,p)\in\man{M}\times\ts{P}$\@.
Then\/ $X\in\JPsections[m+m']{\ts{P};\tb{\man{M}}}$ if and only if\/
$(x,p)\mapsto\lieder{X^p}{f}$ is a jointly parameterised function of class\/
$\C^{m+m'}$ for every\/ $f\in\func[\infty]{\man{M}}$\@.
\begin{proof}
This follows from Corollary~\pldblref{cor:COmm'weak}{pl:COmm'weakcont}\@.
\end{proof}
\end{proposition}

\subsubsection{The holomorphic case}

As with time-varying vector fields, we are not really interested, \emph{per
se}\@, in holomorphic control systems, and in fact we will not even define
the notion.  However, it is possible, and possibly sometimes easier, to
verify that a control system satisfies our rather technical criterion of
being a ``real analytic control system'' by verifying that it possesses an
holomorphic extension.  Thus, in this section, we present the required
holomorphic definitions.  We will consider an holomorphic vector bundle
$\map{\pi}{\man{E}}{\man{M}}$ with an Hermitian fibre metric $\metric$\@.
This defines the seminorms $p^{\hol}_K$\@, $K\subset\man{M}$ compact,
describing the $\CO^{\hol}$-topology for $\sections[\hol]{\man{E}}$ as in
Section~\ref{subsec:COhol-vb}\@.
\begin{definition}
Let $\map{\pi}{\man{E}}{\man{M}}$ be an holomorphic vector bundle and let
$\ts{P}$ be a topological space.  A map
$\map{\xi}{\man{M}\times\ts{P}}{\man{E}}$ such that $\xi(z,p)\in\man{E}_z$
for every $(z,p)\in\man{M}\times\ts{P}$
\begin{compactenum}[(i)]
\item is a \defn{separately parameterised section of class $\C^{\hol}$} if
\begin{compactenum}[(a)]
\item for each $z\in\man{M}$\@, the map $\map{\xi_z}{\ts{P}}{\man{E}}$
defined by $\xi_z(p)=\xi(z,p)$ is continuous and
\item for each $p\in\ts{P}$\@, the map $\map{\xi^p}{\man{M}}{\man{E}}$
defined by $\xi^p(z)=\xi(z,p)$ is of class $\C^\hol$\@,
\end{compactenum}
and
\item is a \defn{jointly parameterised section of class $\C^{\hol}$} if it is
a separately parameterised section of class $\C^{\hol}$ and if the map
$(z,p)\mapsto\xi^p(z)$ is continuous.
\end{compactenum}
By $\SPsections[\hol]{\ts{P};\man{E}}$ we denote the set of separately
parameterised sections of $\man{E}$ of class $\C^{\hol}$ and by
$\JPsections[\hol]{\ts{P};\man{E}}$ we denote the set of jointly
parameterised sections of $\man{E}$ of class $\C^{\hol}$\@.\oprocend
\end{definition}

As in the smooth case, it is possible to give purely topological
characterisations of these classes of sections.
\begin{proposition}\label{prop:paramChol}
Let\/ $\map{\pi}{\man{E}}{\man{M}}$ be an holomorphic vector bundle, let\/
$\ts{P}$ be a topological space, and let\/
$\map{\xi}{\man{M}\times\ts{P}}{\man{E}}$ satisfy\/ $\xi(z,p)\in\man{E}_z$
for every\/ $(z,p)\in\man{M}\times\ts{P}$\@.  Then\/
$\xi\in\JPsections[\hol]{\ts{P};\man{E}}$ if and only if the map\/
$p\mapsto\xi^p\in\sections[\hol]{\man{E}}$ is continuous, where\/
$\sections[\hol]{\man{E}}$ has the\/ $\CO^{\hol}$-topology.
\begin{proof}
We define $\map{\sigma_\xi}{\ts{P}}{\sections[\hol]{\man{E}}}$ by
$\sigma_\xi(p)=\xi^p$\@.

First suppose that $\xi$ is continuous.  Let $K\subset\man{M}$ be compact,
let $\epsilon\in\realp$\@, and let $p_0\in\ts{P}$\@.  Let $z\in K$ and let
$\nbhd{W}_z\subset\man{E}$ be a neighbourhood of $\xi(z,p_0)$ for which
\begin{equation*}
\nbhd{W}_z\subset\setdef{\eta(z')\in\man{E}}
{\dnorm{\eta(z')-\xi(z',p_0)}_{\ol{\metric}}<\epsilon}.
\end{equation*}
By continuity of $\xi$\@, there exist a neighbourhood
$\nbhd{U}_z\subset\man{M}$ of $z$ and a neighbourhood
$\nbhd{O}_z\subset\ts{P}$ of $p_0$ such that
$\xi(\nbhd{U}_z\times\nbhd{O}_z)\subset\nbhd{W}_z$\@.  Now let
$z_1,\dots,z_k\in K$ be such that $K\subset\cup_{j=1}^k\nbhd{U}_{z_j}$ and
let $\nbhd{O}=\cap_{j=1}^k\nbhd{O}_{z_j}$\@.  Then, if $p\in\nbhd{O}$ and
$z\in K$\@, we have $z\in\nbhd{U}_{z_j}$ for some $j\in\{1,\dots,k\}$\@.
Thus $\xi(z,p)\in\nbhd{W}_{z_j}$\@.  Thus
$\dnorm{\xi(z,p)-\xi(z,p_0)}_{\ol{\metric}}<\epsilon$\@.  Therefore, taking
supremums over $z\in K$\@,
$p^{\hol}_K(\sigma_\xi(p)-\sigma_\xi(p_0))\le\epsilon$\@.  As this can be
done for every compact $K\subset\man{M}$\@, we conclude that $\sigma_\xi$ is
continuous.

Next suppose that $\sigma_\xi$ is continuous.  Let
$(z_0,p_0)\in\man{M}\times\ts{P}$ and let $\nbhd{W}\subset\man{E}$ be a
neighbourhood of $\xi(z_0,p_0)$\@.  Let $\nbhd{U}\subset\man{M}$ be a
relatively compact neighbourhood of $z_0$ and let $\epsilon\in\realp$ be such
that
\begin{equation*}
\pi^{-1}(\nbhd{U})\cap\setdef{\eta(z)\in\man{E}}
{\dnorm{\eta(z)-\xi(z,p_0)}_{\ol{\metric}}<\epsilon}\subset\nbhd{W}.
\end{equation*}
By continuity of $\sigma_\xi$\@, let $\nbhd{O}\subset\ts{P}$ be a
neighbourhood of $p_0$ such that
$p^{\hol}_{\closure(\nbhd{U})}(\sigma_\xi(p)-\sigma_\xi(p_0))<\epsilon$ for
$p\in\nbhd{O}$\@.  Therefore,
\begin{equation*}
\dnorm{\sigma_\xi(p)(z)-\sigma_\xi(p_0)(z)}_{\ol{\metric}}<\epsilon,
\qquad(z,p)\in\closure(\nbhd{U})\times\nbhd{O}.
\end{equation*}
Therefore, if $(z,p)\in\nbhd{U}\times\nbhd{O}$\@, we have
$\xi(z,p)\in\nbhd{W}$\@, showing that $\xi$ is continuous at $(z_0,p_0)$\@.
\end{proof}
\end{proposition}

The specialisation of the preceding constructions to vector fields and
functions is immediate.  This gives the spaces $\SPC[\hol]{\ts{P};\man{M}}$
and $\JPC[\hol]{\man{M}}$ of parameterised functions, and the spaces
$\SPsections[\hol]{\ts{P};\tb{\man{M}}}$ and
$\JPsections[\hol]{\ts{P};\tb{\man{M}}}$ of parameterised vector fields.  Let
us verify that we can as well use a weak-$\sL$ version of the preceding
definitions for vector fields in the case when the base manifold is Stein.
\begin{proposition}
Let\/ $\man{M}$ be a Stein manifold, let\/ $\ts{P}$ be a topological space,
and let\/ $\map{X}{\man{M}\times\ts{P}}{\tb{\man{M}}}$ satisfy\/
$X(z,p)\in\tb[z]{\man{M}}$ for every\/ $(z,p)\in\man{M}\times\ts{P}$\@.
Then\/ $X\in\JPsections[\hol]{\ts{P};\tb{\man{M}}}$ if and only if\/
$(x,p)\mapsto\lieder{X^p}{f}$ is a jointly parameterised function of class\/
$\C^{\hol}$ for every\/ $f\in\func[\infty]{\man{M}}$\@.
\begin{proof}
This follows from Corollary~\pldblref{cor:COweak}{pl:COweakcont}\@.
\end{proof}
\end{proposition}

\subsubsection{The real analytic case}

Now we repeat the procedure above for real analytic sections.  We thus will
consider a real analytic vector bundle $\map{\pi}{\man{E}}{\man{M}}$ with
$\nabla^0$ a real analytic linear connection on $\man{E}$\@, $\nabla$ a real
analytic affine connection on $\man{M}$\@, $\metric_0$ a real analytic fibre
metric on $\man{E}$\@, and $\metric$ a real analytic Riemannian metric on
$\man{M}$\@.  This defines the seminorms $p^\omega_{K,\vect{a}}$\@,
$K\subset\man{M}$ compact, $\vect{a}\in\csd(\integernn;\realp)$\@, describing
the $\C^\omega$-topology as in Theorem~\ref{the:Comega-seminorms}\@.
\begin{definition}\label{def:Comega-paramsec}
Let $\map{\pi}{\man{E}}{\man{M}}$ be a real analytic vector bundle and let
$\ts{P}$ be a topological space.  A map
$\map{\xi}{\man{M}\times\ts{P}}{\man{E}}$ such that $\xi(x,p)\in\man{E}_x$
for every $(x,p)\in\man{M}\times\ts{P}$
\begin{compactenum}[(i)]
\item is a \defn{separately parameterised section of class $\C^\omega$} if
\begin{compactenum}[(a)]
\item for each $x\in\man{M}$\@, the map $\map{\xi_x}{\ts{P}}{\man{E}}$ defined
by $\xi_x(p)=\xi(x,p)$ is continuous and
\item for each $p\in\ts{P}$\@, the map $\map{\xi^p}{\man{M}}{\man{E}}$
defined by $\xi^p(x)=\xi(x,p)$ is of class $\C^\omega$\@,
\end{compactenum}
and
\item is a \defn{jointly parameterised section of class $\C^\omega$} if it is
a separately parameterised section of class $\C^\infty$ and if, for each
$(x_0,p_0)\in\man{M}\times\ts{P}$\@, for each
$\vect{a}\in\csd(\integernn,\realp)$\@, and for each $\epsilon\in\realp$\@,
there exist a neighbourhood $\nbhd{U}\subset\man{M}$ of $x_0$ and a
neighbourhood $\nbhd{O}\subset\ts{P}$ of $p_0$ such that
\begin{equation*}
j_m\xi(\nbhd{U}\times\nbhd{O})\subset
\setdef{j_m\eta(x)\in\jet{m}{\man{E}}}
{a_0a_1\cdots a_m\dnorm{j_m\eta(x)-j_m\xi^{p_0}(x)}_{\ol{\metric}_m}
<\epsilon}
\end{equation*}
for every\/ $m\in\integernn$\@, where, of course,
$j_m\xi(x,p)=j_m\xi^p(x)$\@.
\end{compactenum}
By $\SPsections[\omega]{\ts{P};\man{E}}$ we denote the set of separately
parameterised sections of $\man{E}$ of class $\C^\omega$ and by
$\JPsections[\omega]{\ts{P};\man{E}}$ we denote the set of jointly
parameterised sections of $\man{E}$ of class $\C^\omega$\@.\oprocend
\end{definition}

\begin{remark}\label{rem:Comega-param}
The condition that $\xi\in\JPsections[\infty]{\ts{P};\man{E}}$ can be
restated like this: for each $(x_0,p_0)\in\man{M}\times\ts{P}$\@, for each
$m\in\integernn$\@, and for each $\epsilon\in\realp$\@, there exist a
neighbourhood $\nbhd{U}\subset\man{M}$ of $x_0$ and a neighbourhood
$\nbhd{O}\subset\ts{P}$ of $p_0$ such that
\begin{equation*}
j_m\xi(\nbhd{U}\times\nbhd{O})\subset
\setdef{j_m\eta(x)\in\jet{m}{\man{E}}}
{\dnorm{j_m\eta(x)-j_m\xi^{p_0}(x)}_{\ol{\metric}_m}<\epsilon};
\end{equation*}
that this is so is, more or less, the idea of the proof of
Proposition~\ref{prop:paramCinfty}\@.  Phrased this way, one sees clearly the
grammatical similarity between the smooth and real analytic definitions.
Indeed, the grammatical transformation from the smooth to the real analytic
definition is, \emph{put a factor of $a_0a_1\cdots a_m$ before the norm,
precede the condition with ``for every
$\vect{a}\in\csd(\integernn;\realp)$'', and move the ``for every
$m\in\integernn$'' from before the condition to after.}  This was also seen
in the definitions of locally integrally bounded and locally essentially
bounded sections in Section~\ref{sec:time-varying}\@.  Indeed, the
grammatical similarity will be encountered many times in the sequel, and we
shall refer to this to keep ourselves from repeating arguments in the real
analytic case that mirror their smooth counterparts.\oprocend
\end{remark}

The following result records topological characterisations of jointly
parameterised sections in the real analytic case.
\begin{proposition}\label{prop:paramComega}
Let\/ $\map{\pi}{\man{E}}{\man{M}}$ be a real analytic vector bundle, let\/
$\ts{P}$ be a topological space, and let\/
$\map{\xi}{\man{M}\times\ts{P}}{\man{E}}$ satisfy\/ $\xi(x,p)\in\man{E}_x$
for every\/ $(x,p)\in\man{M}\times\ts{P}$\@.  Then\/
$\xi\in\JPsections[\omega]{\ts{P};\man{E}}$ if and only if the map\/
$p\mapsto\xi^p\in\sections[\omega]{\man{E}}$ is continuous, where\/
$\sections[\omega]{\man{E}}$ has the\/ $\C^\omega$-topology.
\begin{proof}
For $\vect{a}\in\csd(\integernn;\realp)$ and $m\in\integernn$\@, given
$\map{\xi}{\man{M}\times\ts{P}}{\man{E}}$ satisfying
$\xi^p\in\sections[\omega]{\man{E}}$\@, we let
$\map{\xi_{\vect{a},m}}{\man{M}\times\ts{P}}{\jet{m}{\man{E}}}$ be the map
\begin{equation*}
\xi_{\vect{a},m}(x,p)=a_0a_1\cdots a_mj_m\xi^p(x).
\end{equation*}
We also denote by $\map{\sigma_\xi}{\ts{P}}{\sections[\omega]{\man{E}}}$ the
map given by $\sigma_\xi(p)=\xi^p$\@.

Suppose that, for every $(x_0,p_0)\in\man{M}\times\ts{P}$\@, for every
$\vect{a}\in\csd(\integernn;\realp)$\@, and for every $\epsilon\in\realp$\@,
there exist a neighbourhood $\nbhd{U}\subset\man{M}$ of $x_0$ and a
neighbourhood $\nbhd{O}\subset\ts{P}$ of $p_0$ such that, if
$(x,p)\in\nbhd{U}\times\nbhd{O}$\@, then
\begin{equation*}
\dnorm{\xi_{\vect{a},m}(x,p)-\xi_{\vect{a},m}(x,p_0)}_{\ol{\metric}_m}<\epsilon,
\qquad m\in\integernn.
\end{equation*}
Let $K\subset\man{M}$ be compact, let $\vect{a}\in\csd(\integernn;\realp)$\@,
let $\epsilon\in\realp$\@, and let $p_0\in\ts{P}$\@.  Let $x\in K$\@.  By
hypothesis, there exist a neighbourhood $\nbhd{U}_x\subset\man{M}$ of $x$
and a neighbourhood $\nbhd{O}_x\subset\ts{P}$ of $p_0$ such that
\begin{equation*}
\xi_{\vect{a},m}(\nbhd{U}_x\times\nbhd{O}_x)\subset
\setdef{j_m\eta(x')\in\jet{m}{\man{E}}}
{\dnorm{a_0a_1\cdots a_mj_m\eta(x')-
\xi_{\vect{a},m}(x',p_0)}_{\ol{\metric}_m}<\epsilon},
\end{equation*}
for each $m\in\integernn$\@.  Now let $x_1,\dots,x_k\in K$ be such that
$K\subset\cup_{j=1}^k\nbhd{U}_{x_j}$ and let
$\nbhd{O}=\cap_{j=1}^k\nbhd{O}_{x_j}$\@.  Then, if $p\in\nbhd{O}$ and $x\in
K$\@, we have $x\in\nbhd{U}_{x_j}$ for some $j\in\{1,\dots,k\}$\@.  Thus
\begin{equation*}
\dnorm{\xi_{\vect{a},m}(x,p)-\xi_{\vect{a},m}(x,p_0)}_{\ol{\metric}_m}<\epsilon,
\qquad m\in\integernn.
\end{equation*}
Therefore, taking supremums over $x\in K$ and $m\in\integernn$\@, we have
$p^\omega_{K,\vect{a}}(\sigma_\xi(p)-\sigma_\xi(p_0))\le\epsilon$\@.  As this
can be done for every compact $K\subset\man{M}$ and every
$\vect{a}\in\csd(\integernn;\realp)$\@, we conclude that $\sigma_\xi$ is
continuous.

Next suppose that $\sigma_\xi$ is continuous, let
$(x_0,p_0)\in\man{M}\times\ts{P}$\@, let
$\vect{a}\in\csd(\integernn;\realp)$\@, and let $\epsilon\in\realp$\@.  Let
$\nbhd{U}$ be a relatively compact neighbourhood of $x_0$\@.  Since
$\sigma_\xi$ is continuous, let $\nbhd{O}$ be a neighbourhood of $p_0$ such
that
\begin{equation*}
p^\omega_{\closure(\nbhd{U}),\vect{a}}
(\sigma_\xi(p)-\sigma_\xi(p_0))<\epsilon,\qquad
p\in\nbhd{O}.
\end{equation*}
Thus, for every $(x,p)\in\nbhd{U}\times\nbhd{O}$\@,
\begin{equation*}
a_0a_1\cdots a_m\dnorm{j_m\xi(x,p)-j_m\xi(x,p_0)}_{\ol{\metric}_m}<\epsilon,
\qquad m\in\integernn,
\end{equation*}
which shows that $\xi\in\JPsections[\omega]{\ts{P};\man{E}}$\@.
\end{proof}
\end{proposition}

As we have done in the smooth and holomorphic cases above, we can specialise
the preceding discussion from sections to vector fields and functions, giving
the spaces $\SPC[\omega]{\ts{P};\man{M}}$ and $\JPC[\omega]{\man{M}}$ of
parameterised functions, and the spaces
$\SPsections[\omega]{\ts{P};\tb{\man{M}}}$ and
$\JPsections[\omega]{\ts{P};\tb{\man{M}}}$ of parameterised vector fields.
We then have the following weak-$\sL$ characterisation for jointly
parameterised vector fields.
\begin{proposition}
Let\/ $\man{M}$ be a real analytic manifold, let\/ $\ts{P}$ be a topological
space, and let\/ $\map{X}{\man{M}\times\ts{P}}{\tb{\man{M}}}$ satisfy\/
$X(x,p)\in\tb[x]{\man{M}}$ for every\/ $(x,p)\in\man{M}\times\ts{P}$\@.
Then\/ $X\in\JPsections[\omega]{\ts{P};\tb{\man{M}}}$ if and only if\/
$(x,p)\mapsto\lieder{X^p}{f}(x)$ is a jointly parameterised function of
class\/ $\C^\omega$ for every\/ $f\in\func[\omega]{\man{M}}$\@.
\begin{proof}
This follows from Corollary~\pldblref{cor:Comegaweak}{pl:Comegaweakcont}\@.
\end{proof}
\end{proposition}

One can wonder about the relationship between sections of jointly
parameterised class $\C^\omega$ and sections that are real restrictions of
sections of jointly parameterised class $\C^{\hol}$\@.  We address this with
a result and an example.  First the result.
\begin{theorem}\label{the:Comega->Cholparam}
Let\/ $\map{\pi}{\man{E}}{\man{M}}$ be a real analytic vector bundle with
holomorphic extension\/ $\map{\ol{\pi}}{\ol{\man{E}}}{\ol{\man{M}}}$ and
let\/ $\ts{P}$ be a topological space.  For a map\/
$\map{\xi}{\man{M}\times\ts{P}}{\man{E}}$ satisfying\/ $\xi(x,p)\in\man{E}_x$
for all\/ $(x,p)\in\man{M}\times\ts{P}$\@, the following statements hold:
\begin{compactenum}[(i)]
\item \label{pl:ComegaCholparam1} if\/
$\xi\in\JPsections[\omega]{\ts{P};\man{E}}$ and if\/ $\ts{P}$ is locally
compact and Hausdorff, then, for each\/ $(x_0,p_0)\in\man{M}\times\ts{P}$\@,
there exist a neighbourhood\/ $\ol{\nbhd{U}}\subset\ol{\man{M}}$ of\/
$x_0$\@, a neighbourhood\/ $\nbhd{O}\subset\ts{P}$ of\/ $p_0$\@, and\/
$\ol{\xi}\in\JPsections[\hol]{\nbhd{O};\ol{\man{E}}|\ol{\nbhd{U}}}$ such
that\/ $\xi(x,p)=\ol{\xi}(x,p)$ for all\/
$(x,p)\in(\man{M}\cap\ol{\nbhd{U}})\times\nbhd{O}$\@;
\item \label{pl:ComegaCholparam2} if there exists a section\/
$\ol{\xi}\in\JPsections[\hol]{\ts{P};\ol{\man{E}}}$ such that\/
$\xi(x,p)=\ol{\xi}(x,p)$ for every\/ $(x,p)\in\man{M}\times\ts{P}$\@, then\/
$\xi\in\JPsections[\omega]{\ts{P};\man{E}}$\@.
\end{compactenum}
\begin{proof}
\eqref{pl:ComegaCholparam1} Let $p_0\in\ts{P}$ and let $\nbhd{O}$ be a
relatively compact neighbourhood of $p_0$\@, this being possible since
$\ts{P}$ is locally compact.  Let $x_0\in\man{M}$\@, let $\nbhd{U}$ be a
relatively compact neighbourhood of $x_0$\@, and let
$\ifam{\ol{\nbhd{U}}_j}_{j\in\integerp}$ be a sequence of neighbourhoods of
$\closure(\nbhd{U})$ in $\ol{\man{M}}$ with the properties that
$\closure(\ol{\nbhd{U}}_j)\subset\ol{\nbhd{U}}_{j+1}$ and that
$\cap_{j\in\integerp}\ol{\nbhd{U}}_j=\closure(\nbhd{U})$\@.  We first note
that
\begin{equation*}
\mappings[0]{\closure(\nbhd{O})}
{\gsections[\hol,\real]{\closure(\nbhd{U})}{\ol{\man{E}}}}\simeq
\func[0]{\closure(\nbhd{O})}\injotimes
\gsections[\hol,\real]{\closure(\nbhd{U})}{\ol{\man{E}}}
\end{equation*}
and
\begin{equation*}
\mappings[0]{\closure(\nbhd{O})}
{\sections[\hol,\real]{\ol{\man{E}}|\ol{\nbhd{U}}_j}}\simeq
\func[0]{\closure(\nbhd{O})}\injotimes
\sections[\hol,\real]{\ol{\man{E}}|\ol{\nbhd{U}}_j},
\end{equation*}
with $\injotimes$ denoting the completed injective tensor product;
see~\cite[Chapter~16]{HJ:81} for the injective tensor product for locally
convex spaces and~\cite[Theorem~1.1.10]{JD/JHF/JS:08} for the preceding
isomorphisms for Banach spaces (the constructions apply more or less verbatim
to locally convex spaces~\cite[Proposition~5.4]{KDB:07}).  One can also
prove, using the argument from the proof
of~\cite[Theorem~1.1.10]{JD/JHF/JS:08} (see top of page~15 of that
reference), that, if $[\ol{\xi}]_K\in\mappings[0]{\closure(\nbhd{O})}
{\gsections[\hol,\real]{\closure(\nbhd{U})}{\ol{\man{E}}}}$\@, then there is
a \emph{sequence} (we know there is a \emph{net})
$\ifam{[\ol{\xi}_k]_{\closure(\nbhd{U})}}_{k\in\integerp}$ in
$\func[0]{\closure(\nbhd{O})}\otimes
\gsections[\hol,\real]{\closure(\nbhd{U})}{\ol{\man{E}}}$ converging to
$[\ol{\xi}]_K$ in the completed injective tensor product topology.  Note that
since $\gsections[\hol,\real]{\closure(\nbhd{U})}{\ol{\man{E}}}$ and
$\sections[\hol,\real]{\ol{\man{E}}|\ol{\nbhd{U}}_j}$\@, $j\in\integerp$\@,
are nuclear, the injective tensor product can be swapped with the projective
tensor product in the above constructions~\cite[Proposition~5.4.2]{AP:69}\@.
One can now reproduce the argument from the proof of
Theorem~\ref{the:Comega->Chol}\@, swapping $\L^1(\tdomain';\real)$ with
$\func[0]{\closure(\nbhd{O})}$ and using the results of \citet{EMM:97}\@, to
complete the proof in this case.

\eqref{pl:ComegaCholparam2} Let $(x_0,p_0)\in\man{M}\times\ts{P}$\@, let
$\vect{a}\in\csd(\integernn;\realp)$\@, and let $\epsilon\in\realp$\@.  Let
$\nbhd{U}\subset\man{M}$ be a relatively compact neighbourhood of $x_0$ and
let $\ol{\nbhd{U}}$ be a relatively compact neighbourhood of
$\closure(\nbhd{U})$\@.  By Proposition~\ref{prop:Cauchyest}\@, there exist
$C,r\in\realp$ such that
\begin{equation*}
p^\infty_{\closure(\nbhd{U}),m}(\sigma_\xi(p)-\sigma_\xi(p_0))\le
Cr^{-m}\sup\setdef{\dnorm{\ol{\xi}(z,p)-\ol{\xi}(z,p_0)}_{\ol{\metric}}}
{z\in\ol{\nbhd{U}}}
\end{equation*}
for all $m\in\integernn$ and $p\in\ts{P}$\@.  Now let $N\in\integernn$ be
such that $a_{N+1}<r$ and let $\nbhd{O}$ be a neighbourhood of $p_0$ such
that
\begin{equation*}
\dnorm{\ol{\xi}(z,p)-\ol{\xi}(z,p_0)}_{\ol{\metric}}<
\frac{\epsilon r^m}{Ca_0a_1\cdots a_m},\qquad m\in\{0,1,\dots,N\},
\end{equation*}
for $(z,p)\in\ol{\nbhd{U}}\times\nbhd{O}$\@.  Then, if
$m\in\{0,1,\dots,N\}$\@, we have
\begin{multline*}
a_0a_1\cdots a_m
\dnorm{j_m\xi^p(x)-j_m\xi^{p_0}(x)}_{\ol{\metric}_m}\\
\le a_0a_1\cdots a_mCr^{-m}
\sup\setdef{\dnorm{\ol{\xi}(z,p)-\ol{\xi}(z,p_0)}_{\ol{\metric}_m}}
{z\in\ol{\nbhd{U}}}
<\epsilon,
\end{multline*}
for $(x,p)\in\nbhd{U}\times\nbhd{O}$\@.  If $m>N$ we also have
\begin{align*}
a_0a_1\cdots a_m
\|j_m\xi^p(x)-&j_m\xi^{p_0}(x)\|_{\ol{\metric}_m}\\
\le&\;a_0a_1\cdots a_Nr^{-N}r^m
\dnorm{j_m\xi^p(x)-j_m\xi^{p_0}(x)}_{\ol{\metric}_m}\\
\le&\;a_0a_1\cdots a_Nr^{-N}r^mCr^{-m}
\sup\setdef{\dnorm{\ol{\xi}(z,p)-\ol{\xi}(z,p_0)}_{\ol{\metric}_m}}
{z\in\ol{\nbhd{U}}}<\epsilon,
\end{align*}
for $(x,p)\in\nbhd{U}\times\nbhd{O}$\@, as desired.
\end{proof}
\end{theorem}

The next example shows that the assumption of local compactness cannot be
generally relaxed.
\begin{example}
Let $\man{M}=\real$\@, let $\ts{P}=\func[\omega]{\real}$\@, and define
$\map{f}{\real\times\ts{P}}{\real}$ by $f(x,g)=g(x)$\@.  Since $g\mapsto f^g$
is the identity map, we conclude from Proposition~\ref{prop:paramComega} that
$f\in\JPC[\omega]{\ts{P};\man{M}}$\@.  Let $x_0\in\real$\@.  We claim that,
for any neighbourhood $\ol{\nbhd{U}}$ of $x_0$ in $\complex$ and any
neighbourhood $\nbhd{O}$ of $0\in\ts{P}$\@, there exists $g\in\nbhd{O}$ such
that $g$\@, and therefore $f^g$\@, does not have an holomorphic extension to
$\ol{\nbhd{U}}$\@.  To see this, let $\sigma\in\realp$ be such that the disk
$\cdisk[]{\sigma}{x_0}$ in $\complex$ is contained in $\ol{\nbhd{U}}$\@.  Let
$K_1,\dots,K_r\subset\real$ be compact, let
$\vect{a}_1,\dots,\vect{a}_r\in\csd(\integernn;\realp)$\@, and let
$\epsilon_1,\dots,\epsilon_r\in\realp$ be such that
\begin{equation*}
\cap_{j=1}^r\setdef{g\in\ts{P}}{p_{K_j,\vect{a}_j}(g)\le\epsilon_j}\subset\nbhd{O}.
\end{equation*}
Now define
\begin{equation*}
g(x)=\frac{\alpha}{1+((x-x_0)/\sigma)^2},\qquad x\in\real,
\end{equation*}
with $\alpha\in\realp$ chosen sufficiently small that
$p_{K_j,\vect{a}_j}(g)<\epsilon_j$\@, $j\in\{1,\dots,r\}$\@, and note that
$g\in\nbhd{O}$ does not have an holomorphic extension to
$\ol{\nbhd{U}}$\@,~\cf~the discussion at the beginning of
Section~\ref{sec:analytic-topology}\@.\oprocend
\end{example}

\subsubsection{Mixing regularity hypotheses}

Just as we discussed with time-varying vector fields in
Section~\ref{subsec:time-mixup}\@, it is possible to consider parameterised
sections with mixed regularity hypotheses.  Indeed, the conditions of
Definitions~\ref{def:Cinfty-paramsec}\@,~\ref{def:Cmm'-paramsec}\@,
and~\ref{def:Comega-paramsec} are joint on state and parameter.  Thus we may
consider the following situation.  Let $m\in\integernn$\@,
$m'\in\{0,\lip\}$\@, $r\in\integernn\cup\{\infty,\omega\}$\@, and
$r'\in\{0,\lip\}$\@.  If $r+r'\ge m+m'$ (with the obvious convention that
$\infty+\lip=\infty$ and $\omega+\lip=\omega$)\@, we may then consider a
parameterised section in
\begin{equation*}
\SPsections[r+r']{\ts{P};\man{E}}\cap\JPsections[m+m']{\ts{P};\man{E}}
\end{equation*}
As with time-varying vector fields, there is nothing wrong with
this\textemdash{}indeed this is often done\textemdash{}as long as one
remembers what is true and what is not in the case when $r+r'>m+m'$\@.

\subsection{Control systems with locally essentially bounded
controls}\label{subsec:Cr-systemsI}

Let us first establish some terminology we will use throughout the remainder
of the paper.
\begin{notation}\label{notn:regularity}
Starting in this section, and continuing throughout the remainder of the
paper, we will simultaneously be considering finitely differentiable,
Lipschitz, smooth, and real analytic hypotheses.  To do this, we will let
$m\in\integernn$ and $m'\in\{0,\lip\}$\@, and consider the regularity classes
$\nu\in\{m+m',\infty,\omega\}$\@.  In such cases we shall require that the
underlying manifold be of class ``$\C^r$\@, $r\in\{\infty,\omega\}$\@, as
required.''  This has the obvious meaning, namely that we consider class
$\C^\omega$ if $\nu=\omega$ and class $\C^\infty$ otherwise.

Proofs will typically break into the four cases $\nu=\infty$\@, $\nu=m$\@,
$\nu=m+\lip$\@, and $\nu=\omega$\@.  In most cases there is a structural
similarity in the way arguments are carried out, so we will oftentimes do all
cases at once.  In doing this, we will, for $K\subset\man{M}$ be compact, for
$k\in\integernn$\@, and for $\vect{a}\in\csd(\integernn;\realp)$\@, denote
\begin{equation*}
p_K=\begin{cases}p^\infty_{K,k},&\nu=\infty,\\p^m_K,&\nu=m,\\
p^{m+\lip}_K,&\nu=m+\lip,\\p^\omega_{K,\vect{a}},&\nu=\omega.\end{cases}
\end{equation*}
Then, using the fact that $\xi\in\LIsections[\nu]{\tdomain;\man{E}}$ if and
only if there exists $g\in\Lloc^1(\tdomain;\realnn)$ such that $p_K(\xi_t)\le
g(t)$ (with a similar sort of assertion for parameterised section), we argue
all cases simultaneously.  The convenience and brevity more than make up for
the slight loss of preciseness in this approach.\oprocend
\end{notation}

With the notions of parameterised sections from the preceding section, we
readily define what we mean by a control system.
\begin{definition}\label{def:Cr-system}
Let $m\in\integernn$ and $m'\in\{0,\lip\}$\@, let\/
$\nu\in\{m+m',\infty,\omega\}$\@, and let $r\in\{\infty,\omega\}$\@, as
required.  A \defn{$\C^\nu$-control system} is a triple
$\Sigma=(\man{M},F,\cs{C})$\@, where
\begin{compactenum}[(i)]
\item $\man{M}$ is a $\C^r$-manifold whose elements are called
\defn{states}\@,
\item $\cs{C}$ is a topological space called the \defn{control set}\@, and
\item $F\in\JPsections[\nu]{\cs{C};\tb{\man{M}}}$\@.\oprocend
\end{compactenum}
\end{definition}

The governing equations for a control system are
\begin{equation*}
\xi'(t)=F(\xi(t),\mu(t)),
\end{equation*}
for suitable functions $t\mapsto\mu(t)\in\cs{C}$ and
$t\mapsto\xi(t)\in\man{M}$\@.  To ensure that these equations make sense, the
differential equation should be shown to have the properties needed for
existence and uniqueness of solutions, as well as appropriate dependence on
initial conditions. We do this by allowing the controls for the system to be
as general as reasonable.
\begin{proposition}\label{prop:open-loopcsI}
Let\/ $m\in\integernn$ and\/ $m'\in\{0,\lip\}$\@, let\/
$\nu\in\{m+m',\infty,\omega\}$\@, and let\/ $r\in\{\infty,\omega\}$\@, as
required.  Let\/ $\Sigma=(\man{M},F,\cs{C})$ be a\/ $\C^\nu$-control system.
If\/ $\mu\in\Lloc^\infty(\tdomain;\cs{C})$ (boundedness here being taking
with respect to the compact bornology) then
$F^\mu\in\LBsections[\nu]{\tdomain,\tb{\man{M}}}$\@, where\/
$\map{F^\mu}{\tdomain\times\man{M}}{\tb{\man{M}}}$ is defined by\/
$F^\mu(t,x)=F(x,\mu(t))$\@.
\begin{proof}
Let us define
$\map{\hat{F}\null^\mu}{\tdomain}{\sections[\nu]{\tb{\man{M}}}}$ by
$\hat{F}\null^\mu(t)=F^\mu_t$\@.  By
Propositions~\ref{prop:paramCinfty}\@,~\ref{prop:paramCmm'}\@,
and~\ref{prop:paramComega}\@, the mapping $u\mapsto F^u$ is continuous.
Since $\hat{F}\null^\mu$ is thus the composition of the measurable function
$\mu$ and the continuous mapping $u\mapsto F^u$\@, it follows that
$\hat{F}\null^\mu$ is measurable.  It follows from
Theorems~\ref{the:Cinfty-tvsec}\@,~\ref{the:Cmm'-tvsec}\@,
and~\ref{the:Comega-tvsec} that $F^\mu$ is a Carath\'eodory vector field of
class $\C^\nu$\@.

Let $\tdomain'\subset\tdomain$ be compact.  Since $\mu$ is locally
essentially bounded, there exists a compact set $K\subset\cs{C}$ such that
\begin{equation*}
\lebmes(\setdef{t\in\tdomain'}{\mu(t)\not\in K})=0.
\end{equation*}
Since the mapping $u\mapsto F^u$ is continuous,
\begin{equation*}
\setdef{F^\mu_t}{t\in\tdomain'}
\end{equation*}
is contained in a compact subset of
$\sections[\nu]{\tb{\man{M}}}$\@,~\ie~$F^\mu$ is locally essentially bounded.
\end{proof}
\end{proposition}

The notion of a trajectory is, of course, well known.  However, we make the
definitions clear for future reference.
\begin{definition}\label{def:cstrajI}
Let $m\in\integernn$ and $m'\in\{0,\lip\}$\@, let\/
$\nu\in\{m+m',\infty,\omega\}$\@, and let $r\in\{\infty,\omega\}$\@, as
required.  Let $\Sigma=(\man{M},F,\cs{C})$ be a $\C^\nu$-control system.  For
an interval $\tdomain\subset\real$\@, a \defn{$\tdomain$-trajectory} is a
locally absolutely continuous curve $\map{\xi}{\tdomain}{\man{M}}$ for which
there exists $\mu\in\Lloc^\infty(\tdomain;\cs{C})$ such that
\begin{equation*}
\xi'(t)=F(\xi(t),\mu(t)),\qquad\ae\ t\in\tdomain.
\end{equation*}
The set of $\tdomain$-trajectories we denote by $\Traj(\tdomain,\Sigma)$\@.
If $\nbhd{U}$ is open, we denote by $\Traj(\tdomain,\nbhd{U},\Sigma)$ those
trajectories taking values in\/ $\nbhd{U}$\@.\footnote{This is not a common
notion in this context, and our introduction of this is for the convenience
of making comparisons in the next section; see
Theorems~\ref{the:cs->gcs-equivI} and~\ref{the:cs->gcs-equivII}\@.}\oprocend
\end{definition}

One may also wish to restrict the class of controls one uses.  Thus we can
consider, for each time-domain $\tdomain$\@, a subset
$\sC(\tdomain)\subset\Lloc^\infty(\tdomain;\cs{C})$\@.  Generally, one will
ask for some compatibility conditions for these subsets, like, for example,
that, if $\tdomain'\subset\tdomain$\@, then $\mu|\tdomain'\in\sC(\tdomain')$
for every $\mu\in\sC(\tdomain)$\@.  For example, one may consider things like
piecewise continuous or piecewise constant controls.  In this case, we denote
by $\Traj(\tdomain,\sC)$ the set of trajectories arising from using controls
from $\sC(\tdomain)$\@.  Similarly, by $\Traj(\tdomain,\nbhd{U},\sC)$ we
denote the trajectories from this set taking values in an open set
$\nbhd{U}$\@.  We shall see in Section~\ref{sec:gcs} that our \gcs{}s provide
a natural means of capturing issues such as this.

\subsection{Control systems with locally integrable controls}\label{subsec:Cr-systemsII}

In this section we specialise the discussion from the preceding section in
one direction, while generalising it in another.  To be precise, we now
consider the case where our control set $\cs{C}$ is a subset of a locally
convex topological vector space, and the system structure is such that the
notion of integrability is preserved (in a way that will be made clear in
Proposition~\ref{prop:open-loopcsII} below).
\begin{definition}
Let $m\in\integernn$ and $m'\in\{0,\lip\}$\@, let\/
$\nu\in\{m+m',\infty,\omega\}$\@, and let $r\in\{\infty,\omega\}$\@, as
required.  A \defn{$\C^\nu$-sublinear control system} is a triple
$\Sigma=(\man{M},F,\cs{C})$\@, where
\begin{compactenum}[(i)]
\item $\man{M}$ is a $\C^r$-manifold whose elements are called
\defn{states}\@,
\item $\cs{C}$ is a subset of a locally convex topological vector space
$\alg{V}$\@, $\cs{C}$ being called the \defn{control set}\@, and
\item \label{pl:sublinear} $\map{F}{\man{M}\times\cs{C}}{\tb{\man{M}}}$ has
the following property: for every continuous seminorm $p$ for
$\sections[\nu]{\tb{\man{M}}}$\@, there exists a continuous seminorm $q$ for
$\alg{V}$ such that
\begin{equation*}\eqoprocend
p(F^{u_1}-F^{u_2})\le q(u_1-u_2),\qquad u_1,u_2\in\cs{C}.
\end{equation*}
\end{compactenum}
\end{definition}

Note that, by Propositions~\ref{prop:paramCinfty}\@,~\ref{prop:paramCmm'}\@,
and~\ref{prop:paramComega}\@, the sublinearity condition~\eqref{pl:sublinear}
implies that a $\C^\nu$-sublinear control system is a $\C^\nu$-control
system.

Let us demonstrate a class of sublinear control systems in which we will be
particularly interested.
\begin{example}\label{eg:control-linear}
The class of sublinear control systems we consider seems quite particular,
but will turn out to be extremely general in our framework.  We let
$m\in\integernn$ and $m'\in\{0,\lip\}$\@, let\/
$\nu\in\{m+m',\infty,\omega\}$\@, and let $r\in\{\infty,\omega\}$\@, as
required.  Let $\alg{V}$ be a locally convex topological vector space, and
let $\cs{C}\subset\alg{V}$\@.  We suppose that we have a continuous linear
map $\Lambda\in\lin{\alg{V}}{\sections[\nu]{\tb{\man{M}}}}$ and we
correspondingly define $\map{F_\Lambda}{\man{M}\times\cs{C}}{\tb{\man{M}}}$
by $F_\Lambda(x,u)=\Lambda(u)(x)$\@.  Continuity of $\Lambda$ immediately
gives that such the control system $(\man{M},F_\Lambda,\cs{C})$ is sublinear,
and we shall call a system such as this a \defn{$\C^\nu$-control-linear
system}\@.

Note that we can regard a control-affine system as a control-linear system as
follows.  For a control-affine system with $\cs{C}\subset\real^k$ and with
\begin{equation*}
F(x,\vect{u})=f_0(x)+\sum_{a=1}^ku^af_a(x),
\end{equation*}
we let $\alg{V}=\real^{k+1}\simeq\real\oplus\real^k$ and take
\begin{equation*}
\cs{C}'=\setdef{(u^0,\vect{u})\in\real\oplus\real^k}
{u^0=1,\ \vect{u}\in\cs{C}},\quad
\Lambda(u^0,\vect{u})=\sum_{a=0}^ku^af_a.
\end{equation*}
Clearly we have $F(x,\vect{u})=F_\Lambda(x,(1,\vect{u}))$ for every
$\vect{u}\in\cs{C}$\@.  Since linear maps from finite-dimensional locally
convex spaces are continuous~\cite[Proposition~2.10.2]{JH:66}\@, we conclude
that control-affine systems are control-linear systems.  Thus they are also
control systems as per Definition~\ref{def:Cr-system}\@.\oprocend
\end{example}

One may want to regard the generalisation from the case where the control set
is a subset of $\real^k$ to being a subset of a locally convex topological
vector space to be mere fancy generalisation, but this is, actually, far from
being the case as we shall see in Section~\ref{sec:gcs}\@.

We also have a version of Proposition~\ref{prop:open-loopcsI} for sublinear control systems.
\begin{proposition}\label{prop:open-loopcsII}
Let\/ $m\in\integernn$ and\/ $m'\in\{0,\lip\}$\@, let\/
$\nu\in\{m+m',\infty,\omega\}$\@, and let\/ $r\in\{\infty,\omega\}$\@, as
required.  Let\/ $\Sigma=(\man{M},F,\cs{C})$ be a $\C^\nu$-sublinear control
system for which\/ $\cs{C}$ is a subset of a locally convex topological
vector space\/ $\alg{V}$\@.  If\/ $\mu\in\Lloc^1(\mathbb{T};\cs{C})$\@,
then\/ $F^\mu\in\LIsections[\nu]{\tdomain;\tb{\man{M}}}$\@, where\/
$\map{F^\mu}{\tdomain\times\man{M}}{\tb{\man{M}}}$ is defined by\/
$F^\mu(t,x)=F(x,\mu(t))$\@.
\begin{proof}
The proof that $F^\mu$ is a Carath\'eodory vector field of class $\C^\nu$
goes exactly as in Proposition~\ref{prop:open-loopcsI}\@.

To prove that $F^\mu\in\LIsections[\nu]{\tdomain;\tb{\man{M}}}$\@, let
$K\subset\man{M}$ be compact, let $k\in\integernn$\@, let
$\vect{a}\in\csd(\integernn;\realp)$\@, and denote
\begin{equation*}
p_K=\begin{cases}p^\infty_{K,k},&\nu=\infty,\\p^m_K,&\nu=m,\\
p^{m+\lip}_K,&\nu=m+\lip,\\p^\omega_{K,\vect{a}},&\nu=\omega.\end{cases}
\end{equation*}
Define $\map{g}{\tdomain}{\realnn}$ by $g(t)=p_K(F^\mu_t)$\@.  We claim that
$g\in\Lloc^\infty(\tdomain;\realnn)$\@.  From the first part of the proof of
Proposition~\ref{prop:open-loopcsI}\@, $t\mapsto F^\mu_t(x)$ is measurable
for every $x\in\man{M}$\@.  By
Theorems~\ref{the:Cinfty-tvsec}\@,~\ref{the:Cmm'-tvsec}\@,
and~\ref{the:Comega-tvsec}\@, it follows that $t\mapsto F^\mu_t$ is
measurable.  Since $p_K$ is a continuous function on
$\sections[\nu]{\tb{\man{M}}}$\@, it follows that $t\mapsto p_K(F^\mu_t)$ is
measurable, as claimed.  We claim that $g\in\Lloc^1(\tdomain;\realnn)$\@.
Note that $X\mapsto p_K(X)$ is a continuous seminorm on
$\sections[\infty]{\tb{\man{M}}}$\@.  By hypothesis, there exists a
continuous seminorm $q$ for the locally convex topology for $\alg{V}$ such
that
\begin{equation*}
p_K(F^{u_1}-F^{u_2})\le q(u_1-u_2)
\end{equation*}
for every $u_1,u_2\in\cs{C}$\@.  Therefore, if $\tdomain'\subset\tdomain$ is
compact and if $u_0\in\cs{C}$\@, we also have
\begin{align*}
\int_{\tdomain'}g(t)\,\d{t}=&\;\int_{\tdomain'}p_K(F^\mu_t)\\
\le&\;\int_{\tdomain'}p_K(F^\mu_t-F^{u_0})\,\d{t}+
\int_{\tdomain'}p_K(F^{u_0})\,\d{t}\\
\le&\int_{\tdomain'}q(\mu(t))\,\d{t}+(q(u_0)+
p_K(F^{u_0}))\lebmes(\tdomain')<\infty,
\end{align*}
the last inequality by the characterisation of Bochner integrability
from~\cite[Theorems~3.2 and~3.3]{RB/AD:11}\@.  Thus $g$ is locally
integrable.  It follows from
Theorems~\ref{the:Cinfty-tvsec}\@,~\ref{the:Cmm'-tvsec}\@,
and~\ref{the:Comega-tvsec} that
$F^\mu\in\LIsections[\nu]{\tdomain;\tb{\man{M}}}$\@, as desired.
\end{proof}
\end{proposition}

There is also a version of the notion of trajectory that is applicable to the
case when the control set is a subset of a locally convex topological space.
\begin{definition}
Let $m\in\integernn$ and $m'\in\{0,\lip\}$\@, let
$\nu\in\{m+m',\infty,\omega\}$\@, and let $r\in\{\infty,\omega\}$\@, as
required.  Let $\Sigma=(\man{M},F,\cs{C})$ be a $\C^\nu$-control system.  For
an interval $\tdomain\subset\real$\@, a \defn{$\tdomain$-trajectory} is a
locally absolutely continuous curve $\map{\xi}{\tdomain}{\man{M}}$ for which
there exists $\mu\in\Lloc^1(\tdomain;\cs{C})$ such that
\begin{equation*}
\xi'(t)=F(\xi(t),\mu(t)),\qquad\ae\ t\in\tdomain.
\end{equation*}
The set of $\tdomain$-trajectories we denote by $\Traj(\tdomain,\Sigma)$\@.
If $\nbhd{U}$ is open, we denote by $\Traj(\tdomain,\nbhd{U},\Sigma)$ those
trajectories taking values in\/ $\nbhd{U}$\@.\oprocend
\end{definition}

\subsection{Differential inclusions}\label{subsec:di}

We briefly mentioned differential inclusions in
Section~\ref{subsubsec:di-intro}\@, but now let us define them properly and
give a few attributes of, and constructions for, differential inclusions of
which we shall subsequently make use.

First the definition.
\begin{definition}
For a smooth manifold $\man{M}$\@, a \defn{differential inclusion} on
$\man{M}$ is a set-valued map $\setmap{\sX}{\man{M}}{\tb{\man{M}}}$ with
nonempty values for which $\sX(x)\subset\tb[x]{\man{M}}$\@.  A
\defn{trajectory} for a differential inclusion $\sX$ is a locally absolutely
continuous curve $\map{\xi}{\tdomain}{\man{M}}$ defined on an interval
$\tdomain\subset\real$ for which $\xi'(t)\in\sX(\xi(t))$ for almost every
$t\in\tdomain$\@.  If $\tdomain\subset\real$ is an interval and if
$\nbhd{U}\subset\man{M}$ is open, by $\Traj(\tdomain,\nbhd{U},\sX)$ we denote
the trajectories of $\sX$ defined on $\tdomain$ and taking values in
$\nbhd{U}$\@.\oprocend
\end{definition}

Of course, differential inclusions will generally not have trajectories, and
to ensure that they do various hypotheses can be made.  Two common attributes
of differential inclusions in this vein are the following.
\begin{definition}
A differential inclusion $\sX$ on a smooth manifold $\man{M}$ is:
\begin{compactenum}[(i)]
\item \defn{lower semicontinuous} at $x_0\in\man{M}$ if, for any
$v_0\in\sX(x_0)$ and any neighbourhood $\nbhd{V}\subset\tb{\man{M}}$ of
$v_0$\@, there exists a neighbourhood $\nbhd{U}\subset\man{M}$ of $x_0$ such
that $\sX(x)\cap\nbhd{V}\not=\emptyset$ for every $x\in\nbhd{U}$\@;
\item \defn{lower semicontinuous} if it is lower semicontinuous at every $x\in\man{M}$\@;
\item \defn{upper semicontinuous} at $x_0\in\man{M}$ if, for every open set
$\tb{\man{M}}\supset\nbhd{V}\supset\sX(x_0)$\@, there exists a neighbourhood
$\nbhd{U}\subset\man{M}$ of $x_0$ such that $\sX(\nbhd{U})\subset\nbhd{V}$\@;
\item \defn{upper semicontinuous} if it is upper semicontinuous at each $x\in\man{M}$\@;
\item \defn{continuous} at $x_0\in\man{M}$ if it is both lower and upper
semicontinuous at $x_0$\@;
\item \defn{continuous} if it is both lower and upper
semicontinuous.\oprocend
\end{compactenum}
\end{definition}

Other useful properties of differential inclusions are the following.
\begin{definition}
A differential inclusion $\sX$ on a smooth manifold $\man{M}$ is:
\begin{compactenum}[(i)]
\item \defn{closed-valued} (\resp~\defn{compact-valued}\@,
\defn{convex-valued}) at $x\in\man{M}$ if $\sX(x)$ is closed (\resp, compact,
convex);
\item \defn{closed-valued} (\resp~\defn{compact-valued}\@,
\defn{convex-valued}) if $\sX(x)$ is closed (\resp, compact, convex) for
every $x\in\man{M}$\@.\oprocend
\end{compactenum}
\end{definition}

Some standard hypotheses for existence of trajectories are then:
\begin{compactenum}
\item $\sX$ is lower semicontinuous with closed and convex values~\cite[Theorem~2.1.1]{JPA/AC:84}\@;
\item $\sX$ is upper semicontinuous with compact and convex values~\cite[Theorem~2.1.4]{JPA/AC:84}\@;
\item $\sX$ is continuous with compact
values~\cite[Theorem~2.2.1]{JPA/AC:84}\@.
\end{compactenum}
These are not matters with which we shall be especially concerned.

A standard operation is to take ``hulls'' of differential inclusions in the
following manner.
\begin{definition}\label{def:dihulls}
Let $r\in\{\infty,\omega\}$\@, let $\man{M}$ be a $\C^r$-manifold, and let
$\setmap{\sX}{\man{M}}{\tb{\man{M}}}$ be a differential inclusion.
\begin{compactenum}[(i)]
\item The \defn{convex hull} of $\sX$ is the differential inclusion
$\cohull(\sX)$ defined by
\begin{equation*}
\cohull(\sX)(x)=\cohull(\sX(x)),\qquad x\in\man{M}.
\end{equation*}
\item The \defn{closure} of $\sX$ is the differential inclusion
$\closure(\sX)$ defined by
\begin{equation*}\eqoprocend
\closure(\sX)(x)=\closure(\sX(x)),\qquad x\in\man{M}.
\end{equation*}
\end{compactenum}
\end{definition}

To close this section, let us make an observation regarding the connection
between control systems and differential inclusions.  Let $m\in\integernn$
and $m'\in\{0,\lip\}$\@, let $\nu\in\{m+m',\infty,\omega\}$\@, and let
$r\in\{\infty,\omega\}$\@, as required.  Let $\Sigma=(\man{M},F,\cs{C})$ be a
$\C^\nu$-control system.  To this system we associate the differential
inclusion $\sX_\Sigma$ by
\begin{equation*}
\sX_\Sigma(x)=\setdef{F^u(x)}{u\in\cs{C}}.
\end{equation*}
Since the differential inclusion $\sX_\Sigma$ is defined by a family of
vector fields, one might try to recover the vector fields $F^u$\@,
$u\in\cs{C}$\@, from $\sX_\Sigma$\@.  The obvious way to do this is to
consider
\begin{equation*}
\sections[\nu]{\sX_\Sigma}\eqdef\setdef{X\in\sections[\nu]{\tb{\man{M}}}}
{X(x)\in\sX_\Sigma(x),\ x\in\man{M}}.
\end{equation*}
Clearly we have $F^u\in\sections[\nu]{\sX_\Sigma}$ for every $u\in\cs{C}$\@.
However, $\sX_\Sigma$ will generally contain vector fields not of the form
$F^u$ for some $u\in\cs{C}$\@.  Let us give an illustration of this.  Let us
consider a smooth control system $(\man{M},F,\cs{C})$ with the following
properties:
\begin{compactenum}
\item $\cs{C}$ is a disjoint union of sets $\cs{C}_1$ and $\cs{C}_2$\@;
\item there exist disjoint open sets $\nbhd{U}_1$ and $\nbhd{U}_2$ such that
$\supp(F^u)\subset\nbhd{U}_1$ for $u\in\cs{C}_1$ and
$\supp(F^u)\subset\nbhd{U}_2$ for $u\in\cs{C}_2$\@.
\end{compactenum}
One then has that
\begin{equation*}
\setdef{c_1F^{u_1}+c_1F^{u_2}}{u_1\in\cs{C}_1,\ u_2\in\cs{C}_2,\
c_1,c_2\in\{0,1\},\ c_1^2+c_2^2\not=0}\subset\sections[\nu]{\sX_\Sigma},
\end{equation*}
showing that there are more sections of $\sX_\Sigma$ than there are control
vector fields.  This is very much related to presheaves and sheaves, to which
we shall now turn our attention.

\section{\Gcs{}s:\ Definitions and fundamental properties}\label{sec:gcs}

In this section we introduce the class of control systems we propose as being
useful mathematical models for the investigation of geometric system
structure.  The reader would do well to remember that this definition makes
no pretences of being simple or user-friendly.  However, we can do some
interesting things with these models, and to illustrate this we present in
Section~\ref{subsec:subriemannian} an elegant formulation of sub-Riemannian
geometry in the framework of \gcs{}s.

\subsection{Presheaves and sheaves of sets of vector fields}\label{subsec:sheaves}

We choose to phrase our notion of control systems in the language of sheaf
theory.  This will seem completely pointless to a reader not used to thinking
in this sort of language.  However, we do believe there are benefits to the
sheaf approach including~(1)~sheaves are the proper framework for
constructing germs of control systems which are often important in the study
of local system structure and~(2)~sheaf theory provides us with a natural
class of mappings between systems that we use to advantage in
Section~\ref{subsec:gcs-category}\@.

We do not even come close to discussing sheaves in any generality; we merely
give the definitions we require, a few of the most elementary consequences
of these definitions, and some representative (for us) examples.
\begin{definition}
Let $m\in\integernn$ and $m'\in\{0,\lip\}$\@, let
$\nu\in\{m+m',\infty,\omega\}$\@, and let $r\in\{\infty,\omega\}$\@, as
required.  Let $\man{M}$ be a manifold of class $\C^r$\@.  A \defn{presheaf
of sets of $\C^\nu$-vector fields} is an assignment to each open set
$\nbhd{U}\subset\man{M}$ a subset $\sF(\nbhd{U})$ of
$\sections[\nu]{\tb{\nbhd{U}}}$ with the property that, for open sets
$\nbhd{U},\nbhd{V}\subset\man{M}$ with $\nbhd{V}\subset\nbhd{U}$\@, the map
\begin{equation*}
\mapdef{r_{\nbhd{U},\nbhd{V}}}{\sF(\nbhd{U})}{\sections[\nu]{\tb{\nbhd{V}}}}
{X}{X|\nbhd{V}}
\end{equation*}
takes values in $\sF(\nbhd{V})$\@.  Elements of $\sF(\nbhd{U})$ are called
\defn{local sections} over $\nbhd{U}$\@.\oprocend
\end{definition}

Let us give some notation to the presheaf of sets of vector fields of which
every other such presheaf is a subset.
\begin{example}
Let $m\in\integernn$ and $m'\in\{0,\lip\}$\@, let
$\nu\in\{m+m',\infty,\omega\}$\@, and let $r\in\{\infty,\omega\}$\@, as
required.  Let $\man{M}$ be a manifold of class $\C^r$\@.  The presheaf of
\emph{all} vector fields of class $\C^\nu$ is denoted by
$\ssections[\nu]{\tb{\man{M}}}$\@.  Thus
$\ssections[\nu]{\tb{\man{M}}}(\nbhd{U})=\sections[\nu]{\tb{\nbhd{U}}}$ for
every open set $\nbhd{U}$\@.  Presheaves such as this are extremely important
in the ``normal'' applications of sheaf theory.  For those with some
background in these more standard applications of sheaf theory, we mention
that our reasons for using the theory are not quite the usual ones.  Such
readers will be advised to be careful not to overlay too much of their past
experience on what we do with sheaf theory here.\oprocend
\end{example}

The preceding notion of a presheaf is intuitively clear, but it does have
some defects.  One of these defects is that one can describe local data that
does not patch together to give global data.  Let us illustrate this with a
few examples.
\begin{examples}\label{eg:!sheaf}
\begin{compactenum}
\item \label{enum:bounded-presheaf} Let $m\in\integernn$ and
$m'\in\{0,\lip\}$\@, let $\nu\in\{m+m',\infty,\omega\}$\@, and let
$r\in\{\infty,\omega\}$\@, as required.  Let us take a manifold $\man{M}$ of
class $\C^r$ with a Riemannian metric $\metric$\@.  Let us define a presheaf
$\subscr{\sF}{bdd}$ by asking that
\begin{equation*}
\subscr{\sF}{bdd}(\nbhd{U})=\setdef{X\in\sections[\nu]{\tb{\man{M}}}}
{\sup\setdef{\dnorm{X(x)}_{\metric}}{x\in\nbhd{U}}<\infty}.
\end{equation*}
Thus $\subscr{\sF}{bdd}$ is comprised of vector fields that are ``bounded.''
This is a perfectly sensible requirement.  However, the following phenomenon
can happen if $\man{M}$ is not compact.  There can exist an open cover
$\ifam{\nbhd{U}_a}_{a\in A}$ for $\man{M}$ and local sections
$X_a\in\subscr{\sF}{bdd}(\nbhd{U}_a)$ that are ``compatible'' in the sense
that $X_a|\nbhd{U}_a\cap\nbhd{U}_b=X_b|\nbhd{U}_a\cap\nbhd{U}_b$\@, for each
$a,b\in A$\@, but such that there is no globally defined section
$X\in\subscr{\sF}{bdd}(\man{M})$ such that $X|\nbhd{U}_a=X_a$ for every $a\in
A$\@.  We leave to the reader the easy job of coming up with a concrete
instance of this.

\item \label{enum:global-sheaf} Let $m\in\integernn$ and $m'\in\{0,\lip\}$\@,
let $\nu\in\{m+m',\infty,\omega\}$\@, and let $r\in\{\infty,\omega\}$\@, as
required.  Let $\man{M}$ be a manifold of class $\C^r$\@.  If
$\sX\subset\sections[\nu]{\tb{\man{M}}}$ is any family of vector fields on
$\man{M}$\@, then we can define an associated presheaf $\sF_{\sX}$ of sets of
vector fields by
\begin{equation*}
\sF_{\sX}(\nbhd{U})=\setdef{X|\nbhd{U}}{X\in\sX}.
\end{equation*}
Note that $\sF(\man{M})$ is necessarily equal to $\sX$\@, and so we shall
typically use $\sF(\man{M})$ to denote the set of globally defined vector
fields giving rise to this presheaf.  A presheaf of this sort will be called
\defn{globally generated}\@.

This sort of presheaf will almost never have nice ``local to global''
properties.  Let us illustrate why this is so.  Let $\man{M}$ be a connected
Hausdorff manifold.  Suppose that the set of globally defined vector fields
$\sF(\man{M})$ has cardinality strictly larger than~$1$ and has the following
property: there exists a disconnected open set $\nbhd{U}\subset\man{M}$ such
that the mapping from $\sF(\nbhd{U})$ to $\sF(\man{M})$ given by
$X|\nbhd{U}\mapsto X$ is injective.  This property will hold for real
analytic families of vector fields, because we can take as $\nbhd{U}$ the
union of a pair of disconnected open sets.  However, the property will also
hold for many reasonable smooth families of vector fields.

We write $\nbhd{U}=\nbhd{U}_1\cup\nbhd{U}_2$ for disjoint open sets
$\nbhd{U}_1$ and $\nbhd{U}_2$\@.  By hypothesis, there exist vector fields
$X_1,X_2\in\sF(\man{M})$ such that $X_1|\nbhd{U}\not=X_2|\nbhd{U}$\@.  Define
local sections $X'_a\in\sF(\nbhd{U}_a)$ by $X'_a=X_a|\nbhd{U}_a$\@,
$a\in\{1,2\}$\@.  The condition
\begin{equation*}
X'_1|\nbhd{U}_1\cap\nbhd{U}_2=X'_2|\nbhd{U}_1\cap\nbhd{U}_2
\end{equation*}
is vacuously satisfied.  But there can be no $X\in\sF(\man{M})$ such that, if
$X'\in\sF(\nbhd{U})$ is given by $X'=X|\nbhd{U}$\@, then $X'|\nbhd{U}_1=X'_1$
and $X'|\nbhd{U}_2=X'_2$\@.

While a globally generated presheaf is unlikely to allow patching from local
to global, this can be easily redressed by undergoing a process known as
``sheafification'' that we will describe below.\oprocend
\end{compactenum}
\end{examples}

The preceding examples suggest that if one wishes to make compatible local
constructions that give rise to a global construction, additional properties
need to be ascribed to a presheaf of sets of vector fields.  This we do as
follows.
\begin{definition}
Let $m\in\integernn$ and $m'\in\{0,\lip\}$\@, let
$\nu\in\{m+m',\infty,\omega\}$\@, and let $r\in\{\infty,\omega\}$\@, as
required.  Let $\man{M}$ be a manifold of class $\C^r$\@.  A presheaf $\sF$
of sets of $\C^\nu$-vector fields is a \defn{sheaf of sets of $\C^\nu$-vector
fields} if, for every open set $\nbhd{U}\subset\man{M}$\@, for every open
cover $\ifam{\nbhd{U}_a}_{a\in A}$ of $\nbhd{U}$\@, and for every choice of
local sections $X_a\in\sF(\nbhd{U}_a)$ satisfying
$X_a|\nbhd{U}_a\cap\nbhd{U}_b=X_b|\nbhd{U}_a\cap\nbhd{U}_b$\@, there exists
$X\in\sF(\nbhd{U})$ such that $X|\nbhd{U}_a=X_a$ for every $a\in
A$\@.\oprocend
\end{definition}

The condition in the definition is called the \defn{gluing condition}\@.
Readers familiar with sheaf theory will note the absence of the other
condition, sometimes called the separation condition, normally placed on a
presheaf in order for it to be a sheaf: it is automatically satisfied for
presheaves of sets of vector fields.

Many of the presheaves that we encounter will not be sheaves, as they will be
globally generated.  Thus let us give some examples of sheaves, just as a
point of reference.
\begin{examples}
\begin{compactenum}
\item Let $m\in\integernn$ and $m'\in\{0,\lip\}$\@, let
$\nu\in\{m+m',\infty,\omega\}$\@, and let $r\in\{\infty,\omega\}$\@, as
required.  Let $\man{M}$ be a $\C^r$-manifold.  The presheaf
$\ssections[\nu]{\tb{\man{M}}}$ of all $\C^\nu$-vector fields is a sheaf.  We
leave the simple and standard working out of this to the reader; it will
provide some facility in working with sheaf concepts for those not already
having this.
\item If instead of considering bounded vector fields as in
part Example~\enumdblref{eg:!sheaf}{enum:bounded-presheaf}\@, we consider the
presheaf of vector fields satisfying a \emph{fixed} bound, then the resulting
presheaf is a sheaf.  Let us be clear.  Let $m\in\integernn$ and
$m'\in\{0,\lip\}$\@, let $\nu\in\{m+m',\infty,\omega\}$\@, and let
$r\in\{\infty,\omega\}$\@, as required.  We let $\man{M}$ be a
$\C^r$-manifold with Riemannian metric $\metric$ and, for $B\in\realp$\@,
define a presheaf $\sF_{\le B}$ by
\begin{equation*}
\sF_{\le B}(\nbhd{U})=\setdef{X\in\sections[\nu]{\tb{\man{M}}}}
{\sup\setdef{\dnorm{X(x)}_{\metric}}{x\in\man{M}}\le B}.
\end{equation*}
The presheaf $\sF_{\le B}$ is a sheaf, as is easily verified.  In this case,
the local constraints for membership are compatible with a global one.
\item Let $m\in\integernn$ and $m'\in\{0,\lip\}$\@, let
$\nu\in\{m+m',\infty,\omega\}$\@, and let $r\in\{\infty,\omega\}$\@, as
required.  Let $\man{M}$ be a $\C^r$-manifold.  Let $A\subset\man{M}$ and
define a presheaf $\sI_A$ of sets of vector fields by
\begin{equation*}
\sI_A(\nbhd{U})=\setdef{X\in\sections[\nu]{\tb{\nbhd{U}}}}{X(x)=0,\ x\in A}.
\end{equation*}
This is a sheaf (again, we leave the verification to the reader) called the
\defn{ideal sheaf} of~$A$\@.\oprocend
\end{compactenum}
\end{examples}

Let us now turn to localising sheaves of sets of vector fields.  Let
$m\in\integernn$ and $m'\in\{0,\lip\}$\@, let
$\nu\in\{m+m',\infty,\omega\}$\@, and let $r\in\{\infty,\omega\}$\@, as
required.  Let $\man{M}$ be a $\C^r$-manifold, let $A\subset\man{M}$\@, and
let $\sN_A$ be the set of neighbourhoods of $A$ in $\man{M}$\@,~\ie~the open
subsets of $\man{M}$ containing $A$\@.  This is a directed set in the usual
way by inclusion,~\ie~$\nbhd{U}\preceq\nbhd{V}$ if
$\nbhd{V}\subset\nbhd{U}$\@.  Let $\sF$ be a sheaf of sets of $\C^\nu$-vector
fields.  The \defn{stalk} of $\sF$ over $A$ is the direct limit
$\dirlim_{\nbhd{U}\in\sN_A}\sF(\nbhd{U})$\@.  Let us be less cryptic about
this.  Let $\nbhd{U},\nbhd{V}\in\sN_A$\@, and let $X\in\sF(\nbhd{U})$ and
$Y\in\sF(\nbhd{V})$\@.  We say $X$ and $Y$ are \defn{equivalent} if there
exists $\nbhd{W}\subset\nbhd{U}\cap\nbhd{V}$ such that
$X|\nbhd{W}=Y|\nbhd{W}$\@.  The \defn{germ} of $X\in\sF(\nbhd{U})$ for
$\nbhd{U}\in\sN_A$ is the equivalence class of $X$ under this equivalence
relation.  If $\nbhd{U}\in\sN_A$ and if $X\in\sF(\nbhd{U})$\@, then we denote
by $[X]_A$ the equivalence class of $X$ in $\sF_A$\@.  The \defn{stalk} of
$\sF$ over $A$ is the set of all equivalence classes.  The stalk of $\sF$
over $A$ is denoted by $\sF_A$\@, and we write $\sF_{\{x\}}$ as $\sF_x$\@.

Let us now describe how a presheaf can be converted in a natural way into a
sheaf.  The description of how to do this for general presheaves is a little
complicated.  However, in the case we are dealing with here, we can be
explicit about this.
\begin{lemma}\label{lem:sheafify}
Let\/ $m\in\integernn$ and\/ $m'\in\{0,\lip\}$\@, let\/
$\nu\in\{m+m',\infty,\omega\}$\@, and let\/ $r\in\{\infty,\omega\}$\@, as
required.  Let\/ $\man{M}$ be a\/ $\C^r$-manifold and let\/ $\sF$ be a
presheaf of sets of\/ $\C^\nu$-vector fields.  For an open set\/
$\nbhd{U}\subset\man{M}$\@, define
\begin{equation*}
\Sh(\sF)(\nbhd{U})=\setdef{X\in\sections[\nu]{\tb{\nbhd{U}}}}{[X]_x\in\sF_x\
\textrm{for every}\ x\in\nbhd{U}}.
\end{equation*}
Then\/ $\Sh(\sF)$ is a sheaf.
\begin{proof}
Let $\nbhd{U}\subset\man{M}$ be open and let $\ifam{\nbhd{U}_a}_{a\in A}$ be
an open cover of $\nbhd{U}$\@.  Suppose that local sections
$X_a\in\Sh(\sF)(\nbhd{U}_a)$\@, $a\in A$\@, satisfy
$X_a|\nbhd{U}_a\cap\nbhd{U}_b=X_b|\nbhd{U}_a\cap\nbhd{U}_b$ for each $a,b\in
A$\@.  Since $\ssections[\nu]{\tb{\man{M}}}$ is a sheaf, there exists
$X\in\sections[\nu]{\tb{\nbhd{U}}}$ such that $X|\nbhd{U}_a=X_a$\@, $a\in
A$\@.  It remains to show that $X\in\Sh(\sF)(\nbhd{U})$\@.  Let
$x\in\nbhd{U}$ and let $a\in A$ be such that $x\in\nbhd{U}_a$\@.  Then we
have $[X]_x=[X_a]_x\in\sF_x$\@, as desired.
\end{proof}
\end{lemma}

With the lemma in mind we have the following definition.
\begin{definition}
Let $m\in\integernn$ and $m'\in\{0,\lip\}$\@, let
$\nu\in\{m+m',\infty,\omega\}$\@, and let $r\in\{\infty,\omega\}$\@, as
required.  Let $\man{M}$ be a $\C^r$-manifold and let $\sF$ be a presheaf of
sets of $\C^\nu$-vector fields.  The \defn{sheafification} of $\sF$ is the
sheaf $\Sh(\sF)$ of sets of vector fields defined by
\begin{equation*}\eqoprocend
\Sh(\sF)(\nbhd{U})=\setdef{X\in\sections[\nu]{\tb{\nbhd{U}}}}
{[X]_x\in\sF_x\ \text{for all}\ x\in\nbhd{U}}.
\end{equation*}
\end{definition}

Let us consider some examples of sheafification.
\begin{examples}\label{eg:sheafify}
\begin{compactenum}
\item \label{enum:shfybdd} Let us consider the presheaf of bounded vector
fields from Example~\enumdblref{eg:!sheaf}{enum:bounded-presheaf}\@.  Let
$m\in\integernn$ and $m'\in\{0,\lip\}$\@, let
$\nu\in\{m+m',\infty,\omega\}$\@, and let $r\in\{\infty,\omega\}$\@, as
required.  Let $\man{M}$ be a $\C^r$-manifold and consider the presheaf
$\subscr{\sF}{bdd}$ of bounded vector fields.  One easily sees that the stalk
of this presheaf at $x\in\man{M}$ is given by
\begin{equation*}
\subscr{\sF}{bdd,$x$}=\setdef{[X]_x}{X\in\sections[\nu]{\tb{\man{M}}}},
\end{equation*}
\ie~there are no restrictions on the stalks coming from the boundedness
restriction on vector fields.  Therefore, $\Sh(\subscr{\sF}{bdd})=\ssections[\nu]{\tb{\man{M}}}$\@.

\item Let us now examine the sheafification of a globally generated presheaf
of sets of vector fields as in Example~\enumdblref{eg:!sheaf}{enum:global-sheaf}\@.
Let $m\in\integernn$ and $m'\in\{0,\lip\}$\@, let
$\nu\in\{m+m',\infty,\omega\}$\@, and let $r\in\{\infty,\omega\}$\@, as
required.  Let $\man{M}$ be a $\C^r$-manifold and let $\sF$ be a globally
generated presheaf of sets of $\C^\nu$-vector fields, with $\sF(\man{M})$ the
global generators.  We will contrast $\sF(\nbhd{U})$ with
$\Sh(\sF)(\nbhd{U})$ to get an understanding of what the sheaf $\Sh(\sF)$
``looks like.''

To do so, for $\nbhd{U}\subset\man{M}$ open and for
$X\in\sections[\nu]{\tb{\nbhd{U}}}$\@, let us define a set-valued map
$\setmap{\kappa_{X,\nbhd{U}}}{\nbhd{U}}{\sF(\man{M})}$ by
\begin{equation*}
\kappa_{X,\nbhd{U}}(x)=\setdef{X'\in\sF(\man{M})}{X'(x)=X(x)}.
\end{equation*}
Generally, since we have asked nothing of the vector field $X$\@, we might
have $\kappa_{X,\nbhd{U}}(x)=\emptyset$ for a chosen $x$\@, or for some
$x$\@, or for every $x$\@.  If, however, we take $X\in\sF(\nbhd{U})$\@, then
$X=X'|\nbhd{U}$ for some $X'\in\sF(\man{M})$\@.  Therefore, there exists a
constant selection of $\kappa_{X,\nbhd{U}}$\@,~\ie~a constant function
$\map{s}{\nbhd{U}}{\sF(\man{M})}$ such that $s(x)\in\kappa_{X,\nbhd{U}}(x)$
for every $x\in\nbhd{U}$\@.  Note that if, for example, $\man{M}$ is
connected and $\nu=\omega$\@, then there will be a \emph{unique} such
constant selection since a real analytic vector field known on an open subset
uniquely determines the vector field on the connected component containing
this open set; this is the Identity Theorem,~\cf~\cite[Theorem~A.3]{RCG:90a}
in the holomorphic case and the same proof applies in the real analytic case.
Moreover, this constant selection in this case will completely characterise
$\kappa_{X,\nbhd{U}}$ in the sense that $\kappa_{X,\nbhd{U}}(x)=\{s(x)\}$\@.

Let us now contrast this with the character of the map $\kappa_{X,\nbhd{U}}$
for a local section $X\in\Sh(\sF)(\nbhd{U})$\@.  In this case, for each
$x\in\nbhd{U}$\@, we have $[X]_x=[X_x]_x$ for some $X_x\in\sF(\man{M})$\@.
Thus there exists a neighbourhood $\nbhd{V}_x\subset\nbhd{U}$ such that
$X|\nbhd{V}_x=X_x|\nbhd{V}_x$\@.  What this shows is that there is a locally
constant selection of $\kappa_{X,\nbhd{U}}$\@,~\ie~a locally constant map
$\map{s}{\nbhd{U}}{\sF(\man{M})}$ such that $s(x)\in\kappa_{X,\nbhd{U}}(x)$
for each $x\in\nbhd{U}$\@.  As above, in the real analytic case when
$\man{M}$ is connected, this locally constant selection is uniquely
determined, and determines $\kappa_{X,\nbhd{U}}$ in the sense that $\kappa_{X,\nbhd{U}}(x)=\{s(x)\}$\@.

Note that locally constant functions are those that are constant on connected
components.  Thus, by passing to the sheafification, we have gained
flexibility by allowing local sections to differ on connected components of
an open set.  While this does not completely characterise the difference
between local sections of the globally generated sheaf $\sF$ and its
sheafification $\Sh(\sF)$\@, it captures the essence of the matter, and
\emph{does} completely characterise the difference when $\nu=\omega$ and
$\man{M}$ is connected.\oprocend
\end{compactenum}
\end{examples}

\subsection{\Gcs{}s}\label{subsec:gcs}

Our definition of a \gcs\ is relatively straightforward, given the
definitions of the preceding section.
\begin{definition}
Let $m\in\integernn$ and $m'\in\{0,\lip\}$\@, let
$\nu\in\{m+m',\infty,\omega\}$\@, and let $r\in\{\infty,\omega\}$\@, as
required.
\begin{compactenum}[(i)]
\item A \defn{$\C^\nu$-\gcs} is a pair $\fG=(\man{M},\sF)$\@, where $\man{M}$
is a manifold of class $\C^r$ whose elements are called \defn{states} and
where $\sF$ is a presheaf of sets of $\C^\nu$-vector fields on $\man{M}$\@.
\item A \gcs\ $\fG=(\man{M},\sF)$ is \defn{complete} if
$\sF$ is a sheaf and is \defn{globally generated} if $\sF$ is globally
generated.
\item The \defn{completion} of $\fG=(\man{M},\sF)$ is the \gcs\
$\Sh(\fG)=(\man{M},\Sh(\sF))$\@.\oprocend
\end{compactenum}
\end{definition}

This is a pretty featureless definition, sorely in need of some connection to
control theory.  Let us begin to build this connection by pointing out the
manner in which more common constructions give rise to \gcs{}s, and vice
versa.
\begin{examples}\label{eg:gcs}
One of the topics of interest to us will be the relationship between our
notion of \gcs{}s and the more common notions of control systems (as in
Sections~\ref{subsec:Cr-systemsI} and~\ref{subsec:Cr-systemsII}) and
differential inclusions (as in Section~\ref{subsec:di}).  We begin here by
making some more or less obvious associations.
\begin{compactenum}
\item \label{enum:cs->gcs} Let $m\in\integernn$ and $m'\in\{0,\lip\}$\@, let
$\nu\in\{m+m',\infty,\omega\}$\@, and let $r\in\{\infty,\omega\}$\@, as
required.  Let $\Sigma=(\man{M},F,\cs{C})$ be a $\C^\nu$-control system.  To
this control system we associate the $\C^\nu$-\gcs\
$\fG_\Sigma=(\man{M},\sF_\Sigma)$ by
\begin{equation*}
\sF_\Sigma(\nbhd{U})=
\setdef{F^u|\nbhd{U}\in\sections[\nu]{\tb{\nbhd{U}}}}{u\in\cs{C}}.
\end{equation*}
The presheaf of sets of vector fields in this case is of the globally
generated variety, as in Example~\enumdblref{eg:!sheaf}{enum:global-sheaf}\@.
According to Example~\enumdblref{eg:!sheaf}{enum:global-sheaf} we should generally
not expect \gcs{}s such as this to be \emph{a priori} complete.  We can,
however, sheafify so that the \gcs\ $\Sh(\fG_\Sigma)$ is complete.

\item \label{enum:gcs->cs} Let us consider a means of going from a large
class of \gcs{}s to a control system.  Let $m\in\integernn$ and
$m'\in\{0,\lip\}$\@, let $\nu\in\{m+m',\infty,\omega\}$\@, and let
$r\in\{\infty,\omega\}$\@, as required.  We suppose that we have a
$\C^\nu$-\gcs\ $\fG=(\man{M},\sF)$ where the presheaf $\sF$ is globally
generated.  We define a $\C^\nu$-control system
$\Sigma_{\fG}=(\man{M},F_{\sF},\cs{C}_{\sF})$ as follows.  We take
$\cs{C}_{\sF}=\sF(\man{M})$\@,~\ie~the control set is our family of globally
defined vector fields and the topology is that induced from
$\sections[\nu]{\tb{\man{M}}}$\@.  We define
\begin{equation*}
\mapdef{F_{\sF}}{\man{M}\times\cs{C}_{\sF}}{\tb{\man{M}}}{(x,X)}{X(x).}
\end{equation*}
(Note that one has to make an awkward choice between writing a vector field
as $u$ or a control as $X$\@, since vector fields are controls.  We have gone
with the latter awkward choice, since it more readily mandates thinking about
what the symbols mean.)  Note that $F_{\sF}^X=X$\@, and so this is somehow
the identity map in disguise.  In order for this construction to provide a
\emph{bona fide} control system, we should check that $F_{\sF}$ is a
parameterised vector field of class $\C^\nu$ according to our
Definitions~\ref{def:Cinfty-paramsec}\@,~\ref{def:Cmm'-paramsec}\@,
and~\ref{def:Comega-paramsec}\@.  According to
Propositions~\ref{prop:paramCinfty}\@,~\ref{prop:paramCmm'}\@,
and~\ref{prop:paramComega}\@, it is sufficient to check that the map
$X\mapsto F_{\sF}^X$ is continuous.  But this is the identity map, which is
obviously continuous!

Note that $\Sigma_{\fG}$ is a control-linear system, according to
Example~\ref{eg:control-linear}\@.

\item \label{enum:di->gcs} Let $m\in\integernn$ and $m'\in\{0,\lip\}$\@, let
$\nu\in\{m+m',\infty,\omega\}$\@, and let $r\in\{\infty,\omega\}$\@, as
required.  Let $\setmap{\sX}{\man{M}}{\tb{\man{M}}}$ be a differential
inclusion.  If $\nbhd{U}\subset\man{M}$ is open, we denote
\begin{equation*}
\sections[\nu]{\sX|\nbhd{U}}=\setdef{X\in\sections[\nu]{\tb{\nbhd{U}}}}
{X(x)\in\sX(x),\ x\in\nbhd{U}}.
\end{equation*}
One should understand, of course, that we may very well have
$\sections[\nu]{\sX|\nbhd{U}}=\emptyset$\@.  This might happen for two
reasons.
\begin{compactenum}[(a)]
\item First, the differential inclusion may lack sufficient regularity to
permit even local sections of the prescribed regularity.
\item Second, even if it permits local sections, there may be be problems
finding sections defined on ``large'' open sets, because there may be global
obstructions.  One might anticipate this to be especially problematic in the
real analytic case, where the specification of a vector field locally
determines its behaviour globally by the Identity
Theorem,~\cf~\cite[Theorem~A.3]{RCG:90a}\@.
\end{compactenum}
This caveat notwithstanding, we can go ahead and define a \gcs\
$\fG_{\sX}=(\man{M},\sF_{\sX})$ with
$\sF_{\sX}(\nbhd{U})=\sections[\nu]{\sX|\nbhd{U}}$\@.

We claim that $\fG_{\sX}$ is complete.  To see this, let
$\nbhd{U}\subset\man{M}$ be open and let $\ifam{\nbhd{U}_a}_{a\in A}$ be an
open cover for $\nbhd{U}$\@.  For each $a\in A$\@, let
$X_a\in\sF_{\sX}(\nbhd{U}_a)$ and suppose that, for $a,b\in A$\@,
\begin{equation*}
X_a|\nbhd{U}_a\cap\nbhd{U}_b=X_b|\nbhd{U}_a\cap\nbhd{U}_b.
\end{equation*}
Since $\ssections[\nu]{\tb{\man{M}}}$ is a sheaf, let
$X\in\sections[\nu]{\tb{\nbhd{U}}}$ be such that $X|\nbhd{U}_a=X_a$ for each
$a\in A$\@.  We claim that $X\in\sF_{\sX}(\nbhd{U})$\@.  Indeed, for
$x\in\nbhd{U}$ we have $X(x)=X_a(x)\in\sX(x)$ if we take $a\in A$ such that $x\in\nbhd{U}_a$\@.

The sheaf $\sF_{\sX}$ is not necessarily globally generated.  Here is a
stupid counterexample.  Let us define $\sX(x)=\tb[x]{\man{M}}$\@,
$x\in\man{M}$\@, so that $\sF_{\sX}=\ssections[\nu]{\tb{\man{M}}}$\@.  For an
open set $\nbhd{U}$\@, there will generally be local sections
$X\in\sections[\nu]{\tb{\nbhd{U}}}$ that are not restrictions to $\nbhd{U}$
of globally defined vector fields; vector fields that ``blow up'' at some
point in the boundary of $\nbhd{U}$ are what one should have in mind.

\item \label{enum:gcs->di} Let $m\in\integernn$ and $m'\in\{0,\lip\}$\@, let
$\nu\in\{m+m',\infty,\omega\}$\@, and let $r\in\{\infty,\omega\}$\@, as
required.  Note that there is also associated to any $\C^\nu$-\gcs\
$\fG=(\man{M},\sF)$ a differential inclusion $\sX_{\fG}$ by
\begin{equation*}
\sX_{\fG}(x)=\setdef{X(x)}{[X]_x\in\sF_x},
\end{equation*}
recalling that $\sF_x$ is the stalk of $\sF$ at $x$\@.\oprocend
\end{compactenum}
\end{examples}

Now note that we can iterate the four constructions and ask to what extent we
end up back where we started.  More precisely, we have the following result.
\begin{proposition}\label{prop:gcs-cs-di}
Let\/ $m\in\integernn$ and\/ $m'\in\{0,\lip\}$\@, let\/
$\nu\in\{m+m',\infty,\omega\}$\@, and let\/ $r\in\{\infty,\omega\}$\@, as
required.  Let\/ $\fG=(\man{M},\sF)$ be a\/ $\C^\nu$-\gcs, let\/
$\Sigma=(\man{M},F,\cs{C})$ be a\/ $\C^\nu$-control system, and let\/ $\sX$
be a differential inclusion.  Then the following statements hold:
\begin{compactenum}[(i)]
\item \label{pl:gcs->cs->gcs} if\/ $\fG$ is globally generated, then\/
$\fG_{\Sigma_{\fG}}=\fG$\@;
\item \label{pl:cs->gcs->cs} if the map\/ $u\mapsto F^u$ from\/ $\cs{C}$ to\/
$\sections[\nu]{\tb{\man{M}}}$ is injective and open onto its image, then\/
$\Sigma_{\fG_{\Sigma}}=\Sigma$\@;
\item \label{pl:gcs->di->gcs} $\sF(\nbhd{U})\subset\sF_{\sX_{\fG}}(\nbhd{U})$
for every open\/ $\nbhd{U}\subset\man{M}$\@;
\item \label{pl:di->gcs->di} $\sX_{\fG_{\sX}}\subset\sX$\@.
\end{compactenum}
\begin{proof}
\eqref{pl:gcs->cs->gcs} Let $\nbhd{U}\subset\man{M}$ be open and let
$X\in\sF(\nbhd{U})$\@.  Then $X=X'|\nbhd{U}$ for $X'\in\sF(\man{M})$\@.  Thus
$X'\in\cs{C}_{\sF}$ and $X'(x)=F(x,X')$ and so
$X\in\sF_{\Sigma_{\fG}}(\nbhd{U})$\@.  Conversely, let
$X\in\sF_{\Sigma_{\fG}}(\nbhd{U})$\@.  Then $X(x)=F(x,X')$\@,
$x\in\nbhd{U}$\@, for some $X'\in\cs{C}_{\sF}$\@.  But this means that
$X(x)=X'(x)$ for $X'\in\sF(\nbhd{U})$ and for all $x\in\nbhd{U}$\@.  In other
words, $X\in\sF(\nbhd{U})$\@.

\eqref{pl:cs->gcs->cs} Note that $\fG_\Sigma$ is globally generated.  Thus we
have
\begin{equation*}
\cs{C}_{\sF_{\Sigma}}=\sF_\Sigma(\man{M})=\setdef{F^u}{u\in\cs{C}}.
\end{equation*}
Since the map $u\mapsto F^u$ is continuous (by
Propositions~\ref{prop:paramCinfty}\@,~\ref{prop:paramCmm'}\@,
and~\ref{prop:paramComega}), and injective and open onto its image (by
hypothesis), it is an homeomorphism onto its image.  Thus
$\cs{C}_{\sF_\Sigma}$ is homeomorphic to $\cs{C}$\@.  Since $u\mapsto F^u$ is
injective we can unambiguously write
\begin{equation*}
F_{\sF_\Sigma}(x,F^u)=F^u(x)=F(x,u).
\end{equation*}

\eqref{pl:gcs->di->gcs} Let $\nbhd{U}\subset\man{M}$ be open.  If
$X\in\sF(\nbhd{U})$\@, then clearly we have $X(x)\in\sX_{\fG}(x)$ for every
$x\in\nbhd{U}$ and so $\sF(\nbhd{U})\subset\sF_{\sX_{\fG}}(\nbhd{U})$\@,
giving the assertion.

\eqref{pl:di->gcs->di} This is obvious.
\end{proof}
\end{proposition}

\begin{remark}
The result establishes the rather surprising correspondence between control
systems $\Sigma=(\man{M},F,\cs{C})$ for which the map $u\mapsto F^u$ is
injective and open onto its image, and the associated control-linear system
$\Sigma_{\fG_\Sigma}=(\man{M},\sF_{\Sigma_{\fG}},\cs{C}_{\sF_\Sigma})$\@.
That is to say, at least at the system level, in our treatment every system
corresponds in a natural way to a control-linear system, albeit with a rather
complicated control set.  This correspondence carries over to trajectories as
well, but one can also weaken these conditions to obtain trajectory
correspondence in more general situations.  These matters we discuss in
detail in Section~\ref{subsec:trajequiv}\@.\oprocend
\end{remark}

Let us make some comments on the hypotheses present in the preceding result.
\begin{remarks}\label{eg:backandforth}
\begin{compactenum}
\item Since $\fG_\Sigma$ is necessarily globally generated for any control
system $\Sigma$\@, the requirement that $\fG$ be globally generated cannot be
dropped in part~\eqref{pl:gcs->cs->gcs}\@.
\item The requirement that the map $u\mapsto F^u$ be injective in
part~\eqref{pl:cs->gcs->cs} cannot be relaxed.  Without this assumption,
there is no way to recover $F$ from $\setdef{F^u}{u\in\cs{C}}$\@.  Similarly,
if this map is not open onto its image, while there may be a bijection
between $\cs{C}$ and $\cs{C}_{\sF_\Sigma}$\@, it will not be an homeomorphism
which one needs for the control systems to be the same.
\item \label{enum:gcs->di->gcs} The converse assertion in
part~\eqref{pl:gcs->di->gcs} does not generally hold, as many counterexamples
show.  Here are two, each of a different character.
\begin{compactenum}[(a)]
\item We take $\man{M}=\real$ and consider the $\C^\omega$-\gcs\
$\fG=(\man{M},\sF)$ where $\sF$ is the globally generated presheaf defined by
the single vector field $x^2\pderiv{}{x}$\@.  Note that
\begin{equation*}
\sX_{\fG}(x)=\begin{cases}\{0\},&x=0,\\\tb[x]{\real},&x\not=0.\end{cases}
\end{equation*}
Therefore,
\begin{equation*}
\sF_{\sX}(\nbhd{U})=\begin{cases}
\setdef{X\in\sections[\omega]{\tb{\nbhd{U}}}}{X(0)=0},&0\in\nbhd{U},\\
\sections[\omega]{\tb{\nbhd{U}}},&0\not\in\nbhd{U}.\end{cases}
\end{equation*}
It holds, therefore, that the vector field $x\pderiv{}{x}$ is a global
section of $\sF_{\sX}$\@, but is not a global section of $\sF$\@.
\item Let us again take $\man{M}=\real$ and now define a smooth \gcs\
$\fG=(\man{M},\sF)$ by asking that $\sF$ be the globally generated presheaf
defined by the vector fields $X_1,X_2\in\sections[\infty]{\real}$\@, where
\begin{equation*}
X_1(x)=\begin{cases}\eul^{-1/x}\pderiv{}{x},&x>0,\\0,&x\le0,\end{cases}
\end{equation*}
and
\begin{equation*}
X_2(x)=\begin{cases}\eul^{-1/x}\pderiv{}{x},&x<0,\\0,&x\ge0.\end{cases}
\end{equation*}
In this case,
\begin{equation*}
\sX_{\fG}(x)=\begin{cases}\{0\},&x=0,\\\{0\}\cup\{\eul^{-1/x}\pderiv{}{x}\},&
x\not=0.\end{cases}
\end{equation*}
Therefore, $\sF_{\sX}$ is the sheafification of the globally generated
presheaf defined by the vector fields $X_1$\@, $X_2$\@, $X_3$\@, and $X_4$\@,
where
\begin{equation*}
X_3(x)=\begin{cases}\eul^{-1/x}\pderiv{}{x},&x\not=0,\\0,&x=0,\end{cases}
\end{equation*}
and $X_4$ is the zero vector field.
\end{compactenum}
\item \label{enum:geomdi} Given the discussion in
Example~\enumdblref{eg:gcs}{enum:di->gcs}\@, one cannot reasonably expect
that we will generally have equality in part~\eqref{pl:di->gcs->di} of the
preceding result.  Indeed, one might even be inclined to say that it is only
differential inclusions satisfying $\sX=\sX_{\fG_{\sX}}$ that are useful in
geometric control theory\ldots\oprocend
\end{compactenum}
\end{remarks}

While we are not yet finished with the task of formulating our
theory\textemdash{}trajectories have yet to appear\textemdash{}it is
worthwhile to make a pause at this point to reflect upon what we have done
and have not done.  After a moments thought, one realises that the difference
between a control system $\Sigma=(\man{M},F,\cs{C})$ and its associated \gcs\
$\fG_\Sigma=(\man{M},\sF_\Sigma)$ is that, in the former case, the control
vector fields are from the \emph{indexed family} $\ifam{F^u}_{u\in\cs{C}}$\@,
while for the \gcs\ we have the \emph{set} $\setdef{F^u}{u\in\cs{C}}$\@.  In
going from the former to the latter we have ``forgotten'' the index $u$ which
we are explicitly keeping track of for control systems.  If the map $u\mapsto
F^u$ is injective, as in
Proposition~\pldblref{prop:gcs-cs-di}{pl:cs->gcs->cs}\@, then there is no
information lost as one goes from the indexed family to the set.  If
$u\mapsto F^u$ is not injective, then this is a signal that the control set
is too large, and perhaps one should collapse it in some way.  In other
words, one can probably suppose injectivity of $u\mapsto F^u$ without loss of
generality.  (Openness of this map is another matter.  As we shall see in
Section~\ref{subsec:trajequiv} below, openness (and a little more) is crucial
for there to be trajectory correspondence between systems and \gcs{}s.)  This
then leaves us with the mathematical semantics of distinguishing between the
indexed family $\ifam{F^u}_{u\in\cs{C}}$ and the subset
$\setdef{F^u}{u\in\cs{C}}$\@.  About this, let us make two observations.
\begin{compactenum}
\item The entire edifice of nonlinear control theory seems, in some sense, to
be built upon the preference of the indexed family over the set.  As we
discuss in the introduction, in applications there are very good reasons for
doing this.  But from the point of view of the general theory, the idea that
one should carefully maintain the labelling of the vector fields from the set
$\setdef{F^u}{u\in\cs{C}}$ seems to be a really unnecessary distraction.
And, moreover, it is a distraction upon which is built the whole notion of
``feedback transformation,'' plus entire methodologies in control theory that
are not
feedback-invariant,~\eg~linearisation,~\cf~Example~\ref{eg:bad-linearise}\@.
So, semantics?  Possibly, but sometimes semantic choices are important.
\item Many readers will probably not be convinced by our attempts to magnify
the distinction between the indexed family $\ifam{F^u}_{u\in\cs{C}}$ and the
set $\setdef{F^u}{u\in\cs{C}}$\@.  As we shall see, however, this distinction
becomes more apparent if one is really dedicated to using sets rather than
indexed families.  Indeed, this deprives one of the notion of ``control,''
and one is forced to be more thoughtful about what one means by
``trajectory.''  It is to this more thoughtful undertaking that we now turn,
slowly.
\end{compactenum}

\subsection{Open-loop systems}\label{subsec:ol}

Trajectories are associated to ``open-loop systems,'' so we first discuss
these.  We first introduce some notation.  Let $m\in\integernn$ and
$m'\in\{0,\lip\}$\@, let $\nu\in\{m+m',\infty,\omega\}$\@, and let
$r\in\{\infty,\omega\}$\@, as required.  For a $\C^\nu$-\gcs\
$\fG=(\man{M},\sF)$\@, we then denote
\begin{equation*}
\LIsections[\nu]{\tdomain;\sF(\nbhd{U})}=
\setdef{\map{X}{\tdomain}{\sF(\nbhd{U})}}
{X\in\LIsections[\nu]{\tdomain;\tb{\nbhd{U}}}},
\end{equation*}
for $\tdomain\subset\real$ an interval and $\nbhd{U}\subset\man{M}$ open.
\begin{definition}
Let $m\in\integernn$ and $m'\in\{0,\lip\}$\@, let
$\nu\in\{m+m',\infty,\omega\}$\@, and let $r\in\{\infty,\omega\}$\@, as
required.  Let $\fG=(\man{M},\sF)$ be a $\C^\nu$-\gcs.  An \defn{open-loop
system} for $\fG$ is a triple $\subscr{\fG}{ol}=(X,\tdomain,\nbhd{U})$ where
\begin{compactenum}[(i)]
\item $\tdomain\subset\real$ is an interval called the \defn{time-domain}\@;
\item $\nbhd{U}\subset\man{M}$ is open;
\item $X\in\LIsections[\nu]{\tdomain;\sF(\nbhd{U})}$\@.\oprocend
\end{compactenum}
\end{definition}

Note that an open-loop system for $\fG=(\man{M},\sF)$ is also an open-loop
system for the completion $\Sh(\fG)$\@, just because
$\sF(\nbhd{U})\subset\Sh(\sF)(\nbhd{U})$\@.  However, of course, there may be
open-loop systems for $\Sh(\fG)$ that are not open-loop systems for $\fG$\@.
This is as it should be, and has no significant ramifications for the theory,
as we shall see as we go along.

In order to see how we should think about an open-loop system, let us
consider this notion in the special case of control systems.
\begin{example}\label{eg:cs->openloop}
Let $m\in\integernn$ and $m'\in\{0,\lip\}$\@, let
$\nu\in\{m+m',\infty,\omega\}$\@, and let $r\in\{\infty,\omega\}$\@, as
required.  Let $\Sigma=(\man{M},F,\cs{C})$ be a $\C^\nu$-control system with
$\fG_\Sigma$ the associated $\C^\nu$-\gcs.  If we let
$\mu\in\Lloc^\infty(\tdomain;\cs{C})$\@, then we have the associated
open-loop system $\fG_{\Sigma,\mu}=(F^\mu,\tdomain,\man{M})$ defined by
\begin{equation*}
F^\mu(t)(x)=F(x,\mu(t)),\qquad t\in\tdomain,\ x\in\man{M}.
\end{equation*}
Proposition~\ref{prop:open-loopcsI} ensures that this is an open-loop system
for the \gcs\ $\fG_\Sigma$\@.

A similar assertion holds if $\cs{C}$ is a subset of a locally convex
topological vector space and $F$ defines a sublinear control system, and if
$\mu\in\Lloc^1(\tdomain;\cs{C})$\@,~%
\cf~Proposition~\ref{prop:open-loopcsII}\@.\oprocend
\end{example}

\begin{notation}
For an open-loop system $\subscr{\fG}{ol}(X,\tdomain,\nbhd{U})$\@, the
notation $X(t)(x)$\@, while accurate, is unnecessarily cumbersome, and we
will often instead write $X(t,x)$ or $X_t(x)$\@, with no loss of clarity and
a gain in aesthetics.\oprocend
\end{notation}

Generally one might wish to place a restriction on the set of open-loop
systems one will use.  This is tantamount to, for usual control systems,
placing restrictions on the controls one might use; one may wish to use
piecewise continuous controls or piecewise constant controls, for example.
For \gcs{}s we do this as follows.
\begin{definition}\label{def:olsubfamily}
Let $m\in\integernn$ and $m'\in\{0,\lip\}$\@, let
$\nu\in\{m+m',\infty,\omega\}$\@, and let $r\in\{\infty,\omega\}$\@, as
required.  Let $\fG=(\man{M},\sF)$ be a $\C^\nu$-\gcs.  An \defn{open-loop
subfamily} for $\fG$ is an assignment, to each interval
$\tdomain\subset\real$ and each open set $\nbhd{U}\subset\man{M}$\@, a subset
$\sO_{\fG}(\tdomain,\nbhd{U})\subset\LIsections[\nu]{\tdomain;\sF(\nbhd{U})}$
with the property that, if $(\tdomain_1,\nbhd{U}_1)$ and
$(\tdomain_2,\nbhd{U}_2)$ are such that $\tdomain_1\subset\tdomain_2$ and
$\nbhd{U}_1\subset\nbhd{U}_2$\@, then
\begin{equation*}\eqoprocend
\setdef{X|\tdomain_1\times\nbhd{U}_1}{X\in\sO_{\fG}(\tdomain_2,\nbhd{U}_2)}
\subset\sO_{\fG}(\tdomain_1,\nbhd{U}_1).
\end{equation*}
\end{definition}

Here are a few common examples of open-loop subfamilies.
\begin{examples}\label{eg:olsf}
Let $m\in\integernn$ and $m'\in\{0,\lip\}$\@, let
$\nu\in\{m+m',\infty,\omega\}$\@, and let $r\in\{\infty,\omega\}$\@, as
required.  Let $\fG=(\man{M},\sF)$ be a $\C^\nu$-\gcs.
\begin{compactenum}
\item The \defn{full subfamily} for $\fG$ is the open-loop subfamily
$\subscr{\sO}{$\fG$,full}$ defined by
\begin{equation*}
\subscr{\sO}{$\fG$,full}(\tdomain,\nbhd{U})=
\LIsections[\nu]{\tdomain;\sF(\nbhd{U})}.
\end{equation*}
Thus the full subfamily contains all possible open-loop systems.  Of course,
every open-loop subfamily will be contained in this one.

\item The \defn{locally essentially bounded subfamily} for $\fG$ is the
open-loop subfamily $\sO_{\fG,\infty}$ defined by asking that
\begin{equation*}
\sO_{\fG,\infty}(\tdomain,\nbhd{U})=
\setdef{X\in\subscr{\sO}{$\fG$,full}(\tdomain,\nbhd{U})}
{X\in\LBsections[\nu]{\tdomain;\tb{\nbhd{U}}}}.
\end{equation*}
Thus, for the locally essentially bounded subfamily, we require that the
condition of being locally integrally $\C^\nu$-bounded be replaced with the
stronger condition of being locally essentially $\C^\nu$-bounded.

\item \label{enum:cptsubfamily} The \defn{locally essentially compact
subfamily} for $\fG$ is the open-loop subfamily $\subscr{\sO}{$\fG$,cpt}$
defined by asking that
\begin{multline*}
\subscr{\sO}{$\fG$,cpt}(\tdomain,\nbhd{U})=
\{X\in\subscr{\sO}{$\fG$,full}(\tdomain,\nbhd{U})|\enspace
\textrm{for every compact subinterval}\ \tdomain'\subset\tdomain\\\
\textrm{there exists a compact}\
K\subset\sections[\nu]{\tdomain;\tb{\nbhd{U}}}\\
\textrm{such that}\ X(t)\subset K\ \textrm{for almost every}\ t\in\tdomain'\}.
\end{multline*}
Thus, for the locally essentially compact subfamily, we require that the
condition of being locally essentially bounded in the von Neumann bornology
(that defines the locally essentially bounded subfamily) be replaced with
being locally essentially bounded in the compact bornology.

We comment that in cases when the compact and von Neumann bornologies agree,
then of course we have $\sO_{\fG,\infty}=\subscr{\sO}{$\fG$,cpt}$\@.  As we
have seen in $\CO^\infty$-\ref{enum:COinfty-nuclear} and
$\C^\omega$-\ref{enum:Comega-nuclear}\@, this is the case when
$\nu\in\{\infty,\omega\}$\@.

\item The \defn{piecewise constant subfamily} for $\fG$ is the open-loop
subfamily $\subscr{\sO}{$\fG$,pwc}$ defined by asking that
\begin{equation*}
\subscr{\sO}{$\fG$,pwc}(\tdomain;\nbhd{U})=
\setdef{X\in\subscr{\sO}{$\fG$,full}(\tdomain,\nbhd{U})}
{t\mapsto X(t)\ \textrm{is piecewise constant}}.
\end{equation*}
Let us be clear what we mean by piecewise constant.  We mean that there is a
partition $\ifam{\tdomain_j}_{j\in J}$ of $\tdomain$ into pairwise disjoint
intervals such that
\begin{compactenum}[(a)]
\item for any compact interval $\mathbb{T}'\subset\tdomain$\@, the set
\begin{equation*}
\setdef{j\in J}{\mathbb{T}'\cap\tdomain_j\not=\emptyset}
\end{equation*}
is finite and such that
\item $X|\tdomain_j$ is constant for each $j\in J$\@.
\end{compactenum}
One might imagine that the piecewise constant open-loop subfamily will be
useful for studying orbits and controllability of \gcs{}s.

\item \label{enum:X-ol} We can associate an open-loop subfamily to an
open-loop system as follows.  Let $m\in\integernn$ and $m'\in\{0,\lip\}$\@,
let $\nu\in\{m+m',\infty,\omega\}$\@, and let $r\in\{\infty,\omega\}$\@, as
required.  Let $\fG=(\man{M},\sF)$ be a $\C^\nu$-\gcs, let $\sO_{\fG}$ be an
open-loop subfamily for $\fG$\@, let $\tdomain$ be a time-domain, let
$\nbhd{U}\subset\man{M}$ be open, and let
$X\in\sO_{\fG}(\tdomain,\nbhd{U})$\@.  We denote by $\sO_{\fG,X}$ the
open-loop subfamily defined as follows.  If $\tdomain'\subset\tdomain$ and
$\nbhd{U}'\subset\nbhd{U}$\@, then we let
\begin{equation*}
\sO_{\fG,X}(\tdomain',\nbhd{U}')=
\setdef{X'\in\sO_{\fG}(\tdomain',\nbhd{U}')}{X'=X|\tdomain'\times\nbhd{U}'}.
\end{equation*}
If $\tdomain'\not\subset\tdomain$ and/or $\nbhd{U}'\not\subset\nbhd{U}$\@,
then we take $\sO_{\fG,X}=\emptyset$\@.  Thus $\sO_{\fG,X}$ is comprised of
those vector fields from $\sO_{\fG}$ that are merely restrictions of $X$ to
smaller domains.  Just why this might be interesting we will only see when we
discuss linearisation about a reference flow in
Section~\ref{subsec:lin-traj-flow}\@.

\item Let $m\in\integernn$ and $m'\in\{0,\lip\}$\@, let
$\nu\in\{m+m',\infty,\omega\}$\@, and let $r\in\{\infty,\omega\}$\@, as
required.  In Proposition~\ref{prop:gcs-cs-di} we saw that there was a pretty
robust correspondence between $\C^\nu$-control systems and $\C^\nu$-\gcs{}s,
\emph{at the system level}\@.  As we make our way towards trajectories, as we
are now doing, this robustness breaks down a little.  To frame this, we can
define an open-loop subfamily for the \gcs\ associated to a $\C^\nu$-control
system $\Sigma=(\man{M},F,\cs{C})$ as follows.  For a time-domain $\tdomain$
and an open $\nbhd{U}\subset\man{M}$\@, we define
\begin{equation*}
\sO_\Sigma(\tdomain,\nbhd{U})=
\setdef{F^\mu|\nbhd{U}}{\mu\in\Lloc^\infty(\tdomain;\cs{C})},
\end{equation*}
recalling that $F^\mu(t,x)=F(x,\mu(t))$\@.  We clearly have
$\sO_\Sigma(\tdomain;\nbhd{U})\subset
\subscr{\sO}{$\fG_\Sigma$,cpt}(\tdomain;\nbhd{U})$ for every time-domain
$\tdomain$ and every open $\nbhd{U}\subset\man{M}$\@; this was proved in the
course of proving Proposition~\ref{prop:open-loopcsI}\@.  Of course, by
virtue of Proposition~\ref{prop:open-loopcsII}\@, we have a corresponding
construction if the control set $\cs{C}$ is a subset of a locally convex
topological vector space, if $F$ is sublinear, and if
$\mu\in\Lloc^1(\tdomain;\cs{C})$\@.  However, we do not generally expect to
have equality of these two open-loop subfamilies.  This, in turn, will have
repercussions on the nature of the trajectories for these subfamilies, and,
therefore, on the relationship of trajectories of a control system and the
corresponding \gcs.  We will consider these matters in
Section~\ref{subsec:trajequiv}\@, and we will see that, for many interesting
classes of control systems, there is, in fact, a natural trajectory
correspondence between the system and its associated \gcs.\oprocend
\end{compactenum}
\end{examples}

Our notion of an open-loop subfamily is very general, and working with the
full generality will typically lead to annoying problems.  There are many
attributes that one may wish for open-loop subfamilies to satisfy in order to
relax some the annoyance.  To illustrate, let us define a typical attribute
that one may require, that of translation-invariance.  Let us define some
notation so that we can easily make the definition.  For a time-domain
$\tdomain$ and for $s\in\real$\@, we denote
\begin{equation*}
s+\tdomain=\setdef{s+t}{t\in\tdomain}
\end{equation*}
and we denote by $\map{\tau_s}{s+\tdomain}{\tdomain}$ the translation map
$\tau_s(t)=t-s$\@.
\begin{definition}
Let $m\in\integernn$ and $m'\in\{0,\lip\}$\@, let
$\nu\in\{m+m',\infty,\omega\}$\@, and let $r\in\{\infty,\omega\}$\@, as
required.  Let $\fG=(\man{M},\sF)$ be a $\C^\nu$-\gcs. An open-loop subfamily
$\sO_{\fG}$ for $\fG$ is \defn{translation-invariant} if, for every
$s\in\real$\@, every time-domain $\tdomain$\@, and every open set
$\nbhd{U}\subset\man{M}$\@, the map
\begin{equation*}
\mapdef{(\tau_s\times\id_{\nbhd{U}})^*}
{\sO_{\fG}(s+\tdomain,\nbhd{U})}{\sO_{\fG}(\tdomain,\nbhd{U})}
{X}{X\scirc(\tau_s\times\id_{\nbhd{U}})}
\end{equation*}
is a bijection.\oprocend
\end{definition}

An immediate consequence of the definition is, of course, that if
$t\mapsto\xi(t)$ is a trajectory (we will formally define the notion of
``trajectory'' in the next section), then so is $t\mapsto\xi(s+t)$ for every
$s\in\real$\@.

Let us now think about how open-loop subfamilies interact with completion.
In order for the definition we are about to make make sense, we should verify
the following lemma.
\begin{lemma}\label{lem:germs->section}
Let\/ $m\in\integernn$\/ and $m'\in\{0,\lip\}$\@, let\/
$\nu\in\{m+m',\infty,\omega\}$\@, and let\/ $r\in\{\infty,\omega\}$\@, as
required.  Let\/ $\man{M}$ be a\/ $\C^r$-manifold, let\/
$\tdomain\subset\real$ be an interval, and let\/
$\map{X}{\tdomain\times\man{M}}{\tb{\man{M}}}$ have the property that\/
$X(t,x)\in\tb[x]{\man{M}}$ for each\/ $(t,x)\in\tdomain\times\man{M}$\@.
Then the following statements hold:
\begin{compactenum}[(i)]
\item \label{pl:locCVF->CVF} if, for each\/ $x\in\man{M}$\@, there exist a
neighbourhood\/ $\nbhd{U}$ of\/ $x$ and\/
$X'\in\CFsections[\nu]{\tdomain;\tb{\nbhd{U}}}$ such that\/
$[X_t]_x=[X'_t]_x$ for every\/ $t\in\tdomain$\@, then\/
$X\in\CFsections[\nu]{\tdomain;\tb{\man{M}}}$\@;
\item \label{pl:locLI->LI} if, for each\/ $x\in\man{M}$\@, there exist a
neighbourhood\/ $\nbhd{U}$ of\/ $x$ and\/
$X'\in\LIsections[\nu]{\tdomain;\tb{\nbhd{U}}}$ such that\/ $[X_t]_x=[X'_t]_x$
for every\/ $t\in\tdomain$\@, then\/
$X\in\LIsections[\nu]{\tdomain;\tb{\man{M}}}$\@;
\item \label{pl:locLB->LB} if, for each\/ $x\in\man{M}$\@, there exist a
neighbourhood\/ $\nbhd{U}$ of\/ $x$ and\/
$X'\in\LBsections[\nu]{\tdomain;\tb{\nbhd{U}}}$ such that\/
$[X_t]_x=[X'_t]_x$ for every\/ $t\in\tdomain$\@, then\/
$X\in\LBsections[\nu]{\tdomain;\tb{\man{M}}}$\@.
\end{compactenum}
\begin{proof}
\eqref{pl:locCVF->CVF} Let $x\in\man{M}$\@.  Since $X$ agrees in some
neighbourhood of $x$ with a Carath\'eodory vector field $X'$\@, it follows
that $t\mapsto X_t(x)=X'_t(x)$ is measurable.  In like manner, let
$t\in\tdomain$ and let $x_0\in\man{M}$\@.  Then $x\mapsto X_t(x)=X'_t(x)$ is
of class $\C^\nu$ in a neighbourhood of $x_0$\@, and so $x\mapsto X_t(x)$ is
of class $\C^\nu$\@.

\eqref{pl:locLI->LI} For $K\subset\man{M}$ be compact, for $k\in\integernn$\@,
and for $\vect{a}\in\csd(\integernn;\realp)$\@, denote
\begin{equation*}
p_K=\begin{cases}p^\infty_{K,k},&\nu=\infty,\\p^m_K,&\nu=m,\\
p^{m+\lip}_K,&\nu=m+\lip,\\p^\omega_{K,\vect{a}},&\nu=\omega.\end{cases}
\end{equation*}
Let $K\subset\man{M}$ be compact, let $x\in K$\@, let $\nbhd{U}_x$ be a
relatively compact neighbourhood of $x$\@, and let
$X_x\in\LIsections[\nu]{\tdomain;\nbhd{U}_x}$ be such that
$[X_t]_x=[X_{x,t}]_x$ for every $t\in\tdomain$\@.  Then there exists
$g_x\in\Lloc^1(\tdomain;\realnn)$ such that
\begin{equation*}
p_{\closure(\nbhd{U}_x)}(X_{x,t})\le g_x(t),\qquad t\in\tdomain.
\end{equation*}
Now let $x_1,\dots,x_k\in K$ be such that
$K\subset\cup_{j=1}^k\nbhd{U}_{x_j}$\@.  Let
$g(t)=\max\{g_{x_1}(t),\dots,g_{x_k}(t)\}$\@, noting that the associated
function $g$ is measurable by~\cite[Proposition~2.1.3]{DLC:80} and is locally
integrable by the triangle inequality, along with the fact that
\begin{equation*}
g(t)\le C(g_{x_1}(t)+\dots+g_{x_k}(t))
\end{equation*}
for some suitable $C\in\realp$ (this is simply the statement of the
equivalence of the $\ell^1$ and $\ell^\infty$ norms for $\real^n$).  We then
have
\begin{equation*}
p_K(X_t)\le g(t),\qquad t\in\tdomain,
\end{equation*}
showing that $X\in\LIsections[\nu]{\tdomain;\tb{\man{M}}}$\@.

\eqref{pl:locLB->LB} This is proved in exactly the same manner, \emph{mutatis
mutandis}\@, as the preceding part of the lemma.
\end{proof}
\end{lemma}

The following definition can now be made.
\begin{definition}
Let $m\in\integernn$ and $m'\in\{0,\lip\}$\@, let
$\nu\in\{m+m',\infty,\omega\}$\@, and let $r\in\{\infty,\omega\}$\@, as
required.  Let $\fG=(\man{M},\sF)$ be a $\C^\nu$-\gcs\ and let $\sO_{\fG}$ be
an open-loop subfamily for $\fG$\@.  The \defn{completion} of $\sO_{\fG}$ is
the open-loop subfamily $\Sh(\sO_{\fG})$ for $\Sh(\fG)$ defined by specifying
that $(X,\tdomain,\nbhd{U})\in\Sh(\sO_{\fG})$ if, for each $x\in\nbhd{U}$\@,
there exist a neighbourhood $\nbhd{U}'\subset\nbhd{U}$ of $x$ and
$(X',\tdomain,\nbhd{U}')\in\sO_{\fG}(\tdomain,\nbhd{U}')$ such that
$[X_t]_x=[X'_t]_x$ for each $t\in\tdomain$\@.\oprocend
\end{definition}

Clearly the completion of an open-loop subfamily is an open-loop subfamily
for the completion.  Moreover, if
$(X,\tdomain,\nbhd{U})\in\sO_{\fG}(\tdomain,\nbhd{U})$\@, then
$(X,\tdomain,\nbhd{U})\in\Sh(\sO_{\fG}(\tdomain,\nbhd{U}))$\@, but one cannot
expect the converse assertion to generally hold.

\subsection{Trajectories}\label{subsec:traj}

With the concept of open-loop system just developed, it is relatively easy to
provide a notion of a trajectory for a \gcs.
\begin{definition}
Let $m\in\integernn$ and $m'\in\{0,\lip\}$\@, let
$\nu\in\{m+m',\infty,\omega\}$\@, and let $r\in\{\infty,\omega\}$\@, as
required.  Let $\fG=(\man{M},\sF)$ be a $\C^\nu$-\gcs\ and let $\sO_{\fG}$ be
an open-loop subfamily for $\fG$\@.
\begin{compactenum}[(i)]
\item For a time-domain $\tdomain$\@, an open set $\nbhd{U}\subset\man{M}$\@,
and for $X\in\sO_{\fG}(\tdomain,\nbhd{U})$\@, an
\defn{$(X,\tdomain,\nbhd{U})$-trajectory} for $\sO_{\fG}$ is a curve
$\map{\xi}{\tdomain}{\nbhd{U}}$ such that $\xi'(t)=X(t,\xi(t))$\@.
\item For a time-domain $\tdomain$ and an open set
$\nbhd{U}\subset\man{M}$\@, a \defn{$(\tdomain,\nbhd{U})$-trajectory} for
$\sO_{\fG}$ is a curve $\map{\xi}{\tdomain}{\nbhd{U}}$ such that
$\xi'(t)=X(t,\xi(t))$ for some $X\in\sO_{\fG}(\tdomain,\nbhd{U})$\@.
\item A \defn{trajectory} for $\sO_{\fG}$ is a curve that is a
$(\tdomain,\nbhd{U})$-trajectory for $\sO_{\fG}$ for some time-domain
$\tdomain$ and some open set $\nbhd{U}\subset\man{M}$\@.\savenum
\end{compactenum}
We denote by:
\begin{compactenum}[(i)]\resumenum
\item $\Traj(X,\tdomain;\nbhd{U})$ the set of
$(X,\tdomain,\nbhd{U})$-trajectories for $\sO_{\fG}$\@;
\item $\Traj(\tdomain,\nbhd{U},\sO_{\fG})$ the set of
$(\tdomain,\nbhd{U})$-trajectories for $\sO_{\fG}$\@;
\item $\Traj(\sO_{\fG})$ the set of trajectories for $\sO_{\fG}$\@.
\end{compactenum}
We shall abbreviate $\Traj(\tdomain,\nbhd{U},\fG)=
\Traj(\tdomain,\nbhd{U},\subscr{\sO}{$\fG$,full})$ and
$\Traj(\fG)=\Traj(\subscr{\sO}{$\fG$,full})$\@.\oprocend
\end{definition}

Sometimes one wishes to keep track of the fact that, associated with a
trajectory is an open-loop system.  The following notion is designed to
capture this.
\begin{definition}
Let $m\in\integernn$ and $m'\in\{0,\lip\}$\@, let
$\nu\in\{m+m',\infty,\omega\}$\@, and let $r\in\{\infty,\omega\}$\@, as
required.  Let $\fG=(\man{M},\sF)$ be a $\C^\nu$-\gcs\ and let $\sO_{\fG}$ be
an open-loop subfamily for $\fG$\@.  A \defn{referenced
$\sO_{\fG}$-trajectory} is a pair $(X,\xi)$ where
$X\in\sO_{\fG}(\tdomain;\nbhd{U})$ and $\xi\in\Traj(X,\tdomain,\nbhd{U})$\@.
By $\Rtraj(\tdomain,\nbhd{U},\sO_{\fG})$ we denote the set of referenced
$\sO_{\fG}$-trajectories for which $X\in\sO_{\fG}(\tdomain;\nbhd{U})$\@.
\end{definition}

In Section~\ref{subsec:trajequiv} below, we shall explore trajectory
correspondences between \gcs{}s, control systems, and differential
inclusions.

The notion of a trajectory immediately gives rise to a certain open-loop
subfamily.  At present it may not be clear why this construction is
interesting, but it will come up in Section~\ref{subsec:lin-traj-flow} when
we talk about linearisations about trajectories.
\begin{example}\label{eg:traj-ol}
Let $m\in\integernn$ and $m'\in\{0,\lip\}$\@, let
$\nu\in\{m+m',\infty,\omega\}$\@, and let $r\in\{\infty,\omega\}$\@, as
required.  Let $\fG=(\man{M},\sF)$ be a $\C^\nu$-\gcs, let $\sO_{\fG}$ be an
open-loop subfamily for $\fG$\@, and let
$\xi\in\Traj(\tdomain,\nbhd{U},\sO_{\fG})$\@.  We denote by $\sO_{\fG,\xi}$
the open-loop subfamily defined as follows.  If $\tdomain'\subset\tdomain$
and $\nbhd{U}'\subset\nbhd{U}$ are such that
$\xi(\tdomain')\subset\nbhd{U}'$\@, then we let
\begin{equation*}
\sO_{\fG,\xi}(\tdomain',\nbhd{U}')=
\setdef{X\in\sO_{\fG}(\tdomain',\nbhd{U}')}{\xi'(t)=X(t,\xi(t)),\
\ae\ t\in\tdomain'}.
\end{equation*}
If $\tdomain'\not\subset\tdomain$ or $\nbhd{U}'\not\subset\nbhd{U}$\@, or if
$\tdomain'\subset\tdomain$ and $\nbhd{U}'\subset\nbhd{U}$ but
$\xi(\tdomain')\not\subset\nbhd{U}'$\@, then we take
$\sO_{\fG,\xi}=\emptyset$\@.  Thus $\sO_{\fG,\xi}$ is comprised of those
vector fields from $\sO_{\fG}$ possessing $\xi$ (restricted to the
appropriate subinterval) as an integral curve.\oprocend
\end{example}

In control theory, trajectories are of paramount importance, often far more
important, say, than systems \emph{per se}\@.  For this reason, one might ask
that completion of a \gcs\ preserve trajectories.  However, this will
generally not be the case, as the following counterexample illustrates.
\begin{example}
We will chat our way through a general example; the reader can very easily
create a specific concrete instance from the general discussion.

Let $m\in\integernn$ and $m'\in\{0,\lip\}$\@, let
$\nu\in\{m+m',\infty,\omega\}$\@, and let $r\in\{\infty,\omega\}$\@, as
required.  We let $\man{M}$ be a $\C^r$-manifold with Riemannian metric
$\metric$\@.  We consider the presheaf $\subscr{\sF}{bdd}$ of bounded
$\C^\nu$-vector fields on $\man{M}$\@, initially discussed in
Example~\enumdblref{eg:!sheaf}{enum:bounded-presheaf}\@.  We let
$\subscr{\fG}{bdd}=(\man{M},\subscr{\sF}{bdd})$ so that, as we saw in
Example~\enumdblref{eg:sheafify}{enum:shfybdd}\@,
$\Sh(\subscr{\sF}{bdd})=\ssections[\nu]{\tb{\man{M}}}$\@.  Let $X$ be a
vector field possessing an integral curve $\map{\xi}{\tdomain}{\man{M}}$ for
which
\begin{equation*}
\limsup_{t\to\sup\tdomain}\dnorm{\xi'(t)}_{\metric}=\infty
\end{equation*}
(this requires that $\tdomain$ be noncompact, of course).

Now let us see how this gives rise to a trajectory for
$\Sh(\subscr{\fG}{bdd})$ that is not a trajectory for $\subscr{\fG}{bdd}$\@.
We let $\tdomain$ be the interval of definition of the integral curve $\xi$
described above.  We consider the open subset $\man{M}\subset\man{M}$\@.  We
then have the open-loop system $(X,\tdomain,\man{M})$ specified by letting
$X(t)=X$ (abusing notation),~\ie~we consider a time-independent open-loop
system.  It is clear, then, that
$\xi\in\Traj(\tdomain,\man{M},\Sh(\subscr{\fG}{bdd}))$ (since
$\Sh(\subscr{\fG}{bdd})=(\man{M},\ssections[\nu]{\tb{\man{M}}})$ as we showed
in Example~\enumdblref{eg:sheafify}{enum:shfybdd}), but that $\xi$ cannot be
a trajectory for $\subscr{\fG}{bdd}$ since any vector field possessing $\xi$
as an integral curve cannot be bounded.\oprocend
\end{example}

Thus we cannot expect sheafification to generally preserve trajectories.
This should be neither a surprise nor a disappointment to us.  It is
gratifying, however, that sheafification \emph{does} preserve trajectories in
at least one important case.
\begin{proposition}\label{prop:sheafify-traj}
Let\/ $m\in\integernn$ and\/ $m'\in\{0,\lip\}$\@, let\/
$\nu\in\{m+m',\infty,\omega\}$\@, and let\/ $r\in\{\infty,\omega\}$\@, as
required.  Let\/ $\fG=(\man{M},\sF)$ be a globally generated\/ $\C^\nu$-\gcs,
let\/ $\tdomain$ be a time-domain, and let\/ $\sO_{\fG}$ be an open-loop
subfamily for\/ $\fG$\@.  For a locally absolutely continuous curve\/
$\map{\xi}{\tdomain}{\man{M}}$ the following statements are equivalent:
\begin{compactenum}[(i)]
\item $\xi\in\Traj(\tdomain,\nbhd{U},\sO_{\fG})$ for some open set\/
$\nbhd{U}\subset\man{M}$\@;
\item $\xi\in\Traj(\tdomain,\nbhd{U}',\Sh(\sO_{\fG}))$ for some open set\/ $\nbhd{U}'\subset\man{M}$\@.
\end{compactenum}
\begin{proof}
Since
$\sO_{\fG}(\tdomain,\nbhd{U})\subset\Sh(\sO_{\fG})(\tdomain,\nbhd{U})$\@, the
first assertion clearly implies the second.  So it is the opposite
implication we need to prove.

Thus let $\nbhd{U}'\subset\man{M}$ be open and suppose that
$\xi\in\Traj(\tdomain,\nbhd{U}',\Sh(\sO_{\fG}))$\@.  Let
$X\in\LIsections[\nu]{\tdomain;\tb{\nbhd{U}'}}$ be such that $\xi$ is an
integral curve for $X$ and such that $X_t\in\Sh(\sF)(\nbhd{U}')$ for every
$t\in\tdomain$\@.  For each \emph{fixed} $\tau\in\tdomain$\@, there exists
$X_\tau\in\LIsections[\nu]{\tdomain;\sF(\man{M})}$ such that
$[X_{\tau,t}]_{\xi(\tau)}=[X_t]_{\xi(\tau)}$ for every $t\in\tdomain$\@.
(This is the definition of $\Sh(\sO_{\fG})$\@, noting that $\sF$ is globally
generated.)  This means that around $\tau$ we have a bounded open interval
$\tdomain_\tau\subset\tdomain$ and a neighbourhood $\nbhd{U}_\tau$ of
$\xi(\tau)$ so that $\xi(\tdomain_\tau)\subset\nbhd{U}_\tau$ and so that
$\xi'(t)=X_\tau(t,\xi(t))$ for almost every $t\in\tdomain_\tau$\@.  By
paracompactness, we can choose a locally finite refinement of these intervals
that also covers $\tdomain$\@.  By repartitioning, we arrive at a locally
finite pairwise disjoint covering $\ifam{\tdomain_j}_{j\in J}$ of $\tdomain$
by subintervals with the following property: the index set $J$ is a finite or
countable subset of $\integer$ chosen so that $t_1<t_2$ whenever
$t_1\in\tdomain_{j_1}$ and $t_2\in\tdomain_{j_2}$ with $j_1<j_2$\@.  That is,
we order the labels for the elements of the partition in the natural way,
this making sense since the cover is locally finite.  By construction, we
have $X_j\in\LIsections[\nu]{\tdomain_j;\sF(\man{M})}$ with the property that
$\xi|\tdomain_j$ is an integral curve for $X_j$\@.  We then define
$\map{\ol{X}}{\tdomain}{\sF(\man{M})}$ by asking that
$\ol{X}|\tdomain_j=X_j$\@.  It remains to show that
$\ol{X}\in\LIsections[\nu]{\tdomain;\sF(\man{M})}$\@.

Because each of the vector fields $X_j$\@, $j\in J$\@, is a Carath\'eodory
vector field, we easily conclude that $\ol{X}$ is also a Carath\'eodory
vector field.

Let $K\subset\man{M}$ be compact, $k\in\integernn$\@,
and $\vect{a}\in\csd(\integernn;\realp)$\@, and denote
\begin{equation*}
p_K=\begin{cases}p^\infty_{K,k},&\nu=\infty,\\p^m_K,&\nu=m,\\
p^{m+\lip}_K,&\nu=m+\lip,\\p^\omega_{K,\vect{a}},&\nu=\omega.\end{cases}
\end{equation*}
For each $j\in J$\@, there then exists $g_j\in\Lloc^1(\tdomain_j;\realnn)$
such that
\begin{equation*}
p_K(X_{j,t})\le g_j(t),\qquad t\in\tdomain_j.
\end{equation*}
Define $\map{g}{\tdomain}{\realnn}$ by asking that $g|\tdomain_j=g_j$\@.  We
claim that $g\in\Lloc^1(\tdomain;\realnn)$\@.  Let
$\mathbb{T}'\subset\tdomain$ be a compact subinterval.  The set
\begin{equation*}
J_{\mathbb{T}'}=\setdef{j\in J}{\mathbb{T}'\cap\tdomain_j\not=\emptyset}.
\end{equation*}
is finite by local finiteness of the cover $\ifam{\tdomain_j}_{j\in J}$\@.
Now we have
\begin{equation*}
\int_{\mathbb{T}'}g(t)\,\d{t}\le
\sum_{j\in J_{\mathbb{T}'}}\int_{\tdomain_j}g_j(t)\,\d{t}<\infty.
\end{equation*}
Since
\begin{equation*}
p_K(\ol{X}_t)\le g(t),\qquad t\in\tdomain,
\end{equation*}
we conclude that $\ol{X}\in\LIsections[\nu]{\tdomain;\tb{\man{M}}}$\@, as
desired.
\end{proof}
\end{proposition}

\subsection{Attributes that can be given to
\gcs{}s}\label{subsec:gcs-attributes}

In this section we show that some typical assumptions that are made for
control systems also can be made for \gcs{}s.  None of this is particularly
earth-shattering, but it does serves as a plausibility check for our
framework, letting us know that it has some common ground with familiar
constructions from control theory.

A construction that often occurs in control theory is to determine a
trajectory as the limit of a sequence of trajectories in some manner.  To
ensure the existence of such limits, the following property for \gcs{}s is
useful.
\begin{definition}
Let $m\in\integernn$ and $m'\in\{0,\lip\}$\@, let
$\nu\in\{m+m',\infty,\omega\}$\@, and let $r\in\{\infty,\omega\}$\@, as
required.  A $\C^\nu$-\gcs\ $\fG=(\man{M},\sF)$ is \defn{closed} if
$\sF(\nbhd{U})$ is closed in the topology of $\sections[\nu]{\tb{\nbhd{U}}}$
for every open set $\nbhd{U}\subset\man{M}$\@.\oprocend
\end{definition}

Here are some examples of control systems that give rise to closed \gcs{}s.
\begin{proposition}\label{prop:closed-gcs}
Let\/ $m\in\integernn$ and\/ $m'\in\{0,\lip\}$\@, let\/
$\nu\in\{m+m',\infty,\omega\}$\@, and let\/ $r\in\{\infty,\omega\}$\@, as
required.  Let\/ $\Sigma=(\man{M},F,\cs{C})$ be a\/ $\C^\nu$-control system
with\/ $\fG_\Sigma$ the associated $\C^\nu$-\gcs\ as in
Example~\enumdblref{eg:gcs}{enum:cs->gcs}\@.  Then\/ $\fG_\Sigma$ is closed if\/
$\Sigma$ has either of the following two attributes:
\begin{compactenum}[(i)]
\item \label{pl:closedgcs1} $\cs{C}$ is compact;
\item \label{pl:closedgcs2} $\cs{C}$ is a closed subset of\/ $\real^k$ and
the system is control-affine,~\ie
\begin{equation*}
F(x,\vect{u})=f_0(x)+\sum_{a=1}^ku^af_a(x),
\end{equation*}
for\/ $f_0,f_1,\dots,f_k\in\sections[\nu]{\tb{\man{M}}}$\@.
\end{compactenum}
\begin{proof}
\eqref{pl:closedgcs1} Let $\nbhd{U}\subset\man{M}$ be open.  By
Propositions~\ref{prop:paramCinfty}\@,~\ref{prop:paramCmm'}\@,
and~\ref{prop:paramComega}\@, the map
\begin{equation*}
\cs{C}\ni u\mapsto F^u\in\sections[\nu]{\tb{\nbhd{U}}}
\end{equation*}
is continuous.  Now let $\nbhd{U}\subset\man{M}$ be open and note that
$\sF_\Sigma(\nbhd{U})$ is the image of $\cs{C}$ under the mapping
\begin{equation*}
\cs{C}\ni u\mapsto F^u|\nbhd{U}\in\sections[\nu]{\tb{\nbhd{U}}}.
\end{equation*}
Thus $\sF_\Sigma(\nbhd{U})$ is compact, and so closed, being the image of a
compact set under a continuous mapping~\cite[Theorem~17.7]{SW:04}\@.

\eqref{pl:closedgcs2} Let $\nbhd{U}\subset\man{M}$ be open.  Just as in the
preceding part of the proof, we consider the mapping $\vect{u}\mapsto
F^{\vect{u}}|\nbhd{U}$\@.  Note that the image of the mapping
\begin{equation*}
\vect{u}\mapsto F^{\vect{u}}=f_0+\sum_{a=1}^ku^af_a
\end{equation*}
is a finite-dimensional affine subspace of the $\real$-vector space
$\sections[\nu]{\tb{\nbhd{U}}}$\@.  Therefore, this image is closed
since~(1)~locally convex topologies are translation invariant (by
construction) and since~(2)~finite-dimensional subspaces of locally convex
spaces are closed~\cite[Proposition~2.10.1]{JH:66}\@.  Moreover, the map
$\vect{u}\mapsto F^{\vect{u}}|\nbhd{U}$ is closed onto its image since any
surjective linear map between finite-dimensional locally convex space is
closed.  We conclude, therefore, that if we restrict this map from all of
$\real^k$ to $\cs{C}$\@, then the image is closed.
\end{proof}
\end{proposition}

Let us next turn to attributes of \gcs{}s arising from the fact, shown in
Example~\enumdblref{eg:gcs}{enum:gcs->di}\@, that \gcs{}s give rise to differential
inclusions in a natural way.
\begin{proposition}
Let\/ $m\in\integernn$ and\/ $m'\in\{0,\lip\}$\@, let\/
$\nu\in\{m+m',\infty,\omega\}$\@, and let\/ $r\in\{\infty,\omega\}$\@, as
required.  If\/ $\fG=(\man{M},\sF)$ is a\/ $\C^\nu$-\gcs, then
\begin{compactenum}[(i)]
\item \label{pl:lsc} $\sX_{\fG}$ is lower semicontinuous and
\item \label{pl:usc} $\sX_{\fG}$ is upper semicontinuous if\/ $\fG$ is
globally generated and\/ $\sF(\man{M})$ is compact.
\end{compactenum}
\begin{proof}
\eqref{pl:lsc} Let $x_0\in\man{M}$ and let $v_{x_0}\in\sX_{\fG}(x_0)$\@.
Then there exist a neighbourhood $\nbhd{W}$ of $x_0$ and $X\in\sF(\nbhd{W})$
such that $X(x_0)=v_{x_0}$\@.  Let $\nbhd{V}\subset\tb{\man{M}}$ be a
neighbourhood of $v_{x_0}$\@.  By continuity of $X$\@, there exists a
neighbourhood $\nbhd{U}\subset\nbhd{W}$ of $x_0$ such that
$X(\nbhd{U})\subset\nbhd{V}$\@.  This implies that $X(x)\in\sX_{\fG}(x)$ for
every $x\in\nbhd{U}$\@, giving lower semicontinuity of $\sX_{\fG}$\@.

\eqref{pl:usc} Let $x_0\in\man{M}$ and let $\nbhd{V}\subset\tb{\man{M}}$ be a
neighbourhood of $\sX_{\fG}(x_0)$\@.  For each $X\in\sF(\man{M})$\@,
$\nbhd{V}$ is a neighbourhood of $X(x_0)$ and so there exist neighbourhoods
$\nbhd{M}_X\subset\man{M}$ of $x_0$ and $\nbhd{C}_X\subset\sF(\man{M})$ of
$X$ such that
\begin{equation*}
\setdef{X'(x)}{x\in\nbhd{M}_X,\ X'\in\nbhd{C}_X}\subset\nbhd{V}.
\end{equation*}
Since $\sF(\man{M})$ is compact, let $X_1,\dots,X_k\in\sF(\man{M})$ be such
that $\sF(\man{M})=\cup_{j=1}^k\cs{C}_{X_j}$\@.  Then the neighbourhood
$\nbhd{U}=\cap_{j=1}^k\nbhd{M}_{X_j}$ of $x_0$ has the property that
$\sX_{\fG}(\nbhd{U})\subset\nbhd{V}$\@.
\end{proof}
\end{proposition}

There are many easy examples to illustrate that compactness of $\sF(\man{M})$
is generally required in part~\eqref{pl:usc} of the preceding result.  Here is
one.
\begin{example}
Let $m\in\integernn$ and $m'\in\{0,\lip\}$\@, let
$\nu\in\{m+m',\infty,\omega\}$\@, and let $r\in\{\infty,\omega\}$\@, as
required.  Let $\man{M}$ be a $\C^r$-manifold and let $x_0\in\man{M}$\@.  Let
$\sF(x_0)$ be the globally generated sheaf of sets of $\C^\nu$-vector fields
defined by
\begin{equation*}
\sF(x_0)(\man{M})=
\setdef{X\in\sections[\nu]{\tb{\man{M}}}}{X(x_0)=0}.
\end{equation*}
We claim that, if we take $\fG=(\man{M},\sF(x_0))$\@, then we have
\begin{equation}\label{eq:F(x0)}
\sX_{\fG}(x)=\begin{cases}\{0_{x_0}\},&x=x_0,\\
\tb[x]{\man{M}},&x\not=x_0.\end{cases}
\end{equation}
In the case $\nu=\infty$ or $\nu=m$\@, this is straightforward.  Let
$\nbhd{U}$ be a neighbourhood of $x\not=x_0$ such that
$x_0\not\in\closure(\nbhd{U})$\@.  By the smooth Tietze Extension
Theorem~\cite[Proposition~5.5.8]{RA/JEM/TSR:88}\@, if
$X\in\sections[\infty]{\tb{\man{M}}}$\@, then there exists
$X'\in\sections[\infty]{\tb{\man{M}}}$ such that $X'|\nbhd{U}=X|\nbhd{U}$ and
such that $X'(x_0)=0_{x_0}$\@.  Thus $[X]_x=[X'_x]$ and so we have
$\sF(x_0)_x=\gsections[\nu]{x}{\man{M}}$ in this case.  From
this,~\eqref{eq:F(x0)} follows.

The case of $\nu=m+\lip$ follows as does the case $\nu=m$\@, noting that a
locally Lipschitz vector field multiplied by a smooth function is still a
locally Lipschitz vector field~\cite[Proposition~1.5.3]{NW:99}\@.

The case of $\nu=\omega$ is a little more difficult, and relies on Cartan's
Theorem~A for coherent sheaves on real analytic manifolds~\cite{HC:57}\@.
Here is the argument for those who know a little about sheaves.  First,
define a sheaf of sets (in fact, submodules) of real analytic vector fields
by
\begin{equation*}
\sI_{x_0}(\nbhd{U})=\begin{cases}
\setdef{X\in\sections[\omega]{\tb{\nbhd{U}}}}{X(x_0)=0_{x_0}},&
x_0\in\nbhd{U},\\
\sections[\omega]{\tb{\nbhd{U}}},&x_0\not\in\nbhd{U}.\end{cases}
\end{equation*}
We note that $\sI_{x_0}$ is a coherent sheaf since it is a finitely generated
subsheaf of the coherent sheaf
$\ssections[\omega]{\tb{\man{M}}}$~\cite[Theorem~3.16]{JPD:12}\@.%
\footnote{This relies on the fact that Oka's Theorem, in the version of ``the
sheaf of sections of a vector bundle is coherent,'' holds in the real
analytic case.  It does, and the proof is the same as for the holomorphic
case~\cite[Theorem~3.19]{JPD:12} since the essential ingredient is the
Weierstrass Preparation Theorem, which holds in the real analytic
case~\cite[Theorem 6.1.3]{SGK/HRP:02}\@.}  Let $x\not=x_0$ and let
$v_x\in\tb[x]{\man{M}}$\@.  By Cartan's Theorem~A, there exist
$X_1,\dots,X_k\in\sI_{x_0}(\man{M})=\sF(x_0)(\man{M})$ such that
$[X_1]_x,\dots,[X_k]_x$ generate
$(\sI_{x_0})_x=\gsections[\omega]{x}{\tb{\man{M}}}$ as a module over the ring
$\gfunc[\omega]{x}{\man{M}}$ of germs of functions at $x$\@.  Let
$[X]_x\in\gsections[\omega]{x}{\tb{\man{M}}}$ be such that $X(x)=v_x$\@.
There then exist $[f^1]_x,\dots,[f^k]_x\in\gfunc[\omega]{x}{\man{M}}$ such
that
\begin{equation*}
[f^1]_x[X_1]_x+\dots+[f^k]_x[X_k]_x=[X]_x.
\end{equation*}
Therefore,
\begin{equation*}
v_x=X(x)=f^1(x)X_1(x)+\dots+f^k(x)X_k(x),
\end{equation*}
and so, taking
\begin{equation*}
X=f^1X_1+\dots+f^kX_k\in\sI_{x_0}(\man{M})=\sF(x_0)(\man{M}),
\end{equation*}
we see that $v_x=X(x)\in\sX_{\fG}(x)$\@, which establishes~\eqref{eq:F(x0)}
in this case.

In any event,~\eqref{eq:F(x0)} holds, and it is easy to see that this
differential inclusion is not upper semicontinuous.\oprocend
\end{example}

We can make the following definitions, rather analogous to those of
Definition~\ref{def:dihulls} for differential inclusions.
\begin{definition}
Let $m\in\integernn$ and $m'\in\{0,\lip\}$\@, let
$\nu\in\{m+m',\infty,\omega\}$\@, and let $r\in\{\infty,\omega\}$\@, as
required.  The $\C^\nu$-\gcs\ $\fG=(\man{M},\sF)$ is:
\begin{compactenum}[(i)]
\item \defn{closed-valued} (\resp~\defn{compact-valued}\@,
\defn{convex-valued}) at $x\in\man{M}$ if $\sX_{\fG}(x)$ is closed (\resp,
compact, convex);
\item \defn{closed-valued} (\resp~\defn{compact-valued}\@,
\defn{convex-valued}) if $\sX_{\fG}(x)$ is closed (\resp, compact, convex)
for every $x\in\man{M}$\@.\oprocend
\end{compactenum}
\end{definition}

One can now talk about taking ``hulls'' under various properties.  Let us
discuss this for the properties of closedness and convexity.  First we need
the definitions we will use.
\begin{definition}
Let $m\in\integernn$ and $m'\in\{0,\lip\}$\@, let
$\nu\in\{m+m',\infty,\omega\}$\@, and let $r\in\{\infty,\omega\}$\@, as
required.  Let $\fG=(\man{M},\sF)$ be a $\C^\nu$-\gcs.
\begin{compactenum}[(i)]
\item The \defn{convex hull} of $\fG$ is the $\C^\nu$-\gcs\
$\cohull(\fG)=(\man{M},\cohull(\sF))$\@, where $\cohull(\sF)$ is the presheaf
of subsets of $\C^\nu$-vector fields given by
\begin{equation*}
\cohull(\sF)(\nbhd{U})=\cohull(\sF(\nbhd{U})),
\end{equation*}
the convex hull on the right being that in the $\real$-vector space
$\sections[\nu]{\tb{\nbhd{U}}}$\@.
\item The \defn{closure} of $\fG$ is the $\C^\nu$-\gcs
\begin{equation*}
\closure(\fG)=(\man{M},\closure(\sF)),
\end{equation*}
where $\closure(\sF)$ is the presheaf of subsets of $\C^\nu$-vector fields
given by $\closure(\sF)(\nbhd{U})=\closure(\sF(\nbhd{U}))$\@, the closure on
the right being that in the $\real$-topological vector space
$\sections[\nu]{\tb{\nbhd{U}}}$\@.  The reader should verify that
$\closure(\sF)$ is indeed a presheaf.\oprocend
\end{compactenum}
\end{definition}

Let us now relate the two different sorts of ``hulls'' we have.
\begin{proposition}
Let\/ $m\in\integernn$ and\/ $m'\in\{0,\lip\}$\@, let\/
$\nu\in\{m+m',\infty,\omega\}$\@, and let\/ $r\in\{\infty,\omega\}$\@, as
required.  Let\/ $\fG=(\man{M},\sF)$ be a\/ $\C^\nu$-\gcs\ with\/ $\sX_{\fG}$
the associated differential inclusion.  Then the following statements hold:
\begin{compactenum}[(i)]
\item \label{pl:gcshull1} $\cohull(\sX_{\fG})=\sX_{\cohull(\fG)}$\@;
\item \label{pl:gcshull2} $\sX_{\closure(\fG)}\subset\closure(\sX_{\fG})$
and\/ $\sX_{\closure(\fG)}=\closure(\sX_{\fG})$ if\/ $\fG$ is globally
generated and\/ $\sF(\man{M})$ is bounded in the compact bornology (or,
equivalently, the von Neumann bornology if $\nu\in\{\infty,\omega\}$).
\end{compactenum}
\begin{proof}
\eqref{pl:gcshull1} Let $x\in\man{M}$\@.  If $v\in\cohull(\sX_{\fG}(x))$\@,
then there exist $v_1,\dots,v_k\in\sX_{\fG}(x)$ and
$c_1,\dots,c_k\in\interval[0,1]$ satisfying $\sum_{j=1}^kc_j=1$ such that
\begin{equation*}
v=c_1v_1+\dots+c_kv_k.
\end{equation*}
Let $\nbhd{U}_1,\dots,\nbhd{U}_k$ be neighbourhoods of $x$ and let
$X_j\in\sF(\nbhd{U}_j)$ be such that $X_j(x)=v_j$\@, $j\in\{1,\dots,k\}$\@.
Then, taking $\nbhd{U}=\cap_{j=1}^k\nbhd{U}_j$\@,
\begin{equation*}
c_1X_1|\nbhd{U}+\dots+c_kX_k|\nbhd{U}\in\cohull(\sF(\nbhd{U})),
\end{equation*}
showing that $\cohull(\sX_{\fG}(x))\subset\sX_{\cohull(\fG)}(x)$\@.

Conversely, let $v\in\sX_{\cohull(\fG)}$\@, let $\nbhd{U}$ be a neighbourhood
of $x$\@, and let $X\in\cohull(\sF(\nbhd{U}))$ be such that $X(x)=v$\@.  Then
\begin{equation*}
X=c_1X_1+\dots+c_kX_k
\end{equation*}
for $X_1,\dots,X_k\in\sF(\nbhd{U})$ and for $c_1,\dots,c_k\in\interval[0,1]$
satisfying $\sum_{j=1}^kc_j=1$\@.  We then have
\begin{equation*}
v=c_1X_1(x)+\dots+c_kX_k(x)\in\cohull(\sX_{\fG})(x),
\end{equation*}
completing the proof of the proposition as concerns convex hulls.

\eqref{pl:gcshull2} Let $x\in\man{M}$\@, let $v\in\sX_{\closure(\fG)}(x)$\@,
let $\nbhd{U}$ be a neighbourhood of $x$\@, and let
$X\in\closure(\sF(\nbhd{U}))$ be such that $X(x)=v$\@.  Let $(I,\preceq)$ be
a directed set and let $\ifam{X_i}_{i\in I}$ be a net in $\sF(\nbhd{U})$
converging to $X$ in the appropriate topology.  Then we have $\lim_{i\in
I}X_i(x)=X(x)$ since the net $\ifam{X_i}_{i\in I}$ converges uniformly in
some neighbourhood of $x$ (this is true for all cases of $\nu$). Thus
$v\in\closure(\sX_{\fG}(x))$\@, as desired.

Suppose that $\sF$ is globally generated with $\sF(\man{M})$ bounded, let
$x\in\man{M}$\@, and let $v\in\closure(\sX_{\fG})(x)$\@.  Thus there exists a
sequence $\ifam{v_j}_{j\in\integerp}$ in $\sX_{\fG}(x)$ converging to $v$\@.
Let $X_j\in\sF(\man{M})$ be such that $X_j(x)=v_j$\@, $j\in\integerp$\@.
Since $\closure(\sF(\man{M}))$ is compact, there is a subsequence
$\ifam{X_{j_k}}_{j_k}$ in $\sF(\man{M})$ converging to
$X\in\closure(\sF(\man{M}))$\@.  Moreover,
\begin{equation*}
X(x)=\lim_{k\to\infty}X_{j_k}(x)=\lim_{j\to\infty}v_j=v
\end{equation*}
since $\ifam{X_{j_k}}_{k\in\integerp}$ converges to $X$ uniformly in some
neighbourhood of $x$ (again, this is true for all $\nu$).  Thus
$v\in\sX_{\closure(\fG)}(x)$\@.

The parenthetical comment in the final assertion of the proof follows since
the compact and von Neumann bornologies agree for nuclear spaces~\cite[Proposition~4.47]{AP:69}\@.
\end{proof}
\end{proposition}

The following example shows that the opposite inclusion stated in the
proposition for closures does not generally hold.
\begin{example}
We will talk our way through a general sort of example, leaving to the reader
the job of instantiating this to give a concrete example.

Let $m\in\integerp$ and $m'\in\{0,\lip\}$\@, let
$\nu\in\{m+m',\infty,\omega\}$\@, and let $r\in\{\infty,\omega\}$\@, as
required.  Let $\man{M}$ be a $\C^r$-manifold.  Let $x\in\man{M}$ and let
$\ifam{X_j}_{j\in\integerp}$ be a sequence of $\C^\nu$-vector fields with the
following properties:
\begin{compactenum}
\item $\ifam{X_j(x)}_{j\in\integerp}$ converges to $0_x$\@;
\item $X_j(x)\not=0_x$ for all $j\in\integerp$\@;
\item there exists a neighbourhood $\nbhd{O}$ of zero in
$\sections[\nu]{\tb{\man{M}}}$ such that, for each $j\in\integerp$\@,
\begin{equation*}
\setdef{k\in\integerp\setminus\{j\}}{X_k-X_j\in\nbhd{O}}=\emptyset.
\end{equation*}
\end{compactenum}
Let $\sF$ be the globally generated presheaf of sets of $\C^\nu$-vector
fields given by $\sF(\man{M})=\setdef{X_j}{j\in\integerp}$\@.  Then
$0_x\in\closure(\sX_{\fG}(x))$\@.  We claim that
$0_x\not\in\sX_{\closure(\fG)}(x)$\@.  To see this, suppose that
$0_x\in\sX_{\closure(\fG)}(x)$\@.  Since $\sF(\man{M})$ is countable, this
implies that there is a subsequence $\ifam{X_{j_k}}_{k\in\integerp}$ that
converges in $\sections[\nu]{\tb{\man{M}}}$\@.  But this is prohibited by the
construction of the sequence $\ifam{X_j}_{j\in\integerp}$\@.\oprocend
\end{example}

\subsection{Trajectory correspondence between \gcs{}s and other sorts of
control systems}\label{subsec:trajequiv}

In Example~\ref{eg:gcs} and Proposition~\ref{prop:gcs-cs-di} we made precise
the connections between various models for control systems: control systems,
differential inclusions, and \gcs{}s.  In order to flesh out these
connections more deeply, in this section we investigate the possible
correspondences between the trajectories for the various models.

We first consider correspondences between trajectories of control systems and
their associated \gcs{}s.  Thus we let $m\in\integernn$ and
$m'\in\{0,\lip\}$\@, let $\nu\in\{m+m',\infty,\omega\}$\@, and let
$r\in\{\infty,\omega\}$\@, as required.  Let $\Sigma=(\man{M},F,\cs{C})$ be a
$\C^\nu$-control system with $\fG_\Sigma$ the associated $\C^\nu$-\gcs, as in
Example~\enumdblref{eg:gcs}{enum:cs->gcs}\@.  As we saw in
Proposition~\pldblref{prop:gcs-cs-di}{pl:cs->gcs->cs}\@, the correspondence
between $\Sigma$ and $\fG_\Sigma$ is perfect, at the system level, when the
map $u\mapsto F^u$ is injective and open onto its image.
Part~\eqref{pl:trajequiv12} of the following result shows that this perfect
correspondence almost carries over at the level of trajectories as well.
Included with this statement we include a few other related ideas concerning
trajectory correspondences.
\begin{theorem}\label{the:cs->gcs-equivI}
Let\/ $m\in\integernn$ and\/ $m'\in\{0,\lip\}$\@, let\/
$\nu\in\{m+m',\infty,\omega\}$\@, and let\/ $r\in\{\infty,\omega\}$\@, as
required.  Let\/ $\Sigma=(\man{M},F,\cs{C})$ be a\/ $\C^\nu$-control system
with\/ $\fG_\Sigma$ the associated\/ $\C^\nu$-\gcs, as in
Example~\enumdblref{eg:gcs}{enum:cs->gcs}\@.  Then the following statements
hold:
\begin{compactenum}[(i)]
\item \label{pl:trajequiv11} $\Traj(\tdomain,\nbhd{U},\Sigma)\subset
\Traj(\tdomain,\nbhd{U},\subscr{\sO}{$\fG_\Sigma$,cpt})$\@;
\item \label{pl:trajequiv12} if the map\/ $u\mapsto F^u$ is injective and
proper, then\/
$\Traj(\tdomain,\nbhd{U},\subscr{\sO}{$\fG_\Sigma$,cpt})\subset
\Traj(\tdomain,\nbhd{U},\Sigma)$\@;
\item \label{pl:trajequiv13} if\/ $\cs{C}$ is a Suslin topological
space\footnote{Recall that this means that $\cs{C}$ is the continuous image
of a complete, separable, metric space.} and if\/ $F$ is proper, then\/
$\Traj(\tdomain,\nbhd{U},\sO_{\fG_\Sigma,\infty})\subset
\Traj(\tdomain,\nbhd{U},\Sigma)$\@.
\item \label{pl:trajequiv14} if, in addition,\/ $\nu\in\{\infty,\omega\}$\@,
then we may replace\/
$\Traj(\tdomain,\nbhd{U},\subscr{\sO}{$\fG_\Sigma$,cpt})$ with
$\Traj(\tdomain,\nbhd{U},\sO_{\fG_\Sigma,\infty})$ in
statements~\eqref{pl:trajequiv11} and~\eqref{pl:trajequiv12}\@.
\end{compactenum}
\begin{proof}
\eqref{pl:trajequiv11} Let $\xi\in\Traj(\tdomain,\nbhd{U},\Sigma)$ and let
$\mu\in\Lloc^\infty(\tdomain;\cs{C})$ be such that
\begin{equation*}
\xi'(t)=F(\xi(t),\mu(t)),\qquad\ae\ t\in\tdomain.
\end{equation*}
Note that, as we saw in Example~\ref{eg:cs->openloop}\@,
$F^\mu|\nbhd{U}\in\sO_{\fG_\Sigma,\infty}(\tdomain,\nbhd{U})$\@, making sure
to note that the conclusions of Proposition~\ref{prop:open-loopcsI} imply
that $F^\mu\in\LBsections[\nu]{\tdomain;\tb{\man{M}}}$\@.  Thus
$\xi\in\Traj(\tdomain,\nbhd{U},\sO_{\fG_\Sigma,\infty})$\@.  To show that, in
fact, $\xi\in\Traj(\tdomain,\nbhd{U},\subscr{\sO}{$\fG_\Sigma$,cpt})$\@, let
$\tdomain'\subset\tdomain$ be a compact subinterval and let $K\subset\cs{C}$
be a compact set such that $\mu(t)\in K$ for almost every $t\in\tdomain'$\@.
Denote
\begin{equation*}
\mapdef{\hat{F}}{\cs{C}}{\sections[\nu]{\tb{\man{M}}}}{u}{F^u.}
\end{equation*}
Since $\hat{F}$ is continuous, $\hat{F}(K)$ is
compact~\cite[Theorem~17.7]{SW:04}\@.  Since $F^\mu_t\in\hat{F}(K)$ for
almost every $t\in\tdomain'$\@, we conclude that
$\xi\in\Traj(\tdomain,\nbhd{U},\subscr{\sO}{$\fG_\Sigma$,cpt})$\@, as
claimed.

\eqref{pl:trajequiv12} Recall from~\cite[Proposition~I.10.2]{NB:89b} that, if
$\hat{F}$ (as defined above) is proper, then it has a closed image, and is a
homeomorphism onto its image.  If
$\xi\in\Traj(\tdomain,\nbhd{U},\subscr{\sO}{$\fG_\Sigma$,cpt})$\@, then there
exists $X\in\subscr{\sO}{$\fG_\Sigma$,cpt}$ such that $\xi'(t)=X(t,\xi(t))$
for almost every $t\in\tdomain$\@.  Note that, since
$X\in\subscr{\sO}{$\fG_\Sigma$,cpt}$\@, we have
$X(t)\in\sF_\Sigma(\man{M})=\image(\hat{F})$\@.  Thus, by hypothesis, there
exists a unique $\map{\mu}{\tdomain}{\cs{C}}$ such that
$\hat{F}\scirc\mu=X$\@.  To show that $\mu$ is measurable, let
$\nbhd{O}\subset\cs{C}$ be open so that $\hat{F}(\nbhd{O})$ is an open subset
of $\image(\hat{F})$\@.  Thus there exists an open set
$\nbhd{O}'\subset\sections[\nu]{\tb{\man{M}}}$ such that
$\hat{F}(\nbhd{O})=\image(\hat{F})\cap\nbhd{O}'$\@.  Then we have
\begin{equation*}
\mu^{-1}(\nbhd{O})=X^{-1}(\hat{F}(\nbhd{O}))=X^{-1}(\nbhd{O}'),
\end{equation*}
giving the desired measurability.  To show that
$\mu\in\Lloc^\infty(\tdomain;\cs{C})$\@, let $\tdomain'\subset\tdomain$ be a
compact subinterval and let $K\subset\sections[\nu]{\tb{\man{M}}}$ be such
that $X(t)\in K$ for almost every $t\in\tdomain'$\@.  Then, since $\hat{F}$
is proper, $\hat{F}\null^{-1}(K)$ is a compact subset of $\cs{C}$\@.  Since
$\mu(t)\in\hat{F}\null^{-1}(K)$ for almost every $t\in\tdomain'$ we conclude
that $\mu\in\Lloc^\infty(\tdomain;\cs{C})$\@.

\eqref{pl:trajequiv13} Let
$\xi\in\Traj(\tdomain,\nbhd{U},\sO_{\fG_\Sigma,\infty})$ and let
$X\in\sO_{\fG_\Sigma,\infty}(\tdomain,\nbhd{U})$ be such that
$\xi'(t)=X(t,\xi(t))$ for almost every $t\in\tdomain$\@.  We wish to
construct $\mu\in\Lloc^\infty(\tdomain,\cs{C})$ such that
\begin{equation*}
\xi'=F(\xi(t),\mu(t)),\qquad\ae\ t\in\tdomain.
\end{equation*}
We fix an arbitrary element $\bar{u}\in\cs{C}$ (it matters not which) and
then define a set-valued map $\setmap{U}{\tdomain}{\cs{C}}$ by
\begin{equation*}
U(t)=\begin{cases}\setdef{u\in\cs{C}}{\xi'(t)=F(\xi(t),u)},&
\xi'(t)\ \textrm{exists},\\\{\bar{u}\},&\textrm{otherwise}.\end{cases}
\end{equation*}
Since $X(t)\in\sF_\Sigma(\man{M})$\@, we conclude that
$X(t)\in\image(\hat{F})$ for every $t\in\tdomain$\@,~\ie~$X(t)=F^u$ for some
$u\in\cs{C}$\@, and so $U(t)\not=\emptyset$ for every $t\in\tdomain$\@.

Properness of $F$ ensures that $U(t)$ is compact for every $t\in\tdomain$\@.
The following lemma shows that \emph{any} selection $\mu$ of $U$ is locally
essentially bounded in the compact bornology.
\begin{prooflemma}\label{plem:muT'}
If\/ $\tdomain'\subset\tdomain$ is a compact subinterval, then the set\/
$\cup\setdef{U(t)}{t\in\tdomain'}$ is contained in a compact subset of\/
$\cs{C}$\@.
\begin{subproof}
Let us define $\map{F_\xi}{\tdomain\times\cs{C}}{\tb{\man{M}}}$ by
$F_\xi(t,u)=F(\xi(t),u)$\@.  We claim that, if $\tdomain'\subset\tdomain$ is
compact, then $F_\xi|\tdomain'\times\cs{C}$ is proper.  To see this, first
define
\begin{equation*}
\mapdef{G_\xi}{\tdomain'\times\cs{C}}{\man{M}\times\cs{C}}
{(t,u)}{(\xi(t),u),}
\end{equation*}
\ie~$G_\xi=\xi\times\id_{\cs{C}}$\@.  With this notation, we have
$F_\xi=F\scirc G_\xi$\@.  Since $F_\xi^{-1}(K)=G_\xi^{-1}(F^{-1}(K))$ and
since $F$ is proper, to show that $F_\xi$ is proper it suffices to show that
$G_\xi$ is proper.  Let $K\subset\man{M}\times\cs{C}$ be compact.  We let
$\map{\pr_1}{\man{M}\times\cs{C}}{\man{M}}$ and
$\map{\pr_2}{\man{M}\times\cs{C}}{\cs{C}}$ be the projections.  Note that
\begin{equation*}
G_\xi^{-1}(K)=(\xi\times\id_{\cs{C}})^{-1}(K)\subset
\xi^{-1}(\pr_1(K))\times\id_{\cs{C}}^{-1}(\pr_2(K)).
\end{equation*}
Since the projections are continuous, $\pr_1(K)$ and $\pr_2(K)$ are
compact~\cite[Theorem~17.7]{SW:04}\@.  Since $\xi$ is a continuous function
whose domain (for our present purposes) is the compact set $\tdomain'$\@,
$\xi^{-1}(\pr_1(K))$ is compact.  Since the identity map is proper,
$\id_{\cs{C}}^{-1}(\pr_2(K))$ is compact.  Thus $G_\xi^{-1}(K)$ is contained
in a product of compact sets.  Since a product of compact sets is
compact~\cite[Theorem~17.8]{SW:04} and $G_\xi^{-1}(K)$ is closed by
continuity of $G_\xi$\@, it follows that $G_\xi^{-1}(K)$ is compact, as
claimed.  Thus $F_\xi|\tdomain'\times\cs{C}$ is proper.

Now, since $\xi$ is a trajectory for the $\sO_{\fG_\Sigma,\infty}$ open-loop
subfamily, there exists a compact set $K'\subset\tb{\man{M}}$ such that
\begin{equation*}
\setdef{\xi'(t)}{t\in\tdomain'}\subset K',
\end{equation*}
adopting the convention that $\xi'(t)$ is taken to satisfy
$\xi'(t)=F(\xi(t),\bar{u})$ when $\xi'(t)$ does not exist; this is an
arbitrary and inconsequential choice.  By our argument above,
$K''\eqdef(F_\xi|\tdomain'\times\cs{C})^{-1}(K')$ is compact.  Therefore, for
each $t\in\tdomain'$\@,
\begin{align*}
\setdef{(t,u)\in\tdomain'\times\cs{C}}
{u\in U(t)}=&\;\setdef{(t,u)\in\tdomain'\times\cs{C}}{F(\xi(t),u)=\xi'(t)}\\
\subset&\;\setdef{(t,u)\in\tdomain'\times\cs{C}}{F(\xi(t),u)\in K'}
\subset K''.
\end{align*}
Defining the compact set (compact by~\cite[Theorem~17.7]{SW:04})
$K=\pr_2(K'')$\@, with $\map{\pr_2}{\tdomain'\times\cs{C}}{\cs{C}}$ the
projection, we then have
\begin{equation*}\eqsubqed
\cup\setdef{U(t)}{t\in\tdomain'}\subset K.
\end{equation*}
\end{subproof}
\end{prooflemma}

We shall now make a series of observations about the set-valued map $U$\@,
using results of \citet{CJH:75} on measurable multi-valued mappings,
particularly with values in Suslin spaces.
\begin{prooflemma}
The set-valued map\/ $U$ is measurable,~\ie~if\/
$\nbhd{O}\subset\cs{C}$ is open, then
\begin{equation*}
U^{-1}(\nbhd{O})=\setdef{t\in\tdomain}{U(t)\cap\nbhd{O}\not=\emptyset}
\end{equation*}
is measurable.
\begin{subproof}
Define
\begin{equation*}
\mapdef{F_\xi}{\tdomain\times\cs{C}}{\tb{\man{M}}}
{(t,u)}{F(\xi(t),u),}
\end{equation*}
noting that $t\mapsto F_\xi(t,u)$ is measurable for each $u\in\cs{C}$ and
that $u\mapsto F_\xi(t,u)$ is continuous for every $t\in\tdomain$\@.  It
follows from~\cite[Theorem~6.4]{CJH:75} that $U$ is measurable as stated.
\end{subproof}
\end{prooflemma}

\begin{prooflemma}
There exists a measurable function\/ $\map{\mu}{\tdomain}{\cs{C}}$ such
that\/ $\mu(t)\in U(t)$ for almost every\/ $t\in\tdomain$\@.
\begin{subproof}
First note that $U(t)$ is a closed subset of $\cs{C}$ since it is either the
singleton $\{\bar{u}\}$ or the preimage of the closed set $\{\xi'(t)\}$ under
the continuous map $u\mapsto F(\xi(t),u)$\@.  It follows
from~\cite[Theorem~3.5]{CJH:75} that
\begin{equation*}
\graph(U)=\setdef{(t,u)\in\tdomain\times\cs{C}}{u\in U(t)}
\end{equation*}
is measurable with respect to the product $\sigma$-algebra of the Lebesgue
measurable sets in $\tdomain$ and the Borel sets in $\cs{C}$\@.  The lemma
now follows from~\cite[Theorem~5.7]{CJH:75}\@.
\end{subproof}
\end{prooflemma}

Now, for $t\in\tdomain$ having the property that $\xi'(t)$ exists and that
$\mu(t)\in U(t)$ (with $\mu$ from the preceding lemma), we have
$\xi'(t)=F(\xi(t),\mu(t))$\@.

\eqref{pl:trajequiv14} This follows by our observation of Example~\enumdblref{eg:olsf}{enum:cptsubfamily}\@.
\end{proof}
\end{theorem}

Let us make some comments on the hypotheses of the preceding theorem.
\begin{remarks}\label{rem:cstraj->gcstrajI}
\begin{compactenum}
\item Part~\eqref{pl:trajequiv12} of the result has assumptions that the map
$u\mapsto F^u$ be injective and proper.  An investigation of the proof shows
that injectivity and openness onto the image of this map is enough to give
trajectories for $\Sigma$ that correspond to measurable controls.  The
additional assumption of properness, which gives the further consequence of
the image of the map $u\mapsto F^u$ being closed, allows us to conclude
boundedness of the controls.  Let us look at these assumptions.
\begin{compactenum}[(a)]
\item By the map $u\mapsto F^u$ being injective, we definitely do not mean
that the map $u\mapsto F(x,u)$ is injective for each $x\in\man{M}$\@; this is
a very strong assumption whose adoption eliminates a large number of
interesting control systems.  For example, if we take $\man{M}=\real$\@,
$\cs{C}=\real$\@, and $F(x,u)=ux\pderiv{}{x}$ to define a $\C^\nu$-control
system for any $\nu\in\{m+m',\infty,\omega\}$ with $m\in\integernn$ and
$m'\in\{0,\lip\}$\@, then the map $u\mapsto F^u$ is injective, but the map
$u\mapsto F(0,u)$ is not.

\item Let us take $\man{M}=\real^2$\@, $\cs{C}=\real$\@, and
\begin{equation*}
F((x_1,x_2),u)=f_1(u)\pderiv{}{x_1}+f_2(u)\pderiv{}{x_2},
\end{equation*}
where $\map{f_1,f_2}{\real}{\real}$ are such that the map
$u\mapsto(f_1(u),f_2(u))$ is injective and continuous, but not a
homeomorphism onto its image.  Such a system may be verified to be a
$\C^\nu$-control system for any $\nu\in\{m+m',\infty,\omega\}$ with
$m\in\integernn$ and $m'\in\{0,\lip\}$ (using
Propositions~\ref{prop:paramCinfty}\@,%
~\ref{prop:paramCmm'}\@, and~\ref{prop:paramComega}).  In this case, we claim
that the map $\hat{F}\colon u\mapsto F^u$ is injective and continuous, but
not a homeomorphism onto its image.  Injectivity of the map is clear and
continuity follows since $F$ is a jointly parameterised vector field of class
$\C^\nu$\@.  Define a linear map
\begin{equation*}
\mapdef{\kappa}{\real^2}{\sections[\nu]{\tb{\man{M}}}}
{(v_1,v_2)}{v_1\pderiv{}{x_1}+v_2\pderiv{}{x_2},}
\end{equation*}
\ie~$\kappa(\vect{v})$ is the constant vector field with components
$(v_1,v_2)$\@.  Using the seminorms for our locally convex topologies the
standard seminorm characterisations of continuous linear maps (as
in~\cite[\S{}III.1.1]{HHS/MPW:99}), we can easily see that $\kappa$ is a
continuous linear map, and so is a homeomorphism onto its closed image
(arguing as in the proof of
Proposition~\pldblref{prop:closed-gcs}{pl:closedgcs2}).  Then
$\hat{F}=\kappa\scirc(f_1\times f_2)$\@, and so we conclude that $\hat{F}$ is
a homeomorphism onto its image if and only if $f_1\times f_2$ is a
homeomorphism onto its image, and this gives our claim.

\item Let us take $\man{M}=\real$\@, $\cs{C}=\real$\@, and
$F(x,u)=\tan^{-1}(u)\pderiv{}{x}$\@.  As with the examples above, we regard
this as a control system of class $\C^\nu$ for any
$\nu\in\{m+m',\infty,\omega\}$\@, for $m\in\integernn$ and
$m'\in\{0,\lip\}$\@.  We claim that $\hat{F}\colon u\mapsto F^u$ is a
homeomorphism onto its image, but is not proper.  This is verified in exactly
the same manner as in the preceding example.

\item If $\cs{C}$ is compact, then $\hat{F}$ is proper because, if
$K\subset\sections[\nu]{\tb{\man{M}}}$ is compact, then
$\hat{F}\null^{-1}(K)$ is closed, and so
compact~\cite[Theorem~17.5]{SW:04}\@.  This gives trajectory correspondence
between a $\C^\nu$-control system and its corresponding \gcs\ for compact
control sets when the map $\hat{F}$ is injective.
\end{compactenum}

\item \label{enum:Suslin-proper} Part~\eqref{pl:trajequiv13} of the result
has two assumptions, that $\cs{C}$ is a Suslin space and that $F$ is proper.
Let us consider some cases where these hypotheses hold.
\begin{compactenum}[(a)]
\item Complete separable metric spaces are Suslin spaces.
\item If $\cs{C}$ is an open or a closed subspace of Suslin space, it is a
Suslin space~\cite[Lemma~6.6.5(ii)]{VIB:07b}\@.
\item For $m\in\integernn$\@, $m'\in\{0,\lip\}$\@, and
$\nu\in\{m+m',\infty,\omega\}$\@, $\sections[\nu]{\tb{\man{M}}}$ is a Suslin
space.  In all except the case of $\nu=\omega$\@, this follows since
$\sections[\nu]{\tb{\man{M}}}$ is a separable, complete, metrisable space.
However, $\sections[\omega]{\tb{\man{M}}}$ is not metrisable.  Nonetheless,
it is Suslin, as argued in Section~\ref{subsec:Comega-props}\@.
\item \label{enum:Fcpt} If $\cs{C}$ is compact, then $F$ is proper.  Indeed,
if $K\subset\tb{\man{M}}$ is compact, then $\tbproj{\man{M}}(K)$ is compact,
and
\begin{equation*}
F^{-1}(K)\subset\tbproj{\man{M}}(K)\times\cs{C},
\end{equation*}
and so the set on the left is compact, being a closed subset of a compact
set~\cite[Theorem~17.5]{SW:04}\@.\oprocend
\end{compactenum}
\end{compactenum}
\end{remarks}

We also have a version of the preceding theorem in the case that the control
set $\cs{C}$ is a subset of a locally convex topological vector
space,~\cf~Proposition~\ref{prop:open-loopcsII}\@.  Here we also specialise
for one of the implications to control-linear systems introduced in
Example~\ref{eg:control-linear}\@.
\begin{theorem}\label{the:cs->gcs-equivII}
Let\/ $m\in\integernn$ and\/ $m'\in\{0,\lip\}$\@, let\/
$\nu\in\{m+m',\infty,\omega\}$\@, and let\/ $r\in\{\infty,\omega\}$\@, as
required.  Let\/ $\Sigma=(\man{M},F,\cs{C})$ be a\/ $\C^\nu$-sublinear
control system for which\/ $\cs{C}$ is a subset of a locally convex
topological vector space\/ $\alg{V}$\@, and let\/ $\fG_\Sigma$ be the
associated\/ $\C^\nu$-\gcs, as in
Example~\enumdblref{eg:gcs}{enum:cs->gcs}\@.  If\/ $\tdomain$ is a
time-domain and if\/ $\nbhd{U}$ is open, then\/
$\Traj(\tdomain,\nbhd{U},\Sigma)\subset
\Traj(\tdomain,\nbhd{U},\subscr{\sO}{$\fG_\Sigma$,full})$\@.

Conversely, if
\begin{compactenum}[(i)]
\item $\Sigma$ is a\/ $\C^\nu$-control-linear system,~\ie~there exists\/
$\Lambda\in\lin{\alg{V}}{\sections[\nu]{\tb{\man{M}}}}$ such that\/
$F(x,u)=\Lambda(u)(x)$\@,
\item $\Lambda$ is injective, and
\item $\Lambda$ is an open mapping onto its image,
\end{compactenum}
then it is also the case that\/
$\Traj(\tdomain,\nbhd{U},\subscr{\sO}{$\fG_\Sigma$,full})\subset
\Traj(\tdomain,\nbhd{U},\Sigma)$\@.
\begin{proof}
We first show that $\Traj(\tdomain,\nbhd{U},\Sigma)\subset
\Traj(\tdomain,\nbhd{U},\subscr{\sO}{$\fG_\Sigma$,full})$\@.  Suppose that
$\xi\in \Traj(\tdomain,\nbhd{U},\Sigma)$\@.  Thus there exists
$\mu\in\Lloc^1(\tdomain;\cs{C})$ such that
\begin{equation*}
\xi'(t)=F(\xi(t),\mu(t)),\qquad \ae\ t\in\tdomain.
\end{equation*}
By Proposition~\ref{prop:open-loopcsII} and Example~\ref{eg:cs->openloop}\@,
$F^\mu|\nbhd{U}\in\subscr{\sO}{$\fG_\Sigma$,full}(\tdomain,\nbhd{U})$ and so
$\xi\in\Traj(\tdomain,\nbhd{U},\subscr{\sO}{$\fG_\Sigma$,full})$\@.

Now let us prove the ``conversely'' assertion of the theorem.  Thus we let
$\xi\in\Traj(\tdomain,\nbhd{U},\subscr{\sO}{$\fG_\Sigma$,full})$ so that
there exists $X\in\subscr{\sO}{$\fG_\Sigma$,full}(\tdomain,\nbhd{U})$ for
which $\xi'(t)=X(t,\xi(t))$ for almost every $t\in\tdomain$\@.  Since
$\Lambda$ is injective and since $X_t\in\Lambda(\cs{C})$ for each
$t\in\tdomain$ (this is the definition of $\fG_\Sigma$), we uniquely define
$\mu(t)\in\cs{C}$ by $\Lambda(\mu(t))=X_t$\@.  We need only show that $\mu$
is locally Bochner integrable.  Let $\Lambda^{-1}$ denote the inverse of
$\Lambda$\@, thought of as a map from $\image(\Lambda)$ to $\alg{V}$\@.  As
$\Lambda$ is open, $\Lambda^{-1}$ is continuous.  From this, measurability of
$\mu$ follows immediately.  To show that $\mu$ is locally Bochner integrable,
let $q$ be a continuous seminorm for the locally convex topology of $\alg{V}$
and, as per~\cite[\S{}III.1.1]{HHS/MPW:99}\@, let $p$ be a continuous
seminorm for the locally convex topology of $\sections[\nu]{\tb{\man{M}}}$
such that $q(\Lambda^{-1}(Y))\le p(Y)$ for every
$Y\in\sections[\nu]{\tb{\man{M}}}$\@.  Then we have, for any compact
subinterval $\tdomain'\subset\tdomain$\@,
\begin{equation*}
\int_{\tdomain'}q(\mu(t))\,\d{t}\le\int_{\tdomain'}p(X_t)\,\d{t}<\infty,
\end{equation*}
giving Bochner integrability of $\mu$ by~\cite[Theorems~3.2 and~3.3]{RB/AD:11}\@.
\end{proof}
\end{theorem}

Let us make some observations about the preceding theorem.
\begin{remarks}\label{rem:cstraj->gcstrajII}
The converse part of Theorem~\ref{the:cs->gcs-equivII} has three
hypotheses:~that the system is control-linear;~that the map from controls to
vector fields is injective;~that the map from controls to vector fields is
open onto its image.  The first hypothesis, linearity of the system, cannot
be weakened except in sort of artificial ways.  As can be seen from the
proof, linearity allows us to talk about the integrability of the associated
control.  Injectivity can be assumed without loss of generality by
quotienting out the kernel if it is not.  Let us consider some cases where
the third hypothesis holds.  Let $m\in\integernn$ and $m'\in\{0,\lip\}$\@,
let $\nu\in\{m+m',\infty,\omega\}$\@, and let $r\in\{\infty,\omega\}$\@, as
required.
\begin{compactenum}
\item Let $\cs{C}\subset\real^k$ and suppose that our system is
$\C^\nu$-control-affine,~\ie
\begin{equation*}
F(x,\vect{u})=f_0(x)+\sum_{a=1}^ku^af_a(x)
\end{equation*}
for $\C^\nu$-vector fields $f_0,f_1,\dots,f_m$\@.  As we pointed out in
Example~\ref{eg:control-linear}\@, this can be regarded as a control-linear
system by taking $\alg{V}=\real\oplus\real^k$
\begin{equation*}
\cs{C}'=\setdef{(u^0,\vect{u})\in\alg{V}}{u^0=1,\ \vect{u}\in\cs{C}},
\end{equation*}
and
\begin{equation*}
\Lambda(u^0,\vect{u})=\sum_{a=0}^ku^af_a.
\end{equation*}
We can assume $\Lambda$ is injective, as mentioned above.  In this case, the
map $\Lambda$ is a homeomorphism onto its image since any map from a
finite-dimensional locally convex space is
continuous~\cite[Proposition~2.10.2]{JH:66}\@.  Thus
Theorem~\ref{the:cs->gcs-equivII} applies to control-affine systems, and
gives trajectory equivalence in this case.

\item \label{enum:Fidentity} The other case of interest to us is that when
$\alg{V}=\sections[\nu]{\tb{\man{M}}}$ and when $\cs{C}\subset\alg{V}$ is
then a family of globally defined vector fields of class $\C^\nu$ on
$\man{M}$\@.  In this case, $\Lambda$ is the identity map on
$\sections[\nu]{\tb{\man{M}}}$\@, so the hypotheses of
Theorem~\ref{the:cs->gcs-equivII} are easily satisfied.  The trajectory
equivalence one gets in this case is that between a globally generated \gcs\
and its corresponding control system as in
Example~\enumdblref{eg:gcs}{enum:gcs->cs}\@.\oprocend
\end{compactenum}
\end{remarks}

One of the conclusions enunciated above is sufficiently interesting to
justify its own theorem.
\begin{theorem}\label{the:gcs->cs-equiv}
Let\/ $m\in\integernn$ and\/ $m'\in\{0,\lip\}$\@, let\/
$\nu\in\{m+m',\infty,\omega\}$\@, and let\/ $r\in\{\infty,\omega\}$\@, as
required.  Let\/ $\fG=(\man{M},\sF)$ be a globally generated\/ $\C^\nu$-\gcs.
As in Example~\enumdblref{eg:gcs}{enum:gcs->cs}\@, let\/
$\Sigma_{\fG}=(\man{M},\Sigma_{\fG},\cs{C}_{\sF})$ be the corresponding
$\C^\nu$-control system.  Then, for each time-domain\/ $\tdomain$ and each
open set\/ $\nbhd{U}\subset\man{M}$\@,\/
$\Traj(\tdomain,\nbhd{U},\subscr{\sO}{$\fG$,full})=
\Traj(\tdomain,\nbhd{U},\Sigma_{\fG})$\@.
\begin{proof}
This is the observation made in
Remark~\enumdblref{rem:cstraj->gcstrajII}{enum:Fidentity}\@.
\end{proof}
\end{theorem}

Now we turn to relationships between trajectories for \gcs{}s and
differential inclusions.  In Example~\enumdblref{eg:gcs}{enum:di->gcs} we
showed how a \gcs\ can be built from a differential inclusion.  However, as
we mentioned in that example, we cannot expect any sort of general
correspondence between trajectories of the differential inclusion and the
\gcs\ constructed from it; differential inclusions are just too irregular.
We can, however, consider the correspondence in the other direction, as the
following theorem indicates.
\begin{theorem}
Let\/ $m\in\integernn$ and\/ $m'\in\{0,\lip\}$\@, let\/
$\nu\in\{m+m',\infty,\omega\}$\@, and let\/ $r\in\{\infty,\omega\}$\@, as
required.  Let\/ $\fG=(\man{M},\sF)$ be a\/ $\C^\nu$-\gcs\ and let\/
$\sX_{\fG}$ be the associated differential inclusion, as in
Example~\enumdblref{eg:gcs}{enum:gcs->di}\@.  For\/ $\tdomain$ a time-domain
and\/ $\nbhd{U}\subset\man{M}$ an open set,\/
$\Traj(\tdomain,\nbhd{U},\fG)\subset\Traj(\tdomain,\nbhd{U},\sX_{\fG})$\@.

Conversely, if\/ $\sF$ is globally generated and if\/ $\sF(\man{M})$ is a
compact subset of\/ $\sections[\nu]{\tb{\man{M}}}$\@, then
$\Traj(\tdomain,\nbhd{U},\sX_{\fG})\subset\Traj(\tdomain,\nbhd{U},\fG)$\@.
\begin{proof}
Since, for an open-loop system $(X,\tdomain,\nbhd{U})$\@,
$X(t)\in\sF(\nbhd{U})$ for every $t\in\tdomain$\@, we have
$X(t,x)\in\sX_{\fG}(x)$ for every $(t,x)\in\tdomain\times\nbhd{U}$\@.  Thus,
if $\xi\in\Traj(\tdomain,\nbhd{U},\fG)$\@, then we have
$\xi'(t)\in\sX_{\fG}(\xi(t))$ for almost every $t\in\tdomain$\@.

For the ``conversely'' part of the theorem, if $\xi$ is a trajectory for the
differential inclusion $\sX_{\fG}$ then, for almost every $t\in\tdomain$\@,
$\xi'(t)=X(\xi(t))$ for some $X\in\sF(\man{M})$\@.  Therefore, let us fix an
arbitrary $\ol{X}\in\sF(\man{M})$ and let us define
$\setmap{U}{\tdomain}{\sF(\man{M})}$ by
\begin{equation*}
U(t)=\begin{cases}\setdef{X\in\sF(\man{M})}{\xi'(t)=X(\xi(t))},&
\xi'(t)\ \textrm{exists},\\\{\ol{X}\},&\textrm{otherwise}.\end{cases}
\end{equation*}
Now we note that
\begin{compactenum}
\item $\cs{C}_{\sF}=\sF(\man{M})$ is a Suslin space, being a closed subset of
a Suslin space, and
\item the map $F_{\sF}$ is proper by
Remark~\enumdblref{rem:cstraj->gcstrajI}{enum:Fcpt}\@.
\end{compactenum}
Thus we are in exactly the right framework to use the proof of
Theorem~\pldblref{the:cs->gcs-equivI}{pl:trajequiv13} to show that there
exists a locally essentially bounded measurable control $t\mapsto X(t)$ for
which
\begin{equation*}
\xi'(t)=F_{\sF}(\xi(t),X(t)),\qquad\ae\ t\in\tdomain,
\end{equation*}
and so $\xi\in\Traj(\tdomain,\nbhd{U},\Sigma_{\fG})$\@, as desired.
\end{proof}
\end{theorem}

Let us comment on the hypotheses of this theorem.
\begin{remark}
The assumption that $\sF(\man{M})$ be compact in the ``conversely'' part of
the preceding theorem is indispensable.  The connection going from
differential inclusion to \gcs\ is too ``loose'' to get any sort of useful
trajectory correspondence, without restricting the class of vector fields
giving rise to the differential inclusion.  Roughly speaking, this is because
a differential inclusion only prescribes the values of vector fields, and the
topologies have to do with derivatives as well.\oprocend
\end{remark}

\subsection{The category of \gcs{}s}\label{subsec:gcs-category}

In our discussion of feedback equivalence in
Section~\ref{subsubsec:parameterised} we indicated that the notion of
equivalence in our framework is not interesting to us.  In this section, we
illustrate why it not interesting by defining a natural notion of
equivalence, and then seeing that it degenerates to something trivial under
natural hypotheses.  We do this in a general way by considering first how one
might define a ``category'' of \gcs{}s with objects and morphisms.  The
problem of equivalence is then the problem of understanding isomorphisms in
this category.  By imposing a naturality condition on morphisms via
trajectories, we prove that isomorphisms are uniquely determined by
diffeomorphisms of the underlying manifolds for the two \gcs{}s.  The notion
of ``direct image'' we use here is common in sheaf theory, and we refer
to~\cite[\eg][Definition~2.3.1]{MK/PS:90} for some discussion.  However, by
far the best presentation that we could find of direct images of presheaves
such as we use here is in the online documentation~\citepalias{Stacks}\@.

Let us first describe how to build maps between \gcs{}s.  This is done first
by making the following definition.
\begin{definition}
Let $m\in\integernn$ and $m'\in\{0,\lip\}$\@, let
$\nu\in\{m+m',\infty,\omega\}$\@, and let $r\in\{\infty,\omega\}$\@, as
required.  Let $\fG=(\man{M},\sF)$ be a $\C^\nu$-\gcs, let $\man{N}$ be
$\C^r$-manifold, and let $\Phi\in\mappings[r]{\man{M}}{\man{N}}$\@.  The
\defn{direct image} of $\fG$ by $\Phi$ is the \gcs\
$\Phi_*\fG=(\man{N},\Phi_*\sF)$ defined by
$\Phi_*\sF(\nbhd{V})=\sF(\Phi^{-1}(\nbhd{V}))$ for $\nbhd{V}\subset\man{N}$
open.\oprocend
\end{definition}

One easily verifies that if $\sF$ is a sheaf, then so too is $\Phi_*\sF$\@.

With the preceding sheaf construction, we can define what we mean by a
morphism of \gcs{}s.
\begin{definition}\label{def:tcs-morphism}
Let $m\in\integernn$ and $m'\in\{0,\lip\}$\@, let
$\nu\in\{m+m',\infty,\omega\}$\@, and let $r\in\{\infty,\omega\}$\@, as
required.  Let $\fG=(\man{M},\sF)$ and $\fH=(\man{N},\sG)$ be
$\C^\nu$-\gcs{}s.  A \defn{morphism} from $\fG$ to $\fH$ is a pair
$(\Phi,\Phi^\sharp)$ such that
\begin{compactenum}[(i)]
\item $\Phi\in\mappings[r]{\man{M}}{\man{N}}$ and
\item $\Phi^\sharp=\ifam{\Phi^\sharp_{\nbhd{V}}}_{\nbhd{V}\,\textrm{open}}$
is a family of mappings
$\map{\Phi^\sharp_{\nbhd{V}}}{\sG(\nbhd{V})}{\Phi_*\sF(\nbhd{V})}$\@,
$\nbhd{V}\subset\man{N}$ defined as follows:
\begin{compactenum}[(a)]
\item there exists a family
$L_{\nbhd{V}}\in\lin{\sections[\nu]{\tb{\nbhd{V}}}}
{\sections[\nu]{\tb{(\Phi^{-1}(\nbhd{V}))}}}$ of continuous linear mappings
satisfying\/ $L_{\nbhd{V}'}=L_{\nbhd{V}}|\sections[\nu]{\tb{\nbhd{V}'}}$ if\/
$\nbhd{V},\nbhd{V}'\subset\man{N}$ are open with\/
$\nbhd{V}'\subset\nbhd{V}$\@;
\item $\Phi^\sharp_{\nbhd{V}}=L_{\nbhd{V}}|\sG(\nbhd{V})$\@.\oprocend
\end{compactenum}
\end{compactenum}
\end{definition}

\begin{definition}
Let $m\in\integernn$ and $m'\in\{0,\lip\}$\@, let
$\nu\in\{m+m',\infty,\omega\}$\@, and let $r\in\{\infty,\omega\}$\@, as
required.  Let $\fG=(\man{M},\sF)$ and $\fH=(\man{N},\sG)$ be
$\C^\nu$-\gcs{}s.  A \defn{morphism} from $\fG$ to $\fH$ is a pair
$(\Phi,\Phi^\sharp)$ where
\begin{compactenum}[(i)]
\item $\Phi\in\mappings[r]{\man{M}}{\man{N}}$ and
\item $\Phi^\sharp=\ifam{\Phi^\sharp_{\nbhd{V}}}_{\nbhd{V}\,\textrm{open}}$
is a family of mappings
$\map{\Phi^\sharp_{\nbhd{V}}}{\sG(\nbhd{V})}{\Phi_*\sF(\nbhd{V})}$\@,
$\nbhd{V}\subset\man{N}$ open, satisfying
\begin{compactenum}[(a)]
\item $\Phi^\sharp_{\nbhd{V}}$ is the restriction to $\sG(\nbhd{V})$ of
$L_{\nbhd{V}}\in\lin{\sections[\nu]{\tb{\nbhd{V}}}}
{\sections[\nu]{\tb{(\Phi^{-1}(\nbhd{V}))}}}$ and
\item
$\Phi^\sharp_{\nbhd{V}'}(Y|\nbhd{V}')=(\Phi^\sharp_{\nbhd{V}}(Y))|\nbhd{V}'$\@,
for $Y\in\sG(\nbhd{V})$ and for open sets $\nbhd{V},\nbhd{V}'\subset\man{N}$
satisfying $\nbhd{V}'\subset\nbhd{V}$\@.\oprocend
\end{compactenum}
\end{compactenum}
\end{definition}

By the preceding definition, we arrive at the ``category of
$\C^\nu$-\gcs{}s'' whose objects are \gcs{}s and whose morphisms are as just
defined.  From the point of view of control theory, one wishes to restrict
these definitions further to account for the fact that morphisms ought to
preserve trajectories.  Therefore, let us see how trajectories come into the
picture.  First we consider open-loop systems.  Thus let $\tdomain$ be a
time-domain and let $\nbhd{V}\subset\man{N}$ be open.  If
$\map{Y}{\tdomain}{\sG(\nbhd{U})}$ then we have
$\Phi^\sharp(Y)_t\eqdef\Phi^\sharp_{\nbhd{V}}(Y_t)\in\sF(\Phi^{-1}(\nbhd{V}))$
for each $t\in\tdomain$\@.  That is, an open-loop system
$(Y,\tdomain,\nbhd{V})$ for $\fH$ gives rise to an open-loop system
$(\Phi^\sharp(Y),\tdomain,\Phi^{-1}(\nbhd{V}))$ for $\fG$\@.  For such a
correspondence to have significance, it must do the more or less obvious
thing to trajectories.
\begin{definition}
Let $m\in\integernn$ and $m'\in\{0,\lip\}$\@, let
$\nu\in\{m+m',\infty,\omega\}$\@, and let $r\in\{\infty,\omega\}$\@, as
required.  Let $\fG=(\man{M},\sF)$ and $\fH=(\man{N},\sG)$ be
$\C^\nu$-\gcs{}s.  A morphism $(\Phi,\Phi^\sharp)$ from $\fG$ to $\fH$ is
\defn{natural} if, for each time-domain $\tdomain$\@, each open
$\nbhd{V}\subset\man{N}$\@, and each
$Y\in\LIsections[\nu]{\tdomain;\sG(\nbhd{V})}$\@, any integral curve
$\map{\xi}{\tdomain'}{\Phi^{-1}(\nbhd{V})}$ for the time-varying vector field
$t\mapsto\Phi^\sharp(Y_t)$ defined on\/ $\tdomain'\subset\tdomain$ has the
property that $\Phi\scirc\xi$ is an integral curve for\/ $Y$\@.\oprocend
\end{definition}

Note that the time-varying vector field $t\mapsto\Phi^\sharp(Y_t)$ from the
definition is locally integrally bounded by~\cite[Lemma~1.2]{RB/AD:11}\@.

We can now characterise these natural morphisms.
\begin{proposition}\label{prop:natural-morphism}
Let\/ $m\in\integernn$ and\/ $m'\in\{0,\lip\}$\@, let\/
$\nu\in\{m+m',\infty,\omega\}$\@, and let\/ $r\in\{\infty,\omega\}$\@, as
required.  Let\/ $\fG=(\man{M},\sF)$ and\/ $\fH=(\man{M},\sG)$ be\/
$\C^\nu$-\gcs{}s.  A morphism\/ $(\Phi,\Phi^\sharp)$ from\/ $\fG$ to\/ $\fH$
is natural if and only if, for each open\/ $\nbhd{V}\subset\man{N}$\@, each\/
$Y\in\sG(\nbhd{V})$\@, each\/ $y\in\nbhd{V}$\@, and each\/
$x\in\Phi^{-1}(y)$\@, we have\/ $\tf[x]{\Phi}(\Phi^\sharp(Y)(x))=Y(y)$\@.
\begin{proof}
First suppose that $(\Phi,\Phi^\sharp)$ is natural, and let
$\nbhd{V}\subset\man{N}$ be open, let $Y\in\sG(\nbhd{V})$\@, let
$y\in\nbhd{V}$\@, and let $x\in\Phi^{-1}(\nbhd{V})$\@.  Let
$\tdomain\subset\real$ be a time-domain for which $0\in\interior(\tdomain)$
and for which the integral curve $\eta$ for $Y$ through $y$ is defined on
$\tdomain$\@.  We consider $Y\in\LIsections[\nu]{\tdomain;\sG(\nbhd{V})}$ by
taking $Y_t=Y$\@,~\ie~$Y$ is a time-independent time-varying vector field.
Note that integral curves of $Y$ can, therefore, be chosen to be
differentiable~\cite[Theorem~1.3]{EEC/NL:84}\@, and \emph{will be}
differentiable if $\nu>0$\@.  Let $\tdomain'\subset\tdomain$ be such that the
differentiable integral curve $\xi$ for $\Phi^\sharp(Y)$ through $x$ is
defined on $\tdomain'$\@.  Since $(\Phi,\Phi^\sharp)$ is natural, we have
$\eta=\Phi\scirc\xi$ on $\tdomain'$\@.  Therefore,
\begin{equation*}
Y(y)=\eta'(0)=\tf[x]{\Phi}(\xi'(0))=\tf[x]{\Phi}(\Phi^\sharp(Y)(x)).
\end{equation*}

Next suppose that, for each open\/ $\nbhd{V}\subset\man{N}$\@, each\/
$Y\in\sG(\nbhd{V})$\@, each\/ $y\in\nbhd{V}$\@, and each\/
$x\in\Phi^{-1}(y)$\@, we have\/ $\tf[x]{\Phi}(\Phi^\sharp(Y)(x))=Y(y)$\@.
Let $\tdomain$ be a time-domain, let $\nbhd{V}\subset\man{N}$ be open, let
$Y\in\LIsections[\nu]{\tdomain;\sG(\nbhd{V})}$\@, and let
$\map{\xi}{\tdomain'}{\Phi^{-1}(\nbhd{V})}$ be an integral curve for the
time-varying vector field $t\mapsto\Phi^\sharp(Y_t)$ defined on
$\tdomain'\subset\tdomain$\@.  Let $\eta=\Phi\scirc\xi$\@.  Then we have
\begin{equation*}
\eta'(t)=\tf[\xi(t)]{\Phi}(\Phi^\sharp(Y_t)(\xi(t)))=Y_t(\eta(t))
\end{equation*}
for almost every $t\in\tdomain'$\@, showing that $\eta$ is an integral curve
for $Y$\@.
\end{proof}
\end{proposition}

Note that the condition $\tf[x]{\Phi}(\Phi^\sharp(Y)(x))=Y(y)$ is consistent
with the regularity conditions for $X$ and $Y$\@.  In the cases
$\nu\in\{m,\infty,\omega\}$\@, this is a consequence of the Chain Rule
(see~\cite[Proposition~2.2.8]{SGK/HRP:02} for the real analytic case).  In
the Lipschitz case this is a consequence of \cite[Example~1.4(c)]{MG:99}
combined with \cite[Proposition~1.2.2]{NW:99}\@.

To make a connection with more common notions of mappings between control
systems, let us do the following.  Let $m\in\integernn$ and
$m'\in\{0,\lip\}$\@, let $\nu\in\{m+m',\infty,\omega\}$\@, and let
$r\in\{\infty,\omega\}$\@, as required.  Suppose that we have two
$\C^\nu$-control systems $\Sigma_1=(\man{M}_1,F_1,\cs{C}_1)$ and
$\Sigma_2=(\man{M}_2,F_2,\cs{C}_2)$\@.  As \gcs{}s, these are globally
generated, so let us not fuss with general open sets for the purpose of this
illustrative discussion.  We then suppose that we have a mapping
$\Phi\in\mappings[r]{\man{M}_1}{\man{M}_2}$ and a mapping
$\map{\kappa}{\man{M}_1\times\cs{C}_2}{\cs{C}_1}$\@, which gives rise to a
correspondence between the system vector fields by
\begin{equation*}
\Phi^\sharp(F_2^{u_2})(x_1)=F_1^{\kappa(x_1,u_2)}(x_1).
\end{equation*}
The condition of naturality means that a trajectory $\xi_1$ for $\Sigma_1$
satisfying
\begin{equation*}
\xi'_1(t)=F_1(\xi_1(t),\kappa(\xi_1(t),\mu_2(t)))
\end{equation*}
gives rise to a trajectory $\xi_2=\Phi\scirc\xi_1$ for $\Sigma_2$\@, implying
that
\begin{equation*}
\xi'_2=\tf[\xi_1(t)]{\Phi}(\xi'_1(t))=\tf[\xi_1(t)]{\Phi}\scirc
F_1(\xi_1(t),\kappa(\xi_1(t),\mu_2(t))).
\end{equation*}
Thus
\begin{equation*}
F_2(x_2,u_2)=\tf[x_1]{\Phi}\scirc F_1(x_1,\kappa(x_1,u_2))
\end{equation*}
for every $x_1\in\Phi^{-1}(x_2)$\@.

There may well be some interest in studying general natural morphisms, but we
will not pursue this right at the moment.  Instead, let us simply think about
isomorphisms in the category of \gcs{}s.
\begin{definition}
Let $m\in\integernn$ and $m'\in\{0,\lip\}$\@, let
$\nu\in\{m+m',\infty,\omega\}$\@, and let $r\in\{\infty,\omega\}$\@, as
required.  Let $\fG=(\man{M},\sF)$ and $\fH=(\man{N},\sG)$ be
$\C^\nu$-\gcs{}s.  An \defn{isomorphism} from $\fG$ to $\fH$ is a morphism
$(\Phi,\Phi^\sharp)$ such that\/ $\Phi$ is a diffeomorphism and
$L_{\nbhd{V}}$ is an isomorphism (in the category of locally convex
topological vector spaces) for every open $\nbhd{V}\subset\man{N}$\@, where
$L_{\nbhd{V}}$ is such that\/
$\Phi^\sharp_{\nbhd{V}}=L_{\nbhd{V}}|\sG(\nbhd{V})$ as in
Definition~\ref{def:tcs-morphism}\@.\oprocend
\end{definition}

It is now easy to describe the natural isomorphisms.
\begin{proposition}\label{prop:natural-isomorphism}
Let\/ $m\in\integernn$ and\/ $m'\in\{0,\lip\}$\@, let\/
$\nu\in\{m+m',\infty,\omega\}$\@, and let\/ $r\in\{\infty,\omega\}$\@, as
required.  Let\/ $\fG=(\man{M},\sF)$ and\/ $\fH=(\man{N},\sG)$ be\/
$\C^\nu$-\gcs{}s.  A morphism\/ $(\Phi,\Phi^\sharp)$ from\/ $\fG$ to\/ $\fH$
is a natural isomorphism if and only if\/ $\Phi$ is a diffeomorphism and
\begin{equation*}
\sG(\Phi(\nbhd{U}))=\setdef{(\Phi|\nbhd{U})_*X}{X\in\sF(\nbhd{U})}
\end{equation*}
for every open set\/ $\nbhd{U}\subset\man{M}$\@.
\begin{proof}
According to Proposition~\ref{prop:natural-morphism}\@, if
$\nbhd{V}\subset\man{N}$ is open and if $Y\in\sG(\nbhd{V})$\@, we have
$(\Phi|\Phi^{-1}(\nbhd{V}))_*(\Phi^\sharp(Y))=Y$ or
$\Phi^\sharp(Y)=(\Phi|\Phi^{-1}(\nbhd{V}))^*Y$\@.  Since $\Phi^\sharp$ is a
bijection from $\sG(\nbhd{V})$ to $\sF(\Phi^{-1}(\nbhd{V}))$\@, we conclude
that
\begin{equation*}
\sF(\Phi^{-1}(\nbhd{V}))=\setdef{(\Phi|\Phi^{-1}(\nbhd{V}))^*Y}
{Y\in\sG(\nbhd{V})}.
\end{equation*}
This is clearly equivalent to the assertion of the theorem since $\Phi$ must
be a diffeomorphism.
\end{proof}
\end{proposition}

In words, natural isomorphisms simply amount to the natural correspondence of
vector fields under the push-forward $\Phi_*$\@.  (One should verify that
push-forward is continuous as a mapping between locally convex spaces.  This
amounts to proving continuity of composition, and for this we point to places
in the literature from which this can be deduced.  In the smooth and finitely
differentiable cases this can be shown using an argument fashioned after that
from~\cite[Proposition~1]{JNM:69}\@.  In the Lipschitz case, this follows
because the Lipschitz constant of a composition is bounded by the product of
the Lipschitz constants~\cite[Proposition~1.2.2]{NW:99}\@.  In the real
analytic case, this follows from Sublemma~\ref{psublem:compcont} from the
proof of Lemma~\ref{lem:pissy-estimate}\@.)  In particular, if one wishes to
consider only the identity diffeomorphism,~\ie~only consider the ``feedback
part'' of a feedback transformation, we see that the only natural isomorphism
is simply the identity morphism.  In this way we see that the notion of
equivalence for \gcs{}s is either very trivial (it is easy to understand when
systems are equivalent) or very difficult (the study of equivalence classes
contains as a special case the classification of vector fields up to
diffeomorphism), depending on your tastes.  It is our view that the
triviality (or impossibility) of equivalence is a virtue of the formulation
since all structure except that of the manifold and the vector fields has
been removed; there is no extraneous structure.  We refer to
Section~\ref{subsubsec:parameterised} for further discussion.

\subsection{A tautological control system formulation of sub-Riemannian
geometry}\label{subsec:subriemannian}

In our preceding discussion of \gcs{}s, we strove to make connections between
\gcs{}s and standard control models.  We do not wish to give the impression,
however, that \gcs{}s are mere fancy reformulations of standard control
systems.  In this section we give an application, sub-Riemannian geometry,
that illustrates the \emph{per se} value of \gcs{}s.

Let us define the basic structure of sub-Riemannian geometry.
\begin{definition}
Let $m\in\integernn$ and $m'\in\{0,\lip\}$\@, let
$\nu\in\{m+m',\infty,\omega\}$\@, and let $r\in\{\infty,\omega\}$\@, as
required.  A \defn{$\C^\nu$-sub-Riemannian manifold} is a pair
$(\man{M},\metric)$ where $\man{M}$ is a $\C^r$-manifold and $\metric$ is a
$\C^\nu$-tensor field of type $(2,0)$ such that $\metric(x)$ is
positive-semidefinite as a quadratic function on
$\ctb[x]{\man{M}}$\@.\oprocend
\end{definition}

Associated with a sub-Riemannian structure $\metric$ on $\man{M}$ is a
distribution that we now describe.  First of all, we have a map
$\map{\metric^\sharp}{\ctb{\man{M}}}{\tb{\man{M}}}$ defined by
\begin{equation*}
\natpair{\beta_x}{\metric^\sharp(\alpha_x)}=\metric(\beta_x,\alpha_x).
\end{equation*}
We then denote by $\dist{D}_{\metric}=\image(\metric^\sharp)$ the associated
distribution.  Note that $\dist{D}_{\metric}$ is a distribution of class
$\C^\nu$ since, for each $x\in\man{M}$\@, there exist a neighbourhood
$\nbhd{U}$ of $x$ and a family of $\C^\nu$-vector fields $\ifam{X_a}_{a\in
A}$ on $\nbhd{U}$ (namely the images under $\metric^\sharp$ of the coordinate
basis vector fields, if we choose $\nbhd{U}$ to be a coordinate chart domain)
such that
\begin{equation*}
\dist{D}_{\metric,y}=\dist{D}_{\metric}\cap\tb[y]{\man{M}}=
\vecspan[\real]{X_a(y)|\ a\in A}
\end{equation*}
for every $y\in\nbhd{U}$\@.  There is also an associated \defn{sub-Riemannian
metric} for $\dist{D}_{\metric}$\@,~\ie~an assignment to each $x\in\man{M}$
an inner product $\metric(x)$ on $\dist{D}_{\metric,x}$\@.  This is denoted
also by $\metric$ and defined by
\begin{equation*}
\metric(u_x,v_x)=\metric(\alpha_x,\beta_x),
\end{equation*}
where $u_x=\metric^\sharp(\alpha_x)$ and $v_x=\metric^\sharp(\beta_x)$\@, and
where we joyously abuse notation.

An absolutely continuous curve $\map{\gamma}{\interval[a,b]}{\man{M}}$ is
\defn{$\dist{D}_{\metric}$-admissible} if
$\gamma'(t)\in\dist{D}_{\metric,\gamma(t)}$ for almost every
$t\in\interval[a,b]$\@.  The \defn{length} of a
$\dist{D}_{\metric}$-admissible curve $\map{\gamma}{\interval[a,b]}{\man{M}}$
is
\begin{equation*}
\ell_{\metric}(\gamma)=\int_a^b\sqrt{\metric(\gamma'(t),\gamma'(t))}\,\d{t}.
\end{equation*}
As in Riemannian geometry, the length of a $\dist{D}_{\metric}$-admissible
curve is independent of parameterisation, and so curves can be considered to
be defined on $\interval[0,1]$\@.  We can then define the
\defn{sub-Riemannian distance} between $x_1,x_2\in\man{M}$ by
\begin{multline*}
\d_{\metric}(x_1,x_2)=\inf\{\ell_{\metric}(\gamma)|\enspace
\map{\gamma}{\interval[0,1]}{\man{M}}\ \textrm{is an absolutely}\\
\textrm{continuous curve for which}\ \gamma(0)=x_1\ \textrm{and}\
\gamma(1)=x_2\}.
\end{multline*}
One of the problems of sub-Riemannian geometry is to determine length
minimising curves,~\ie~sub-Riemannian geodesics.

A common means of converting sub-Riemannian geometry into a standard control
problem is to choose a $\metric$-orthonormal basis $\ifam{X_1,\dots,X_k}$ for
$\dist{D}_{\metric}$ and so consider the control-affine system with dynamics
prescribed by
\begin{equation*}
F(x,\vect{u})=\sum_{a=1}^ku^aX_a(x),\qquad x\in\man{M},\ \vect{u}\in\real^k.
\end{equation*}
Upon doing this, $\dist{D}_{\metric}$-admissible curves are evidently
trajectories for this control-affine system.  Moreover, for a trajectory
$\map{\xi}{\interval[0,1]}{\man{M}}$ satisfying
\begin{equation*}
\xi'(t)=\sum_{a=1}^ku^a(t)X_a(\xi(t)),
\end{equation*}
we have
\begin{equation*}
\ell_{\metric}(\xi)=\int_0^1\dnorm{\vect{u}(t)}\,\d{t}.
\end{equation*}

The difficulty, of course, with the preceding approach to sub-Riemannian
geometry is that there may be no $\metric$-orthonormal basis for
$\dist{D}_{\metric}$\@.  This can be the case for at least two
reasons:~(1)~the distribution $\dist{D}_{\metric}$ may not have locally
constant rank;~(2)~when the distribution $\dist{D}_{\metric}$ has locally
constant rank, the global topology of $\man{M}$ may prohibit the existence of
a global basis,~\eg~on even-dimensional spheres there is no global basis for
vector fields, orthonormal or otherwise.  However, one can formulate
sub-Riemannian geometry in terms of a \gcs\ in a natural way.  Indeed,
associated to $\dist{D}_{\metric}$ is the \gcs\
$\fG_{\metric}=(\man{M},\sF_{\metric})$\@, where, for an open subset
$\nbhd{U}\subset\man{M}$\@,
\begin{equation*}
\sF_{\metric}(\nbhd{U})=\setdef{X\in\sections[\nu]{\tb{\nbhd{U}}}}
{X(x)\in\dist{D}_{\metric,x},\ x\in\nbhd{U}}.
\end{equation*}
One readily verifies that $\sF_{\metric}$ is a sheaf.

Let us see how we can regard our tautological control system formulation as
that for an ``ordinary'' control system, with a suitable control set, as per
Example~\enumdblref{eg:gcs}{enum:gcs->cs}\@.  First of all, note that the
sheaf $\sF_{\metric}$ is not globally generated; this is because it is a
sheaf,~\cf~Example~\enumdblref{eg:!sheaf}{enum:global-sheaf}\@.  However, it
can be regarded as the sheafification of the globally generated sheaf with
global generators $\sF_{\metric}(\man{M})$\@.
\begin{lemma}
The sheaf\/ $\sF_{\metric}$ is the sheafification of the globally generated
presheaf with generators\/ $\sF_{\metric}(\man{M})$\@.
\begin{proof}
This is a result about sheaf cohomology, and we will not give all details
here.  Instead we will simply point to the main facts from which the
conclusion follows.  First of all, to prove the assertion, it suffices by
Lemma~\ref{lem:sheafify} to show that $\sF_{\metric,x}$ is generated, as a
module over the ring $\gfunc[\nu]{x}{\man{M}}$\@, by germs of global
sections.  In the cases $\nu\in\{m,m+\lip,\infty\}$\@, the fact that the
sheaf of rings of smooth functions admits partitions of unity implies that
the sheaf $\sfunc[\nu]{\man{M}}$ is a fine sheaf of
rings~\cite[Example~3.4(d)]{ROW:08}\@.  It then follows
from~\cite[Example~3.4(e)]{ROW:08} that the sheaf $\sF_{\metric}$ is also
fine and so soft~\cite[Proposition~3.5]{ROW:08}\@.  Because of this, the
cohomology groups of positive degree for this sheaf
vanish~\cite[Proposition~3.11]{ROW:08}\@, and this ensures that germs of
global sections generate all stalks (more or less by definition of cohomology
in degree $1$).  In the case $\nu=\omega$\@, the result is quite nontrivial.
First of all, by a real analytic adaptation
of~\cite[Corollary~H9]{RCG:90b}\@, one can show that $\sF_{\metric}$ is
locally finitely generated.  Then, $\sF_{\metric}$ being a finitely generated
subsheaf of the coherent sheaf $\ssections[\omega]{\tb{\man{M}}}$\@, it is
itself coherent~\cite[Theorem~3.16]{JPD:12}\@.  Then, by Cartan's
Theorem~A~\cite{HC:57}\@, we conclude that $\sF_{\metric,x}$ is generated by
germs of global sections.
\end{proof}
\end{lemma}

By the preceding lemma and Proposition~\ref{prop:sheafify-traj}\@, we can as
well consider the globally generated presheaf with global generators
$\sF_{\metric}(\man{M})$\@, and so trajectories are those of the associated
``ordinary'' control system
$\Sigma_{\metric}=(\man{M},F_{\metric},\cs{C}_{\metric})$\@, where
$\cs{C}_{\metric}=\sF_{\metric}(\man{M})$ and $F_{\metric}(x,X)=X(x)$\@.

Let us next formulate the sub-Riemannian geodesic problem in the framework of
tautological control systems.  First of all, it is convenient when performing
computations to work with energy rather than length as the quantity we are
minimising.  To this end, for an absolutely continuous
$\dist{D}_{\metric}$-admissible curve
$\map{\gamma}{\interval[a,b]}{\man{M}}$\@, we define the \defn{energy} of
this curve to be
\begin{equation*}
E_{\metric}(\gamma)=\frac{1}{2}\int_a^b\metric(\gamma'(t),\gamma'(t))\,\d{t}.
\end{equation*}
A standard argument shows that curves that minimise energy are in 1--1
correspondence with curves that minimise length and are parameterised to have
an appropriate constant speed~\cite[Proposition~1.4.3]{RM:02}\@.  We can and
do, therefore, consider the energy minimisation problem.  We let
$x_1,x_2\in\man{M}$ and let $\sO_{x_1,x_2}$ be the open-loop subfamily for
which the members of $\sO_{x_1,x_2}(\tdomain,\nbhd{U})$ are those vector
fields $X\in\LIsections[\nu]{\tdomain;\sF_{\metric}(\nbhd{U})}$ having the
property that there exist $t_1,t_2\in\tdomain$ with $t_1<t_2$\@,
$\nbhd{U}'\subset\nbhd{U}$\@, and
$\xi\in\Traj(\interval[{t_1},{t_2}],\nbhd{U}',\sO_{\fG_{\metric},X})$ (see
Example~\enumdblref{eg:olsf}{enum:X-ol} for notation) such that
$\xi(t_1)=x_1$ and $\xi(t_2)=x_2$\@.  If
$X\in\sO_{q_1,q_2}(\tdomain,\nbhd{U})$\@, let us denote by $\Traj(X,x_1,x_2)$
those integral curves $\map{\xi}{\interval[{t_1},{t_2}]}{\man{M}}$ for $X$
with the property that $\xi(t_1)=x_1$ and $\xi(t_2)=x_2$\@.  We can then
define
\begin{equation*}
\sC_{\metric}(X)=\inf\setdef{E_{\metric}(\xi)}{\xi\in\Traj(X,x_1,x_2)}.
\end{equation*}
The goal, then, is to find an interval $\tdomain_*\subset\real$\@, an open
set $\nbhd{U}_*$\@, and $X_*\in\sO_{x_1,x_2}(\tdomain_*,\nbhd{U}_*)$ such that
\begin{equation*}
\sC_{\metric}(X_*)\le\sC_{\metric}(X),\qquad
X\in\sO_{x_1,x_2}(\tdomain,\nbhd{U}),\ \tdomain\ \textrm{an interval},
\nbhd{U}\subset\man{M}\ \textrm{open}.
\end{equation*}

Let us apply the classical Maximum Principle of \citet{LSP/VGB/RVG/EFM:61}\@,
leaving aside the technicalities caused by the complicated topology of the
control set.  The dealing with of these technicalities will be the subject of
future work.  We thus suppose that we have a length minimising trajectory
$\xi_*\in\Traj(X_*,x_1,x_2)$ for
$X_*\in\sO_{x_1,x_2}(\tdomain_*,\nbhd{U}_*)$\@.  The \defn{Hamiltonian} for
the system has the form
\begin{equation*}
\mapdef{H_{\metric}}{\ctb{\nbhd{U}_*}\times\sF_{\metric}(\nbhd{U}_*)}
{\real}{(\alpha_x,X)}{\natpair{\alpha_x}{X(x)}+
\lambda_0\tfrac{1}{2}\metric(X(x),X(x)),}
\end{equation*}
where $\lambda_0\in\{0,-1\}$\@.  If we consider only normal
extremals,~\ie~supposing that $\lambda_0=-1$\@, then the Maximum Principle
prescribes that $\map{X_*}{\ctb{\nbhd{U}_*}}{\tb{\nbhd{U}_*}}$ should be a
bundle map over $\id_{\nbhd{U}_*}$ chosen so that $X_*(\alpha_x)$ maximises
the function
\begin{equation*}
v_x\mapsto\natpair{\alpha_x}{v_x}-\tfrac{1}{2}\metric(v_x,v_x).
\end{equation*}
Standard finite-dimensional optimisation gives
$X_*(x)=\metric^\sharp(\alpha_x)$\@.  The \defn{maximum Hamiltonian} is then
obtained by substituting this value of the ``control'' into the Hamiltonian:
\begin{equation*}
\mapdef{\supscr{H}{max}_{\metric}}{\ctb{\man{M}}}{\real}
{\alpha_x}{\tfrac{1}{2}\metric(\alpha_x,\alpha_x).}
\end{equation*}
The normal extremals are then integral curves of the Hamiltonian vector field
associated with the Hamiltonian $\supscr{H}{max}_{\metric}$\@.

The preceding computations, having banished the usual parameterisation by
control, are quite elegant when compared to manner in which one applies the
Maximum Principle to the ``usual'' control formulation of sub-Riemannian
geometry.  The calculations are also more general and global.  However, to
make sense of them, one has to prove an appropriate version of the Maximum
Principle, something which will be forthcoming.  For the moment, we mention
that a significant r\^ole in this will be played by appropriate needle
variations constructed by dragging variations along a trajectory to the final
endpoint.  The manner in which one drags these variations has to do with
linearisation, to which we now turn our attention.

\section{Linearisation of \gcs{}s}\label{sec:linearisation}

As an illustration of the fact that it is possible to do non-elementary
things in the framework of \gcs{}s, we present a fully developed theory for
the linearisation of these systems.  This theory is both satisfying and
revealing.  It is satisfying because it is very simple (if one knows a little
tangent bundle geometry) and it is revealing because, for example, it
clarifies and rectifies the hiccup with classical linearisation theory that
was revealed in Example~\ref{eg:bad-linearise}\@.

Before we begin, it is worth pointing out that, apart from the problem
revealed in Example~\ref{eg:bad-linearise}\@, there are other difficulties
with the very idea of classical Jacobian linearisation to which blind eyes
seem to be routinely turned in practice.  First of all, for models of the
form ``$F(x,u)$\@,'' one must assume that differentiation with respect to $u$
can be done.  For models of this sort, there is no reason to assume the
control set to be a subset of $\real^m$\@, and so one runs into a problem
right away.  Even so, if one restricts to control-affine systems, where the
notion of differentiation with respect to $u$ seems not to be problematic,
one must ignore the fact that the control set is generally not an open set,
and so these derivatives are not so easily made sense of.  Therefore, even
for the typical models one studies in control theory, there are good reasons
to revisit the notion of linearisation.

We point out that geometric linearisation of control-affine systems, and a
Linear Quadratic Regulator theory in this framework, has been carried out by
\citet{ADL/DRT:10}\@.  But even the geometric approach in that work is
refined and clarified by what we present here.

In this section we work with systems of general regularity, only requiring
that they be at least once differentiable so that we can easily define their
linearisation.  For dealing with Lipschitz systems, we will use the following
result.
\begin{lemma}\label{lem:mlip->m-1lip}
For a smooth vector bundle\/ $\map{\pi}{\man{E}}{\man{M}}$ and for\/
$m\in\integerp$\@, if\/ $\xi\in\sections[m+\lip]{\man{E}}$\@, then\/
$j_1\xi\in\sections[(m-1)+\lip]{\jet{1}{\man{E}}}$\@.  Moreover,\/
$\dil{j_{m-1}(j_1\xi)}(x)=\dil{j_m\xi}(x)$ for every\/ $x\in\man{M}$\@.
\begin{proof}
We need to show that $j_{m-1}(j_1\xi)$ is locally Lipschitz.  This, however,
is clear since $j_{m-1}j_1\xi$ is the image of $j_m\xi$ under the injection
of $\jet{m}{\man{E}}$ in
$\jet{m-1}{\jet{1}{\man{E}}}$~\cite[Definition~6.2.25]{DJS:89}\@, and since
$j_m\xi$ is Lipschitz by hypothesis.

The last formula in the statement of the lemma requires us to make sense of
$\dil j_{m-1}(j_1\xi)$\@.  This is made sense of using the fact that, by
Lemma~\ref{lem:Jrdecomp}\@, one has
$\jet{1}{\man{E}}\simeq\ctb{\man{M}}\otimes\man{E}$\@, and so the Riemannian
metric $\metric$ on $\man{M}$\@, the fibre metric $\metric_0$\@, the
Levi-Civita connection $\nabla$ on $\man{M}$\@, and the
$\metric_0$-orthogonal linear connection $\nabla^0$ induce a fibre metric and
linear connection in the vector bundle $\jet{1}{\man{E}}$ as in
Sections~\ref{subsec:Jmdecomp} and~\ref{subsec:olGm}\@.  Now let us examine
the inclusion of $\jet{m}{\man{E}}$ in $\jet{m-1}{\jet{1}{\man{E}}}$ to
verify the final assertion of the lemma.  We use Lemma~\ref{lem:Jrdecomp} to
write
\begin{equation*}
\jet{m}{\man{E}}\simeq\oplus_{j=0}^m\Symalg[j]{\ctb{\man{M}}}
\otimes\man{E}.
\end{equation*}
In this case, the inclusion of $\jet{m}{\man{E}}$ in
$\jet{1}{\jet{m-1}{\man{E}}}$ becomes identified with the natural inclusions
\begin{equation*}
\Symalg[j]{\ctb{\man{M}}}\otimes\man{E}\rightarrow 
\Symalg[j-1]{\ctb{\man{M}}}\otimes\ctb{\man{M}}\otimes\man{E},
\qquad j\in\{0,1,\dots,m-1\},
\end{equation*}
given by 
\begin{equation*}
\alpha^1\odot\dots\odot\alpha^j\otimes e\mapsto
\sum_{k=1}^j\alpha^1\odot\dots\odot 
\alpha^{k-1}\odot\alpha^{k+1}\odot\dots\odot\alpha^j\otimes\alpha^k\otimes e.
\end{equation*}
The fibre metric on $\Symalg[j]{\ctb{\man{M}}}$ is the restriction of that
on $\tensor[j]{\ctb{\man{M}}}$\@.  Thus the preceding inclusion preserves
the fibre metrics since these are defined componentwise on the tensor
product.  Similarly, since the connection in the symmetric and tensor
products is defined so as to satisfy the Leibniz rule for the tensor product,
the injection above commutes with parallel translation.  It now follows from
the definition of dilatation that the final formula in the statement of the
lemma holds.
\end{proof}
\end{lemma}

\subsection{Tangent bundle geometry}

To make the constructions in this section, we recall a little tangent bundle
geometry.  Throughout this section, we let $m\in\integerp$\@,
$m'\in\{0,\lip\}$\@, and let $\nu\in\{m+m',\infty,\omega\}$\@.  We take
$r\in\{\infty,\omega\}$\@, as required.  The meaning of ``$\nu-1$'' is
obvious for all $\nu$\@.  But, to be clear, $\infty-1=\infty$\@,
$\omega-1=\omega$\@, and, given Lemma~\ref{lem:mlip->m-1lip}\@,
$(m+\lip)-1=(m-1)+\lip$\@.

Let $X\in\sections[\nu]{\tb{\man{M}}}$\@.  We will lift $X$ to a vector field
on $\tb{\man{M}}$ in two ways.  The first is the vertical lift, and is
described first by a vector bundle map
$\map{\vlft}{\tbproj{\man{M}}^*\tb{\man{M}}}{\tb{\tb{\man{M}}}}$ as follows.
Let $x\in\man{M}$ and let $v_x,w_x\in\tb[x]{\man{M}}$\@.  The
\textbf{vertical lift} of $u_x$ to $v_x$ is given by
\begin{equation*}
\vlft(v_x,u_x)=\derivatzero{}{t}(v_x+tu_x).
\end{equation*}
Now, if $X\in\sections[\nu]{\tb{\man{M}}}$\@, we define
$\vlift{X}\in\sections[\nu]{\tb{\tb{\man{M}}}}$ by
$\vlift{X}(v_x)=\vlft(v_x,X(x))$\@.  In coordinates $(x^1,\dots,x^n)$ for
$\man{M}$ with $((x^1,\dots,x^n),(v^1,\dots,v^n))$ the associated natural
coordinates for $\tb{\man{M}}$, if $X=X^j\pderiv{}{x^j}$\@, then
$\vlift{X}=X^j\pderiv{}{v^j}$\@.  The vertical lift is a very simple vector
field.  It is tangent to the fibres of $\tb{\man{M}}$\@, and is in fact
constant on each fibre.

The other lift of $X\in\sections[\nu]{\tb{\man{M}}}$ that we shall use is the
\defn{tangent lift}\footnote{This is also frequently called the
\defn{complete lift}\@.  However, ``tangent lift'' so much better captures
the essence of the construction, that we prefer our terminology.  Also, the
dual of the tangent lift is used in the Maximum Principle, and this is then
conveniently called the ``cotangent lift.''} which is the vector field
$\tlift{X}$ on $\tb{\man{M}}$ of class $\C^{\nu-1}$ whose flow is given by
$\flow{\tlift{X}}{t}(v_x)=\tf[x]{\flow{X}{t}}(v_x)$\@.  Therefore,
explicitly,
\begin{equation*}
\tlift{X}(v_x)=\derivatzero{}{t}\tf[x]{\flow{X}{t}}(v_x).
\end{equation*}
In coordinates as above, if $X=X^j\pderiv{}{x^j}$\@, then
\begin{equation}\label{eq:XTcoor}
\tlift{X}=X^j\pderiv{}{x^j}+\pderiv{X^j}{x^k}v^k\pderiv{}{v^j}.
\end{equation}
One recognises the ``linearisation'' of $X$ in this expression, but one
should understand that the second term in this coordinate expression
typically has no meaning by itself.  The flow for $\tlift{X}$ is related to
that for $X$ according to the following commutative diagram:
\begin{equation}\label{eq:XTflow}
\xymatrix{{\tb{\man{M}}}\ar[r]^{\flow{\tlift{X}}{t}}\ar[d]_{\tbproj{\man{M}}}&
{\tb{\man{M}}}\ar[d]^{\tbproj{\man{M}}}\\{\man{M}}\ar[r]_{\flow{X}{t}}&{\man{M}}}
\end{equation}
Thus $\tlift{X}$ projects to $X$ in the sense that
$\tf[v_x]{\tbproj{\man{M}}}(\tlift{X}(v_x))=X(x)$\@.  Moreover, $\tlift{X}$
is a ``linear'' vector field (as befits its appearance in ``linearisation''
below), which means that the diagram
\begin{equation}\label{eq:XTlinear}
\xymatrix{{\tb{\man{M}}}\ar[r]^{\tlift{X}}\ar[d]_{\tbproj{\man{M}}}&
{\tb{\tb{\man{M}}}}\ar[d]^{\tf{\tbproj{\man{M}}}}\\
{\man{M}}\ar[r]_X&{\tb{\man{M}}}}
\end{equation}
defines $\tlift{X}$ as a vector bundle map over $X$\@.

We will be interested in the flow of the tangent lift in the time-varying
case, and the next lemma indicates how this works.
\begin{lemma}\label{lem:XTflow}
Let $m\in\integerp$ and $m'\in\{0,\lip\}$\@, let
$\nu\in\{m+m',\infty,\omega\}$\@, and let $r\in\{\infty,\omega\}$\@, as
required.  Let\/ $\man{M}$ be a\/ $\C^r$-manifold and let\/
$\tdomain\subset\real$ be a time-domain.  For\/
$X\in\LIsections[\nu]{\tdomain;\tb{\man{M}}}$ define
$\map{\tlift{X}}{\tdomain\times\tb{\man{M}}}{\tb{\tb{\man{M}}}}$ by\/
$\tlift{X}(t,v_x)=\tlift{(X(t))}(v_x)$\@.  Then
\begin{compactenum}[(i)]
\item \label{pl:tvtlift1} $\tlift{X}\in\LIsections[\nu-1]{\tdomain;\tb{\tb{\man{M}}}}$\@,
\item \label{pl:tvtlift2} if\/ $(t,t_0,x_0)\in D_X$\@, then\/
$(t,t_0,v_{x_0})\in D_{\tlift{X}}$ for every\/
$v_{x_0}\in\tb[x_0]{\man{M}}$\@, and
\item \label{pl:tvtlift3} $\tlift{X}(t,v_x)=\tderivatzero{}{\tau}\tf[x]{\flow{X}{t+\tau}[t]}(v_x)$\@.
\end{compactenum}
\begin{proof}
\eqref{pl:tvtlift1} Since differentiation with respect to $x$ preserves
measurability in $t$\@,\footnote{Derivatives are limits of sequences of
difference quotients, each of which is measurable, and limits of sequences of
measurable functions are measurable~\cite[Proposition~2.1.4]{DLC:80}\@.} and
since the coordinate expression for $\tlift{X}$ involves differentiating the
coordinate expression for $X$\@, we conclude that $\tlift{X}$ is a
Carath\'eodory vector field.  To show that
$\tlift{X}\in\LIsections[\nu-1]{\tdomain;\tb{\tb{\man{M}}}}$ requires,
according to our definitions of Section~\ref{sec:time-varying}\@, an affine
connection on $\tb{\man{M}}$ and a Riemannian metric on $\tb{\man{M}}$\@.  We
suppose, of course, that we have an affine connection $\nabla$ and a
Riemannian metric $\metric$ on $\man{M}$\@.  For simplicity of some of the
computations below, and without loss of generality, we shall suppose that
$\nabla$ is torsion-free.  In case $r=\omega$\@, we suppose these are real
analytic, according to Lemma~\ref{lem:analytic-conn}\@.  In case $\nu=m+\lip$
for some $m\in\integerp$\@, we assume that $\nabla$ is the Levi-Civita
connection associated with $\metric$\@.

Let us first describe the Riemannian metric on $\tb{\man{M}}$ we shall use.
The affine connection $\nabla$ gives a splitting
$\tb{\tb{\man{M}}}\simeq\tbproj{\man{M}}^*\tb{\man{M}}\oplus
\tbproj{\man{M}}^*\tb{\man{M}}$~\cite[\S11.11]{IK/PWM/JS:93}\@.  We adopt the
convention that the second component of this decomposition is the vertical
component so $\tf[v_x]{\tbproj{\man{M}}}$ restricted to the first component
is an isomorphism onto $\tb[x]{\man{M}}$\@,~\ie~the first component is
``horizontal.''  If $X\in\sections[\nu]{\tb{\man{M}}}$ we denote by
$\hlift{X}\in\sections[\nu]{\tb{\tb{\man{M}}}}$ the unique horizontal vector
field for which $\tf[v_x]{\tbproj{\man{M}}}(\hlift{X}(v_x))=X(x)$ for every
$v_x\in\tb{\man{M}}$\@,~\ie~$\hlift{X}$ is the ``horizontal lift'' of $X$\@.
Let us denote by
$\map{\hor,\ver}{\tb{\tb{\man{M}}}}{\tbproj{\man{M}}^*\tb{\man{M}}}$ the
projections onto the first and second components of the direct sum
decomposition.  This then immediately gives a Riemannian metric
$\tlift{\metric}$ on $\tb{\man{M}}$ by
\begin{equation*}
\tlift{\metric}(X_{v_x},Y_{v_x})=\metric(\hor(X_{v_x}),\hor(Y_{v_x}))+
\metric(\ver(X_{v_x}),\ver(Y_{v_x})).
\end{equation*}
This is called the \defn{Sasaki metric}~\cite{SS:58} in the case that
$\nabla$ is the Levi-Civita connection associated with $\metric$\@.

Now let us determine how an affine connection on $\tb{\man{M}}$ can be
constructed from $\nabla$\@.  There are a number of ways to lift an affine
connection from $\man{M}$ to one on $\tb{\man{M}}$\@, many of these being
described by \citet{KY/SI:73}\@.  We shall use the so-called ``tangent lift''
of $\nabla$\@, which is the unique affine connection $\tlift{\nabla}$ on
$\tb{\man{M}}$ satisfying
$\tlift{\nabla}_{\tlift{X}}\tlift{Y}=\tlift{(\nabla_XY)}$ for
$X,Y\in\sections[\nu]{\tb{\man{M}}}$~\cite[\S7]{KY/SK:66}\@,
\cite[page~30]{KY/SI:73}\@.

We have the following sublemma.
\begin{proofsublemma}\label{psublem:nablaT}
If\/ $X\in\sections[\nu]{\tb{\man{M}}}$\@, if\/ $v_x\in\tb{\man{M}}$\@, if\/
$k\in\integernn$ satisfies\/ $k\le\nu$\@, if\/
$X_1,\dots,X_k\in\sections[\infty]{\tb{\man{M}}}$\@, and if\/
$Z_a\in\{\tlift{X_a},\vlift{X_a}\}$\@,\/ $a\in\{1,\dots,k\}$\@, then the
following formula holds:
\begin{equation*}
(\tlift{\nabla})^{(k)}\tlift{X}(Z_1,\dots,Z_k)=
\begin{cases}\vlift{(\nabla^{(k)}X(X_1,\dots,X_k))},&Z_a\
\textrm{is vertical for some}\ a\in\{1,\dots,k\},\\
\tlift{(\nabla^{(k)}X(X_1,\dots,X_k))},&\textrm{otherwise}.\end{cases}
\end{equation*}
\begin{subproof}
By~\cite[Proposition~7.2]{KY/SK:66}\@, we have
\begin{equation*}
\tlift{\nabla}\tlift{X}(\tlift{X_1})=\tlift{(\nabla X(X_1))},\quad
\tlift{\nabla}\tlift{X}(\vlift{X_1})=\vlift{(\nabla X(X_1))},
\end{equation*}
giving the result when $k=1$\@.  Suppose the result is true for
$k\in\integerp$\@, and let $Z_a\in\{\tlift{X_a},\vlift{X_a}\}$\@,
$a\in\{1,\dots,m+1\}$\@.  First suppose that
$Z_{k+1}=\tlift{X_{k+1}}(v_x)$\@.  We then compute, using the fact that
covariant differentiation commutes with
contraction~\cite[Theorem~7.03(F)]{CTJD/TP:91}\@,
\begin{multline}\label{eq:nablaTm1}
(\tlift{\nabla})^{(k+1)}\tlift{X}(Z_1,\dots,Z_m,Z_{k+1})=
\tlift{\nabla}_{\tlift{X_{k+1}}}((\tlift{\nabla})^{(k)}
\tlift{X})(Z_1,\dots,Z_k)\\
-\sum_{j=1}^k(\tlift{\nabla})^{(k)}\tlift{X}(Z_1,\dots,
\tlift{\nabla}_{\tlift{X_{k+1}}}Z_j,\dots,Z_k).
\end{multline}
We now consider two cases.
\begin{compactenum}
\item \emph{None of $Z_1,\dots,Z_k$ are vertical:} In this case, by the
induction hypothesis,
\begin{equation*}
((\tlift{\nabla})^{(k)}\tlift{X})(Z_1,\dots,Z_k)=
\tlift{(\nabla^{(k)}X)(U_1,\dots,U_k))},
\end{equation*}
and~\cite[Proposition~7.2]{KY/SK:66} gives
\begin{equation*}
\tlift{\nabla}_{\tlift{X_{k+1}}}((\tlift{\nabla})^{(k)}
\tlift{X})(Z_1,\dots,Z_k)=
\tlift{(\nabla_{X_{k+1}}(\nabla^{(k)}X)(U_1,\dots,U_k))}.
\end{equation*}
Again using~\cite[Proposition~7.2]{KY/SK:66} and also using the induction
hypothesis, we have, for $j\in\{1,\dots,k\}$\@,
\begin{equation*}
(\tlift{\nabla})^{(k)}\tlift{X}(Z_1,\dots,\tlift{\nabla}_{\tlift{X_{k+1}}}Z_j,
\dots,Z_k)=\tlift{(\nabla^{(k)}X(U_1,\dots,\nabla_{X_{k+1}}U_j,\dots,U_k))}.
\end{equation*}
Combining the preceding two formulae with~\eqref{eq:nablaTm1} gives the
desired conclusion for $k+1$ in this case.
\item \emph{At least one of $Z_1,\dots,Z_k$ is vertical:} In this case, we
have
\begin{equation*}
((\tlift{\nabla})^{(k)}\tlift{X})(Z_1,\dots,Z_k)=
\vlift{(\nabla^{(k)}X)(U_1,\dots,U_k))}
\end{equation*}
by the induction hypothesis.  Applications
of~\cite[Proposition~7.2]{KY/SK:66} and the induction hypothesis give the
formulae
\begin{equation*}
\tlift{\nabla}_{\tlift{X_{k+1}}}((\tlift{\nabla})^{(k)}
\tlift{X})(Z_1,\dots,Z_k)=
\vlift{(\nabla_{X_{k+1}}(\nabla^{(k)}X)(U_1,\dots,U_k))}.
\end{equation*}
and, for $j\in\{1,\dots,k\}$\@,
\begin{equation*}
(\tlift{\nabla})^{(k)}\tlift{X}(Z_1,\dots,\tlift{\nabla}_{\tlift{X_{k+1}}}Z_j,
\dots,Z_k)=\vlift{(\nabla^{(k)}X(U_1,\dots,\nabla_{X_{k+1}}U_j,\dots,U_k))}.
\end{equation*}
Combining the preceding two formulae with~\eqref{eq:nablaTm1} again gives the
desired conclusion for $k+1$ in this case.
\end{compactenum}

If we take $Z_{k+1}=\vlift{X_{k+1}}$\@, an entirely similar argument gives
the result for this case for $k+1$\@, and so completes the proof of the
sublemma.
\end{subproof}
\end{proofsublemma}

To complete the proof of the lemma, let us for the moment simply regard $X$
as a vector field of class $\C^\nu$\@, not depending on time.  We will make
use of the fact that, for every $v_x\in\tb{\man{M}}$\@,
$\tb[v_x]{\tb{\man{M}}}$ is spanned by vector fields of the form
$\tlift{X_1}+\vlift{Y_1}$ since vertical lifts obviously span the vertical
space and since tangent lifts of nonzero vector fields are complementary to
the vertical space.  Therefore, for a fixed $v_x$\@, we can choose
$X_1,\dots,X_n,Y_1,\dots,Y_n\in\sections[\infty]{\man{M}}$ so that
$\ifam{\tlift{X_1}(v_x),\dots,\tlift{X_n}(v_x)}$ and
$\ifam{\vlift{Y_1}(v_x),\dots,\vlift{Y_n}(v_x)}$ comprise
$\tlift{\metric}$-orthonormal bases for the horizontal and vertical
subspaces, respectively, of $\tb[v_x]{\tb{\man{M}}}$\@.  Note that these
vector fields depend on $v_x$\@, but for the moment we will fix $v_x$\@.  We
use the following formula given by \citet[Lemma~4.5]{MBL/ADL:12} for any
vector field $W$ of class $\C^\nu$ on $\man{M}$\@:
\begin{equation}\label{eq:XT-XH}
\tlift{W}(v_x)=\hlift{W}(v_x)+\vlft(v_x,\nabla_{v_x}W(x)),
\end{equation}
keeping in mind that we are supposing $\nabla$ to be torsion-free.

By the sublemma, if $Z_a=\tlift{X_{j_a}}$\@, $a\in\{1,\dots,k\}$\@, then we
have
\begin{multline}\label{eq:nablaThor}
(\tlift{\nabla})^{(k-1)}\tlift{X}(v_x)(Z_1(v_x),\dots,Z_k(v_x))=
\hlift{(\nabla^{(k-1)}X(x)(X_{j_1}(x),\dots,X_{j_k}(x)))}\\
+\vlft(v_x,\nabla_{v_x}(\nabla^{(k-1)}X(X_{j_1},\dots,X_{j_k}))(x)),
\end{multline}
using~\eqref{eq:XT-XH} with $W=\nabla^{(k-1)}X(X_{j_1},\dots,X_{j_k})$\@.
Again using~\eqref{eq:XT-XH}\@, now with $W=X_{j_a}$\@, we have
\begin{equation*}
\tlift{X_{j_a}}(v_x)=\hlift{X_{j_a}}(v_x)+
\vlft(v_x,\nabla_{v_x}X_{j_a}(x)).
\end{equation*}
Since $\tlift{X_{j_a}}$ was specified so that it is horizontal at $v_x$\@,
its vertical part must be zero, whence $\nabla_{v_x}X_{j_a}(x)=0$\@.
Therefore, expanding the second term on the right in~\eqref{eq:nablaThor}\@,
we get
\begin{multline}\label{eq:nablaThor-expand}
(\tlift{\nabla})^{(k-1)}\tlift{X}(v_x)(Z_1(v_x),\dots,Z_k(v_x))=
\hlift{(\nabla^{(k-1)}X(x)(X_{j_1}(x),\dots,X_{j_k}(x)))}\\
+\vlft(v_x,\nabla^{(k)}X(x)(X_{j_1}(x),\dots,X_{j_k}(x),v_x)).
\end{multline}
Symmetrising this formula with respect to $\{1,\dots,k\}$ gives
\begin{multline}\label{eq:PknablaT1}
P^k_{\tlift{\nabla}}(\tlift{X})(v_x)(Z_1(v_x),\dots,Z_k(v_x))=
\hlift{(P^k_\nabla(X)(x)(X_{j_1}(x),\dots,X_{j_k}(x)))}\\
+\vlft\Bigl(v_x,\nabla_{v_x}P^k_\nabla(X)(x)
(X_{j_1},\dots,X_{j_k})\Bigr),
\end{multline}
where, adopting the notation from Section~\ref{subsec:Jmdecomp}\@,
$P^k_\nabla(X)=\Sym_k\otimes\id_{\tb{\man{M}}}(\nabla^{(k-1)}X)$\@.  Now
consider $Z_a\in\{\tlift{X_{j_a}},\vlift{Y_{j_a}}\}$\@,
$a\in\{1,\dots,k\}$\@, and suppose that at least one of these vector fields
is vertical.  Then, by the sublemma, we immediately have the estimate
\begin{equation}\label{eq:PknablaT2}
P^k_{\tlift{\nabla}}(\tlift{X})(v_x)(Z_1(v_x),\dots,Z_k(v_x))=
\vlift{(P^k_\nabla(X_{j_1}(x),\dots,X_{j_k}(x)))},
\end{equation}
where $\hat{X}_{j_1},\dots,\hat{X}_{j_k}$ are chosen from $X_1,\dots,X_n$ and 
$Y_1,\dots,Y_n$\@, corresponding to the way that $Z_1,\dots,Z_k$ are defined.

Now let us use these formulae in the various regularity classes to obtain the lemma.

$\nu=\infty$\@: Let $K\subset\tb{\man{M}}$ be compact and let
$m\in\integernn$\@.  For the moment, suppose that $X$ is time-independent.
Combining~\eqref{eq:PknablaT1} and~\eqref{eq:PknablaT2}\@, and noting that
they hold as we evaluate $P^m_{\tlift{\nabla}}(\tlift{X})(v_x)$ on a
$\tlift{\metric}$-orthonormal basis for $\tb[v_x]{\tb{\man{M}}}$\@, we obtain
the estimate
\begin{equation*}
\dnorm{P^m_{\tlift{\nabla}}(\tlift{X})(v_x)}_{\tlift{\metric}_m}
\le C(\dnorm{P^m_\nabla(X)(x)}_{\metric_m}+
\dnorm{P^{m+1}_\nabla(X)(x)}_{\metric_{m+1}}\dnorm{v_x}_{\metric}),\qquad
v_x\in K,
\end{equation*}
for some $C\in\realp$\@.  Now, if we make use of the fibre norms induced on
jet bundles as in Section~\ref{subsec:olGm}\@, we have
\begin{equation*}
\dnorm{j_m\tlift{X}(v_x)}_{\ol{\tlift{\metric}}_m}\le
C(\dnorm{j_mX(x)}_{\ol{\metric}_m}+
\dnorm{j_{m+1}X(x)}_{\ol{\metric}_{m+1}}\dnorm{v_x}_{\metric}),\qquad v_x\in K,
\end{equation*}
for some possibly different $C\in\realp$\@.  Since
$v_x\mapsto\dnorm{v_x}_{\metric}$ is bounded on $K$\@, the previous estimate
gives
\begin{equation}\label{eq:jmXTest}
\dnorm{j_m\tlift{X_t}(v_x)}_{\ol{\tlift{\metric}}_m}\le
C\dnorm{j_{m+1}X_t(x)}_{\ol{\metric}_{m+1}},\qquad v_x\in K,\ t\in\tdomain,
\end{equation}
for some appropriate $C\in\realp$\@.

Now we consider time-dependence, supposing that
$X\in\LIsections[\infty]{\tdomain;\tb{\man{M}}}$\@.  Then there exists
$f\in\Lloc^1(\tdomain;\realnn)$ such that
\begin{equation*}
\dnorm{j_{m+1}X_t(x)}_{\ol{\metric}_{m+1}}\le f(t),\qquad x\in K,\ t\in\tdomain.
\end{equation*}
We then immediately have
\begin{equation*}
\dnorm{j_m\tlift{X_t}(v_x)}_{\ol{\tlift{\metric}}_m}\le Cf(t),\qquad
x\in K,\ t\in\tdomain,
\end{equation*}
showing that $\tlift{X}\in\LIsections[\infty]{\tdomain;\tb{\tb{\man{M}}}}$\@,
as desired.

$\nu=m$\@: This case follows directly from the computations in the smooth case.

$\nu=m+\lip$\@: Here we take $m=1$ as the general situation follows by
combining this with the previous case.  We consider $X$ to be
time-independent for the moment.  We let $K\subset\tb{\man{M}}$ be compact.
By Lemma~\ref{lem:dil-deriv} we have
\begin{multline*}
\dil{\tlift{X}}(v_x)=
\inf\{\sup\setdef{\dnorm{\tlift{\nabla}_{Y_{v_y}}
\tlift{X}}_{\tlift{\metric}}}
{v_y\in\closure(\nbhd{W}),\ \dnorm{Y_{v_y}}_{\tlift{\metric}}=1,\ \tlift{X}\
\textrm{differentiable at}\ v_y}|\\
\nbhd{W}\ \textrm{is a relatively compact neighbourhood of}\ v_x\}.
\end{multline*}
Now we make use of Lemma~\ref{lem:Jrdecomp}\@,~\eqref{eq:jmXTest}\@, and the
fact that $K$ is compact, to reduce this to an estimate
\begin{multline*}
\dil{\tlift{X}}(v_x)\le C\inf\{
\sup\setdef{\dnorm{j_2X(y)}_{\ol{\metric}_1}}
{y\in\closure(\nbhd{U}),\ j_1X\ \textrm{differentiable at}\ y}|\\
{\nbhd{U}\ \textrm{a relatively compact neighbourhood of}\ x}\}
\end{multline*}
for some $C\in\realp$ and for every $x\in K$\@.  By Lemma~\ref{lem:dil-deriv}
then gives $\dil{\tlift{X}}(v_x)\le C\dil{j_1X}(x)$ for $x\in K$\@.  From
this we obtain the estimate
\begin{equation*}
\lambda^\lip_K(\tlift{X})\le Cp^{1+\lip}_{\tbproj{\man{M}}(K)}(X).
\end{equation*}
From the proof above in the smooth case, we have
\begin{equation*}
p^0_K(\tlift{X})\le C'p^1_{\tbproj{\man{M}}(K)}(X).
\end{equation*}
Combining these previous two estimates gives
\begin{equation*}
p^\lip_K(\tlift{X})\le Cp^{1+\lip}_{\tbproj{\man{M}}(K)}(X)
\end{equation*}
for some $C\in\realp$\@, and from this, this part of the result follows
easily after adding the appropriate time-dependence.

$\nu=\omega$\@: For the moment, we take $X$ to be time-independent.  The
following sublemma will allow us to estimate the last term
in~\eqref{eq:PknablaT1}\@.
\begin{proofsublemma}
Let\/ $\man{M}$ be a real analytic manifold, let\/ $\nabla$ be a real
analytic affine connection on\/ $\man{M}$\@, let\/ $\metric$ be a real
analytic Riemannian metric on\/ $\man{M}$\@, and let\/ $K\subset\man{M}$ be
compact.  Then there exist\/ $C,\sigma\in\realp$ such that
\begin{equation*}
\dnorm{\nabla^kP^k_\nabla(X)(x)}_{\ol{\metric}_{k+1}}\le
2\dnorm{j_{k+1}X(x)}_{\ol{\metric}_{k+1}}
\end{equation*}
for every\/ $x\in K$ and\/ $k\in\integernn$\@.
\begin{subproof}
We use Lemma~\ref{lem:Jrdecomp} to represent elements of
$\jet{k}{\tb{\man{M}}}$\@.  Following~\cite[\S17.1]{IK/PWM/JS:93}\@, we think
of a connection $\tilde{\nabla}\null^k$ on $\jet{k}{\tb{\man{M}}}$ as being
defined by a vector bundle mapping
\begin{equation*}
\xymatrix{{\jet{k}{\tb{\man{M}}}}\ar[r]^{\tilde{S}_k}\ar[d]&
{\jet{1}{\jet{k}{\tb{\man{M}}}}}\ar[d]\\
{\man{M}}\ar[r]&{\man{M}}}
\end{equation*}
The connection $\nabla^{[k]}$\@, thought of in this way and using the
decomposition of Lemma~\ref{lem:Jrdecomp}\@, gives the associated vector
bundle mapping as zero.  Now, with our identifications, we see that
$P^k_\nabla(X)=j_kX-j_{k-1}X$\@, noting that $\jet{k-1}{\tb{\man{M}}}$ is a
subbundle of $\jet{k}{\tb{\man{M}}}$ with our identification.  Therefore, by
definition of $\nabla^{[k]}$\@,
\begin{equation*}
\nabla^k(P^k_\nabla(X))=\nabla^{[k]}(j_kX-j_{k-1}X)=j_1(j_kX-j_{k-1}X).
\end{equation*}
As we pointed out in the proof of Lemma~\ref{lem:mlip->m-1lip} above, the
inclusion of $\jet{k+1}{\tb{\man{M}}}$ in $\jet{1}{\jet{k}{\tb{\man{M}}}}$
preserves the fibre metric.  Therefore,
\begin{equation*}
\dnorm{\nabla^k(P^k_\nabla(X))(x)}_{\metric_k}\le
\dnorm{j_{k+1}X(x)}_{\ol{\metric}_{k+1}}+
\dnorm{j_kX(x)}_{\ol{\metric}_k}\le2\dnorm{j_{k+1}X(x)}_{\ol{\metric}_{k+1}},
\end{equation*}
as desired.
\end{subproof}
\end{proofsublemma}

Let $K\subset\tb{\man{M}}$ be compact and let
$\vect{a}\in\csd(\integernn;\realp)$\@.  As in the smooth case, but now using
the preceding sublemma, we obtain an estimate
\begin{equation*}
\dnorm{j_m\tlift{X}(v_x)}_{\ol{\tlift{\metric}}_m}\le
C\dnorm{j_{m+1}X(x)}_{\ol{\metric}_{m+1}},
\qquad x\in K,\ m\in\integernn,
\end{equation*}
for some suitable $C\in\realp$\@.

Now, taking $X\in\LIsections[\omega]{\tdomain;\tb{\man{M}}}$\@, there exists
$f\in\Lloc^1(\tdomain;\realnn)$ such that
\begin{equation*}
a'_0a'_1\cdots a'_{m+1}\dnorm{j_{m+1}X_t(x)}_{\ol{\metric}_{m+1}}
\le f(t),\qquad x\in K,\ t\in\tdomain,\ m\in\integernn,
\end{equation*}
where $a'_{j+1}=a_j$\@, $j\in\{1,\dots,m\}$\@, and $a'_0=C$\@.  We then
immediately have
\begin{equation*}
a_0a_1\cdots a_m\dnorm{j_m\tlift{X_t}(v_x)}_{\ol{\tlift{\metric}}_m}
\le f(t),\qquad x\in K,\ t\in\tdomain,\ m\in\integernn,
\end{equation*}
showing that $\tlift{X}\in\LIsections[\omega]{\tdomain;\tb{\tb{\man{M}}}}$\@,
as desired. 

\eqref{pl:tvtlift3} We now prove the third assertion.  It is local, so we
work in a chart.  Thus we assume that we are working in an open subset
$\nbhd{U}\subset\real^n$\@.  We let
$\map{\vect{X}}{\tdomain\times\nbhd{U}}{\real^n}$ be the principal part of
the vector field so that a trajectory for $\vect{X}$ is a curve
$\map{\vect{\xi}}{\tdomain}{\nbhd{U}}$ satisfying
\begin{equation*}
\deriv{}{t}\vect{\xi}(t)=\vect{X}(t,\vect{\xi}(t)),\qquad\ae\ t\in\tdomain.
\end{equation*}
The solution with initial condition $\vect{x}_0$ and $t_0$ we denote by
$t\mapsto\flow*{\vect{X}}(t,t_0,\vect{x}_0)$\@.  For fixed
$(t_0,\vect{x}_0)\in\tdomain\times\nbhd{U}$ and for $t$ sufficiently close to
$t_0$\@, let us define a linear map
$\mat{\Psi}(t)\in\Hom_{\real}(\real^n;\real^n)$ by
\begin{equation*}
\mat{\Psi}(t)\cdot\vect{w}=
\plinder{3}{\flow*{\vect{X}}}(t,t_0,\vect{x}_0)\cdot\vect{w}.
\end{equation*}
We have
\begin{equation*}
\deriv{}{t}\flow*{\vect{X}}(t,t_0,\vect{x}_0)=
\vect{X}(t,\flow*{\vect{X}}(t,t_0,\vect{x}_0)),\qquad\ae\ t,
\end{equation*}
for $t$ sufficiently close to $t_0$\@.  Therefore,
\begin{align*}
\deriv{}{t}\plinder{3}{\flow*{\vect{X}}}(t,t_0,\vect{x}_0)=&\;
\plinder{3}{(\tderiv{}{t}\flow*{\vect{X}}(t,t_0,\vect{x}_0))}\\
=&\;\plinder{2}{\vect{X}(t,\flow*{\vect{X}}(t,t_0,\vect{x}_0))}\cdot
\plinder{3}{\flow*{\vect{X}}(t,t_0,\vect{x}_0)}.
\end{align*}
In the preceding expression, we have used~\cite[Corollary~2.2]{FS/HvdM:00} to
swap the time and spatial derivatives.  This shows that
$t\mapsto\mat{\Psi}(t)$ satisfies the initial value problem
\begin{equation*}
\deriv{}{t}\mat{\Psi}(t)=
\plinder{2}{\vect{X}(t,\flow*{\vect{X}}(t,t_0,\vect{x}_0))}\cdot\mat{\Psi}(t),
\qquad\mat{\Psi}(t_0)=\mat{I}_n.
\end{equation*}
By~\cite[Proposition~C.3.8]{EDS:98}\@, $t\mapsto\mat{\Psi}(t)$ can be defined
for all $t$ such that $(t,t_0,\vect{x}_0)\in D_{\vect{X}}$\@.  Moreover, for
$\vect{v}_0\in\real^n$ (which we think of as being the tangent space at
$\vect{x}_0$), the curve
$t\mapsto\vect{v}(t)\eqdef\mat{\Psi}(t)\cdot\vect{v}_0$ satisfies
\begin{equation*}
\deriv{}{t}\vect{v}(t)=
\plinder{2}{\vect{X}}(t,\flow*{\vect{X}}(t,t_0,\vect{x}_0))\cdot\vect{v}(t).
\end{equation*}

Returning now to geometric notation, the preceding chart computations, after
sifting through the notation, show that
\begin{equation*}
\flow*{\tlift{X}}(t,t_0,v_{x_0})=\tf[x]{\flow*{X}(t,t_0,x_0)}(v_{x_0}),
\end{equation*}
and differentiation with respect to $t$ at $t_0$ gives this part of the
lemma.

\eqref{pl:tvtlift2} This was proved along the way to
proving~\eqref{pl:tvtlift3}\@.
\end{proof}
\end{lemma}

We will also use some features of the geometry of the double tangent
bundle,~\ie~$\tb{\tb{\man{M}}}$\@.  This is an example of what is known as a
``double vector bundle,'' and we refer to~\cite[Chapter~9]{KCHM:05} as a
comprehensive reference.  A review of the structure we describe here can be
found~\cite{MBL/ADL:12}\@, along with an interesting application of this
structure.  We begin by noting that the double tangent bundle possesses two
natural vector bundle structures over
$\map{\tbproj{\man{M}}}{\tb{\man{M}}}{\man{M}}$\@:
\begin{equation*}
\xymatrix{{\tb{\tb{\man{M}}}}\ar[r]^{\tbproj{\tb{\man{M}}}}
\ar[d]_{\tf{\tbproj{\man{M}}}}&{\tb{\man{M}}}\ar[d]^{\tbproj{\man{M}}}\\
{\tb{\man{M}}}\ar[r]_{\tbproj{\man{M}}}&{\man{M}}}\qquad\qquad
\xymatrix{{\tb{\tb{\man{M}}}}\ar[r]^{\tf{\tbproj{\man{M}}}}
\ar[d]_{\tbproj{\tb{\man{M}}}}&{\tb{\man{M}}}\ar[d]^{\tbproj{\man{M}}}\\
{\tb{\man{M}}}\ar[r]_{\tbproj{\man{M}}}&{\man{M}}}
\end{equation*}
The left vector bundle structure is called the \defn{primary vector bundle}
and the right the \defn{secondary vector bundle}\@.  We shall denote vector
addition in the vector bundles as follows.  If $u,v\in\tb{\tb{\man{M}}}$
satisfy $\tbproj{\tb{\man{M}}}(u)=\tbproj{\tb{\man{M}}}(v)$\@, then the sum
of $u$ and $v$ in the primary vector bundle is denoted by $u+_1v$\@.  If
$u,v\in\tb{\tb{\man{M}}}$ satisfy
$\tf{\tbproj{\man{M}}}(u)=\tf{\tbproj{\man{M}}}(v)$\@, then the sum of $u$
and $v$ in the secondary vector bundle is denoted by $u+_2v$\@.

The two vector bundle structures admit a naturally defined isomorphism
between them, described as follows.  Let $\rho$ be a smooth map from a
neighbourhood of $(0,0)\in\real^2$ to $\man{M}$\@.  We shall use coordinates
$(s,t)$ for $\real^2$\@.  For fixed $s$ and $t$ define
$\rho_s(t)=\rho^t(s)=\rho(s,t)$\@.  We then denote
\begin{equation*}
\pderiv{}{t}\rho(s,t)=\deriv{}{t}\rho_s(t)\in\tb[\rho(s,t)]{\man{M}},\quad
\pderiv{}{s}\rho(s,t)=\deriv{}{s}\rho^t(s)\in\tb[\rho(s,t)]{\man{M}}.
\end{equation*}
Note that $s\mapsto\pderiv{}{t}\rho(s,t)$ is a curve in $\tb{\man{M}}$ for
fixed $t$\@.  The tangent vector field to this curve we denote by
\begin{equation*}
s\mapsto\pderiv{}{s}\pderiv{}{t}\rho(s,t)\in
\tb[\pderiv{}{t}\rho(s,t)]{\tb{\man{M}}}.
\end{equation*}
We belabour the development of the notation somewhat since these partial
derivatives are not the usual partial derivatives from calculus, although the
notation might make one think they are.  For example, we do not generally
have equality of mixed partials,~\ie~generally we have
\begin{equation*}
\pderiv{}{s}\pderiv{}{t}\rho(s,t)\not=\pderiv{}{t}\pderiv{}{s}\rho(s,t).
\end{equation*}

Now let $\rho_1$ and $\rho_2$ be smooth maps from a neighbourhood of
$(0,0)\in\real^2$ to $\man{M}$\@.  We say two such maps are
\textbf{equivalent} if
\begin{gather*}
\pderiv{}{s}\pderiv{}{t}\rho_1(0,0)=
\pderiv{}{s}\pderiv{}{t}\rho_2(0,0).
\end{gather*}
To the equivalence classes of this equivalence relation, we associate points
in $\tb{\tb{\man{M}}}$ by
\begin{equation*}
[\rho]\mapsto\pderiv{}{s}\pderiv{}{t}\rho(0,0).
\end{equation*}
The set of equivalence classes is easily seen to be exactly the double
tangent bundle $\tb{\tb{\man{M}}}$\@.  We easily verify that
\begin{equation}\label{eq:rho-project}
\tbproj{\tb{\man{M}}}([\rho])=\pderiv{}{t}\rho(0,0),\quad
\tf{\tbproj{\man{M}}}([\rho])=\pderiv{}{s}\rho(0,0).
\end{equation}

Next, using the preceding representation of points in $\tb{\tb{\man{M}}}$\@,
we relate the two vector bundle structures for $\tb{\tb{\man{M}}}$ by
defining a canonical involution of $\tb{\tb{\man{M}}}$\@.  If $\rho$ is a
smooth map from a neighbourhood of $(0,0)\in\real^2$ into $M$\@, define
another such map by $\bar{\rho}(s,t)=\rho(t,s)$\@.  We then define the
\defn{canonical tangent bundle involution} as the map
$\map{I_{\man{M}}}{\tb{\tb{\man{M}}}}{\tb{\tb{\man{M}}}}$ given by
$I_{\man{M}}([\rho])=[\bar{\rho}]$\@.  Clearly $I_{\man{M}}\circ
I_{\man{M}}=\id_{\tb{\tb{\man{M}}}}$\@.  In a natural coordinate chart for
$\tb{\tb{\man{M}}}$ associated to a natural coordinate chart for
$\tb{\man{M}}$\@, the local representative of $I_{\man{M}}$ is
\begin{equation*}
((\vect{x},\vect{v}),(\vect{u},\vect{w}))\mapsto
((\vect{x},\vect{u}),(\vect{v},\vect{w})).
\end{equation*}
One readily verifies that $I_{\man{M}}$ is a vector bundle isomorphism from
$\tb{\tb{\man{M}}}$ with the primary (\resp~secondary) vector bundle
structure to $\tb{\tb{\man{M}}}$ with the secondary (\resp~primary) vector
bundle structure~\cite[Lemma~A.4]{MBL/ADL:12}\@.

The following technical lemma is Lemma~A.5 from~\cite{MBL/ADL:12}\@.
\begin{lemma}\label{lem:IMvlft}
If\/ $w\in\tb{\tb{\man{M}}}$ satisfies\/ $\tbproj{\tb{\man{M}}}(w)=v$ and\/
$\tf{\tbproj{\man{M}}}=u$ and if\/ $z\in\tb[x]{\man{M}}$\@, then
\begin{equation*}
w+_2I_{\man{M}}\circ\vlft(u,z)=w+_1\vlft(v,z).
\end{equation*}
\end{lemma}

The final piece of tangent bundle geometry we will consider concerns
presheaves and sheaves of sets of vector fields on tangent bundles.  We shall
need the following natural notion of such a presheaf.
\begin{definition}
Let $m\in\integernn$ and $m'\in\{0,\lip\}$\@, let
$\nu\in\{m+m',\infty,\omega\}$\@, and let $r\in\{\infty,\omega\}$\@, as
required.  Let $\man{M}$ be a $\C^r$-manifold and let $\sG$ be a presheaf of
sets of vector fields of class\/ $\C^\nu$ on $\tb{\man{M}}$\@.  The presheaf
$\sG$ is \defn{projectable} if
\begin{equation*}\eqoprocend
\sG(\nbhd{W})=\setdef{Z|\nbhd{W}}
{Z\in\sG(\tbproj{\man{M}}^{-1}(\tbproj{\man{M}}(\nbhd{W})))}.
\end{equation*}
\end{definition}

The idea is that a projectable sheaf is determined by the local sections over
the open sets $\tbproj{\man{M}}^{-1}(\nbhd{U})$ for $\nbhd{U}\subset\man{M}$
open.

\subsection{Linearisation of systems}

Throughout this section, unless stated otherwise, we let $m\in\integerp$\@,
$m'\in\{0,\lip\}$\@, and let $\nu\in\{m+m',\infty,\omega\}$\@.  We take
$r\in\{\infty,\omega\}$\@, as required.

When linearising, one typically does so about a trajectory.  We will do this
also.  But before we do so, let us provide the notion of the linearisation of
a \emph{system}\@.  The result, gratifyingly, is a system on the tangent
bundle.  Before we produce the definition, let us make a motivating
computation.  We let $\fG=(\man{M},\sF)$ be a globally generated \gcs\ of
class $\C^\nu$\@.  By Example~\enumdblref{eg:gcs}{enum:gcs->cs}\@, we have
the corresponding $\C^\nu$-control system
$\Sigma_{\fG}=(\man{M},F_{\sF},\cs{C}_{\sF})$ with
$\cs{C}_{\sF}=\sF(\man{M})$ and $F_{\sF}(x,X)=X(x)$\@.  This is a control
system whose control set is a vector space, and so is a candidate for
classical Jacobian linearisation, provided one is prepared to overlook
technicalities of differentiation in locally convex spaces\ldots and we are
for the purposes of this motivational computation.  In Jacobian linearisation
one considers perturbations of state and control.  In our framework, we
linearise about a state/control $(x,X)$\@.  We perturb the state by
considering a $\C^1$-curve $\map{\gamma}{J}{\man{M}}$ defined on an interval
$J$ for which $0\in\interior(J)$ and with $\gamma'(0)=v_x$\@.  Thus we
perturb the state in the direction of $v_x$\@.  We perturb the control from
$X$ in the direction of $Y\in\sF(\man{M})$ by considering a curve of controls
$s\mapsto X+sY$\@.  Let us then define $\map{\sigma}{\nbhd{N}}{\man{M}}$ on a
neighbourhood $\nbhd{N}$ of $(0,0)\in\real^2$ by
\begin{equation*}
\sigma(t,s)=\flow{X+sY}{t}(\gamma(s));
\end{equation*}
thus $\sigma(t,s)$ gives the flow at time $t$ corresponding to the
perturbation at parameter $s$\@.  Now we compute
\begin{align*}
\pderiv{}{t}\pderiv{}{s}\sigma(t,s)=&\;
\pderiv{}{t}\pderiv{}{s}\flow{X+sY}{t}(\gamma(s))\\
=&\;\pderiv{}{t}\pderiv{}{s}\flow{X}{t}(\gamma(s))+
\pderiv{}{t}\pderiv{}{s}\flow{X+sY}{t}(x)\\
=&\;\pderiv{}{t}\tf[x]{\flow{X}{t}}(\gamma'(s))+
I_{\man{M}}\Bigr(\pderiv{}{s}\pderiv{}{t}\flow{X+sY}{t}(x)\Bigr)\\
=&\;\pderiv{}{t}\tf[x]{\flow{X}{t}}(\gamma'(s))+
I_{\man{M}}\Bigr(\pderiv{}{s}(X+sY)(\flow{X+sY}{t}(x))\Bigr),
\end{align*}
from which we have
\begin{equation}\label{eq:linear-derivation}
\pderiv{}{t}\pderiv{}{s}\sigma(0,0)=
\tlift{X}(v_x)+I_{\man{M}}(\vlft(X(x),Y(x)))=\tlift{X}(v_x)+\vlift{Y}(v_x),
\end{equation}
using Lemma~\ref{lem:IMvlft}\@.

The formula clearly suggests what the linearisation of a \gcs\ should be.
However, we need the following lemma to make a sensible definition in our
sheaf framework.
\begin{lemma}\label{lem:linear-presheaf}
Let\/ $m\in\integerp$ and\/ $m'\in\{0,\lip\}$\@, let\/
$\nu\in\{m+m',\infty,\omega\}$\@, and let\/ $r\in\{\infty,\omega\}$\@, as
required.  Let\/ $\sF$ be a presheaf of sets of\/ $\C^\nu$-vector fields on
a\/ $\C^r$-manifold\/ $\man{M}$\@.  Then there exist unique projectable
presheaves\/ $\tlift{\sF}$ and\/ $\vlift{\sF}$ of\/ $\C^{\nu-1}$-vector
fields and\/ $\C^\nu$-vector fields on\/ $\tb{\man{M}}$ with the property
that
\begin{equation*}
\tlift{\sF}(\tbproj{\man{M}}^{-1}(\nbhd{U}))=
\setdef{\tlift{X}}{X\in\sF(\nbhd{U})}
\end{equation*}
and
\begin{equation*}
\vlift{\sF}(\tbproj{\man{M}}^{-1}(\nbhd{U}))=
\setdef{\vlift{X}}{X\in\sF(\nbhd{U})}
\end{equation*}
for every open set\/ $\nbhd{U}\subset\man{M}$\@.  Moreover,
\begin{compactenum}[(i)]
\item \label{pl:linear-sheaf1} $\tlift{\sF}$ is a sheaf if and only if\/ $\sF$
is a sheaf,
\item \label{pl:linear-sheaf2} $\vlift{\sF}$ is a sheaf if and only if\/
$\sF$ is a sheaf,
\item \label{pl:linear-sheaf3} $\Sh(\tlift{\sF})=\tlift{\Sh(\sF)}$\@, and
\item \label{pl:linear-sheaf4} $\Sh(\vlift{\sF})=\vlift{\Sh(\sF)}$\@.
\end{compactenum}
\begin{proof}
Let $\nbhd{W}\subset\tb{\man{M}}$ be open and note that
$\nbhd{U}_{\nbhd{W}}=\tbproj{\man{M}}(\nbhd{W})$ is open.  For
$\nbhd{W}\subset\tb{\man{M}}$ open we define
\begin{equation*}
\tlift{\sF}(\nbhd{W})=
\setdef{\tlift{X}|\nbhd{W}}{X\in\sF(\nbhd{U}_{\nbhd{W}})}
\end{equation*}
and
\begin{equation*}
\vlift{\sF}(\nbhd{W})=
\setdef{\vlift{X}|\nbhd{W}}{X\in\sF(\nbhd{U}_{\nbhd{W}})}.
\end{equation*}
If $\nbhd{W},\nbhd{W}'\subset\tb{\man{M}}$ are open with
$\nbhd{W}'\subset\nbhd{W}$ and if
$\tlift{X}|\nbhd{W}\in\tlift{\sF}(\nbhd{W})$\@, then, for
$v_x\in\nbhd{W}'$\@, we have
\begin{equation*}
(\tlift{X}(v_x)|\nbhd{W}')(v_x)=
(\tlift{(X|\nbhd{U}_{\nbhd{W}'})})(v_x),
\end{equation*}
this making sense since $\tlift{X}(v_x)$ depends only on the values of $X$ in
a neighbourhood of $x$\@, and since $\nbhd{U}_{\nbhd{W}'}$ contains a
neighbourhood of $x$ if $v_x\in\nbhd{W}'$\@.  In any case, we have that
\begin{equation*}
\tlift{X}|\nbhd{W}'\in\tlift{\sF}(\nbhd{W}'),
\end{equation*}
which shows that $\tlift{\sF}$ is a presheaf.  A similar argument, of course,
works for $\vlift{\sF}$\@.  This gives the existence assertion of the lemma.
Uniqueness follows immediately from the requirement that $\tlift{\sF}$ and
$\vlift{\sF}$ be projectable.

\eqref{pl:linear-sheaf1} Suppose that $\sF$ is a sheaf.  We shall first show
that $\tlift{\sF}$ is a sheaf.  Let $\nbhd{W}\subset\tb{\man{M}}$ be open,
and let $\ifam{\nbhd{W}_a}_{a\in A}$ be an open cover of $\nbhd{W}$\@.  Let
$Z_a\in\tlift{\sF}(\nbhd{W}_a)$\@, supposing that
\begin{equation*}
Z_a|\nbhd{W}_a\cap\nbhd{W}_b=Z_b|\nbhd{W}_a\cap\nbhd{W}_b
\end{equation*}
for $a,b\in A$\@.  For each $a\in A$\@, we have, by our definition of
$\tlift{\sF}$ above, $Z_a=\tlift{X_a}|\nbhd{W}_a$ for
$X_a\in\sF(\nbhd{U}_{\nbhd{W}_a})$\@.  Using the fact that
$\sections[\nu-1]{\tb{\tb{\man{M}}}}$ is a sheaf, we infer that there exists
$Z\in\sections[\nu-1]{\tb{\tb{\man{M}}}}$ such that
$Z|\nbhd{W}_a=\tlift{X_a}|\nbhd{W}_a$ for each $a\in A$\@.  Now, for each
$x\in\nbhd{U}_{\nbhd{W}}$\@, let us fix $a_x\in A$ such that
$x\in\tbproj{\man{M}}(\nbhd{W}_a)$\@.  Note that
$Z|\nbhd{W}_{a_x}=\tlift{X_{a_x}}|\nbhd{W}_{a_x}$ and so there is a
neighbourhood $\nbhd{U}_x\subset\nbhd{U}_{\nbhd{W}_{a_x}}$ of $x$ and
$X_x\in\sections[\nu-1]{\tb{\nbhd{U}_x}}$ such that
$X_x=X_{a_x}|\nbhd{U}_x$\@.  In particular, $X_x\in\sF(\nbhd{U}_x)$\@.
Moreover, since $\tlift{\sF}$ is projectable, we can easily see that
$[X_x]_x$ is independent of the rule for choosing $a_x$\@.  Now let
$x_1,x_2\in\man{M}$ and let $x\in\nbhd{U}_{x_1}\cap\nbhd{U}_{x_2}$\@.  By
projectability of $\tlift{\sF}$\@, there exist a neighbourhood
$\nbhd{V}_x\subset\nbhd{U}_{x_1}\cap\nbhd{U}_{x_2}$ and
$X'_x\in\sF(\nbhd{V}_x)$ such that
\begin{equation*}
\tlift{X_{a_{x_j}}}|\nbhd{W}_{a_{x_j}}\cap\tbproj{\man{M}}^{-1}(\nbhd{V}_x)=
\tlift{(X'_x)}|\nbhd{W}_{a_{x_j}},\qquad j\in\{1,2\}.
\end{equation*}
We conclude, therefore, that $X_{x_1}(x)=X_{x_2}(x)$\@.  Thus we have an open
covering $\ifam{\nbhd{U}_x}_{x\in\nbhd{U}_{\nbhd{W}}}$ of
$\nbhd{U}_{\nbhd{W}}$ and local sections $X_x\in\sF(\nbhd{U}_x)$ pairwise
agreeing on intersections.  Since $\sF$ is a sheaf, there exists
$X\in\sF(\nbhd{U}_{\nbhd{W}})$ such that $X|\nbhd{U}_x=X_x$ for each
$x\in\nbhd{U}_{\nbhd{W}}$\@.  Since
\begin{equation*}
\tlift{X}|\nbhd{W}_{a_x}\cap\tbproj{\man{M}}^{-1}(\nbhd{U}_x)=
\tlift{X_x}|\nbhd{W}_{a_x}\cap\tbproj{\man{M}}^{-1}(\nbhd{U}_x)=
\tlift{X_{a_x}}|\nbhd{W}_{a_x}\cap\tbproj{\man{M}}^{-1}(\nbhd{U}_x),
\end{equation*}
projectability of $\tlift{\sF}$ allows us to conclude that
$Z=\tlift{X}|\nbhd{W}$\@.

Now suppose that $\tlift{\sF}$ is a sheaf and let $\nbhd{U}\subset\man{M}$ be
open, let $\ifam{\nbhd{U}_a}_{a\in A}$ be an open covering of $\nbhd{U}$\@,
and let $X_a\in\sF(\nbhd{U}_a)$\@, $a\in A$ be such that
$X_a|\nbhd{U}_a\cap\nbhd{U}_b=X_b|\nbhd{U}_a\cap\nbhd{U}_b$\@.  This implies
that
\begin{equation*}
\tlift{X_a}|\tbproj{\man{M}}^{-1}(\nbhd{U}_a\cap\nbhd{U}_b)=
\tlift{X_b}|\tbproj{\man{M}}^{-1}(\nbhd{U}_a\cap\nbhd{U}_b).
\end{equation*}
Therefore, by hypothesis, there exists $X\in\sF(\nbhd{U})$ such that
$\tlift{X}|\tbproj{\man{M}}^{-1}(\nbhd{U}_a)=\tlift{X_a}$ for each $a\in
A$\@.  Projecting to $\man{M}$ gives $X|\nbhd{U}_a=X_a$ for each $a\in A$\@,
showing that $\sF$ is a sheaf.

\eqref{pl:linear-sheaf2} To show that $\vlift{\sF}$ is a sheaf can be made
with an identically styled argument as above in showing that $\tlift{\sF}$ is
a sheaf.  The argument, indeed, is even easier since vertical lifts do not
depend on the value of their projections in a neighbourhood of a point in
$\tb{\man{M}}$\@, only on the projection at the point.

\eqref{pl:linear-sheaf3} Let $\nbhd{W}\subset\tb{\man{M}}$ be open and let
$Z\in\Sh(\tlift{\sF})(\nbhd{W})$\@.  This means that, for each
$v_x\in\nbhd{W}$\@, $[Z]_{v_x}\in\tlift{\sF_{0,v}}$\@.  Therefore, there
exist a neighbourhood $\nbhd{W}_{v_x}$ of $v_x$ and
$X_x\in\sF(\nbhd{U}_{\nbhd{W}_{v_x}})$ such that
$Z|\nbhd{W}_{v_x}=\tlift{X_x}|\nbhd{W}_{v_x}$\@.  We now proceed as in the
preceding part of the proof.  Thus, for each $x\in\nbhd{U}_{\nbhd{W}}$ let us
fix $v_x\in\nbhd{W}$\@.  Note that
$Z|\nbhd{W}_{v_x}=\tlift{X_{v_x}}|\nbhd{W}_{v_x}$ and so there is a
neighbourhood $\nbhd{U}_x\subset\nbhd{U}_{\nbhd{W}_{v_x}}$ of $x$ and
$X_x\in\sections[\nu-1]{\tb{\nbhd{U}_x}}$ such that
$X_x=X_{v_x}|\nbhd{U}_x$\@.  In particular, $X_x\in\sF(\nbhd{U}_x)$\@.
Moreover, since $\tlift{\sF}$ is projectable, we can easily see that
$[X_x]_x$ is independent of the rule for choosing $v_x\in\nbhd{W}$\@.  Now
let $x_1,x_2\in\man{M}$ and let $x\in\nbhd{U}_{x_1}\cap\nbhd{U}_{x_2}$\@.  By
projectability of $\tlift{\sF}$\@, there exist a neighbourhood
$\nbhd{V}_x\subset\nbhd{U}_{x_1}\cap\nbhd{U}_{x_2}$ and
$X'_x\in\sF(\nbhd{V}_x)$ such that
\begin{equation*}
\tlift{X_{v_{x_j}}}|\nbhd{W}_{v_{x_j}}\cap\tbproj{\man{M}}^{-1}(\nbhd{V}_x)=
\tlift{(X'_x)}|\nbhd{W}_{v_{x_j}},\qquad j\in\{1,2\}.
\end{equation*}
We conclude, therefore, that $X_{x_1}(x)=X_{x_2}(x)$\@.  Thus we have an open
covering $\ifam{\nbhd{U}_x}_{x\in\nbhd{U}_{\nbhd{W}}}$ and local sections
$X_x\in\sF(\nbhd{U}_x)$ pairwise agreeing on intersections.  Thus there
exists $X\in\Sh(\sF(\nbhd{U}_{\nbhd{W}}))$ such that
$X|\nbhd{U}_x=X_x$ for each $x\in\nbhd{U}_{\nbhd{W}}$\@.  Since
\begin{equation*}
\tlift{X}|\nbhd{W}_{v_x}\cap\tbproj{\man{M}}^{-1}(\nbhd{U}_x)=
\tlift{X_x}|\nbhd{W}_{v_x}\cap\tbproj{\man{M}}^{-1}(\nbhd{U}_x)=
\tlift{X_{v_x}}|\nbhd{W}_{v_x}\cap\tbproj{\man{M}}^{-1}(\nbhd{U}_x),
\end{equation*}
projectability of $\Sh(\tlift{\sF})$ allows us to conclude that
$Z=\tlift{X}|\nbhd{W}$\@,~\ie~$Z\in\tlift{\Sh(\sF)}(\nbhd{W})$\@.

\eqref{pl:linear-sheaf4} A similar argument as in the preceding part of the
proof works to give this part of the proof as well.
\end{proof}
\end{lemma}

With the preceding computations and sheaf lemma as motivation, we make the
following definition.
\begin{definition}
Let $m\in\integerp$ and $m'\in\{0,\lip\}$\@, let
$\nu\in\{m+m',\infty,\omega\}$\@, and let $r\in\{\infty,\omega\}$\@, as
required.  Let $\fG=(\man{M},\sF)$ be a $\C^\nu$-\gcs.  The
\defn{linearisation} of $\fG$ is the $\C^{\nu-1}$-\gcs\
$\tf{\fG}=(\tb{\man{M}},\tf{\sF})$\@, where the projectable presheaf of sets
of vector fields $\tf{\sF}$ is characterised uniquely by the requirement
that, for every open subset $\nbhd{U}\subset\man{M}$\@,
\begin{equation*}\eqoprocend
\tf{\sF}(\tbproj{\man{M}}^{-1}(\nbhd{U}))=
\setdef{\tlift{X}+\vlift{Y}}{X,Y\in\sF(\nbhd{U})}.
\end{equation*}
\end{definition}

This definition may look a little strange at a first glance.  However, as we
go along, we shall use the definition in more commonplace settings, and we
will see then that it connects to more familiar constructions.

\subsection{Trajectories for linearisations}\label{subsec:linear-traj}

As a \gcs, $\tf{\fG}$ provides a forum for all of the constructions of
Sections~\ref{subsec:gcs}\@,~\ref{subsec:ol}\@, and~\ref{subsec:traj}
concerning such systems.  In particular, the linearisation has trajectories,
so let us look at these.

Let us first think about open-loop systems.  By definition, an open-loop
system for $\tf{\fG}$ is a triple $(Z,\tdomain,\nbhd{W})$ with
$\tdomain\subset\real$ an interval, $\nbhd{W}\subset\tb{\man{M}}$ an open
set, and $Z\in\LIsections[\nu-1]{\tdomain;\tf{\sF}(\nbhd{W})}$\@.  Thus
$Z(t)=\tlift{X(t)}+\vlift{Y(t)}$ for
$\map{X,Y}{\tdomain}{\sF(\tbproj{\man{M}}(\nbhd{W}))}$\@.  We will write
$Z=\tlift{X}+\vlift{Y}$ with the understanding that this means precisely what
we have just written.  We should, however, verify that $X$ and $Y$ have
useful properties.
\begin{lemma}\label{lem:open-loop-lin}
Let\/ $m\in\integerp$ and\/ $m'\in\{0,\lip\}$\@, let\/
$\nu\in\{m+m',\infty,\omega\}$\@, and let\/ $r\in\{\infty,\omega\}$\@, as
required.  Let\/ $\fG=(\man{M},\sF)$ be a\/ $\C^\nu$-\gcs\ with
linearisation\/ $\tf{\fG}=(\tb{\man{M}},\tf{\sF})$\@.  Let\/ $\tdomain$ be a
time-domain and let\/ $\nbhd{W}\subset\tb{\man{M}}$ be open.  If\/
$Z\in\LIsections[\nu-1]{\tdomain;\tf{\sF}(\nbhd{W})}$ is given by
\begin{equation*}
Z(t,v_x)=\tlift{X}(t,v_x)+\vlift{Y}(t,v_x)
\end{equation*}
for maps\/
$\map{X,Y}{\tdomain\times\tbproj{\man{M}}(\nbhd{W})}{\tb{\man{M}}}$ for
which\/ $X_t,Y_t\in\sections[\nu]{\tdomain;\tbproj{\man{M}}(\nbhd{W})}$ for
every\/ $t\in\tdomain$\@, then
$X\in\LIsections[\nu]{\tdomain;\sF(\tbproj{\man{M}}(\nbhd{W}))}$ and\/ $Y\in\LIsections[\nu-1]{\tdomain;\sF(\tbproj{\man{M}}(\nbhd{W}))}$\@.
\begin{proof}
It is possible to make oneself believe the lemma by a coordinate computation.
However, we shall give a coordinate-free proof.  To do this, we will use the
Riemannian metric $\tlift{\metric}$ and the affine connection
$\tlift{\nabla}$ on $\tb{\man{M}}$ defined by a Riemannian metric $\metric$
and affine connection $\nabla$ on $\man{M}$\@, as described in the proof of
Lemma~\ref{lem:XTflow}\@.  For simplicity, and since we will make use of some
formulae derived in the proof of Lemma~\ref{lem:XTflow} where this assumption
was made, we suppose that $\nabla$ is torsion-free.

Since we will be calculating iterated covariant differentials as in
Section~\ref{subsec:COinfty-vb}\@, only now using the affine connection
$\tlift{\nabla}$ on $\tb{\man{M}}$\@, we should also think about the
character of $\tensor[k]{\ctb{\tb{\man{M}}}}$\@.  For
$v_x\in\tb[x]{\man{M}}$\@, $\tf[v_x]{\tbproj{\man{M}}}$ is a surjective
linear mapping from $\tb[v_x]{\tb{\man{M}}}$ to $\tb[x]{\man{M}}$\@.  Thus
its dual, $\dual{(\tf[v_x]{\tbproj{\man{M}}})}$\@, is an injective linear
mapping from $\ctb[x]{\man{M}}$ to $\ctb[v_x]{\tb{\man{M}}}$\@.  It induces,
therefore, an injective linear mapping from $\tensor[k]{\ctb[x]{\man{M}}}$
to
$\tensor[k]{\ctb[v_x]{\tb{\man{M}}}}$~\cite[Proposition~III.5.2.2]{NB:89a}\@.
\citet{KY/SK:66} call this the \defn{vertical lift} of
$\tensor[k]{\ctb{\man{M}}}$ into $\tensor[k]{\ctb{\tb{\man{M}}}}$\@.  Note
that vertically lifted tensors, thought of as multilinear maps, vanish if
they are given a vertical vector as one of their arguments,~\ie~they are
``semi-basic'' (in fact, they are even ``basic'').  Note that
$\ctb[v_x]{\tb{\man{M}}}\simeq\ctb[x]{\man{M}}\oplus\ctb[x]{\man{M}}$ by
dualising the splitting of the tangent bundle.  So as to notationally
distinguish between the two components of the direct sum, let us denote the
first component by $(\ctb[x]{\man{M}})_1$ and the second component by
$(\ctb[x]{\man{M}})_2$\@, noting that the first component is defined to be
the image of the canonical injection from $\ctb[x]{\man{M}}$ to
$\ctb[v_x]{\tb{\man{M}}}$\@.  We then have
\begin{equation*}
\tensor[k]{(\ctb[x]{\man{M}})_1\oplus(\ctb[x]{\man{M}})_2}\simeq
\bigoplus_{a_1,\dots,a_k\in\{1,2\}}(\ctb[x]{\man{M}})_{a_1}\otimes\dots\otimes
(\ctb[x]{\man{M}})_{a_k}
\end{equation*}
by~\cite[\S{}III.5.5]{NB:89a}\@.  Let
\begin{equation*}
\map{\pi_k}{\tensor[k]{\ctb[v_x]{\tb{\man{M}}}}}
{(\ctb[x]{\man{M}})_1\otimes\dots\otimes(\ctb[x]{\man{M}})_1}
\end{equation*}
be the projection onto the component of the direct sum decomposition.

With all of the preceding, we can now make sense of the following sublemma.
We adopt the notation~\eqref{eq:nabla(j)} introduced in the proof of
Theorem~\ref{the:COinfty-weak}\@.
\begin{proofsublemma}
If, for\/ $X,Y\in\sections[\nu]{\tb{\man{M}}}$\@, we have\/
$Z=\tlift{X}+\vlift{Y}$\@, then we have
\begin{equation*}
\pi_k\otimes\id_{\tb{\tb{\man{M}}}}((\tlift{\nabla})^{(k)}Z(0_x))=
\nabla^{(k)}X(x)\oplus(\nabla^{(k)}Y(x))
\end{equation*}
for\/ $k\in\integernn$ satisfying\/ $k\le\nu$\@.
\begin{subproof}
Obviously we can consider two special cases, the first where $Y=0$ and the
second where $X=0$\@. When $Y=0$\@, the result follows from
Sublemma~\ref{psublem:nablaT} from the proof of Lemma~\ref{lem:XTflow}\@,
especially the formula~\eqref{eq:nablaThor-expand} we derived from the
sublemma.  When $X=0$ the result immediately follows from the same sublemma.
\end{subproof}
\end{proofsublemma}

By the preceding sublemma, $Z(t,0_x)=X(t,x)\oplus Y(t,x)$\@.  Since the
projections onto the first and second component of the direct sum
decomposition of $\tb{\tb{\man{M}}}$ are continuous, we immediately conclude
that $X,Y\in\CFsections[\nu]{\tdomain;\tb{(\tbproj{\man{M}}(\nbhd{W}))}}$\@.

The remainder of the proof breaks into the various cases of regularity.

$\nu=\infty$\@: Let $K\subset\man{M}$ be compact and let $m\in\integernn$\@.
Since $K$ is also a compact subset of $\tb{\man{M}}$\@, there exists
$g\in\Lloc^1(\tdomain;\realnn)$ such that
\begin{equation*}
\dnorm{j_mZ(t,0_x)}_{\ol{\tlift{\metric}}_m}\le g(t),\qquad
t\in\tdomain,\ x\in K.
\end{equation*}
Let $\map{\ol{\pi}_m}{\jet{m}{\tb{\tb{\man{M}}}}}
{\oplus_{j=0}^m\tensor[j]{\tbproj{\man{M}}^*\tb{\man{M}}}
\otimes\tb{\tb{\man{M}}}}$ be defined by
\begin{equation*}
\ol{\pi}_m(j_mZ'(v_x))=\sum_{j=0}^m\pi_j\otimes\id_{\tb{\tb{\man{M}}}}
((\tlift{\nabla})^{(j-1)}Z'(v_x)),
\end{equation*}
this making sense by virtue of Lemma~\ref{lem:Jrdecomp}\@.  By the sublemma,
by the definition of $\tlift{\metric}$\@, and by the definition of the fibre
metrics on $\jet{m}{\tb{\man{M}}}$ and $\jet{m}{\tb{\tb{\man{M}}}}$ induced
by the decomposition of Lemma~\ref{lem:Jrdecomp}\@, we have
\begin{equation*}
\dnorm{\ol{\pi}_m(j_mZ(t,0_x))}^2_{\ol{\tlift{\metric}}_m}=
\dnorm{j_mX(t,x)}^2_{\ol{\metric}_m}+\dnorm{j_mY(t,x)}^2_{\ol{\metric}_m}.
\end{equation*}
This gives
\begin{equation*}
\dnorm{j_mX(t,x)}_{\ol{\metric}_m}\le g(t),\enspace
\dnorm{j_mY(t,x)}_{\ol{\metric}_m}\le g(t),\qquad t\in\tdomain,\ x\in K,
\end{equation*}
which gives the lemma in this case.

$\nu=m$\@: From the computations above in the smooth case we have that $X$
and $Y$ are locally integrally $\C^{m-1}$-bounded.  To show $X$ is, in fact,
locally integrally $\C^m$-bounded, we will use the computations from the
proof of Lemma~\ref{lem:XTflow}\@.  Let $K\subset\man{M}$ and let
\begin{equation*}
K_1=\setdef{v_x\in\tb{\man{M}}}{x\in K,\ \dnorm{v_x}_{\metric}\le1}
\end{equation*}
so $K_1$ is a compact subset of $\tb{\man{M}}$\@.  For the moment, let us fix
$t\in\tdomain$\@.  We now recall equation~\eqref{eq:PknablaT1} which gives a
formula for $P^m_{\tlift{\nabla}}(\tlift{X_t})$ when all arguments are
horizontal.  Since, in the expression~\eqref{eq:PknablaT1}\@, $v_x$ is
arbitrary, by letting it vary over vectors of unit length we get an estimate
\begin{equation*}
\dnorm{P^m_\nabla(X_t)(x)}_{\metric_m}\le
C(p^{m-1}_K(X_t)+p^{m-1}_{K_1}(\tlift{X_t}))
\end{equation*}
for some $C\in\realp$\@.  Since $X,Y\in\LIsections[m-1]{\tdomain;\man{M}}$
and since
$\tlift{X}=Z-\vlift{Y}\in\LIsections[m-1]{\tdomain;\tb{\man{M}}}$\@, by
Lemma~\ref{lem:Jrdecomp} there exists $g\in\Lloc^1(\tdomain;\realnn)$ such
that
\begin{equation*}
\dnorm{j_mX_t(x)}\le g(t),\qquad (t,x)\in\tdomain\times K,
\end{equation*}
which gives $X\in\LIsections[m]{\tb{\man{M}}}$\@.  

$\nu=m+\lip$\@: This follows from the computations above, using
Lemma~\ref{lem:dil-deriv}\@,~\cf~the proof of the Lipschitz part of the proof
of Lemma~\ref{lem:XTflow}\@.

$\nu=\omega$\@: Let $K\subset\man{M}$ be compact and let\/
$\vect{a}\in\csd(\integernn;\realp)$\@.  Since $K$ is also a compact subset
of $\tb{\man{M}}$\@, there exists $g\in\Lloc^1(\tdomain;\realnn)$ such that
\begin{equation*}
a_0a_1\cdots a_m\dnorm{j_mZ(t,0_x)}_{\ol{\tlift{\metric}}_m}
\le g(t),\qquad t\in\tdomain,\ x\in K,\ m\in\integernn.
\end{equation*}
As in the smooth case we have
\begin{equation*}
\dnorm{\ol{\pi}_m(j_mZ(t,0_x))}^2_{\ol{\tlift{\metric}}_m}=
\dnorm{j_mX(t,x)}^2_{\ol{\metric}_m}+\dnorm{j_mY(t,x)}^2_{\ol{\metric}_m}.
\end{equation*}
This gives
\begin{equation*}
a_0a_1\cdots a_m\dnorm{j_mX(t,x)}_{\ol{\metric}_m}\le g(t),\quad
a_0a_1\cdots a_m\dnorm{j_mY(t,x)}_{\ol{\metric}_m}\le g(t),
\end{equation*}
for $t\in\tdomain$\@, $x\in K$\@, and $m\in\integernn$\@, which gives the
lemma.
\end{proof}
\end{lemma}

Next let us think about open-loop subfamilies for linearisations.  Generally
speaking, one may wish to consider different classes of open-loop systems for
the ``tangent lift part'' and the ``vertical lift part'' of a linearised
system.  The open-loop systems for the tangent lift part will be those giving
rise to reference trajectories and reference flows.  On the other hand, the
open-loop systems for the vertical lift part will be those that we will allow
as perturbing the reference flow.  There is no reason that these should be the
same.  While this proliferation of open-loop subfamilies will lead to some
notational complexity, the freedom to carefully account for these
possibilities is one of the strengths of our theory.  Indeed, in standard
Jacobian linearisation, it is difficult to keep track of how the
controls\textemdash{}constraints on them and attributes of
them\textemdash{}are carried over to the linearisation.  In our theory, this
is natural.

We first make tangent and vertical lift constructions for open-loop subfamilies.
\begin{definition}
Let $m\in\integerp$ and $m'\in\{0,\lip\}$\@, let
$\nu\in\{m+m',\infty,\omega\}$\@, and let $r\in\{\infty,\omega\}$\@, as
required.  Let $\fG=(\man{M},\sF)$ be a $\C^\nu$-\gcs\ with linearisation
$\tf{\fG}=(\tb{\man{M}},\tf{\sF})$\@, and let $\sO_{\fG}$ be an open-loop
subfamily for $\fG$\@.
\begin{compactenum}[(i)]
\item The \defn{tangent lift} of $\sO_{\fG}$ is the open-loop subfamily
$\tlift{\sO}_{\fG}$ for\/ $(\tb{\man{M}},\tlift{\sF})$ defined by
\begin{equation*}
\tlift{\sO}_{\fG}(\tdomain,\nbhd{W})=\setdef{\tlift{X}|\nbhd{W}}
{X\in\sO_{\fG}(\tdomain,\tbproj{\man{M}}(\nbhd{W}))}
\end{equation*}
for a time-domain $\tdomain$ and for $\nbhd{W}\subset\tb{\man{M}}$ open.
\item The \defn{vertical lift} of $\sO_{\fG}$ is the open-loop subfamily
$\vlift{\sO}_{\fG}$ for\/ $(\tb{\man{M}},\vlift{\sF})$ defined by
\begin{equation*}
\vlift{\sO}_{\fG}(\tdomain,\nbhd{W})=\setdef{\vlift{Y}|\nbhd{W}}
{Y\in\sO_{\fG}(\tdomain,\tbproj{\man{M}}(\nbhd{W}))}
\end{equation*}
for a time-domain $\tdomain$ and for $\nbhd{W}\subset\tb{\man{M}}$ open.\oprocend
\end{compactenum}
\end{definition}

\begin{definition}
Let $m\in\integerp$ and $m'\in\{0,\lip\}$\@, let
$\nu\in\{m+m',\infty,\omega\}$\@, and let $r\in\{\infty,\omega\}$\@, as
required.  Let $\fG=(\man{M},\sF)$ be a $\C^\nu$-\gcs\ with linearisation
$\tf{\fG}=(\tb{\man{M}},\tf{\sF})$\@.  An \defn{open-loop subfamily} for
$\tf{\fG}$ defined by a pair $(\sO_{\fG,0},\sO_{\fG,1})$ of open-loop
subfamilies for $\fG$ is the open-loop subfamily
$\tlift{\sO}_{\fG,0}+\vlift{\sO}_{\fG,1}$ defined by:
\begin{equation*}\eqoprocend
\tlift{X}+\vlift{Y}\in(\tlift{\sO}_{\fG,0}+
\vlift{\sO}_{\fG,1})(\tdomain,\nbhd{W})\enspace\iff\enspace \tlift{X}\in\tlift{\sO}_{\fG,0}(\tdomain,\tbproj{\man{M}}(\nbhd{W})),\
\vlift{Y}\in\vlift{\sO}_{\fG,1}(\tdomain,\tbproj{\man{M}}(\nbhd{W})).
\end{equation*}
\end{definition}

Note that the restriction properties of open-loop subfamilies as per
Definition~\ref{def:olsubfamily} are satisfied by our construction above, so
the result is indeed an open-loop subfamily for $\tf{\fG}$\@.

Next we can define what we mean by trajectories for the linearisation in the
more or less obvious way.
\begin{definition}
Let $m\in\integerp$ and $m'\in\{0,\lip\}$\@, let
$\nu\in\{m+m',\infty,\omega\}$\@, and let $r\in\{\infty,\omega\}$\@, as
required.  Let $\fG=(\man{M},\sF)$ be a $\C^\nu$-\gcs\ with linearisation
$\tf{\fG}=(\tb{\man{M}},\tf{\sF})$\@.  Let $\sO_{\fG,0}$ and $\sO_{\fG,1}$ be
open-loop subfamilies for $\fG$\@.
\begin{compactenum}[(i)]
\item For a time-domain $\tdomain$\@, an open set
$\nbhd{W}\subset\tb{\man{M}}$\@, and for $X\in\sO_{\fG,0}(\tdomain,\nbhd{U})$
and $Y\in\sO_{\fG,1}$\@, an \defn{$(X,Y,\tdomain,\nbhd{W})$-trajectory} for
$(\sO_{\fG,0},\sO_{\fG,1})$ is a curve $\map{\Upsilon}{\tdomain}{\nbhd{W}}$
such that $\Upsilon'(t)=\tlift{X}(t,\Upsilon(t))+\vlift{Y}(t,\Upsilon(t))$\@.
\item For a time-domain $\tdomain$ and an open set
$\nbhd{W}\subset\tb{\man{M}}$\@, a \defn{$(\tdomain,\nbhd{W})$-trajectory}
for the pair $(\sO_{\fG,0},\sO_{\fG,1})$ is a
$(\tdomain,\nbhd{W})$-trajectory for
$\tlift{\sO}_{\fG,0}+\vlift{\sO}_{\fG,1}$\@.
\item A plain \defn{trajectory} for the pair $(\sO_{\fG,0},\sO_{\fG,1})$ is a
curve that is a $(\tdomain,\nbhd{W})$-trajectory for
$(\sO_{\fG,0},\sO_{\fG,1})$ for some time-domain $\tdomain$ and some open
$\nbhd{W}\subset\tb{\man{M}}$\@.\savenum
\end{compactenum}
We denote by:
\begin{compactenum}[(i)]\resumenum
\item $\Traj(X,Y,\tdomain;\nbhd{W})$ the set of
$(X,Y,\tdomain,\nbhd{U})$-trajectories for $(\sO_{\fG,0},\sO_{\fG,1})$\@;
\item $\Traj(\tdomain,\nbhd{W},(\sO_{\fG,0},\sO_{\fG,1}))$ the set of
$(\tdomain,\nbhd{U})$-trajectories for $(\sO_{\fG,0},\sO_{\fG,1})$\@;
\item $\Traj(\sO_{\fG,0},\sO_{\fG,1})$ the set of trajectories for
$(\sO_{\fG,0},\sO_{\fG,1})$\@.
\end{compactenum}
We shall abbreviate
\begin{equation*}
\Traj(\tdomain,\nbhd{W},(\subscr{\sO}{$\fG$,full},\subscr{\sO}{$\fG$,full}))=
\Traj(\tdomain,\nbhd{W},\tf{\fG})
\end{equation*}
and $\Traj(\subscr{\sO}{$\fG$,full},\subscr{\sO}{$\fG$,full})=
\Traj(\tf{\fG})$\@.\oprocend
\end{definition}

Now that we have been clear about what we mean by the trajectory of a
linearised system, let us say some things about these trajectories.
\begin{proposition}
Let\/ $m\in\integerp$ and\/ $m'\in\{0,\lip\}$\@, let\/
$\nu\in\{m+m',\infty,\omega\}$\@, and let\/ $r\in\{\infty,\omega\}$\@, as
required.  Let\/ $\fG=(\man{M},\sF)$ be a $\C^\nu$-\gcs\ with linearisation\/
$\tf{\fG}$\@, and let\/ $\sO_{\fG,0}$ and\/ $\sO_{\fG,1}$ be open-loop
subfamilies for\/ $\fG$\@.  Let\/ $\tdomain\subset\real$ be a time-domain and
let\/ $\nbhd{W}\subset\tb{\man{M}}$ be open.  If\/
$\tlift{\xi}\in\Traj(\tdomain,\nbhd{W},(\sO_{\fG,0},\sO_{\fG,1}))$ then the
following statements hold:
\begin{compactenum}[(i)]
\item \label{pl:lintraj1} there exist\/ $X\in\sO_{\fG,0}(\tdomain,\tbproj{\man{M}}(\nbhd{W}))$ and\/
$Y\in\sO_{\fG,1}(\tdomain,\tbproj{\man{M}}(\nbhd{W}))$ such that
\begin{equation*}
(\tlift{\xi})'(t)=\tlift{X}(t,\xi^T(t))+\vlift{Y}(t,\xi^T(t));
\end{equation*}
\item \label{pl:lintraj2} there exists\/
$\xi\in\Traj(\tdomain,\tbproj{\man{M}}(\nbhd{W}),\sO_{\fG,0})$ such that the
diagram
\begin{equation*}
\xymatrix{{\tdomain}\ar[r]^{\tlift{\xi}}\ar[rd]_{\xi}&
{\tb{\man{M}}}\ar[d]^{\tbproj{\nbhd{U}}}\\&{\man{M}}}
\end{equation*}
commutes,~\ie~$\tlift{\xi}$ is a vector field along\/ $\xi$\@.
\end{compactenum}
\begin{proof}
The first assertion follows from Lemma~\ref{lem:open-loop-lin}\@.  The second
assertion follows by taking $\xi=\tbproj{\man{M}}\scirc\tlift{\xi}$\@, and
noting that
\begin{equation*}
\xi'(t)=\tf[\tlift{\xi}(t)]{\tbproj{\man{M}}}((\tlift{\xi})'(t))\\
=\tf[\tlift{\xi}(t)]{\tbproj{\man{M}}}(\tlift{X}(t,\tlift{\xi}(t))+
\vlift{Y}(t,\xi^T(t)))=X(t,\xi(t))
\end{equation*}
and $X$ is an open-loop system for $\sO_{\fG,0}$\@.
\end{proof}
\end{proposition}

\subsection{Linearisation about reference trajectories and reference flows}\label{subsec:lin-traj-flow}

Let us now slowly begin to pull back our general notion of linearisation to
something more familiar.  In this section we will linearise about two sorts
of things, trajectories and flows.  We will see in the next section that it
is the distinction between these two things that accounts for the problems
observed in Example~\ref{eg:bad-linearise}\@.

But for now, we proceed in general.  We let $\fG$ be a \gcs\ and $\sO_{\fG}$
an open-loop subfamily.  We recall from Example~\ref{eg:traj-ol} that, if
$\tdomain$ is a time-domain, if $\nbhd{U}\subset\man{M}$ is open, and if
$\xi\in\Traj(\tdomain,\nbhd{U},\sO_{\fG})$\@, then $\sO_{\fG,\xi}$ is the
open-loop subfamily associated to the trajectory $\xi$\@,~\ie~all open-loop
systems from $\sO_{\fG}$ possessing $\xi$ as a trajectory.  Having made this
recollection, we make the following definition.
\begin{definition}
Let $m\in\integerp$ and $m'\in\{0,\lip\}$\@, let
$\nu\in\{m+m',\infty,\omega\}$\@, and let $r\in\{\infty,\omega\}$\@, as
required.  Let $\fG=(\man{M},\sF)$ be a $\C^\nu$-\gcs\ with linearisation
$\tf{\fG}$\@.  Let $\sO_{\fG,0}$ and $\sO_{\fG,1}$ be open-loop subfamilies
for $\fG$\@, let $\tdomain$ be a time-domain, let $\nbhd{U}\subset\man{M}$ be
open, and let $\subscr{\xi}{ref}\in\Traj(\tdomain,\nbhd{U},\sO_{\fG,0})$\@.
The \defn{$(\sO_{\fG,0},\sO_{\fG,1})$-linearisation of $\fG$ about
$\subscr{\xi}{ref}$} is the open-loop subfamily
$\tlift{\sO}_{\fG,0,\subscr{\xi}{ref}}+\vlift{\sO}_{\fG,1}$ for $\tf{\fG}$\@.
A \defn{trajectory} for this linearisation is a
$(\tdomain',\nbhd{W})$-trajectory $\Upsilon$ for
$(\sO_{\fG,0,\subscr{\xi}{ref}},\sO_{\fG,1})$ satisfying
$\tbproj{\man{M}}\scirc\Upsilon=\subscr{\xi}{ref}$\@, and where
$\tdomain'\subset\tdomain$ and
$\nbhd{W}\subset\tbproj{\man{M}}^{-1}(\nbhd{U})$\@.\oprocend
\end{definition}

By definition, a trajectory for the linearisation about the reference
trajectory $\subscr{\xi}{ref}$ is a curve
$\map{\Upsilon}{\tdomain'}{\nbhd{W}}$ satisfying
\begin{equation*}
\Upsilon'(t)=\tlift{X}(t,\Upsilon(t))+\vlift{Y}(t,\Upsilon(t)),
\end{equation*}
for $X\in\sO_{\fG,0,\subscr{\xi}{ref}}(\tdomain',\tbproj{\man{M}}(\nbhd{W}))$
and for $Y\in\sO_{\fG,1}(\tdomain',\tbproj{\man{M}}(\nbhd{W}))$\@, and where
$\Upsilon$ is a tangent vector field along $\subscr{\xi}{ref}$\@.  Note that
there may well be trajectories for
$(\sO_{\fG,0,\subscr{\xi}{ref}},\sO_{\fG,1})$ that are \emph{not} vector
fields along $\subscr{\xi}{ref}$\@; we just do not call these trajectories
for the linearisation \emph{about} $\subscr{\xi}{ref}$\@.

Let us now talk about linearisation, not about a trajectory, but about a
flow.  Here we recall the notion of the open-loop subfamily associated to an
open-loop system in Example~\enumdblref{eg:olsf}{enum:X-ol}\@.
\begin{definition}
Let $m\in\integerp$ and $m'\in\{0,\lip\}$\@, let
$\nu\in\{m+m',\infty,\omega\}$\@, and let $r\in\{\infty,\omega\}$\@, as
required.  Let $\fG=(\man{M},\sF)$ be a $\C^\nu$-\gcs\ with linearisation
$\tf{\fG}$\@.  Let $\sO_{\fG,0}$ and $\sO_{\fG,1}$ be open-loop subfamilies
for $\fG$\@, let $\tdomain$ be a time-domain, let $\nbhd{U}\subset\man{M}$ be
open, and let $\subscr{X}{ref}\in\sO_{\fG,0}(\tdomain,\nbhd{U})$\@.  The
\defn{$\sO_{\fG,1}$-linearisation of $\fG$ about $\subscr{X}{ref}$} is the
open-loop subfamily $\tlift{\sO}_{\fG,0,\subscr{X}{ref}}+\vlift{\sO}_{\fG,1}$
for $\tf{\fG}$\@.  A \defn{trajectory} for this linearisation is a
$(\tdomain',\nbhd{W})$-trajectory for
$(\sO_{\fG,0,\subscr{X}{ref}},\sO_{\fG,1})$\@, where
$\tdomain'\subset\tdomain$ and where
$\nbhd{W}\subset\tbproj{\man{M}}^{-1}(\nbhd{U})$\@.\oprocend
\end{definition}

By definition, a trajectory for the linearisation about the reference flow
$\subscr{X}{ref}$ is a curve $\map{\Upsilon}{\tdomain'}{\nbhd{W}}$ satisfying
\begin{equation*}
\Upsilon'(t)=\tlift{\subscr{X}{ref}}(t,\Upsilon(t))+\vlift{Y}(t,\Upsilon(t)),
\end{equation*}
for $Y\in\sO_{\fG,1}(\tdomain',\tbproj{\man{M}}(\nbhd{W}))$\@.  Note that the
definition of $\sO_{\fG,0,\subscr{X}{ref}}$ necessarily implies that
$\tbproj{\man{M}}\scirc\Upsilon$ is an integral curve for
$\subscr{X}{ref}$\@.  Unlike the case of linearisation about a reference
trajectory, we do not specify that the trajectories for the linearisation
about a reference flow follow a specific trajectory for $\fG$\@, although one
can certainly do this as well.

\subsection{Linearisation about an equilibrium point}\label{subsec:equil-linearise}

Continuing to make things concrete, let us consider linearising about trivial
reference trajectories and reference flows.  We begin by considering what an
equilibrium point is in our framework.
\begin{definition}
Let $m\in\integerp$ and $m'\in\{0,\lip\}$\@, let
$\nu\in\{m+m',\infty,\omega\}$\@, and let $r\in\{\infty,\omega\}$\@, as
required.  Let $\fG=(\man{M},\sF)$ be a $\C^\nu$-\gcs\ and let
$x_0\in\man{M}$\@.
\begin{compactenum}[(i)]
\item The \defn{\gcs\ for $\fG$ at $x_0$} is the
$\C^\nu$-\gcs\ $\fG_{x_0}=(\man{M},\Eq_{\sF,x_0})$\@, where
\begin{equation*}
\Eq_{\sF,x_0}(\nbhd{U})=\setdef{X\in\sF(\nbhd{U})}{X(x_0)=0_{x_0}}.
\end{equation*}
\item If there exists an open set $\nbhd{U}\subset\man{M}$ for which
$\Eq_{\sF,x_0}(\nbhd{U})\not=\emptyset$\@, then $x_0$ is an
\defn{equilibrium point} for $\fG$\@.\oprocend
\end{compactenum}
\end{definition}

Of course, by properties of presheaves, if $X\in\Eq_{\sF,x_0}(\nbhd{U})$\@,
then $X|\nbhd{V}\in\Eq_{\sF,x_0}(\nbhd{V})$ for every open set
$\nbhd{V}\subset\nbhd{U}$\@.  Thus $\fG_{x_0}$ is indeed a \gcs.

Let us examine the nature of \gcs{}s at $x_0$\@.  This amounts to
understanding any particular structure that one can associate to vector
fields that vanish at a point.  This is the content of the following lemma.
\begin{lemma}
Let\/ $\man{M}$ be a smooth manifold, let\/ $x_0\in\man{M}$\@, and let\/
$X\in\sections[1]{\man{M}}$\@.  If\/ $X(x_0)=0_{x_0}$\@, then there exists a
unique\/ $A_{X,x_0}\in\End_{\real}(\tb[x_0]{\man{M}})$ satisfying either of
the following equivalent characterisations:
\begin{compactenum}[(i)]
\item \label{pl:AXx01} noting that\/ $\map{\tlift{X}|\tb[x_0]{\man{M}}}{\tb[x_0]{\man{M}}}
{\vb[0_{x_0}]{\tb{\man{M}}}\simeq\tb[x_0]{\man{M}}}$\@, $A_{X,x_0}=\tlift{X}|\tb[x_0]{\man{M}}$\@;
\item \label{pl:AXx02} $A_{X,x_0}(v_{x_0})=[V,X](x_0)$ where\/
$V\in\sections[\infty]{\man{M}}$ satisfies\/ $V(x_0)=v_{x_0}$\@.
\end{compactenum}
\begin{proof}
We will show that the characterisation from part~\eqref{pl:AXx01} makes sense,
and that it agrees with the second characterisation.

First, note that, since $X(x_0)=0_{x_0}$\@,
$\tf[v_{x_0}]{\tbproj{\man{M}}}(\tlift{X}(v_{x_0}))=0_{x_0}$ for every
$v_{x_0}\in\tb[x_0]{\man{M}}$\@.  Thus
$\tlift{X}(v_{x_0})\in\vb[0_{x_0}]{\tb{\man{M}}}$\@, as claimed.  That
$\tlift{X}|\tb[x_0]{\man{M}}$ is linear is a consequence of the fact that
$\tlift{X}$ is a linear vector field,~\ie~that the diagram~\eqref{eq:XTlinear}
commutes.  In the particular case that $X(x_0)=0_{x_0}$\@, the diagram
implies that $\tlift{X}$ is a linear map from $\tb[x_0]{\man{M}}$ to
$\tb[0_{x_0}]{\tb{\man{M}}}$\@.  As we already know that
$\tlift{X}|\tb[x_0]{\man{M}}$ is $\vb[0_{x_0}]{\tb{\man{M}}}$-valued, the
characterisation from part~\eqref{pl:AXx01} does indeed uniquely define an
endomorphism of $\tb[x_0]{\man{M}}$\@.

Let us now show that the characterisation of part~\eqref{pl:AXx02} agrees
with that of part~\eqref{pl:AXx01}\@.
By~\cite[Theorem~4.2.19]{RA/JEM/TSR:88}\@, we have
\begin{align*}
\vlft(0_{x_0},[V,X](x_0))=&\;
\derivatzero{}{t}\tf[\flow{X}{-t}(x_0)]{\flow{X}{t}}\scirc
V\scirc\flow{X}{-t}(x_0)\\
=&\;\derivatzero{}{t}\tf[x_0]{\flow{X}{t}}\scirc V(x_0)=\tlift{X}(V(x_0)),
\end{align*}
as desired.
\end{proof}
\end{lemma}

According to the lemma, we can make the following definitions.
\begin{definition}
Let $m\in\integerp$ and $m'\in\{0,\lip\}$\@, let
$\nu\in\{m+m',\infty,\omega\}$\@, and let $r\in\{\infty,\omega\}$\@, as
required.  Let $\fG=(\man{M},\sF)$ be a $\C^\nu$-\gcs.  For an equilibrium
point $x_0\in\man{M}$ for $\fG$\@, we define
\begin{equation*}
\sL_{\sF,x_0}=\setdef{A_{X,x_0}}{[X]_{x_0}\in(\Eq_{\sF,x_0})_{x_0}}
\end{equation*}
(where $(\Eq_{\sF,x_0})_{x_0}$ denotes the stalk of the presheaf
$\Eq_{\sF,x_0}$ at $x_0$) and
\begin{equation*}\eqoprocend
\sF(x_0)=\setdef{X(x_0)}{[X]_{x_0}\in\sF_{x_0}}.
\end{equation*}
\end{definition}

Associated to an equilibrium point are natural notions of open-loop systems
that preserve the equilibrium point.
\begin{definition}
Let $m\in\integerp$ and $m'\in\{0,\lip\}$\@, let
$\nu\in\{m+m',\infty,\omega\}$\@, and let $r\in\{\infty,\omega\}$\@, as
required.  Let $\fG=(\man{M},\sF)$ be a $\C^\nu$-\gcs.  If $x_0\in\man{M}$
and if $\sO_{\fG}$ is an open-loop subfamily for $\fG$\@, the open-loop
subfamily $\sO_{\fG,x_0}$ is defined by specifying that, for a time-domain
$\tdomain$ and an open set $\nbhd{U}\subset\man{M}$\@,
\begin{equation*}\eqoprocend
\sO_{\fG,x_0}(\tdomain,\nbhd{U})=
\setdef{X\in\sO_{\fG}(\tdomain;\nbhd{U})}{X(t)\in\Eq_{\sF,x_0}(\nbhd{U}),\
t\in\tdomain}.
\end{equation*}
\end{definition}

Note that the only trajectory of $\sO_{\fG,x_0}$ passing through $x_0$ is the
constant trajectory $t\mapsto x_0$\@, as it should be.

It is now more or less obvious how one should define linearisations about an
equilibrium point.  This can be done for trajectories and flows.  We start
with trajectories.
\begin{definition}
Let $m\in\integerp$ and $m'\in\{0,\lip\}$\@, let
$\nu\in\{m+m',\infty,\omega\}$\@, and let $r\in\{\infty,\omega\}$\@, as
required.  Let $\fG=(\man{M},\sF)$ be a $\C^\nu$-\gcs\ with linearisation
$\tf{\fG}$\@.  Let $\sO_{\fG,0}$ and $\sO_{\fG,1}$ be open-loop subfamilies
for $\fG$ and let $x_0\in\man{M}$\@.  The
\defn{$(\sO_{\fG,0},\sO_{\fG,1})$-linearisation of $\fG$ about $x_0$} is the
open-loop subfamily $\tlift{\sO}_{\fG,0,x_0}+\vlift{\sO}_{\fG,1}$ for
$\tf{\fG}$\@.  A \defn{trajectory} for this linearisation is a
$(\tdomain,\nbhd{W})$-trajectory for the
$(\sO_{\fG,0,x_0},\sO_{\fG,1})$-linearisation about the trivial reference
trajectory $t\mapsto x_0$\@, where $\tdomain$ is a time-domain and where
$\nbhd{W}$ is a neighbourhood of $\tb[x_0]{\man{M}}$\@.\oprocend
\end{definition}

By definition and by the characterisation of $\tlift{X}$ at equilibria, a
trajectory for the linearisation about $x_0$ will be a curve
$\map{\Upsilon}{\tdomain}{\tb[x_0]{\man{M}}}$ satisfying
\begin{equation*}
\Upsilon'(t)=A_{X(t),x_0}(\Upsilon(t))+b(t),
\end{equation*}
where $t\mapsto X(t)$ is a curve in $\sL_{\sF,x_0}$ whose nature is
determined by the open-loop subfamily $\sO_{\fG,0}$\@,~\eg~it may be locally
integrable, locally essentially bounded, piecewise constant, \etc, and where
$t\mapsto b(t)$ is a curve in $\sF(x_0)\subset\tb[x_0]{\man{M}}$\@, again
whose nature is determined by the open-loop subfamily $\sO_{\fG,1}$\@.  Note
that the linearisation about $x_0$ will, therefore, generally be a family of
\emph{time-dependent} linear systems on $\tb[x_0]{\man{M}}$\@.  This may come
as a surprise to those used to Jacobian linearisation, but we will see in
Example~\ref{eg:good-linearise} below how this arises in practice.

Let us now talk about linearisation about an equilibrium point, not about a 
trajectory, but about a flow.
\begin{definition}
Let $m\in\integerp$ and $m'\in\{0,\lip\}$\@, let
$\nu\in\{m+m',\infty,\omega\}$\@, and let $r\in\{\infty,\omega\}$\@, as
required.  Let $\fG=(\man{M},\sF)$ be a $\C^\nu$-\gcs\ with linearisation
$\tf{\fG}$\@.  Let $\sO_{\fG,0}$ and $\sO_{\fG,1}$ be open-loop subfamilies
for $\fG$\@, let $\tdomain$ be a time-domain, let $x_0\in\man{M}$\@, let
$\nbhd{U}\subset\man{M}$ be a neighbourhood of $x_0$\@, and let
$\subscr{X}{ref}\in\sO_{\fG,0,x_0}(\tdomain,\nbhd{U})$\@.  The
\defn{$\sO_{\fG,1}$-linearisation of $\fG$ about $(\subscr{X}{ref},x_0)$} is
the open-loop subfamily
$\tlift{\sO}_{\fG,0,\subscr{X}{ref}}+\vlift{\sO}_{\fG,1}$ for $\tf{\fG}$\@.
A \defn{trajectory} for this linearisation is a
$(\tdomain',\nbhd{W})$-trajectory for
$(\sO_{\fG,0,\subscr{X}{ref}},\sO_{\fG,1})$\@, where
$\tdomain'\subset\tdomain$ and where
$\nbhd{W}\subset\tbproj{\man{M}}^{-1}(\nbhd{U})$\@.\oprocend
\end{definition}

In this case, we have a prescribed curve $t\mapsto\subscr{X}{ref}(t)$ such
that $\subscr{X}{ref}(t,x_0)=0_{x_0}$ for every $t$\@.  Thus this defines a
curve $A_{\subscr{X}{ref}(t),x_0}$ in $\sL_{\sF,x_0}$\@.  By definition, a
trajectory for the linearisation about the pair $(\subscr{X}{ref},x_0)$ is a
curve $\map{\Upsilon}{\tdomain'}{\tb[x_0]{\man{M}}}$ satisfying
\begin{equation*}
\Upsilon'(t)=A_{\subscr{X}{ref}(t),x_0}(\Upsilon(t))+b(t),
\end{equation*}
where $t\mapsto b(t)$ is a curve in $\sF(x_0)$ having properties determined
by the open-loop subfamily $\sO_{\fG,1}$\@.  Note that this linearisation
will still generally be time-dependent, but it is now a \emph{single}
time-dependent linear system, not a family of them, as with linearisation
about a trajectory.  Moreover, if $\subscr{X}{ref}$ is chosen to be
time-independent, then the linearisation will also be time-invariant.  But
there is no reason in the general theory to do this.

The above comments about the possibility of time-varying linearisations
notwithstanding, there is one special case where we can be sure that
linearisations will be time-independent, and this is when $\sL_{\sF,x_0}$
consists of a single vector field.  The following result gives a common case
where this happens.  Indeed, the ubiquity of this situation perhaps explains
the neglect of the general situation that has led to the seeming
contradictions in the standard treatments, such as are seen in Example~\ref{eg:bad-linearise}\@.
\begin{proposition}\label{prop:nice-cas}
Let\/ $\Sigma=(\man{M},F,\cs{C})$ be a\/ $\C^1$-control-affine system with\/
$\cs{C}\subset\real^k$ and
\begin{equation*}
F(x,\vect{u})=f_0(x)+\sum_{a=1}^ku^af_a(x).
\end{equation*}
For\/ $x_0\in\man{M}$\@, suppose that
\begin{compactenum}[(i)]
\item there exists\/ $\vect{u}_0\in\cs{C}$ such that
\begin{equation*}
f_0(x_0)=\sum_{a=1}^ku^a_0f_a(x_0)
\end{equation*}
and
\item $\ifam{f_1(x_0),\dots,f_k(x_0)}$ is linearly independent.
\end{compactenum}
Then\/ $x_0$ is an equilibrium point for\/ $\fG_\Sigma$ and\/
$\sL_{\sF_\Sigma,x_0}$ consists of a single linear map.
\begin{proof}
Let us define
\begin{equation*}
f'_0=f_0-\sum_{a=1}^ku^a_0f_a,
\end{equation*}
noting that $f'_0\in\sF_\Sigma$\@.  Since $f'_0(x_0)=0_{x_0}$\@, we conclude
that $x_0$ is an equilibrium point.  Now suppose that
$F(x_0,\vect{u})=0_{x_0}$\@.  Thus
\begin{equation*}
f_0(x_0)+\sum_{a=1}^ku^af_a(x_0)=0_{x_0}\enspace\implies\enspace
f_0(x_0)=-\sum_{a=1}^ku^af_a(x_0).
\end{equation*}
This last equation has a solution for $\vect{u}$\@, namely
$\vect{u}=-\vect{u}_0$\@, and since $\ifam{f_1(x_0),\dots,f_m(x_0)}$ is
linearly independent, this solution is unique.  Thus, for any neighbourhood
$\nbhd{U}$ of $x_0$\@,
\begin{equation*}
\Eq_{\sF_\Sigma,x_0}(\nbhd{U})=\Bigl\{f_0-\sum_{a=1}^ku^a_0f_a(x_0)\Bigr\}=
\{f'_0(x_0)\}.
\end{equation*}
This shows that $\sL_{\sF_\Sigma,x_0}=\{A_{f'_0,x_0}\}$\@, as claimed.
\end{proof}
\end{proposition}

While we are definitely not giving a comprehensive account of controllability
in this paper\textemdash{}see Section~\ref{subsec:future-cont} for a
discussion of controllability\textemdash{}in order to ``close the loop'' on
Example~\ref{eg:bad-linearise}\@, let us consider here how one talks about
linear controllability in our framework.  First we introduce some general
notation.
\begin{definition}
Let $\alg{F}$ be a field, let $\alg{V}$ be an $\alg{F}$-vector space, let
$\sL\subset\End_{\alg{F}}(\alg{V})$\@, and let $S\subset\alg{V}$\@.  By
$\inprod{\sL}{S}$ we denote the smallest subspace of\/ $\alg{V}$
that~(i)~contains $S$ and~(ii)~is invariant under $L$ for every
$L\in\sL$\@.\oprocend
\end{definition}

We can give a simple description of this subspace.
\begin{lemma}\label{lem:lincont-form}
If\/ $\alg{F}$ is a field, if\/ $\alg{V}$ is an $\alg{F}$-vector space, if\/
$\sL\subset\End_{\alg{F}}(\alg{V})$\@, and if\/ $S\subset\alg{V}$\@, then\/
$\inprod{\sL}{S}$ is spanned by elements of\/ $\alg{V}$ of the form
\begin{equation}\label{eq:ULS}
L_1\scirc\dots\scirc L_k(v),\qquad k\in\integernn,\ L_1,\dots,L_k\in\sL,\
v\in S.
\end{equation}
\begin{proof}
Let $\alg{U}_{\sL,S}$ be the subspace spanned by elements of the
form~\eqref{eq:ULS}\@.  Clearly $S\subset\alg{U}_{\sL,S}$ (taking the
convention that $L_1\scirc\dots\scirc L_k(v)=v$ if $k=0$) and, if
$L\in\sL$\@, then $L(\alg{U}_{\sL,S})\subset\alg{U}_{\sL,S}$ since an
endomorphism from $\sL$ maps a generator of the form~\eqref{eq:ULS} to another
generator of this form.  Therefore,
$\inprod{\sL}{S}\subset\alg{U}_{\sL,S}$\@.  Now, if $v\in S$\@, then clearly
$v\in\inprod{\sL}{S}$\@.  Since $\inprod{\sL}{S}$ is invariant under
endomorphisms from $\sL$\@, $L(v)\in\inprod{\sL}{S}$ for every $v\in S$ and
$L\in\sL$\@.  Recursively, we see that all generators of the
form~\eqref{eq:ULS} are in $\inprod{\sL}{S}$\@, whence
$\alg{U}_{\sL,S}\subset\inprod{\sL}{S}$ since $\alg{U}_{\sL,S}$ is a subspace.
\end{proof}
\end{lemma}

With the preceding as setup, let us make the following definition.
\begin{definition}\label{def:lincont}
Let $m\in\integerp$ and $m'\in\{0,\lip\}$\@, let
$\nu\in\{m+m',\infty,\omega\}$\@, and let $r\in\{\infty,\omega\}$\@, as
required.  Let $\fG=(\man{M},\sF)$ be a $\C^\nu$-\gcs\ with linearisation
$\tf{\fG}$\@, and let $x_0\in\man{M}$ be an equilibrium point for $\fG$\@.
The system $\fG$ is \defn{linearly controllable} at $x_0$ if there exists
$S\subset\sF(x_0)$ such that~(i)~$0_{x_0}\in\cohull(S)$
and~(ii)~$\inprod{\sL_{\sF,x_0}}{S}=\tb[x_0]{\man{M}}$\@.\oprocend
\end{definition}

\begin{remark}
For readers who may not recognise the relationship between our definition of
linear controllability and the classical Kalman rank
test~\cite[Theorem~13.3]{RWB:70}\@, we make the following comments.  Consider
the linear system
\begin{equation*}
\dot{\vect{x}}(t)=\mat{A}\vect{x}(t)+\mat{B}\vect{u}(t),
\end{equation*}
with $\vect{x}\in\real^n$\@, $\vect{u}\in\real^m$\@, and for appropriately
sized matrices $\mat{A}$ and $\mat{B}$\@.  Using Lemma~\ref{lem:lincont-form}
and the Cayley\textendash{}Hamilton Theorem, it is easy to check that the
smallest $\mat{A}$-invariant subspace containing $\image(\mat{B})$ is exactly
the columnspace of the Kalman controllability matrix,
\begin{equation*}
\left[\begin{array}{c|c|c|c}\mat{B}&\mat{A}\mat{B}&\cdots&
\mat{A}^{n-1}\mat{B}\end{array}\right].
\end{equation*}
For the more geometric approach to topics in linear system theory, we refer
to the excellent book of \citet{WMW:85}\@.\oprocend
\end{remark}

We state linear controllability as a definition, not a theorem, because we do
not want to develop the definitions required to state a theorem.  However, it
\emph{is} true that a system that is linearly controllable according to our
definition is small-time locally controllable in the usual sense of the word.
This is proved by \citet[Theorem~5.14]{COA:10}\@.  The setting of
\citeauthor{COA:10} is not exactly that of our paper.  However, it is easy to
see that this part of \citeauthor{COA:10}'s development easily translates to
what we are doing here.

Let us close this section, and the technical part of the paper, by revisiting
Example~\ref{eg:bad-linearise} where we saw that the classical picture of
Jacobian linearisation presents some problems.
\begin{example}\label{eg:good-linearise}
We work with the system
\begin{align*}
\dot{x}_1(t)=&\;x_2(t),\\
\dot{x}_2(t)=&\;x_3(t)u_1(t),\\
\dot{x}_3(t)=&\;u_2(t).
\end{align*}
We could as well work with the other representation for the system from
Example~\ref{eg:bad-linearise}\@, but since the family of vector fields is
the same (what changes between the two representations is the
parameterisation of the set of vector fields!), we will get the same
conclusions; this, after all, is the point of our feedback-invariant
approach.

This, of course, is a control-affine system, and the resulting \gcs\ is
$\fG=(\real^3,\sF)$ where $\sF$ is the globally generated presheaf with
\begin{equation*}
\sF(\real^3)=\setdef{f_0+u^1f_1+u^2f_2}{(u^1,u^2)\in\real^2},
\end{equation*}
with
\begin{equation*}
f_0=x_2\pderiv{}{x_1},\quad f_1=x_3\pderiv{}{x_2},\quad f_2=\pderiv{}{x_3}.
\end{equation*}
We have an equilibrium point at $(0,0,0)$\@.
\begin{prooflemma}
$\Eq_{\sF,(0,0,0)}(\real^3)=f_0+\vecspan[\real]{f_1}$\@.
\begin{subproof}
It is clear that $f_0(0,0,0)=f_1(0,0,0)=0$\@, and, therefore, any linear
combination of $f_0$ and $f_1$ will also vanish at $(0,0,0)$\@, and
particularly those from the affine subspace $f_0+\vecspan[\real]{f_1}$\@.
Conversely, if
\begin{equation*}
f_0(0,0,0)+u^1f_1(0,0,0)+u^2f_2(0,0,0)=0,
\end{equation*}
then $u^2=0$ and so the resulting vector field is in the asserted affine
subspace.
\end{subproof}
\end{prooflemma}

We, therefore, have
\begin{equation*}
\sL_{\sF,(0,0,0)}=
\asetdef{\begin{bmatrix}0&1&0\\0&0&a\\0&0&0\end{bmatrix}}{a\in\real}.
\end{equation*}
We also have
\begin{equation*}
\sF((0,0,0))=\setdef{bf_2(0,0,0)}{b\in\real}=
\asetdef{\begin{bmatrix}0\\0\\b\end{bmatrix}}{b\in\real}.
\end{equation*}
Thus a curve in $\sL_{\sF,(0,0,0)}$ has the form
\begin{equation*}
t\mapsto\begin{bmatrix}0&1&0\\0&0&a(t)\\0&0&0\end{bmatrix}
\end{equation*}
for a function $a$ having whatever properties might be induced from the
open-loop subfamily $\sO_{\fG,0}$ one is using,~\eg~locally integrable,
locally essentially bounded.  A curve in $\sF((0,0,0))$ has the form
\begin{equation*}
t\mapsto\begin{bmatrix}0\\0\\b(t)\end{bmatrix}
\end{equation*}
for a function $b$ having whatever properties might be induced
from the open-loop subfamily $\sO_{\fG,1}$ one is using.  Trajectories for
the linearisation about $(0,0,0)$ then satisfy
\begin{equation*}
\begin{bmatrix}\dot{v}_1(t)\\\dot{v}_2(t)\\\dot{v}_3(t)\end{bmatrix}=
\begin{bmatrix}0&1&0\\0&0&a(t)\\0&0&0\end{bmatrix}
\begin{bmatrix}v_1(t)\\v_2(t)\\v_3(t)\end{bmatrix}+
\begin{bmatrix}0\\0\\b(t)\end{bmatrix}.
\end{equation*}
Note that this is not a fixed time-varying linear system, but a family of
these, since the function $a$ is not \emph{a priori} specified, but is
variable.

Next let us look at two instances of linearisation about a reference flow by
choosing the two reference flows $X_1=f_0$ and $X_2=f_0+f_1$\@.  We use
coordinates $((x_1,x_2,x_3),(v_1,v_2,v_3))$ for $\tb{\real^3}$ and we compute
\begin{equation*}
\tlift{X_1}=x_2\pderiv{}{x_1}+v_2\pderiv{}{v_1},\qquad
\tlift{X_2}=x_2\pderiv{}{x_1}+x_3\pderiv{}{x_2}+
v_2\pderiv{}{v_2}+v_3\pderiv{}{x_2}.
\end{equation*}
If $t\mapsto Y(t)$ is a time-dependent vector field with values in $\sF(\real^3)$\@, then
\begin{equation*}
Y_t=f_0+\nu_1(t)f_1+\nu_2(t)f_2=x_2\pderiv{}{x_1}+\nu_1(t)x_3\pderiv{}{x_2}
+\nu_2(t)\pderiv{}{x_3},
\end{equation*}
for functions $\nu_1$ and $\nu_2$ whose character is determined by the
open-loop subfamily $\sO_{\fG,1}$\@.  The linearisation about the two
reference flows are described by the differential equations
\begin{equation*}
\begin{aligned}
\dot{x}_1(t)=&\;x_2(t),\\
\dot{x}_2(t)=&\;0;\\
\dot{x}_3(t)=&\;0,\\
\dot{v}_1(t)=&\;v_2(t)+x_2(t),\\
\dot{v}_2(t)=&\;\nu_1(t)x_3(t),\\
\dot{v}_3(t)=&\;\nu_2(t),
\end{aligned}\qquad
\begin{aligned}
\dot{x}_1(t)=&\;x_2(t),\\
\dot{x}_2(t)=&\;x_3(t);\\
\dot{x}_3(t)=&\;0,\\
\dot{v}_1(t)=&\;v_2(t)+x_2(t),\\
\dot{v}_2(t)=&\;v_3(t)+\nu_1(t)x_3(t),\\
\dot{v}_3(t)=&\;\nu_2(t),
\end{aligned}
\end{equation*}
respectively.  The linearisations about $(X_1,(0,0,0))$ and $(X_2,(0,0,0))$
will be time-independent since the vector fields $X_1$ and $X_2$ are
time-independent, and we easily determine that these linearisations are given
by
\begin{equation*}
\begin{bmatrix}\dot{v}_1(t)\\\dot{v}_2(t)\\\dot{v}_3(t)\end{bmatrix}=
\begin{bmatrix}0&1&0\\0&0&0\\0&0&0\end{bmatrix}
\begin{bmatrix}v_1(t)\\v_2(t)\\v_3(t)\end{bmatrix}+
\begin{bmatrix}0\\0\\\nu_2(t)\end{bmatrix}
\end{equation*}
and
\begin{equation*}
\begin{bmatrix}\dot{v}_1(t)\\\dot{v}_2(t)\\\dot{v}_3(t)\end{bmatrix}=
\begin{bmatrix}0&1&0\\0&0&1\\0&0&0\end{bmatrix}
\begin{bmatrix}v_1(t)\\v_2(t)\\v_3(t)\end{bmatrix}+
\begin{bmatrix}0\\0\\\nu_2(t)\end{bmatrix},
\end{equation*}
respectively.  These are exactly the two distinct linearisations we
encountered in Example~\ref{eg:bad-linearise}\@.  Thus we can see here what
was going on in Example~\ref{eg:bad-linearise}\@: we were linearising about
two different reference flows.  This also highlights the dangers of explicit
and fixed parameterisations by control:~one can unknowingly make choices that
affect conclusions.

We comment that the reason this example does not meet the conditions of
Proposition~\ref{prop:nice-cas} is that the vector fields $f_1$ and $f_2$ are
not linearly independent at $(0,0,0)$\@.  The distribution generated by these
vector fields has $(0,0,0)$ as a singular point.  These sorts of matters will
doubtless be interesting in subsequent studies of geometric control systems
in our framework.

Finally, using Lemma~\ref{lem:lincont-form}\@, we can easily conclude that
this system is linearly controllable.\oprocend
\end{example}

\section{Future work}\label{sec:future-work}

There is a lot of control theory that has yet to be done in our framework of
\gcs{}s.  In this closing section, we discuss a few avenues for future work,
and provide a few preliminary ideas related to these directions.

\subsection{Controllability}\label{subsec:future-cont}

The controllability of nonlinear systems comprises a vast and difficult
component of the geometric control theory literature.  A number of papers
have been published addressing the seemingly impenetrable nature of the
problems of
controllability~\cite{AAA:99,RMB/MK:02,MK:90b,MK:90a,MK:06,EDS:88}\@.
Despite this, there has been substantial effort dedicated to determining
sufficient or necessary conditions for
controllability~\cite{AAA/RVG:93a,AB/GS:83,RMB/GS:84,RMB/GS:86,RMB/GS:93,
GWH/HH:70,HH:74,HH:76a,HH:76b,HH:77,HH:82a,HH/MK:87,MK:87,MK:88c,MK:92,
MK:98,MK:99,GS:86,HJS:73,HJS:78,HJS:83a,HJS:87,HJS/VJ:72}\@.  The problem of
controllability has a certain lure that attracts researchers in geometric
control theory.  The problem is such a natural one that it feels as if it
should be possible to obtain complete results, at least in some quite general
situations.  However, this objective remains to be fulfilled.

Our view is that one of the reasons for this is that many of the approaches
to controllability are not feedback-invariant.  An extreme example of this
are methods for studying controllability of control-affine systems,
\emph{fixing} a drift vector field $f_0$ and control vector fields
$f_1,\dots,f_m$\@, and using these as generators of a free Lie algebra.  In
this sort of analysis, Lie series are truncated, leading to the notion of
``nilpotent approximation'' of control systems.  These ideas are reflected in
a great many of the papers cited above.  The difficulty with this approach is
that it will behave very badly under feedback
transformations,~\cf~Example~\ref{eg:bad-linearise}\@.  This is discussed by \citet{ADL:12a}\@.

One approach is then to attempt to find feedback-invariant conditions for
local controllability.  In the first-order case,~\ie~the more or less linear
case, this leads to Definition~\ref{def:lincont}\@; see
also~\cite{RMB/GS:84}\@.  Second-order feedback-invariant conditions are
considered in~\cite{JBG:98,RMH/ADL:02c}\@.  Any attempts to determine
higher-order feedback-invariant controllability conditions have, as far as we
know, met with no success.  Indeed, the likelihood of this approach leading
anywhere seems very small, given the extremely complicated manner in which
feedback transformations interact with controllability conditions.

Therefore, the most promising idea would appear to be to develop a framework
for control theory that has feedback-invariance ``built in.''  It is this
that we have done in this paper.  In his PhD thesis, \citet{COA:10} provides
a class of control variations that is well-suited to our feedback-invariant
approach.  \citet{COA/ADL:11} have used these control variations to
completely characterise controllability of a class of homogeneous systems.
It will be an interesting project to apply the variations of
\citeauthor{COA:10} in our framework to see what sorts of conditions
for controllability naturally arise.

\subsection{Optimal control and the Maximum Principle}

It should be a fairly straightforward exercise to formulate optimal control
problems in our framework.  Also, our approach to linearisation in
Section~\ref{sec:linearisation} already provides us with the natural means by
which needle variations can be transported along reference trajectories, and
so one expects that an elegant version of the Maximum Principle of
\citet{LSP/VGB/RVG/EFM:61} will be possible.  There will be a resemblance in
this to the work of \citet{HJS:97b}\@, who provides already a satisfying
formulation of the Maximum Principle on manifolds.  In the same way as the
natural feedback-invariance of our formulation should aid in the study of
controllability, it should also aid in the study of higher-order conditions
for optimality.  In geometric control theory, the study of so-called singular
extremals (those not characterised by the Maximum Principle) is problematic
for multi-input systems, so hopefully our approach can shed light on this.

As outlined in Section~\ref{subsec:subriemannian}\@, problems in
sub-Riemannian geometry fit naturally into the \gcs\ framework, and can
likely be handled well by a theory of optimal control for \gcs{}s.

\subsection{Feedback and stabilisation theory}

There are, one could argue, three big problems in control theory.  Two,
controllability and optimal control, are discussed above.  The third is
stabilisation.  This problem, being one of enormous practical importance, has
been comprehensively studied, mainly from the point of view of Lyapunov
theory, where the notion of a ``control-Lyapunov function'' provides a useful
device for characterising when a system is
stabilisable~\cite{FHC/YSL/EDS/AIS:97} and for stabilisation if one is
known~\cite{EDS:89b}\@.  Our view is that Lyapunov characterisations for
stabilisability are important from a practical point of view, but, from a
fundamental point of view, merely replace one impenetrable notion,
``stabilisability,'' with another, ``existence of a control-Lyapunov
function.''  This is expressed succinctly by \citeauthor{EDS:98}\@.
\begin{quote}\itshape
In any case, all converse Lyapunov results are purely existential, and are of
no use in guiding the search for a Lyapunov function.  The search for such
functions is more of an art than a science, and good physical insight into a
given system plus a good amount of trial and error is typically the only way
to proceed.\textemdash\citet[page~259]{EDS:98}
\end{quote}
As \citeauthor{EDS:98} goes on to explain, there are many heuristics for
guessing control-Lyapunov functions.  However, this is unsatisfying if one is
seeking a general understanding of the problem of stabilisability, and not
just a means of designing stabilising controllers for classes of systems.

It is also the case that there has been virtually no work on stabilisability
from a geometric perspective.  Topological characterisations of
stabilisability such as those of \citet{RWB:83} (refined by
\citet{RO/LP/IM:03,JZ:89}) and \citet{JMC:90} are gratifying when they are
applicable, but they are far too coarse to provide anything even close to a
complete characterisation of the problem.  Indeed, the extremely detailed and
intricate analysis of controllability, as reflected by the work we cite
above, is simply not present for stabilisability.  It is fair to say that,
outside the control-Lyapunov framework, very little work has been done in
terms of really understanding the structural obstructions to stabilisability.
Moreover, it is also fair to say that almost none of the published literature
on stabilisation and stabilisability passes the ``acid test'' for feedback
invariance that we discuss in Section~\ref{subsec:modelcomp}\@.  For
researchers such as ourselves interested in structure, this in an
unsatisfying state of affairs.

Our framework provides a natural means of addressing problems like this, just
as with controllability and optimal control, because of the
feedback-invariance of the framework.  Indeed, upon reflection, one sees that
the problem of stabilisability should have some relationships with that of
controllability, although little work has been done along these lines (but
see the PhD thesis of \citet{PI:12b}).  This area of research is wide open~\cite{ADL:12a}\@.

\subsection{Linear system theory}

Our definition of linearisation suggests an immediate generalisation from
tangent bundles to vector bundles.  Let us quickly see how it will work,
making no pretence to the level of generality of the main body of the paper.
\begin{definition}
Let $r\in\{\infty,\omega\}$ and let $\map{\pi}{\man{E}}{\man{M}}$ be vector
bundle of class $\C^r$\@.  A vector field $X\in\sections[r]{\tb{\man{E}}}$ is
\defn{linear} if
\begin{compactenum}[(i)]
\item $X$ is $\pi$-projectable,~\ie~there exists a vector field $\pi
X\in\sections[r]{\tb{\man{M}}}$ such that $\tf[e_x]{\pi}(X(e_x))=\pi X(x)$
for every $x\in\man{M}$ and $e_x\in\man{E}_x$\@, and
\item $X$ is a vector bundle mapping for which the diagram
\begin{equation*}
\xymatrix{{\man{E}}\ar[r]^{X}\ar[d]_{\pi}&{\tb{\man{E}}}\ar[d]^{\tf{\pi}}\\
{\man{M}}\ar[r]_{\pi X}&\tb{\man{M}}}
\end{equation*}
commutes.\oprocend
\end{compactenum}
\end{definition}

The prototypical linear vector field is the tangent lift $\tlift{X}$\@, which
is a linear vector field on the vector bundle
$\map{\tbproj{\man{M}}}{\tb{\man{M}}}{\man{M}}$ according to the preceding
definition.  One may show that flows of linear vector fields are such that
the diagram
\begin{equation*}
\xymatrix{{\man{E}}\ar[r]^{\flow{X}{t}}\ar[d]_{\pi}&
{\man{E}}\ar[d]^{\pi}\\{\man{M}}\ar[r]_{\flow{\pi X}{t}}&{\man{M}}}
\end{equation*}
commutes and $\flow{X}{t}|\man{E}_x$ is an isomorphism of $\man{E}_x$ with
$\man{E}_{\flow{\pi X}{t}(x)}$~\cite[Proposition~47.9]{IK/PWM/JS:93}\@.

Vertical lifts are also easily defined for vector bundles.  We first define
the vector bundle map $\map{\vlft}{\pi^*\man{E}}{\tb{\man{E}}}$ as follows.
Let $x\in\man{M}$ and let $e_x,f_x\in\man{E}_x$\@.  The \textbf{vertical
lift} of $f_x$ to $e_x$ is given by
\begin{equation*}
\vlft(e_x,f_x)=\tderivatzero{}{t}(e_x+tf_x).
\end{equation*}
Now, if $\xi\in\sections[\infty]{\man{E}}$\@, we define
$\vlift{\xi}\in\sections[\infty]{\tb{\man{E}}}$ by
$\vlift{\xi}(e_x)=\vlft(e_x,\xi(x))$\@.

One also has the notion of a projectable presheaf of vector fields on a
vector bundle.
\begin{definition}
Let $r\in\{\infty,\omega\}$\@, let $\map{\pi}{\man{E}}{\man{M}}$ be a vector
bundle of class $\C^r$\@, and let $\sG$ be a presheaf of sets of vector
fields of class $\C^r$ on $\man{E}$\@.  The presheaf $\sG$ is
\defn{projectable} if
\begin{equation*}\eqoprocend
\sG(\nbhd{W})=\setdef{Z|\nbhd{W}}{Z\in\sG(\pi^{-1}(\pi(\nbhd{W})))}.
\end{equation*}
\end{definition}

One also has the more or less obvious notion of presheaves of sets of
sections of $\man{E}$\@.
\begin{definition}
Let $r\in\{\infty,\omega\}$ and let $\map{\pi}{\man{E}}{\man{M}}$ be a vector
bundle of class $\C^r$\@.  A \defn{presheaf of sets of $\C^r$-sections} of
$\man{E}$ is an assignment, to each open set $\nbhd{U}\subset\man{M}$\@, a
subset $\sF(\nbhd{U})$ of $\sections[r]{\man{E}|\nbhd{U}}$ with the property
that, for open sets $\nbhd{U},\nbhd{V}\subset\man{M}$ with
$\nbhd{V}\subset\nbhd{U}$\@, the map
\begin{equation*}
\mapdef{r_{\nbhd{U},\nbhd{V}}}{\sF(\nbhd{U})}{\sections[r]{\tb{\nbhd{V}}}}
{\xi}{\xi|\nbhd{V}}
\end{equation*}
takes values in $\sF(\nbhd{V})$\@.  Elements of $\sF(\nbhd{U})$ are called
\defn{local sections} over $\nbhd{U}$\@.\oprocend
\end{definition}

One also has an analogue of Lemma~\ref{lem:linear-presheaf} for vector
bundles, which makes sense of the following, final, definition.
\begin{definition}
Let $r\in\{\infty,\omega\}$ and let $\map{\pi}{\man{E}}{\man{M}}$ be a vector
bundle of class $\C^r$\@.  A \defn{$\C^r$-linear system} on $\man{E}$ is a
$\C^r$-\gcs\ $\fG=(\man{E},\sF)$\@, where the projectable presheaf of sets of
vector fields $\sF$ is characterised uniquely by the requirement that, for
every open subset $\nbhd{U}\subset\man{M}$\@,
\begin{equation*}
\sF(\pi^{-1}(\nbhd{U}))=
\setdef{X+\vlift{Y}}{X\in\sF_0(\pi^{-1}(\nbhd{U}),\ Y\in\sF_1(\nbhd{U})},
\end{equation*}
where $\sF_0$ is a projectable presheaf of sets of linear vector fields
on $\man{E}$ and $\sF_1$ is a presheaf of sets of sections of $\man{E}$\@.\oprocend
\end{definition}

This is then a class of \gcs{}s containing linearisations of \gcs{}s as a
special case.  One is then interested in what one can say about problems of
control\textemdash{}controllability, optimal control theory,
stabilisation\textemdash{}for these systems.  An approach to this is
presented in~\cite{ADL/DRT:10} for control-affine systems.
In~\cite[Chapter~5]{FC/WK:00} one can find a setup along these lines, but
with a decidedly different perspective.

\subsection{The category of \gcs{}s}

In Section~\ref{subsec:gcs-category} we introduced morphisms between \gcs{}s
with the objective of showing that our framework is feedback-invariant.  The
notion of morphism we present is one that is natural and possibly easy to
work with.  It would be, therefore, interesting to do all of the exercises of
category theory with the category of \gcs{}s.  That is, one would like to
study epimorphisms, monomorphisms, subobjects, quotient objects, products,
coproducts, pull-backs, push-outs, and various functorial operations in this
category.  Many of these may not be interesting or useful, or even exist.
But probably some of it would be of interest.  For example, \citet{PT/GJP:05}
study quotients of control systems, and \citet{VIE:99} studies various
categorical constructions for control-affine systems.

\subsection{Real analytic chronological calculus}

As we have mentioned a few times, the treatment of real analytic time-varying
vector fields by \citet{AAA/RVG:78} is carried out under very restrictive
hypotheses, namely that the real analytic vector fields are required to admit
bounded holomorphic extensions to a fixed neighbourhood in the
complexification whose width is bounded uniformly from below.  Even in the
case of compact real analytic manifolds, this is a severe restriction.  With
the theory of real analytic time-varying vector fields presented in this
paper, a fully functioning real analytic chronological calculus ought to be
feasible.

Moreover, the results that we have proved above allow a strengthening of the
existing results of \citet{AAA/RVG:78}\@, even in the smooth case, in the
following way.  \citeauthor{AAA/RVG:78} do everything ``weakly.''  By this we
mean the following.  Vector fields are characterised by
\citeauthor{AAA/RVG:78} by what they do to functions,~\ie~they use what we
call the weak-$\sL$ topology.  In
Theorems~\ref{the:COinfty-weak}\@,~\ref{the:COmm'-weak}\@,
and~\ref{the:Comega-weak} we see that this is equivalent to working directly
with the appropriate topologies for vector fields.  Probably this is well
understood in the finitely differentiable and smooth cases, but in this paper
we have understood that this is also true in the real analytic case.  Also,
when dealing with matters such as measurability, integrability, and absolute
continuity, \citeauthor{AAA/RVG:78} reduce to the scalar case by first
composing all objects with the evaluation functionals $\ev_x$ as in the proof
of Theorem~\ref{the:Cinfty-tvsec} (and by implication, in the proofs of
Theorems~\ref{the:Cmm'-tvsec} and~\ref{the:Comega-tvsec}), and defining and
computing with the scalar versions of these notions.  However,
Theorems~\ref{the:smooth-td-summary}\@,~\ref{the:mm'-td-summary}\@,
and~\ref{the:ra-td-summary} ensure that this is equivalent to doing
computations in the spaces of finitely differentiable, smooth, or real
analytic vector fields.  Again, perhaps this is understood in the finitely
differentiable and smooth cases, but we have shown that this is also true in
the real analytic case.

Thus, combining the preceding two paragraphs, one should be able to develop
the chronological calculus of \citet{AAA/RVG:78} into a more powerful and
broadly applicable tool.

\end{document}